\documentclass[a4paper, reqno]{amsart}

\usepackage{comment}
\usepackage{amssymb} 
\usepackage[T1]{fontenc}
\usepackage[latin1]{inputenc}
\usepackage{amsmath}
\usepackage{amsfonts}
\usepackage{amsxtra}
\usepackage{ae}
\usepackage[all, cmtip]{xy}
\usepackage{enumerate}
\usepackage{url}
\usepackage{mathrsfs}
\usepackage{verbatim}
\usepackage{graphicx}
\usepackage[pageanchor]{hyperref}
\usepackage{cancel}
\usepackage{etoolbox}

\usepackage[refpage]{nomencl}  

\renewcommand\nomgroup[1]{%
	\bigskip
	\item[\bfseries
	 \ifstrequal{#1}{A}{General notation}{%
	 	\ifstrequal{#1}{C}{$q$-calculus}{%
	 	\ifstrequal{#1}{G}{Actions of algebras and groups}{%
		 \ifstrequal{#1}{O}{Other symbols}{%
		  \ifstrequal{#1}{S}{Subscripts, superscripts, and decorations}{}%
	}}}}%
  ]
  \medskip
	}

\usepackage[usenames,dvipsnames]{xcolor}

\newcommand{\lie}{\mathfrak}

\newcommand{\NN}{\mathbb{N}}
\newcommand{\ZZ}{\mathbb{Z}}

\newcommand{\RR}{\mathbb{R}}
\newcommand{\CC}{\mathbb{C}}

\newcommand{\Poly}{\mathcal{O}}
\newcommand{\dual}{\vee}
\newcommand{\precon}{{}^\mathrm{c}}
\newcommand*{\ket}{\rangle}
\newcommand*{\bra}{\langle}
\newcommand*{\ad}{\mathsf{ad}}

\newcommand*{\A}{\mathcal{A}}
\newcommand*{\B}{\mathcal{B}}

\newcommand*{\D}{\mathcal{D}}

\newcommand*{\M}{\mathcal{M}}
\newcommand*{\N}{\mathcal{N}}
\newcommand*{\E}{\mathcal{E}}
\newcommand*{\F}{\mathcal{F}}

\renewcommand*{\H}{\mathcal{H}}

\renewcommand*{\L}{\mathcal{L}}
\renewcommand*{\P}{\mathcal{P}}
\newcommand*{\R}{\mathcal{R}}
\renewcommand*{\S}{\mathcal{S}}
\renewcommand*{\O}{\mathcal{O}}
\newcommand*{\T}{\mathcal{T}}

\newcommand*{\X}{\mathcal{X}}

\newcommand*{\Z}{\mathcal{Z}}

\newcommand*{\CF}{\mathfrak{C}}
\newcommand*{\DF}{\mathfrak{D}}

\newcommand*{\mx}{\mathsf{f}}
\newcommand*{\red}{\mathsf{r}}

\newcommand*{\res}{\mathsf{res}}
\newcommand*{\hit}{\!\rightharpoonup\!}
\newcommand*{\hitby}{\!\leftharpoonup\!}

\newcommand*{\CH}{\mathbb{C}}

\newcommand*{\KH}{\mathbb{K}}
\newcommand*{\LH}{\mathbb{L}}

\newcommand*{\opp}{\mathsf{opp}}
\newcommand*{\cop}{\mathsf{cop}}

\newcommand*{\ext}{\mathsf{ext}}

\newcommand{\half}{{\frac12}}

\newcommand{\roots}{\mathbf{Q}}
\newcommand{\simpleroots}{\Sigma}
\newcommand{\weights}{\mathbf{P}}

\renewcommand{\Re}{\mathrm{Re}}
\renewcommand{\Im}{\mathrm{Im}}

\DeclareMathOperator{\Ph}{Ph}
\DeclareMathOperator{\HC}{HC}

\DeclareMathOperator{\End}{End}
\DeclareMathOperator{\Soc}{Soc}
\DeclareMathOperator{\tr}{tr}
\DeclareMathOperator{\ch}{ch}

\DeclareMathOperator{\Aut}{Aut}
\DeclareMathOperator{\Hom}{Hom}

\DeclareMathOperator{\ind}{ind}
\DeclareMathOperator{\im}{im}

\DeclareMathOperator{\id}{id}

\DeclareMathOperator{\ann}{ann}

\DeclareMathOperator{\sign}{sign}
\DeclareMathOperator{\Sh}{Sh}
\DeclareMathOperator{\Quot}{Quot}
\DeclareMathOperator{\Spm}{Spm}
\DeclareMathOperator{\Det}{Det}

\newenvironment{bnum}
{\begin{list}{}
    {\setlength{\labelwidth}{15pt}
     \setlength{\leftmargin}{\labelwidth}
    }
}
{\end{list}}

\numberwithin{equation}{section}

\theoremstyle{plain}
\newtheorem{theorem}{Theorem}[section]
\newtheorem{prop}[theorem]{Proposition}
\newtheorem{lemma}[theorem]{Lemma}
\newtheorem{cor}[theorem]{Corollary}

\theoremstyle{definition}
\newtheorem{definition}[theorem]{Definition}

\theoremstyle{plain}
\newtheorem*{theorem*}{Theorem}
\newtheorem*{prop*}{Proposition}
\newtheorem*{lemma*}{Lemma}
\newtheorem*{cor*}{Corollary}

\theoremstyle{definition}
\newtheorem*{definition*}{Definition}
\newtheorem*{example*}{Example}

\makenomenclature

\begin{document}

\title{Complex semisimple quantum groups and representation theory}

\author{Christian Voigt}
\address{School of Mathematics and Statistics \\
         University of Glasgow \\
         15 University Gardens \\
         Glasgow G12 8QW \\
         UK
}
\email{christian.voigt@glasgow.ac.uk}

\author{Robert Yuncken}
\address{Laboratoire de Math\'ematiques \\ 
         Universit\'e Blaise Pascal \\
         Complexe universitaire des C\'ezeaux \\
         63177 Aubi\`ere Cedex \\
         France 
}

\email{Robert.Yuncken@math.univ-bpclermont.fr}

\thanks{The first author would like to thank the Isaac Newton Institute for Mathematical Sciences, Cambridge, for support
and hospitality during the programme Operator Algebras: Subfactors and their applications,
where work on this paper was undertaken. This work was supported by EPSRC grant no EP/K032208/1.
The first author was supported by the Polish National Science Centre grant no. 2012/06/M/ST1/00169.
This paper was partially supported by the grant H2020-MSCA-RISE-2015-691246-QUANTUM DYNAMICS.
The second author was supported by the project SINGSTAR of the Agence Nationale de la Recherche, ANR-14-CE25-0012-01 and by the CNRS PICS project OpPsi.  Both authors would also like to thank the Erwin Schr\"odinger Institute in Vienna, where some of the present research was undertaken.
}

\subjclass[2010]
{
20G42, 
17B37, 
16T05, 
46L67, 
81R50, 
46L65. 
}

\maketitle

\newpage

\section*{Introduction}

The theory of quantum groups, originating from the study of integrable systems, has seen a rapid development from the mid 1980's with far-reaching connections to various branches of mathematics, including knot theory, representation theory and operator algebras, see \cite{Drinfeldicm}, \cite{Lusztigbook}, \cite{CPbook}. While the term \emph{quantum group} itself has no precise definition, it is used to denote a number of related constructions, including in particular quantized universal enveloping algebras of semisimple Lie algebras and, dually, deformations of the algebras of polynomial functions on the corresponding semisimple groups. In operator algebras, the theory of locally compact quantum groups \cite{KVLCQG} is a powerful framework which allows one to extend Pontrjagin duality to a fully noncommutative setting.

In both the algebraic and the analytic theory of quantum groups, an important role is played by the Drinfeld double, also known as the quantum double, which is designed to produce solutions to the quantum Yang-Baxter equation.  The algebraic version of this construction, due to Drinfeld, appears already in \cite{Drinfeldicm}.  The analytical analogue of the Drinfeld double was defined and studied first by Podle\'s and Woronowicz in \cite{PWlorentz}.

If one applies the Drinfeld double construction to the Hopf algebra of functions on the $ q $-deformation of a compact semisimple Lie group, then, in accordance with the quantum duality principle \cite{Drinfeldicm}, \cite{Semenov-Tyan-Shanskii}, the resulting quantum group can be viewed as a quantization of the corresponding complex semisimple group. This fact, which was exhibited in \cite{PWlorentz} in the case of the quantum Lorentz group $SL_q(2,\mathbb{C})$, is key to understanding the structure of these Drinfeld doubles. More precisely, one can transport techniques from the representation theory of classical complex semisimple groups to the quantum situation, and it turns out that the main structural results carry over, albeit with sometimes quite different proofs.

These notes contain an introduction to the theory of complex semisimple quantum groups, that is, Drinfeld doubles of compact quantum groups arising from $ q $-deformations. Our main aim is to present the classification of irreducible Harish-Chandra modules for these quantum groups, or equivalently the irreducible Yetter-Drinfeld modules of $ q $-deformations of compact semisimple Lie groups. We also treat some operator algebraic aspects of these constructions, but we put our main emphasis on the algebraic considerations based on quantized universal enveloping algebras.

The main reason for going into a considerable amount of detail on the algebraic side is that the existing literature does not quite contain the results in the form needed for this purpose. Many authors work over the field $ \mathbb{Q}(q) $ of rational functions in $ q $, while we are mainly interested in the case that $ q \in \mathbb{C}^\times $ is not a root of unity. Although it is folklore that this does not affect the general theory in a serious way, there are subtle differences which are easily overlooked, in particular when it comes to Verma modules and Harish-Chandra bimodules. In addition, different conventions are used in the literature, which can make it cumbersome to combine results from different sources.

In the first part of these notes, consisting of Chapters \ref{chhopf} and \ref{chuqg}, we work over a general ground field $ \mathbb{K} $ and a deformation parameter $ q \in \mathbb{K}^\times $ which is not a root of unity. Technically, this means that one has to start from an element $ s \in \mathbb{K} $ such that $ s^L = q $ for a certain number $ L \in \mathbb{N} $ depending on the type of the underlying semisimple Lie algebra. No assumptions on the characteristic of $ \mathbb{K} $ are made, and in particular we do not rely on specialization at $ q = 1 $ in order to transport results from the classical situation to the quantum case. However, we shall freely use general constructions and facts from classical Lie theory as can be found, for instance, in \cite{HumphreysLie}.

In the second part of the notes, consisting of Chapters \ref{chcomplexqg}, \ref{chcato} and \ref{chreptheory}, we restrict ourselves to the case $ \mathbb{K} = \mathbb{C} $ and assume that $ q = e^h \neq 1 $ is a positive real number. Some parts of the material could be developed in greater generality, but for certain arguments the specific properties of the ground field $ \mathbb{C} $ and its exponential map do indeed play a role. This is the case, for instance, for the characterization of dominant and antidominant weights in Chapter \ref{chcato}.

We have used a number of textbooks as the basis of our presentation, let us mention in particular the books by Lusztig \cite{Lusztigbook}, Jantzen \cite{Jantzenqg}, Klimyk and Schm\"udgen \cite{KS}, Chari and Pressley \cite{CPbook}, Brown and Goodearl \cite{BrownGoodearllectures}, and the book by Humphreys on category $ \O $ in the classical setting \cite{HumphreysO}. The key source for the second half of these notes is the book by Joseph \cite{Josephbook}, which in turn builds to a large extent on work of Joseph and Letzter. Our account is largely self-contained, the only notable exception being the section on separation of variables in Chapter \ref{chuqg}, which relies on the theory of canonical bases. Nonetheless, we assume that the reader has some previous acquaintance with quantum groups, and we have kept short those arguments which can be easily found elsewhere. Let us also point out that our bibliography is far from complete, and the reader should consult the above mentioned sources for a historically accurate attribution of the results covered here, along with much more background and motivation.

\medskip

Let us now explain in detail how this text is organized.

In Chapter \ref{chhopf} we review the theory of Hopf algebras and multiplier Hopf algebras. While the basics of Hopf algebras can be found in numerous textbooks, the non-unital version, developed by Van Daele \cite{vDmult}, \cite{vDadvances} under the name of multiplier Hopf algebras, is not as widely known.  Multiplier Hopf algebras provide a natural setting for the study of some aspects of quantized enveloping algebras, like the universal $ R $-matrix. Moreover, they are a very useful tool for studying the link between the algebraic and analytic theory of quantum groups \cite{KvD}. Since only a limited amount of the general theory of multiplier Hopf algebras is needed for our purposes, we refer to the original sources for most of the proofs.

Chapter \ref{chuqg} contains an exposition of the basic theory of quantized universal enveloping algebras. Most of the material is standard, but we have made an effort to establish uniform conventions and notation. We work throughout with what is often called the \emph{simply connected} version of $ U_q(\mathfrak{g}) $. This is crucial for some of the more advanced parts of the theory, notably in connection with the $l$-functionals.
Our discussion of the braid group action on $ U_q(\mathfrak{g}) $ and its modules is significantly more detailed than what can be found in the standard textbooks. We have also tried to simplify and streamline various arguments in the literature, in particular we avoid some delicate filtration arguments in \cite{Josephbook} in the proof of Noetherianity of the locally finite part of $ U_q(\mathfrak{g}) $.

In Chapter \ref{chcomplexqg} we introduce our main object of study, namely complex semisimple quantum groups. These quantum groups can be viewed as quantizations of complex semisimple Lie groups considered as \emph{real} groups. As indicated above, they are obtained by applying the quantum double construction to compact semisimple quantum groups. We discuss the structure of complex quantum groups as locally compact quantum groups, including a description of their Haar weights and dual Haar weights.

The remaining two chapters are devoted to representation theory. Chapter \ref{chcato} contains a discussion of category $ \O $ for quantized universal enveloping algebras. This is parallel to the theory for classical universal enveloping algebras, but some peculiar new features arise in the quantum situation due to the periodicity in the space of purely imaginary weights.

Finally, Chapter \ref{chreptheory}, covers the representation theory of complex semisimple quantum groups and their associated Harish-Chandra modules. We also provide a detailed account of intertwining operators, by taking advantage of the relation between the category of Harish-Chandra modules and category $ \O $.
This both simplifies and extends the work by Pusz and Woronowicz \cite{Puszunitaryrepresentations}, \cite{PWunitaryrepresentations} on the quantum Lorentz group $ SL_q(2, \mathbb{C}) $.

Let us conclude with some remarks on notation and conventions. By default, an algebra is a unital associative algebra over a commutative ground ring, which will typically be a field $ \mathbb{K} $. In some situations we have to deal with non-unital algebras and their multiplier algebras, but it should always be clear from the context when a non-unital algebra appears. Unadorned tensor products are over the ground field $ \mathbb{K} $. In some situations we also use $ \otimes $ to denote the tensor product of Hilbert spaces, the spatial tensor product of von Neumann algebras, or the minimal tensor product of $ C^* $-algebras. 
Again, the meaning should be clear from the context. 
We will use the standard leg numbering notation for operators on multiple tensor products. 
A full list of notation is included at the end of these notes.

We have aimed to make the transition between the algebraic and analytic points of view as convenient as possible. This results in some slightly unconventional choices for the algebraically-minded reader, in particular in terms of the pairing between quantized enveloping algebras and their dual function algebras, which we define to be skew-pairings by default.

Last but not least, it is a pleasure to thank a number of people with whom we have discussed aspects of quantum groups and representation theory over the past few years.
Let us mention in particular Y.\ Arano, P.\ Baumann, K.\ A.\ Brown, I.\ Heckenberger, N.\ Higson, U.\ Kr\"ahmer, S.\ Neshveyev, S.\ Riche and D.\ Vogan; we are grateful to all of them for sharing their insight with us.

\newpage

\tableofcontents

\newpage

\section{Multiplier Hopf algebras} \label{chhopf}

In this chapter we collect definitions and basic results regarding Hopf algebras and multiplier Hopf algebras. Algebras are 
not assumed to have identities in general. Throughout we shall work over an arbitrary base field $ \mathbb{K} $, and 
all tensor products are over $ \mathbb{K} $.

\subsection{Hopf algebras} 

The language of Hopf algebras is the starting point for the study of noncommutative symmetries. In this section we collect some basic definitions and 
facts from the theory. There is a large variety of textbooks devoted to Hopf algebras and quantum groups where more information can be found, let us mention in 
particular \cite{CPbook}, \cite{Jantzenqg}, \cite{Kasselbook}, \cite{KS}. 

An algebra is a vector space $ H $ together with a linear map $ m: H \otimes H \rightarrow H $ such that the diagram 
\begin{align*}
\xymatrix{
H \otimes H \otimes H \ar@{->}[r]^{\;\;\; m \otimes \id} \ar@{->}[d]^{\id \otimes m} & H \otimes H \ar@{->}[d]^{m} \\
H \otimes H \ar@{->}[r]^{m} & H
     }
\end{align*}
is commutative. An algebra $ H $ is called unital if there exists a linear map $ u: \mathbb{K} \rightarrow H $ such that the diagram 
\begin{align*}
\xymatrix{
\mathbb{K} \otimes H \ar@{-}[dr]^{\cong} 
\ar[r]^{u \otimes \id} & H \otimes H \ar[d]^{m} & H \otimes \mathbb{K} \ar[l]_{\id \otimes u} \\
& H \ar@{-}[ur]^{\cong} &
}
\end{align*}
is commutative. Here the isomorphisms are induced by scalar multiplication.  As usual, we write $fg$ for $m(f\otimes g)$, and $1_H$ or simply $1$ for $u(1)$ in the unital case.

A (unital) algebra homomorphism between (unital) algebras $ A $ and $ B $ 
is a linear map $ \varphi: A \rightarrow B $ such that $ \varphi m_A = m_B(\varphi \otimes \varphi) $ (and $ \varphi u_A = u_B $ in the unital case). 
Here $ m_A, m_B $ denote the multiplication maps of $ A $ and $ B $, respectively, and $ u_A, u_B $ the unit maps in the unital case. 
In other words, $ \varphi $ is an algebra homomorphism iff $ \varphi(fg) = \varphi(f) \varphi(g) $ for all $ f,g \in A $,
 and additionally $ \varphi(1_A) = 1_B $ in the unital case.

By definition, a coalgebra is a vector space $ H $ together with a linear map $ \Delta: H \rightarrow H \otimes H $%
\label{nom:coalgebra_coproduct}
such that the diagram
\begin{align*}
\xymatrix{
H \ar@{->}[r]^{\Delta} \ar@{->}[d]^{\Delta} & H \otimes H \ar@{->}[d]^{\id \otimes \Delta} \\
H \otimes H \ar@{->}[r]^{\!\!\!\! \Delta \otimes \id} & H \otimes H \otimes H
     }
\end{align*}
is commutative. 
A coalgebra $ H $ is called counital if there exists a linear map $ \epsilon: H \rightarrow \mathbb{K} $
\nomenclature{$\epsilon$}{counit}%
such that the diagram 
\begin{align*}
\xymatrix{
\mathbb{K} \otimes H \ar@{-}[dr]^{\cong} 
\ar@{<-}[r]^{\epsilon \otimes \id} & H \otimes H \ar@{<-}[d]^{\Delta} & H \otimes \mathbb{K} \ar@{<-}[l]_{\id \otimes \epsilon} \\
& H \ar@{-}[ur]^{\cong} &}
\end{align*}
is commutative. 

We shall use the Sweedler notation 
$$
\Delta(c) = c_{(1)} \otimes c_{(2)} 
$$
when performing computations in coalgebras, where we are suppressing a summation on the right-hand side. For instance, the counit axiom reads $ \epsilon(c_{(1)}) c_{(2)} = c = c_{(1)} \epsilon(c_{(2)}) $ in this notation. 

In analogy to the case of algebras, a (counital) coalgebra homomorphism between (counital) coalgebras $ C $ and $ D $ 
is a linear map $ \varphi: C \rightarrow D $ such that $ (\varphi \otimes \varphi) \Delta_C = \Delta_D \varphi $ (and $ \epsilon_D \varphi = \epsilon_C $ in the counital case). 

If $ C $ is a (counital) coalgebra then the linear dual space $ C^* = \Hom(C, \mathbb{K}) $ is a (unital) algebra with multiplication 
$$
(fg, c) = (f, c_{(1)}) (g, c_{(2)})
$$
(and unit element $ \epsilon_C $). 
Here $ (\;, \;) $ denotes the canonical pairing between $ C $ and $ C^* $. 

Conversely, if $ A $ is a finite dimensional (unital) algebra, then the linear dual space $ A^* $ becomes a (counital) coalgebra using the transpose 
$ m^*: A^* \rightarrow (A \otimes A)^* \cong A^* \otimes A^* $ of the 
multiplication map $ m: A \otimes A \rightarrow A $ (and counit $ \epsilon(f) = f(1) $). 
We point out that if $ A $ is infinite dimensional this will typically break down since then 
$ (A \otimes A)^* \neq A^* \otimes A^* $.

However, at least in the finite dimensional situation, we have a complete duality between algebras and coalgebras. This can be used to 
transport concepts from algebras to coalgebras and vice versa. For instance, a coalgebra $ C $ is called cosimple if it does not 
admit any proper subcoalgebras. Here a subcoalgebra of $ C $ is a linear subspace $ D \subset C $ such that $ \Delta(D) \subset D \otimes D $. 
This notion corresponds to simplicity for the dual algebra $ A $. A coalgebra is called cosemisimple if it is a direct sum of its simple 
subcoalgebras.

A basic example of a cosimple coalgebra is the dual coalgebra $ C = M_n(\mathbb{K})^* $ of the algebra of $ n \times n $-matrices over $ \mathbb{K} $. 
Explicitly, the coproduct of $ C $ is given by 
$$
\Delta(u_{ij}) = \sum_{k = 1}^n u_{ik} \otimes u_{kj} 
$$
for $ 1 \leq i,j \leq n $, where the elements $ u_{ij} \in C $ are dual to the basis of standard matrix units for $ M_n(\mathbb{K}) $. 
All cosemisimple coalgebras we will encounter later on are direct sums of such cosimple matrix coalgebras.

A Hopf algebra structure is the combination of the structures of a unital algebra and a counital coalgebra as follows. 

\begin{definition} 
\label{def:Hopf_algebra}
A bialgebra is a unital algebra $ H $ which is at the same time a counital coalgebra such that 
the comultiplication $ \Delta: H \rightarrow H \otimes H $ and the counit $ \epsilon: H \rightarrow \mathbb{K} $ are 
algebra homomorphisms. \\
A bialgebra is a Hopf algebra if there exists a linear map $ S: H \rightarrow H $,
\label{nom:S}
called the antipode, such that the diagrams 
\begin{align*}
\xymatrix{
H \ar@{->}[r]^{u \epsilon} \ar@{->}[d]^{\Delta} & H \ar@{<-}[d]^{m} \\
H \otimes H \ar@{->}[r]^{\id \otimes S} & H \otimes H
     }
\qquad 
\xymatrix{
H \ar@{->}[r]^{u \epsilon} \ar@{->}[d]^{\Delta} & H \ar@{<-}[d]^{m} \\
H \otimes H \ar@{->}[r]^{S \otimes \id} & H \otimes H
     }
\end{align*}
are commutative. 
\end{definition} 
In the definition of a bialgebra, one can equivalently require that the multiplication map $ m $ and the unit map $ u $ are coalgebra homomorphisms. 

For many examples of Hopf algebras the antipode $ S $ is an invertible linear map. We write $ S^{-1} $ for the inverse of $ S $ in this situation. 
We remark that $ S $ is always invertible if $ H $ is a finite dimensional Hopf algebra. 

If $ H $ is a finite dimensional Hopf algebra then the dual space $ \hat{H} = H^* = \Hom(H, \mathbb{K}) $ can be naturally equipped with a 
Hopf algebra structure with comultiplication $ \hat{\Delta}: \hat{H} \rightarrow \hat{H} \otimes \hat{H} $ and counit $ \hat{\epsilon} $ such that 
\begin{align*}
(x y, f) &= (x \otimes y, \Delta(f)) \\
(\hat{\Delta}(x), f \otimes g) &= (x, gf) \\
\hat{\epsilon}(x) &= (x, 1) \\ 
(\hat{S}(x), f) &= (x, S^{-1}(f)) 
\end{align*}
for all $ x, y \in \hat{H} $ and $ f,g \in H $.  Here $ (\;, \;) $ denotes the canonical pairing between $ H $ and $ \hat{H} $, 
or $ H \otimes H $ and $ \hat{H} \otimes \hat{H} $ respectively, given by evaluation.
Note that we are choosing the convention whereby $\hat{H}$ and $H$ are skew-paired, in the sense that the comultiplication $\hat\Delta$ is dual to the opposite multiplication $ f\otimes g \mapsto gf$.

In order to generalize this construction to certain infinite dimensional examples we shall work with the theory of multiplier Hopf algebras 
discussed in the next section.

\subsection{Multiplier Hopf algebras} \label{secmha}

The theory of multiplier Hopf algebras developed by van Daele and his coauthors \cite{vDadvances}, \cite{DvD}, \cite{DvDZ} is 
an extension of the theory of Hopf algebras to the case where the underlying algebras do not necessarily have identity elements. 

\subsubsection{Essential algebras} 

In order to obtain a reasonable theory, it is necessary to impose some conditions on the multiplication of an algebra. 
We will work with algebras that are essential in the following sense. 
\begin{definition} 
An algebra $ H $ is called essential if $ H \neq 0 $ and the multiplication map induces an isomorphism $ H \otimes_H H \cong H $.
\end{definition}
Clearly, every unital algebra is essential. 

More generally, assume that $ H $ has local units in the sense that for every finite set of elements $ h_1, \dots, h_n \in H $ there 
exists $ u, v \in H $ such that $ u h_j = h_j $ and $ h_j v = h_j $ for all $ j $. Then $ H $ is essential. 
Regular multiplier Hopf algebras, to be defined below, automatically have local units, see \cite{DvDZ}. 

The prototypical example of a non-unital essential algebra to keep in mind is an algebra of the form 
$$
H = \bigoplus_{i \in I} A_i, 
$$
where $ (A_i)_{i \in I} $ is a family of unital algebras and multiplication is componentwise. In fact, all non-unital essential algebras 
that we will encounter later on are of this form. 

Let $ H $ be an algebra. A left $ H $-module $ V $ is called essential if the canonical map $ H \otimes_H V \rightarrow V $ is 
an isomorphism. An analogous definition can be given for right modules. In particular, an essential algebra $ H $ is an essential
left and right module over itself.

\subsubsection{Algebraic multiplier algebras} 

To proceed further we need to discuss multipliers. A left multiplier for an algebra $ H $ is a linear map 
$ L: H \rightarrow H $ such that $ L(fg) = L(f)g $ for all $ f,g \in H $. Similarly, a right multiplier 
is a linear map $ R: H \rightarrow H $ such that $ R(fg) = fR(g) $ for all $ f,g \in H $. 
We let $ \M_l(H) $ and $ \M_r(H) $ be the spaces of left and right multipliers, respectively. 
These spaces become algebras with multiplication given by composition of maps. 

\begin{definition}
The algebraic multiplier algebra $ \M(H) $
\nomenclature[o$M$]{$\M(H)$}{multiplier algebra of an algebra $H$}%
of an algebra $ H $ is the space of all 
pairs $ (L,R) $ where $ L $ is a left multiplier and $ R $ is a right multiplier 
for $ H $ such that $ fL(g) = R(f) g $ for all $ f,g \in H $. 
\end{definition} 

The algebra structure of $ \M(H) $ is inherited from $ \M_l(H) \oplus \M_r(H) $. 
There is a natural homomorphism $ \iota: H \rightarrow \M(H) $. By construction, $ H $ is a left and right $ \M(H) $-module in a natural way. 

In the case that $ H = \bigoplus_{i \in I} A_i
$
is a direct sum of unital algebras $ A_i $ for $ i \in I $, it is straightforward to check that 
$$
\M(H) \cong \prod_{i \in I} A_i
$$
is the direct product of the algebras $ A_i $. 

Let $ H $ and $ K $ be algebras and let $ \varphi: H \rightarrow \M(K) $ be a 
homomorphism. Then $ K $ is a left and right $ H $-module in an obvious way. 
We say that the homomorphism $ \varphi: H \rightarrow \M(K) $ is essential if it turns 
$ K $ into an essential left and right $ H $-module, that is, we require 
$ H \otimes_H K \cong K \cong K \otimes_H H $. Note that the identity map $ \id: H \rightarrow H $ 
defines an essential homomorphism $ H \rightarrow \M(H) $ iff the algebra $ H $ is essential. 

\begin{lemma}
Let $ H $ be an algebra and let $ \varphi: H \rightarrow \M(K) $ be an essential homomorphism into the multiplier algebra of an essential 
algebra $ K $. Then there exists a unique unital homomorphism $ \Phi: \M(H) \rightarrow \M(K) $ such that $ \Phi \iota = \varphi $ where 
$ \iota: H \rightarrow \M(H) $ is the canonical map.  
\end{lemma}

\proof We obtain a linear map $ \Phi_l: \M_l(H) \rightarrow \M_l(K) $ by 
\begin{equation*} 
\M_l(H) \otimes K \cong
\xymatrix{
\M_l(H) \otimes H \otimes_H K \; \ar@{->}[r]^{\;\;\;\;\;\;\;\; m \otimes \id}
& H \otimes_H K }
\cong K 
\end{equation*}
and accordingly a linear map $ \Phi_r: \M_r(H) \rightarrow \M_r(K) $ by 
\begin{equation*} 
K \otimes \M_r(H) \cong
\xymatrix{
K \otimes_H H \otimes \M_r(H) \; \ar@{->}[r]^{\;\;\;\;\;\;\;\; \id \otimes m}
& K \otimes_H H }
\cong K.
\end{equation*}
It is straightforward to check that $ \Phi((L,R)) = (\Phi_l(L),\Phi_r(R)) $ defines a unital homomorphism 
$ \Phi: \M(H) \rightarrow \M(K) $ such that $ \Phi \iota = \varphi $. Uniqueness of $ \Phi $ follows from the 
fact that $ \varphi(H) \cdot K = K $ and $ K \cdot \varphi(H) = K $. \qed 

We note that essential homomorphisms behave well under tensor products. 
More precisely, assume that $ H_1, H_2 $ are essential algebras and let $ \varphi_1: H_1 \rightarrow \M(K_1) $ and $ \varphi_2: H_2 \rightarrow \M(K_2) $ be essential 
homomorphisms into the multiplier algebras of algebras $ K_1 $ and $ K_2 $. Then the induced homomorphism 
$ \varphi_1 \otimes \varphi_2 : H_1 \otimes H_2 \rightarrow \M(K_1 \otimes K_2) $ is essential. 

Following the terminology of van Daele, we say that an algebra $ H $ is nondegenerate if the multiplication map 
$ H \times H \rightarrow H $ defines a nondegenerate bilinear pairing. That is, $ H $ is nondegenerate iff 
$ fg = 0 $ for all $ g \in H $ implies $ f = 0 $ and $ fg = 0 $ for all $ f $ implies $ g = 0 $. 
These conditions can be reformulated by saying that the natural maps 
$$
H \rightarrow \M_l(H), \quad H \rightarrow \M_r(H)
$$
are injective. In particular, for a nondegenerate algebra the canonical map $ H \rightarrow \M(H) $ is injective. 

Nondegeneracy of an algebra is a consequence of the existence of a faithful linear functional in the following sense. 
\begin{definition} \label{defflf}
Let $ H $ be an algebra. A linear functional $ \omega: H \rightarrow \mathbb{K} $ is called faithful if
$ \omega(fg) = 0 $ for all $ g $ implies $ f = 0 $ and $ \omega(fg) = 0 $ for all $ f $ implies $ g = 0 $. 
\end{definition}

\subsubsection{Multiplier Hopf algebras} 

In this subsection we introduce the notion of a multiplier Hopf algebra. 

Let $ H $ be an essential algebra and let $ \Delta: H \rightarrow \M(H \otimes H) $ 
be a homomorphism. The left Galois maps $ \gamma_l, \gamma_r: H \otimes H \rightarrow \M(H \otimes H) $ for $ \Delta $ are defined by 
$$
\gamma_l(f\otimes g) = \Delta(f)(g \otimes 1),\qquad \gamma_r(f \otimes g) = \Delta(f) (1 \otimes g). 
$$
Similarly, the right Galois maps $ \rho_l, \rho_r: H \otimes H \rightarrow \M(H \otimes H) $ for $ \Delta $ are defined by 
$$
\rho_l(f \otimes g) = (f \otimes 1) \Delta(g), \qquad \rho_r(f \otimes g) = (1 \otimes f) \Delta(g). 
$$
These maps, or rather their appropriate analogues on the Hilbert space level, play an important role in the analytical study of 
quantum groups \cite{BSUM}, \cite{KVLCQG}. 
Our terminology originates from the theory of Hopf-Galois extensions, see for instance \cite{Montgomery}. 

We say that an essential homomorphism $ \Delta: H \rightarrow \M(H \otimes H) $%
\label{nom:multiplier_coproduct}
is a comultiplication 
if it is coassociative, that is, if 
$$
(\Delta \otimes \id) \Delta = (\id \otimes \Delta)\Delta 
$$
holds. Here both sides are viewed as maps from $ H $ to $ \M(H \otimes H \otimes H) $. 
An essential algebra homomorphism $ \varphi: H \rightarrow \M(K) $ between algebras with comultiplications is called a 
coalgebra homomorphism if $ \Delta \varphi = (\varphi \otimes \varphi) \Delta $. 

We need some more terminology. The opposite algebra $ H^{\opp} $
\nomenclature[s$H$1]{$H^\opp$}{opposite algebra}%
of $ H $ is the space $ H $ equipped with 
the opposite multiplication. That is, the multiplication $ m^{\opp} $ in $ H^{\opp} $ is defined by $ m^{\opp} = m \sigma $ where 
$ m: H \otimes H \rightarrow H $ is the multiplication in $ H $ and $ \sigma: H \otimes H \rightarrow H \otimes H $ 
is the flip map given by $ \sigma(f \otimes g) = g \otimes f $.
\label{nom:flip1}%
An algebra antihomomorphism between $ H $ and $ K $ 
is an algebra homomorphism from $ H $ to $ K^{\opp} $. Equivalently, an algebra antihomomorphism can be viewed as 
an algebra homomorphism $ H^{\opp} \rightarrow K $. 
If $ \Delta: H \rightarrow \M(H \otimes H) $ is a comultiplication then $ \Delta $ also defines a comultiplication 
$ H^{\opp} \rightarrow \M(H^{\opp} \otimes H^{\opp}) $.

Apart from changing the order of multiplication we may also reverse the order of a comultiplication. 
If $ \Delta: H \rightarrow \M(H \otimes H) $ is a comultiplication then the opposite comultiplication $ \Delta^{\cop} $
\nomenclature[s$H$2]{$H^\cop$}{opposite coalgebra}%
is the essential homomorphism 
from $ H $ to $ \M(H \otimes H) $ defined by $ \Delta^{\cop} = \sigma \Delta $. Here $ \sigma: \M(H \otimes H) \rightarrow \M(H \otimes H) $ 
is the extension of the flip map to multipliers. 
We write $ H^{\cop} $ for $ H $ equipped with the opposite comultiplication. 
Using opposite comultiplications we obtain the notion of a coalgebra antihomomorphism. 
That is, a coalgebra antihomomorphism between $ H $ and $ K $ 
is a coalgebra homomorphism from $ H $ to $ K^{\cop} $, or equivalently, from $ H^{\cop} $ to $ K $.

Let us now give the definition of a multiplier Hopf algebra \cite{vDmult}. 

\begin{definition} \label{defmhopf}
A multiplier Hopf algebra is an essential algebra $ H $ together with a comultiplication $ \Delta: H \rightarrow \M(H \otimes H) $ 
such that the Galois maps $ \gamma_r, \rho_l $ are isomorphisms
from $H\otimes H$ to $H\otimes H \subset \M(H\otimes H)$.

A regular multiplier Hopf algebra is an essential algebra $ H $ together with a comultiplication $ \Delta: H \rightarrow \M(H \otimes H) $ 
such that all Galois maps $ \gamma_l, \gamma_r, \rho_l, \rho_r $ are isomorphisms from $H\otimes H$ to $H\otimes H$. 
\end{definition}
A morphism between multiplier Hopf algebras $ H $ and $ K $ is an essential algebra homomorphism 
$ \alpha: H \rightarrow \M(K) $ such that $ (\alpha \otimes \alpha)\Delta = \Delta \alpha $. 

Note that a multiplier Hopf algebra $ H $ is regular iff $ H^{\opp} $ is a multiplier Hopf algebra, or equivalently, iff $ H^{\cop} $ 
is a multiplier Hopf algebra. 

We have the following fundamental result due to van Daele \cite{vDmult}. 
\begin{theorem} \label{thmcharmha} 
Let $ H $ be a multiplier Hopf algebra. Then there exists an essential algebra homomorphism  $ \epsilon: H \rightarrow \mathbb{K} $ and 
an algebra antihomomorphism $ S: H \rightarrow \M(H) $ such that 
\begin{equation*}
(\epsilon \otimes \id)\Delta = \id = (\id \otimes \epsilon)\Delta
\end{equation*}
and
\begin{equation*}
 m(S \otimes \id) \gamma_r = \epsilon \otimes \id, \qquad m(\id \otimes S) \rho_l = \id \otimes \epsilon.
\end{equation*}
If $ H $ is a regular multiplier Hopf algebra, then $ S $ is a linear isomorphism from $ H $ to $ H $. 
\end{theorem} 
If follows from Theorem \ref{thmcharmha} that Definition \ref{defmhopf} reduces to the definition of Hopf algebras in the unital 
case. A unital regular multiplier Hopf algebra is the same thing as a Hopf algebra with invertible antipode. 

We will be exclusively interested in regular multiplier Hopf algebras. 
For a regular multiplier Hopf algebra, the antipode $ S $ is a coalgebra antihomomorphism, that is, it satisfies $ (S \otimes S)\Delta = \Delta^{\cop} S $, 
see \cite{vDmult}.

Given a multiplier Hopf algebra we will again use the Sweedler notation 
$$
\Delta(f) = f_{(1)} \otimes f_{(2)} 
$$
for the comultiplication of $ H $. This is useful for writing calculations formally, although in principle one 
always has to reduce everything to manipulations with Galois maps.

\subsection{Integrals} \label{secintegrals}

In this section we discuss integrals for multiplier Hopf algebras. For a detailed treatment we refer to \cite{vDadvances}.  The 
case of ordinary Hopf algebras can be found in \cite{Sweedler}.  

\subsubsection{The definition of integrals}
\label{sec:defintegrals}

Let $ H $ be a regular multiplier Hopf algebra. For technical reasons we shall suppose throughout that $ H $ admits a faithful 
linear functional, see Definition \ref{defflf}. 
As noted before that definition, we may therefore view $H$
as a subalgebra of the algebraic multiplier algebra $ \M(H) $. 
An analogous statement applies to tensor powers of $ H $.

Assume that $ \omega $ is a linear functional on $ H $. Then we define for any $ f \in H $ a multiplier $ (\id \otimes \omega)\Delta(f) \in \M(H) $ by 
\begin{align*}
(\id \otimes \omega)\Delta(f) \cdot g &= (\id \otimes \omega)\gamma_l(f \otimes g) \\
g \cdot (\id \otimes \omega)\Delta(f) &= (\id \otimes \omega)\rho_l(g \otimes f). 
\end{align*}
To check that this is indeed a two-sided multiplier observe that  
$$
(f \otimes 1) \gamma_l(g \otimes h) = \rho_l(f \otimes g)(h \otimes 1)
$$
for all $ f,g,h \in H $. In a similar way we define $ (\omega \otimes \id)\Delta(f) \in \M(H) $ by 
\begin{align*}
(\omega \otimes \id)\Delta(f) \cdot g &= (\omega \otimes \id)\gamma_r(f \otimes g) \\
g \cdot (\omega \otimes \id)\Delta(f) &= (\omega \otimes \id)\rho_r(g \otimes f).
\end{align*}
\begin{definition}
Let $ H $ be a regular multiplier Hopf algebra. A linear functional $ \phi: H \rightarrow \mathbb{K} $ is called a left invariant integral%
\label{nom:left_integral}
 if 
\begin{equation*}
(\id \otimes \phi)\Delta(f) = \phi(f) 1
\end{equation*}
for all $ f \in H $. Similarly, a linear functional $ \psi: H \rightarrow \mathbb{K} $ is called a right invariant integral%
\label{nom:right_integral}
if 
\begin{equation*}
(\psi \otimes \id)\Delta(f) = \psi(f) 1
\end{equation*}
for all $ f \in H $. 
\end{definition}
\begin{definition} 
	\label{defmhopfintegrals}
A regular multiplier Hopf algebra with integrals is a regular multiplier Hopf algebra $ H $ together with 
a faithful left invariant functional $ \phi: H \rightarrow \mathbb{K} $ and a faithful right invariant 
functional $ \psi: H \rightarrow \mathbb{K} $. 
\end{definition} 
Note that the existence of one of $ \phi $ or $ \psi $ is sufficient. More precisely, one obtains a faithful right/left invariant functional 
by taking a faithful left/right invariant functional and precomposing with the antipode. 

It can be shown that left/right invariant integrals are always unique up to a scalar, see Section 3 in \cite{vDadvances}. 
Moreover, they are related by a modular element, that is, there exists an invertible multiplier $ \delta \in M(H) $
\nomenclature{$\delta$}{modular element}%
such that 
$$
(\phi \otimes \id) \Delta(f) = \phi(f) \delta 
$$
for all $ f \in H $.
The multiplier $ \delta $ is a group-like element in the sense that 
$$
\Delta(\delta) = \delta \otimes \delta, \quad \epsilon(\delta) = 1, \quad S(\delta) = \delta^{-1}. 
$$
It follows that a right invariant integral $\psi$ is defined by
\[
 \psi(f) = \phi(f\delta)= \phi(\delta f).
\]

We say that $ H $ is unimodular if $ \delta = 1 $, or equivalently if we can choose $ \psi = \phi $.

For ordinary (unital) Hopf algebras, the existence of integrals is closely related to cosemisimplicity. 

\begin{prop} \label{cosemisimpleintegral}
Let $ H $ be a Hopf algebra. Then the following conditions are equivalent. 
\begin{bnum}
\item[a)] $ H $ is cosemisimple. 
\item[b)] $ H $ admits a left and right invariant integral $ \phi $ such that $ \phi(1) = 1 $. 
\end{bnum} 
\end{prop} 

\begin{proof} $ a) \Rightarrow b) $ Assume that $ H = \bigoplus_{\lambda \in \Lambda} C_\lambda $ is written as a direct sum of its simple coalgebras, 
and without loss of generality let us write $ C_0 \subset H $ for the one-dimensional simple subcoalgebra spanned by $ 1 \in H $. 
Then we can define $ \phi: H \rightarrow \mathbb{K} $ to be the projection onto $ C_0 \cong \mathbb{K} $.  
Since $ \Delta(C_\lambda) \subset C_\lambda \otimes C_\lambda $ it is straightforward to check that $ \phi $ is a left and right invariant integral
such that $ \phi(1) = 1 $. 

$ b) \Rightarrow a) $. Let us only sketch the argument. Firstly, using the invariant integral $ \phi $ one proves that every comodule for $ H $ 
is a direct sum of simple comodules. This is then shown to be equivalent to 
$ H $ being cosemisimple. We refer to chapter 14 of \cite{Sweedler} for the details. 
\end{proof}

\subsubsection{The dual multiplier Hopf algebra} 
\label{sec:dual_multiplier_Hopf_algebra}

Given a regular multiplier Hopf algebra with integrals we shall introduce the dual multiplier Hopf algebra 
and discuss the Biduality Theorem. These constructions and results are due to van Daele \cite{vDadvances}. 

We define $ \hat{H} $
\nomenclature[s$H$3]{$\hat{H}$}{dual of a multiplier Hopf algebra, with operations $\hat{\Delta}$, $\hat{\epsilon}$, $\hat{S}$, \textit{etc.}}%
as the linear subspace of the dual space $ H^* = \Hom(H, \mathbb{K}) $ given by all functionals of the 
form $ \F(f) $ for $ f \in H $ where
$$ 
(\F(f), h) = \F(f)(h) = \phi(h f). 
$$
\label{nom:F_Hopf}%
It can be shown that one obtains the same space of linear functionals upon replacing $ \phi $ by $ \psi $, or reversing the 
order of multiplication under the integral in the above formula, see \cite{vDadvances}. 
Note that these other choices correspond to reversing the multiplication 
or comultiplication of $ H $, respectively.

Using the evaluation of linear functionals, we obtain the canonical pairing $ \hat{H} \times H \rightarrow \mathbb{K} $ by 
$$ 
(x, f) = x(f) 
$$ 
\label{nom:pairing}
for $ x \in \hat{H} $ and $ h \in H $. 
Note that the formula for $\F(f)$ above can be used to extend this to a pairing $\hat{H} \times \M(H) \to \mathbb{K}$.
In a similar way one obtains pairings between tensor powers of $ H $ and $ \hat{H} $. 

Let us now explain how the space $ \hat{H} $ can be turned into a regular multiplier Hopf algebra with integrals. We point out, however, that we work with 
the \emph{opposite comultiplication} on $ \hat{H} $ compared to the conventions in \cite{vDadvances}. 

\begin{theorem} \label{dualmulthopfconstruct}
Let $ H $ be a regular multiplier Hopf algebra with integrals. 
The vector space $ \hat{H} $ becomes a regular multiplier Hopf algebra with 
comultiplication $ \hat{\Delta}: \hat{H} \rightarrow \M(\hat{H} \otimes \hat{H}) $, 
counit $ \hat{\epsilon}: \hat{H} \rightarrow \mathbb{K} $, and antipode $ \hat{S}: \hat{H} \rightarrow \hat{H} $, such that 
\begin{align*}
(x y, f) &= (x \otimes y, \Delta(f)) \\
(\hat{\Delta}(x), f \otimes g) &= (x, gf) \\
\hat{\epsilon}(x) &= (x, 1) \\ 
(\hat{S}(x), f) &= (x, S^{-1}(f)) 
\end{align*}
for $ x, y \in \hat{H} $ and $ f, g \in H $. A left invariant integral for $ \hat{H} $ is given by 
$$
\hat{\phi}(\F(f)) = \epsilon(f) 
$$%
\label{nom:dual_left_integral}
for $ f \in H $. 
\end{theorem} 

\begin{proof} We shall only sketch two parts of the argument; a central part of the proof consists in showing that the above formulas are well-defined. 

Firstly, for $ f, g \in H $ the product $ \F(f) \F(g) $ in $ \hat{H} $ is defined by
$$
\F(f) \F(g) = \phi(S^{-1}(g_{(1)}) f) \F(g_{(2)}) = (\F \otimes \phi) \gamma_l^{-1}(g \otimes f). 
$$
Then formally using 
\begin{align*}
(\F(f) \F(g))(h) &= (\F(f) \otimes \F(g), \Delta(h)) \\ 
&= \phi(h_{(1)} f) \phi(h_{(2)} g) \\
&= \phi(h_{(1)} g_{(2)} S^{-1}(g_{(1)}) f) \phi(h_{(2)} g_{(3)}) \\ 
&= \phi(S^{-1}(g_{(1)}) f) \phi(h g_{(2)}) \\ 
&= \phi(S^{-1}(g_{(1)}) f) \F(g_{(2)})(h)
\end{align*} 
we see that this agrees with the transposition of the coproduct of $ H $. 
From the fact that the product can be described in this way it follows easily that multiplication in $ \hat{H} $ is associative. 

Secondly, assuming that $ \hat{\Delta} $ is well-defined let us check that $ \hat{\phi} $ is left invariant. For this we compute 
\begin{align*}
((\id \otimes \hat{\phi}&) \hat{\Delta}(\F(f)), h) 
= \hat{\phi}(\F(hf)) = \epsilon(f) \epsilon(h) = (1 \hat{\phi}(\F(f)), h)
\end{align*}
for any $ h \in H $. 

For the construction of the remaining structure maps and the verification of the axioms see \cite{vDadvances}. 
\end{proof}

The linear isomorphism $\F:H \to \hat{H}$ is often referred to as the \emph{Fourier transform}. In fact, if $ H = C(G) $ is the algebra of functions on a finite group $ G $, this map corresponds to the canonical identification between $ C(G) $ and the group algebra $ \mathbb{C}[G] $, equipping $ C(G) $ with the convolution product. Of course, if $G$ is abelian, then the convolution algebra of $ G $ is isomorphic to the algebra of functions on the Pontrjagin dual group $ \hat{G} $ with pointwise multiplication. This recovers the usual notion of Fourier transform.

We shall again use the Sweedler notation
$$ 
 \hat{\Delta}(x) = x_{(1)} \otimes x_{(2)} 
$$ 
for the comultiplication of $ \hat{H} $. The compatibility between $ H $ and $ \hat{H} $ can be summarized as follows. 
\begin{prop} \label{multhopfpairingprop}
Let $ H $ be a regular multiplier Hopf algebra with integrals. Then the canonical evaluation pairing between $ \hat{H} $ and $ H $ satisfies 
\begin{align*}
(xy, f) &= (x, f_{(1)}) (y, f_{(2)}), \qquad (x, fg) = (x_{(2)}, f) (x_{(1)}, g) 
\end{align*}
and 
\begin{align*}
(\hat{S}(x), f) &= (x, S^{-1}(f)) 
\end{align*} 
for all $ x,y \in \hat{H} $ and $ f,g \in H $. 
\end{prop} 
Note that the final relation implies 
$$
(\hat{S}^{-1}(x), f) = (x, S(f)) 
$$
for all $ x \in \hat{H} $ and $ f \in H $.

The next major result is the Biduality Theorem, which states that $\hat{\hat{H}} \cong H$ as multiplier Hopf algebras. 
It is important to note, however, that under this isomorphism the canonical pairing $\hat{\hat{H}} \times \hat{H} \to \mathbb{K}$ will not correspond to evaluation between $H$ and $\hat{H}$.  
This is due to our choice to work with the coopposite comultiplication in the definition of the dual. 
We therefore introduce the following convention, which will apply throughout this work.

\begin{definition}
\label{def:reverse_pairing}
 If $H$ is a regular multiplier Hopf algebra, we define the reverse pairing $H \times \hat{H} \to \mathbb{K}$ by
 \[
  (f, x) = (x, S^{-1}(f)), \qquad\qquad
 \]
 for all $f\in H$, $x\in \hat{H}$.
\end{definition}

Let us stress that with these conventions, we have $ (f, x) \neq (x, f) $ in general.

Since $ \hat{H} $ is again a regular multiplier Hopf algebra with integrals, we can construct its dual. 
As in the case of finite dimensional Hopf algebras one has a duality result, see \cite{vDadvances}.

\begin{theorem}[Biduality Theorem] \label{bidualitymh}
Let $ H $ be a regular multiplier Hopf algebra with integrals. Then the dual of $ \hat{H} $ is isomorphic to $ H $.  The isomorphism is given by
\[
 B : H \to \hat{\hat{H}}; \qquad B(f)(x) = (f,x),
\]
for $f\in H$, $x\in\hat{H}$.
\end{theorem}

\begin{proof} Again we shall only give a brief sketch of the argument. 
The key step in the proof---which we will not carry out---is to show that $ B $ is well-defined. Injectivity 
then follows from the nondegeneracy of the pairing $ ( \;, \;) $. The fact that $ B $ is compatible with multiplication and comultiplication 
is a consequence of Proposition \ref{multhopfpairingprop}. 
\end{proof}

Let us conclude this section by writing down a formula for the inverse of the Fourier transform $\F$.

\begin{prop}
 \label{prop:Fourier_inverse}
 Let $H$ be a regular multiplier Hopf algebra.  For any $x,y \in \hat{H}$ we have
 \[
  (y, \F^{-1}(x)) = \hat{\phi}(\hat{S}^{-1}(y)x) .
 \]
\end{prop}

\begin{proof}
 Let us first observe that for any $x\in \hat{H}$, $f\in H$ we have the formula
 \[
  x \F(f)  = (x, S^{-1}(f_{(1)}))\,\F(f_{(2)}),
 \]
 which can be confirmed by pairing each side with some $g\in H$ and using the left invariance of $\phi$.
 Then, using the proposed formula for $\F^{-1}$ from the proposition, we compute
 \begin{align*}
  (y, \F^{-1} \F(f)) &= \hat{\phi}(\hat{S}^{-1}(y) \F(f)) = \hat{\phi}(\F(f_{(2)})) (y, f_{(1)}) = (y, f), 
 \end{align*}
 for all $ f \in H, y \in \hat{H} $.  Since $\F$ is surjective by definition, we see that $\F$ is invertible with inverse $\F^{-1}$ as proposed.
\end{proof}

If, as above, we identify the dual of $\hat{H}$ with $H$ via the reverse pairing $H \times \hat{H} \to \mathbb{K}$, then we can reinterpret the relation $ \F^{-1} \F = \id $ by saying that $\hat{\F} = \F^{-1}$.
Indeed, notice that 
$$
(\hat{\F}(x), y) = \hat{\phi}(yx) = (\hat{S}(y), \F^{-1}(x)) = (\F^{-1}(x), y) 
$$
\nomenclature[o$F$2]{$\hat{\F}$}{inverse linear isomorphism $\hat{H}\to H$}%
for all $ x, y \in \hat{H} $. 

This simple relation between $ \F $ and $ \hat{\F} $ is a feature of our conventions to work with the coopposite comultiplication on the dual.

\subsection{The Drinfeld double - algebraic level} \label{secddalgebraic}

In this section we discuss the algebraic version of the Drinfeld double construction. This construction produces a regular multiplier 
Hopf algebra $ L \bowtie K $ out of two regular multiplier Hopf algebras $ K, L $ equipped with an invertible skew-pairing. 
For ordinary Hopf algebras---that is, in the situation that $ K $ and $ L $ are unital---this is covered in all standard textbooks, see 
for instance \cite{KS}. For a more detailed account in the non-unital setting we refer to \cite{DvD}.

This construction will be crucial for the later chapters, since the quantized convolution algebra of a complex semisimple group is defined as the Drinfeld double of the quantized convolution algebra of the associated compact form, see Section \ref{seccqg} for a precise statement.

Let us first discuss skew-pairings in the unital case. If $ K $ and $ L $ are Hopf algebras then a skew-pairing between $ L $ and $ K $ is simply a bilinear map 
$ \tau: L \times K \rightarrow \mathbb{K} $ satisfying
\begin{align*}
\tau(xy, f) &= \tau(x, f_{(1)}) \tau(y, f_{(2)}) \\ 
\tau(x, fg) &= \tau(x_{(1)}, g) \tau(x_{(2)}, f) \\
\tau(1, f) &= \epsilon_K(f) \\ 
\tau(x, 1) &= \epsilon_L(x)
\end{align*}
for all $ f,g \in K $ and $ x,y \in L $. 

In the non-unital case this needs to be phrased more carefully. Assume that $ K, L $ are regular multiplier Hopf algebras 
and let $ \tau: L \times K \rightarrow \mathbb{K} $ be a bilinear map. For $ x \in L $ and $ f \in K $ we define linear
maps $ \tau_x: K \rightarrow \mathbb{K} $ and $ _f\tau: L \rightarrow \mathbb{K} $ by 
$$
\tau_x(g) = \tau(x, g), \qquad _f\tau(y) = \tau(y, f). 
$$ 
We can then define multipliers 
$ \tau(x_{(1)}, f) x_{(2)} = (_f \tau \otimes \id) \Delta_L(x), \tau(x_{(2)}, f) x_{(1)} = (\id \otimes _f\tau) \Delta_L(x) \in \M(L) $ 
by first multiplying the outer leg of the coproduct of $ x $ with elements of $ L $, and then applying $ _f\tau $ to the first leg. 
In the same way one obtains multipliers 
$ \tau(x, f_{(1)}) f_{(2)} = (\tau_x \otimes \id) \Delta_K(f), \tau(x, f_{(2)}) f_{(1)} = (\id \otimes \tau_x) \Delta_K(f) \in \M(K) $.  

We will say that the pairing $ \tau $ is regular if the following two conditions are satisfied. 
Firstly, we require that all the multipliers defined above are in fact contained in $ L $ and $ K $, respectively, and not just in $ \M(L) $ 
and $ \M(K) $. Secondly, we require that the linear span of all $ \tau(x_{(1)}, f) x_{(2)} $ is equal to $ L $, 
and the linear span of all $ \tau(x, f_{(1)}) f_{(2)} $ is equal to $ K $, as well as the analogous conditions for the 
flipped comultiplications. It is shown in \cite{DvD} that this implies that the multiplier 
$$
\tau(x_{(1)}, f_{(1)}) x_{(2)} \otimes f_{(2)} 
$$ 
of $ \M(L \otimes K) $ is in fact contained in $ L \otimes K $. 
The same holds for the multipliers obtained by flipping the comultiplications in this formula in one or both factors.

\begin{definition} \label{defskewpairing}
Let $ K, L $ be regular multiplier Hopf algebras. A skew-pairing between $ L $ and $ K $ is a regular bilinear map 
$ \tau: L \times K \rightarrow \mathbb{K} $ which 
satisfies 
\begin{bnum}
\item[a)] $ \tau_{xy}(f) = \tau_x(\id \otimes \tau_y) \Delta_K(f) = \tau_y(\tau_x \otimes \id) \Delta_K(f) $
\item[b)] $ _{fg}\tau(x) = _g\tau(\id \otimes _f\tau) \Delta_L(x) = _f\tau(_g\tau \otimes \id) \Delta_L(x) $ 
\end{bnum}
for all $ f,g \in K $ and $ x,y \in L $. 
\end{definition} 

Occasionally we will identify $ \tau $ with its associated linear map $ L \otimes K \rightarrow \mathbb{K} $ and 
write $ \tau(x \otimes f) $ instead of $ \tau(x,f) $.

Note that the conditions in Definition \ref{defskewpairing} are well-defined 
by the regularity assumption on $ \tau $. Informally, these conditions can be written as 
$$ 
\tau(xy, f) = \tau(x, f_{(1)}) \tau(y, f_{(2)}), \quad \tau(x, fg) = \tau(x_{(1)}, g) \tau(x_{(2)}, f) 
$$ 
for $ f, g \in K, x, y \in L $, as in the unital case explained above. 

Regularity implies that the skew-pairing $ \tau $ induces essential module structures of $ K $ acting on $ L $ and vice versa. 
More precisely, $ K $ becomes an essential left $ L $-module and an essential right $ L $-module by 
$$ 
x \triangleright f = \tau(x, f_{(2)}) f_{(1)}, \qquad f \triangleleft x = \tau(x, f_{(1)}) f_{(2)}
$$ 
and $ L $ becomes an essential left and right $ K $-module by 
$$ 
f \triangleright x = \tau(x_{(1)}, f) x_{(2)}, \qquad x \triangleleft f = \tau(x_{(2)}, f) x_{(1)}.
$$ 
This means in particular that $ \tau $ extends canonically to bilinear pairings $ \M(L) \times K \rightarrow \mathbb{K} $ 
and $ L \times \M(K) \rightarrow \mathbb{K} $ by setting 
$$
\tau(x, f) = \epsilon_K(x \triangleright f) = \epsilon_K(f \triangleleft x), \qquad 
\tau(y, g) = \epsilon_L(g \triangleright y) = \epsilon_L(y \triangleleft g)
$$
for $ x \in \M(L), f \in K $ and $ y \in L, g \in \M(K) $, respectively. 
With this notation, we get in particular 
$$ 
\tau(x, 1) = \epsilon_L(x), \qquad \tau(1, f) = \epsilon_K(f) 
$$ 
for all $ x \in L, f \in K $, again reproducing the formulas in the unital case.

The convolution product $ \tau \eta $ of two regular pairings $ \tau, \eta: L \otimes K \rightarrow \mathbb{K} $ is defined by 
$$
(\tau \eta)(x \otimes f) = \eta(x_{(1)} \otimes f_{(1)}) \tau(x_{(2)} \otimes f_{(2)}).
$$ 
This defines an associative multiplication on the linear space of all regular pairings $ L \otimes K \rightarrow \mathbb{K} $, 
with unit element $ \epsilon_L \otimes \epsilon_K $. 
Any skew-pairing $\tau: L \otimes K \to \mathbb{K}$ between regular multiplier Hopf algebras is convolution invertible, with inverse given by
\[
 \tau^{-1}(x, f) = \tau(S_L(x),f) = \tau(x, S_K^{-1}(f)).
\]

Let us now introduce the Drinfeld double of skew-paired regular multiplier Hopf algebras. 

\begin{definition} 
\label{def:Drinfeld_double}
Let $ K, L $ be regular multiplier Hopf algebras and let $ \tau: L \otimes K \rightarrow \mathbb{K} $ be a skew-pairing between them. 
The Drinfeld double $ L \bowtie K $ is the regular multiplier Hopf algebra with underlying vector 
space $ L \otimes K $, equipped with the multiplication 
$$
(x \bowtie f)(y \bowtie g) = x \tau(y_{(1)}, f_{(1)}) \,y_{(2)}  \bowtie f_{(2)} \, \tau^{-1}(y_{(3)}, f_{(3)}) g, 
$$
the coproduct 
$$ 
\Delta_{L \bowtie K}(x \bowtie f) = (x_{(1)} \bowtie f_{(1)}) \otimes (x_{(2)} \bowtie f_{(2)}),  
$$
antipode 
$$
S_{L \bowtie K}(x \bowtie f) = \tau^{-1}(x_{(1)}, f_{(1)}) S_L(x_{(2)}) \bowtie S_K(f_{(2)}) \tau(x_{(3)}, f_{(3)}),  
$$
and the counit 
$$
 \epsilon_{L \bowtie K}(x \bowtie f) = \epsilon_L(x) \epsilon_K(f) .
$$ 
\end{definition} 
It is proved in \cite{DvD}  that these structures indeed turn $ L \bowtie K $ into a regular multiplier Hopf algebra such 
that $ \M(L \bowtie K) $ contains both $ L $ and $ K $ as multiplier Hopf subalgebras in a natural way. 
Note that the formula for the antipode is forced by requiring that it be compatible with the antipodes $S_L$ and $S_K$ on the subalgebras $L$ and $K$ in $ \M(L \bowtie K) $, thanks to the formula
\[
 S_{L\bowtie K}(x\bowtie f) = S_{L\bowtie K}((x\bowtie 1)(1\bowtie f)) = (1 \bowtie S_K(f))\,(S_L(x)\bowtie 1),
\]
for $x\in L$, $f\in K$.

Formally, associativity of the multiplication in $ L \bowtie K $ follows from 
\begin{align*}
&((x \bowtie f)(y \bowtie g))(z \bowtie h) = \left(x \tau(y_{(1)}, f_{(1)}) y_{(2)} \bowtie f_{(2)} \tau^{-1}(y_{(3)}, f_{(3)}) g\right)(z \bowtie h) \\ 
&= x \tau(y_{(1)}, f_{(1)}) y_{(2)} \tau(z_{(1)}, f_{(2)} g_{(1)}) z_{(2)} \bowtie 
f_{(3)} g_{(2)} \tau^{-1}(z_{(3)}, f_{(4)} g_{(3)}) \tau^{-1}(y_{(3)}, f_{(5)}) h \\ 
&= x \tau(y_{(1)} z_{(2)}, f_{(1)}) \tau(z_{(1)}, g_{(1)}) y_{(2)} z_{(3)} \bowtie 
f_{(2)} \tau^{-1}(y_{(3)} z_{(4)}, f_{(3)}) g_{(2)} \tau^{-1}(z_{(5)}, g_{(3)}) h \\  
&= (x \bowtie f)\left(y \tau(z_{(1)}, g_{(1)}) z_{(2)} \bowtie g_{(2)} \tau^{-1}(z_{(3)}, g_{(3)}) h \right)\\  
&= (x \bowtie f)((y \bowtie g)(z \bowtie h)) 
\end{align*} 
for $ x \bowtie f, y \bowtie g, z \bowtie h \in L \bowtie K $. To check that $ \Delta_{L \bowtie K} $ is an algebra homomorphism one may 
formally compute 
\begin{align*}
&\Delta_{L \bowtie K}((x \bowtie f)(y \bowtie g)) 
= \Delta_{L \bowtie K}(x \tau(y_{(1)}, f_{(1)}) y_{(2)} \bowtie f_{(2)} \tau^{-1}(y_{(3)}, f_{(3)}) g) \\ 
&= x_{(1)} \tau(y_{(1)}, f_{(1)}) y_{(2)} \bowtie f_{(2)}g_{(1)} \otimes 
x_{(2)} y_{(3)} \bowtie \tau^{-1}(y_{(4)}, f_{(4)}) f_{(3)} g_{(2)} \\ 
&= (x_{(1)} \tau(y_{(1)}, f_{(1)}) y_{(2)} \bowtie f_{(2)} \tau^{-1}(y_{(3)}, f_{(3)}) g_{(1)})  \\
&\hspace{3cm} \otimes (x_{(2)} \tau(y_{(4)}, f_{(4)}) y_{(5)} \bowtie f_{(5)} \tau^{-1}(y_{(6)}, f_{(6)}) g_{(2)}) \\  
&= (x_{(1)} \bowtie f_{(1)} \otimes x_{(2)} \bowtie f_{(2)}) (y_{(1)} \bowtie g_{(1)} \otimes y_{(2)} \bowtie g_{(2)}) \\  
&= \Delta_{L \bowtie K}(x \bowtie f) \Delta_{L \bowtie K}(y \bowtie g)
\end{align*} 
for $ x \bowtie f, y \bowtie g \in L \bowtie K $. It is clear from the definition that $ \Delta_{L \bowtie K} $ is coassociative and that 
$ \epsilon_{L \bowtie K} $ is a counit for $ \Delta_{L \bowtie K} $. 
For the antipode axiom we compute 
\begin{align*}
&m_{L \bowtie K}(S_{L \bowtie K} \otimes \id)\Delta_{L \bowtie K}(x \bowtie f) 
= m_{L \bowtie K}(S_{L \bowtie K} \otimes \id)(x_{(1)} \bowtie f_{(1)} \otimes x_{(2)} \bowtie f_{(2)}) \\
&= m_{L \bowtie K}(\tau^{-1}(x_{(1)}, f_{(1)}) S_L(x_{(2)}) \bowtie S_K(f_{(2)}) \tau(x_{(3)}, f_{(3)}) \otimes x_{(4)} \bowtie f_{(4)}) \\ 
&
= \tau^{-1}(x_{(1)}, f_{(1)}) S_L(x_{(2)}) \tau(x_{(4)}, S_K(f_{(4)})) x_{(5)}  
\\
&
\hspace{5cm} \bowtie S_K(f_{(3)}) \times \tau^{-1}(x_{(6)}, S_K(f_{(2)})) \tau(x_{(3)}, f_{(5)}) f_{(6)} 
 \\ 
&= \tau^{-1}(x_{(1)}, f_{(1)}) S_L(x_{(2)}) x_{(3)} \bowtie 
\tau^{-1}(x_{(4)}, S_K(f_{(2)})) S_K(f_{(3)}) f_{(4)} \\ 
&= \tau^{-1}(x_{(1)}, f_{(1)}) 1\bowtie 
\tau^{-1}(x_{(2)}, S_K(f_{(2)})) 1 \\ 
&= \epsilon_{L \bowtie K}(x \bowtie f) 1, 
\end{align*} 
where $ m_{L \bowtie K} $ denotes the multiplication of $ L \bowtie K $. In a similar way one verifies the other antipode condition. 

Let $ K, L $ be regular multiplier Hopf algebras and let $ \rho : L \otimes K \rightarrow \mathbb{K} $ be an invertible skew-pairing. 
Then $ \tau: K \otimes L \rightarrow \mathbb{K} $ given by 
$$
\tau(f, x) = \rho(x, S_K^{-1}(f)) = \rho (S_L(x), f) 
$$ 
is also an invertible skew-pairing. 

\begin{definition}
\label{def:Rosso_form}
With the above notation, the Rosso form is the bilinear form $\kappa$ on the Drinfeld double $ L \bowtie K $ defined by 
$$
\kappa(x \bowtie f, y \bowtie g) = \rho(y, f) \rho(S_L(x), S_K^{-1}(g)) 
  = \tau(S_K(f), y) \tau(g, S_L(x)).
$$
\end{definition}

The left adjoint action of a regular multiplier Hopf algebra $ H $ on itself is defined by 
$$ 
\ad(f)(g) = f \rightarrow g = f_{(1)} g \, S(f_{(2)}). 
$$
\label{nom:adjoint_action1}
A bilinear form $ \kappa $ on $ H $ is called $ \ad $-invariant if 
$$
\kappa(f \rightarrow g, h) = \kappa(g, S(f) \rightarrow h) 
$$
for all $ f,g,h \in H $. 

For the following result compare Section 8.2.3 in \cite{KS}. 

\begin{prop} \label{rossoadinvariant}
Let $ K, L $ be regular multiplier Hopf algebras equipped with an invertible skew-pairing $ \rho: L \times K \rightarrow \mathbb{K} $. 
Then the Rosso form on the Drinfeld double $ L \bowtie K $ is $ \ad $-invariant. 
\end{prop} 

\proof Since $ L \bowtie K $ is generated as an algebra by $ L $ and $ K $ it suffices to consider the adjoint action of $ x \in L $ 
and $ f \in K $ on an element $ y \bowtie g $ of the double. We compute 
\begin{align*}
\ad(x)(y \bowtie g) &= (x_{(1)} \bowtie 1)(y \bowtie g)(S_L(x_{(2)}) \bowtie 1) \\ 
&= x_{(1)} y\, S_L(x_{(3)}) \rho(S_L(x_{(4)}), g_{(1)})  \bowtie  \rho^{-1}(S_L(x_{(2)}), g_{(3)}) g_{(2)}.
\end{align*} 
Noting that $ \rho(S_L(x), S_K(f)) = \rho(x, f)$ we obtain
\begin{align*}
&\kappa(x \rightarrow (y \bowtie g), z \bowtie h) \\
&= \kappa(x_{(1)} y S_L(x_{(3)}) \rho^{-1}(x_{(4)}, g_{(1)}) \bowtie \rho^{-1}(S_L(x_{(2)}), g_{(3)}) g_{(2)}, z \bowtie h) \\
&= \rho^{-1}(S_L(x_{(1)} y S_L(x_{(3)})), h) \rho^{-1}(x_{(4)}, g_{(1)}) \rho^{-1}(S_L(x_{(2)}), g_{(3)}) \rho(z, g_{(2)}) \\
&= \rho^{-1}(S_L(x_{(1)}), h_{(1)}) \rho^{-1}(S_L(y), h_{(2)}) \rho^{-1}(S^2_L(x_{(3)})), h_{(3)})  \\
&\hspace{4cm} \times \rho^{-1}(x_{(4)}, g_{(1)}) \rho^{-1}(S_L(x_{(2)}), g_{(3)}) \rho(z, g_{(2)}) \\
&= \rho^{-1}(S_L(x_{(1)}), h_{(1)}) \rho^{-1}(S_L(y), h_{(2)}) \rho^{-1}(S_L^2(x_{(3)}), h_{(3)}) \rho(S_L(x_{(4)}) z S_L^2(x_{(2)}), g) \\
&= \kappa\left (y \bowtie g, S_L(x_{(4)}) z S_L^2(x_{(2)}) \rho^{-1}(S_L(x_{(1)}), h_{(1)}) \bowtie \rho^{-1}(S_L^2(x_{(3)}), h_{(3)}) h_{(2)}\right) \\
&= \kappa\left(y \bowtie g, S_L(x) \rightarrow (z \bowtie h)\right).  
\end{align*}
Similarly, we have 
\begin{align*}
\ad(f)(y \bowtie g) &= (1 \bowtie f_{(1)})(y \bowtie g)(1 \bowtie S_K(f_{(2)})) \\ 
&= \rho(y_{(1)}, f_{(1)}) y_{(2)} \bowtie \rho^{-1}(y_{(3)}, f_{(3)}) f_{(2)} g S_K(f_{(4)})
\end{align*} 
and so
\begin{align*}
&\kappa(f \rightarrow (y \bowtie g), z \bowtie h) \\
&= \kappa(\rho(y_{(1)}, f_{(1)}) y_{(2)} \bowtie \rho^{-1}(y_{(3)}, f_{(3)}) f_{(2)} g S_K(f_{(4)}), z \bowtie h) \\
&= \rho(y_{(1)}, f_{(1)}) \rho^{-1}(S_L(y_{(2)}), h) \rho^{-1}(y_{(3)}, f_{(3)}) \rho(z, f_{(2)} g S_K(f_{(4)})) \\
& =\rho(y_{(1)}, f_{(1)}) \rho(y_{(2)},S_K^{-2}(h)) \rho(y_{(3)},S_K^{-1}(f_{(3)}))
   \rho(z_{(1)},S_K(f_{(4)})) \rho(z_{(2)}, g) \rho(z_{(3)},f_{(2)}) 
 \\
&= \rho(y,  S_K^{-1}(f_{(3)})   S_K^{-2}(h)  f_{(1)})  \rho(z_{(1)}, S_K(f_{(4)})) \rho(z_{(2)}, g) \rho^{-1}(z_{(3)}, S_K(f_{(2)})) \\
&= \kappa(y \bowtie g, \rho(z_{(1)}, S_K(f_{(4)})) z_{(2)} \bowtie \rho^{-1}(z_{(3)}, S_K(f_{(2)})) S_K(f_{(3)}) h S_K^2(f_{(1)})) \\
&= \kappa(y \bowtie g, S_K(f) \rightarrow (z \bowtie h)).  
\end{align*}
This yields the claim. \qed 

\section{Quantized universal enveloping algebras} \label{chuqg}

In this chapter we collect background material on quantized universal enveloping algebras. We give in particular a detailed account of 
the construction of the braid group action and PBW-bases, and discuss the finite dimensional representation theory in the setting  
that the base field $ \mathbb{K} $  is an arbitrary field and the deformation parameter $ q \in \mathbb{K}^\times $ is not a root of unity. 
Our presentation mainly follows the textbooks \cite{Jantzenqg}, \cite{Josephbook} and \cite{Lusztigbook}. 

\subsection{$ q $-calculus} \label{secqnumbers}

Let $ \mathbb{K} $ be a field.  We will assume $q\in\mathbb{K}^\times$ is not a root of unity.

For $ n \in \mathbb{Z} $ we write 
$$
[n]_q = \frac{q^n - q^{-n}}{q - q^{-1}} = q^{-n + 1} + q^{-n + 3} + \cdots + q^{n - 1} 
$$
\nomenclature[c,q1]{$[n]_q$}{$q$-number $\displaystyle[n]_q=\frac{q^n-q^{-n}}{q-q^{-1}}$}%
for the corresponding $ q $-number. For $ n \in \mathbb{N} $ we set 
$$
[n]_q! = \prod_{k = 1}^n [k]_q, 
$$
\nomenclature[c,q2]{$[n]_q"!$}{$q$-factorial}%
and in addition we define $ [0]_q! = 1 $. 

The $ q $-binomial coefficients are defined by 
$$
\begin{bmatrix} n \\ m \end{bmatrix}_q = \frac{[n]_q [n - 1]_q \cdots [n - m + 1]_q}{[m]_q!} 
$$
\nomenclature[c,q3]{$\displaystyle\begin{bmatrix} n \\ m \end{bmatrix}_q$}{$q$-binomial coefficient}%
for $ n \in \mathbb{Z} $ and $ m \in \mathbb{N} $. In addition we declare $ \begin{bmatrix} n \\ 0 \end{bmatrix}_q = 1 $ for all $ n \in \mathbb{Z} $ 
and $ \begin{bmatrix} n \\ m \end{bmatrix}_q  = 0 $ for $ m < 0 $. 
Note that if $ n \in \mathbb{N}_0 $ we have 
$$
\begin{bmatrix} n \\ m \end{bmatrix}_q = \frac{[n]_q!}{[n - m]_q! [m]_q!} 
$$
for $ 0 \leq m \leq n $. We shall often omit the subscripts $ q $ if no confusion is likely. 

\begin{lemma} \label{lemqcalc1}
We have
$$
q^m 
\begin{bmatrix}
n \\ 
m 
\end{bmatrix}
+ q^{m - n - 1} 
\begin{bmatrix}
n \\ 
m - 1
\end{bmatrix}
= \begin{bmatrix}
n + 1\\ 
m 
\end{bmatrix}
$$
and
$$
q^{-m} 
\begin{bmatrix}
n \\ 
m 
\end{bmatrix}
+ q^{n + 1 - m} 
\begin{bmatrix}
n \\ 
m - 1
\end{bmatrix}
= \begin{bmatrix}
n + 1\\ 
m 
\end{bmatrix}
$$
for any $n\in\mathbb{Z}$ and $m\in\mathbb{N}$. 
\end{lemma} 

\begin{proof} 
We calculate 
\begin{align*}
 \left(
 q^m 
\begin{bmatrix}
n \\ 
m 
\end{bmatrix} \right.
&
\left.
+ q^{m - n - 1} 
\begin{bmatrix}
n 
\\ 
m - 1
\end{bmatrix} 
\right) \times [m]
\\
&= q^m [n]\ldots[n-m+2][n-m+1] + q^{m-n-1}[n]\ldots[n-m+2][m] \\ 
&= \left( q^m[n-m+1] + q^{m-n-1}[m] \right) [n]\ldots[n-m+2] \\ 
&= \left(
{\textstyle \frac{q^{n + 1} - q^{2m - n - 1} + q^{2m - n - 1} - q^{-n - 1}}{q - q^{-1}} }
\right) [n]\ldots[n-m+2] \\
&= [n+1][n]\ldots[n-m+2] \\
&= \begin{bmatrix}
n + 1\\ 
m 
\end{bmatrix}\times [m].
\end{align*}
This yields the first claim. The second follows by replacing $q$ by $q^{-1}$.
\end{proof}

\begin{lemma} \label{lembinomialvanishing}
For any $n\in\mathbb{N}_0$ we have 
$$
\sum_{m = 0}^n (-1)^m 
\begin{bmatrix}
n \\ 
m 
\end{bmatrix}
q^{m(n-1)} =
\begin{cases}
 1 & \text{if }n=0,\\
 0 & \text{if }n>0
\end{cases} 
$$
and
$$
\sum_{m = 0}^n (-1)^m 
\begin{bmatrix}
n \\ 
m 
\end{bmatrix}
q^{-m(n-1)} =
\begin{cases}
 1 & \text{if }n=0,\\
 0 & \text{if }n>0.
\end{cases} 
$$
\end{lemma} 

\begin{proof} For $n=0$ or $1$ the claim can be checked directly. For $n\geq2$ we compute, using Lemma \ref{lemqcalc1},
\begin{align*}
\sum_{m = 0}^n (-1)^m &
\begin{bmatrix}
n \\ 
m 
\end{bmatrix}
q^{m(n-1)} \\
&= 
\sum_{m = 0}^{n - 1} (-1)^m q^m 
\begin{bmatrix}
n - 1 \\ 
m 
\end{bmatrix}
q^{nm-m}
+ \sum_{m = 1}^n 
(-1)^m q^{m - n} 
\begin{bmatrix}
n - 1 \\ 
m - 1
\end{bmatrix}
q^{nm-m} \\ 
&= 
\sum_{m = 0}^{n - 1} (-1)^m q^{nm}
\begin{bmatrix}
n - 1 \\ 
m 
\end{bmatrix}
+ \sum_{m = 1}^n
(-1)^m q^{n(m - 1)} 
\begin{bmatrix}
n - 1 \\ 
m - 1
\end{bmatrix}
= 0.
\end{align*}
For $ n = 1 $ the claim can be checked directly.
The second equation follows by replacing $q$ by $q^{-1}$.
\end{proof}

Consider the special case $ \mathbb{K} = \mathbb{Q}(q) $. It is clear by definition 
that $ [n] \in \mathbb{Z}[q, q^{-1}] \subset \mathbb{Q}(q) $. 
Moreover, from Lemma \ref{lemqcalc1} one obtains by induction that all $ q $-binomial coefficients 
are contained in $ \mathbb{Z}[q, q^{-1}] $. That is, the above formulas, suitably interpreted, are valid in $ \mathbb{Z}[q, q^{-1}] $. 
With this in mind, one can dispense of our initial assumption that $q$ is not a root of unity.

\subsection{The definition of $ U_q(\mathfrak{g}) $} \label{secdefuqg}

In this section we explain the construction of the quantized universal enveloping algebra of a semisimple Lie algebra.

\subsubsection{Semisimple Lie algebras}
\label{sec:semisimple_Lie_algebras}

Let $ \mathfrak{g} $ a semisimple Lie algebra over $ \mathbb{C} $ of rank $ N $. 
\nomenclature[o$g$]{$\lie{g}$}{semisimple Lie algebra over $\mathbb{C}$}%
\nomenclature{$N$}{rank of $\lie{g}$}%
We fix a Cartan subalgebra $ \mathfrak{h} \subset \mathfrak{g} $.
\nomenclature[o$h$]{$\lie{h}$}{Cartan subalgebra of $\lie{g}$}%
We write $\bf\Delta$
\nomenclature[o$\Delta$5]{$\bf\Delta$}{set of roots of $\lie{g}$}%
for the set of simple roots and fix a set $\bf\Delta^+$
\nomenclature[o$\Delta$5]{$\bf\Delta^+$}{set of positive roots of $\lie{g}$}%
of positive roots and denote by $ \simpleroots = \{\alpha_1, \dots, \alpha_N\} $
\nomenclature[o$\Sigma$2]{$\simpleroots$}{set of simple roots of $\lie{g}$}%
the associated set of simple roots.
We write $ (\;,\;) $ for the bilinear form on $ \mathfrak{h}^* $ obtained by rescaling the Killing form such that the 
shortest root $ \alpha $ of $ \mathfrak{g} $ satisfies $ (\alpha,\alpha) = 2 $. 
Moreover we set 
$$ 
d_i = (\alpha_i, \alpha_i)/2
$$
\nomenclature{$d_i$}{$=(\alpha_i,\alpha_i)/2$}
\nomenclature[o$d_i$2]{$d_\alpha$}{$=(\alpha,\alpha)/2$}%
for all $ i = 1, \dots, N $ and let 
$$
\alpha_i^\vee = d_i^{-1} \alpha_i
$$
\nomenclature{$\alpha^\vee$}{coroot $d_\alpha^{-1}\alpha$}%
\nomenclature[s$\alpha^\vee$]{$\alpha^\vee$}{coroot $d_\alpha^{-1}\alpha$}%
be the simple coroot corresponding to $ \alpha_i $. 

Denote by $ \varpi_1, \dots, \varpi_N $
\nomenclature[o$\pi$1]{$\varpi_i$}{fundamental weight}%
the fundamental weights of $ \mathfrak{g} $, 
satisfying the relations $ (\varpi_i, \alpha_j^\vee) = \delta_{ij} $. 
We write 
$$ 
\weights = \bigoplus_{j = 1}^N \mathbb{Z} \varpi_j,
 \qquad \roots = \bigoplus_{j = 1}^N \mathbb{Z} \alpha_j,
 \qquad \roots^\vee = \bigoplus_{j = 1}^N \mathbb{Z} \alpha_j^\vee,
$$ 
\nomenclature[o$P$]{$\weights$}{weight lattice}%
\nomenclature[o$Q$5]{$\roots$}{root lattice}%
\nomenclature[o$Q$8]{$\roots^\vee$}{coroot lattice}%
for the weight, root and coroot lattices of $ \mathfrak{g} $, respectively. Note that $ \roots \subset \weights \subset \mathfrak{h}^* $ and that $\roots^\vee$ identifies with the $\ZZ$-dual of $\weights$ under the pairing.
 
The set $ \weights^+ $
\nomenclature[o$P$1]{$\weights^+$}{set of non-negative integer combinations of fundamental weights}%
of dominant integral weights is the set of all non-negative integer combinations of the fundamental weights. 
We also write $ \roots^+ $
\nomenclature[o$Q$6]{$\roots^+$}{set of non-negative integer combinations of simple roots}%
for the non-negative integer combinations of the simple roots.

The Cartan matrix for $ \mathfrak{g} $ is the matrix $(a_{ij})_{1\leq i,j\leq N}$ with coefficients 
$$ 
a_{ij} = (\alpha_i^\vee, \alpha_j) = \frac{2 (\alpha_i, \alpha_j)}{(\alpha_i, \alpha_i)}. 
$$ 
\nomenclature{$a_{ij}$}{$= (\alpha_i^\vee, \alpha_j) = \frac{2 (\alpha_i, \alpha_j)}{(\alpha_i, \alpha_i)}$, coefficient of Cartan matrix}%

We write $W$
\nomenclature[o$W$21]{$W$}{Weyl group}%
for the Weyl group, that is, the finite group of automorphisms of $\weights$ generated by the reflections $s_i$ in the hyperplanes orthogonal to the simple roots $\alpha_i\in\bf\Delta^+$, namely
\[
 s_i(\mu) = \mu - (\alpha_i^\vee,\mu)\alpha_i
\]
\nomenclature{$s_i$}{reflection in hyperplane orthogonal to $\alpha_i$}%
for all  $\mu\in \weights$.

Throughout we fix the smallest positive integer $ L $
\nomenclature{$L$}{smallest positive integer such that $(\varpi_i, \varpi_j)\in\frac{1}{L}\mathbb{Z}$ for all $i,j$}%
such that the numbers $ (\varpi_i, \varpi_j) $ take values 
in $ \frac{1}{L} \mathbb{Z} $ for all $ 1 \leq i,j \leq N $, see Section 1 in \cite{Sawinrootsofunity} for the explicit values of $ L $ 
in all types and more information.

In the sequel we shall also work over more general ground fields $ \mathbb{K} $. We will always keep
the notions of weights and roots from the classical case as above. We remark that associated to the complex semisimple Lie algebra $ \mathfrak{g} $ we 
can construct a Lie algebra over $ \mathbb{K} $ using the Serre presentation of $ \mathfrak{g} $. We will sometimes work implicitly with 
this Lie algebra and, by slight abuse of notation, denote it by $ \mathfrak{g} $ again. This apparent ambiguity is resolved by observing that 
the starting point for all constructions in the sequel is in fact a finite Cartan matrix, rather than a semisimple Lie algebra.

\subsubsection{The quantized universal enveloping algebra without Serre relations}

For practical purposes it is convenient to construct the quantized universal enveloping algebra in two steps, mimicking the 
standard approach to defining Kac-Moody algebras. In this subsection we present the first step, namely the 
definition of a Hopf algebra which will admit the quantized universal enveloping algebra as a quotient. 

Let $ \mathfrak{g} $ be a semisimple Lie algebra. We keep the notation introduced above. 

\begin{definition} 
\label{def:no_Serre}
Let $ \mathbb{K} $ be a field and $ q = s^L \in \mathbb{K}^\times $
be an invertible element such that $ q_i \neq \pm 1 $ for 
all $ 1 \leq i \leq N $ where $ q_i = q^{d_i} $. The algebra $ \tilde{U}_q(\mathfrak{g}) $%
over $ \mathbb{K} $ has 
generators $ K_\lambda $ for $ \lambda \in \weights $, and $ E_i, F_i $ for $ i = 1, \dots, N $, and the defining relations 
\begin{align*}
K_0 &= 1 \\
K_\lambda K_\mu &= K_{\lambda + \mu} \\
K_\lambda E_j K_\lambda^{-1} &= q^{(\lambda, \alpha_j)} E_j \\
K_\lambda F_j K_\lambda^{-1} &= q^{-(\lambda, \alpha_j)} F_j \\ 
[E_i, F_j] &= \delta_{ij} \frac{K_i - K_i^{-1}}{q_i - q_i^{-1}} 
\end{align*} 
for all $ \lambda, \mu \in \weights $ and all $ 1 \leq i, j \leq N $. Here we abbreviate $ K_i = K_{\alpha_i} $ for all simple roots. 
\end{definition} 

Our hypothesis on $ q $ means that we always require $ q^2 \neq 1 $, and in fact $ q^4 \neq 1 $ or $ q^6 \neq 1 $ in the 
case that $ \mathfrak{g} $ contains a component of type $ B_n, C_n, F_4 $ or $ G_2 $, respectively.

Notice that $ q^{(\alpha_i, \alpha_j)} = q_i^{a_{ij}} $.  Therefore we have 
$$
K_i E_j K_i^{-1} = q_i^{a_{ij}} E_j, \qquad K_i F_j K_i^{-1} = q_i^{-a_{ij}} F_j 
$$
for all $ 1 \leq i,j \leq N $.  In particular
\[
 K_iE_iK_i^{-1} = q_i^2E_i, \qquad K_iF_iK_i^{-1} = q_i^{-2}F_i.
\]

The algebra $ \tilde{U}_q(\mathfrak{g}) $ admits a Hopf algebra structure as follows. 

\begin{lemma} \label{lemuqhopf}
The algebra $ \tilde{U}_q(\mathfrak{g}) $ is a Hopf algebra with 
comultiplication $ \hat{\Delta}: \tilde{U}_q(\mathfrak{g}) \rightarrow \tilde{U}_q(\mathfrak{g}) \otimes \tilde{U}_q(\mathfrak{g}) $ given by 
\begin{align*}
\hat{\Delta}(K_\mu) &= K_\mu \otimes K_\mu, \\
\hat{\Delta}(E_i) &= E_i \otimes K_i + 1 \otimes E_i \\
\hat{\Delta}(F_i) &= F_i \otimes 1 + K_i^{-1} \otimes F_i,  
\end{align*}
counit $ \hat{\epsilon}: \tilde{U}_q(\mathfrak{g}) \rightarrow \mathbb{K} $ given by 
\begin{align*}
\hat{\epsilon}(K_\lambda) = 1, \qquad \hat{\epsilon}(E_j) = 0, \qquad \hat{\epsilon}(F_j) = 0, 
\end{align*} 
and antipode $ \hat{S}: \tilde{U}_q(\mathfrak{g}) \rightarrow \tilde{U}_q(\mathfrak{g}) $ given by 
\begin{align*}
\hat{S}(K_\lambda) = K_{-\lambda}, \qquad \hat{S}(E_j) = -E_j K_j^{-1}, \qquad \hat{S}(F_j) = -K_j F_j. 
\end{align*}
\end{lemma} 

\begin{proof} These are straightforward calculations. In order to verify that $ \hat{\Delta} $ extends to a 
homomorphism as stated we compute, for instance, 
\begin{align*}
[\hat{\Delta}(E_i), \hat{\Delta}(F_j)] &= [E_i \otimes K_i + 1 \otimes E_i, F_j \otimes 1 + K_j^{-1} \otimes F_j] \\
&= K_j^{-1} \otimes [E_i, F_j] + [E_i, F_j] \otimes K_i \\ 
&= \frac{\delta_{ij}}{q_i - q_i^{-1}} \hat{\Delta}(K_i - K_i^{-1}).  
\end{align*} 
To see that $ \hat{S} $ defines an antihomomorphism we check for instance 
\begin{align*}
\hat{S}([E_i, F_j]) &= \frac{\delta_{ij}}{q_i - q_i^{-1}} \hat{S}(K_i - K_i^{-1}) \\ 
&= \frac{\delta_{ij}}{q_i - q_i^{-1}} (K_i^{-1} - K_i) \\ 
&= K_j F_j E_i K_i^{-1} - K_j E_i F_j K_i^{-1} \\
&= K_j F_j E_i K_i^{-1} - E_i K_i^{-1} K_j F_j \\ 
&= [-K_j F_j, -E_i K_i^{-1}] = [\hat{S}(F_j), \hat{S}(E_i)].
\end{align*} 
The Hopf algebra axioms are easily verified on generators. \end{proof} 

We will continue to use the Sweedler notation 
$$ 
\hat{\Delta}(X) = X_{(1)} \otimes X_{(2)} 
$$ 
for the coproduct of $ X \in \tilde{U}_q(\mathfrak{g}) $. 

The following fact is frequently useful. 

\begin{lemma} \label{lemflippingEandF}
There exists an algebra automorphism $ \omega: \tilde{U}_q(\mathfrak{g}) \rightarrow \tilde{U}_q(\mathfrak{g}) $ 
given by 
$$
\omega(K_\lambda) = K_{-\lambda}, \qquad \omega(E_j) = F_j, \qquad \omega(F_j) = E_j 
$$
on generators. Moreover $ \omega $ is a coalgebra antihomomorphism. 
\end{lemma} 

\begin{proof} It is straightforward to verify that the defining relations of $ \tilde{U}_q(\mathfrak{g}) $ are preserved by the above 
assignment. For the coalgebra antihomomorphism property notice that 
$$ 
\sigma \hat{\Delta}(\omega(E_j)) = F_j \otimes K_j^{-1} + 1 \otimes F_j = (\omega \otimes \omega)(\hat{\Delta}(E_j)) 
$$
for all $ j = 1, \dots, N $. Here $ \sigma $ denotes the flip map. \end{proof} 

Note that, as a consequence of Lemma \ref{lemflippingEandF}, we have $\omega\circ\hat{S} = \hat{S}^{-1}\circ\omega$.

Let us record the following explicit formulas for certain coproducts in $ \tilde{U}_q(\mathfrak{g}) $. 

\begin{lemma} \label{lemEiFicomult}
We have 
\begin{align*}
\hat{\Delta}(E_j^r) &= \sum_{i = 0}^r q_j^{i (r - i)} 
\begin{bmatrix}
r \\ 
i
\end{bmatrix}_{q_j}
E_j^i \otimes E_j^{r - i} K_j^i \\ 
\hat{\Delta}(F_j^r) &= \sum_{i = 0}^r q_j^{i (r - i)} 
\begin{bmatrix}
r \\ 
i
\end{bmatrix}_{q_j}
F_j^{r - i} K_j^{-i} \otimes F_j^i
\end{align*}
for all $ r \in \mathbb{N} $ and $ j = 1, \dots, N $. 
\end{lemma} 

\begin{proof} We use induction on $ r $. The claim for $ \hat{\Delta}(E_j^1) = \hat{\Delta}(E_j) $ is clear. Assume the formula for $ \hat{\Delta}(E_j^r) $ 
is proved for $ r \in \mathbb{N} $. Using Lemma \ref{lemqcalc1}, we compute
\begin{align*}
\hat{\Delta}&(E_j^{r + 1}) = \hat{\Delta}(E_j) \hat{\Delta}(E_j^r) \\
&= (E_j \otimes K_j + 1 \otimes E_j) \sum_{i = 0}^r q_j^{i (r - i)} 
\begin{bmatrix}
r \\ 
i
\end{bmatrix}_{q_j}
E_j^i \otimes E_j^{r - i} K_j^i \\
&= \sum_{i = 0}^r  \left( q_j^{2(r - i)} q_j^{i (r - i)} 
\begin{bmatrix}
r \\ 
i
\end{bmatrix}_{q_j}
E_j^{i + 1} \otimes E_j^{r - i} K_j^{i + 1} + 
q_j^{i (r - i)}
\begin{bmatrix}
r \\ 
i
\end{bmatrix}_{q_j}
E_j^i \otimes E_j^{r + 1 - i} K_j^i \right) \allowdisplaybreaks \\
&= \sum_{i = 1}^{r + 1} q_j^{2(r + 1 - i)} q_j^{(i - 1) (r + 1 - i)} 
\begin{bmatrix}
r \\ 
i - 1
\end{bmatrix}_{q_j}
E_j^i \otimes E_j^{r + 1 - i} K_j^i \\
&\qquad + \sum_{i = 0}^r q_j^{i (r - i)} 
\begin{bmatrix}
r \\ 
i
\end{bmatrix}_{q_j}
E_j^i \otimes E_j^{r + 1 - i} K_j^i \\
&= \sum_{i = 0}^{r + 1} q_j^{i (r + 1 - i)} 
\begin{bmatrix}
r + 1 \\ 
i
\end{bmatrix}_{q_j}
E_j^i \otimes E_j^{r + 1 - i} K_j^i.
\end{align*} 
The formula for $ \hat{\Delta}(F_j^r) $ is proved in a similar way, or by applying the automorphism $ \omega $ from Lemma \ref{lemflippingEandF}
to the first formula. \end{proof} 

Let 
$$ 
\rho = \sum_{i = 1}^N \varpi_i = \frac{1}{2} \sum_{\alpha \in {\bf \Delta}^+} \alpha 
$$ 
\nomenclature{$\rho$}{half-sum of positive roots}%
be the half-sum of all positive roots. 
Note that for all $1\leq i \leq N$, we have
\begin{align*}
 (\rho,\alpha_i^\vee) &= 1, 
 & (\rho,\alpha_i) &= d_i.
\end{align*}

\begin{lemma}
	\label{lem:S2}
	For all $X\in U_q(\mathfrak{g})$ we have
	\[
	 \hat{S}^2(X) = K_{2\rho} X K_{-2\rho}.
	\]
\end{lemma}

\begin{proof}
  This is easily checked on generators. 
\end{proof}

In particular, the antipode $ \hat{S} $ is invertible, so that $ \tilde{U}_q(\mathfrak{g}) $ is a regular Hopf algebra, compare the remarks after Theorem \ref{thmcharmha}.
The 
inverse of $ \hat{S} $ is given on generators by 
\begin{align*}
\hat{S}^{-1}(K_\lambda) = K_{-\lambda}, \qquad \hat{S}^{-1}(E_j) = -K_j^{-1} E_j , \qquad \hat{S}^{-1}(F_j) = -F_j K_j. 
\end{align*}

Let $ \tilde{U}_q(\mathfrak{n}_+) $%
be the subalgebra of $ \tilde{U}_q(\mathfrak{g}) $ generated by the elements 
$ E_1, \dots, E_N $, and let $ \tilde{U}_q(\mathfrak{n}_-) $
be the subalgebra generated by $ F_1, \dots, F_N $. 
Moreover we let $ \tilde{U}_q(\mathfrak{h}) $
be the subalgebra generated by the elements $ K_\lambda $ 
for $ \lambda \in \weights $. We write $ \tilde{U}_q(\mathfrak{b}_+) $
for the subalgebra of $ \tilde{U}_q(\mathfrak{g}) $ generated by 
$ E_1, \dots, E_N $ and all $ K_\lambda $ for $ \lambda \in \weights $, and similarly we write $ \tilde{U}_q(\mathfrak{b}_-) $
for the 
subalgebra generated by the elements $ F_1, \dots, F_N, K_\lambda $ for $ \lambda \in \weights $. The algebras $ \tilde{U}_q(\mathfrak{h}) $ 
and $ \tilde{U}_q(\mathfrak{b}_\pm) $ are Hopf subalgebras. It is often convenient to use the 
automorphism $ \omega $ from Lemma \ref{lemflippingEandF} to transport results for $ \tilde{U}_q(\mathfrak{n}_+) $ 
and $ \tilde{U}_q(\mathfrak{b}_+) $ to $ \tilde{U}_q(\mathfrak{n}_-) $ and $ \tilde{U}_q(\mathfrak{b}_-) $, and vice versa.
\label{nom:Uqb_tilde} 

\begin{prop} \label{tildeuqgtriangular}
Multiplication in $ \tilde{U}_q(\mathfrak{g}) $ induces a linear isomorphism 
$$
\tilde{U}_q(\mathfrak{n}_-) \otimes \tilde{U}_q(\mathfrak{h}) \otimes \tilde{U}_q(\mathfrak{n}_+) \cong \tilde{U}_q(\mathfrak{g}). 
$$
\end{prop} 

\begin{proof} In order to prove the claim it suffices to show that the elements $ F_I K_\mu E_J $, 
for $ F_I = F_{i_1} \cdots F_{i_k}, E_J = E_{j_1} \cdots E_{j_l} $ finite sequences of simple root vectors and 
$ \mu \in \weights $, form a linear basis of $ \tilde{U}_q(\mathfrak{g}) $. This in turn is 
an easy consequence of the Diamond Lemma \cite{Bergmandiamond}. Indeed, from the definition of $ \tilde{U}_q(\mathfrak{g}) $ we see 
that there are only overlap ambiguities, and all of these turn out to be resolvable. \end{proof} 

In a similar way we obtain linear isomorphisms 
$$
\tilde{U}_q(\mathfrak{n}_\pm) \otimes \tilde{U}_q(\mathfrak{h}) \cong \tilde{U}_q(\mathfrak{b}_\pm); 
$$
moreover $ \tilde{U}_q(\mathfrak{h}) $ identifies with the group algebra of $ \weights $, whereas $ \tilde{U}_q(\mathfrak{n}_\pm) $ 
is isomorphic to the free algebra generated by the $ E_i $ or $ F_j $, respectively. 

The algebra $ \tilde{U}_q(\mathfrak{g}) $ is equipped with a $ \roots $-grading such that the generators $ K_\lambda $ have degree $ 0 $, 
the elements $ E_i $ degree $ \alpha_i $, and the elements $ F_i $ degree $ -\alpha_i $. We shall also refer to the degree of a homogeneous element with 
respect to this grading as its weight, and denote by $ \tilde{U}_q(\mathfrak{g})_\beta \subset \tilde{U}_q(\mathfrak{g}) $ the subspace of all elements 
of weight $ \beta $. The weight grading induces a direct sum decomposition 
$$
\tilde{U}_q(\mathfrak{g}) = \bigoplus_{\beta \in \roots} \tilde{U}_q(\mathfrak{g})_\beta 
$$
of $ \tilde{U}_q(\mathfrak{g}) $. 
Notice that $ \tilde{U}_q(\mathfrak{g})_\alpha \tilde{U}_q(\mathfrak{g})_\beta \subset \tilde{U}_q(\mathfrak{g})_{\alpha + \beta} $ for 
all $ \alpha, \beta \in \roots $. 

If $ q \in \mathbb{K}^\times $ is not a root of unity we can describe the weight grading equivalently in terms of the adjoint action. 
More precisely, we have that $ X \in \tilde{U}_q(\mathfrak{g}) $ is of weight $ \lambda $ iff 
$$ 
\ad(K_\mu)(X) = K_\mu X K_\mu^{-1} = q^{(\mu, \lambda)} X 
$$ 
for all $ \mu \in \weights $. 

As for any Hopf algebra, we may consider the adjoint action of $ \tilde{U}_q(\mathfrak{g}) $ on itself, given by 
$$
\ad(X)(Y) = X \rightarrow Y = X_{(1)} Y \hat{S}(X_{(2)}). 
$$
\label{nom:left_adjoint_action1}
Explicitly, we obtain  
\begin{align*} 
E_j \rightarrow Y &= Y \hat{S}(E_j) + E_j Y \hat{S}(K_j) = - Y E_j K_j^{-1} + E_j Y K_j^{-1} = [E_j, Y] K_j^{-1} \\ 
F_j \rightarrow Y &= K_j^{-1} Y \hat{S}(F_j) + F_j Y = - K_j^{-1} Y K_j F_j + F_j Y \\ 
K_\lambda \rightarrow Y &= K_\lambda Y \hat{S}(K_\lambda) = K_\lambda Y K_\lambda^{-1}
\end{align*} 
for all $ 1 \leq j \leq N $ and $ \lambda \in \weights $. 

We shall occasionally also consider the right adjoint action of $ \tilde{U}_q(\mathfrak{g}) $ on itself, given by 
$$
Y \leftarrow X = \hat{S}(X_{(1)}) Y X_{(2)}. 
$$
\label{nom:right_adjoint_action1}
These actions are linked via 
$$
\omega(X \rightarrow Y) 
= \omega(Y) \leftarrow \hat{S}^{-1}(\omega(X)), 
$$
where $ \omega $ is the automorphism from Lemma \ref{lemflippingEandF}. 

\begin{lemma} \label{adjointpowers}
We have 
\begin{align*} 
E_j^r \rightarrow Y &= \sum_{i = 0}^r (-1)^{r - i} q_j^{-(r - i)(r - 1)} 
\begin{bmatrix}
r \\ 
i
\end{bmatrix}_{q_j}
E_j^i Y E_j^{r - i} K_j^{-r} \\ 
F_j^r \rightarrow Y &= \sum_{i = 0}^r (-1)^i q_j^{i (r - 1)} 
\begin{bmatrix}
r \\ 
i
\end{bmatrix}_{q_j}
F_j^{r - i} K_j^{-i} Y K_j^i F_j^i
\end{align*} 
for all $ r \in \mathbb{N} $. 
\end{lemma} 

\begin{proof} We apply Lemma \ref{lemEiFicomult} to obtain 
\begin{align*}
E_j^r \rightarrow Y &= \sum_{i = 0}^r q_j^{i (r - i)} 
\begin{bmatrix}
r \\ 
i
\end{bmatrix}_{q_j}
E_j^i Y \hat{S}(E_j^{r - i} K_j^i) \\ 
&= \sum_{i = 0}^r q_j^{i (r - i)} 
\begin{bmatrix}
r \\ 
i
\end{bmatrix}_{q_j}
(-1)^{r - i} E_j^i Y K_j^{-i} (E_j K_j^{-1})^{r - i} \\ 
&= \sum_{i = 0}^r q_j^{i (r - i)}  
\begin{bmatrix}
r \\ 
i
\end{bmatrix}_{q_j}
(-1)^{r - i} q_j^{-(r - i)(r - i - 1) - 2i(r - i)} E_j^i Y E_j^{r - i} K_j^{-r} \\
&= \sum_{i = 0}^r (-1)^{r - i} q_j^{-(r - i)(r - 1)} 
\begin{bmatrix}
r \\ 
i
\end{bmatrix}_{q_j}
E_j^i Y E_j^{r - i} K_j^{-r}. 
\end{align*}
Similarly, 
\begin{align*}
F_j^r \rightarrow Y &= \sum_{i = 0}^r q_j^{i (r - i)} 
\begin{bmatrix}
r \\ 
i
\end{bmatrix}_{q_j}
F_j^{r - i} K_j^{-i} Y \hat{S}(F_j)^i \\ 
&= \sum_{i = 0}^r q_j^{i (r - i)} 
\begin{bmatrix}
r \\ 
i
\end{bmatrix}_{q_j}
(-1)^i q_j^{(i - 1) i} 
F_j^{r - i} K_j^{-i} Y K_j^i F_j^i \\ 
&= \sum_{i = 0}^r (-1)^i q_j^{i (r - 1)} 
\begin{bmatrix}
r \\ 
i
\end{bmatrix}_{q_j}
F_j^{r - i} K_j^{-i} Y K_j^i F_j^i. 
\end{align*}
This yields the claim. \end{proof}

\subsubsection{The Serre elements}
\label{sec:Serre_elements}

For the construction of $ U_q(\mathfrak{g}) $ we shall be interested in 
the \emph{Serre elements} of $ \tilde{U}_q(\mathfrak{g}) $, defined by 
\begin{align*}
u^+_{ij} &= \sum_{k = 0}^{1 - a_{ij}} (-1)^k \begin{bmatrix} 1 - a_{ij} \\ k \end{bmatrix}_{q_i}
E_i^{1 - a_{ij} - k} E_j E_i^k \\ 
u^-_{ij} &= \sum_{k = 0}^{1 - a_{ij}} (-1)^k \begin{bmatrix} 1 - a_{ij} \\ k \end{bmatrix}_{q_i}
F_i^{1 - a_{ij} - k} F_j F_i^k, 
\end{align*}
\nomenclature{$u^\pm_{ij}$}{Serre elements}%
for $ 1 \leq i,j \leq N $ with $ i \neq j $. We note that 
$$
u^-_{ij} = F_i^{1 - a_{ij}} \rightarrow F_j, \qquad u^+_{ij} = E_j \leftarrow \hat{S}^{-1}(E_i^{1 - a_{ij}})
$$
due to Lemma \ref{adjointpowers}.  Note also that $ \omega(u^+_{ij}) = u^-_{ij} $.

\begin{prop} \label{serrecomultiplication}
We have 
\begin{align*}
\hat{\Delta}(u^+_{ij}) &= u^+_{ij} \otimes K_i^{1 - a_{ij}} K_j + 1 \otimes u^+_{ij} \\ 
\hat{\Delta}(u^-_{ij}) &= u^-_{ij} \otimes 1 + K_i^{-(1 - a_{ij})} K_j^{-1}  \otimes u^-_{ij}
\end{align*}
and 
\begin{align*}
\hat{S}(u^+_{ij}) &= -u^+_{ij} K_i^{a_{ij} - 1} K_j^{-1} \\
\hat{S}(u^-_{ij}) &= -K_i^{1 - a_{ij}} K_j u^-_{ij}
\end{align*}
for all $ 1 \leq i,j \leq N $ with $ i \neq j $. 
\end{prop} 

\begin{proof} Let us prove the first claim for $ u^+_{ij} $, following Lemma 4.10 in \cite{Jantzenqg}. We note that 
$$
\hat{\Delta}(u^+_{ij}) = \sum_{k = 0}^{1 - a_{ij}} (-1)^{1 - a_{ij} - k} \begin{bmatrix} 1 - a_{ij} \\ k \end{bmatrix}_{q_i}
\hat{\Delta}(E_i^k) \hat{\Delta}(E_j) \hat{\Delta}(E_i^{1 - a_{ij} - k}), 
$$
and by Lemma \ref{lemEiFicomult} we have 
$$
\hat{\Delta}(E_i^r) = E_i^r \otimes K_i^r + 1 \otimes E_i^r + \sum_{m = 1}^{r - 1} q_i^{m (r - m)} 
\begin{bmatrix}
r \\ 
m
\end{bmatrix}_{q_i}
E_i^m \otimes E_i^{r - m} K_i^m.  
$$
It follows that $ \hat{\Delta}(u^+_{ij}) $ has the form 
$$
\hat{\Delta}(u^+_{ij}) = u^+_{ij} \otimes K_i^{1 - a_{ij}} K_j + 1 \otimes u^+_{ij} + 
\sum_{m = 1}^{1 - a_{ij}} E^m_i \otimes X_m + \sum_{m,n} E_i^m E_j E_i^n \otimes Y_{mn} 
$$
for suitable elements $ X_k, Y_{kl} $, where the second sum is over all $ m, n \geq 0 $ such that $ m + n \leq 1 - a_{ij} $. 

We shall show that all $ X_k $ and $ Y_{kl} $ are zero. Collecting terms in the binomial expansions we obtain 
\begin{align*}
Y_{mn} &= \sum_{k = m}^{1 - a_{ij} - n} (-1)^{1 - a_{ij} - k} \begin{bmatrix} 1 - a_{ij} \\ k \end{bmatrix}_{q_i} q_i^{m(k - m)} \\
&\qquad \times 
\begin{bmatrix} k \\ m \end{bmatrix}_{q_i} E_i^{k - m} K_i^m K_j q_i^{n(1 - a_{ij} - k - n)} 
\begin{bmatrix} 1 - a_{ij} - k \\ n \end{bmatrix}_{q_i} E_i^{1 - a_{ij} - k - n} K_i^n. 
\end{align*}

We have 
\begin{align*}
\begin{bmatrix} r \\ k \end{bmatrix}_{q_i}
\begin{bmatrix} k \\ m \end{bmatrix}_{q_i} 
\begin{bmatrix} r - k \\ n \end{bmatrix}_{q_i}
&= \frac{[r]! [k]! [r - k]!}{[k]! [m]! [n]! [r - k]! [k - m]! [r - k - n]!} \\
&= \frac{[r]!}{[m]! [n]! [k - m]! [r - k - n]!} \\
&= \frac{[r - m - n]! [r - n]! [r]!}{[k - m]! [m]! [n]! [r - k - n]! [r - m - n]! [r - n]!} \\
&= \begin{bmatrix} r - m - n \\ k - m \end{bmatrix}_{q_i}
\begin{bmatrix} r - n \\ m \end{bmatrix}_{q_i}
\begin{bmatrix} r \\ n \end{bmatrix}_{q_i}. 
\end{align*}
This shows
$$
Y_{mn} = (-1)^{1 - a_{ij}} y_{mn} \begin{bmatrix} 1 - a_{ij} - n \\ m \end{bmatrix}_{q_i} 
\begin{bmatrix} 1 - a_{ij} \\ n \end{bmatrix}_{q_i} E_i^{1 - a_{ij} - m - n} K_i^{m + n} K_j 
$$
where 
\begin{align*}
y_{mn} &= \sum_{k = m}^{1 - a_{ij} - n} (-1)^k \begin{bmatrix} 1 - a_{ij} - m - n \\ k  - m \end{bmatrix}_{q_i} 
q_i^{m(k - m) + (n + 2m + a_{ij})(1 - a_{ij} - k - n)} \\
&= (-1)^m q_i^{(n + 2m + a_{ij})(1 - a_{ij} - m - n)} \\
&\qquad \times \sum_{l = 0}^{1 - a_{ij} - m - n} (-1)^l \begin{bmatrix} 1 - a_{ij} - m - n \\ l \end{bmatrix}_{q_i} 
q_i^{(1 - a_{ij} - m - n)l} q_i^{-l} 
= 0 
\end{align*}
according to Lemma \ref{lembinomialvanishing}. 

For $ X_m $ we get 
\begin{align*}
X_m &= \sum_{k = 0}^{1 - a_{ij}} (-1)^{1 - a_{ij} - k} \begin{bmatrix} 1 - a_{ij} \\ k \end{bmatrix}_{q_i} 
\sum_l q_i^{l(k - l)} \begin{bmatrix} k \\ l \end{bmatrix}_{q_i} E_i^{k - l} K_i^l E_j \\
&\qquad \times q_i^{(m - l)(1 - a_{ij} - k - m + l)} \begin{bmatrix} 1 - a_{ij} - k  \\ m - l \end{bmatrix}_{q_i} 
E_i^{1 - a_{ij} - k - m + l} K_i^{m - l}, 
\end{align*}
where $ l $ runs over $ \max(0, m + k + a_{ij} - 1) \leq l \leq \min(k, m) $. 
We get 
\begin{align*}
X_m &= (-1)^{1 - a_{ij}} \sum_{s = 0}^{1 - a_{ij} - m} a_{ms} E_i^s E_j E_i^{1 - a_{ij} - m - s} K_i^m, 
\end{align*}
where 
\begin{align*}
a_{ms} &= \sum_{k = s}^{m + s} (-1)^k \begin{bmatrix} 1 - a_{ij} \\ k \end{bmatrix}_{q_i} 
\begin{bmatrix} k \\ s \end{bmatrix}_{q_i} 
\begin{bmatrix} 1 - a_{ij} - k \\ m + s - k \end{bmatrix}_{q_i} \\
&\qquad \times q_i^{(k - s)s + (m + s - k)(1 - a_{ij} - m - s) + a_{ij} (k - s) + 2(k - s)(1 - a_{ij} - m - s)}. 
\end{align*}
Using 
\begin{align*}
\begin{bmatrix} r \\ k \end{bmatrix}_{q_i}
\begin{bmatrix} k \\ s \end{bmatrix}_{q_i} 
\begin{bmatrix} r - k \\ m + s - k \end{bmatrix}_{q_i}
&= \frac{[r]! [k]! [r - k]!}{[k]! [s]! [m + s - k]! [r - k]! [k - s]! [r - m - s]!} \\
&= \frac{[r]!}{[s]![m + s - k]! [k - s]! [r - m - s]!} \\
&= \frac{[r]! [r - m]! [m]!}{[m]! [s]! [k - s]! [r - m]! [r - m - s]! [m + s - k]!} \\
&= \begin{bmatrix} r \\ m \end{bmatrix}_{q_i}
\begin{bmatrix} r - m \\ s \end{bmatrix}_{q_i}
\begin{bmatrix} m \\ k - s \end{bmatrix}_{q_i},
\end{align*}
this simplifies to 
\begin{align*}
a_{ms} &= \sum_{k = s}^{m + s} (-1)^k \begin{bmatrix} 1 - a_{ij} \\ m \end{bmatrix}_{q_i} 
\begin{bmatrix} 1 - a_{ij} - m \\ s \end{bmatrix}_{q_i} 
\begin{bmatrix} m \\ k - s \end{bmatrix}_{q_i} \\
&\qquad \times q_i^{(k - s)s + (m + s - k)(1 - a_{ij} - m - s) + a_{ij} (k - s) + 2(k - s)(1 - a_{ij} - m - s)} \\
&= \sum_{k = s}^{m + s} (-1)^k \begin{bmatrix} 1 - a_{ij} \\ m \end{bmatrix}_{q_i} 
\begin{bmatrix} 1 - a_{ij} - m  \\ s \end{bmatrix}_{q_i} 
\begin{bmatrix} m \\ k - s \end{bmatrix}_{q_i} \\
&\qquad \times q_i^{(k - s)s + (m + k - s)(1 - a_{ij} - m - s) + a_{ij} (k - s)} \\
&= \begin{bmatrix} 1 - a_{ij} \\ m \end{bmatrix}_{q_i} 
\begin{bmatrix} 1 - a_{ij} - m  \\ s \end{bmatrix}_{q_i} 
\sum_{k = s}^{m + s} (-1)^k \begin{bmatrix} m \\ k - s \end{bmatrix}_{q_i} 
q_i^{m (1 - a_{ij} - m - s) + (k - s)(1 - m)} \\
&= (-1)^s q_i^{m (1 - a_{ij} - m - s)} \begin{bmatrix} 1 - a_{ij} \\ m \end{bmatrix}_{q_i} 
\begin{bmatrix} 1 - a_{ij} - m  \\ s \end{bmatrix}_{q_i} 
\sum_{l = 0}^m (-1)^l \begin{bmatrix} m \\ l \end{bmatrix}_{q_i} 
q_i^{l (1 - m)} = 0, 
\end{align*}
taking into account Lemma \ref{lembinomialvanishing}. 

The claim for the coproduct of $ u^-_{ij} $ is obtained by applying the automorphism $ \omega $ from Lemma \ref{lemflippingEandF}. 

The remaining assertions follow from the formulas for $ \hat{\Delta}(u^\pm_{ij}) $ and the antipode relations. \end{proof}

\begin{lemma} \label{EFpowercommutationhelp} 
For any $ 1 \leq k \leq N $ and $ r \in \mathbb{N} $ we have 
$$
[E_k, F_k^r] = [r]_{q_k} F_k^{r - 1} \frac{1}{q_k - q_k^{-1}} (q_k^{-(r - 1)} K_k  - q_k^{r - 1}K_k^{-1}). 
$$
\end{lemma} 

\begin{proof} For $ r = 1 $ the formula clearly holds. Assuming the claim for $ r $, we compute 
\begin{align*}
[E_k, &F_k^{r + 1}] = [E_k, F_k^r] F_k + F_k^r [E_k , F_k] \\
&= [r]_{q_k} F_k^{r - 1} \frac{1}{q_k - q_k^{-1}} (q_k^{-(r - 1)} K_k  - q_k^{r - 1}K_k^{-1}) F_k 
+ F_k^r \frac{1}{q_k - q_k^{-1}}(K_k - K_k^{-1}) \\ 
&= \frac{1}{q_k - q_k^{-1}} [r]_{q_k} F_k^r (q_k^{-(r + 1)} K_k  - q_k^{r + 1} K_k^{-1}) 
+ F_k^r \frac{1}{q_k - q_k^{-1}}(K_k - K_k^{-1}) \\ 
&= \frac{1}{q_k - q_k^{-1}} (q_k^{r - 1} + q_k^{r - 3} + \cdots + q_k^{-r + 1}) F_k^r (q_k^{-(r + 1)} K_k  - q_k^{r + 1} K_k^{-1}) \\
&\qquad + F_k^r \frac{1}{q_k - q_k^{-1}}(K_k - K_k^{-1}) \\ 
&= (q_k^{r} + q_k^{r - 2} + \cdots + q_k^{-r}) F_k^r \frac{1}{q_k - q_k^{-1}} (q_k^{-r} K_k  - q_k^{r} K_k^{-1}) \\
&= [r + 1]_{q_k} F_k^r \frac{1}{q_k - q_k^{-1}} (q_k^{-r} K_k  - q_k^{r} K_k^{-1}).  
\end{align*}
This finishes the proof. \end{proof} 

\begin{lemma} \label{lemucommutator}
For $ 1 \leq i,j \leq N $ and $ i \neq j $ we have 
$$
[E_k, u^-_{ij}] = 0 = [F_k, u^+_{ij}] 
$$
for all $ k = 1, \dots, N $. 
\end{lemma} 

\begin{proof} Let us verify $ [E_k, u^-_{ij}] = 0 $. If $ k \neq i,j  $ then the elements $ E_k $ commute with $ F_i, F_j $, 
so that the assertion is obvious. Assume now $ k = i $. Note that 
$$
[E_k, u^-_{ij}] = \ad(E_k)(u^-_{ij}) K_k,  
$$
so it suffices to show $ \ad(E_k)(u^-_{kj}) = 0 $. As observed at the beginning of this subsection, we have
$ \ad(F_k^{1 - a_{kj}})(F_j) = u^-_{kj} $. Combining this relation with Lemma \ref{EFpowercommutationhelp} we obtain 
\begin{align*}
\ad(E_k)(u^-_{kj}) &= \ad(E_k) \ad(F_k^{1 - a_{kj}})(F_j) \\
&= \ad(F_k^{1 - a_{kj}}) \ad(E_k)(F_j) + \ad([E_k, F_k^{1 - a_{kj}}])(F_j) \\
&= \ad(F_k^{1 - a_{kj}}) \ad(E_k)(F_j) \\
&\qquad + [1 - a_{kj}]_{q_k} \ad(F_k^{- a_{kj}}) {\textstyle \frac{1}{q_k - q_k^{-1}}} (q_k^{a_{kj}} \ad(K_k) - q_k^{-a_{kj}} \ad(K_k^{-1}))(F_j) \\
&= 0,  
\end{align*} 
by using $ [E_k, F_j] = 0 $ and $ \ad(K_k)(F_j) = q_k^{-a_{kj}} F_j $. 

Now let $ k = j $. Then $ [E_k, F_i] = 0 $ 
and $ K_k F_i 
= q_i^{a_{ik}} F_i K_k $,  
so that 
\begin{align*}
[E_k, u^-_{ik}] &= \sum_{l = 0}^{1 - a_{ik}} (-1)^l \begin{bmatrix} 1 - a_{ik} \\ l \end{bmatrix}_{q_i}
F_i^{1 - a_{ik} - l} [E_k, F_k] F_i^l \\
&= \frac{1}{q_k - q_k^{-1}} \sum_{l = 0}^{1 - a_{ik}} (-1)^l \begin{bmatrix} 1 - a_{ik} \\ l \end{bmatrix}_{q_i}
F_i^{1 - a_{ik} - l} (K_k - K_k^{-1}) F_i^l \\
&= \frac{1}{q_k - q_k^{-1}} \sum_{l = 0}^{1 - a_{ik}} (-1)^l \begin{bmatrix} 1 - a_{ik} \\ l \end{bmatrix}_{q_i} F_i^{1 - a_{ik}}
(q_i^{l a_{ik}} K_k - q_i^{l (-a_{ik})} K_k^{-1}) = 0 
\end{align*}
by Lemma \ref{lembinomialvanishing}. 

Finally, the equality $ [F_k, u^+_{ij}] = 0 $ follows by applying the automorphism $ \omega $ 
from Lemma \ref{lemflippingEandF}. \end{proof}

\subsubsection{The quantized universal enveloping algebra}

Let us now give the definition of the quantized universal enveloping algebra of $ \mathfrak{g} $. 
We fix a ground field $ \mathbb{K} $. 

\begin{definition} \label{defuqg}
 \nomenclature{$U_q(\lie{g})$}{quantized enveloping algebra}%
Let $ \mathbb{K} $ be a field and $ q = s^L \in \mathbb{K}^\times $ be an invertible element such that $ q_i \neq \pm 1 $ for all $ i $. 
The algebra $ U_q(\mathfrak{g}) $ over $ \mathbb{K} $ has generators $ K_\lambda $ 
for $ \lambda \in \weights $, and $ E_i, F_i $ for $ i = 1, \dots, N $, and the defining relations for $ U_q(\mathfrak{g}) $ are 
\begin{align*}
K_0 &= 1 \\
K_\lambda K_\mu &= K_{\lambda + \mu} \\
K_\lambda E_j K_\lambda^{-1} &= q^{(\lambda, \alpha_j)} E_j \\
K_\lambda F_j K_\lambda^{-1} &= q^{-(\lambda, \alpha_j)} F_j \\ 
[E_i, F_j] &= \delta_{ij} \frac{K_i - K_i^{-1}}{q_i - q_i^{-1}} 
\end{align*} 
for all $ \lambda, \mu \in \weights $ and all $ i, j $, together with the quantum Serre relations 
\begin{align*}
&\sum_{k = 0}^{1 - a_{ij}} (-1)^k \begin{bmatrix} 1 - a_{ij} \\ k \end{bmatrix}_{q_i}
E_i^{1 - a_{ij} - k} E_j E_i^k = 0 \\ 
&\sum_{k = 0}^{1 - a_{ij}} (-1)^k \begin{bmatrix} 1 - a_{ij} \\ k \end{bmatrix}_{q_i}
F_i^{1 - a_{ij} - k} F_j F_i^k = 0. 
\end{align*}
In the above formulas we abbreviate $ K_i = K_{\alpha_i} $ for all simple roots, and we use the notation $ q_i = q^{d_i} $. 
\end{definition} 

We occasionally write $E_\alpha$ and $F_\alpha$ to denote $E_i$ and $F_i$ when $\alpha=\alpha_i$ is a simple root.

One often finds a slightly different version of the quantized universal enveloping algebra in the literature, 
only containing elements $ K_\lambda $ for $ \lambda \in \roots $. For our purposes it will be crucial to work with the algebra 
as in Definition \ref{defuqg}.

It follows from Lemma \ref{lemuqhopf} and Proposition \ref{serrecomultiplication} that $ U_q(\mathfrak{g}) $ is a Hopf algebra 
with comultiplication $ \hat{\Delta}: U_q(\mathfrak{g}) \rightarrow U_q(\mathfrak{g}) \otimes U_q(\mathfrak{g}) $ given by 
\label{nom:Uqg_Hopf}
\begin{align*}
\hat{\Delta}(K_\lambda) &= K_\lambda \otimes K_\lambda, \\
\hat{\Delta}(E_i) &= E_i \otimes K_i + 1 \otimes E_i \\
\hat{\Delta}(F_i) &= F_i \otimes 1 + K_i^{-1} \otimes F_i,  
\end{align*}
counit $ \hat{\epsilon}: U_q(\mathfrak{g}) \rightarrow \mathbb{K} $ given by 
\begin{align*}
\hat{\epsilon}(K_\lambda) = 1, \qquad \hat{\epsilon}(E_j) = 0, \qquad \hat{\epsilon}(F_j) = 0, 
\end{align*} 
and antipode $ \hat{S}: U_q(\mathfrak{g}) \rightarrow U_q(\mathfrak{g}) $ given by 
\begin{align*}
\hat{S}(K_\lambda) = K_{-\lambda}, \qquad \hat{S}(E_j) = -E_j K_j^{-1}, \qquad \hat{S}(F_j) = -K_j F_j. 
\end{align*}

Let $ U_q(\mathfrak{n}_+) $
\nomenclature{$U_q(\mathfrak{n}_+)$, $U_q(\mathfrak{n}_-)$}{upper and lower triangular subalgebras of $U_q(\lie{g})$}%
be the subalgebra of $ U_q(\mathfrak{g}) $ generated by the elements 
$ E_1, \dots, E_N $, and let $ U_q(\mathfrak{n}_-) $ be the subalgebra generated by $ F_1, \dots, F_N $. 
Moreover we let $ U_q(\mathfrak{h}) $
\nomenclature{$U_q(\lie{h})$}{Cartan subalgebra of $U_q(\lie{g})$}%
 be the subalgebra generated by the elements $ K_\lambda $ 
for $ \lambda \in \weights $.  

\begin{prop} \label{uqgtriangular}
Multiplication in $ U_q(\mathfrak{g}) $ induces a linear isomorphism 
$$
U_q(\mathfrak{n}_-) \otimes U_q(\mathfrak{h}) \otimes U_q(\mathfrak{n}_+) \cong U_q(\mathfrak{g}), 
$$
the quantum analogue of the triangular decomposition. 
\end{prop} 

\begin{proof} Let us write $ I^+ $ for the ideal in $ \tilde{U}_q(\mathfrak{n}_+) $ 
generated by the elements $ u^+_{ij} $. We claim that under the isomorphism of Proposition \ref{tildeuqgtriangular} 
the ideal in $ \tilde{U}_q(\mathfrak{g}) $ generated by the $ u^+_{ij} $ identifies with 
$ \tilde{U}_q(\mathfrak{n}_-) \otimes \tilde{U}_q(\mathfrak{h}) \otimes I^+ $. The left ideal property is obvious, 
and the right ideal property follows using Lemma \ref{lemucommutator}. 
An analogous claim holds for the ideal $ I^- $ in  $ \tilde{U}_q(\mathfrak{n}_-) $ 
generated by the elements $ u^-_{ij} $. 

It follows that 
$$
I = I^- \otimes \tilde{U}_q(\mathfrak{h}) \otimes \tilde{U}_q(\mathfrak{n}_+) + 
\tilde{U}_q(\mathfrak{n}_-) \otimes \tilde{U}_q(\mathfrak{h}) \otimes I^+
$$ 
is an ideal in $ \tilde{U}_q(\mathfrak{g}) $, and the quotient algebra $ \tilde{U}_q(\mathfrak{g})/I $ is canonically 
isomorphic to $ U_q(\mathfrak{g}) $. The assertion now follows from Proposition \ref{tildeuqgtriangular}. Note in particular that 
we obtain a canonical isomorphism $ U_q(\mathfrak{h}) \cong \tilde{U}_q(\mathfrak{h}) $. \end{proof}

We write $ U_q(\mathfrak{b}_+) $
\nomenclature{$U_q(\mathfrak{b}_+)$, $U_q(\mathfrak{b}_-)$}{Borel subalgebras of $U_q(\lie{g})$}%
for the subalgebra of $ U_q(\mathfrak{g}) $ generated by 
$ E_1, \dots, E_N $ and all $ K_\lambda $ for $ \lambda \in \weights $, and similarly we write $ U_q(\mathfrak{b}_-) $ for the 
subalgebra generated by the elements $ F_1, \dots, F_N, K_\lambda $ for $ \lambda \in \weights $. These algebras are  Hopf subalgebras. It follows from Proposition \ref{uqgtriangular} that these algebras, as well 
as the algebras $ U_q(\mathfrak{n}_\pm) $ and $ U_q(\mathfrak{h}) $, are canonically isomorphic to 
the universal algebras with the appropriate generators and relations from Definition \ref{defuqg}. 

As a consequence of Proposition \ref{uqgtriangular} we obtain linear isomorphisms  
$$ 
U_q(\mathfrak{h}) \otimes U_q(\mathfrak{n}_\pm) \cong U_q(\mathfrak{b}_\pm).
$$ 
Also note that $ U_q(\mathfrak{h}) $ is isomorphic to the group algebra of $ \weights $. 

We observe that the automorphism $ \omega: \tilde{U}_q(\mathfrak{g}) \rightarrow \tilde{U}_q(\mathfrak{g}) $ from 
Lemma \ref{lemflippingEandF} induces an algebra automorphism of $ U_q(\mathfrak{g}) $, which we will again denote by $ \omega $. 
This automorphism allows us to interchange the upper and lower triangular parts of $ U_q(\mathfrak{g}) $. 
We state this explicitly in the following Lemma.

\begin{lemma} \label{defOmega}
There exists an algebra automorphism $ \omega: U_q(\mathfrak{g}) \rightarrow U_q(\mathfrak{g}) $ given by 
$$
\omega(K_\lambda) = K_{-\lambda}, \qquad \omega(E_j) = F_j, \qquad \omega(F_j) = E_j 
$$
on generators. Moreover $ \omega $ is a coalgebra antihomomorphism. 
\\
Similarly, there exists an algebra anti-automorphism $ \Omega: U_q(\mathfrak{g}) \rightarrow U_q(\mathfrak{g}) $ given by 
$$
\Omega(K_\lambda) = K_\lambda, \qquad \Omega(E_j) = F_j, \qquad \Omega(F_j) = E_j 
$$
on generators. 
\end{lemma} 

\begin{proof} It is straightforward to verify that the defining relations of $ U_q(\mathfrak{g}) $ are preserved by these assignments. \end{proof} 

We note that the automorphism $ \omega $ can be used to translate formulas which depend on our convention for the comultiplication of $ U_q(\mathfrak{g}) $ 
to the convention using the coopposite comultiplication. 

Let us introduce another symmetry of $ U_q(\mathfrak{g}) $. 

\begin{lemma} \label{deftau} 
There is a unique algebra anti-automorphism $ \tau: U_q(\mathfrak{g}) \rightarrow U_q(\mathfrak{g}) $ such that 
$$
\tau(E_j) = K_j F_j, \qquad \tau(F_j) = E_j K_j^{-1}, \qquad \tau(K_\lambda) = K_\lambda 
$$
\nomenclature{$\tau$}{algebra anti-automorphism, coalgebra automorphism of $U_q(\lie{g})$}%
for $ j = 1, \dots, N $ and $ \lambda \in \weights $. 
Moreover $ \tau $ is involutive, that is $ \tau^2 = \id $, and a coalgebra homomorphism. 
\end{lemma}

\begin{proof}
 Observe first that all of the Hopf algebra relations of $U_q(\mathfrak{g})$ are preserved by the transformation
 \[
  E_i \mapsto -E_i, \qquad F_i\mapsto -F_i, \qquad K_\lambda \mapsto K_\lambda,
 \]
 so that these generate an involutive Hopf algebra automorphism of $U_q(\mathfrak{g})$.  Now $\tau$ agrees on generators with the composition of this map with the involution $\hat{S}\circ\omega$, which is an algebra anti-automorphism and coalgebra automorphism.  The result follows.
\end{proof}

In the same way as for $ \tilde{U}_q(\mathfrak{g}) $ one obtains a $ \roots $-grading on $ U_q(\mathfrak{g}) $ such that the 
generators $ K_\lambda $ have degree $ 0 $, the elements $ E_i $ degree $ \alpha_i $, and the elements $ F_i $ degree $ -\alpha_i $. 
The degree of a homogeneous element with respect to this grading will again be referred to as its weight, and 
we write $ U_q(\mathfrak{g})_\beta \subset U_q(\mathfrak{g}) $
\nomenclature{$U_q(\lie{g})_\beta$}{subspace of $U_q(\lie{g})$ of weight $\beta\in\roots$}%
for the subspace of all elements 
of weight $ \beta $. The weight grading induces a direct sum decomposition 
$$
U_q(\mathfrak{g}) = \bigoplus_{\beta \in \roots} U_q(\mathfrak{g})_\beta 
$$
of $ U_q(\mathfrak{g}) $, and we have $ U_q(\mathfrak{g})_\alpha U_q(\mathfrak{g})_\beta \subset U_q(\mathfrak{g})_{\alpha + \beta} $ for 
all $ \alpha, \beta \in \roots $. 

If $ q \in \mathbb{K}^\times $ is not a root of unity the weight grading is determined by the adjoint action of $ U_q(\mathfrak{g}) $ 
on itself in the same way as for $ \tilde{U}_q(\mathfrak{g}) $.

\subsubsection{The restricted integral form}
\label{sec:restricted_integral_form}

Instead of working over a field $ \mathbb{K} $, it is sometimes necessary to consider more general coefficient rings.  

Let us fix a semisimple Lie algebra $ \mathfrak{g} $ and put $ q = s^L \in \mathbb{Q}(s) $. 
We consider the ring $ \A = \mathbb{Z}[s,s^{-1}] \subset \mathbb{Q}(s) $.
\nomenclature[o$A$]{$\A$}{$=\mathbb{Z}[s,s^{-1}]$}%
With this notation in place, we define the restricted integral form of the quantized universal enveloping algebra, compare \cite{Lusztigrootsofone}.

\begin{definition} 
\label{def:integral_form}
The restricted integral form $ U_q^\A(\mathfrak{g}) $
\nomenclature{$U_q^\A(\mathfrak{g})$}{restricted integral form of $U_q(\lie{g})$}
of $ U_q(\mathfrak{g}) $ is the $\A$-subalgebra of the quantized universal enveloping algebra $ U_q(\mathfrak{g}) $ 
over $ \mathbb{Q}(s) $ generated by the elements $ K_\lambda $ for $ \lambda \in \weights $ and
$$ 
[K_i; 0] = \frac{K_i - K_i^{-1}}{q_i - q_i^{-1}} 
$$ 
\label{nom:K:0}%
for $ i = 1, \dots, N $, together with the divided powers 
$$ 
E_i^{(r)} = \frac{1}{[r]_{q_i}!} E_i^r, \qquad F_i^{(r)} = \frac{1}{[r]_{q_i}!} F_i^r
$$
\nomenclature{$E_i^{(r)}$}{divided power $\frac{1}{[r]_{q_i}"!} E_i^r$}%
\nomenclature{$F_i^{(r)}$}{divided power $\frac{1}{[r]_{q_i}"!} F_i^r$}%
for $ i = 1, \dots, N $. 
\end{definition}

Let us show that $ U_q^\A(\mathfrak{g}) $ is a Hopf algebra over the commutative ring $ \A $ in a natural way. 

\begin{lemma} \label{integralhopfstructure}
The comultiplication, counit and antipode of $ U_q(\mathfrak{g}) $ induce on $ U_q^\A(\mathfrak{g}) $ the structure of a Hopf algebra over $ \A $. 
\end{lemma} 

\begin{proof} Let us verify that the formulas for the coproducts of all generators of $ U_q^\A(\mathfrak{g}) $ make sense 
as elements in $ U_q^\A(\mathfrak{g}) \otimes_A U_q^\A(\mathfrak{g}) $. For the generators $ K_\lambda $ this is obvious. 
Moreover we observe 
$$
\hat{\Delta}([K_i;0]) = [K_i;0] \otimes K_i + K_i^{-1} \otimes [K_i;0]. 
$$
According to Lemma \ref{lemEiFicomult} we have 
\begin{align*}
\hat{\Delta}(E_j^{(r)}) &= \sum_{i = 0}^r q_j^{i (r - i)} E_j^{(i)} \otimes E_j^{(r - i)} K_j^i \\ 
\hat{\Delta}(F_j^{(r)}) &= \sum_{i = 0}^r q_j^{i (r - i)} F_j^{(r - i)} K_j^{-i} \otimes F_j^{(i)}
\end{align*}
for all $ j $, so the assertion also holds for the divided powers. 

It is clear that the antipode and counit of $ U_q(\mathfrak{g}) $ induce corresponding maps on the level of $ U_q^\A(\mathfrak{g}) $. 
The Hopf algebra axioms are verified in the same way as in Lemma \ref{lemuqhopf}. \end{proof}

We next define the bar involution of $U_q^\A(\mathfrak{g})$.
The field automorphism $ \beta $ of $\mathbb{Q}(s) $ determined by $ \beta(s) = s^{-1} $
\label{nom:beta_field_automorphism}%
restricts to a ring automorphism of $ \A $.  This extends to $U_q^\A(\mathfrak{g})$ as follows.

\begin{lemma} \label{defbarinvolution}
The quantized universal enveloping algebra $ U_q(\mathfrak{g}) $ over $ \mathbb{Q}(s) $ admits an 
automorphism $ \beta: U_q(\mathfrak{g}) \rightarrow U_q(\mathfrak{g}) $
of $ \mathbb{Q} $-algebras such that $ \beta(s) = s^{-1} $ and 
\begin{align*}
\beta(K_\mu) = K_{-\mu}, \qquad \beta(E_i) = E_i, \qquad \beta(F_i) = F_i
\end{align*}
for all $ \mu \in \weights $ and $ 1 \leq i \leq N $.  This restricts to an automorphism of $U_q^\A(\mathfrak{g})$.
\end{lemma} 

\begin{proof} We can view $ U_q(\mathfrak{g}) $ as a $ \mathbb{Q}(s) $-algebra using the scalar action $ c \bullet X = \beta(c) X $ for $ c \in \mathbb{Q}(s) $. 
Let us write $ U_q(\mathfrak{g})^\beta $ for the resulting $ \mathbb{Q}(s) $-algebra. Sending $ K_\mu $ to $ K_{-\mu} $ and fixing the generators $ E_i, F_j $ 
determines a homomorphism $ U_q(\mathfrak{g}) \rightarrow U_q(\mathfrak{g})^\beta $ of $ \mathbb{Q}(s) $-algebras. By slight abuse of notation 
we can view this as the desired automorphism $ \beta: U_q(\mathfrak{g}) \rightarrow U_q(\mathfrak{g}) $ of $ \mathbb{Q} $-algebras. 

The fact that $\beta$ restricts to $U_q^\A(\mathfrak{g})$ follows directly from the definition of the integral form.
\end{proof} 

In the literature, the bar involution $ \beta $ of $ U_q(\mathfrak{g}) $ is usually denoted by a bar. 
We shall write $ \beta(X) $ instead of $ \overline{X} $ for $ X \in U_q(\mathfrak{g}) $ in order 
to avoid confusion with $ * $-structures later on. Note that $ \beta: U_q(\mathfrak{g}) \rightarrow U_q(\mathfrak{g}) $ is not a Hopf algebra automorphism.

For $l\in\mathbb{Z}$ and $1\leq i\leq N$ we will use the notation
\[
 [K_i;l]_{q_i} = [K_i; l] = \frac{q_i^lK_i - q_i^{-l}K_i^{-1}}{q_i - q_i^{-1}},
\]
\label{nom:K:l}%
and for any $m\in\mathbb{N}_0$ we define
\[
 \begin{bmatrix} K_i; l \\ m \end{bmatrix}_{q_i}
   = \frac{[K_i;l]_{q_i} [K_i;l-1]_{q_i} \ldots [K_i;l-m+1]_{q_i}}{[m]_{q_i}!}
   = \prod_{j = 1}^m \frac{q_i^{l + 1 - j} K_i - q_i^{-(l + 1 - j)} K_i^{-1}}{q_i^j - q_i^{-j}},
\]
\nomenclature[o$(K_i)$3]{$\begin{bmatrix} K_i; l \\ m \end{bmatrix}_{q_i}$}{$= \displaystyle \frac{[K_i;l]_{q_i} [K_i;l-1]_{q_i} \ldots [K_i;l-m+1]_{q_i}}{[m]_{q_i}"!}$}%
where, as usual, we interpret this as $1$ when $m=0$.  Note that 
\[
 \begin{bmatrix} K_i; l \\ 1 \end{bmatrix}_{q_i} = [K_i;l]_{q_i}.
\]
These expressions should be considered as elements in $U_q(\mathfrak{g})$. We will see in Proposition \ref{PBWcartanintegral} that they in fact belong to the restricted integral form $U_q^\A(\mathfrak{g})$.

\begin{lemma} \label{uqgintegralKlemma}
Let $ m \in \mathbb{N}_0, l \in \mathbb{Z} $ and $ 1 \leq i \leq N $. Then the following relations hold. 
\begin{bnum} 
\item[a)] We have 
$$
q_i^{-(m + 1)} \begin{bmatrix}
K_i; l \\
m + 1
\end{bmatrix}_{q_i}
= K_i^{-1} q_i^{-l} \begin{bmatrix}
K_i; l - 1 \\
m 
\end{bmatrix}_{q_i}
+ \begin{bmatrix}
K_i; l - 1\\
m + 1
\end{bmatrix}_{q_i}
$$
and 
$$
q_i^{m + 1} \begin{bmatrix}
K_i; l \\
m + 1
\end{bmatrix}_{q_i}
= K_i q_i^l \begin{bmatrix}
K_i; l - 1 \\
m 
\end{bmatrix}_{q_i}
+ \begin{bmatrix}
K_i; l - 1\\
m + 1
\end{bmatrix}_{q_i}. 
$$
\item[b)] If $ l \geq 0 $ then
$$
\begin{bmatrix}
K_i; l \\
m
\end{bmatrix}_{q_i}
= \sum_{j = 0}^m
q_i^{l(m - j)} 
\begin{bmatrix}
l \\
j
\end{bmatrix}_{q_i}
K_i^{-j}
\begin{bmatrix}
K_i; 0 \\
m - j
\end{bmatrix}_{q_i}.
$$
\item[c)] If $ l < 0 $ then
$$
\begin{bmatrix}
K_i; l \\
m
\end{bmatrix}_{q_i}
= \sum_{j = 0}^m (-1)^j
q_i^{l(j - m)} 
\begin{bmatrix}
- l + j - 1 \\
j
\end{bmatrix}_{q_i}
K_i^j
\begin{bmatrix}
K_i; 0 \\
m - j
\end{bmatrix}_{q_i}. 
$$
\end{bnum}
\end{lemma} 

\begin{proof} $ a) $ To verify the first equation we use induction on $ m $. For $ m = 0 $ one computes
\begin{align*}
q_i^{-1} \begin{bmatrix}
K_i; l \\
1
\end{bmatrix}_{q_i}
&= \frac{q_i^{l - 1} K_i - q_i^{-l - 1} K_i^{-1}}{q_i - q_i^{-1}} \\
&= \frac{q_i^{l - 1} K_i - q_i^{-(l - 1)} K_i^{-1} + q_i^{-l} (q_i - q_i^{-1}) K_i^{-1}}{q_i - q_i^{-1}} \\
&= q_i^{-l} K_i^{-1} + \frac{q_i^{l - 1} K_i - q_i^{-(l - 1)} K_i^{-1}}{q_i - q_i^{-1}} \\
&= q_i^{-l} K_i^{-1} \begin{bmatrix}
K_i; l - 1 \\
0
\end{bmatrix}_{q_i}
+ \begin{bmatrix}
K_i; l - 1\\
1
\end{bmatrix}_{q_i}
\end{align*}
as desired. Assume that the assertion is proved for $ m - 1 \geq 0 $ and compute 
\begin{align*}
&q_i^{-(m + 1)} 
\begin{bmatrix}
K_i; l \\
m + 1
\end{bmatrix}_{q_i}
= q_i^{-m} 
\begin{bmatrix}
K_i; l \\
m 
\end{bmatrix}_{q_i}
\frac{q_i^{l - m - 1} K_i - q_i^{-(l - m + 1)} K_i^{-1}}{q_i^{m + 1} - q_i^{-(m + 1)}} \\
&= \biggl(K_i^{-1} q_i^{-l} 
\begin{bmatrix}
K_i; l - 1 \\
m - 1
\end{bmatrix}_{q_i}
+ \begin{bmatrix}
K_i; l - 1\\
m
\end{bmatrix}_{q_i}\biggr) \frac{q_i^{l - m - 1} K_i - q_i^{-(l - m + 1)} K_i^{-1}}{q_i^{m + 1} - q_i^{-(m + 1)}} \\
&= K_i^{-1} q_i^{-l} 
\begin{bmatrix}
K_i; l - 1 \\
m - 1
\end{bmatrix}_{q_i} \frac{q_i^{l - m - 1} K_i - q_i^{-(l - m + 1)} K_i^{-1}}{q_i^{m + 1} - q_i^{-(m + 1)}} \\
&\qquad  + 
\begin{bmatrix}
K_i; l - 1\\
m 
\end{bmatrix}_{q_i}
\frac{q_i^{-(l - m - 1)} K_i^{-1} - q_i^{-(l - m + 1)} K_i^{-1}}{q_i^{m + 1} - q_i^{-(m + 1)}} 
+ 
\begin{bmatrix}
K_i; l - 1\\
m + 1
\end{bmatrix}_{q_i}
\\
&
= 
K_i^{-1} q_i^{-l} 
\begin{bmatrix}
K_i; l - 1 \\
m 
\end{bmatrix}_{q_i} \biggl(\frac{q_i^{- 1}}{q_i^{m + 1} - q_i^{-(m + 1)}}(q_i^m - q_i^{-m}) 
+ 
\frac{q_i^{m + 1} - q_i^{m - 1}}{q_i^{m + 1} - q_i^{-(m + 1)}} \biggr) 
\\
&\qquad 
+ 
\begin{bmatrix}
K_i; l - 1\\
m + 1
\end{bmatrix}_{q_i} \\
&= K_i^{-1} q_i^{-l} \begin{bmatrix}
K_i; l - 1 \\
m  
\end{bmatrix}_{q_i}
+ \begin{bmatrix}
K_i; l - 1\\
m + 1
\end{bmatrix}_{q_i}. 
\end{align*}
This yields the first equality. Applying the automorphism $ \beta $ to this equality yields the second claim, where we note that the 
terms $ \begin{bmatrix}
K_i; l \\
m
\end{bmatrix}_{q_i}
$ are fixed under $ \beta $ for all $ l \in \mathbb{Z} $ and $ m \in \mathbb{N}_0 $. 

$ b) $ We use induction on $ l $. Assume without loss of generality that $ m > 0 $. 
For $ l = 0 $ the assertion clearly holds. To check the inductive step we use $ a) $ and the inductive hypothesis to calculate
\begin{align*}
\begin{bmatrix}
K_i; l + 1 \\
m
\end{bmatrix}_{q_i}
&= K_i^{-1} q_i^{m - l - 1} 
\begin{bmatrix}
K_i; l \\
m - 1
\end{bmatrix}_{q_i}
+ q_i^m \begin{bmatrix}
K_i; l \\
m 
\end{bmatrix}_{q_i} \\
&= K_i^{-1} q_i^{m - l - 1} 
\sum_{j = 0}^{m - 1}
q_i^{l(m - 1 - j)} 
\begin{bmatrix}
l \\
j
\end{bmatrix}_{q_i}
K_i^{-j}
\begin{bmatrix}
K_i; 0 \\
m - 1 - j
\end{bmatrix}_{q_i} \\
&\qquad + q_i^m \sum_{j = 0}^m
q_i^{l(m - j)} 
\begin{bmatrix}
l \\
j
\end{bmatrix}_{q_i}
K_i^{-j}
\begin{bmatrix}
K_i; 0 \\
m - j
\end{bmatrix}_{q_i} \\
&= q_i^{m - l - 1} 
\sum_{j = 1}^m
q_i^{l(m - j)} 
\begin{bmatrix}
l \\
j - 1
\end{bmatrix}_{q_i}
K_i^{-j}
\begin{bmatrix}
K_i; 0 \\
m - j
\end{bmatrix}_{q_i} \\
&\qquad + q_i^m \sum_{j = 0}^m
q_i^{l(m - j)} 
\begin{bmatrix}
l \\
j
\end{bmatrix}_{q_i}
K_i^{-j}
\begin{bmatrix}
K_i; 0 \\
m - j
\end{bmatrix}_{q_i} \\
&= \sum_{j = 0}^m
q_i^{(l + 1)(m - j)} 
\begin{bmatrix}
l + 1 \\
j
\end{bmatrix}_{q_i}
K_i^{-j}
\begin{bmatrix}
K_i; 0 \\
m - j
\end{bmatrix}_{q_i}, 
\end{align*}
taking into account Lemma \ref{lemqcalc1} in the last step.

$ c) $ We first show 
\begin{align*}
\prod_{k = 1}^m \frac{q_i^{- k} K_i - q_i^{k} K_i^{-1}}{q_i^k - q_i^{-k}} 
&= \sum_{j = 0}^m (-1)^j
q_i^{-(j - m)} 
K_i^j
\prod_{k = 1}^{m - j} \frac{q_i^{1 - k} K_i - q_i^{-(1 - k)} K_i^{-1}}{q_i^k - q_i^{-k}} 
\end{align*} 
by induction on $ m $. For $ m = 0 $ this relation clearly holds. 
If it holds for $ m \geq 0 $ then using the inductive hypothesis we get 
\begin{align*}
\prod_{k = 1}^{m + 1} &\frac{q_i^{- k} K_i - q_i^{k} K_i^{-1}}{q_i^k - q_i^{-k}} 
= \prod_{k = 1}^{m + 1} \frac{q_i^{- k} K_i - q_i^{k} K_i^{-1}}{q_i^k - q_i^{-k}} + 
q_i^{m + 1} \prod_{k = 1}^{m + 1} \frac{q_i^{1 - k} K_i - q_i^{-(1 - k)} K_i^{-1}}{q_i^k - q_i^{-k}} \\
&- q_i^{m + 1} \frac{K_i - K_i^{-1}}{q_i^{m + 1} - q_i^{-(m + 1)}} 
\prod_{k = 1}^{m} \frac{q_i^{- k} K_i - q_i^{k} K_i^{-1}}{q_i^k - q_i^{-k}} \\
&= q_i^{m + 1} \prod_{k = 1}^{m + 1} \frac{q_i^{1 - k} K_i - q_i^{-(1 - k)} K_i^{-1}}{q_i^k - q_i^{-k}} 
+ \biggl(\frac{q_i^{-(m + 1)} K_i - q_i^{m + 1} K_i^{-1}}{q_i^{m + 1} - q_i^{-(m + 1)}} \\
&\qquad - \frac{q_i^{m + 1} K_i - q_i^{m + 1} K_i^{-1}}{q_i^{m + 1} - q_i^{-(m + 1)}} \biggr) 
\sum_{j = 0}^{m} (-1)^j q_i^{-(j - m)} K_i^j
\prod_{k = 1}^{m - j} \frac{q_i^{1 - k} K_i - q_i^{-(1 - k)} K_i^{-1}}{q_i^k - q_i^{-k}} \\
&= q_i^{m + 1} \prod_{k = 1}^{m + 1} \frac{q_i^{1 - k} K_i - q_i^{-(1 - k)} K_i^{-1}}{q_i^k - q_i^{-k}} \\
&\qquad - K_i \sum_{j = 0}^{m} (-1)^j
q_i^{-(j - m)} K_i^j
\prod_{k = 1}^{m - j} \frac{q_i^{1 - k} K_i - q_i^{-(1 - k)} K_i^{-1}}{q_i^k - q_i^{-k}} \allowdisplaybreaks \\
&= q_i^{m + 1} \prod_{k = 1}^{m + 1} \frac{q_i^{1 - k} K_i - q_i^{-(1 - k)} K_i^{-1}}{q_i^k - q_i^{-k}} \\
&\qquad + K_i \sum_{j = 1}^{m + 1} (-1)^j
q_i^{-(j - m - 1)} K_i^{j - 1}
\prod_{k = 1}^{m + 1 - j} \frac{q_i^{1 - k} K_i - q_i^{-(1 - k)} K_i^{-1}}{q_i^k - q_i^{-k}} \\
&= \sum_{j = 0}^{m + 1} (-1)^j
q_i^{-(j - m - 1)} 
K_i^j
\prod_{k = 1}^{m + 1 - j} \frac{q_i^{1 - k} K_i - q_i^{-(1 - k)} K_i^{-1}}{q_i^k - q_i^{-k}} 
\end{align*} 
as desired. 

Let us now prove the assertion. For $ l = -1 $ we obtain 
\begin{align*}
\begin{bmatrix}
K_i; -1 \\
m
\end{bmatrix}_{q_i}
&= \prod_{k = 1}^m \frac{q_i^{- k} K_i - q_i^{k} K_i^{-1}}{q_i^k - q_i^{-k}} \\
&= \sum_{j = 0}^m (-1)^j
q_i^{-(j - m)} 
K_i^j
\prod_{k = 1}^{m - j} \frac{q_i^{1 - k} K_i - q_i^{-(1 - k)} K_i^{-1}}{q_i^k - q_i^{-k}} \\
&= \sum_{j = 0}^m (-1)^j
q_i^{-(j - m)} 
K_i^j
\begin{bmatrix}
K_i; 0 \\
m - j
\end{bmatrix}_{q_i}
\end{align*} 
using the above calculation. 
Assume now that the claim holds for some $ l < 0 $ and all $ m \in \mathbb{N}_0 $. For $ l - 1 $ and $ m = 0 $ the assertion clearly 
holds too. Using induction on $ m $ and applying the second formula from $ a) $ we get 
\begin{align*}
\begin{bmatrix}
K_i; l - 1 \\
m
\end{bmatrix}_{q_i}
&= 
q_i^m \begin{bmatrix}
K_i; l \\
m 
\end{bmatrix}_{q_i}
- K_i q_i^l \begin{bmatrix}
K_i; l - 1 \\
m - 1
\end{bmatrix}_{q_i}
\\
&= q_i^m \sum_{j = 0}^m (-1)^j
q_i^{l(j - m)} 
\begin{bmatrix}
- l + j - 1 \\
j
\end{bmatrix}_{q_i}
K_i^j
\begin{bmatrix}
K_i; 0 \\
m - j
\end{bmatrix}_{q_i} \\
&\qquad - K_i q_i^l 
\sum_{j = 0}^{m - 1} (-1)^j
q_i^{(l - 1)(j - m + 1)} 
\begin{bmatrix}
- l + j \\
j
\end{bmatrix}_{q_i}
K_i^j
\begin{bmatrix}
K_i; 0 \\
m - 1 - j
\end{bmatrix}_{q_i} \\
&= q_i^m \sum_{j = 0}^m (-1)^j
q_i^{l(j - m)} 
\begin{bmatrix}
- l + j - 1 \\
j
\end{bmatrix}_{q_i}
K_i^j
\begin{bmatrix}
K_i; 0 \\
m - j
\end{bmatrix}_{q_i} \\
&\qquad + q_i^l
\sum_{j = 1}^{m} (-1)^j
q_i^{(l - 1)(j - m)} 
\begin{bmatrix}
- l + j - 1 \\
j - 1
\end{bmatrix}_{q_i}
K_i^j
\begin{bmatrix}
K_i; 0 \\
m - j
\end{bmatrix}_{q_i} \\
&= \sum_{j = 0}^m (-1)^j
q_i^{(l - 1)(j - m)} 
\begin{bmatrix}
- l + j \\
j
\end{bmatrix}_{q_i}
K_i^j
\begin{bmatrix}
K_i; 0 \\
m - j
\end{bmatrix}_{q_i}
\end{align*} 
as claimed, again taking into account Lemma \ref{lemqcalc1}. \end{proof}

We collect some more formulas. 

\begin{lemma} \label{uqgintegralmultiplication}
Assume $ m \in \mathbb{N} $ and $ n \in \mathbb{N}_0 $. Then we have 
\begin{align*}
\sum_{j = 0}^n 
(-1)^j q_i^{m(n - j)} 
\begin{bmatrix}
m + j - 1 \\
j
\end{bmatrix}_{q_i}
K_i^j 
\begin{bmatrix}
K_i; 0 \\
m 
\end{bmatrix}_{q_i}
\begin{bmatrix}
K_i; 0 \\
n - j
\end{bmatrix}_{q_i} 
&= \begin{bmatrix}
K_i; 0 \\
m 
\end{bmatrix}_{q_i}
\begin{bmatrix}
K_i; -m \\
n 
\end{bmatrix}_{q_i} \\
&= \begin{bmatrix}
m + n \\
n 
\end{bmatrix}_{q_i}
\begin{bmatrix}
K_i; 0 \\
m + n
\end{bmatrix}_{q_i} 
\end{align*}
for all $ i = 1, \dots, N $. 
\end{lemma} 

\begin{proof} The first equality follows from Lemma \ref{uqgintegralKlemma} $ c) $. 
For the second equality we compute
\begin{align*}
\begin{bmatrix}
K_i; 0 \\
m
\end{bmatrix}_{q_i}
&\begin{bmatrix}
K_i; -m \\
n
\end{bmatrix}_{q_i} \\
&= \prod_{j = 1}^m \frac{q_i^{1 - j} K_i - q_i^{-(1 - j)} K_i^{-1}}{q_i^j - q_i^{-j}} 
\prod_{k = 1}^n \frac{q_i^{-m + 1 - k} K_i - q_i^{-(-m + 1 - k)} K_i^{-1}}{q_i^k - q_i^{-k}} \\
&= \prod_{j = 1}^m \frac{q_i^{1 - j} K_i - q_i^{-(1 - j)} K_i^{-1}}{q_i^j - q_i^{-j}} 
\prod_{k = m + 1}^{m + n} \frac{q_i^{1 - k} K_i - q_i^{-(1 - k)} K_i^{-1}}{q_i^{k - m} - q_i^{-(k - m)}} \\
&= \begin{bmatrix}
m + n \\
n 
\end{bmatrix}_{q_i} 
\prod_{j = 1}^{m + n} \frac{q_i^{1 - j} K_i - q_i^{-(1 - j)} K_i^{-1}}{q_i^j - q_i^{-j}} \\
&= \begin{bmatrix}
m + n \\
n 
\end{bmatrix}_{q_i}
\begin{bmatrix}
K_i; 0\\
m + n
\end{bmatrix}_{q_i}
\end{align*}
as desired. \end{proof} 

\begin{lemma} \label{uqgintegralcommutation}
Given any $ l \in \mathbb{Z}, m \in \mathbb{N}_0 $ we have 
\begin{align*}
\begin{bmatrix}
K_i; l \\
m
\end{bmatrix}_{q_i} 
E_i^{(r)} 
&= E_i^{(r)}
\begin{bmatrix}
K_i; l + 2r \\
m
\end{bmatrix}_{q_i} \\
\begin{bmatrix}
K_i; l \\
m
\end{bmatrix}_{q_i} 
F_i^{(r)} 
&= F_i^{(r)}
\begin{bmatrix}
K_i; l - 2r \\
m
\end{bmatrix}_{q_i} 
\end{align*}
for all $ r \in \mathbb{N}_0 $. 
\end{lemma} 

\begin{proof} We compute 
\begin{align*}
\begin{bmatrix}
K_i; l \\
m
\end{bmatrix}_{q_i} 
E_i^{(r)} &= \prod_{j = 1}^m \frac{q_i^{l + 1 - j} K_i - q_i^{-(l + 1 - j)} K_i^{-1}}{q_i^j - q_i^{-j}} E_i^{(r)} \\
&= E_i^{(r)} \prod_{j = 1}^m \frac{q_i^{l + 1 - j + 2r} K_i - q_i^{-(l + 1 - j + 2r)} K_i^{-1}}{q_i^j - q_i^{-j}} \\
&= E_i^{(r)}
\begin{bmatrix}
K_i; l + 2r \\
m
\end{bmatrix}_{q_i} 
\end{align*}
as desired. The proof of the second formula is analogous. \end{proof} 

According to Lemma \ref{EFpowercommutationhelp} We have 
\begin{align*}
[E_i, F_i^m] &= [m]_{q_i} F_i^{m - 1} [K_i; 1 - m]_{q_i} 
\end{align*}
for $ m \geq 1 $.
This can also be written in the form 
\begin{align*}
[E_i, F_i^{(m)}] &= F_i^{(m - 1)} [K_i; 1 - m]_{q_i}.  
\end{align*}
Let us record a more general commutation relation. 

\begin{prop} \label{EFpowercommutation}
We have 
\begin{align*}
E_i^{(r)} F_i^{(s)} &= 
\sum_{r,s \geq t \geq 0} F_i^{(s - t)} E_i^{(r - t)} 
\begin{bmatrix}
K_i; r - s \\
t
\end{bmatrix}_{q_i} \\
&= \sum_{r,s \geq t \geq 0} F_i^{(s - t)} 
\begin{bmatrix}
K_i; 2t - r - s \\
t
\end{bmatrix}_{q_i} 
E_i^{(r - t)} 
\end{align*}
for all $ r, s \in \mathbb{N}_0 $.  
\end{prop} 

\begin{proof} Let us verify the first equality. 
For $ r = 0 $ and $ s \in \mathbb{N}_0 $ arbitrary the formula clearly holds. For $ r = 1 $ and $ s \in \mathbb{N}_0 $ arbitrary the 
formula follows from the computation just before the proposition. 

Assume that the equality is true for some $ r \geq 1 $ and all $ s \in \mathbb{N}_0 $. Then we compute 
\begin{align*}
E_i^{(r + 1)} &F_i^{(s)} = \frac{1}{[r + 1]_{q_i}} E_i^{(r)} E_i F_i^{(s)} \\
&= \frac{1}{[r + 1]_{q_i}} (E_i^{(r)} F_i^{(s - 1)} [K_i; 1 - s] + E_i^{(r)} F_i^{(s)} E_i) \\
&= \frac{1}{[r + 1]_{q_i}} \sum_{r,s - 1 \geq t \geq 0} F_i^{(s - 1 - t)}
E_i^{(r - t)} 
\begin{bmatrix}
K_i; r - s + 1 \\
t
\end{bmatrix}_{q_i}
[K_i; 1 - s] \\
&\qquad + \frac{1}{[r + 1]_{q_i}} \sum_{r,s \geq t \geq 0} F_i^{(s - t)}
E_i^{(r - t)} 
\begin{bmatrix}
K_i; r - s \\
t
\end{bmatrix}_{q_i}
E_i \\
&= \frac{1}{[r + 1]_{q_i}} \sum_{r + 1,s \geq t + 1 \geq 1} F_i^{(s - 1 - t)}
E_i^{(r - t)} 
\begin{bmatrix}
K_i; r - s + 1 \\
t
\end{bmatrix}_{q_i}
[K_i; 1 - s] \\
&\qquad + \frac{1}{[r + 1]_{q_i}} \sum_{r,s \geq t \geq 0} [r + 1 - t]_{q_i} F_i^{(s - t)}
E_i^{(r + 1 - t)} 
\begin{bmatrix}
K_i; r - s + 2 \\
t
\end{bmatrix}_{q_i}
\\
&= \frac{1}{[r + 1]_{q_i}} \sum_{r + 1,s \geq t \geq 1} F_i^{(s - t)}
E_i^{(r + 1 - t)} 
\begin{bmatrix}
K_i; r - s + 1 \\
t - 1
\end{bmatrix}_{q_i}
[K_i; 1 - s] \\
&\qquad + \frac{1}{[r + 1]_{q_i}} \sum_{r+1,s \geq t \geq 0} [r + 1 - t]_{q_i} F_i^{(s - t)}
E_i^{(r + 1 - t)} 
\begin{bmatrix}
K_i; r - s + 2 \\
t
\end{bmatrix}_{q_i}, 
\end{align*}
using Lemma \ref{uqgintegralcommutation}. Therefore it suffices to show 
\begin{align*}
\begin{bmatrix}
K_i; r - s + 1 \\
t - 1
\end{bmatrix}_{q_i}
[K_i; 1 - s] &+ 
[r + 1 - t]_{q_i} 
\begin{bmatrix}
K_i; r - s + 2 \\
t
\end{bmatrix}_{q_i} \\
&= 
[r + 1]_{q_i}
\begin{bmatrix}
K_i; r + 1 - s \\
t
\end{bmatrix}_{q_i}
\end{align*}
for $ r,s \geq t \geq 1 $; note that the contribution corresponding to $ t = 0 $ is covered by the second term in the 
previous expression.

We compute the left-hand side,
\begin{align*}
 &\begin{bmatrix} K_i; r - s + 1 \\ t - 1 \end{bmatrix}_{q_i}  [K_i; 1 - s] + 
   [r + 1 - t]_{q_i} \begin{bmatrix} K_i; r - s + 2 \\ t \end{bmatrix}_{q_i} \\
 & = \frac{[K_i;r-s+1]\ldots[K_i;r-s-t+3]}{[t]_{q_i}!}
   \Big(
    [t]_{q_i} [K_i;1-s] + [r+1-t]_{q_i} [K_i;r-s+2]
   \Big)
\end{align*}
The expression in parentheses here is
\begin{align*}
&([t]_{q_i} [K_i; 1 - s] + [r + 1 - t]_{q_i} [K_i; r - s + 2]) \\
&= {\textstyle \frac{1}{ (q_i - q_i^{-1})^2}} 
    \bigg( (q_i^{t} - q_i^{-t})(q_i^{1 - s} K_i - q_i^{-(1 - s)} K_i^{-1}) \\
    &\hspace{3cm} + (q_i^{r + 1 - t} - q_i^{-(r + 1 - t)})(q_i^{r - s + 2} K_i - q_i^{-(r - s + 2)} K_i^{-1}) 
   \bigg)\\
&= {\textstyle \frac{1}{ (q_i - q_i^{-1})^2}}
    \bigg( (q_i^{2r + 3 - t - s} - q_i^{-t + 1 - s}) K_i + (q_i^{-2r - 3 + t + s} - q_i^{t - 1 + s}) K_i^{-1} 
   \bigg)\\ 
&= {\textstyle \frac{1}{ (q_i - q_i^{-1})^2}}(q_i^{r + 1} - q_i^{-(r + 1)})(q_i^{r - s + 2 - t} K_i - q_i^{-(r - s + 2 - t)} K_i^{-1}) \\
&= [r + 1]_{q_i} [K_i; r - s + 2 - t],
\end{align*}
which yields the claim.

The second equality follows from Lemma \ref{uqgintegralcommutation}. \end{proof} 

We denote by $ U_q^\A(\mathfrak{n}_+) $
\nomenclature{$U_q^\A(\mathfrak{n}_\pm)$}{subalgebras of $U_q^\A(\lie{g})$ generated by divided powers of $E_i$, $F_i$, respectively}%
the subalgebra of $ U_q^\A(\mathfrak{g}) $ generated 
by the divided powers $ E_1^{(t_1)}, \dots, E_N^{(t_N)} $ with $ t_j \in \mathbb{N}_0 $, and 
let $ U_q^\A(\mathfrak{n}_-) $ be the subalgebra generated by the corresponding divided powers $ F_1^{(t_1)}, \dots, F_N^{(t_N)} $. 
Moreover we let $ U_q^\A(\mathfrak{h}) $
\nomenclature{$U_q^\A(\mathfrak{h})$}{Cartan subalgebra $U_q^\A(\mathfrak{g})\cap U_q(\mathfrak{h})$}%
be the intersection of $ U_q^\A(\mathfrak{g}) $ 
with $ U_q(\mathfrak{h}) $. 

The algebras $ U_q^\A(\mathfrak{n}_\pm) $ will be studied in greater 
detail later on. Here we shall determine the structure of $ U_q^\A(\mathfrak{h}) $, see Theorem 6.7 in \cite{Lusztigrootsofone}.

\begin{prop} \label{PBWcartanintegral}
The algebra $ U_q^\A(\mathfrak{h}) $ is free as an $ \A $-module, with basis given by all elements of the form
$$
K_\lambda \prod_{i = 1}^N K_i^{s_i} 
\begin{bmatrix}
K_i; 0 \\
m_i
\end{bmatrix}_{q_i}
$$
where  $s_i \in \{0,1\} $, $ m_i \in \mathbb{N}_0 $ and $ \lambda $ belongs to a set of representatives of the cosets of $\roots$ in $\weights$. 
\end{prop}

\begin{proof} Let us first show by induction on $ m $ that $ U_q^\A(\mathfrak{h}) $ contains all elements 
$ \begin{bmatrix}
K_i;l \\
m
\end{bmatrix}_{q_i}
$ 
for $ l \in \mathbb{Z} $ and $ m \in \mathbb{N}_0 $. The assertion clearly holds for $ m = 0 $. 
Assume that it holds for some $ m \geq 0 $. Then Proposition \ref{EFpowercommutation} applied to $ r = s = m + 1 $ and induction shows that 
$ \begin{bmatrix}
K_i; 0 \\
m + 1
\end{bmatrix}_{q_i}
$ 
is contained in $ U_q^\A(\mathfrak{h}) $. Finally, parts $ b) $ and $ c) $ of Lemma \ref{uqgintegralKlemma} and induction allow us to conclude that $ \begin{bmatrix}
K_i; l \\
m + 1
\end{bmatrix}_{q_i}
$ 
is in $ U_q^\A(\mathfrak{h}) $ for all positive and negative integers $ l $, respectively.

For $ \lambda \in \weights $ and $ {\bf s} = (s_1, \dots, s_N) \in \mathbb{Z}^N$, ${\bf m} = (m_1, \dots, m_N) \in \mathbb{N}_0^N $ let us abbreviate 
$$ 
K(\lambda, {\bf s}, {\bf m}) = K_\lambda \prod_{i = 1}^N K_i^{s_i} \begin{bmatrix}
K_i; 0 \\
m_i
\end{bmatrix}_{q_i}. 
$$ 
According to the definition of $ U_q^\A(\mathfrak{h}) $ and our above considerations, these elements are all contained in $ U_q^\A(\mathfrak{h}) $. 

Let us write $H$ for the $ \A $-linear span of all 
elements of the form $K_\lambda \prod_{i=1}^N \begin{bmatrix}
K_i; l_i \\
m_i
\end{bmatrix}_{q_i}
$
with $\lambda\in\weights$,  $ l_i \in \mathbb{Z} $ and $ m_i \in \mathbb{N}_0 $.  
Then $ H \subset U_q^\A(\mathfrak{h}) $.
We claim that $H$ is closed under multiplication.  For this, it suffices to consider products of the form
$ \begin{bmatrix}
K_i; l \\
m 
\end{bmatrix}_{q_i}
\begin{bmatrix}
K_i; l' \\
m'
\end{bmatrix}_{q_i}
$
with $l,l'\in\ZZ$, $m,m'\in\NN_0$.  Using Lemma \ref{uqgintegralKlemma}, we can reduce to the case $l=l'=0$, which follows from Lemma \ref{uqgintegralmultiplication} by induction on $m'$.
Next, note that
according to Proposition \ref{EFpowercommutation} and Lemma \ref{uqgintegralcommutation}, 
the $ \A $-module $ U_q^\A(\mathfrak{n}_-) H U_q^\A(\mathfrak{n}_+) $ is closed under multiplication, and therefore agrees with $ U_q^\A(\mathfrak{g}) $. 
This implies $ H = U_q^\A(\mathfrak{h}) $. 

Now fix a set of coset representatives of $\roots$ in $\weights$, which we will denote temporarily by $\mathbf{R}$. Let us write $ L $ for the $ \A $-linear span of all $ K(\lambda, {\bf s}, {\bf m}) $ 
with $ \lambda \in \mathbf{R} $ and $ {\bf s} = (s_1, \dots, s_N) \in \{0,1\}^N, {\bf m} = (m_1, \dots, m_N) \in \mathbb{N}_0^N $.  
We will prove that $L$ contains $K(\lambda,{\bf t}, {\bf m})$ for all $\lambda\in\mathbf{P}$, ${\bf t} \in \mathbb{Z}^N$, ${\bf m} \in \NN_0^N$ by induction on ${\bf t}$.  
For this, it is enough to show that $L$ contains $\prod_{i=1}^N K_i^{t_i} \begin{bmatrix}
K_i; 0 \\
m_i
\end{bmatrix}_{q_i} $ for all ${\bf t}\in\ZZ^N$ and ${\bf m}\in \mathbb{N}_0^N$.
This is trivially true for  ${\bf t } \in \{0,1\}^N$ and any ${\bf m}$.  
From the definition of $\begin{bmatrix} K_i; -m_i \\ t_i+1 \end{bmatrix}_{q_i}$ we can write
\begin{align*}
 K_i^{t_i+2}
   &= a K_i \begin{bmatrix} K_i; -m_i \\ t_i+1 \end{bmatrix}_{q_i} + \sum_{k=-t_i}^{t_i} b_k K_i^k
\end{align*}
for some $a, b_{-t_i},\ldots, b_{t_i} \in \A$, and likewise
\begin{align*}
 K_i^{-t_i-1}
   = a' \begin{bmatrix} K_i; -m_i \\ t_i+1 \end{bmatrix}_{q_i} + \sum_{k=-t_i+1}^{t_i+1} b'_k K_i^k
\end{align*}
for some $a', b'_{-t_i+1}, \ldots, b'_{t_i+1} \in \A$.  Using the second equality from Lemma \ref{uqgintegralmultiplication} we therefore obtain
\begin{align*}
  K_i^{t_i+2} \begin{bmatrix} K_i; 0 \\ m_i \end{bmatrix}_{q_i}
   &= a \begin{bmatrix} m_i+t_i+1 \\ t_i+1 \end{bmatrix}_{q_i} K_i 
      \begin{bmatrix} K_i; 0 \\ m_i+t_i+1 \end{bmatrix}_{q_i} + \sum_{k=-t_i}^{t_i} b_k K_i^k 
      \begin{bmatrix} K_i; 0 \\ m_i \end{bmatrix}_{q_i}, \\
  K_i^{-t_i-1} \begin{bmatrix} K_i; 0 \\ m_i \end{bmatrix}_{q_i} 
   &= a' \begin{bmatrix} m_i+t_i+1 \\ t_i+1 \end{bmatrix}_{q_i}
      \begin{bmatrix} K_i; 0 \\ m_i+t_i+1 \end{bmatrix}_{q_i} + \sum_{k=-t_i+1}^{t_i+1} b'_k K_i^k 
      \begin{bmatrix} K_i; 0 \\ m_i \end{bmatrix}_{q_i}.
\end{align*}
Performing an induction on $t_i$ successively for $i=1,\ldots,N$, we see that every $K(\lambda,{\bf t},{\bf m})$ is a linear combination of the elements $K(\lambda,{\bf s}, {\bf m})$ with ${\bf s} \in \{0,1\}^N$.  
Taking into account Lemma \ref{uqgintegralKlemma} $ b), c)$ we conclude $ L = H = U_q^\A(\mathfrak{h}) $ as desired. 

It remains to show that the elements $ K(\lambda, {\bf s}, {\bf m}) $ in $ L $ are linearly independent over $ \A $. For this it suffices to observe that the elements $ K(0, {\bf s}, {\bf m}) $ with ${\bf s} \in \{0,1\}^N$ and $ {\bf m} = (m_1, \dots, m_N) $ satisfying $ m_j \leq m $ for all $ j = 1, \dots, N $ 
span the $ 2^N (m + 1)^N $-dimensional subspace of $ U_q(\mathfrak{g}) $ over $ \mathbb{Q}(s) $ generated by all elements 
of the form $ K_1^{r_1} \cdots K_N^{r_N} $ such that $ -m \leq r_i \leq m + 1 $ for all $ i $.  
\end{proof}

It will be shown later on that $ U_q^\A(\mathfrak{g}) $ is indeed an integral form of the quantized universal enveloping algebra $ U_q(\mathfrak{g}) $ 
over $ \mathbb{K} = \mathbb{Q}(s) $ in the sense that there is a canonical isomorphism $ \mathbb{Q}(s) \otimes_\A U_q^\A(\mathfrak{g}) \cong U_q(\mathfrak{g}) $. 

Using the integral form one can also show that $ U_q(\mathfrak{g}) $ tends to the classical universal enveloping algebra $ U(\mathfrak{g}) $ 
as $ q $ tends to $ 1 $. More precisely, let $ \mathbb{K} $ be a field and let $ U_1(\mathfrak{g}) = U_q^\A(\mathfrak{g}) \otimes_\A \mathbb{K} $ 
be the algebra obtained by extension of scalars from $ \A $ to $ \mathbb{K} $ such that $ s $ is mapped to $ 1 \in \mathbb{K} $. 

Recall that the classical universal enveloping algebra $ U(\mathfrak{g}) $ of $ \mathfrak{g} $ over $ \mathbb{K} $ is the $ \mathbb{K} $-algebra with 
generators $ E_i, F_i, H_i $ for $ 1 \leq i \leq N $ satisfying 
\begin{align*}
[H_i, H_j] &= 0 \\ 
[H_i, E_j] &= a_{ij} E_j \\ 
[H_i, F_j] &= -a_{ij} F_j \\ 
[E_i, F_j] &= \delta_{ij} H_j 
\end{align*}
for all $ 1 \leq i, j \leq N $ and the Serre relations 
\begin{align*}
& \sum_{k = 0}^{1 - a_{ij}} (-1)^k \begin{pmatrix} 1 - a_{ij} \\ k \end{pmatrix}
E_i^k E_j E_i^{1 - a_{ij} - k} = 0 \\ 
& \sum_{k = 0}^{1 - a_{ij}} (-1)^k \begin{pmatrix} 1 - a_{ij} \\ k \end{pmatrix}
F_i^k F_j F_i^{1 - a_{ij} - k} = 0
\end{align*} 
for $ i \neq j $. 

\begin{prop} \label{quantumtoclassical}
Let $ U(\mathfrak{g}) $ be the universal enveloping algebra of $ \mathfrak{g} $. Then there 
exists a canonical surjective homomorphism $ U_1(\mathfrak{g}) \rightarrow U(\mathfrak{g}) $ of Hopf algebras. 
\end{prop} 

\begin{proof} Let us write $ e_i, f_i, h_i $ for the generators $ E_i, F_i, [K_i; 0] $ viewed as elements of $ U_1(\mathfrak{g}) $, and similarly write 
$ k_\lambda $ for $ K_\lambda $. Moreover let $ I_q \subset U_q^\A(\mathfrak{g}) $ be the ideal generated by all elements $ K_\lambda - 1 $. 
We claim that $ U_1(\mathfrak{g})/I_1 $ is isomorphic to $ U(\mathfrak{g}) $. 

For this we need to verify that the generators $ e_i, f_i, h_i \in U_1(\mathfrak{g})/I $ satisfy the defining relations of $ U(\mathfrak{g}) $. 
Clearly the elements $ h_i $ commute among themselves. We compute 
\begin{align*}
[[K_i;0], E_j] &= \frac{K_i - K_i^{-1}}{q_i - q_i^{-1}} E_j - E_j \frac{K_i - K_i^{-1}}{q_i - q_i^{-1}} \\
&= \frac{(1 - q_i^{-a_{ij}}) K_i - (1 - q_i^{a_{ij}}) K_i^{-1}}{q_i - q_i^{-1}} E_j \\
&= [a_{ij}]_{q_i} E_j 
\end{align*}
in $ U_q(\mathfrak{g})/I_q $, and hence deduce $ [h_i, e_j] = a_{ij} e_j $ in $ U_1(\mathfrak{g})/I_1 $. 
Similarly one obtains $ [h_i, f_j] = -a_{ij} f_j $, and the relation $ [e_i, f_j] = \delta_{ij} h_j $ follows from the definitions. 
Finally, to check the Serre relations for the $ e_i $ and $ f_j $ it suffices to observe that quantum binomial coefficients 
become ordinary binomial coefficients at $ q = 1 $. 

Hence we obtain an algebra homomorphism $ \iota: U_1(\mathfrak{g}) \rightarrow U(\mathfrak{g}) $ by sending $ e_i, f_i, h_i $ to $ E_i, F_i, H_i $, respectively. 
This map is clearly surjective. It induces an isomorphism $ U_1(\mathfrak{g})/I_1 \rightarrow U(\mathfrak{g}) $ whose inverse is obtained by using the universal property of $ U(\mathfrak{g}) $. 

To check that this map is compatible with coproducts we use Lemma \ref{integralhopfstructure} to calculate 
\begin{align*}
\hat{\Delta}([K_i; 0]) = [K_i; 0] \otimes 1 + 1 \otimes [K_i; 0] 
\end{align*} 
in $ U_q(\mathfrak{g})/I_q $, and the corresponding formulas for $ \hat{\Delta}(E_i), \hat{\Delta}(F_j) $ follow from the definitions. 
This yields the claim. \end{proof}

Let us finally note that the integral form $ U_q^\A(\mathfrak{g}) $ plays a prominent role in the construction of canonical bases. More recently, it has 
appeared in connection with categorification of quantum groups.

\subsection{Verma modules} \label{secvermauqg}

This section contains some basic definitions and results on weight modules and Verma modules for $ U_q(\mathfrak{g}) $. 
Throughout, our basic assumption is that $ \mathbb{K} $ is a field 
and $ q = s^L \in \mathbb{K}^\times $ is not a root of unity. At this stage, some constructions and results work for more general parameters $ q $, 
but our interest is in the case that $ q $ is not a root of unity anyway.

\subsubsection{The parameter space $ \mathfrak{h}^*_q $} 
\label{sec:weights}

The notion of a (non-integral) weight for a representation of the quantum enveloping algebra $ U_q(\mathfrak{g}) $ is more subtle than in the classical case. 	To illustrate this, let us start with the case where the ground field is $\mathbb{K}=\CC$ and $q$ is a strictly positive real.  This will be the main case of interest in later chapters.  

Let $h\in\RR^\times$
\nomenclature{$h$}{element of $\RR$ (or $\mathbb{K}$) such that $e^h=q$}%
with $q=e^h$ and put $\hbar=h/2\pi$.
\nomenclature[o$h$2]{$\hbar$}{$=h/2\pi$}%
One says that a vector $v$ in a representation of $ U_q(\lie{g}) $ is a weight vector of weight $\lambda\in\mathfrak{h}^* = \Hom_{\mathbb{C}}(\mathfrak{h},\mathbb{C})$ if it is a common eigenvector for the action of $U_q(\lie{h})$ with
\[
K_\mu \cdot v = q^{(\lambda,\mu)} v, \qquad\qquad \text{for all }\mu\in\weights.
\]
Note, however, that the characters
\[
  U_q(\mathfrak{h}) \to \mathbb{C}; 
   \qquad   K_\mu \mapsto q^{(\lambda,\mu)}
\]
are not all distinct.  Specifically, $q^{(\lambda,\cdot)} \equiv 1$ as a function on $\weights$ whenever $\lambda \in \frac{2\pi i}{h} \roots^\vee = i\hbar^{-1}\roots^\vee$.  Therefore, the appropriate parameter space for weights is 
\[
\lie{h}^*_q  = \lie{h}^*/i\hbar^{-1}\roots^\vee.
\]
\label{nom:hq-star-C}	
	
For more general ground fields $\mathbb{K}$, weights should be understood as  algebra characters of the Cartan part $U_q(\mathfrak{h})$ of $U_q(\mathfrak{g})$.   The set of such characters is an abelian group with product dual to the coproduct on $U_q(\mathfrak{h})$, namely for any characters $\chi_1,\chi_2$ of $U_q(\mathfrak{h})$, 
\[
  \chi_1 \chi_2 (K) = (\chi_1 \otimes \chi_2) (\hat{\Delta}K) , 
    \qquad \text{for all } K\in U_q(\mathfrak{h}).
\]
This character group is isomorphic to $\Hom(\weights,\mathbb{K}^\times)$, upon observing that $U_q(\mathfrak{h})$ is the group algebra $\mathbb{K}[\weights]$.
In order to be consistent with the notation when $\mathbb{K} = \mathbb{C}$, we will use a formal $q$-exponential notation for elements of $\Hom(\weights,\mathbb{K}^\times)$, as follows.

\begin{definition}
	\label{def:hq-star}
	The parameter space $\mathfrak{h}_q^*$ is the character group of $U_q(\mathfrak{h})$, written additively.  For an element $\lambda\in \mathfrak{h}_q^*$ we denote the corresponding character  of $U_q(\mathfrak{h})$ by $\chi_\lambda$.  
	\nomenclature{$\chi_\lambda$}{character of $U_q(\lie{h})$ associated to $\lambda\in\lie{h}_q^*$}
	The associated elements of $\Hom(\weights,\mathbb{K}^\times)$ will be written as formal $q$-exponentials,
	\[
	   \weights \to \mathbb{K}^\times;
	   \qquad \mu \mapsto q^{(\lambda,\mu)},
	\]
	\nomenclature[o$q$4]{$q^{(\lambda,\cdot)}$}{formal notation for element of $\Hom(\weights,\mathbb{K}^\times)$ associated to $\lambda\in \lie{h}_q^*$}%
	so that $\chi_\lambda(K_\mu) = q^{(\lambda,\mu)}$ for all $\mu\in\weights$.
\end{definition}

To justify this notation, let us make some observations according to the properties of the base field $\mathbb{K}$.

Firstly, consider the case where $\mathbb{K}$ is a field satisfying no further assumptions.  Let $\lambda\in\weights$ be an integral weight.  Note that $(\lambda,\mu)\in \frac{1}{L}\mathbb{Z}$ for all $\mu\in\weights$, so it makes sense to define a character of $U_q(\mathfrak{h}_q)$ by
\[
  \chi_\lambda(K_\mu)= q^{(\lambda,\mu)} = s^{L(\lambda,\mu)}.
\]
The map $\lambda \mapsto \chi_\lambda$ is a homomorphism of $\weights$ into the characters of $U_q(\mathfrak{h})$, and we claim that it is injective.  Indeed, if $\lambda$ is in the kernel of this map, then we have $q^{(\lambda,\mu)} = 1$ for all $\mu\in\weights$.  Since $q$ is not a root of unity, this implies $\lambda=0$.  Therefore, we obtain an embedding
\[
 \weights \subset \mathfrak{h}_q^*
\]
which is compatible with the formal $q$-exponential notation $q^{(\lambda,\mu)}$ above.

Occasionally, we will consider the stronger condition that $q=s^{2L}$, in which case we have an embedding $\frac{1}{2}\weights \subset \mathfrak{h}_q^*$.

The other case of interest to us is that the field $ \mathbb{K} $ is exponential and the parameter $ q $ is obtained by exponentiating an 
element $ h \in \mathbb{K} $. Recall that an exponential field is a field $ \mathbb{K} $ of characteristic zero 
together with a group homomorphism
$ e^\bullet: \mathbb{K} \rightarrow \mathbb{K}^\times $, which we will denote by $ e^\bullet(h) = e^h $. 
We shall assume that $ q = e^h $ for some element $ h \in \mathbb{K} $, 
and write $ q^\bullet $ for the group homomorphism $ \mathbb{K} \rightarrow \mathbb{K}^\times $ given by $ q^\bullet(x) = e^{hx} $.  

In this case we set $ \mathfrak{h}^* = \mathbb{K} \otimes_\mathbb{Z} \weights $.  We can associate to any $\lambda \in \mathfrak{h}^*$ the character $\chi_\lambda$ given by
\[
\chi_\lambda (K_\mu) = q^{(\lambda,\mu)} =  e^{h(\lambda,\mu)}.
\]
Again, the map $\lambda \mapsto \chi_\lambda$ is a homomorphism from $\mathfrak{h}^*$ to the characters of $U_q(\mathfrak{h})$, and therefore induces a homomorphism
\[
 \mathfrak{h}^* \to \mathfrak{h}_q^*; \qquad \lambda \mapsto \chi_\lambda,
\]
which is compatible with the formal $q$-exponential notation.  Under the continuing assumption that $q$ is not a root of unity, the kernel of this map is $\ker(q^\bullet) \roots^\vee$, so we obtain an embedding
\[
 \mathfrak{h}^* / \ker(q^\bullet)\roots^\vee \subset \lie{h}_q^*.
\]

Of course, the prototypical example of an exponential field is $ \mathbb{K} = \mathbb{C} $ with the standard exponential function. In this case, $q^\bullet:\mathbb{C} \to \mathbb{C}^\times$ is surjective, with kernel
$\ker(q^\bullet) = (2\pi i/h)\mathbb{Z} = i \hbar^{-1}\mathbb{Z}$.
Using the fact that any character of $U_q(\mathfrak{h})$ is uniquely determined by its values on the generators $K_{\varpi_i}$ for $i=1,\ldots,n$, one can check that the map $\mathfrak{h}^* \to \mathfrak{h}_q^*$ is a surjection, and we have an isomorphism
\[
  \mathfrak{h}_q^* \cong \mathfrak{h}^* / i\hbar^{-1}\roots^\vee.
\]
In other words, Definition \ref{def:hq-star} is consistent with our definition of $\mathfrak{h}_q^*$ at the start of this subsection in the case $\mathbb{K}=\mathbb{C}$.

There is an action of the Weyl group $W$ on the quantum weight space $\mathfrak{h}_q^*$, induced by the action on $\Hom(\weights,\mathbb{K}^\times)$.  Specifically, if $w\in W$ and $\lambda\in\mathfrak{h}_q^*$ then $w\lambda$ is defined by
\[
 q^{(w\lambda,\mu)} = q^{(\lambda,w^{-1}\mu)}
\]
for all $\mu\in\weights$.
This extends the usual action of $W$ on $\weights \subset \mathfrak{h}_q^*$.

We can extend the usual ordering on $\roots$ to $\lie{h}^*_q$ as follows.  Recall that we write $\roots^+$
for the set of non-negative integer combinations of the simple roots.

\begin{definition}
	\label{def:hq-star-order}
	We define a relation $ \geq $ on $ \mathfrak{h}^*_q $ by saying that
	$$ 
	\lambda \geq \mu \qquad \text{if}\qquad \lambda - \mu \in \roots^+ .
	$$ 
	Here we are identifying $\roots^+$ with its image in $ \mathfrak{h}^*_q $. 
\end{definition}

Since $\roots^+$ is embedded in $\mathfrak{h}_q^*$, it is a simple matter to check that this is indeed a partial order.

\subsubsection{Weight modules and highest weight modules}

Let $ V $ be a left module over $ U_q(\mathfrak{g}) $. For any $ \lambda \in \mathfrak{h}_q^* $ we define the weight space 
$$
V_\lambda = \{v \in V \mid K_\mu \cdot v = q^{(\mu, \lambda)} v \; \text{for all} \; \mu \in \weights \}. 
$$
A vector $ v \in V $ is said to have weight $ \lambda $ iff $ v \in V_\lambda $. 
\nomenclature[s]{$V_\lambda$}{$\lambda$-weight space of a $U_q(\lie{g})$-module $V$}%

\begin{definition} \label{weightmodule}
A module $ V $ over $ U_q(\mathfrak{g}) $ is called a weight module if it is the direct sum of its weight spaces $ V_\lambda $ 
for $ \lambda \in \mathfrak{h}_q^* $.  We say that $\lambda$ is a weight of $V$ if $V_\lambda\neq 0$.
\end{definition} 

Every submodule of a weight module is again a weight module. This is a consequence of Artin's Theorem on linear independence 
of characters. 

\begin{definition}
	\label{def:integrable_module}
A $U_q(\mathfrak{g})$-module $V$ is called integrable if it is a weight module whose weights all belong to $\weights$, and
the operators $ E_i, F_j $ are locally nilpotent on $ V $ for all $ 1 \leq i,j \leq N $. 
\end{definition}

Such modules are sometimes referred to as type 1 modules.

\begin{definition} \label{defprimitive}
A vector $ v $ in a weight module $ V $ is called primitive if 
$$
E_i \cdot v = 0 \quad \text{for all} \; i = 1, \dots, N. 
$$
A module of highest weight $ \lambda \in \mathfrak{h}_q^* $ is a weight module $ V $ with a primitive cyclic vector $ v_\lambda \in V $ of 
weight $ \lambda $. 
\end{definition} 

If $ V $ is a weight module such that all weight spaces are finite dimensional, we define the restricted dual of $ V $ 
to be the $ U_q(\mathfrak{g}) $-module 
$$ 
V^\vee = \bigoplus_{\lambda \in \mathfrak{h}_q^*} \Hom(V_\lambda, \mathbb{K}), 
$$ 
\label{nom:restricted_dual1}%
with the left $ U_q(\mathfrak{g}) $-module structure given by 
$$
(X \cdot f)(v) = f(\tau(X) \cdot v).  
$$
Here $ \tau $ is the automorphism from Lemma \ref{deftau}. It is clear that $ V^\vee $ is again a weight module with $ (V^\vee)_\lambda = \Hom(V_\lambda, \mathbb{K}) $. Notice that we have a canonical isomorphism $ V^{\vee \vee} \cong V $ 
since $ \tau $ is involutive. 
If $ 0 \rightarrow K \rightarrow M \rightarrow Q \rightarrow 0 $ is an exact sequence of weight modules with finite dimensional 
weight spaces then the dual sequence $ 0 \rightarrow Q^\vee \rightarrow M^\vee \rightarrow K^\vee \rightarrow 0 $ is again exact.

\subsubsection{The definition of Verma modules} 

Let us now come to the definition of Verma modules. Recall that $ U_q(\mathfrak{b}_+)$ denotes the subalgebra of $ U_q(\mathfrak{g}) $ generated by 
the elements $ E_i $ and $ K_\lambda $.  The projection
\[
 U_q(\lie{b}_+) \cong  U_q(\lie{h}) \otimes U_q(\lie{n}_+) \stackrel{\id\otimes\epsilon}{\longrightarrow} U_q(\lie{h}) 
\]
which kills the generators $E_i$ is an algebra homomorphism, and we can use this to extend any character $ \chi: U_q(\mathfrak{h}) \rightarrow \mathbb{K} $ to a character of $ U_q(\mathfrak{b}_+)$ which sends the generators $ E_i $ to zero.
 We write again $ \chi $ for the resulting 
homomorphism $ U_q(\mathfrak{b}_+) \rightarrow \mathbb{K} $. 

Conversely, every algebra homomorphism $ U_q(\mathfrak{b}_+) \rightarrow \mathbb{K} $ vanishes on the elements $ E_i $, and therefore 
is determined by a homomorphism $ U_q(\mathfrak{h}) \rightarrow \mathbb{K} $. 
In particular, for each $\lambda\in\mathfrak{h}_q^*$ we have the character $\chi_\lambda$ of $U_q(\lie{b}_+)$ determined, as above, by
\[
 \chi_\lambda(K_\mu) = q^{(\lambda,\mu)}, \qquad \text{for all }\mu\in\weights.
\]

\begin{definition}
The Verma module $ M(\lambda) $ associated to $ \lambda \in \mathfrak{h}^*_q $ is the 
induced $ U_q(\mathfrak{g}) $-module 
$$
M(\lambda) = U_q(\mathfrak{g}) \otimes_{U_q(\mathfrak{b}_+)} \mathbb{K}_\lambda
$$
\nomenclature[o$K_\lambda$2]{$\mathbb{K}_\lambda$}{one-dimensional $U_q(\mathfrak{h})$-module or $U_q(\mathfrak{b})$-module}%
\label{nom:C_lambda}%
\nomenclature{$M(\lambda)$}{Verma module for $U_q(\lie{g})$}%
where $ \mathbb{K}_\lambda $ denotes the one-dimensional $ U_q(\mathfrak{b}_+) $-module $ \mathbb{K} $ with the action induced from the character $ \chi_\lambda $. The vector $ v_\lambda = 1 \otimes 1 \in U_q(\mathfrak{g}) \otimes_{U_q(\mathfrak{b}_+)} \mathbb{K}_\lambda $
\nomenclature{$v_\lambda$}{highest weight vector of the Verma module $M(\lambda)$ or its simple quotient $V(\lambda)$}
is called the highest weight vector of $ M(\lambda) $. 
\end{definition}

By construction, $ M(\lambda) $ for $ \lambda \in \mathfrak{h}_q^* $ is a highest 
weight module of highest weight $ \lambda $, and every other highest weight module of highest weight $ \lambda $ is a quotient of $ M(\lambda) $. 

Using Proposition \ref{uqgtriangular} one checks that $ M(\lambda) $ is free as a $ U_q(\mathfrak{n}_-) $-module. 
More specifically, the natural map $ U_q(\mathfrak{n}_-) \rightarrow M(\lambda) $ given by $Y \mapsto Y\cdot v_\lambda$ sends
$ U_q(\mathfrak{n}_-)_\mu $ bijectively on
to $ M(\lambda)_{\lambda - \mu} $. 

As in the classical case, the Verma module $ M(\lambda) $ contains a unique maximal submodule. 

\begin{lemma} \label{vermaquotient}
The Verma module $ M(\lambda) $ contains a unique maximal proper 
submodule $ I(\lambda) $,
\nomenclature{$I_\lambda$}{maximal proper submodule of the Verma module $M(\lambda)$}%
namely the linear span of all submodules not containing the highest weight 
vector $ v_\lambda $. 
\end{lemma} 
\begin{proof} Since every submodule $ U $ of $ M(\lambda) $ is a weight module, it is a proper submodule iff it does not contain the highest 
weight vector. That is, $ U $ is contained in the sum of all weight spaces different from $ M(\lambda)_{\lambda} $. 
In particular, the sum $ I(\lambda) $ of all proper submodules does not contain $ v_\lambda $, which means that $ I(\lambda) $ is the unique maximal proper 
submodule of $ M(\lambda) $. \end{proof}

The resulting simple quotient module $ M(\lambda)/I(\lambda) $ will be denoted $ V(\lambda) $.
\label{nom:simple_quotient}%
It is again a weight module, 
and all its weight spaces are finite dimensional. We may therefore form the restricted dual 
$ V(\lambda)^\vee $ of $ V(\lambda) $. Note that $ V(\lambda)^\vee $ is simple 
for any $ \lambda \in \mathfrak{h}_q^* $ because $ V(\lambda)^{\vee \vee} \cong V(\lambda) $ is simple. Since $ V(\lambda)^\vee $ is 
again a highest weight module of highest weight $ \lambda $ we conclude $ V(\lambda)^\vee \cong V(\lambda) $.

\subsection{Characters of $U_q(\lie{g})$}
\label{sec:Uqg_characters}

Unlike the classical case, the quantized enveloping algebras $U_q(\lie{g})$ admit several one-dimensional representations, with a corresponding multiplicity for all finite dimensional representations of $U_q(\lie{g})$.  We will only really be interested in the integrable (type $1$) finite dimensional modules, but at certain points we will need to acknowledge the existence of their non-integrable analogues.

\begin{definition}
 \label{def:Uqg_characters}
 We write $\mathbf{X}_q$ for the set of weights $\omega\in\lie{h}_q^*$ satisfying
 \[
  q^{(\omega,\alpha)} = \pm1, \qquad \text{for all } \alpha\in\roots.
 \]
 \nomenclature[o$X_q$]{$\mathbf{X}_q$}{set of weights $\omega\in\lie{h}_q^*$ satisfying $q^{(\omega,\alpha)} = \pm1$ for all $\alpha\in\roots$}%
\end{definition}

\begin{prop}
 \label{prop:Uqg_characters}
 The algebra characters of $U_q(\mathfrak{g})$ are in one-to-one correspondence with elements of $\mathbf{X}_q$.  Specifically, for every $\omega\in\mathbf{X}_q$, there is an algebra character
 $\chi_\omega : U_q(\lie{g}) \to \mathbb{K}$ defined on generators by
 \[
  \chi_\omega(K_\mu) = q^{(\omega,\mu)}, \quad
  \chi_\omega(E_i) = \chi_\omega(F_i) = 0,
 \]
 and every algebra character of $U_q(\lie{g})$ is of this form.
\end{prop}

\begin{proof}
 One can directly check that the formulas for $\chi_\omega$ respect the defining relations of $U_q(\mathfrak{g})$.  Conversely, suppose $\chi$ is an algebra character of $U_q(\lie{g})$.  Considering the restriction of $\chi$ to $U_q(\lie{h})$, there is $\omega\in\lie{h}_q^*$ with $\chi(K_\mu) = q^{(\omega,\mu)}$ for all $\mu\in\weights$.  Since
 \[
  0 = [\chi(E_i),\chi(F_i)] = \frac{q^{(\omega,\alpha_i)} - q^{-(\omega,\alpha_i)}}{q_i - q_i^{-1}},
 \]
 we have $q^{(\omega,\alpha_i)} = \pm1$ for all $i$, so $\omega \in \mathbf{X}_q$.  Finally, from the relation $K_\mu E_j K_\mu^{-1} = q^{(\alpha_j,\mu)}E_j$ and the fact that $q$ is not a root of unity, we deduce that $\chi(E_j)=0$, and similarly $\chi(F_j)=0$.
\end{proof}

\begin{definition}
 \label{def:non-type_1_weights}
 We define the extended integral weight lattice by
 \[
  \weights_q = \weights + \mathbf{X}_q.
 \]
 \nomenclature[o$P_q$]{$\weights_q$}{$=\weights+\mathbf{X}_q$}%
 We also put $\weights_q^+ = \weights^+ + \mathbf{X}_q$.
 \nomenclature[o$P_q^+$]{$\weights_q^+$}{$=\weights^++\mathbf{X}_q$}%
\end{definition}

Note that $\weights \cap \mathbf{X}_q = \{0\}$ by the assumption that $q$ is not a root of unity.

\subsection{Finite dimensional representations of $ U_q(\mathfrak{sl}(2, \mathbb{K})) $} 
\label{sec:fin_dim_sl2}

In this section we discuss the finite dimensional representation theory of $ U_q(\mathfrak{sl}(2, \mathbb{K})) $. 
We work over a field $ \mathbb{K} $ containing the non-root of unity $ q = s^2 \in \mathbb{K}^\times $, and we 
write $ s = q^{\frac{1}{2}} $ for simplicity.

When discussing $ U_q(\mathfrak{sl}(2, \mathbb{K})) $, we will use the following notation for the generators.  Let $\alpha$
\nomenclature{$\alpha$}{simple root of $\mathfrak{sl}(2,\mathbb{C})$, identified with $1\in\half\mathbb{Z}\cong\weights$}%
be the unique simple root and $\varpi = \half\alpha$
\nomenclature[o$\pi$2]{$\varpi$}{fundamental weight of $\mathfrak{sl}(2,\mathbb{C})$, identified with $\half\in\half\mathbb{Z}\cong\weights$}%
be the fundamental weight.  We write
\[
 E = E_1, \qquad F = F_1, \qquad K=K_\varpi, \qquad K^2=K_\alpha = K_1.
\]
We warn the reader that with this notation one must replace $K_i$ with $K^2$ when specializing formulas for $U_q(\lie{g})$ in terms of $E_i,F_i,K_i$ to formulas for $ U_{q_i}(\mathfrak{sl}(2, \mathbb{K})) $.

We will often identify $ \weights $ with $ \frac{1}{2} \mathbb{Z} $ in the sequel, so that $\varpi$ corresponds to $\frac12$ and $\alpha$ corresponds to $1$.  Note that $K$ acts on any vector $v$ of weight $m\in\frac12\mathbb{Z}$ by
\[
 K\cdot v = q^m v.
\]

It is not hard to check that the formulas 

\begin{align*}
\pi_{\frac{1}{2}}(E) = 
\begin{pmatrix}
0 & q^{-\frac{1}{2}} \\ 
0 & 0 
\end{pmatrix}, \quad
\pi_{\frac{1}{2}}(K) = 
\begin{pmatrix}
q^{\frac{1}{2}} & 0 \\ 
0 & q^{-\frac{1}{2}} 
\end{pmatrix}, \quad
\pi_{\frac{1}{2}}(F) = 
\begin{pmatrix}
0 & 0 \\ 
q^{\frac{1}{2}} & 0 
\end{pmatrix}
\end{align*}
define a $ 2 $-dimensional irreducible representation $ \pi_{\frac{1}{2}}: U_q(\mathfrak{sl}(2, \mathbb{K})) \rightarrow M_2(\mathbb{K}) $.
We will write $ V(\frac12) $ for this module.  This is consistent with our notation for simple quotients of Verma modules, see Proposition \ref{sl2irreps} below.

In order to construct further simple modules we shall use
the following specific instance of Lemma \ref{EFpowercommutationhelp}.  We recall the notation
\[
 [K^2; m] = \frac{q^m K^2 - q^{-m} K^{-2}}{q - q^{-1}}
\]
for $ m \in \mathbb{Z} $.

\begin{lemma} \label{EFbasic}
In $ U_q(\mathfrak{sl}(2, \mathbb{K})) $ we have 
$$
[E, F^{m + 1}] = \frac{1}{q - q^{-1}} [m + 1]_q F^m (q^{-m} K^2 - q^m K^{-2}) = [m + 1]_q F^m [K^2; -m]  
$$
for all $ m \in \mathbb{N}_0 $.  
\end{lemma}

\begin{prop} \label{sl2irreps}
	Let $\lambda\in\lie{h}_q^*$.  
	\begin{bnum}
		\item[a)] The Verma module $M(\lambda)$ is simple if and only if $\lambda \notin \weights_q^+$.
		\item[b)] If $\lambda\in\weights_q^+$, so that $\lambda=n\alpha+\omega$  with $n\in\half\mathbb{N}_0$, $\omega \in \mathbf{X}_q$ then $M(\lambda)$ has a unique simple quotient $V(\lambda)$ of highest weight $\lambda$ and dimension $2n+1$.  Up to isomorphism, these are the only finite dimensional simple weight modules of $U_q(\mathfrak{sl}(2,\mathbb{K}))$.
		\item[c)] If $\lambda\in\weights^+$ then $V(\lambda)$ is integrable, and these are the only simple integrable $U_q(\mathfrak{sl}(2,\mathbb{K}))$-modules up to isomorphism.
	\end{bnum}
\end{prop}

\begin{proof}
 \emph{a)} The Verma module $M(\lambda)$ is non-simple if and only if it contains a primitive vector of the form $F^{m+1}\cdot v_\lambda$ for some $m \geq 0$.  According to Lemma \ref{EFbasic}, this occurs if and only if $q^{-m} q^{(\lambda,\alpha)} - q^m q^{-(\lambda,\alpha)} = 0$, or equivalently $q^{(\lambda-m\varpi,\alpha)} = \pm1$.  Putting $\omega=\lambda - m\varpi$ we have $\omega \in\mathbf{X}_q$ and so $\lambda = m\varpi+\omega \in \weights_q^+$. 
 
 \emph{b)} In the case $\lambda = m\varpi+\omega \in \weights_q^+$, the simple quotient is spanned by $\{F^i\cdot v_\lambda \mid i=0,\ldots m\}$.  Putting $m=2n$  we obtain the dimension formula.
 Moreover, any finite dimensional simple $ U_q(\mathfrak{sl}(2, \mathbb{K})) $-module is necessarily a highest weight module, and so a quotient of some Verma module $M(\lambda)$. 
 
 \emph{c)} Likewise, any simple integrable module is necessarily a highest weight module with highest weight in $\weights$.  Since $F$ must act locally nilpotently, \emph{a)} implies that $\lambda\in\weights^+$. 
\end{proof}

 In the sequel we will write $V(n)$ for the integrable module $V(n\alpha)$ with $n\in\half\mathbb{N}_0$. 
Let us give an explicit description of $ V(n) $. 

\begin{lemma} \label{vnexplicit}
Let $ n \in \frac{1}{2} \mathbb{N}_0 $. The module $ V(n) $ over $ U_q(\mathfrak{sl}(2, \mathbb{K})) $ has a $ \mathbb{K} $-linear 
basis $ v_n, v_{n - 1}, \dots, v_{-n + 1}, v_{-n} $ such that 
\begin{align*}
K \cdot v_j &= q^j v_j, \\
F \cdot v_j &= v_{j - 1}, \\
E \cdot v_j &= [n - j] [n + j + 1] v_{j + 1}. 
\end{align*}
Here we interpret $ v_k = 0 $ if $ k $ is not contained in the set $ \{n, n - 1, \dots, -n\} $. If we rescale this basis as
$$ 
v_{(j)} = \frac{1}{[n - j]!} v_j 
$$ 
we have 
\begin{align*}
K \cdot v_{(j)} &= q^j v_{(j)}, \\ 
F \cdot v_{(j)} &= [n - j + 1] v_{(j - 1)}, \\
E \cdot v_{(j)} &= [n + j + 1] v_{(j + 1)}. 
\end{align*}
\end{lemma} 
\begin{proof} This is an explicit description of the basis obtained from applying the operator $ F $ to the highest 
weight vector $ v_n $. Since $ FK = q KF $ we obtain the action of $ K $ on the vectors $ v_j $ by induction, for the inductive step use 
$$ 
K \cdot v_{j - 1} = K F \cdot v_j = q^{-1} FK \cdot v_j = q^{j - 1} F \cdot v_j = q^{j - 1} v_{j - 1}. 
$$
Moreover, Lemma \ref{EFbasic} yields 
\begin{align*}
E F^{r + 1} \cdot v_n &= [E, F^{r + 1}] \cdot v_n = [r + 1][2n - r] F^r \cdot v_n 
\end{align*}
for all $ r \in \mathbb{N}_0 $, hence setting $ r = n - j - 1 $ we get 
\begin{align*}
E \cdot v_j = E F^{n - j} \cdot v_n &= [n - j][n + j + 1] F^{n - j - 1} \cdot v_n \\
&= [n - j][n + j + 1] v_{j + 1} 
\end{align*}
The second set of equations follows immediately from these formulas, indeed, we have 
\begin{align*}
E \cdot v_{(j)} &= \frac{1}{[n - j]!} E \cdot v_j = \frac{[n + j + 1]}{[n - j - 1]!} v_{j + 1} = [n + j + 1] v_{(j + 1)}, \\ 
F \cdot v_{(j)} &= \frac{1}{[n - j]!} F \cdot v_j = \frac{[n - j + 1]}{[n - j + 1]!} v_{j - 1} = [n - j + 1] v_{(j - 1)} 
\end{align*}
as desired. \end{proof}  

We point out that the labels of the basis vectors in Lemma \ref{vnexplicit} run over half-integers if $ n $ is a half-integer, and 
over integers if $ n \in \mathbb{N}_0 $.

The non-integrable simple modules are obtained by twisting the integrable simple  modules by a character.  Specifically, if $n\in\half\mathbb{N}$ and $\omega\in\mathbf{X}_q$ then
\[
 V(n\alpha  +\omega) \cong V(n) \otimes V(\omega),
\]
since both are simple modules of highest weight $n\alpha + \omega$.  More generally, we have the following.

\begin{lemma}
	\label{lem:reduction_to_integrable_for_sl2}
	All weights of a finite dimensional weight module $V$ for $U_q(\mathfrak{sl}(2,\mathbb{K}))$ belong to $\weights_q$.  Moreover, $V$ admits a direct sum decomposition
	\[
	  V \cong \bigoplus_{\omega\in\mathbf{X}_q} W^\omega \otimes V(\omega)
	\]
	where the $W^\omega$ are finite dimensional integrable $U_q(\mathfrak{sl}(2,\mathbb{K}))$-modules. 
\end{lemma}

\begin{proof}
	Let $\lambda\in\lie{h}_q^*$ be a weight of $V$.  If we fix $v\in V_\lambda$ nonzero then $E^k\cdot v$ is primitive vector for some $k\in\mathbb{N}_0$.  By Proposition \ref{sl2irreps} we must have $\lambda+k\alpha \in \weights_q^+$ and so $\lambda\in\weights_q$.  
	
	Now fix $\omega\in\mathbf{X}_q$.  The subspace
	\[
	V_{\weights + \omega} = \bigoplus_{\mu\in\weights} V_{\mu+\omega}
	\]
	is invariant under the $U_q(\lie{g})$-action.  We obtain a direct sum decomposition $V=\bigoplus_{\omega\in\mathbf{X}_q} V_{\weights+\omega}$.  Putting $W^\omega = V_{\weights+\omega}\otimes V(-\omega)$ yields the result.
\end{proof}

We can now prove complete reducibility for finite dimensional weight modules of $U_q(\mathfrak{sl}(2, \mathbb{K})) $.

\begin{prop} \label{slq2completeirreducibility}
Every finite dimensional weight module of $ U_q(\mathfrak{sl}(2, \mathbb{K})) $ is completely reducible. 
\end{prop} 

\begin{proof} 
By Lemma \ref{lem:reduction_to_integrable_for_sl2}, it suffices to proof the proposition for finite dimensional integrable modules $V$.  In this case, the weights of $V$ belong to $\weights$.  

Let $\mu\in\weights$ be maximal among the weights of $V$.
Any vector $ v \in V_\mu $ generates a finite dimensional $ U_q(\mathfrak{sl}(2, \mathbb{K})) $-module isomorphic to $ V(\mu) $. 
Since the intersection of this submodule with $ V_\mu $ is $ \mathbb{K} v $, one checks that $ V $ contains a direct sum of $ \dim(V_\mu) $ 
copies of $ V(\mu) $ as a submodule.  

Writing $ K $ for this submodule and $ Q = V/K $ for the corresponding quotient, we obtain a short exact sequence 
$ 0 \rightarrow K \rightarrow V \rightarrow Q \rightarrow 0 $ 
of finite dimensional weight modules. Using induction on dimension we may assume that $ Q $ is a direct sum 
of simple modules. Moreover the dual sequence 
$ 0 \rightarrow Q^\vee \rightarrow V^\vee \rightarrow K^\vee \rightarrow 0 $, 
see the constructions after Definition \ref{weightmodule}, is split by construction because the highest weight subspace of $ V^\vee $ maps 
isomorphically onto the highest weight subspace of $ K^\vee $. Applying duality again shows that the original sequence 
is split exact, which means that $ V $ is a direct sum of simple modules. 
\end{proof}

\subsection{Finite dimensional representations of $ U_q(\mathfrak{g}) $}

In this section we begin the study of finite dimensional representations of $U_q(\mathfrak{g})$.  Ultimately, in Section \ref{sec:Casimir}, we will prove complete reducibility of finite dimensional weight modules and obtain a classification the finite dimensional simple modules, although those results require considerably more machinery than in the case $ \mathfrak{g} = \mathfrak{sl}(2, \mathbb{K}) $ discussed above.  For more information we refer to \cite{Jantzenqg}, \cite{Lusztigbook}. 

Throughout this section we assume that $ \mathfrak{g} $ is a semisimple Lie algebra and $ q = s^L \in \mathbb{K}^\times $ is not a root of unity.

\subsubsection{Rank-one quantum subgroups}

For each $1 \leq i \leq N$, we write $ U_{q_i}(\mathfrak{g}_i) \subset U_q(\mathfrak{g}) $
\nomenclature[o$g_i$]{$\lie{g}_i$}{rank-one Lie subalgebra of $\lie{g}$ generated by $E_i$, $F_i$}%
\nomenclature{$U_{q_i}(\mathfrak{g}_i)$}{subalgebra of $U_q(\mathfrak{g})$ generated by $ E_i, F_i, K_i^{\pm 1} $}%
for the subalgebra 
generated by $ E_i, F_i, K_i^{\pm 1} $.  It is isomorphic as a Hopf algebra to the subalgebra of $U_{q_i}(\mathfrak{sl}(2,\mathbb{K}))$ generated by $E,F,K^{\pm2}$. 

The classification of finite dimensional simple modules for $U_{q_i}(\lie{g}_i)$ is essentially identical to that of $U_{q_i}(\mathfrak{sl}(2,\mathbb{K}))$: such modules are indexed by highest weights of the form $\mu=n\alpha_i + \omega$ where $n\in\half\mathbb{N}_0$ and $q_i^{(\omega,\cdot)}$ denotes a character of $\ZZ\alpha_i$ with values in $\{\pm1\}$.

Note that if $v\in V_\lambda$ is a vector of weight $\lambda \in \lie{h}_q^*$ for some $U_q(\mathfrak{g})$-module $V$, then $K_i$ acts on $v$ by
\[
 K_i\cdot v = q_i^{(\lambda,\alpha_i^\vee)} v.
\]
In particular, if $V$ is an integrable $U_q(\mathfrak{g})$-module then it becomes an integrable $U_q(\mathfrak{g}_i)$-module by restriction, wherein vectors in $V_\lambda$ have weight
\[ 
 \textstyle\half(\lambda,\alpha_i^\vee) \in \half\ZZ
\]
under the usual identification of half-integers with integral weights for $U_{q_i}(\mathfrak{sl}(2,\mathbb{K}))$.

\subsubsection{Finite dimensional modules}

As a consequence of the above remarks, we obtain some basic structural results for finite dimensional $U_q(\lie{g})$-modules.  
We begin with the analogue of Lemma \ref{lem:reduction_to_integrable_for_sl2}, which will let us restrict our attention to integrable modules.

\begin{lemma}
 \label{lem:reduction_to_integrable}
 All the weights of a finite dimensional weight module $V$ for $U_q(\lie{g})$ belong to $\weights_q$.  Moreover, $V$ admits a direct sum decomposition 
 \[
  V \cong \bigoplus_{\omega\in\mathbf{X}_q} W^\omega \otimes V(\omega)
 \]
 where each $W^\omega$ is a finite dimensional integrable $U_q(\lie{g})$-module and $V(\omega)$ is the one dimensional representation with weight $\omega\in\mathbf{X}_q$. 
\end{lemma}

\begin{proof}
 Let $\lambda\in \lie{h}^*_q$ and $v\in V_\lambda$ be a nonzero weight vector.  For each $i=1,\ldots,N$, we have
 \[
  K_i\cdot v = q^{(\lambda,\alpha_i)}v.
 \]
 Considering $V$ as a $U_{q_i}(\mathfrak{sl}(2,\mathbb{K}))$-module by restriction, Lemma \ref{lem:reduction_to_integrable_for_sl2} implies that $q^{(\lambda,\alpha_i)} \in \pm q^\ZZ$.  It follows that $\lambda\in\weights_q$.
 The direct sum decomposition then follows exactly as in the proof of Lemma \ref{lem:reduction_to_integrable_for_sl2}. 
\end{proof}

\begin{prop}
 \label{prop:weights_Weyl_invariant}
 The set of weights of any integrable $U_q(\mathfrak{g})$-module $V$, and their multiplicities, are invariant under the Weyl group.
\end{prop}

\begin{proof}
 Let $1\leq i\leq N$.
 By Proposition \ref{slq2completeirreducibility} $V$ decomposes as a direct sum of irreducible integrable $U_{q_i}(\mathfrak{g}_i)$-modules.  It follows that the set of weights of $V$ and their multiplicities are invariant under the simple reflection $s_i$.  Since these generate $W$, the result follows.
\end{proof}

Let us fix a dominant integral weight $ \mu \in \weights^+ $ and $ 1 \leq i \leq N $. Then $ (\alpha^\vee_i, \mu) \in \mathbb{N}_0 $, 
and we claim that $ v = F_i^{(\alpha_i^\vee, \mu) + 1} \cdot v_\mu $ 
is a primitive vector in $ M(\mu) $. Indeed, Lemma \ref{EFpowercommutationhelp} shows that 
$ E_i \cdot v = 0 $, and $ E_j \cdot v = 0 $ for $ j \neq i $ follows from the fact that $ [E_j, F_i] = 0 $ 
and the primitivity of $ v_\lambda $. The vector $ v $ has weight $ s_i . \mu $, where $ s_i \in W $ denotes 
the simple reflection corresponding to $ \alpha_i $ and 
$$
w . \lambda = w (\lambda + \rho) - \rho
$$
denotes the shifted Weyl group action. Hence $ v $ generates a homomorphic image of 
$ M(s_i . \mu) $ inside $ M(\mu) $. By slight abuse of notation we write $ M(s_i . \mu) \subset M(\mu) $ 
for this submodule, since it will be shown later that the homomorphism $ M(s_i . \mu) \rightarrow M(\mu) $ is indeed injective.

We are now ready to describe the finite dimensional integrable quotients of Verma modules, compare Section 5.9. in \cite{Jantzenqg}.  

\begin{theorem} \label{thmfdintegrable}
Let $ \lambda \in \weights^+ $. Then the largest integrable quotient $ L(\lambda) $ of $ M(\lambda) $ is finite dimensional and given by 
$$ 
L(\lambda) = M(\lambda)\bigg/\sum_{i = 1}^N M(s_i . \lambda). 
$$ 
\nomenclature{$L(\lambda)$}{integrable quotient of the Verma module $M(\lambda)$ with $\lambda\in\weights^+$}%
If $ \lambda \in \mathfrak{h}^*_q \setminus \weights^+ $ then the Verma module $ M(\lambda) $ does not admit any nontrivial integrable quotients. 
\end{theorem} 

\begin{proof} 
It is clear that any integrable quotient of $ M(\lambda) $ must annihilate the sum of the modules $ M(s_i . \lambda) $, 
because otherwise the action of $ F_i $ on $ v_\lambda $ fails to be locally finite. 

In order to show that the action of $ F_i $ on any $ v \in L(\lambda) $ is locally nilpotent
we use that $ v $ can be written as a sum of terms $ F \cdot v_\lambda $ where $ F $ is a monomial in $ F_1, \dots, F_N $. 
Let us show by induction on the degree $ r $ of $ F = F_{i_1} \cdots F_{i_r} $ that $ F_i^k \cdot F \cdot v_\lambda = 0 $ for some $ k \in \mathbb{N} $. 
For $ r = 0 $ the claim is obvious from the definition of $ L(\lambda) $. Assume now $ F = F_j Y \cdot v_\lambda $ where $ Y $ is a monomial of 
degree $ r $ and $ 1 \leq j \leq N $. If $ j = i $ then the claim follows from our inductive hypothesis, so let us assume $ j \neq i $. 
Writing $ u = Y \cdot v_\lambda $ we have 
$$ 
F_i \cdot (F_j \cdot u) = (F_i \rightarrow F_j) \cdot u + K_i^{-1} F_j K_i F_i \cdot u. 
$$ 
By applying this relation iteratively we obtain $ F_i^k \cdot (F_j \cdot u) = 0 $ in $ L(\lambda) $ for a suitable $ k \in \mathbb{N} $, using the 
inductive hypothesis and the fact that the Serre relations imply $ F_i^{1 - a_{ij}} \rightarrow F_j = 0 $. 
It follows that $ L(\lambda) $ is integrable, and it is indeed the largest integrable quotient of $ M(\lambda) $. 

By Proposition \ref{prop:weights_Weyl_invariant}, 
 the set of weights of $ L(\lambda) $ is 
invariant under the Weyl group action. Since each weight space is finite dimensional this implies that $ L(\lambda) $ is finite dimensional.

Finally, suppose $M(\lambda)$ admits a nontrivial integrable quotient for some $\lambda \in \lie{h}_q^* \setminus \weights^+$.  The definition of integrability implies $\lambda\in\weights$.  But since $\lambda\notin\weights^+$,  
we have $ (\alpha_i^\vee, \lambda) \notin \mathbb{N}_0 $ for some $ i $. 
It follows that the $ U_{q_i}(\mathfrak{g}_i) $-module generated by $ v_\lambda $ is infinite dimensional and simple by Proposition \ref{sl2irreps}, and hence 
$F_i$ does not act locally nilpotently on any quotient of $M(\lambda)$. 
\end{proof} 

Note that we will ultimately observe, in Theorem \ref{thmfdirreducible}, that when $\mu\in\weights$, the finite dimensional quotients $L(\mu)$ are irreducible, and hence $L(\mu) = V(\mu)$.

\begin{cor}
 The Verma module $M(\lambda)$ admits a finite dimensional quotient if and only if $\lambda\in\weights_q^+$.
\end{cor}

\begin{proof}
 If $M(\lambda)$ admits a nontrivial finite dimensional quotient then Lemma \ref{lem:reduction_to_integrable} implies $\lambda\in\weights_q$, so $\lambda = \mu+\omega$ for some $\mu\in\weights$, $\omega\in\mathbf{X}_q$.  It is easy to check that $M(\mu) \cong M(\lambda)\otimes V(-\omega)$.  It follows that $M(\mu)$ has a nontrivial finite dimensional hence integrable quotient.  Theorem \ref{thmfdintegrable} implies $\mu\in\weights^+$, so $\lambda\in\weights_q^+$.
 
 Conversely, if $\lambda = \mu+\omega \in\weights_q^+$ then $M(\lambda)\cong M(\mu)\otimes V(\omega)$ admits the finite dimensional quotient $L(\lambda)\otimes V(\omega)$.
\end{proof}

The following technical lemma says roughly that in weight spaces sufficiently close to the highest weight, the finite dimensional module $L(\mu)$ resembles the Verma module $M(\mu)$.

\begin{lemma}
 \label{lem:L_looks_like_Uqn}
 Let $\mu=\sum_i m_i \varpi_i\in \weights^+$.  If $\lambda = \sum_i l_i \alpha_i \in \roots^+$ with $l_i\leq m_i$ for all $i=1,\ldots, N$ then $ L(\mu)_{\mu - \lambda} \cong M(\mu)_{\mu - \lambda} \cong U_q(\mathfrak{n}_-)_{-\lambda} $ as vector spaces, via the map which sends $F\in U_q(\mathfrak{n}_-)_{-\lambda} $ to $F\cdot v_\mu \in L(\mu)_{\mu - \lambda}$.
\end{lemma}

\begin{proof}
 Note that $s_i.\mu = \mu - m_i\alpha_i$.  Therefore, with the given conditions on $\lambda$ we have $\mu-\lambda \not\leq s_i.\mu$ for every $i=1,\ldots,N$, and the first claim follows from Theorem \ref{thmfdintegrable}.  The second claim then follows from the definition of $\omega$.
\end{proof}

If $\nu = \sum_i n_i\varpi_i \in \weights^+$, we will use the notation $L(-\nu)$
\nomenclature[o$L(\nu)$]{$L(-\nu)$}{finite dimensional integrable module with lowest weight $\nu$, where $\nu\in\weights^+$}%
to denote the finite dimensional module $L(\nu)$ twisted  by the automorphism $ \omega $ from Lemma \ref{defOmega}, so that the action of $ X \in U_q(\mathfrak{g}) $ on $ v \in L(-\nu) = L(\nu) $ is given by $ \omega(X) \cdot v $. In this case we have $ L(-\nu)_{-\nu + \lambda} \cong U_q(\mathfrak{n}_+)_{\lambda} $ provided $ l_i \leq n_i $ for all $ i = 1, \dots, N $.

An important fact is that finite dimensional integrable representations 
separate the points of $ U_q(\mathfrak{g}) $, compare for instance Section 7.1.5 of \cite{KS}. 

\begin{theorem} \label{uqgfdseparation}
The finite dimensional integrable  representations of $ U_q(\mathfrak{g}) $ separate points. That is, 
if $ X \in U_q(\mathfrak{g}) $ satisfies $ \pi(X) = 0 $ for all finite dimensional integrable representations $ \pi: U_q(\mathfrak{g}) \rightarrow \End(V) $, 
then $ X = 0 $.  
\end{theorem} 

\begin{proof}
Let us write $ v_\mu $ for the highest weight vector in $ L(\mu) $ and $ v_{-\nu} $ for the lowest weight vector in $ L(-\nu) $. 
Assume $ X \neq 0 $ and write this element as a finite sum
$$
X = \sum_{i, \eta, j} f_{ij} K_\eta e_{ij}, 
$$
where $ f_{ij} \in U_q(\mathfrak{n}_-)_{-\beta_i}, e_{ij} \in U_q(\mathfrak{n}_+)_{\gamma_j} $ 
for pairwise distinct weights $ \beta_1, \dots, \beta_m $ and $ \gamma_1, \dots, \gamma_n $, and $ \eta \in \weights $. 
We shall show that $ X $ acts in a nonzero fashion on $ L(-\nu) \otimes L(\mu) $ for suitable $ \mu, \nu \in \weights^+ $. 

Assume this is not the case, so that $ X $ acts by zero on all $ L(-\nu) \otimes L(\mu) $ for $ \mu, \nu \in \weights^+ $. 
Let $ \gamma_t $ be maximal among the $ \gamma_j $ such that some $ f_{it} K_\eta e_{it} $ is nonzero. 
For each $ i, j $ we have $ \hat{\Delta}(e_{ij}) = e_{ij} \otimes K_{\gamma_j} + b $ where 
$ b $ is a sum of terms in $ U_q(\mathfrak{n}_+) \otimes (U_q(\mathfrak{b}_+) \cap \ker(\hat{\epsilon})) $ which 
vanish on $ v_{-\nu} \otimes v_\mu $. 
We thus get 
$$
e_{ij} \cdot (v_{-\nu} \otimes v_\mu) = q^{(\gamma_j, \mu)} e_{ij} \cdot v_{-\nu} \otimes v_\mu.  
$$
Moreover, since $ \hat{\Delta}(K_\eta) = K_\eta \otimes K_\eta $ we obtain
$$
K_\eta e_{ij} \cdot (v_{-\nu} \otimes v_\mu) = q^{(\eta, \mu - \nu + \gamma_j)} q^{(\gamma_j, \mu)} e_{ij} \cdot v_{-\nu} \otimes v_\mu. 
$$ 
Finally, using that $ \hat{\Delta}(f_{ij}) = K_{-\beta_i} \otimes f_{ij} + c $ 
where $ c \in (U_q(\mathfrak{b}_-) \cap \ker(\hat{\epsilon})) \otimes (U_q(\mathfrak{n}_-) $ 
we get 
$$
f_{ij} K_\eta e_{ij} \cdot (v_{-\nu} \otimes v_\mu) = q^{(\beta_i, \nu - \gamma_j)} 
q^{(\eta, \mu - \nu + \gamma_j)} q^{(\gamma_j, \mu)} e_{ij} \cdot v_{-\nu} \otimes f_{ij} \cdot v_\mu + r, 
$$ 
where $ r $ is a sum of terms in $ L(-\nu)_{-\nu + \delta} \otimes L(\mu) $ with $ \delta \neq \gamma_t $. 
It follows that the component of $ X \cdot (v_{-\nu} \otimes v_\mu) $ contained in $ L(-\nu)_{-\nu + \gamma_t} \otimes L(\mu) $ is 
$$
0 = \sum_{i, \eta} q^{(\beta_i, \nu - \gamma_t)} 
q^{(\eta, \mu - \nu + \gamma_t)} q^{(\gamma_t, \mu)} e_{it} \cdot v_{-\nu} \otimes f_{it} \cdot v_\mu.  
$$
Assume without loss of generality that all $ e_{it}, f_{it} $ for $ 1 \leq i \leq m $ are nonzero, and 
using Lemma \ref{lem:L_looks_like_Uqn}, 
fix $ \nu $ large enough such 
that $ e_{it} \cdot v_{-\nu} $ is nonzero for all $ i $. 
Choosing $ \mu $ large enough we get that the vectors $ e_{it} \cdot v_{-\nu} \otimes f_{it} \cdot v_\mu $ are linearly independent. 
Hence 
$$
0 = \sum_{\eta} q^{(\beta_i, \nu - \gamma_t)} q^{(\eta, \mu - \nu + \gamma_t)} q^{(\gamma_t, \mu)}
$$
for all $ i $ in this case. Cancelling the common nonzero factors $ q^{(\beta_i, \nu - \gamma_t)} q^{(\gamma_t, \mu)} $ and 
writing $ c_\eta = q^{(\eta, \gamma_t - \nu)} $ we deduce 
$$
0 = \sum_{\eta} q^{(\eta, \mu - \nu + \gamma_t)} = \sum_\eta c_\eta q^{(\eta, \mu)} = \sum_\eta c_\eta \chi_\eta(K_\mu).  
$$
Consider the subsemigroup $ P \subset \weights^+ $ of all $ \mu \in \weights^+ $ 
such that the vectors $ e_{it} \cdot v_{-\nu} \otimes f_{it} \cdot v_\mu $ are linearly independent. 
Since the characters $ \chi_\eta $ on $ P $ 
are pairwise distinct,
Artin's Theorem on the linear independence of characters 
implies $ q^{(\eta, \gamma_t - \nu)} = c_\eta = 0 $ for all $ \eta $, which is clearly impossible. 

Hence our initial assumption that $ X $ acts by zero on all modules of the form $ V(-\nu) \otimes V(\mu) $ was wrong. This finishes the proof. \end{proof}

\subsection{Braid group action and PBW basis} \label{secbraid}

In this section we explain the construction of the PBW basis for $ U_q(\mathfrak{g}) $ and its integral form. 
A detailed exposition is given in part VI of \cite{Lusztigbook} and chapter 8 of \cite{Jantzenqg}; we shall follow the discussion in \cite{Jantzenqg}. 
In order to define the PBW basis one first constructs an action of the braid 
group $ B_\mathfrak{g} $ of $ \mathfrak{g} $ on $ U_q(\mathfrak{g}) $ and its integrable modules. 
Throughout this section we assume that $ q = s^L \in \mathbb{K}^\times $ is not a root of unity. 

Let $ \mathfrak{g} $ be a semisimple Lie algebra of rank $ N $. The braid group $ B_\mathfrak{g} $
\nomenclature{$B_{\mathfrak{g}}$}{braid group}%
is obtained from the Cartan 
matrix $ (a_{ij}) $ of $ \mathfrak{g} $ as follows, see \cite{CPbook} Section 8.1. 
Let $ m_{ij} $
\nomenclature{$m_{ij}$}{index $2,3,4,6$ according to $a_{ij}a_{ji}=0,1,2,3$, respectively}%
be equal to $ 2,3,4,6 $ iff $ a_{ij} a_{ji} $ equals $ 0,1,2,3 $, respectively. 
By definition, $ B_\mathfrak{g} $ is the group with generators $ T_j $
\nomenclature{$T_j$}{generator of the braid group $B_{\mathfrak{g}}$}%
for $ j = 1, \dots, N $ and relations  
$$
T_i T_j T_i \cdots = T_j T_i T_j \cdots 
$$
for all $ 1 \leq i, j \leq N $ with $ i \neq j $, where there are $ m_{ij} $ terms on each side of the equation.

For $ \mathfrak{g} $ of type $ A $ one obtains the classical braid groups in this way. 
In general, there is a canonical quotient homomorphism $ B_\mathfrak{g} \rightarrow W $, sending the generators $ T_j $ to 
the simple reflections $ s_j $. Indeed, the Weyl group $ W $ is the quotient of $ B_\mathfrak{g} $ by the relations $ T_j^2 = 1 $ 
for $ j = 1, \dots, N $. 

\subsubsection{The case of $ \mathfrak{sl}(2, \mathbb{K}) $} 

We use the notation for $U_q(\mathfrak{sl}(2,\mathbb{K}))$ from Section \ref{sec:fin_dim_sl2}, in particular we write $K=K_\varpi$ and $K^2=K_\alpha$.
Let $ V $ be an integrable module over $ U_q(\mathfrak{sl}(2, \mathbb{K})) $. 
We label the weights by $ \frac{1}{2} \mathbb{Z} $, so that $ v \in V $ is of weight $ m $ if $ K \cdot v = q^m v $. 
The corresponding weight space is denoted by $ V_m $. 

We define operators $ \T_\pm: V \rightarrow V $ by 
\begin{align*}
\T_+(v) &= \sum_{\substack{r,s,t \geq 0 \\ r - s + t = 2m}} (-1)^s q^{s - rt} F^{(r)} E^{(s)} F^{(t)} \cdot v \\
\T_-(v) &= \sum_{\substack{r,s,t \geq 0 \\ r - s + t = 2m}} (-1)^s q^{rt - s} F^{(r)} E^{(s)} F^{(t)} \cdot v
\end{align*} 
\label{nom:T_pm}%
for $ v \in V_m $, where $ E^{(i)} = E^i/[i]! $ and $ F^{(j)} = F^j/[j]! $ are the divided powers. 
Note that these operators are well-defined because the action of $ E $ and $ F $ on $ V $ is locally finite. 
It is also clear that $ \T_\pm $ maps $ V_m $ to $ V_{-m} $ for any weight $ m $.

Recall the basis $ v_{(n)}, v_{(n - 1)}, \dots, v_{(-n)} $ of $ V(n) $ from Lemma \ref{vnexplicit}. 
According to Lemma \ref{vnexplicit} we obtain the formulas 
\begin{align*}
E^{(r)} \cdot v_{(j)} &= \frac{[n + j + 1][n + j + 2] \cdots [n + j + r]}{[r]!}  \cdot v_{(j)} = 
\begin{bmatrix}
n + j + r \\ 
r
\end{bmatrix}
v_{(j + r)} \\
F^{(r)} \cdot v_{(j)} &= \frac{[n - j + 1][n - j + 2] \cdots [n - j + r]}{[r]!}  \cdot v_{(j)} = 
\begin{bmatrix}
n - j + r \\ 
r
\end{bmatrix}
v_{(j - r)} 
\end{align*} 
for the action of divided powers on these basis vectors, where we adopt the convention $ v_{(n + 1)} = 0 = v_{(-n - 1)} $. 

Using these formulas we shall derive explicit expressions for the action of the operators $ \T_\pm $ on the basis vectors $ v_{(j)} $. 
As a preparation we state some properties of $ q $-binomial coefficients. 

\begin{lemma} \label{tlemma1}
Let $ a \in \mathbb{Z} $ and $ b,m \in \mathbb{N}_0 $. Then 
$$
\sum_{i = 0}^m q^{ai - b(m - i)} 
\begin{bmatrix}
a \\ 
m - i
\end{bmatrix}
\begin{bmatrix}
b \\ 
i
\end{bmatrix}
= \begin{bmatrix}
a + b \\ 
m
\end{bmatrix}. 
$$
\end{lemma}

\begin{proof} We proceed using induction on $ b $. For $ b = 0 $ the sum reduces to the single term 
$$
q^{a0 - 0(m - 0)} 
\begin{bmatrix}
a \\ 
m - 0
\end{bmatrix}
\begin{bmatrix}
0 \\ 
0
\end{bmatrix}
= \begin{bmatrix}
a + 0 \\ 
m
\end{bmatrix}. 
$$
Suppose the assertion holds for some $ b $. Using Lemma \ref{lemqcalc1} with $ q $ replaced by $ q^{-1} $ 
and the inductive hypothesis we compute 
\begin{align*}
\sum_{i = 0}^m q^{ai - (b + 1)(m - i)} &
\begin{bmatrix}
a \\ 
m - i
\end{bmatrix}
\begin{bmatrix}
b + 1 \\ 
i
\end{bmatrix} \\ 
&= \sum_{i = 0}^m q^{ai - b(m - i) - m + i} 
\begin{bmatrix}
a \\ 
m - i
\end{bmatrix}
\biggl(q^{-i}
\begin{bmatrix}
b \\ 
i
\end{bmatrix} + q^{b - i + 1} 
\begin{bmatrix}
b \\ 
i - 1
\end{bmatrix}\biggr) \\
&= \sum_{i = 0}^m q^{ai - b(m - i) - m} 
\begin{bmatrix}
a \\ 
m - i
\end{bmatrix}
\begin{bmatrix}
b \\ 
i
\end{bmatrix} \\
&\qquad + \sum_{i = 1}^m q^{ai - b(m - i - 1) - m + 1} 
\begin{bmatrix}
a \\ 
m - i
\end{bmatrix}
\begin{bmatrix}
b \\ 
i - 1
\end{bmatrix} \\
&= q^{-m} 
\begin{bmatrix}
a + b \\ 
m 
\end{bmatrix}
+ 
\sum_{i = 0}^{m - 1} q^{ai + a - b(m - 1 - i) + b - m + 1} 
\begin{bmatrix}
a \\ 
m - i - 1
\end{bmatrix}
\begin{bmatrix}
b \\ 
i 
\end{bmatrix} \\
&= q^{-m} 
\begin{bmatrix}
a + b\\ 
m 
\end{bmatrix}
+ q^{a + b - m + 1} 
\begin{bmatrix}
a + b \\ 
m - 1
\end{bmatrix} \\
&= \begin{bmatrix}
a + b + 1 \\ 
m
\end{bmatrix}
\end{align*}
as desired. \end{proof} 

\begin{lemma} \label{tlemma2}
For any $ k,l \in \mathbb{N}_0 $ we have 
$$
\sum_{a = 0}^k \sum_{c = 0}^k (-1)^{a + c + k} q^{\pm(a + c - ac - k - kl)} 
\begin{bmatrix}
k \\ 
a
\end{bmatrix}
\begin{bmatrix}
a + l\\ 
k - c
\end{bmatrix} 
\begin{bmatrix}
c + l \\ 
c
\end{bmatrix}
= 1.  
$$
\end{lemma}
\begin{proof} Using $ 
\begin{bmatrix}
c + l \\ 
c
\end{bmatrix} = 
(-1)^c
\begin{bmatrix}
-l - 1 \\ 
c
\end{bmatrix} 
$ 
and $ - (a + l) c -(l + 1)(k - c) = c - ac -k - kl $ we obtain 
\begin{align*}
\sum_{a = 0}^k \sum_{c = 0}^k &(-1)^{a + c + k} q^{\pm(a + c - ac - k - kl)} 
\begin{bmatrix}
k \\ 
a
\end{bmatrix} 
\begin{bmatrix}
a + l\\ 
k - c
\end{bmatrix}
\begin{bmatrix}
c + l \\ 
c
\end{bmatrix} \\
&= \sum_{a = 0}^k \sum_{c = 0}^k (-1)^{a + k} q^{\pm(a + c - ac - k - kl)} 
\begin{bmatrix}
k \\ 
a
\end{bmatrix} 
\begin{bmatrix}
a + l \\ 
k - c
\end{bmatrix}
\begin{bmatrix}
- l - 1 \\ 
c
\end{bmatrix} \\
&= 
\sum_{a = 0}^k (-1)^{a + k} q^{\pm a} \begin{bmatrix}
k \\ 
a
\end{bmatrix} 
\sum_{c = 0}^k q^{\pm(- (a + l) c -(l + 1)(k - c))} 
\begin{bmatrix}
a + l \\ 
k - c
\end{bmatrix} 
\begin{bmatrix}
- l - 1 \\ 
c
\end{bmatrix}
\\
&= 
\sum_{a = 0}^k (-1)^{a + k} q^{\pm a} \begin{bmatrix}
k \\ 
a
\end{bmatrix} 
\begin{bmatrix}
a - 1 \\ 
k
\end{bmatrix}, 
\end{align*}
using Lemma \ref{tlemma1} for the sum over $ c $ in the final step. 
The last factor in this expression vanishes for $ 0 \leq a - 1 < k $, so only the summand for $ a = 0 $ survives. 
We conclude 
\begin{align*}
\sum_{a = 0}^k \sum_{c = 0}^k (-1)^{a + c + k} q^{\pm(a + c - ac - k - kl)} 
\begin{bmatrix}
k \\ 
a
\end{bmatrix} 
\begin{bmatrix}
a + l \\ 
k - c
\end{bmatrix}
\begin{bmatrix}
c + l \\ 
c
\end{bmatrix}
&= (-1)^k
\begin{bmatrix}
k \\ 
0
\end{bmatrix} 
\begin{bmatrix}
- 1 \\ 
k
\end{bmatrix} \\
&= (-1)^k \cdot 1 \cdot (-1)^k = 1 
\end{align*}
as claimed. \end{proof} 

We are now ready to obtain explicit formulas for the action of $ \T_\pm $ on the simple module $ V(n) $. 

\begin{prop} \label{texplicit}
For $ V = V(n) $ and any $ j \in \{n, n - 1, \dots, -n\} $ we have 
\begin{align*}
\T_\pm(v_{(j)}) &= (-1)^{n - j} q^{\pm (n - j)(n + j + 1)} v_{(-j)}.
\end{align*}
In particular, the operators $ \T_\pm $ are invertible with inverses given by 
\begin{align*}
\T_\pm^{-1}(v_{(j)}) &= (-1)^{n + j} q^{\mp (n + j)(n - j + 1)} v_{(-j)} = (-q)^{\mp 2j} \T_\mp(v_{(j)}).  
\end{align*}
\end{prop} 

\begin{proof} The formulas for the action of divided powers before Lemma \ref{tlemma1} imply 
\begin{align*}
&\T_\pm(v_{(j)}) = \sum_{\substack{r,s,t \geq 0 \\ r - s + t = 2j}} (-1)^s q^{\pm(s - rt)} F^{(r)} E^{(s)} F^{(t)} \cdot v_{(j)} \\
&= \sum_{\substack{r,s,t \geq 0 \\ r - s + t = 2j}} (-1)^s q^{\pm(s - rt)} 
\begin{bmatrix}
n - j + t \\ 
t
\end{bmatrix}
F^{(r)} E^{(s)} \cdot v_{(j - t)} \\ 
&= \sum_{\substack{r,s,t \geq 0 \\ r - s + t = 2j}} (-1)^s q^{\pm(s - rt)} 
\begin{bmatrix}
n + j - t + s \\ 
s
\end{bmatrix}
\begin{bmatrix}
n - j + t \\ 
t
\end{bmatrix}
F^{(r)} \cdot v_{(j + s - t)} \\
&= \sum_{\substack{r,s,t \geq 0 \\ r - s + t = 2j}} (-1)^s q^{\pm(s - rt)} 
\begin{bmatrix}
n - j + t - s + r \\ 
r
\end{bmatrix}
\begin{bmatrix}
n + j - t + s \\ 
s
\end{bmatrix}
\begin{bmatrix}
n - j + t \\ 
t
\end{bmatrix}
v_{(-j)}. 
\end{align*}
Note that all terms in this sum with $t>n+j$ or $r>n+j$ must vanish, since the vectors $v_{(j-t)}$ or $v_{(j+s-t)}=v_{(r-j)}$ vanish in these cases.  Also, when $t\leq n+j$, the middle binomial coefficient can be rewritten as
\[
 \begin{bmatrix} n + j - t + s \\ s \end{bmatrix}
  =  \begin{bmatrix} n + j - t + s \\ n+j-t \end{bmatrix},
\]
and the latter expression vanishes when $s < 0$.
Using $ s = r + t - 2j $ we therefore have to show 
\begin{align*}
 \sum_{r,t = 0}^{n + j} (-1)^{r + t - 2j} &q^{\pm(r + t - 2j - rt)} 
\begin{bmatrix}
n + j \\ 
r
\end{bmatrix}
\begin{bmatrix}
n - j + r \\ 
n + j - t
\end{bmatrix}
\begin{bmatrix}
n - j + t \\ 
t
\end{bmatrix} \\
&= (-1)^{n - j} q^{\pm (n - j)(n + j + 1)}.
\end{align*}
Setting $ k = n + j, l = n - j, a = r, c = t $ this amounts to 
\begin{align*}
\sum_{a,c = 0}^k (-1)^{l - k + a + c} q^{\pm(l - k + a + c - ac)} 
\begin{bmatrix}
k \\ 
a
\end{bmatrix}&
\begin{bmatrix}
l + a \\ 
k - c
\end{bmatrix}
\begin{bmatrix}
l + c \\ 
c
\end{bmatrix} 
= (-1)^l q^{\pm l(k + 1)}, 
\end{align*}
hence the first claim follows from Lemma \ref{tlemma2}. 

In order to verify the formula for $ \T_\pm^{-1} $ we compute 
\begin{align*}
(-1)^{n + j} &q^{\mp (n + j)(n - j + 1)} \T_\pm(v_{(-j)}) \\
&= (-1)^{n + j} q^{\mp (n + j)(n - j + 1)} (-1)^{n + j} q^{\pm (n + j)(n - j + 1)} v_{(j)} = v_{(j)}, 
\end{align*} 
and to check $ \T_\pm^{-1}(v_{(j)}) = (-q)^{\mp 2j} \T_\mp(v_{(j)}) $ we observe 
\begin{align*}
q^{\mp (n + j)(n - j + 1)} &= q^{\mp ((n + j)(n - j) + n + j)} \\
&= q^{\mp 2j} q^{\mp ((n + j)(n - j) + n - j)} = q^{\mp 2j}  q^{\mp (n - j)(n + j + 1)}.  
\end{align*}
This finishes the proof. \end{proof} 

The inverses of $ \T_\pm $ can be described in a similar fashion to the way we defined the operators $ \T_\pm $ themselves, as we discuss next.  

\begin{cor} \label{corformulatinverse}
Let $ V $ be an integrable module over $ U_q(\mathfrak{sl}(2, \mathbb{K})) $. Then we have 
\begin{align*}
\T_+^{-1}(v) &= \sum_{\substack{r,s,t \geq 0 \\ -r + s - t = 2m}} (-1)^s q^{rt - s} E^{(r)} F^{(s)} E^{(t)} \cdot v, \\
\T_-^{-1}(v) &= \sum_{\substack{r,s,t \geq 0 \\ -r + s - t = 2m}} (-1)^s q^{s - rt} E^{(r)} F^{(s)} E^{(t)} \cdot v
\end{align*} 
for $ v \in V_m $. 
\end{cor}

\begin{proof} It suffices to consider $ V = V(n) $. Consider the involutive linear automorphism $ \omega_V: V \rightarrow V $ given by 
$ \omega_V(v_{(j)}) = v_{(-j)} $. Using the automorphism $ \omega $ from Lemma \ref{defOmega} and the formulas before Lemma \ref{tlemma1} we obtain 
$$ 
\omega(X) \cdot \omega_V(v) = \omega_V(X \cdot v) 
$$ 
for all $ X \in U_q(\mathfrak{sl}(2, \mathbb{K})) $ and $ v \in V $. 
Accordingly, the claim becomes $ \T_\pm^{-1} \omega_V = \omega_V \T_\mp $. 
This in turn follows immediately from Proposition \ref{texplicit}. \end{proof} 

Let us also record the following commutation relations. 

\begin{lemma} \label{tcommutationlemma}
Let $ V $ be an integrable $ U_q(\mathfrak{sl}(2, \mathbb{K})) $-module and $ v \in V $. Then we have 
\begin{align*}
\T_\pm(E \cdot v) &= - K^{\pm 2} F \cdot \T_\pm(v) \\ 
\T_\pm(K \cdot v) &= K^{-1} \cdot \T_\pm(v) \\ 
\T_\pm(F \cdot v) &= - EK^{\mp 2} \cdot \T_\pm(v)
\end{align*}
and 
\begin{align*}
E \cdot \T_\pm(v) &= - \T_\pm(FK^{\mp 2}  \cdot v) \\ 
K \cdot \T_\pm(v) &= \T_\pm(K^{-1} \cdot v) \\ 
F \cdot \T_\pm(v) &= - \T_\pm(K^{\pm 2}E \cdot v). 
\end{align*}
\end{lemma} 

\begin{proof} It suffices to consider $ V = V(n) $ for $ n \in \frac{1}{2} \mathbb{N}_0 $. 
Since $ \T_\pm $ interchanges the weight spaces with weights $ \pm j $ for all $ j $ the formulas for the action of $ K $ are obvious. 
Using Lemma \ref{vnexplicit} and Lemma \ref{texplicit} we compute 
\begin{align*}
\T_\pm(E \cdot v_{(j)}) &= [n + j + 1] \T_\pm(v_{(j + 1)}) \\
&= (-1)^{n - j - 1} [n + j + 1] q^{\pm (n - j - 1)(n + j + 2)} v_{(-j - 1)} \\
&= (-1)^{n - j - 1} [n + j + 1] q^{\pm (n - j)(n + j) \pm(-n - j - 2 + 2n -2j)} v_{(-j - 1)} \\
&= -q^{\pm (-2j - 2)} (-1)^{n - j} q^{\pm(n - j)(n + j) \pm(n - j)} [n + j + 1] v_{(-j - 1)} \\ 
&= -q^{\pm (-2j - 2)} (-1)^{n - j} q^{\pm(n - j)(n + j + 1)} F \cdot v_{(-j)} \\ 
&= - q^{\pm (-2j - 2)} F \cdot \T_\pm(v_{(j)}) \\ 
&= - K^{\pm 2} F \cdot \T_\pm(v_{(j)}), 
\end{align*}
and similarly
\begin{align*}
\T_\pm(F \cdot v_{(j)}) &= [n - j + 1] \T_\pm(v_{(j - 1)}) \\
&= (-1)^{n - j + 1} [n - j + 1] q^{\pm(n - j + 1)(n + j)} v_{(-j + 1)} \\
&= (-1)^{n - j + 1} [n - j + 1] q^{\pm(n - j)(n + j) \pm(n + j)} v_{(-j + 1)} \\
&= -q^{\pm 2j} (-1)^{n - j} q^{\pm(n - j)(n + j) \pm(n - j)} [n - j + 1] v_{(-j + 1)} \\ 
&= -q^{\pm 2j} (-1)^{n - j} q^{\pm(n - j)(n + j + 1)} E \cdot v_{(-j)} \\ 
&= -q^{\pm 2j} E \cdot \T_\pm(v_{(j)}) \\ 
&= -E K^{\mp 2} \cdot \T_\pm(v_{(j)}). 
\end{align*}
The remaining formulas can be easily deduced from these relations. \end{proof}

\subsubsection{The case of general $ \mathfrak{g} $}  

Now let $ \mathfrak{g} $ be an arbitrary semisimple Lie algebra. 
Based on the constructions in the previous subsection we will define operators $ \T_i $ for $ i = 1, \dots, N $ acting on every 
integrable $ U_q(\mathfrak{g}) $-module. 

Specifically, let $ V $ be an integrable $ U_q(\mathfrak{g}) $-module and recall that $ U_{q_i}(\mathfrak{g}_i) $ for $ 1 \leq i \leq N $ denotes 
the subalgebra of $ U_q(\mathfrak{g}) $ generated by $ E_i, F_i, K_i^{\pm 1} $. Using the canonical embedding 
$ U_{q_i}(\mathfrak{g}_i) \subset U_q(\mathfrak{g}) $ we can view $ V $ as an 
integrable $ U_{q_i}(\mathfrak{g}_i) $-module, and we define $ \T_i: V \rightarrow V $ to be 
the operator corresponding to $ \T = \T_+ $ in the notation of the previous subsection. 
Explicitly, 
\begin{align*}
\T_i(v) &= \sum_{\substack{r,s,t \geq 0 \\ r - s + t = m}} (-1)^s q_i^{s - rt} F_i^{(r)} E_i^{(s)} F_i^{(t)} \cdot v 
\end{align*} 
\label{nom:T_i}%
for $ v \in V_\lambda $ and $ m = (\alpha_i^\vee, \lambda) $. According to Proposition \ref{texplicit}, the operators $ \T_i $ are bijective, 
and $ \T_i $ maps $ V_\lambda $ onto $ V_{s_i \lambda} $. 
From this observation it follows that 
$$
\T_i(K_\mu \cdot v) = K_{s_i \mu} \cdot \T_i(v) 
$$
for all $ \mu \in \weights $. Let us also record the following formulas. 

\begin{lemma} \label{Ticommutationeasy}
Let $ V $ be an integrable $ U_q(\mathfrak{g}) $-module and $ v \in V $. Then we have 
\begin{align*}
\T_i(E_i \cdot v) &= - K_i F_i \cdot \T_i(v) \\ 
\T_i(F_i \cdot v) &= - E_i K_i^{-1} \cdot \T_i(v) \\ 
E_i \cdot \T_i(v) &= - \T_i(F_i K_i^{-1}\cdot v) \\ 
F_i \cdot \T_i(v) &= - \T_i(K_i E_i \cdot v)
\end{align*}
for all $ i = 1, \dots, N $. 
\end{lemma} 

\begin{proof} This is an immediate consequence of Lemma \ref{tcommutationlemma}. \end{proof} 

Our first aim is to obtain explicit formulas for $ \T_i(E_j \cdot v) $ and $ E_j \cdot \T_i(v) $ also in the case $ i \neq j $, and similarly 
for $ F_j $ instead of $ E_j $. 

Rewriting the results from Lemma \ref{adjointpowers} in terms of divided powers we obtain 
\begin{align*}
E^{(m)}_i \rightarrow Y &= 
\sum_{k = 0}^m (-1)^{m - k} q_i^{-(m - k)(m - 1)} E_i^{(k)} Y E_i^{(m - k)} K_i^{-m} \\
F_i^{(m)} \rightarrow Y &= \sum_{k = 0}^m (-1)^k q_i^{k (m - 1)} 
F_i^{(m - k)} K_i^{-k} Y K_i^k F_i^{(k)}. 
\end{align*}
In particular, we have 
\begin{align*}
F_i^{(m)} \rightarrow F_j &= \sum_{k = 0}^m (-1)^k q_i^{k (m - 1 + a_{ij})} F_i^{(m - k)} F_j F_i^{(k)}, 
\end{align*}
using $ K_i^{-k} F_j K_i^k = q_i^{k a_{ij}} F_j $. 
\begin{lemma} \label{Ticommutationlemma}
If $ i \neq j $ then for all $ m, l \in \mathbb{N}_0 $ we have 
\begin{align*}
(F_i^{(m)} \rightarrow F_j) E_i^{(l)} &= \sum_{k = 0}^l (-1)^k 
\begin{bmatrix}
-a_{ij} - m + k \\ 
k
\end{bmatrix}_{q_i}
q_i^{k (l - 1)} E_i^{(l - k)} (F_i^{(m - k)} \rightarrow F_j) K_i^k \\
(F_i^{(m)} \rightarrow F_j) F_i^{(l)} &= \sum_{k = 0}^l (-1)^k 
\begin{bmatrix}
m + k \\ 
k
\end{bmatrix}_{q_i}
q_i^{-l (a_{ij} + 2m) - k(l - 1)} F_i^{(l - k)} (F_i^{(m + k)} \rightarrow F_j). 
\end{align*}
\end{lemma} 
\begin{proof} We use induction on $ l $ for both formulas. For $ l = 0 $ there is nothing to show. 

Assume the first formula holds for some $ l $. Note that since $ E_i $ and $ F_j $ commute we have $ E_i \rightarrow F_j = 0 $, and hence 
\begin{align*}
(E_i F_i^{(m)}) \rightarrow F_j &= [-a_{ij} + 1 - m]_{q_i} F_i^{(m - 1)} \rightarrow F_j
\end{align*}
due to Lemma \ref{EFbasic}.  
Writing $ a(m) = F^{(m)}_i \rightarrow F_j $ we deduce 
\begin{align*}
a(m) E_i K_i^{-1} &= E_i a(m) K_i^{-1} - E_i \rightarrow a(m) \\
&= E_i a(m) K_i^{-1} - [-a_{ij} + 1 - m]_{q_i} a(m - 1), 
\end{align*}
which implies 
\begin{align*}
a(m - k) E_i &= E_i a(m - k) - [-a_{ij} + 1 + k - m]_{q_i} a(m - k - 1) K_i 
\end{align*}
for all $ 0 \leq k \leq m $. 
Now from the induction hypothesis we obtain 
\begin{align*}
a(m) E_i^{(l + 1)} &= a(m) \frac{E_i^{(l)} E_i}{[l + 1]_{q_i}} \\
&= \sum_{k = 0}^l \frac{(-1)^k}{[l + 1]_{q_i}} 
\begin{bmatrix}
-a_{ij} - m + k \\ 
k
\end{bmatrix}_{q_i}
q_i^{k (l - 1)} E_i^{(l - k)} (F_i^{(m - k)} \rightarrow F_j) K_i^k E_i. 
\end{align*}
Since $ K_i^k E_i = q_i^{2k} E_i K_i^k $ our above computation yields 
\begin{align*}
a(m) E_i^{(l + 1)} &= \sum_{k = 0}^{l + 1} b_k E_i^{(l + 1 - k)} a(m - k) K_i^k, 
\end{align*}
where 
\begin{align*}
b_k &= \frac{(-1)^k}{[l + 1]_{q_i}}
\begin{bmatrix}
-a_{ij} - m + k \\ 
k
\end{bmatrix}_{q_i}
q_i^{k(l - 1) + 2k} [l + 1 - k]_{q_i} \\
&\qquad - \frac{(-1)^{k - 1}}{[l + 1]_{q_i}}
\begin{bmatrix}
-a_{ij} - m + k - 1 \\ 
k - 1
\end{bmatrix}_{q_i}
q_i^{(k - 1)(l - 1) + 2k - 2} [-a_{ij} - m + k]_{q_i} \\
&= \frac{(-1)^k}{[l + 1]_{q_i}}
\begin{bmatrix}
-a_{ij} - m + k \\ 
k
\end{bmatrix}_{q_i}
q_i^{kl + k} [l + 1 - k]_{q_i} \\
&\qquad + \frac{(-1)^k}{[l + 1]_{q_i}}
\begin{bmatrix}
-a_{ij} - m + k \\ 
k 
\end{bmatrix}_{q_i}
q_i^{kl + k - l - 1} [k]_{q_i}, 
\end{align*}
note here that the first summand vanishes for $ k = l + 1 $, 
and the second summand vanishes for $ k = 0 $. 
Using 
\begin{align*}
\frac{q_i^{k}[l + 1 - k]_{q_i} + q_i^{k - l - 1}[k]_{q_i}}{[l + 1]_{q_i}} 
&= \frac{q_i^{k}(q_i^{l + 1 - k} - q^{-l - 1 + k}) + q_i^{k - l - 1}(q_i^k - q_i^{-k})}{q_i^{l + 1} - q_i^{-l - 1}} = 1
\end{align*}
we conclude 
$$
b_k = (-1)^k
\begin{bmatrix}
-a_{ij} - m + k \\ 
k 
\end{bmatrix}_{q_i}
q_i^{kl}
$$
as desired. 

For the second formula we proceed in the same way. Assume that the formula holds for some $ l $. 
Writing again $ a(m) = F^{(m)}_i \rightarrow F_j $ we compute 
\begin{align*}
a(m) F_i^{(l + 1)} &= a(m) \frac{F_i^{(l)} F_i}{[l + 1]_{q_i}} \\
&= \sum_{k = 0}^l \frac{(-1)^k}{[l + 1]_{q_i}}
\begin{bmatrix}
m + k \\ 
k
\end{bmatrix}_{q_i}
q_i^{-l (a_{ij} + 2m) - k(l - 1)} F_i^{(l - k)} a(m + k) F_i
\end{align*}
according to the induction hypothesis. Note also that 
\begin{align*}
a(m) F_i &= q_i^{-2m - a_{ij}} K_i^{-1} a(m) K_i F_i \\
&= q_i^{-2m - a_{ij}} (F_i a(m) - F_i \rightarrow a(m)) \\
&= q_i^{-2m - a_{ij}} (F_i a(m) - [m + 1]_{q_i} a(m + 1)), 
\end{align*}
using $ [m + 1]_{q_i} F_i^{(m + 1)} = F_i F_i^{(m)} $ in the last step. We thus have 
\begin{align*}
a(m + k) F_i &= q_i^{-2m - 2k - a_{ij}} (F_i a(m + k) - [m + k + 1]_{q_i} a(m + k + 1)),  
\end{align*}
and combining these formulas we obtain 
\begin{align*}
a(m) F_i^{(l + 1)} &= \sum_{k = 0}^{l + 1} c_k F_i^{(l + 1 - k)} a(m + k) 
\end{align*}
where 
\begin{align*}
c_k &= \frac{(-1)^k}{[l + 1]_{q_i}}
\begin{bmatrix}
m + k \\ 
k
\end{bmatrix}_{q_i}
q_i^{-(l + 1) (a_{ij} + 2m) - 2k - k(l - 1)} [l + 1 - k]_{q_i} \\
&\qquad - \frac{(-1)^{k - 1}}{[l + 1]_{q_i}}
\begin{bmatrix}
m + k - 1\\ 
k - 1
\end{bmatrix}_{q_i}
q_i^{-(l + 1) (a_{ij} + 2m) - 2k + 2 - (k - 1)(l - 1)} [m + k]_{q_i} \\
&= \frac{(-1)^k}{[l + 1]_{q_i}}
\begin{bmatrix}
m + k \\ 
k
\end{bmatrix}_{q_i}
q_i^{-(l + 1) (a_{ij} + 2m) - kl} q^{-k} [l + 1 - k]_{q_i} \\
&\qquad + \frac{(-1)^k}{[l + 1]_{q_i}}
\begin{bmatrix}
m + k \\ 
k
\end{bmatrix}_{q_i}
q_i^{-(l + 1) (a_{ij} + 2m) - kl - k + l + 1} [k]_{q_i}, 
\end{align*}
note in particular that the first summand vanishes for $ k = l + 1 $, and the second summand vanishes for $ k = 0 $. 
Using 
\begin{align*}
\frac{q_i^{-k}[l + 1 - k]_{q_i} + q_i^{l + 1 - k}[k]_{q_i}}{[l + 1]_{q_i}} 
&= \frac{q_i^{-k}(q_i^{l + 1 - k} - q^{-l - 1 + k}) + q_i^{l + 1 - k}(q_i^k - q_i^{-k})}{q_i^{l + 1} - q_i^{-l - 1}} = 1
\end{align*}
we conclude 
$$
c_k = (-1)^k
\begin{bmatrix}
m + k \\ 
k
\end{bmatrix}_{q_i}
q_i^{-(l + 1) (a_{ij} + 2m) - kl}. 
$$
This finishes the proof. \end{proof} 

\begin{prop} \label{Tigeneratorscommutation}
Let $ V $ be an integrable $ U_q(\mathfrak{g}) $-module. If $ i \neq j $ we have 
\begin{align*}
\T_i(E_j \cdot v) &= (-q_i)^{-a_{ij}} (E_j \leftarrow \hat{S}^{-1}(E^{(-a_{ij})}_i)) \cdot \T_i(v) \\ 
\T_i(F_j \cdot v) &= (F^{(-a_{ij})}_i \rightarrow F_j) \cdot \T_i(v) 
\end{align*}
for all $ v \in V $. 
\end{prop}

\begin{proof} Let us first consider the second formula. In the same way as in the proof of Lemma \ref{Ticommutationlemma} we shall 
abbreviate $ a(m) = F_i^{(m)} \rightarrow F_j $. Since $ F_i^{-a_{ij} + 1} \rightarrow F_j = u^-_{ij} = 0 $, the second part of 
Lemma \ref{Ticommutationlemma} implies
$$
a(-a_{ij}) F_i^{(l)} = (F_i^{(-a_{ij})} \rightarrow F_j) F_i^{(l)} = 
q_i^{l a_{ij}} F_i^{(l)} (F_i^{(-a_{ij})} \rightarrow F_j) = q_i^{l a_{ij}} F_i^{(l)} a(-a_{ij})
$$
for any $ l $. 
Moreover we obtain 
\begin{align*}
a(-a_{ij}) E_i^{(l)} &= (F_i^{(-a_{ij})} \rightarrow F_j) E_i^{(l)} \\
&= \sum_{k = 0}^l (-1)^k q_i^{k (l - 1)} E_i^{(l - k)} (F_i^{(-a_{ij} - k)} \rightarrow F_j) K_i^k \\
&= \sum_{k = 0}^l (-1)^k q_i^{k (l - 1)} E_i^{(l - k)} a(-a_{ij} - k) K_i^k 
\end{align*}
from the first part of Lemma \ref{Ticommutationlemma}. 

Assume $ v \in V_\lambda $ for some $ \lambda \in \weights $. Writing $ p = (\alpha_i^\vee, \lambda) $, the above formulas imply
\begin{align*}
a(-a_{ij}) &\T_i(v) = \sum_{\substack{r,s,t \geq 0 \\ r - s + t = p}} (-1)^s q_i^{s - rt} a(-a_{ij}) F_i^{(r)} E_i^{(s)} F_i^{(t)} \cdot v \\ 
&= \sum_{\substack{r,s,t \geq 0 \\ r - s + t = p}} (-1)^s q_i^{s - rt + r a_{ij}} F_i^{(r)} a(-a_{ij}) E_i^{(s)} F_i^{(t)} \cdot v \\ 
&= \sum_{\substack{r,s,t \geq 0 \\ r - s + t = p}} \sum_{k = 0}^s (-1)^{k + s} 
q_i^{s - rt + r a_{ij} + k (s - 1)} F_i^{(r)} E_i^{(s - k)} a(-a_{ij} - k) K_i^k F_i^{(t)} \cdot v \\ 
&= \sum_{\substack{r,s,t \geq 0 \\ r - s + t = p}} \sum_{k = 0}^s (-1)^{k + s} 
q_i^{s - r(t - a_{ij}) + k (s - 1 + p - 2t)} F_i^{(r)} E_i^{(s - k)} a(-a_{ij} - k) F_i^{(t)} \cdot v \\
&= \sum_{\substack{r,s,t \geq 0 \\ r - s + t = p}} \sum_{k = 0}^s \sum_{l = 0}^t (-1)^{k + s} 
q_i^{s - r(t - a_{ij}) + k (s - 1 + p - 2t)} F_i^{(r)} E_i^{(s - k)} \times \\ 
&\qquad (-1)^l \begin{bmatrix}
-a_{ij} - k + l \\ 
l
\end{bmatrix}_{q_i}
q_i^{-t (-a_{ij} -2k) - l(t - 1)} F_i^{(t - l)} a(-a_{ij} - k + l) \cdot v \\
&= \sum_{\substack{r,s,t \geq 0 \\ r - s + t = p}} \sum_{k = 0}^s \sum_{l = 0}^t (-1)^{k + s + l} 
\begin{bmatrix}
-a_{ij} - k + l \\ 
l
\end{bmatrix}_{q_i} \times \\
&\qquad q_i^{s - r(t - a_{ij}) + k (s - 1 + p) + t a_{ij} - l(t - 1)} F_i^{(r)} E_i^{(s - k)} F_i^{(t - l)} a(-a_{ij} - k + l) \cdot v,  
\end{align*} 
using 
$$
K_i^k F_i^{(t)} \cdot v = q_i^{-2 k t} F_i^{(t)} K_i^k \cdot v = q_i^{kp - 2kt} F_i^{(t)} \cdot v 
$$
for $ k, t \in \mathbb{N}_0 $ and the second part of Lemma \ref{Ticommutationlemma}. 
Setting $ b = s - k, c = t - l $ we have $ -r + b - c = -p - k + l $, and since $ a(-a_{ij} - k + l) = 0 $ for $ l > k $ we see that 
$ a(-a_{ij}) \T_i(v) $ is a linear combination of terms $ F_i^{(r)} E_i^{(b)} F_i^{(c)} a(-a_{ij} - h) \cdot v $ with 
$ -r + b - c = -p - h $ and $ h \geq 0 $. Explicitly, the coefficient of $ F_i^{(r)} E_i^{(b)} F_i^{(c)} a(-a_{ij} - h) \cdot v $ is 
\begin{align*} 
\sum_{l = 0}^{-a_{ij} - h} (-1)^{b + l} 
\begin{bmatrix}
-a_{ij} - h \\ 
l
\end{bmatrix}_{q_i} 
q_i^{b + k - r(c + l - a_{ij}) + k (b + k - 1 + p) + (c + l)a_{ij} - l(c + l - 1)} ,  
\end{align*}  
where $ k = l + h $. If we insert $ k = l + h $ in this expression we obtain 
\begin{align*} 
(-1)^b &q_i^{b + h - r(c - a_{ij}) + h(b + h - 1 + p) + c a_{ij}}  
\sum_{l = 0}^{-a_{ij} - h} (-1)^l 
\begin{bmatrix}
-a_{ij} - h \\ 
l
\end{bmatrix}_{q_i} \times \\
&\qquad q_i^{l(1 - r + b + l + h - 1 + p + h + a_{ij} - c - l + 1)} \\
&= (-1)^b q_i^{b - r(c - a_{ij}) + h(b + h + p) + c a_{ij}}  
\sum_{l = 0}^{-a_{ij} - h} (-1)^l 
\begin{bmatrix}
-a_{ij} - h \\ 
l
\end{bmatrix}_{q_i} \times \\
&\qquad q_i^{l(1 - r + b + 2h + p + a_{ij} - c)}. 
\end{align*}  
Using $ -r + b - c = -p - h $, the sum over $ l $ in the previous formula becomes 
\begin{align*} 
\sum_{l = 0}^{-a_{ij} - h} (-1)^l 
\begin{bmatrix}
-a_{ij} - h \\ 
l
\end{bmatrix}_{q_i} 
q_i^{l(1 + h + a_{ij})} = \delta_{-a_{ij} - h, 0} 
\end{align*}  
according to Lemma \ref{lembinomialvanishing}. Hence only the term with $ h = -a_{ij} $ survives, and since 
$$ 
q_i^{b - r(c - a_{ij}) + h(b + h + p) + c a_{ij}} = q_i^{b - rc + h(-r + b - c + h + p)} = q_i^{b - rc} 
$$ 
for $ h = -a_{ij} $, we conclude 
\begin{align*} 
a(-a_{ij}) \T_i(v) &= \sum_{\substack{r,b,c \geq 0 \\ r - b + c = p - a_{ij}}} (-1)^b
q_i^{b - rc} F_i^{(r)} E_i^{(b)} F_i^{(c)} a(0) v = \T_i(F_j \cdot v) 
\end{align*}
as desired.

Let us now consider the first relation. 
Fix $ v \in V_\lambda $ and let $ \omega_V $ denote the automorphism of $ V $ defined as in the proof of Corollary \ref{corformulatinverse} by viewing 
$ V $ as a $ U_{q_i}(\mathfrak{q}_i) $-module. Then combining 
Proposition \ref{texplicit} and the proof Corollary \ref{corformulatinverse} yields 
$$ 
\omega_V \T_i \omega_V(v) = (-q_i)^{(\alpha_i^\vee, \lambda)} \T_i(v).
$$ 
Moreover we recall that 
\begin{align*}
\omega(X \rightarrow Y) &= \omega(X_{(1)} Y \hat{S}(X_{(2)})) \\
&= \omega(X_{(1)}) \omega(Y) \hat{S}^{-1}(\omega(X_{(2)})) \\
&= \omega(X)_{(2)} \omega(Y) \hat{S}^{-1}(\omega(X)_{(1)}) = \omega(Y) \leftarrow \hat{S}^{-1}(\omega(X)) 
\end{align*}
for $ X, Y \in U_q(\mathfrak{g}) $, where $ Y \leftarrow X = \hat{S}(X_{(1)}) Y X_{(2)} $. 
Using that 
the vector $ F_j \cdot \omega_V(v) $ has weight $ -\lambda - \alpha_j $ we therefore obtain 
\begin{align*}
\T_i(E_j \cdot v) &= \T_i(\omega(F_j) \cdot v) \\ 
&= \T_i(\omega_V (F_j \cdot \omega_V(v))) \\ 
&= (-q_i)^{(\alpha_i^\vee, -\lambda - \alpha_j)} \omega_V (\T_i(F_j \cdot \omega_V(v))) \\ 
&= (-q_i)^{-(\alpha_i^\vee, \lambda + \alpha_j)} \omega_V((F_i^{(-a_{ij})} \rightarrow F_j) \cdot \T_i(\omega_V(v))) \\ 
&= (-q_i)^{-(\alpha_i^\vee, \lambda) - a_{ij}} \omega(F_i^{(-a_{ij})} \rightarrow F_j) \cdot \omega_V(\T_i(\omega_V(v))) \\ 
&= (-q_i)^{-(\alpha_i^\vee, \lambda) - a_{ij}} (-q_i)^{(\alpha_i^\vee, \lambda)} \omega(F_i^{(-a_{ij})} \rightarrow F_j) \cdot \T_i(v) \\ 
&= (-q_i)^{- a_{ij}} (\omega(F_j) \leftarrow \hat{S}^{-1}(\omega(F_i^{(-a_{ij})}))) \cdot \T_i(v) \\ 
&= (-q_i)^{- a_{ij}} (E_j \leftarrow \hat{S}^{-1}(E^{(-a_{ij})}_i)) \cdot \T_i(v). 
\end{align*}
This finishes the proof. \end{proof} 

We obtain an algebra antiautomorphism $ \gamma: U_q(\mathfrak{g}) \rightarrow U_q(\mathfrak{g}) $ by defining 
$$
\gamma(E_i) = E_i, \qquad \gamma(F_i) = F_i, \qquad \gamma(K_\mu) = K_{-\mu} 
$$
\nomenclature{$\gamma$}{algebra anti-automorphism of $U_q(\lie{g})$}%
on generators. Notice that $ \gamma $ is involutive, that is, $ \gamma^2 = \id $.

Let us define the $ \gamma $-twisted adjoint action by 
$$
X \overset{\gamma}{\rightarrow} Y = \gamma(X \rightarrow \gamma(Y)) 
$$
\label{nom:gamma-twisted}%
for $ X, Y \in U_q(\mathfrak{g}) $. Note that 
\begin{align*}
E_i \overset{\gamma}{\rightarrow} Y &= \gamma([E_i, \gamma(Y)] K_i^{-1}) = K_i [Y, E_i] \\
F_i \overset{\gamma}{\rightarrow} Y &= \gamma(F_i \gamma(Y) - K_i^{-1} \gamma(Y) K_i F_i) = Y F_i - F_i K_i^{-1} Y K_i \\ 
K_\mu \overset{\gamma}{\rightarrow} Y &= \gamma(K_\mu \gamma(Y) K_\mu^{-1}) = K_\mu Y K_\mu^{-1}
\end{align*} 
for all $ Y \in U_q(\mathfrak{g}) $, where $ 1 \leq i \leq N $ and $ \mu \in \weights $. 

\begin{lemma} \label{tauadcommutation1}
Suppose $ X, Y \in U_q(\mathfrak{g}) $ satisfy $ X \cdot \T_i(v) = \T_i(Y \cdot v) $ for any vector $ v $ in an integrable $ U_q(\mathfrak{g}) $-module. Then 
\begin{align*}
(E_i \rightarrow X) \cdot \T_i(v) &= \T_i((F_i \overset{\gamma}{\rightarrow} Y) \cdot v) \\
(F_i \rightarrow X) \cdot \T_i(v) &= \T_i((E_i \overset{\gamma}{\rightarrow} Y) \cdot v).  
\end{align*}
\end{lemma}

\begin{proof} Using the formula $ E_i \cdot \T_i(v) = - \T_i(F_iK_i^{-1} \cdot v) $ from Lemma \ref{Ticommutationeasy} we compute 
\begin{align*} 
(E_i \rightarrow X) \cdot \T_i(v) &= [E_i, X] K_i^{-1} \cdot \T_i(v) \\
&= \T_i((-F_i K_i^{-1} Y + Y F_i K_i^{-1}) K_i \cdot v) \\
&= \T_i((Y F_i - F_i K_i^{-1} Y K_i) \cdot v) \\
&= \T_i((F_i \overset{\gamma}{\rightarrow} Y) \cdot v) 
\end{align*} 
for $ v \in V $. Similarly, using $ F_i \cdot \T_i(v) = - \T_i(K_i E_i \cdot v) $ we obtain 
\begin{align*} 
(F_i \rightarrow X) \cdot \T_i(v) &= (F_i X - K_i^{-1} X K_i F_i) \cdot \T_i(v) \\
&= \T_i((-K_i E_i Y + K_i Y K_i^{-1} K_i E_i) \cdot v) \\
&= \T_i(K_i [Y, E_i] \cdot v) \\
&= \T_i((E_i \overset{\gamma}{\rightarrow} Y) \cdot v). 
\end{align*} 
This yields the claim. \end{proof} 

\begin{lemma} \label{Ticommutationlemmahigherterms}
Let $ V $ be an integrable $ U_q(\mathfrak{g}) $-module and $ v \in V $. If $ i \neq j $ we have 
\begin{align*}
(F_i^{(-a_{ij} - l)} \rightarrow F_j) \cdot \T_i(v) &= \T_i((F_i^{(l)} \overset{\gamma}{\rightarrow} F_j) \cdot v)
\end{align*}
for any $ 0 \leq l \leq -a_{ij} $. 
\end{lemma} 

\begin{proof} Using induction on $ l $, Lemma \ref{tauadcommutation1} and Proposition \ref{Tigeneratorscommutation} imply 
\begin{align*}
\T_i((F_i^{(l)} \overset{\gamma}{\rightarrow} F_j) \cdot v) &= (E_i^{(l)} F^{(-a_{ij})}_i \rightarrow F_j) \cdot \T_i(v). 
\end{align*} 
As already observed in the proof of Lemma \ref{Ticommutationlemma}, we have
\begin{align*}
(E_i F_i^{(m)}) \rightarrow F_j &= [-a_{ij} + 1 - m]_{q_i} F_i^{(m - 1)} \rightarrow F_j. 
\end{align*}
Applying this formula iteratively we obtain
\begin{align*}
(E_i^{(l)} F_i^{(-a_{ij})}) \rightarrow F_j &= F_i^{(-a_{ij} - l)} \rightarrow F_j, 
\end{align*}
which yields the desired formula. \end{proof} 

\begin{prop} \label{tiformulas}
Let $ V $ be an integrable $ U_q(\mathfrak{g}) $-module and $ v \in V $. If $ i \neq j $ we have 
\begin{align*}
\sum_{k = 0}^{-a_{ij} - l} (-1)^k q_i^{k (l + 1)} &E_i^{(k)} E_j E_i^{(-a_{ij} - l - k)} \cdot \T_i(v) \\
&= \sum_{k = 0}^l (-1)^k q_i^{-k (l - 1 + a_{ij})} \T_i(E_i^{(l - k)} E_j E_i^{(k)} \cdot v)
\end{align*} 
and 
\begin{align*}
\sum_{k = 0}^{-a_{ij} - l} (-1)^k q_i^{-k (l + 1)} &F_i^{(-a_{ij} - l - k)} F_j F_i^{(k)} \cdot \T_i(v) \\
&= \sum_{k = 0}^l (-1)^k q_i^{k (l - 1 + a_{ij})} \T_i(F_i^{(k)} F_j F_i^{(l - k)} \cdot v) 
\end{align*} 
for any $ 0 \leq l \leq -a_{ij} $. 
\end{prop} 

\begin{proof} Let us first consider the second formula. Recall from the calculations before Lemma \ref{Ticommutationlemma} that 
\begin{align*}
F_i^{(m)} \rightarrow F_j &= \sum_{k = 0}^m (-1)^k q_i^{k (m - 1 + a_{ij})} F_i^{(m - k)} F_j F_i^{(k)}. 
\end{align*}
Since $ \gamma $ is an anti-automorphism fixing the generators $ F_k $ for all $ k $ we therefore obtain 
\begin{align*}
F_i^{(m)} \overset{\gamma}{\rightarrow} F_j &= \gamma(F_i^{(m)} \rightarrow F_j) = \sum_{k = 0}^m (-1)^k q_i^{k (m - 1 + a_{ij})} F_i^{(k)} F_j F_i^{(m - k)}. 
\end{align*}
Combining these formulas with Lemma \ref{Ticommutationlemmahigherterms} yields 
\begin{align*}
\sum_{k = 0}^{-a_{ij} - l} (-1)^k q_i^{k (-l - 1)} &F_i^{(-a_{ij} - l - k)} F_j F_i^{(k)} \cdot \T_i(v) \\
&= \sum_{k = 0}^l (-1)^k q_i^{k (l - 1 + a_{ij})} \T_i(F_i^{(k)} F_j F_i^{(l - k)} \cdot v)
\end{align*} 
as desired. 

Next observe that $ \omega(X) \cdot \T_i(v) = (-q_i)^{(\alpha_i^\vee, \mu)} \T_i(\omega(Y) \cdot v) $ if $ X $ has weight $ \mu $ with respect 
to the adjoint action and $ X \cdot \T_i(v) = \T_i(Y \cdot v) $ for all vectors $ v $ in integrable $ U_q(\mathfrak{g}) $-modules. 
Hence the formula proved above implies 
\begin{align*}
\sum_{k = 0}^{-a_{ij} - l} &(-1)^k q_i^{- k (-l - 1)} E_i^{(k)} E_j E_i^{(-a_{ij} - l - k)} \cdot \T_i(v) \\
&= q_i^{(-a_{ij} - l)(l + 1)} \sum_{k = 0}^{-a_{ij} - l} (-1)^k q_i^{(-a_{ij} - l - k) (-l - 1)} E_i^{(k)} E_j E_i^{(-a_{ij} - l - k)} \cdot \T_i(v) \\
&= (-1)^{a_{ij} + l} q_i^{(-a_{ij} - l)(l + 1)} \sum_{k = 0}^{-a_{ij} - l} (-1)^k q_i^{k (-l - 1)} E_i^{(-a_{ij} - l - k)} E_j E_i^{(k)} \cdot \T_i(v) \\
&= (-1)^{a_{ij} + l} q_i^{(-a_{ij} - l)(l + 1)} (-q_i)^{2l + a_{ij}} \sum_{k = 0}^l (-1)^k q_i^{k (l - 1 + a_{ij})} \T_i(E_i^{(k)} E_j E_i^{(l - k)} \cdot v) \\
&= q_i^{(-a_{ij} - l)(l + 1)} q_i^{2l + a_{ij}} q_i^{l (l - 1 + a_{ij})} \sum_{k = 0}^l (-1)^k q_i^{-k (l - 1 + a_{ij})} \T_i(E_i^{(l - k)} E_j E_i^{(k)} \cdot v) \\
&= \sum_{k = 0}^l (-1)^k q_i^{-k (l - 1 + a_{ij})} \T_i(E_i^{(l - k)} E_j E_i^{(k)} \cdot v), 
\end{align*} 
upon reindexing $ k $ to $ -a_{ij} - l - k $, and further below reindexing $ k $ to $ l - k $, respectively. 
Note here that $ F_i^{(-a_{ij} - l)} \rightarrow F_j $ has weight $ (a_{ij} + l)\alpha_i - \alpha_j $, so 
that 
$$ 
(\alpha_i^\vee, (a_{ij} + l)\alpha_i - \alpha_j) = 2l + a_{ij}. 
$$
This finishes the proof. \end{proof} 

Note that putting $l=0$ in Proposition \ref{tiformulas} yields in particular 
\begin{align*}
\T_i(E_j \cdot v) &= \sum_{k = 0}^{-a_{ij}} (-1)^k q_i^k E_i^{(k)} E_j E_i^{(-a_{ij} - k)} \cdot \T_i(v) \\
\T_i(F_j \cdot v) &= \sum_{k = 0}^{-a_{ij}} (-1)^k q_i^{-k} F_i^{(-a_{ij} - k)} F_j F_i^{(k)} \cdot \T_i(v).  
\end{align*}

\subsubsection{The braid group action on $ U_q(\mathfrak{g}) $}

We are now ready to construct automorphisms $ \T_i $ of $ U_q(\mathfrak{g}) $ for $ i = 1, \dots, N $. 

\begin{theorem} \label{deftiautomorphism}
For $ i = 1, \dots, N $ there exist algebra automorphisms $ \T_i: U_q(\mathfrak{g}) \rightarrow U_q(\mathfrak{g}) $ 
satisfying 
\begin{align*}
\T_i(K_\mu) &= K_{s_i \mu},
\quad \T_i(E_i) = - K_i F_i, \quad \T_i(F_i) = - E_i K_i^{-1}\\ 
\T_i(E_j) &= \sum_{k = 0}^{-a_{ij}} (-1)^k q_i^k E_i^{(k)} E_j E_i^{(-a_{ij} - k)}, \quad i \neq j \\ 
\T_i(F_j) &= \sum_{k = 0}^{-a_{ij}} (-1)^k q_i^{-k} F_i^{(-a_{ij} - k)} F_j F_i^{(k)}, \quad i \neq j. 
\end{align*}
\end{theorem} 

\begin{proof} 
We shall define $ \T_i: U_q(\mathfrak{g}) \rightarrow U_q(\mathfrak{g}) $ using conjugation with the operators $ \T_i $ on integrable 
$ U_q(\mathfrak{g}) $-modules. More precisely, we consider the canonical map $ U_q(\mathfrak{g}) \rightarrow \prod \End(V) \subset \End(\bigoplus V) $ 
where $ V $ runs over all integrable $ U_q(\mathfrak{g}) $-modules. According to Theorem \ref{uqgfdseparation} this map is an embedding, 
and we shall identify $ X \in U_q(\mathfrak{g}) $ with its image in $ \prod \End(V) $. 
We then define $ \T_i(X) \in \prod \End(V) $ by setting  
$$ 
\T_i(X) = \T_i X \T_i^{-1}, 
$$ 
where the operators $ \T_i \in \prod \End(V) $ are  defined as in the previous subsection. 
In order to show that this prescription yields a well-defined automorphism of $ U_q(\mathfrak{g}) $ it is enough to check that conjugation 
with $ \T_i \in \prod \End(V) $  preserves the action of $ U_q(\mathfrak{g}) $. In fact, it is sufficient to verify this for the standard 
generators of $ U_q(\mathfrak{g}) $. 

For the generators $ K_\mu $ for $ \mu \in \weights $ the claim follows from elementary weight considerations, indeed the corresponding 
relation $ \T_i K_\mu \T_i^{-1} = K_{s_i \mu} $ was already pointed out in the discussion before Lemma \ref{Ticommutationeasy}. 
For the generators $ E_i, F_i $ the assertion follows from Proposition \ref{tiformulas}. \end{proof}

Next we shall verify the braid group relations. 

\begin{theorem} \label{uqgbraidrelations}
The automorphisms $ \T_i $ of $ U_q(\mathfrak{g}) $ satisfy the braid relations in $ B_\mathfrak{g} $, that is, 
$$
\T_i \T_j \T_i \cdots = \T_j \T_i \T_j \cdots 
$$
for all $ 1 \leq i, j \leq N $ such that $ i \neq j $, with $ m_{ij} $ operators on each side of the equation. 
In other words, the maps $ \T_i: U_q(\mathfrak{g}) \rightarrow U_q(\mathfrak{g}) $ induce an action 
of the braid group $ B_{\mathfrak{g}} $ on $ U_q(\mathfrak{g}) $ by algebra automorphisms.
\end{theorem} 

\begin{proof} Recall that $ m_{ij} = 2,3,4,6 $ iff $ a_{ij} a_{ji} = 0,1,2,3 $, respectively. 
We have $ a_{ij} = 0 = a_{ji} $ if $ m_{ij} = 2 $, and $ a_{ij} = -1 = a_{ji} $ if $ m_{ij} = 3 $. Note that 
we may assume that $ a_{ij} = -2, a_{ji} = -1 $ and $ a_{ij} = -3, a_{ji} = -1 $ in the cases $ m_{ij} = 4, 6 $, respectively.
The cases $ m = m_{ij} = 2,3,4,6 $ are verified separately by explicit calculations, compare \cite{Lusztigbook} and \cite{Jantzenqg}.
\begin{enumerate} 
\item[$ m = 2 $:] This corresponds to $ (\alpha_i, \alpha_j) = 0 $. In this case we clearly 
have $ \T_i \T_j = \T_j \T_i $ on modules because the operators $ E_i, F_i $ commute pairwise with $ E_j, F_j $. 
Since the action of $ \T_i $ on $ U_q(\mathfrak{g}) $ is obtained by conjugating with the operator $ \T_i $ on modules, 
it follows that $ \T_i \T_j = \T_j \T_i $ as automorphisms of $ U_q(\mathfrak{g}) $. 
\item[$ m = 3 $:] Since $ a_{ij} = - 1 = a_{ji} $ we have $ q_i = q_j $, 
and the formulas in Proposition \ref{tiformulas} give 
\begin{align*}
\T_i(E_j) &= E_j E_i - q_i E_i E_j = \T_j^{-1}(E_i) \\ 
\T_i(F_j) &= F_i F_j - q_i^{-1} F_j F_i = \T_j^{-1}(F_i).  
\end{align*}
Symmetrically, we obtain 
\begin{align*}
\T_j(E_i) &= E_i E_j - q_i E_j E_i = \T_i^{-1}(E_j) \\ 
\T_j(F_i) &= F_j F_i - q_i^{-1} F_i F_j = \T_i^{-1}(F_j). 
\end{align*}
It suffices to verify $ \T_i \T_j \T_i(X) = \T_j \T_i \T_j(X) $ for $ X = E_k, F_k $ and $ k = 1, \dots, N $ as well as $ X = K_\lambda $ 
for $ \lambda \in \weights $. For the generators $ K_\lambda $ the claim follows from the fact that the action of the Weyl group on $ \mathfrak{h}^* $ 
satisfies the braid relations. 

Let us consider $ X = E_j $. The above relations imply $ \T_j \T_i(E_j) = E_i $, and hence we obtain 
$$
\T_i \T_j \T_i(E_j) = \T_i(E_i) = -K_i F_i.  
$$
On the other hand we get 
$$
\T_j \T_i \T_j(E_j) = \T_j \T_i(-K_j F_j) = -\T_j \T_i(K_j) \T_j \T_i(F_j)  = - K_i F_i 
$$
since the above relations imply $ \T_j \T_i(F_j) = F_i $, and moreover we have $ \T_j \T_i(K_j) = K_i $. For the latter use  
\begin{align*}
s_j s_i(\alpha_j) &= s_j(\alpha_j - (\alpha_i^\vee, \alpha_j) \alpha_i) \\
&= s_j(\alpha_j + \alpha_i) = - \alpha_j + \alpha_j + \alpha_i = \alpha_i. 
\end{align*}
For $ X = F_j $ analogous considerations give 
$$
\T_i \T_j \T_i(F_j) = \T_i(F_i) = - E_i K_i^{-1}
$$
and 
$$
\T_j \T_i \T_j(F_j) = \T_j \T_i(-E_j K_j^{-1}) = - \T_j \T_i(E_j) \T_j \T_i(K_j^{-1}) = - E_i K_i^{-1}. 
$$
Since the situation is symmetric in $ i $ and $ j $ we also get the claim for $ X = E_i, F_i $. 

Therefore it remains to consider the case $ X = E_k, F_k $ for $ k \neq i, j $. 
For $ k \neq i,j $ we have $ a_{ik} = 0 $ or $ a_{jk} = 0 $. Again by symmetry we may assume $ a_{jk} = 0 $. Then $ [E_k, E_j] = 0 = [E_k, F_j] $. 
From the construction of the automorphism $ \T_j $ using conjugation by the operator $ \T_j $ we therefore get $ \T_j(E_k) = E_k $. 
Moreover recall that $ \T_j \T_i(E_j) = E_i $ and $ \T_j \T_i(F_j) = F_i $. Hence we also get 
\begin{align*}
[\T_j \T_i(E_k), E_i] &= \T_j \T_i([E_k, E_j]) = 0, \\
[\T_j \T_i(E_k), F_i] &= \T_j \T_i([E_k, F_j]) = 0. 
\end{align*} 
Again by the construction of $ \T_i $ we therefore get $ \T_i \T_j \T_i(E_k) = \T_j \T_i(E_k) $. Combining these considerations yields 
$$
\T_i \T_j \T_i(E_k) = \T_j \T_i(E_k) = \T_j \T_i \T_j(E_k)
$$
as desired. Upon replacing $ E_k $ by $ F_k $, the above argument shows 
$$ 
\T_i \T_j \T_i(F_k) = \T_j \T_i(F_k) = \T_j \T_i \T_j(F_k)
$$
as well. 
\item[$ m = 4 $:] Let us assume $ a_{ij} = -2, a_{ji} = -1 $. In this case we have $ q_i = q $ and $ q_j = q^2 $. 
As in the case $ m = 3 $, Proposition \ref{tiformulas} implies
\begin{align*}
\T_j(E_i) &= E_i E_j - q^2 E_j E_i, \\
\T_j^{-1}(E_i) &= E_j E_i - q^2 E_i E_j. 
\end{align*}
Moreover, taking the first formula of Proposition \ref{tiformulas} for $ l = 1 $ gives 
$$
\T_i(E_i E_j - q^2 E_j E_i) = E_j E_i - q^2 E_i E_j.
$$
Hence $ \T_i \T_j(E_i) = \T_j^{-1}(E_i) $, or equivalently $ \T_j \T_i \T_j(E_i) = E_i $. 

Swapping the roles of $ i $ and $ j $, Proposition \ref{tiformulas} gives 
\begin{align*}
\T_i(E_j) &= E_j E_i^{(2)} - q E_i E_j E_i + q^2 E_i^{(2)} E_j, \\
\T_i^{-1}(E_j) &= E_i^{(2)} E_j - q E_i E_j E_i + q^2 E_j E_i^{(2)}.  
\end{align*}
Multiplying these equations with $ [2]_{q_i} = q_i + q_i^{-1} = q + q^{-1} $ we obtain 
\begin{align*}
[2]_{q_i} \T_i(E_j) &= E_j E_i^2 - (1 + q^2) E_i E_j E_i + q^2 E_i^2 E_j \\
&= \T_j^{-1}(E_i) E_i - E_i \T_j^{-1}(E_i)
\end{align*}
and 
\begin{align*}
[2]_{q_i} \T_i^{-1}(E_j) &= E_i^2 E_j - (1 + q^2) E_i E_j E_i + q^2 E_j E_i^2 \\ 
&= E_i \T_j(E_i) - \T_j(E_i) E_i.   
\end{align*}
Combining these equations yields $ \T_j \T_i(E_j) = \T_i^{-1}(E_j) $, or equivalently $ \T_i \T_j \T_i(E_j) = E_j $. 

In a completely analogous fashion we obtain 
$$ 
\T_i \T_j \T_i(F_j) = F_j, \qquad \T_j \T_i \T_j(F_i) = F_i. 
$$
Let us now show $ \T_i \T_j \T_i \T_j(X) = \T_j \T_i \T_j \T_i(X) $ for $ X = E_k, F_k $ where $ k = 1, \dots, N $. Since 
\begin{align*}
s_j s_i s_j(\alpha_i) &= s_j s_i(\alpha_i - a_{ji} \alpha_j) \\
&= s_j(-\alpha_i + \alpha_j - a_{ij} \alpha_i) \\
&= s_j(\alpha_i + \alpha_j) = \alpha_i + \alpha_j - \alpha_j = \alpha_i 
\end{align*}
we have $ \T_j \T_i \T_j(K_i) = K_i $, and hence our above computations imply
\begin{align*}
\T_i \T_j \T_i \T_j(E_i) = \T_i(E_i) 
&= - K_i F_i \\
&= - \T_j \T_i \T_j(K_i) \T_j \T_i \T_j(F_i) \\ 
&= \T_j \T_i \T_j (-K_i F_i) \\
&= \T_j \T_i \T_j \T_i(E_i). 
\end{align*}
In a similar way one obtains the assertion for $ X = E_j $ and $ X = F_i, F_j $. 

Consider now $ k \neq i,j $. Then we have $ (\alpha_i^\vee, \alpha_k) = 0 $ or $ (\alpha_j^\vee, \alpha_k) = 0 $. 
If $ (\alpha_j^\vee, \alpha_k) = 0 $ we obtain $ [E_k, E_j] = 0 = [E_k, F_j] $ and hence $ \T_j(E_k) = E_k $. Moreover, our above computations yield 
\begin{align*}
[\T_i \T_j \T_i(E_k), E_j] &= \T_i \T_j \T_i([E_k, E_j]) = 0, \\
[\T_i \T_j \T_i(E_k), F_j] &= \T_i \T_j \T_i([E_k, F_j]) = 0. 
\end{align*} 
This implies $ \T_j \T_i \T_j \T_i(E_k) = \T_i \T_j \T_i(E_k) $, and thus 
\begin{align*}
\T_j \T_i \T_j \T_i(E_k) = \T_i \T_j \T_i(E_k) = \T_i \T_j \T_i \T_j(E_k). 
\end{align*} 
Replacing $ E_k $ by $ F_k $ in this argument yields the claim for $ X = F_k $. Finally, if $ (\alpha_i^\vee, \alpha_k) = 0 $ we swap the roles of $ i $ and $ j $ 
in the above argument, and obtain the assertion in the same way. 
\item[$ m = 6 $:] 
Let us assume $ a_{ij} = -3, a_{ji} = -1 $. In this case we have $ q_i = q $ and $ q_j = q^3 $. 
Given $ 0 \leq l \leq 3 $, let us abbreviate 
\begin{align*}
E^-_{ij}(l) &= \sum_{k = 0}^l (-1)^k q_i^{k (4 - l)} E_i^{(k)} E_j E_i^{(l - k)}, \\
E^+_{ij}(l) &= \sum_{k = 0}^l (-1)^k q_i^{k (4 - l)} E_i^{(l - k)} E_j E_i^{(k)}. 
\end{align*} 
As in the case $ m = 3 $, Proposition \ref{tiformulas} implies
\begin{align*}
\T_j(E_i) &= E_i E_j - q^3 E_j E_i = E^+_{ij}(1) \\
\T_j^{-1}(E_i) &= E_j E_i - q^3 E_i E_j = E^-_{ij}(1). 
\end{align*}
Moreover the relations from Proposition \ref{tiformulas} give
\begin{align*}
E_j &= \T_i(E_{ij}^+(3)), \\
E_{ij}^-(1) &= \T_i(E_{ij}^+(2)), \\
E_{ij}^-(2) &= \T_i(E_{ij}^+(1)), \\
E_{ij}^-(3) &= \T_i(E_j). 
\end{align*} 
Explicitly, we have 
\begin{align*}
E_{ij}^+(3) &= E_i^{(3)} E_j - q E_i^{(2)} E_j E_i + q^2 E_i E_j E_i^{(2)} - q^3 E_j E_i^{(3)}, \\ 
E_{ij}^-(3) &= E_j E_i^{(3)} - q E_i E_j E_i^{(2)} + q^2 E_i^{(2)} E_j E_i - q^3 E_i^{(3)} E_j, \\
E_{ij}^+(2) &= E_i^{(2)} E_j - q^2 E_i E_j E_i + q^4 E_j E_i^{(2)}, \\ 
E_{ij}^-(2) &= E_j E_i^{(2)} - q^2 E_i E_j E_i + q^4 E_i^{(2)} E_j.
\end{align*}
Hence we obtain
\begin{align*}
[3]_q E_{ij}^+(3) &= [2]_q^{-1} E_i^3 E_j - q [3]_q E_i^{(2)} E_j E_i \\
&\qquad + q^2 [3]_q E_i E_j E_i^{(2)} - q^3 [2]_q^{-1} E_j E_i^3 \\
&= E_i(E_i^{(2)} E_j - q^2 E_i E_j E_i + q^4 E_j E_i^{(2)}) \\
&\qquad - q^{-1} (E_i^{(2)} E_j - q^2 E_i E_j E_i + q^4 E_j E_i^{(2)}) E_i \\
&= E_i E_{ij}^+(2) - q^{-1} E_{ij}^+(2) E_i
\end{align*}
using $ q(q^2 + 1 + q^{-2}) = q^2(q + q^{-1}) + q^{-1} $ and $ q^2(q^2 + 1 + q^{-2}) = q^4 + q (q + q^{-1}) $. \\ 
Similarly, 
\begin{align*}
[3]_q E_{ij}^-(3) &= [2]_q^{-1} E_j E_i^3 - q [3]_q E_i E_j E_i^{(2)} \\
&\qquad + q^2 [3]_q E_i^{(2)} E_j E_i - q^3 [2]_q^{-1} E_i^3 E_j \\
&= (E_j E_i^{(2)} - q^2 E_i E_j E_i + q^4 E_i^{(2)} E_j) E_i \\
&\qquad - q^{-1} E_i(E_j E_i^{(2)} - q^2 E_i E_j E_i + q^4 E_i^{(2)} E_j) \\
&= E_{ij}^-(2) E_i - q^{-1} E_i E_{ij}^-(2). 
\end{align*}
Let us also record
\begin{align*}
[2]_q E_{ij}^+(2) &= E_i^2 E_j - q^2 (q + q^{-1}) E_i E_j E_i + q^4 E_j E_i^2 \\
&= E_i E_{ij}^+(1) - q E_{ij}^+(1) E_i
\end{align*}
and 
\begin{align*}
[2]_q E_{ij}^-(2) &= E_j E_i^2 - q^2 (q + q^{-1}) E_i E_j E_i + q^4 E_i^2 E_j \\
&= E_{ij}^-(1) E_i - q E_i E_{ij}^-(1).
\end{align*}
Using $ E^\pm_{ij}(1) = \T_j^{\pm 1}(E_i) $ we obtain 
\begin{align*}
E_{ij}^+(2) &= [2]_q^{-1}(E_i \T_j(E_i) - q \T_j(E_i) E_i) \\
&= [2]_q^{-1} \T_j(E^-_{ij}(1) E_i - q E_i E^-_{ij}(1)) \\
&= \T_j(E_{ij}^-(2)). 
\end{align*}
Now we compute 
\begin{align*}
\T_j \T_i(E_j) &= \T_j(E^-_{ij}(3)) \\
&= [3]_q^{-1}\T_j(E_{ij}^-(2) E_i - q^{-1} E_i E_{ij}^-(2)) \\
&= [3]_q^{-1}(\T_j(E_{ij}^-(2)) E^+_{ij}(1) - q^{-1} E^+_{ij}(1) \T_j^{-1}(E_{ij}^-(2))) \\
&= [3]_q^{-1}(E_{ij}^+(2) E^+_{ij}(1)  - q^{-1} E^+_{ij}(1) E_{ij}^+(2))
\end{align*}
and 
\begin{align*}
\T_j^{-1} \T_i^{-1}(E_j) &= \T_j^{-1}(E^+_{ij}(3)) \\
&= [3]_q^{-1}\T_j^{-1}(E_i E_{ij}^+(2) - q^{-1} E_{ij}^+(2) E_i) \\
&= [3]_q^{-1}(E^-_{ij}(1) \T_j^{-1}(E_{ij}^+(2)) - q^{-1} \T_j^{-1}(E_{ij}^+(2)) E^-_{ij}(1)) \\
&= [3]_q^{-1}(E^-_{ij}(1) E_{ij}^-(2) - q^{-1} E_{ij}^-(2) E^-_{ij}(1)). 
\end{align*}
We thus obtain 
\begin{align*}
\T_i^{-1} \T_j^{-1} \T_i^{-1}(E_j) &= [3]_q^{-1}\T_i^{-1}(E^-_{ij}(1) E_{ij}^-(2) - q^{-1} E_{ij}^-(2) E^-_{ij}(1)) \\
&= [3]_q^{-1}(E_{ij}^+(2) E^+_{ij}(1)  - q^{-1} E^+_{ij}(1) E_{ij}^+(2)) \\
&= \T_j \T_i(E_j), 
\end{align*} 
or equivalently, $ \T_i \T_j \T_i \T_j \T_i(E_j) = E_j $. 

We also have 
\begin{align*}
\T_i \T_j(E_i) &= \T_i(E^+_{ij}(1)) = E^-_{ij}(2)
\end{align*}
and 
\begin{align*}
\T_i^{-1} \T_j^{-1}(E_i) &= \T_i^{-1}(E^-_{ij}(1)) = E^+_{ij}(2). 
\end{align*}
Therefore
\begin{align*}
\T_j \T_i \T_j(E_i) &= \T_j(E^-_{ij}(2)) = E^+_{ij}(2) = \T_i^{-1} \T_j^{-1}(E_i), 
\end{align*}
and hence $ \T_j \T_i \T_j \T_i \T_j(E_i) = E_i $. 

In a completely analogous fashion one calculates 
$$ 
\T_i \T_j \T_i \T_j \T_i(F_j) = F_j, \qquad \T_j \T_i \T_j \T_i \T_j(F_i) = F_i. 
$$
Let us now verify $ \T_i \T_j \T_i \T_j \T_i \T_j(X) = \T_j \T_i \T_j \T_i \T_j \T_i(X) $ for $ X = E_k, F_k $ where $ k = 1, \dots, N $. 
Since 
\begin{align*}
s_j s_i s_j s_i s_j(\alpha_i) &= s_j s_i s_j s_i(\alpha_i - a_{ji} \alpha_j) \\
&= s_j s_i s_j(-\alpha_i + \alpha_j - a_{ij} \alpha_i) \\
&= s_j s_i s_j(2 \alpha_i + \alpha_j) \\
&= s_j(4 \alpha_i + 2\alpha_j - \alpha_j + a_{ij} \alpha_i) \\
&= s_j(\alpha_i + \alpha_j) = \alpha_i + \alpha_j - \alpha_j = \alpha_i 
\end{align*}
we have $ \T_j \T_i \T_j \T_i \T_j(K_i) = K_i $, and hence our above computations imply
\begin{align*}
\T_i \T_j \T_i \T_j \T_i \T_j(E_i) = \T_i(E_i) 
&= - K_i F_i\\
&= - \T_j \T_i \T_j \T_i \T_j(K_i) \T_j \T_i \T_j \T_i \T_j(F_i) \\
&= \T_j \T_i \T_j \T_i \T_j(-K_i F_i) \\ 
&= \T_j \T_i \T_j \T_i \T_j \T_i(E_i). 
\end{align*}
In a similar way one obtains the assertion for $ X = E_j $ and $ X = F_i, F_j $. 

Consider now $ k \neq i,j $. Then we have $ (\alpha_i^\vee, \alpha_k) = 0 $ or $ (\alpha_j^\vee, \alpha_k) = 0 $. 
If $ (\alpha_j^\vee, \alpha_k) = 0 $ we obtain $ [E_k, E_j] = 0 = [E_k, F_j] $ and hence $ \T_j(E_k) = E_k $. Moreover, our above computations yield 
\begin{align*}
[\T_i \T_j \T_i \T_j \T_i(E_k), E_j] &= \T_i \T_j \T_i \T_j \T_i([E_k, E_j]) = 0, \\
[\T_i \T_j \T_i \T_j \T_i(E_k), F_j] &= \T_i \T_j \T_i \T_j \T_i([E_k, F_j]) = 0. 
\end{align*} 
This implies $ \T_j \T_i \T_j \T_i \T_j \T_i(E_k) = \T_i \T_j \T_i \T_j \T_i(E_k) $, and thus 
\begin{align*}
\T_j \T_i \T_j \T_i \T_j \T_i(E_k) = \T_i \T_j \T_i \T_j \T_i(E_k) = \T_i \T_j \T_i \T_j \T_i \T_j(E_k). 
\end{align*} 
Replacing $ E_k $ by $ F_k $ in this argument yields the claim for $ X = F_k $. Finally, if $ (\alpha_i^\vee, \alpha_k) = 0 $ we swap $ i $ and $ j $ 
in the above computations and obtain the assertion in the same way. 
\end{enumerate} 
This finishes the proof. \end{proof}

We remark that the algebra automorphisms $\T_i$ in Theorem \ref{uqgbraidrelations}	are not coalgebra automorphisms.

The action of the braid group $B_\mathfrak{g}$ on $U_q(\mathfrak{g})$ does not descend to an action of the Weyl group $W$, since $\T_i^2 \neq \id$.  Nevertheless, there is a standard way to associate to each $w\in W$ an automorphism $\T_w$, as we now explain.

If $ w = s_{i_1} \cdots s_{i_k} $ is a reduced expression of $ w \in W $, then 
any other reduced expression is obtained by applying a sequence of \emph{elementary moves} to $ s_{i_1} \cdots s_{i_k} $. Here 
an elementary move is the replacement of a substring $ s_i s_j s_i \cdots $ of length $ m_{ij} $ in the given expression by the string $ s_j s_i s_j \cdots $ 
of the same length. We refer to Section 8.1 in \cite{HumphreysReflectiongroups} for the details. 
As a consequence, according to Theorem \ref{uqgbraidrelations} the automorphism 
$$ 
\T_w = \T_{i_1} \cdots \T_{i_k} 
$$ 
\nomenclature[o$T_w$]{$\T_w$}{braid automorphism associated to a Weyl group element $w$}%
of $ U_q(\mathfrak{g}) $ depends only on $ w $ and not on the reduced expression for $w$. 

Note that we have $ \T_w(K_\mu) = K_{w\mu} $ for 
all $ \mu \in \weights $ by Theorem \ref{deftiautomorphism}.

\subsubsection{The Poincar\'e-Birkhoff-Witt Theorem}

Let $ w \in W $ and let us fix a reduced expression $ w = s_{i_1} \cdots s_{i_t} $. If $ w \neq e $ then the roots 
$ \beta_r = s_{i_1} \cdots s_{i_{r - 1}} \alpha_{i_r} $ for $ 1 \leq r \leq t $ are pairwise distinct and positive, see 5.6 in \cite{HumphreysReflectiongroups}. 
We define the root vectors of $ U_q(\mathfrak{g}) $ associated to the reduced expression $ w = s_{i_1} \cdots s_{i_t} $ by 
$$
E_{\beta_r} = \T_{i_1} \cdots \T_{i_{r - 1}}(E_{i_r}), \quad F_{\beta_r} = \T_{i_1} \cdots \T_{i_{r - 1}}(F_{i_r})
$$
\label{nom:quantum_root_vector_w}%
for $ 1 \leq r \leq t $; in the special case $ w = e $ we declare $ 1 $ to be the unique associated root vector. 
Let us point out that these vectors do depend on the choice of the reduced expression for $ w $ in general, see for instance the discussion 
in Section 8.15 of \cite{Jantzenqg}.

Our first aim is to show that the root vectors $ E_{\beta_r} $ associated to a reduced expression for $ w $ are contained in $ U_q(\mathfrak{n}_+) $. 
For this we need some preparations. 

\begin{lemma} \label{PBWpositivehelplemma}
Let $ i \neq j $ and assume $ w $ is contained in the subgroup of $ W $ generated by $ s_i $ and $ s_j $. If $ w \alpha_i \in {\bf \Delta}^+ $ 
then $ \T_w(E_i) $ is contained in the subalgebra of $ U_q(\mathfrak{n}_+) $ generated by $ E_i $ and $ E_j $, and if $ w \alpha_i = \alpha_k $ 
for some $ 1 \leq k \leq N $ then $ \T_w(E_i) = E_k $. 
\end{lemma} 

\begin{proof} Let us remark first that $ w \alpha_i = \alpha_k $ can only happen for $ k = i,j $. \\
Without loss of generality we may assume $ w \neq e $. The proof will proceed case-by-case depending on the order $ m $ 
of $ s_i s_j $ in $ W $. Let us write $ \bra s_i, s_j \ket $ for the subgroup of $ W $ generated by $ s_i, s_j $. 
\begin{enumerate} 
\item[$ m = 2 $:] In this case $ s_i $ and $ s_j $ commute, and the only nontrivial element $ w \in \bra s_i, s_j \ket $ 
satisfying $ w \alpha_i \in {\bf \Delta}^+ $ is $ w = s_j $; note that $ s_j \alpha_i = \alpha_i $. 
The formula in Theorem \ref{deftiautomorphism} gives $ \T_{s_j}(E_i) = \T_j(E_i) = E_i $. 
\item[$ m = 3 $:] In this case we have $ s_i s_j s_i = s_j s_i s_j $, and the nontrivial elements $ w \in \bra s_i, s_j \ket $ satisfying 
$ w \alpha_i \in {\bf \Delta}^+ $ are $ w = s_j $ and $ w = s_i s_j $; note that in this case $ \bra s_i, s_j \ket \cong S_3 $. 
Using the computations in the proof of Theorem \ref{uqgbraidrelations} we obtain $ \T_{s_j}(E_i) = \T_j(E_i) = E_i E_j - q_i E_j E_i $ and 
$ \T_{s_i s_j}(E_i) = \T_i \T_j(E_i) = E_j $. 
\item[$ m = 4 $:] In this case we have $ s_i s_j s_i s_j = s_j s_i s_j s_i $, and the nontrivial elements $ w \in \bra s_i, s_j \ket $ satisfying 
$ w \alpha_i \in {\bf \Delta}^+ $ are $ w = s_j, w = s_i s_j, w = s_j s_i s_j $; in fact applying these elements to $ \alpha_i $ yields
two distinct positive roots. 
Note that we have here $ \bra s_i, s_j \ket \cong D_4 = S_2 \wr S_2 $, the dihedral group of order $ 8 $. 
This is the Weyl group of type $ B_2 $. 

If $ (\alpha_i, \alpha_i) < (\alpha_j, \alpha_j) $, then as in the proof of Theorem \ref{uqgbraidrelations} we obtain 
\begin{align*} 
\T_{s_j}(E_i) &= \T_j(E_i) = E_i E_j - q^2 E_j E_i, \\
\T_{s_i s_j}(E_i) &= \T_i \T_j(E_i) = \T_j^{-1}(E_i) = E_j E_i - q^2 E_i E_j, \\
\T_{s_j s_i s_j}(E_i) &= \T_j \T_i \T_j(E_i) = E_i 
\end{align*}
If $ (\alpha_i, \alpha_i) > (\alpha_j, \alpha_j) $, we have to swap the roles of $ i $ and $ j $ in the proof of Theorem \ref{uqgbraidrelations} and obtain 
\begin{align*} 
\T_{s_j}(E_i) &= \T_j(E_i) = E_i E_j^{(2)} - q E_j E_i E_j + q^2 E_j^{(2)} E_i, \\
\T_{s_i s_j}(E_i) &= \T_i \T_j(E_i) = \T_j^{-1}(E_i) = E_j^{(2)} E_i - q E_j E_i E_j + q^2 E_i E_j^{(2)}, \\
\T_{s_j s_i s_j}(E_i) &= \T_j \T_i \T_j(E_i) = E_i 
\end{align*}
In both cases the required properties hold. 
\item[$ m = 6 $:] In this case we have $ s_i s_j s_i s_j s_i s_j = s_j s_i s_j s_i s_j s_i $, and for the nontrivial elements $ w \in \bra s_i, s_j \ket $ satisfying 
$ w \alpha_i \in {\bf \Delta}^+ $ one finds $ w = s_j, w = s_i s_j, w = s_j s_i s_j, w = s_i s_j s_i s_j, w = s_j s_i s_j s_i s_j $; in fact 
applying these elements to $ \alpha_i $ yields three distinct positive roots. 
Note that we have $ \bra s_i, s_j \ket \cong D_6 $, the dihedral group 
of order $ 12 $. This is the Weyl group of type $ G_2 $. 

If $ (\alpha_i, \alpha_i) < (\alpha_j, \alpha_j) $, then using the notation in the proof of Theorem \ref{uqgbraidrelations} we obtain 
\begin{align*} 
\T_{s_j}(E_i) &= \T_j(E_i) = E^+_{ij}(1), \\
\T_{s_i s_j}(E_i) &= \T_i \T_j(E_i) = E^-_{ij}(2), \\
\T_{s_j s_i s_j}(E_i) &= \T_j \T_i \T_j(E_i) = E^+_{ij}(2), \\
\T_{s_i s_j s_i s_j}(E_i) &= \T_i \T_j \T_i \T_j(E_i) = \T_j^{-1}(E_i) = E^-_{ij}(1), \\
\T_{s_j s_i s_j s_i s_j}(E_i) &= \T_j \T_i \T_j \T_i \T_j(E_i) = E_i. 
\end{align*}
If $ (\alpha_i, \alpha_i) > (\alpha_j, \alpha_j) $, we have to swap the roles of $ i $ and $ j $ in the proof of Theorem \ref{uqgbraidrelations} and obtain 
\begin{align*} 
\T_{s_j}(E_i) &= \T_j(E_i) = E^-_{ji}(3), \\
\T_{s_i s_j}(E_i) &= \T_i \T_j(E_i) = [3]_q^{-1}(E^+_{ji}(2) E_{ji}^+(1) - q^{-1} E_{ji}^+(1) E^+_{ji}(2)), \\
\T_{s_j s_i s_j}(E_i) &= \T_j \T_i \T_j(E_i) = [3]_q^{-1}(E_{ji}^-(1) E^-_{ji}(2)  - q^{-1} E^-_{ji}(2) E_{ji}^-(1)) \\
\T_{s_i s_j s_i s_j}(E_i) &= \T_i \T_j \T_i \T_j(E_i) = \T_j^{-1}(E_i) = E^+_{ji}(3) \\
\T_{s_j s_i s_j s_i s_j}(E_i) &= \T_j \T_i \T_j \T_i \T_j(E_i) = E_i. 
\end{align*}
In both cases the required properties hold. 
\end{enumerate}
This finishes the proof. \end{proof} 

\begin{lemma} \label{PBWpositivelemma}
Let $ w \in W $. If $ w \alpha_i \in {\bf \Delta}^+ $ then $ \T_w(E_i) \in U_q(\mathfrak{n}_+) $, and if $ w \alpha_i = \alpha_k $ 
for some $ 1 \leq k \leq N $ then $ \T_w(E_i) = E_k $. 
\end{lemma} 

\begin{proof} We use induction on the length $ l(w) $ of $ w $, the case $ l(w) = 0 $ being trivial. Hence suppose $ l(w) > 0 $. 
Then according to Section 10.2 in \cite{HumphreysLie} there exists $ 1 \leq j \leq N $ such that $ w \alpha_j \in {\bf \Delta}^- $, 
note in particular that we must have $ i \neq j $. Let us write again $ \bra s_i, s_j \ket \subset W $ for the subgroup generated by $ s_i $ and $ s_j $. 
According to Section 1.10 in \cite{HumphreysReflectiongroups} we find $ u \in W $ and $ v \in \bra s_i, s_j \ket $ such 
that $ w = uv $ and $ u \alpha_i, u \alpha_j \in {\bf \Delta}^+ $, with lengths satisfying $ l(w) = l(u) + l(v) $. 
In particular, since $ u $ maps both $ \alpha_i $ and $ \alpha_j $ into $ {\bf \Delta}^+ $ we have $ u \neq w $, and we conclude $ l(u) < l(w) $. 
Applying the inductive hypothesis to $ u $ we get $ \T_u(E_i) \in U_q(\mathfrak{n}_+) $ and $ \T_u(E_j) \in U_q(\mathfrak{n}_+) $. 
A similar reasoning as above shows $ v \alpha_i \in {\bf \Delta}^+ $ since $ w \alpha_i \in {\bf \Delta}^+ $. 
Hence Lemma \ref{PBWpositivehelplemma} implies that $ \T_v(E_i) $ is contained in the subalgebra of $ U_q(\mathfrak{n}_+) $ generated by $ E_i $ and $ E_j $. 
We conclude $ \T_w(E_i) = \T_u \T_v(E_i) \in U_q(\mathfrak{n}_+) $ as desired. 

For the second claim we proceed again by induction on $ l(w) $. Using the above considerations, 
it suffices to show in the inductive step that $ v \alpha_i $ is a simple root. Indeed, the second part of 
Lemma \ref{PBWpositivehelplemma} will then yield the claim. However, if $ v \alpha_i $ is not simple we can write 
$ v \alpha_i = m \alpha_i + n \alpha_j $ with $ m, n \in \mathbb{N} $ because $ v \in \bra s_i, s_j \ket $ and 
$ v \alpha_i \in {\bf \Delta}^+ $. 
Then $ w \alpha_i = m u \alpha_i + n u \alpha_j $ is a sum of positive roots because $ u \alpha_i, u \alpha_j \in {\bf \Delta}^+ $. 
This contradicts the assumption that $ w \alpha_i = \alpha_k $ is simple. Hence $ v \alpha_i $ is simple, and this finishes the proof. \end{proof}

If $ w \in W $ and $ w = s_{i_1} \cdots s_{i_t} $ is a reduced expression of $ w $ then we shall call the vectors $ E^{a_1}_{\beta_1} \cdots E^{a_t}_{\beta_t} $ 
with $ a_j \in \mathbb{N}_0 $ for all $ j $ the associated PBW-vectors. We will show next that these vectors are always contained in $ U_q(\mathfrak{n}_+) $. 

\begin{prop} \label{PBWpositive}
Let $ w \in W $ and $ w = s_{i_1} \cdots s_{i_t} $ be a reduced expression of $ w $. Then the associated 
PBW-vectors $ E^{a_1}_{\beta_1} \cdots E^{a_t}_{\beta_t} $ are contained in $ U_q(\mathfrak{n}_+) $.
\end{prop}

\begin{proof} Since $ U_q(\mathfrak{n}_+) $ is a subalgebra of $ U_q(\mathfrak{g}) $ it suffices to show that $ E_{\beta_r} $ is 
contained in $ U_q(\mathfrak{n}_+) $ for all $ 1 \leq r \leq t $. As indicated earlier on, $ s_{i_1} \cdots s_{i_{r - 1}} \alpha_{i_r} $ 
is a positive root. Considering the element $ s_{i_1} \cdots s_{i_{r - 1}} \in W $ and $ i = i_r $ in Lemma \ref{PBWpositivelemma} 
we see that $ \T_{s_{i_1} \cdots s_{i_{r - 1}}}(E_{i_r}) $ is contained in $ U_q(\mathfrak{n}_+) $. 
Moreover, observe that $ s_{i_1} \cdots s_{i_{r - 1}} $ is a reduced expression since it is part of a reduced expression. It follows that 
$ \T_{s_{i_1} \cdots s_{i_{r - 1}}} = \T_{s_{i_1}} \cdots \T_{s_{i_{r - 1}}} $, and hence  
$ \T_{s_{i_1} \cdots s_{i_{r - 1}}}(E_{i_r}) = E_{\beta_r} $. This finishes the proof. \end{proof} 

Next we discuss linear independence. 

\begin{prop} \label{PBWindependence}
Let $ w \in W $ and $ w = s_{i_1} \cdots s_{i_t} $ be a reduced expression of $ w $. Then the associated PBW-vectors $ E^{a_1}_{\beta_1} \cdots E^{a_t}_{\beta_t} $ are linearly 
independent.
\end{prop} 

\begin{proof} As a preparation, assume first that $ X_1, \dots, X_m \in U_q(\mathfrak{n}_+) $ satisfy $ \T_i^{-1}(X_j) \in U_q(\mathfrak{n}_+) $ 
for all $ j $ and $ \sum_j X_j E_i^j = 0 $ or $ \sum_j E_i^j X_j = 0 $. We claim that $ X_j = 0 $ for $ j = 1, \dots, m $. 

In order to verify this consider the case $ \sum_j E_i^j X_j = 0 $. Applying $ \T_i^{-1} $ yields 
\begin{align*}
0 = \T_i^{-1}\biggl(\sum_{j = 1}^m E_i^j X_j \biggr) = \sum_{j = 1}^m \T_i^{-1}(E_i)^j \T_i^{-1}(X_j) 
= \sum_{j = 1}^m (-F_i K_i^{-1})^j \T_i^{-1}(X_j).  
\end{align*} 
The vectors $ (-F_i K_i^{-1})^j \in U_q(\mathfrak{b}_-) $ have pairwise distinct weights and $ \T_i^{-1}(X_j) \in U_q(\mathfrak{n}_+) $ 
by assumption. Hence the triangular decomposition of Proposition \ref{uqgtriangular} yields the claim. The case $ \sum_j X_j E_i^j = 0 $ is analogous. 

Let us now use induction on $ t $ to prove the Proposition. For $ t = 0 $, there is only one PBW-vector, namely the element $ 1 $, corresponding 
to the empty word. The assertion clearly holds in this case. Assume now that $ t > 0 $. Then the PBW-vectors have the 
form $ E_{i_1}^k \T_{i_1}(X_{kl})  $ where $ k \in \mathbb{N}_0 $ and $ X_{kl} $ is a PBW-vector for the reduced 
expression $ s_{i_2} \cdots s_{i_t} $ of $ s_{i_1} w $. 
If $ c_{kl} $ are scalars with $ \sum_{k,l} c_{kl} E_{i_1}^k \T_{i_1}(X_{kl}) = 0 $, then the first part of our proof 
combined with Proposition \ref{PBWpositive} implies $ \sum_l c_{kl} X_{kl} = 0 $ 
for all $ k $. By the inductive hypothesis we get $ c_{kl} = 0 $ for all $ k,l $, and this finishes the proof. \end{proof} 

Let $ w \in W $ and $ w = s_{i_1} \cdots s_{i_t} $ be a reduced expression of $ w $. We show next that the subspace of $ U_q(\mathfrak{n}_+) $ 
spanned by the associated PBW-vectors is independent of the reduced expression of $ w $. 

\begin{lemma} \label{PBWexpressionhelplemma}
If $ i \neq j $ and $ w $ is the longest element in the subgroup of $ W $ generated by $ s_i $ and $ s_j $, then the PBW-vectors associated to 
a reduced expression of $ w $ span the subalgebra of $ U_q(\mathfrak{n}_+) $ generated by $ E_i $ and $ E_j $. 
\end{lemma} 

\begin{proof} Note that there are precisely two reduced expressions for $ w $, namely $ w = s_i s_j s_i \cdots $ and $ w = s_j s_i s_j \cdots $. Both 
consist of $ m $ factors, where $ m $ is the order of $ s_i s_j $. We shall proceed case-by-case depending on $ m $. 
\begin{enumerate} 
\item[$ m = 2 $:] In this case we have $ \T_j(E_i) = E_i, \T_i(E_j) = E_j $, the two reduced expressions 
for $ w $ have associated PBW-vectors $ E_i^{r_1} E_j^{r_2} $ and $ E_j^{t_1} E_i^{t_2} $, respectively, with $ r_k, t_k \in \mathbb{N}_0 $. 
Since $ E_i $ and $ E_j $ commute, the span is in both cases equal to the subalgebra of $ U_q(\mathfrak{n}_+) $ generated by $ E_i $ and $ E_j $. 
\item[$ m = 3 $:] In this case the situation is symmetric in $ i $ and $ j $. Hence it suffices to show that the PBW-vectors 
associated with $ s_j s_i s_j $ span the subalgebra generated by $ E_i $ and $ E_j $. 
From the proof of Theorem \ref{uqgbraidrelations} we have $ \T_{s_j}(E_i) = \T_j(E_i) = E_i E_j - q_i E_j E_i $ 
and $ \T_{s_j s_i}(E_j) = \T_j \T_i(E_j) = E_i $. 
Then the claim is that the vectors $ E_j^{r_1} \T_j(E_i)^{r_2} E_i^{r_3} $ with $ r_k \in \mathbb{N}_0 $ 
span the subalgebra generated by $ E_i $ and $ E_j $. 

It suffices to show that the linear span of these vectors is closed under left multiplication by $ E_i $ and $ E_j $. For $ E_j $ 
this is trivial, so let us consider $ E_i $. Using the quantum Serre relations we get 
\begin{align*}
E_i \T_j(E_i) &= E_i^2 E_j - q_i E_i E_j E_i \\
&= q_i^{-1} E_i E_j E_i - E_j E_i^2 = q_i^{-1} \T_j(E_i) E_i,  
\end{align*}
and we have 
$$ 
E_i E_j = q_i E_j E_i + \T_j(E_i) 
$$
by definition. Applying these relations iteratively we see that any term 
of the form $ E_i E_j^{r_1} \T_j(E_i)^{r_2} E_i^{r_3} $ can again be written as a linear combination of such terms. 
\item[$ m = 4 $:] As in the proof of Theorem \ref{uqgbraidrelations} we shall use the convention $ a_{ij} = -2, a_{ji} = -1 $. 
We have $ q_i = q $ and $ q_j = q^2 $ in this case. 

For the reduced expression $ w = s_j s_i s_j s_i $ we obtain the associated root 
vectors $ E_j, \T_j(E_i), \T_j \T_i(E_j) = \T_i^{-1}(E_j), \T_j \T_i \T_j(E_i) = E_i $. 
As before, it suffices to check that left multiplication by $ E_i $ preserves the linear span of all terms 
of the form 
$$ 
E_j^{r_1} \T_j(E_i^{r_2}) \T_i^{-1}(E_j^{r_3}) E_i^{r_4} 
$$
with $ r_k \in \mathbb{N}_0 $. \\
To this end we observe 
$$
E_i E_j = q^2 E_j E_i + \T_j(E_i), 
$$
and 
\begin{align*}
E_i \T_j(E_i) &= \T_j(E_i) E_i + [2]_q \T_i^{-1}(E_j).  
\end{align*}
The quantum Serre relations give 
\begin{align*}
[2]_q E_i \T_i^{-1}(E_j) &= E_i^3 E_j - (q + q^{-1})q E_i^2 E_j E_i + q^2 E_i E_j E_i^2 \\
&= q^{-2} E_i^2 E_j E_i - (q + q^{-1})q^{-1} E_i E_j E_i^2 + E_j E_i^3 \\
&= q^{-2} [2]_q \T_i^{-1}(E_j) E_i.  
\end{align*}
Moreover we have 
$$
E_i \T_j^{-1}(E_j) = -E_i F_j K_j^{-1} = - q_j^{a_{ji}} F_j K_j^{-1} E_i = q^{-2} \T_j^{-1}(E_j) E_i 
$$
and hence 
$$
\T_j(E_i) E_j = q^{-2} E_j \T_j(E_i). 
$$
Finally, the quantum Serre relations give 
\begin{align*}
[2]_q \T_i(E_j) E_i &= E_j E_i^3 - (q + q^{-1})q E_i E_j E_i^2 + q^2 E_i^2 E_j E_i \\
&= q^{-2} E_i E_j E_i^2 - (q + q^{-1})q^{-1} E_i^2 E_j E_i + E_i^3 E_j \\
&= q^{-2} [2]_q E_i \T_i(E_j)  
\end{align*}
and hence using $ \T_j^{-1} \T_i^{-1}(E_j) = \T_i(E_j) $ we get 
$$
\T_i^{-1}(E_j) \T_j(E_i) = q^{-2} \T_j(E_i) \T_i^{-1}(E_j).  
$$
Combining these commutation relations yields the assertion for $ w = s_j s_i s_j s_i $. 

For the reduced expression $ w = s_i s_j s_i s_j $ we obtain the associated root 
vectors $ E_i, \T_i(E_j), \T_i \T_j(E_i) = \T_j^{-1}(E_i), \T_i \T_j \T_i(E_j) = E_j $. 
We check that right multiplication by $ E_i $ preserves the linear span of all terms of the form 
$$ 
E_i^{r_1} \T_i(E_j^{r_2}) \T_j^{-1}(E_i^{r_3}) E_j^{r_4} 
$$
with $ r_k \in \mathbb{N}_0 $. 

To this end we observe 
$$
E_j E_i = q^2 E_i E_j + \T_j^{-1}(E_i), 
$$
and 
\begin{align*}
\T_j^{-1}(E_i) E_i &= E_i \T_j^{-1}(E_i) + [2]_q \T_i(E_j).  
\end{align*}
As above, the quantum Serre relations give 
\begin{align*}
[2]_q \T_i(E_j) E_i &= q^{-2} [2]_q E_i \T_i(E_j)
\end{align*}
and 
\begin{align*}
[2]_q E_i \T_i^{-1}(E_j) &= q^{-2} [2]_q \T_i^{-1}(E_j) E_i, 
\end{align*}
the second of which implies 
$$
\T_j^{-1}(E_i) \T_i(E_j) = q^{-2} \T_i(E_j) \T_j^{-1}(E_i). 
$$
Finally, we have the relation 
$$
\T_j(E_j) E_i = - F_j K_j E_i = - q_j^{a_{ji}} E_i F_j K_j = q^{-2} E_i \T_j(E_j) 
$$
which implies 
$$
E_j \T_j^{-1}(E_i) = q^{-2} \T_j^{-1}(E_i) E_j. 
$$
Again, combining these relations yields the assertion. 
\item[$ m = 6 $:] As in the proof of Theorem \ref{uqgbraidrelations} we shall use the convention $ a_{ij} = -3, a_{ji} = -1 $. 
We have $ q_i = q $ and $ q_j = q^3 $ in this case. 

Consider first the reduced expression $ w = s_i s_j s_i s_j s_i s_j $.  
We obtain the associated root vectors $ E_i, \T_i(E_j) = E^-_{ij}(3), \T_i \T_j(E_i) = E^-_{ij}(2), \T_i \T_j \T_i(E_j) = \T_j^{-1} \T_i^{-1}(E_j), 
\T_i \T_j \T_i \T_j(E_i) = \T_j^{-1}(E_i) = E^-_{ij}(1), E_j $. 
We check that right multiplication by $ E_i $ preserves the linear span of the vectors 
$$ 
E_i^{r_1} E^-_{ij}(3)^{r_2} E^-_{ij}(2)^{r_4} \T_j^{-1} \T_i^{-1}(E_j) E^-_{ij}(1)^{r_5} E_j^{r_6}
$$
with $ r_j \in \mathbb{N}_0 $. 

We first note 
$$
E_j E_i = q^3 E_i E_j + E^-_{ij}(1)
$$
and 
$$
E^-_{ij}(1) E_i = q E_i E^-_{ij}(1) + [2]_q E^-_{ij}(2). 
$$
Using 
\begin{align*}
[2]_q E_j E^-_{ij}(2) &= E_j E_{ij}^-(1) E_i - q E_j E_i E_{ij}^-(1) \\
&= q^{-3} E_{ij}^-(1) E_j E_i - q E_{ij}^-(1)^2 - q^4 E_i E_j E^-_{ij}(1) \\
&= q^{-3} E^-_{ij}(1)^2 + E^-_{ij}(1) E_i E_j - q E_{ij}^-(1)^2 - q E_i E^-_{ij}(1) E_j \\
&= (q^{-3} - q) E^-_{ij}(1)^2 + E_{ij}^-(1) E_i E_j - q E_i E_{ij}^-(1) E_j \\
&= (q^{-3} - q) E^-_{ij}(1)^2 + [2]_q E^-_{ij}(2) E_j
\end{align*}
and $ \T_j^{-1} \T_i^{-1}(E^-_{ij}(1)) = \T_j^{-1}(E^+_{ij}(2)) = E^-_{ij}(2) $ we get 
$$
[2]_q \T_j^{-1} \T_i^{-1}(E_j) E_i = 
(q^{-3} - q) E^-_{ij}(2)^2 + [2]_q E_i \T_j^{-1} \T_i^{-1}(E_j). 
$$
In addition, 
\begin{align*}
E^-_{ij}(2) E_i - &q^{-1} E_i E^-_{ij}(2) \\
&= [3]_q E_j E_i^{(3)} - q^2 [2]_q E_i E_j E_i^{(2)} + q^4 E_i^{(2)} E_j E_i \\ 
&\qquad - q^{-1} E_i E_j E_i^{(2)} + q [2]_q E_i^{(2)} E_j E_i - q^3 [3]_q E_i^{(3)} E_j \\
&= [3]_q E^-_{ij}(3), 
\end{align*}
using $ q [3]_q = q^{-1} + q + q^3 = q^2[2]_q + q^{-1} $ and $ q^2[3]_q = 1 + q^2 + q^4 = q^4 + q [2]_q $. 
Finally, we recall $ E^-_{ij}(3) E_i = -q^{-3} E_i E^-_{ij}(3) $. 

Moreover $ E_j E^-_{ij}(1) = q^{-3} E^-_{ij}(1) E_j $ and 
$$
E^-_{ij}(1) E^-_{ij}(2) = q E^-_{ij}(2) E^-_{ij}(1) + [3]_q \T_j^{-1} \T_i^{-1}(E_j).  
$$
We have $ E_i \T_i^{-1}(E_j) = E_i E^+_{ij}(3) = -q^{-3} E^+_{ij}(3) E_i $, and since $ \T_j^{-1}(E_i) = E^-_{ij}(1) $ we get 
$$
E^-_{ij}(1) \T_j^{-1} \T_i^{-1}(E_j) = - q^{-3} \T_j^{-1} \T_i^{-1}(E_j) E^-_{ij}(1). 
$$  
In addition, $ E^+_{ij}(3) E^+_{ij}(2) = q^{-3} E^+_{ij}(2) E^+_{ij}(3) $ implies 
$$
\T_j^{-1} \T_i^{-1}(E_j) E^-_{ij}(2) = q^{-3} E^-_{ij}(2) \T_j^{-1} \T_i^{-1}(E_j). 
$$
Finally, we recall $ E_j E^+_{ij}(1) = q^3 E^+_{ij}(1) E_j $, which gives
$$ 
E^-_{ij}(3) E^-_{ij}(2) = q^3 E^-_{ij}(2) E^-_{ij}(3).
$$ 
Now consider the reduced expression $ w = s_j s_i s_j s_i s_j s_i $. The associated root vectors are
$ 
E_j, \T_j(E_i) = E^+_{ij}(1), \T_j \T_i(E_j) = \T_i^{-1} \T_j^{-1} \T_i^{-1}(E_j), \T_j \T_i \T_j(E_i) = E^+_{ij}(2), 
\T_j \T_i \T_j \T_i(E_j) = \T_i^{-1}(E_j) = E^+_{ij}(3), E_i. 
$
We have to check that left multiplication by $ E_i $ preserves the linear span of all terms of the form 
$$ 
E_j^{r_1} E^+_{ij}(1)^{r_2} \T_i^{-1}\T_j^{-1}\T_i^{-1}(E_j^{r_3}) E^+_{ij}(2)^{r_4} E^+_{ij}(3)^{r_5} E_i^{r_6}
$$
with $ r_j \in \mathbb{N}_0 $.

To this end we observe 
$$
E_i E_j = q^3 E_j E_i + E^+_{ij}(1), 
$$
and 
\begin{align*}
E_i E^+_{ij}(1) &= q E^+_{ij}(1) E_i + [2]_q E^+_{ij}(2), 
\end{align*}
Next observe that $ \T_i^{-1} \T_j^{-1}(E_i) = \T_i^{-1}(E^-_{ij}(1)) = E^+_{ij}(2) $ and $ \T_i \T_j \T_i \T_j \T_i(E_j) = E_j $.  
We obtain 
\begin{align*}
[2]_q E^+_{ij}(2) E_j &= E_i E_{ij}^+(1) E_j - q E_{ij}^+(1) E_i E_j \\
&= q^{-3} E_i E_j E_{ij}^+(1) - q E_{ij}^+(1)^2 + q^4 E^+_{ij}(1)E_j E_i \\
&= q^{-3} E^+_{ij}(1)^2 - E_j E_i E^+_{ij}(1) - q E_{ij}^+(1)^2 + q E_j E^+_{ij}(1) E_i \\
&= (q^{-3} - q) E^+_{ij}(1)^2 - E_j E_i E_{ij}^+(1) + q E_j E_{ij}^+(1) E_i \\
&= (q^{-3} - q) E^+_{ij}(1)^2 - [2]_q E_j E^+_{ij}(2) 
\end{align*}
Using $ \T_j \T_i(E^+_{ij}(1)) = E^+_{ij}(2) $ we get 
\begin{align*}
[2]_q E_i \T_i^{-1} \T_j^{-1} \T_i^{-1}(E_j) &= (q + q^{-3}) E^+_{ij}(2)^2 - [2]_q \T_i^{-1} \T_j^{-1} \T_i^{-1}(E_j) E_i.  
\end{align*}
In addition, 
\begin{align*}
E_i E^+_{ij}(2) - &q^{-1} E^+_{ij}(2) E_i \\
&= [3]_q E_i^{(3)} E_j - q^2 [2]_q E_i^{(2)} E_j E_i + q^4 E_i E_j E_i^{(2)} \\ 
&\qquad - q^{-1} E_i^{(2)} E_j E_i + q [2]_q E_i E_j E_i^{(2)} - q^3 [3]_q E_j E_i^{(3)} \\
&= [3]_q E^+_{ij}(3), 
\end{align*}
using $ q [3]_q = q^{-1} + q + q^3 = q^2[2]_q + q^{-1} $ and $ q^2[3]_q = 1 + q^2 + q^4 = q^4 + q [2]_q $. 
Using the quantum Serre relations we obtain 
\begin{align*}
[3]_q! &E_i E^+_{ij}(3) \\
&= E_i^4 E_j -q [3]_q E_i^3 E_j E_i + q^2 [3]_q E_i^2 E_j E_i^2 - q^3 E_i E_j E_i^3 \\
&= -q^{-3} E_i^3 E_j E_i + q^{-2} [3]_q E_i^2 E_j E_i^2 - q^{-1} [3]_q E_i E_j E_i^3 + E_j E_i^4 \\
&= -q^{-3} [3]_q! E^+_{ij}(3) E_i,  
\end{align*}
taking into account $ [4]_q = q^{-3} + q[3]_q = q^3 + q^{-1} [3]_q $ and $ [4]_q = (q^2 + q^{-2})[2]_q $. \\
Moreover, 
\begin{align*}
E^+_{ij}(1) E_j &= E_i E_j^2 - q^3 E_j E_i E_j \\
&= q^{-3} E_j E_i E_j - E_j^2 E_i \\ 
&= q^{-3} E_j E^+_{ij}(1) 
\end{align*}
using the quantum Serre relation
$$
E_i E_j^2 - (q^3 + q^{-3}) E_j E_i E_j + E_j^2 E_i = 0. 
$$
We have 
$$ 
E^+_{ij}(2) E^+_{ij}(1) = [3]_q \T_i^{-1} \T_j^{-1} \T_i^{-1}(E_j) + q^{-1} E^+_{ij}(1) E^+_{ij}(2), 
$$ 
and from $ E^+_{ij}(1) E_j = q^{-3} E_j E^+_{ij}(1) $ we get 
\begin{align*}
E^+_{ij}(2) \T_i^{-1} &\T_j^{-1} \T_i^{-1}(E_j) = \T_j \T_i(E^+_{ij}(1) E_j) \\
&= q^{-3} \T_j \T_i(E_j E^+_{ij}(1)) = q^{-3} \T_i^{-1} \T_j^{-1} \T_i^{-1}(E_j) E^+_{ij}(2).  
\end{align*}
Using again the quantum Serre relations we compute 
\begin{align*}
[3]_q! & E^-_{ij}(3) E_i \\
&= E_j E_i^4 -q [3]_q E_i E_j E_i^3 + q^2 [3]_q E_i^2 E_j E_i^2 - q^3 E_i^3 E_j E_i \\
&= -q^{-3} E_i E_j E_i^3 + q^{-2} [3]_q E_i^2 E_j E_i^2 - q^{-1} [3]_q E_i^3 E_j E_i + E_i^4 E_j \\
&= -q^{-3} [3]_q! E_i E^-_{ij}(3),  
\end{align*}
taking into account $ [4]_q = q^{-3} + q[3]_q = q^3 + q^{-1} [3]_q $ and $ [4]_q = (q^2 + q^{-2})[2]_q $. 
This yields 
\begin{align*}
\T_i^{-1} \T_j^{-1} &\T_i^{-1}(E_j) E^+_{ij}(1) = \T_j(E^-_{ij}(3) E_i) \\
&= -q^{-3} \T_j(E_i E^-_{ij}(3)) = -q^{-3} E^+_{ij}(1) \T_i^{-1} \T_j^{-1} \T_i^{-1}(E_j).  
\end{align*}
We have 
\begin{align*}
E_j E^-_{ij}(1) &= E_j^2 E_i - q^3 E_j E_i E_j \\
&= q^{-3} E_j E_i E_j - E_i E_j^2 \\ 
&= q^{-3} E^-_{ij}(1) E_j,  
\end{align*}
whence 
\begin{align*}
E^+_{ij}(3) E^+_{ij}(2) &= \T_i^{-1}(E_j E^-_{ij}(1)) \\
&= q^{-3} \T_i^{-1}(E^-_{ij}(1) E_j) = q^{-3} E^+_{ij}(2) E^+_{ij}(3). 
\end{align*}
Combining these relations yields the assertion. 
\end{enumerate}
This finishes the proof. \end{proof} 

\begin{prop} \label{PBWspanindependent}
Let $ w \in W $ and $ w = s_{i_1} \cdots s_{i_t} $ be a reduced expression of $ w $. The linear subspace of $ U_q(\mathfrak{n}_+) $ spanned by the associated 
PBW-vectors $ E^{a_1}_{\beta_1} \cdots E^{a_t}_{\beta_t} $ depends on $ w $, but not on the choice of the reduced expression for $ w $. 
\end{prop} 

\begin{proof} We use induction on $ l(w) $, the cases $ l(w) = 0 $ and $ l(w) = 1 $ being trivial. 
Now assume $ l(w) > 1 $ and let $ w = s_{i_1} \cdots s_{i_r} $ 
be a reduced expression. Any other reduced expression $ s_{j_1} \cdots s_{j_r} $ is obtained from a sequence of elementary moves applied 
to $ s_{i_1} \cdots s_{i_r} $, so it suffices to show that the subspace spanned by the associated PBW-vectors does not change if 
$ s_{j_1} \cdots s_{j_r} $ is obtained by applying a single elementary move to $ s_{i_1} \cdots s_{i_r} $. Set $ \alpha = \alpha_{i_1} $ 
and $ \beta = \alpha_{j_1} $. 

Assume first $ \alpha = \beta $ write $ u = s_{i_1} w = s_{j_1} w $. The the inductive hypothesis applied to $ u $ shows that the subspace $ U $ 
associated to $ u $ is the same for the reduced expressions $ s_{i_2} \cdots s_{i_r} $ and $ s_{j_2} \cdots s_{j_r} $ of $ u $. It follows 
from the definitions that the space spanned by the PBW-vectors associated to both $ s_{i_1} \cdots s_{i_r} $ and $ s_{j_1} \cdots s_{j_r} $ is the sum of the 
subspaces $ E^k_{i_1} \T_{\alpha_i}(U) $ for $ k \in \mathbb{N}_0 $. In particular, these spaces agree for both reduced expressions of $ w $. 

Assume now $ \alpha \neq \beta $. Then the elementary move changes the first letter of the expression $ s_{i_1} \cdots s_{i_r} $, and therefore 
the corresponding elementary move has to take place at the beginning of the expression. Let $ u $ be the longest word in the subgroup of $ W $ 
generated by $ s_{i_1} $ and $ s_{j_1} $. Then $ s_{i_1} \cdots s_{i_r} $ and $ s_{j_1} \cdots s_{j_r} $ start with the two distinct reduced 
expressions for $ u $ and agree after that. Let us write $ w = uv $ with $ l(w) = l(u) + l(v) $. According to Lemma \ref{PBWexpressionhelplemma} 
the subspace $ U $ spanned by the PBW-vectors associated to $ u $ is independent of the two expressions. Let us write $ V $ for 
the span of PBW-vectors associated to the reduced expression of $ v $ induced by $ s_{i_1} \cdots s_{i_r} $, or equivalently $ s_{j_1} \cdots s_{j_r} $. 
Then the span of the PBW-vectors associated to the given expressions for $ w $ is equal to the $ U \T_u(V) $, 
in particular, it does not depend on the reduced expressions. \end{proof} 

In the sequel we will write $ U_q(\mathfrak{n}_+)[w] \subset U_q(\mathfrak{n}_+) $
\nomenclature{$U_q(\mathfrak{n}_+)[w]$}{subspace of $U_q(\mathfrak{n}_+)$ spanned by PBW-vectors associated to $w\in W$}%
for the subspace spanned by the PBW-vectors associated with any reduced expression of $ w $. 

Let us fix a reduced expression $ w_0 = s_{i_1} \cdots s_{i_n} $ for the longest element $ w_0 $ of $ W $.
\nomenclature{$w_0$}{longest element of the Weyl group}%
Then the positive 
roots of $ \mathfrak{g} $ can be uniquely written in the form $ \beta_r = s_{i_1} \cdots s_{i_{r - 1}} \alpha_{i_r} $ for $ 1 \leq r \leq n $.

\begin{definition}
\label{def:quantum_root_vectors}
Fix a reduced expression $w_0 = s_{i_1} \cdots s_{i_n} $ for the longest element of $W$.
We define the \emph{quantum root vectors} of $ U_q(\mathfrak{g}) $ to be the associated PBW-vectors
$$
E_{\beta_r} = \T_{i_1} \cdots \T_{i_{r - 1}}(E_{i_r}), \quad F_{\beta_r} = \T_{i_1} \cdots \T_{i_{r - 1}}(F_{i_r}).
$$ 
\end{definition}

Let us now formulate and prove the Poincar\'e-Birkhoff-Witt Theorem for $ U_q(\mathfrak{g}) $. 

\begin{theorem}[PBW-basis - Non-root of unity case] \label{PBWgeneralfield}
Assume $ q \in \mathbb{K}^\times $ is not a root of unity. Then the elements 
$$
F^{b_1}_{\beta_1} \cdots F^{b_n}_{\beta_n} K_\lambda E^{a_1}_{\beta_1} \cdots E^{a_n}_{\beta_n}
$$
where $ b_j, a_j \in \mathbb{N}_0 $ for all $ j $,  $ \lambda \in \weights $ and $ \beta_1, \dots, \beta_n $ are the 
positive roots of $ \mathfrak{g} $, form a vector space basis of $ U_q(\mathfrak{g}) $ over $ \mathbb{K} $. 
\end{theorem} 

\begin{proof}
According to the triangular decomposition of $ U_q(\mathfrak{g}) $, see Proposition \ref{uqgtriangular}, it suffices to show that 
the vectors $ E^{a_1}_{\beta_1} \cdots E^{a_n}_{\beta_n} $ with $ a_j \in \mathbb{N}_0 $ form a basis of $ U_q(\mathfrak{n}_+) $. 
The corresponding claim for $ U_q(\mathfrak{n}_-) $ can then be obtained by applying the automorphism $ \omega $. Indeed, one checks that  $\T_i(\omega(E_j))$ and $\omega(\T_i(E_j))$ agree up to invertible elements in $\mathbb{K}$, for every $1\leq i,j \leq N$.

Firstly, according to Proposition \ref{PBWindependence} the PBW-vectors $ E^{a_1}_{\beta_1} \cdots E^{a_n}_{\beta_n} $ 
with $ a_j \in \mathbb{N}_0 $ are linearly independent. 
To show that these vectors span $ U_q(\mathfrak{n}_+) $ we need an auxiliary consideration. 
If $ w \in W $ and $ \alpha_i \in \Sigma $ satisfy $ w^{-1} \alpha_i \in {\bf \Delta}^- $ then we 
have $ E_i U_q(\mathfrak{n}_+)[w] \subset U_q(\mathfrak{n}_+)[w] $. Indeed, in this case we find a reduced decomposition 
$ w = s_{i_1} \cdots s_{i_t} $ with $ i_1 = i $, hence the PBW-vectors start with powers of $ E_i $. Clearly, the resulting 
space is invariant under left multiplication by $ E_i $. 

Let us now show that the PBW-vectors span $ U_q(\mathfrak{n}_+) $. By definition, these vectors span the subspace $ U_q(\mathfrak{n}_+)[w_0] $. 
Since $ w_0 $ is the longest word of $ W $ we have $ w_0^{-1}(\alpha_i) \in {\bf \Delta}^- $ and thus 
$ E_i U_q(\mathfrak{n}_+)[w_0] \subset U_q(\mathfrak{n}_+)[w_0] $ for all $ i $. Hence $ U_q(\mathfrak{n}_+)[w_0] $ is closed 
under left multiplication by the generators $ E_1, \dots, E_N $. By construction, the space $ U_q(\mathfrak{n}_+)[w_0] $ 
also contains $ 1 $. It follows that $ U_q(\mathfrak{n}_+)[w_0] = U_q(\mathfrak{n}_+) $ as desired. \end{proof}

We point out that one can reverse the ordering in the basis obtained in Theorem \ref{PBWgeneralfield}. 
That is, the vectors 
$$
F^{b_n}_{\beta_n} \cdots F^{b_1}_{\beta_1} K_\lambda E^{a_n}_{\beta_n} \cdots E^{a_1}_{\beta_1}
$$
where $ b_j, a_j \in \mathbb{N}_0 $ for all $ j $, $ \lambda \in \weights $ 
and $ \beta_1, \dots, \beta_n $ are the positive roots of $ \mathfrak{g} $,
form again a vector space basis of $ U_q(\mathfrak{g}) $. This can be seen using the anti-automorphism $\delta$
of $ U_q(\mathfrak{g}) $
which preserves all generators but exchanges $ q $ with $ q^{-1} $. 
One checks that the automorphisms $ \T_i^{-1} \delta \T_i \delta $ act by certain scalars on the generators, and hence 
diagonally on weight spaces of $ U_q(\mathfrak{g}) $. Hence we have for all $ 1 \leq j \leq n $ that 
$ \delta(E_{\beta_j}) $ is equal to $ E_{\beta_j} $ up to a scalar,
and since $ \delta $ is anti-multiplicative it follows from Theorem \ref{PBWgeneralfield} that the vectors $ E^{a_n}_{\beta_n} \cdots E^{a_1}_{\beta_1} $ 
form indeed a basis. We will frequently make use of this fact in the sequel, and refer to the vectors $ E^{a_n}_{\beta_n} \cdots E^{a_1}_{\beta_1} $ 
again as PBW-vectors. In a similar way we proceed for $ F^{b_n}_{\beta_n} \cdots F^{b_1}_{\beta_1} $ of course. 

Let us also remark that the rank $ 1 $ case of Theorem \ref{PBWgeneralfield} follows immediately from Proposition \ref{uqgtriangular}. 

As a consequence of the PBW-basis Theorem, we see that the dimension of the weight spaces $ U_q(\mathfrak{n}_\pm)_\nu $ for $ \nu \in \roots^\pm $ 
coincide with their classical counterparts. Here we are writing $ \roots^- = - \roots^+ $.  
\nomenclature[o$Q$7]{$\roots^-$}{$=-\roots^+$}%
These dimensions can be described in terms of the Kostant partition function 
$P:\roots \to \NN_0$, defined by
\[
 P(\nu) = \left| \{ (r_1,\ldots,r_n)\in\NN_0^n \mid r_1\beta_1 + \cdots + r_n\beta_n = \nu \}\right|.
\]
\label{nom:Kostant_partition_function1}%

\begin{prop} 
	\label{prop:Uqn_weight_multiplicities}
For all $\nu\in\roots^+$ we have $ \dim U_q(\mathfrak{n}_\pm)_{\pm\nu} = P(\nu)$, where $P$ denotes Kostant's partition function.
\end{prop} 

\begin{proof} Let us only consider the case of $ U_q(\mathfrak{n}_+) $. 
According to Theorem \ref{PBWgeneralfield}, a linear basis for $ U_q(\mathfrak{n}_+)_{\nu} $ for $ \nu \in \roots^+ $ is given by 
all PBW-vectors $ E^{a_1}_{\beta_1} \cdots E^{a_n}_{\beta_n} $ such that $ a_1 \beta_1 + \cdots a_n \beta_N = \nu $. 
The number of these vectors is precisely $ P(\nu) $. 
\end{proof} 

Let us finish this subsection by the following result, originally due to Levendorskii,
which will help us to analyze 
the structure of $ U_q(\mathfrak{n}_\pm) $ further below. 

\begin{prop} \label{uqnplusqcommutationrelations}
If $ 0 \leq r < s \leq n $ then we have 
$$
E_{\beta_r} E_{\beta_s} - q^{(\beta_r, \beta_s)} E_{\beta_s} E_{\beta_r} = 
\sum_{t_{r + 1}, \dots t_{s - 1} = 0}^\infty c_{t_{r + 1}, \dots, t_{s - 1}} E_{\beta_{s - 1}}^{t_{s - 1}} \cdots E_{\beta_{r + 1}}^{t_{r + 1}}, 
$$
with only finitely many coefficients $ c_{t_{r + 1}, \dots, t_{s - 1}} \in \mathbb{K} $ in the sum on the right hand side being nonzero. 
 Likewise,
$$
F_{\beta_r} F_{\beta_s} - q^{(\beta_r, \beta_s)} F_{\beta_s} F_{\beta_r} = 
\sum_{t_{r + 1}, \dots t_{s - 1} = 0}^\infty d_{t_{r + 1}, \dots, t_{s - 1}} F_{\beta_{s - 1}}^{t_{s - 1}} \cdots F_{\beta_{r + 1}}^{t_{r + 1}}, 
$$
with only finitely many $d_{t_{r + 1}, \dots, t_{s - 1}} \in \mathbb{K}$ nonzero.
\end{prop} 

\begin{proof} According to the PBW-Theorem \ref{PBWgeneralfield} we can write 
$$
E_{\beta_r} E_{\beta_s} = 
\sum_{t_1, \dots t_n \in \mathbb{N}_0} c_{t_1, \dots, t_n} E_{\beta_n}^{t_n} \cdots E_{\beta_1}^{t_1}, 
$$
with only finitely many nonzero terms on the right hand side. Assume that there exists $ a < r $ such that some $ c_{0, \dots, 0, t_a, \dots t_n} \neq 0 $ with $t_a>0$, and 
without loss of generality pick the smallest such $ a $. If we apply $ \T_{i_a}^{-1} \cdots \T_{i_1}^{-1} $ to both sides of the above equation then
the left hand side, and all terms on the right hand side with $ t_a = 0 $ are contained in $ U_q(\mathfrak{n}_+) $. 
Indeed, all these terms are products of root vectors of the form $ \T_{i_{a + 1}} \cdots \T_{i_{l - 1}}(E_{i_l}) $ such that $ a + 1 \leq l $, 
recall that $ E_{\beta_l} = \T_{i_1} \cdots \T_{i_{l - 1}}(E_{i_l}) $. The word $ s_{i_{a + 1}} \cdots s_{i_l} $ is a reduced expression 
and hence $ s_{i_{a + 1}} \cdots s_{i_{l - 1}} \alpha_{i_l} $ is a positive root. 
According to Lemma \ref{PBWpositivelemma} this implies $ E_{\beta_l} = \T_{i_1} \cdots \T_{i_{l - 1}}(E_{i_l}) \in U_q(\mathfrak{n}_+) $.

The terms with $ t_a > 0 $ are of the form $ X \T_{i_a}^{-1}(E_{i_a})^{t_a} $ for $ X \in U_q(\mathfrak{n}_+) $. 
Since $ \T_{i_a}^{-1}(E_{i_a}) = - F_{i_a} K_{i_a}^{-1} $ this is a contradiction to the PBW-Theorem \ref{PBWgeneralfield}.

In a similar way one checks that $ c_{t_1, \dots, t_b, 0, \dots, 0} \neq 0 $ for $ b > s $ is impossible. Indeed, assume that 
a nonzero coefficient of this form exists and pick $ b $ maximal with this property. 
Applying $ \T_{i_{b - 1}}^{-1} \cdots \T_{i_1}^{-1} $ to both sides of the equation leads to expressions in $ U_q(\mathfrak{n}_-) $, 
except for the terms with $ t_b \neq 0 $. 

Next, consider the nonzero coefficients $ c_{0, \dots, 0, t_r, \dots, t_s, 0, \dots} $ with $ t_r $ or $ t_s $ nonzero.  In a similar way as above we apply $ \T_{i_r}^{-1} \cdots \T_{i_1}^{-1} $ and $ \T_{i_{s - 1}}^{-1} \cdots \T_{i_1}^{-1} $, respectively, 
and compare both sides of the equation. The highest powers of $ F_{i_r} $ or $ E_{i_s} $ appearing in the resulting expressions must have the same degree 
on both sides, and this forces $ t_r = 1 $ and $ t_s = 1 $, respectively. 
As a consequence, for any nonzero coefficient of this form we obtain $ t_j = 0 $ for all $ r < j < s $ by comparing weights. 

It follows that we can write 
$$
E_{\beta_r} E_{\beta_s} - c E_{\beta_s} E_{\beta_r} =
\sum_{t_{r + 1}, \dots t_{s - 1} = 0}^\infty c_{t_{r + 1}, \dots, t_{s - 1}} E_{\beta_{s - 1}}^{t_{s - 1}} \cdots E_{\beta_{r + 1}}^{t_{r + 1}}, 
$$
with an as yet unspecified constant $ c $. In order to compute $ c $ we apply $ \T_{i_r}^{-1} \cdots \T_{i_1}^{-1} $ to both sides of the previous equation. 
Then the left hand side becomes 
$$
- F_{i_r} K_{i_r}^{-1} E + c E F_{i_r} K_{i_r}^{-1}
$$
where $ E = \T_{i_{r + 1}} \cdots \T_{i_{s - 1}}(E_{i_s}) $. Since $ E $ has weight $ s_{i_r} \cdots s_{i_1}(\beta_s) $ we get 
$ K_{i_r}^{-1} E = q^{-(\alpha_{i_r}, s_{i_r} \cdots s_{i_1}(\beta_s))} E K_{i_r}^{-1} $. Moreover 
$$
-(\alpha_{i_r}, s_{i_r} \cdots s_{i_1}(\beta_s)) = (s_{i_1} \cdots s_{i_{r - 1}}(\alpha_{i_r}), \beta_s) = (\beta_r, \beta_s), 
$$
using that $ s_{i_r}(\alpha_{i_r}) = -\alpha_{i_r} $. Accordingly, we obtain 
$$
- F_{i_r} K_{i_r}^{-1} E + c E F_{i_r} K_{i_r}^{-1} = (- q^{(\beta_r, \beta_s)} F_{i_r} E + c E F_{i_r}) K_{i_r}^{-1}. 
$$
Writing $ E $ as a sum of monomials in the generators $ E_1, \dots, E_N $ we get $ F_{i_r} E - E F_{i_r} \in U_q(\mathfrak{b}_+) $. Comparing 
with the right hand side of the formula and using once more the PBW-Theorem \ref{PBWgeneralfield} yields $ c = q^{(\beta_r, \beta_s)} $ as desired.

The claim for the PBW-vectors of $ U_q(\mathfrak{n}_-) $
follows by applying the automorphism $\omega$ to the formula in $ U_q(\mathfrak{n}_+) $ and using the fact that  $\T_i(\omega(E_j))$ and $\omega(\T_i(E_j))$ agree up to an invertible element in $\mathbb{K}$,  for every $i,j$.
\end{proof}

\subsubsection{The integral PBW-basis}

Recall that the integral version $ U_q^\A(\mathfrak{g}) $ of the quantized universal enveloping algebra is defined over 
$ \A = \mathbb{Z}[s, s^{-1}] $ where $ q = s^L $. Let us fix a reduced expression $ w_0 = w_{i_1} \cdots w_{i_n} $ for the longest 
element $ w_0 $ of $ W $. We define the \emph{restricted PBW-elements} for $ U_q(\mathfrak{n}_+) $ by 
$$
E^{(a_1)}_{\beta_1} \cdots E^{(a_n)}_{\beta_n},
$$
where 
$$ 
E_{\beta_r}^{(a_r)} = \frac{1}{[a_r]_{q_{\beta_r}}!} E_\beta^{a_r} 
$$
and $ q_{\beta_r} = q_{i_r} $ if $ E_{\beta_r} = \T_{i_1} \cdots \T_{i_{r - 1}}(E_{i_r}) $.
\label{nom:q_beta_r}
That is, the difference to the ordinary PBW-elements is that we multiply by the inverse of $ [a_1]_{q_{\beta_1}}! \cdots [a_n]_{q_{\beta_n}}! $. 
In a similar way one obtains restricted PBW-elements for $ U_q(\mathfrak{n}_-) $. 

Before we can state the integral version of the PBW-Theorem we need a refinement of Proposition \ref{Tigeneratorscommutation}. 

\begin{lemma} \label{Tigeneratorscommutationhigherpowers}
Let $ V $ be an integrable $ U_q(\mathfrak{g}) $-module. For any $ n \in \mathbb{N} $ we have 
\begin{align*}
\T_i(E_j^{(n)} \cdot v) &= (-q_i)^{-na_{ij}} (E_j^{(n)} \leftarrow \hat{S}^{-1}(E^{(-na_{ij})}_i)) \cdot \T_i(v) \\ 
\T_i(F_j^{(n)} \cdot v) &= (F^{(-na_{ij})}_i \rightarrow F_j^{(n)}) \cdot \T_i(v) 
\end{align*}
for all $ v \in V $. 
\end{lemma}

\begin{proof} Let us first consider the second formula. According to Proposition \ref{Tigeneratorscommutation} we have the claim for $ n = 1 $. 
In order to prove the assertion for general $ n $ it suffices to 
show $ (F^{(-a_{ij})}_i \rightarrow F_j)^n = F^{(-na_{ij})}_i \rightarrow (F_j^n) $ by induction. 
Assuming that the claim holds for $ n - 1 $, 
we compute, using Lemma \ref{lemEiFicomult},
\begin{align*}
F^{(-na_{ij})}_i \rightarrow (F_j^n)
 &
 =
 \left( (F^{(-na_{ij})}_i )_{(1)} \rightarrow F_j^{n-1} \right)
 \left( (F^{(-na_{ij})}_i )_{(2)} \rightarrow F_j \right)
\\
 &= \sum_{k = 0}^{-n a_{ij}} q_i^{k(-na_{ij} - k)} (F_i^{(-n a_{ij} - k)} K_i^{-k} \rightarrow F_j^{n - 1}) 
(F_i^{(k)} \rightarrow F_j) .
\end{align*}
Using that $ F_i^{(-a_{ij} + 1)} \rightarrow F_j = 0 $ and the inductive hypothesis, all of the terms in this sum vanish except when $k=-a_{ij}$, leaving
\begin{align*}
F^{(-na_{ij})}_i \rightarrow (F_j^n)
&= q_i^{-a_{ij}(-(n - 1) a_{ij})} (F_i^{(-(n - 1)a_{ij})} K_i^{a_{ij}} \rightarrow F_j^{n - 1}) 
(F_j^{(-a_{ij})} \rightarrow F_j) \\
&= (F_i^{(-(n - 1)a_{ij})} \rightarrow F_j^{n - 1}) (F_j^{(-a_{ij})} \rightarrow F_j) \\ 
&= (F_j^{(-a_{ij})} \rightarrow F_j)^{n - 1} (F_j^{(-a_{ij})} \rightarrow F_j) = (F_j^{(-a_{ij})} \rightarrow F_j)^n 
\end{align*}
as desired. 

Using this result, one obtains the first formula in the same way as in the proof of Proposition \ref{Tigeneratorscommutation}. \end{proof} 

Let us now present an integral version of the PBW Theorem \ref{PBWgeneralfield}. Recall that $ U_q^\A(\mathfrak{n}_\pm) \subset U_q^\A(\mathfrak{g}) $ 
is the subalgebra generated by all divided powers $ E_i^{(a_i)} $ and $ F_j^{(b_j)} $, respectively. 

\begin{theorem}[PBW-basis - Integral Case] \label{PBWintegral}
The elements 
$$
E^{(a_1)}_{\beta_1} \cdots E^{(a_n)}_{\beta_n}
$$
where $ a_j \in \mathbb{N}_0 $ and $ \beta_1, \dots, \beta_n $ are the 
positive roots of $ \mathfrak{g} $, form a free basis of $ U_q^\A(\mathfrak{n}_+) $. Similarly, the elements 
$$
F^{(b_1)}_{\beta_1} \cdots F^{(b_n)}_{\beta_n}
$$
where $ b_j \in \mathbb{N}_0 $ and $ \beta_1, \dots, \beta_n $ are the 
positive roots of $ \mathfrak{g} $, form a free basis of $ U_q^\A(\mathfrak{n}_-) $. 
\end{theorem} 

\begin{proof} Recall that we can view $ U_q^\A(\mathfrak{g}) $ as a subring of the quantized universal enveloping algebra $ U_q(\mathfrak{g}) $ 
defined over $ \mathbb{K} = \mathbb{Q}(s) $. 
Using the formulas in Theorem \ref{deftiautomorphism} and Lemma \ref{Tigeneratorscommutationhigherpowers}  
we see that the operators $ \T_i $ map $ U_q^\A(\mathfrak{g}) $ into itself. 
It follows in particular that the PBW-vectors are indeed contained in $ U_q^\A(\mathfrak{g}) $. 

In the same way as in the proof of Theorem \ref{PBWgeneralfield} one checks that left multiplication with 
$ E_i^{(r)} $ preserves the $ \A $-span of the PBW-basis vectors $ E^{(a_1)}_{\beta_1} \cdots E^{(a_n)}_{\beta_n} $ 
for all $ 1 \leq i \leq N $ and $ r \in \mathbb{N}_0 $. 
We point out in particular that the commutation relations obtained in Lemma \ref{PBWexpressionhelplemma} only involve coefficients 
from $ \A $, so that the integral version of Proposition \ref{PBWspanindependent} holds true. That is, the $ \A $-linear span 
of the restricted PBW-vectors $ E^{(a_1)}_{\beta_1} \cdots E^{(a_n)}_{\beta_n} $ inside $ U_q(\mathfrak{n}_+) $ 
is independent of the reduced expression of $ w_0 $. 
Linear independence of the restricted PBW-elements over $ \A $ follows from linear independence over $ \mathbb{Q}(s) $. 

The assertion for $ U_q(\mathfrak{n}_-) $ is obtained in an analogous way. \end{proof} 

Using Theorem \ref{PBWintegral} we see that the triangular decomposition of $ U_q(\mathfrak{g}) $ carries over to the integral form as well. 
Recall that $ U_q^\A(\mathfrak{h}) $ denotes 
the intersection of $U_q^\A(\mathfrak{g})$ with $U_q(\mathfrak{h})$.

\begin{prop} \label{uqgtriangularintegral}
Multiplication in $ U_q^\A(\mathfrak{g}) $ induces an isomorphism 
$$
U_q^\A(\mathfrak{n}_-) \otimes_\A U_q^\A(\mathfrak{h}) \otimes_\A U_q^\A(\mathfrak{n}_+) \cong U_q^\A(\mathfrak{g})
$$
of $ \A $-modules, and we have 
$$
\mathbb{Q}(s) \otimes_\A U_q^\A(\mathfrak{g}) \cong U_q(\mathfrak{g}), 
$$
where $ U_q(\mathfrak{g}) $ is the quantized universal enveloping algebra over $ \mathbb{Q}(s) $. 
\end{prop} 

\begin{proof} It is easy to see that the image in $ U_q^\A(\mathfrak{g}) $ under the multiplication map is a subalgebra containing all 
generators, and therefore yields all of $ U_q^\A(\mathfrak{g}) $. According to Theorem \ref{PBWintegral} and Theorem \ref{PBWgeneralfield} 
we have $ \mathbb{Q}(s) \otimes_\A U_q^\A(\mathfrak{n}_\pm) \cong U_q(\mathfrak{n}_\pm) $. 
From Proposition \ref{PBWcartanintegral} we know $ \mathbb{Q}(s) \otimes_\A U_q^\A(\mathfrak{h}) \cong U_q(\mathfrak{h}) $. 
Hence the multiplication map $ U_q^A(\mathfrak{n}_-) \otimes_\A U_q^\A(\mathfrak{h}) \otimes_\A U_q^\A(\mathfrak{n}_+) \rightarrow U_q^\A(\mathfrak{g}) $ is 
injective. The remaining claim then follows taking into account Proposition \ref{uqgtriangular}. \end{proof} 

Proposition \ref{uqgtriangularintegral} shows in particular that $ U_q^\A(\mathfrak{g}) $ is indeed an integral version of $ U_q(\mathfrak{g}) $ 
in a natural sense.

\subsection{The Drinfeld pairing and the quantum Killing form} \label{seckilling}

In this section we discuss the Drinfeld pairing and the quantum Killing form, compare \cite{Tanisakikilling}. 

\subsubsection{Drinfeld pairing }

Let us start by defining the Drinfeld pairing.   
We assume that $ \mathbb{K} $ is a field and that $ q = s^L \in \mathbb{K}^\times $ is not a root of unity.

We recall that a skew-pairing of Hopf algebras $L$ and $K$ means a bilinear map $\rho : L\times K \to \mathbb{K}$ such that
\begin{align*}
 \rho(xy,f) &= \rho(x,f_{(1)}) \rho(y,f_{(2)}),
 & \rho(x,fg) &= \rho(x_{(1)},g) \rho(x_{(2)},f), \\
 \rho(1,f) &= \epsilon_K(f), 
 & \rho(x,1) &= \epsilon_L(x).
\end{align*}

\begin{prop} 
There exists a unique skew-pairing $ \rho: U_q(\mathfrak{b}_-) \otimes U_q(\mathfrak{b}_+) \rightarrow \mathbb{K} $ determined by 
$$
\rho(K_\alpha, K_\beta) = q^{(\alpha, \beta)}, \qquad \rho(K_\alpha, E_i) = 0 = \rho(F_j, K_\beta), \qquad 
\rho(F_j, E_i) = \frac{\delta_{ij}}{q_i - q_i^{-1}}
$$
\nomenclature{$\rho$}{variant of the Drinfeld pairing, $\rho:U_q(\lie{b}_-)\otimes U_q(\lie{b}_+)\to\mathbb{K}$}%
for all $ \alpha, \beta \in \weights $ and $ i,j = 1, \dots, N $. 
\end{prop} 

\begin{proof} Uniqueness of the form is clear since the elements $ K_\alpha, E_i $ and $ K_\beta, F_j $ generate 
$ U_q(\mathfrak{b}_+) $ and $ U_q(\mathfrak{b}_-) $, respectively.

We define linear functionals $ \kappa_\lambda, \phi_i \in U_q(\mathfrak{b}_+)^* $ for $ \lambda \in \weights $ and $ i = 1, \dots, N $ by
$$
\kappa_\lambda(K_\beta X) = q^{(\lambda, \beta)} \hat{\epsilon}(X) 
$$
for $ X \in U_q(\mathfrak{n}_+) $ and 
$$
\phi_i(K_\beta E_i) = \frac{1}{q_i - q_i^{-1}}, 
\qquad \phi_i(K_\beta X) = 0 \quad \text{if } X \in U_q(\mathfrak{n}_+)_\gamma, \gamma \neq \alpha_i. 
$$

Consider the convolution product on $U_q(\mathfrak{b}_+)^*$ dual to comultiplication on  $U_q(\mathfrak{b}_+)$, namely
\[
 \varphi\varphi' = (\varphi\otimes\varphi')\circ\hat\Delta, \qquad\qquad \text{for }\varphi, \varphi' \in U_q(\mathfrak{b}_+)^*.
\]
With this product, $\kappa_\lambda$ is convolution invertible with inverse 
$ \kappa_\lambda^{-1} = \kappa_{-\lambda} $.  Moreover, we claim that
\begin{align*}
\kappa_\lambda \phi_j &= q^{-(\alpha_j, \lambda)} \phi_j \kappa_\lambda.
\end{align*}
Indeed, for $ X \in U_q(\mathfrak{n}_+) $ we have $ (\kappa_\lambda \phi_j \kappa_\lambda^{-1})(K_\mu X) = 0 = \phi_j(K_\mu X) $ 
if $ X \in U_q(\mathfrak{n}_+)_\beta $ for $ \beta \neq \alpha_j $ and  
\begin{align*}
(\kappa_\lambda \phi_j)(K_\mu E_j) &= \kappa_\lambda(K_\mu) \phi_j(K_\mu E_j ) + \kappa_\lambda(K_\mu E_j) \phi_j(K_\mu K_j) \\
&= q^{(\lambda, \mu)} \phi_j(K_\mu E_j) \\ 
&= q^{-(\alpha_j, \lambda)}(\phi_j(K_\mu) \kappa_\lambda(K_\mu E_j) + \phi_j(K_\mu E_j) \kappa_\lambda(K_\mu K_j)) \\
&= q^{-(\alpha_j, \lambda)}(\phi_j \kappa_\lambda)(K_\mu E_j).
\end{align*}
It follows that we obtain an algebra homomorphism $ \tilde{\gamma}: \tilde{U}_q(\mathfrak{b}_-) \rightarrow U_q(\mathfrak{b}_+)^* $ by 
setting $ \tilde{\gamma}(K_\lambda) = \kappa_\lambda $ and $ \tilde{\gamma}(F_j) = \phi_j $. Here $ \tilde{U}_q(\mathfrak{b}_-) $ 
denotes the algebra without Serre relations discussed in Section \ref{secdefuqg}.  

We define $ \tilde{\rho}(Y,X) = \tilde{\gamma}(Y)(X) $. 
Then we have $ \tilde{\rho}(YZ, X) = \tilde{\rho}(Y, X_{(1)}) \tilde{\rho}(Z, X_{(2)}) $ by construction. In order to 
show $ \tilde{\rho}(Y, WX) = \tilde{\rho}(Y_{(2)}, W) \tilde{\rho}(Y_{(1)}, X) $ it suffices to consider generators in the 
first variable, since the relation holds for $ \tilde{\rho}(Y_1 Y_2, WX) $ for arbitrary $ W, X $ iff it holds for 
all $ \tilde{\rho}(Y_1, WX) $ and $ \tilde{\rho}(Y_2, WX) $, for arbitrary $ W, X $.

Let $ E_{\bf i} = E_{i_1} \cdots E_{i_k}, E_{\bf j} = E_{j_1} \cdots E_{j_l} $ and $ \mu, \nu \in \weights $ be given, and consider $\tilde{\rho}(K_\lambda, K_\mu E_{\bf i} K_\nu E_{\bf j})$.  If either one of $ E_{\bf i} $ or $ E_{\bf j} $ is different from $ 1 $ we get
\[
 \tilde{\rho}(K_\lambda, K_\mu E_{\bf i} K_\nu E_{\bf j}) = 0  = \tilde{\rho}(K_\lambda, K_\nu E_{\bf j}) 
\tilde{\rho}(K_\lambda, K_\mu E_{\bf i}).
\]
Otherwise we get
\begin{align*} 
\tilde{\rho}(K_\lambda, K_\mu K_\nu) &= \kappa_\lambda(K_\mu K_\nu) = q^{(\lambda, \mu + \nu)} 
= \kappa_\lambda(K_\nu) \kappa_\lambda(K_\mu) = \tilde{\rho}(K_\lambda, K_\nu) \tilde{\rho}(K_\lambda, K_\mu) .
\end{align*}

Moreover we have 
\begin{align*} 
\tilde{\rho}(F_k, K_\mu E_{\bf i} K_\nu E_{\bf j}) &= \phi_k(K_\mu E_{\bf i} K_\nu E_{\bf j}) = 0 
= (\phi_k \otimes 1 + \kappa_k^{-1} \otimes \phi_k)(K_\nu E_{\bf j} \otimes K_\mu E_{\bf i}) 
\end{align*} 
unless precisely one of the terms $ E_{\bf i} $ or $ E_{\bf j} $ has length one. 
In these cases we obtain 
\begin{align*}
\tilde{\rho}(F_k, K_\mu E_i K_\nu) = \phi_k(K_\mu E_i K_\nu) &= q^{-(\alpha_i, \nu)} \phi_k(K_\mu K_\nu E_i) \\
&= q^{-(\alpha_i, \nu)} \delta_{ik} \frac{1}{q_i - q_i^{-1}} \\
&= (\phi_k \otimes 1 + \kappa_k^{-1} \otimes \phi_k)(K_\nu \otimes K_\mu E_i) \\
&= \tilde{\rho}((F_k)_{(1)}, K_\nu) \tilde{\rho}((F_k)_{(2)}, K_\mu E_i) 
\end{align*}
and
\begin{align*}
\tilde{\rho}(F_k, K_\mu K_\nu E_j) &= \phi_k(K_\mu K_\nu E_j) \\
&= \delta_{jk} \frac{1}{q_j - q_j^{-1}} \\
&= (\phi_k \otimes 1 + \kappa_k^{-1} \otimes \phi_k)(K_\nu E_j \otimes K_\mu) \\
&= \tilde{\rho}((F_k)_{(1)}, K_\nu E_j) \tilde{\rho}((F_k)_{(2)}, K_\mu), 
\end{align*}
respectively.

Next we claim that
\begin{align*}
\tilde{\gamma}(u^-_{ij}) = \sum_{k = 0}^{1 - a_{ij}} (-1)^k \begin{bmatrix} 1 - a_{ij} \\ k \end{bmatrix}_{q_i}
&\phi_i^{1 - a_{ij} - k} \phi_j \phi_i^k = 0,  
\end{align*}
where we recall that the Serre elements $ u^-_{ij} $ were introduced in Subsection \ref{sec:Serre_elements}. 
For this, notice that the individual summands vanish on all monomials in the $ E_k $ except on terms of the form 
$ K_\nu E_i^r E_j E_i^{1 - a_{ij} - r} $. According to Proposition \ref{serrecomultiplication} and the relation
$ \tilde{\rho}(Y, WX) = \tilde{\rho}(Y_{(2)}, W) \tilde{\rho}(Y_{(1)}, X) $ for all $ W, X \in U_q(\mathfrak{b}_+) $ established 
above we have 
$$ 
\tilde{\rho}(u^-_{ij}, WX) = \tilde{\rho}(u^-_{ij}, W) \tilde{\rho}(\kappa_i^{-(1 - a_{ij})} \kappa_j^{-1}, X) 
+ \tilde{\rho}(1, W) \tilde{\rho}(u^-_{ij}, X)
$$ 
for all $ W, X \in U_q(\mathfrak{b}_+) $, therefore the claim follows from degree considerations. We conclude that the map $ \tilde{\gamma} $ 
factorises to an algebra homomorphism $ \gamma: U_q(\mathfrak{b}_-) \rightarrow U_q(\mathfrak{b}_+)^* $ such 
that $ \gamma(K_\lambda) = \kappa_\lambda $ and $ \gamma(F_j) = \phi_j $. We define 
$ \rho(Y,X) = \gamma(Y)(X) $. 

By construction, we have $ \rho(YZ, X) = \rho(Y, X_{(1)}) \rho(Z, X_{(2)}) $ for all $ Y,Z \in U_q(\mathfrak{b}_-) $ and $ X \in U_q(\mathfrak{b}_+) $, 
and the relation $ \rho(Y, WX) = \rho(Y_{(2)}, W) \rho(Y_{(1)}, X) $ for $ Y \in U_q(\mathfrak{b}_-) $ and $ W, X \in U_q(\mathfrak{b}_+) $ 
is inherited from $ \tilde{\rho} $. 
Moreover, we notice that $ \rho(1, X) = \hat{\epsilon}(X) $ because $ \gamma $ is an algebra homomorphism. Finally, 
$ \rho(K_\lambda, 1) = \kappa_\lambda(1)
= 1 = \hat{\epsilon}(K_\lambda) $ and $ \rho(F_i, 1) = \phi_i(1) = 0 = \hat{\epsilon}(F_i) $, which 
implies $ \rho(Y, 1) = \hat{\epsilon}(Y) $ for all $ Y \in U_q(\mathfrak{b}_-) $. This finishes the proof. \end{proof} 

Note that the form $ \rho $ satisfies 
$$
\rho(Y, X) = 0 \qquad \text{ for } X \in U_q(\mathfrak{n}_+)_\alpha, Y \in U_q(\mathfrak{n}_-)_\beta, \alpha \neq -\beta. 
$$

In the sequel we will work both with $\rho$ and the closely related form $\tau$ defined as follows.

\begin{definition}
 \label{def:Drinfeld_pairing}
 We shall refer to the skew-pairing $ \tau: U_q(\mathfrak{b}_+) \otimes U_q(\mathfrak{b}_-) \rightarrow \mathbb{K} $ given by 
 $$
  \tau(X, Y) = \rho(\hat{S}(Y), X). 
 $$
 \nomenclature{$\tau$}{Drinfeld pairing, $\tau:U_q(\lie{b}_+)\otimes U_q(\lie{b}_-) \to \mathbb{K}$}%
 as the Drinfeld pairing.
\end{definition}

\begin{lemma} \label{tauskewpairingproperties}
The Drinfeld pairing $ \tau $ satisfies 
$$
\tau(K_\alpha, K_\beta) = q^{-(\alpha, \beta)}, \qquad \tau(E_i, K_\alpha) = 0 = \tau(K_\beta, F_j), \qquad 
\tau(E_i, F_j) = -\frac{\delta_{ij}}{q_i - q_i^{-1}}
$$
for all $ \alpha, \beta \in \weights $ and $ i,j = 1, \dots, N $. 
Moreover 
\begin{align*}
\tau(X K_\mu, Y K_\nu) = q^{-(\mu, \nu)} \tau(X, Y) \qquad \text{ for } X \in U_q(\mathfrak{n}_+), Y \in U_q(\mathfrak{n}_-). 
\end{align*}
\end{lemma} 

\begin{proof} The formulas in the first part of the Lemma are verified by direct computation from the definition of $ \rho $. 

For $ Y = 1 \in U_q(\mathfrak{n}_-) $ the second claim follows immediately from $ \tau(X K_\mu, K_\nu) = \rho(K_\nu^{-1}, X K_\mu) $. For $ Y = F_j $ we 
see that $ \tau(X K_\mu, F_j) $ and $ \tau(X, F_j) $ vanish unless $ X $ is a multiple of $ E_j $, and we compute 
\begin{align*}
\tau(E_j K_\mu, F_j) &= \tau(E_j, F_j) \tau(K_\mu, 1) + \tau(E_j, K_j^{-1}) \tau(K_\mu, F_j) \\ 
&= -\frac{1}{q_j - q_j^{-1}} = \tau(E_j, F_j). 
\end{align*}

We get the formula for general $ Y $ by induction, using 
\begin{align*}
\tau(X K_\mu, F_j Y K_\nu) &= \tau(X_{(1)} K_\mu, Y K_\nu) \tau(X_{(2)} K_\mu, F_j) \\ 
&= q^{-(\mu, \nu)} \tau(X_{(1)}, Y) \tau(X_{(2)}, F_j) \\ 
&= q^{-(\mu, \nu)} \tau(X, F_j Y), 
\end{align*}
in the inductive step. 
Notice here that $ \hat{\Delta}(X) $ is contained in $ U_q(\mathfrak{n}_+) \otimes U_q(\mathfrak{b}_+) $. \end{proof} 

A crucial property of the pairings $ \rho $ and $ \tau $ is that they can be used to exhibit $ U_q(\mathfrak{g}) $ as a quotient of the Drinfeld 
double $ U_q(\mathfrak{b}_-) \bowtie U_q(\mathfrak{b}_+) $. 

\begin{lemma} \label{lemdrinfeldformcommutation}
For all $ X \in U_q(\mathfrak{b}_+) $ and $ Y \in U_q(\mathfrak{b}_-) $ we have 
\begin{align*}
XY &= \rho(Y_{(1)}, X_{(1)}) Y_{(2)} X_{(2)} \rho(\hat{S}(Y_{(3)}), X_{(3)}) \\
&= \tau(\hat{S}(X_{(1)}), Y_{(1)}) Y_{(2)} X_{(2)} \tau(X_{(3)}, Y_{(3)}) 
\end{align*}
in $ U_q(\mathfrak{g}) $.
\end{lemma} 

\begin{proof} It suffices to check this on generators. We compute 
\begin{align*}
K_\mu K_\nu &= K_\nu K_\mu = q^{(\mu, \nu)} K_\nu K_\mu q^{-(\mu, \nu)}, \\ 
E_i K_\mu &= q^{-(\mu, \alpha_i)} K_\mu E_i \\
&= \rho(K_\mu, 1) K_\mu E_i \rho(\hat{S}(K_\mu), K_i) \\ 
&= \rho(K_\mu, 1) K_\mu E_i \rho(\hat{S}(K_\mu), K_i) 
+ \rho(K_\mu, E_i) K_\mu K_i \rho(\hat{S}(K_\mu), K_i) \\
&\qquad + \rho(K_\mu, 1) K_\mu \rho(\hat{S}(K_\mu), E_i), \\ 
E_i F_j &= F_j E_i - \delta_{ij} \frac{K_i^{-1}}{q_i - q_i^{-1}} + \delta_{ij} \frac{K_i}{q_i - q_i^{-1}} \\
&= \rho(K_j^{-1}, 1) F_j E_i \rho(\hat{S}(1), K_i) + \rho(K_j^{-1}, 1) K_j^{-1} \rho(\hat{S}(F_j), E_i) \\
&\qquad + \rho(F_j, E_i) K_i \rho(\hat{S}(1), K_i). 
\end{align*} 
This yields the claim. \end{proof} 

We note that we can rephrase the commutation relation of Lemma \ref{lemdrinfeldformcommutation} as
\begin{align*}
	X_{(1)}Y_{(1)} \rho(Y_{(2)},X_{(2)}) = \rho(Y_{(1)},X_{(1)}) Y_{(2)}X_{(2)}
\end{align*}
for $ X \in U_q(\mathfrak{b}_+) $ and $ Y \in U_q(\mathfrak{b}_-) $.

In the double $ U_q(\mathfrak{b}_-) \bowtie U_q(\mathfrak{b}_+) $ we shall denote the Cartan generators 
of $ U_q(\mathfrak{b}_+) $ by $ K^+_\lambda $ and the Cartan generators of $ U_q(\mathfrak{b}_-) $ by $ K^-_\lambda $. 

\begin{cor}
The two-sided ideal $ I $
of the Drinfeld double $ U_q(\mathfrak{b}_-) \bowtie U_q(\mathfrak{b}_+) $ generated by the elements 
$ K^-_\lambda - K^+_\lambda $ 
for $ \lambda \in \weights $ is a Hopf ideal, and there is a canonical isomorphism 
$$
(U_q(\mathfrak{b}_-) \bowtie U_q(\mathfrak{b}_+))/I \cong U_q(\mathfrak{g})
$$
of Hopf algebras. 
\end{cor} 

\begin{proof} The canonical embedding maps $ U_q(\mathfrak{b}_+) \rightarrow U_q(\mathfrak{g}) $ and 
$ U_q(\mathfrak{b}_-) \rightarrow U_q(\mathfrak{g}) $ define a Hopf algebra homomorphism 
$ f: (U_q(\mathfrak{b}_-) \bowtie U_q(\mathfrak{b}_+))/I \rightarrow U_q(\mathfrak{g}) $ according to Lemma 
\ref{lemdrinfeldformcommutation}. Conversely, we obtain an algebra homomorphism 
$ g: U_q(\mathfrak{g}) \rightarrow (U_q(\mathfrak{b}_-) \bowtie U_q(\mathfrak{b}_+))/I $ by setting 
$$ 
g(K_\lambda) = K^-_\lambda = K^+_\lambda, \quad g(E_i) = 1 \bowtie E_i,  \quad g(F_j) = F_j \bowtie 1, 
$$ 
since the defining relations of $ U_q(\mathfrak{g}) $ are satisfied by these elements. It is evident 
that $ f $ and $ g $ define inverse isomorphisms. \end{proof}

\subsubsection{The quantum Killing form}

Using the above constructions we can define the quantum Killing form on $ U_q(\mathfrak{g}) $ in terms of the Rosso form for the 
Drinfeld double.  In this subsection, we will need a slightly stronger assumption on $q$, namely $q=s^{2L} \in \mathbb{K}^\times$, again not a root of unity.

According to Definition \ref{def:Rosso_form}, the Rosso form on the 
double $ U_q(\mathfrak{b}_-) \bowtie U_q(\mathfrak{b}_+) $ is given by 
$$
\kappa(Y_1 \bowtie X_1, Y_2 \bowtie X_2) = \tau(\hat{S}(X_1), Y_2) \tau(X_2, \hat{S}(Y_1)), 
$$
for $ X_1, X_2 \in U_q(\mathfrak{b}_+), Y_1, Y_2 \in U_q(\mathfrak{b}_-) $, or equivalently 
$$
\kappa(Y_1 \bowtie \hat{S}^{-1}(X_1), Y_2 \bowtie \hat{S}(X_2)) = \tau(X_1, Y_2) \tau(X_2, Y_1). 
$$
These formulas inspire the following definition. 

\begin{definition} \label{defquantumkilling}
Assume that $ q = s^{2L} \in \mathbb{K}^\times $. 
The quantum Killing form for $ U_q(\mathfrak{g}) $ is the bilinear form $ \kappa: U_q(\mathfrak{g}) \times U_q(\mathfrak{g}) \rightarrow \mathbb{K} $ given by 
$$
\kappa(Y_1 K_\mu \hat{S}^{-1}(X_1), Y_2 K_\nu \hat{S}(X_2)) = q^{(\mu, \nu)/2} \tau(X_1, Y_2) \tau(X_2, Y_1). 
$$
for $ X_1, X_2 \in U_q(\mathfrak{n}_+), Y_1, Y_2 \in U_q(\mathfrak{n}_-), \mu, \nu \in \weights $. 
\end{definition} 

We remark that we need here a slightly stronger requirement on $ q $ than usual in order for the terms $ q^{(\lambda, \mu)/2} $ 
to be well-defined. The factor of $ \tfrac{1}{2} $ is needed to compensate for the 
doubling of the Cartan part in $ U_q(\mathfrak{b}_-) \bowtie U_q(\mathfrak{b}_+) $, see also the remarks further below. 

Heuristically, the formula in Definition \ref{defquantumkilling} can be explained by formally writing 
\begin{align*} 
\kappa(Y_1 K_{\mu/2} &\bowtie K_{\mu/2} \hat{S}^{-1}(X_1), Y_2 K_{\nu/2} \bowtie K_{\nu/2} \hat{S}(X_2)) \\
&= \kappa(Y_1 K_{\mu/2} \bowtie \hat{S}^{-1}(X_1 K_{-\mu/2}), Y_2 K_{\nu/2} \bowtie \hat{S}(X_2 K_{-\nu/2})) \\
&= q^{(\mu, \nu)/2} \tau(X_1, Y_2) \tau(X_2, Y_1), 
\end{align*}
in the Rosso form. This amounts to ``splitting up'' the Cartan part of $ U_q(\mathfrak{g}) $ 
evenly among $ U_q(\mathfrak{b}_+) $ and $ U_q(\mathfrak{b}_-) $.

\begin{prop} \label{killinginvariance} 
The quantum Killing form $ \kappa: U_q(\mathfrak{g}) \times U_q(\mathfrak{g}) \rightarrow \mathbb{K} $ is $ \ad $-invariant, that is, 
$$
\kappa(Z \rightarrow X, Y) = \kappa(X, \hat{S}(Z) \rightarrow Y) 
$$
for all $ X, Y, Z \in U_q(\mathfrak{g}) $. 
\end{prop} 

\begin{proof}  
Let us consider $ X = Y_1 K_\mu \hat{S}^{-1}(X_1) $ and $ Y = Y_2 K_\nu \hat{S}(X_2) $. 
For $ Z = K_\lambda $ and $ X_j \in U_q(\mathfrak{n}_+)_{\beta_j} $ and $ Y_j \in U_q(\mathfrak{n}_-)_{\gamma_j} $ we get 
\begin{align*}
\kappa(K_\lambda \rightarrow (Y_1 K_\mu \hat{S}^{-1}(X_1))&, K_\lambda \rightarrow (Y_2 K_\beta \hat{S}(X_2))) \\
&= q^{(\lambda, \beta_1 + \gamma_1)} q^{(\lambda, \beta_2 + \gamma_2)} \kappa(Y_1 K_\mu \hat{S}^{-1}(X_1), Y_2 K_\nu \hat{S}(X_2)) \\
&= \delta_{\beta_1, -\gamma_2} \delta_{\beta_2, - \gamma_1} \kappa(Y_1 K_\mu \hat{S}^{-1}(X_1), Y_2 K_\nu \hat{S}(X_2)) \\
&= \kappa(Y_1 K_\alpha \hat{S}^{-1}(X_1), Y_2 K_\beta \hat{S}(X_2)). 
\end{align*}
Now consider $ Z = \hat{S}^{-1}( E_i ) $.
We compute 
\begin{align*}
E_i &\rightarrow (Y_2 K_\nu \hat{S}(X_2)) = E_i Y_2 K_\nu \hat{S}(X_2) \hat{S}(K_i) + Y_2 K_\nu \hat{S}(X_2) \hat{S}(E_i) \\
&= q^{(\alpha_i, \beta_2)} E_i Y_2 K_{\nu - \alpha_i} \hat{S}(X_2) + Y_2 K_\nu \hat{S}(E_i X_2) \\
&= q^{(\alpha_i, \beta_2)} \biggl(\tau(\hat{S}(E_i), (Y_2)_{(1)}) (Y_2)_{(2)} K_i \tau(K_i, (Y_2)_{(3)}) K_{\nu - \alpha_i} \hat{S}(X_2) \\ 
&\qquad + \tau(1, (Y_2)_{(1)}) (Y_2)_{(2)} E_i \tau(K_i, (Y_2)_{(3)}) K_{\nu - \alpha_i} \hat{S}(X_2) \\
&\qquad + \tau(1, (Y_2)_{(1)}) (Y_2)_{(2)} \tau(E_i, (Y_2)_{(3)}) K_{\nu - \alpha_i} \hat{S}(X_2) \biggr) 
+ Y_2 K_\nu \hat{S}(E_i X_2) \\
&= q^{(\alpha_i, \beta_2)} \biggl(\tau(\hat{S}(E_i), (Y_2)_{(1)}) (Y_2)_{(2)} K_\nu \hat{S}(X_2) \\ 
&\qquad + Y_2 E_i K_{\nu - \alpha_i} \hat{S}(X_2) + \tau(E_i, (Y_2)_{(2)}) (Y_2)_{(1)} K_{\nu - \alpha_i} \hat{S}(X_2) \biggr) 
+ Y_2 K_\nu \hat{S}(E_i X_2) \\
&= q^{(\alpha_i, \beta_2)} \tau(\hat{S}(E_i), (Y_2)_{(1)}) (Y_2)_{(2)} K_\nu \hat{S}(X_2) - q^{(\alpha_i, \beta_2)} Y_2 \hat{S}(E_i) K_\nu \hat{S}(X_2) \\
&\qquad + q^{(\alpha_i, \beta_2)} \tau(E_i, (Y_2)_{(2)}) (Y_2)_{(1)} K_{\nu - \alpha_i} \hat{S}(X_2) + Y_2 K_\nu \hat{S}(E_i X_2) \\
&= -q^{(\alpha_i, \beta_2)} \tau(E_i K_i^{-1}, (Y_2)_{(1)}) (Y_2)_{(2)} K_\nu \hat{S}(X_2) \\
&\qquad - q^{(\alpha_i, \beta_2)} q^{-(\alpha_i, \nu)} Y_2 K_\nu \hat{S}(X_2 E_i) \\
&\qquad + q^{(\alpha_i, \beta_2)} \tau(E_i, (Y_2)_{(2)}) (Y_2)_{(1)} K_{\nu - \alpha_i} \hat{S}(X_2) \\
&\qquad + Y_2 K_\nu \hat{S}(E_i X_2).
\end{align*}
In the third term, note that $\tau(E_i,\cdot)$ annihilates all PBW vectors in $U_q(\mathfrak{b}_-)$ except those of the form $F_i K_\lambda$ for some $\lambda \in \weights$.   We therefore obtain
\begin{align*}
\kappa(Y_1 K_\mu \hat{S}^{-1}(X_1), E_i &\rightarrow (Y_2 K_\nu \hat{S}(X_2))) \\
&= -q^{(\alpha_i, \beta_2)} q^{(\mu, \nu)/2} \tau(E_i K_i^{-1}, (Y_2)_{(1)}) \tau(X_1, (Y_2)_{(2)}) \tau(X_2, Y_1) \\
&\qquad - q^{(\alpha_i, \beta_2)} q^{-(\alpha_i, \nu)} q^{(\mu, \nu)/2} \tau(X_1, Y_2) \tau(X_2 E_i, Y_1) \\
&\qquad + q^{(\alpha_i, \beta_2)} q^{(\mu, - \alpha_i)/2} q^{(\mu, \nu - \alpha_i)/2} \tau(E_i, (Y_2)_{(2)}) \tau(X_1, (Y_2)_{(1)}) \tau(X_2, Y_1) \\
&\qquad + q^{(\mu, \nu)/2} \tau(X_1, Y_2) \tau(E_i X_2, Y_1) \\
&= -q^{(\alpha_i, \beta_2 - \beta_1)} q^{(\mu, \nu)/2} \tau(E_i X_1, Y_2) \tau(X_2, Y_1) \\
&\qquad - q^{(\alpha_i, \beta_2)} q^{-(\alpha_i, \nu)} q^{(\mu, \nu)/2} \tau(X_1, Y_2) \tau(X_2 E_i, Y_1) \\
&\qquad + q^{(\alpha_i, \beta_2)} q^{-(\mu, \alpha_i)} q^{(\mu, \nu)/2} \tau(X_1 E_i, Y_2) \tau(X_2, Y_1) \\
&\qquad + q^{(\mu, \nu)/2} \tau(X_1, Y_2) \tau(E_i X_2, Y_1). 
\end{align*}
On the other hand, using 
\begin{align*}
E_i &\rightarrow (Y_1 K_\mu \hat{S}^{-1}(X_1)) = E_i Y_1 K_\mu \hat{S}^{-1}(X_1) \hat{S}(K_i) + Y_1 K_\mu \hat{S}^{-1}(X_1) \hat{S}(E_i) \\
&= q^{(\alpha_i, \beta_1)} E_i Y_1 K_{\mu - \alpha_i} \hat{S}^{-1}(X_1) + q^{(\alpha_i, \alpha_i)} Y_1 K_\mu \hat{S}^{-1}(X_i) \hat{S}^{-1}(E_i) \\
&= q^{(\alpha_i, \beta_1)} \biggl(\tau(\hat{S}(E_i), (Y_1)_{(1)}) (Y_1)_{(2)} K_i \tau(K_i, (Y_1)_{(3)}) K_{\mu - \alpha_i} \hat{S}^{-1}(X_1) \\ 
&\qquad + \tau(1, (Y_1)_{(1)}) (Y_1)_{(2)} E_i \tau(K_i, (Y_1)_{(3)}) K_{\mu - \alpha_i} \hat{S}^{-1}(X_1) \\
&\qquad + \tau(1, (Y_1)_{(1)}) (Y_1)_{(2)} \tau(E_i, (Y_1)_{(3)}) K_{\mu - \alpha_i} \hat{S}^{-1}(X_1) \biggr) 
+ q^{(\alpha_i, \alpha_i)} Y_1 K_\mu \hat{S}^{-1}(E_i X_1) \\
&= q^{(\alpha_i, \beta_1)} \biggl(\tau(\hat{S}(E_i), (Y_1)_{(1)}) (Y_1)_{(2)} K_\mu \hat{S}^{-1}(X_1) \\ 
&\qquad + Y_1 E_i K_{\mu - \alpha_i} \hat{S}^{-1}(X_1) + \tau(E_i, (Y_1)_{(2)}) (Y_1)_{(1)} K_{\mu - \alpha_i} \hat{S}^{-1}(X_1) \biggr) 
+ q^{(\alpha_i, \alpha_i)} Y_1 K_\mu \hat{S}^{-1}(E_i X_1) \\
&= q^{(\alpha_i, \beta_1)} \tau(\hat{S}(E_i), (Y_1)_{(1)}) (Y_1)_{(2)} K_\mu \hat{S}^{-1}(X_1) 
- q^{(\alpha_i, \beta_1 + \alpha_i)} Y_1 \hat{S}^{-1}(E_i) K_\mu \hat{S}^{-1}(X_1) \\
&\qquad + q^{(\alpha_i, \beta_1)} \tau(E_i, (Y_1)_{(2)}) (Y_1)_{(1)} K_{\mu - \alpha_i} \hat{S}^{-1}(X_1) 
+ q^{(\alpha_i, \alpha_i)} Y_1 K_\mu \hat{S}^{-1}(E_i X_1) \\
&= -q^{(\alpha_i, \beta_1)} \tau(E_i K_i^{-1}, (Y_1)_{(1)}) (Y_1)_{(2)} K_\mu \hat{S}^{-1}(X_1) \\
&\qquad - q^{(\alpha_i, \beta_1 + \alpha_i)} q^{-(\alpha_i, \mu)} Y_1 K_\mu \hat{S}^{-1}(X_1 E_i) \\
&\qquad + q^{(\alpha_i, \beta_1)} \tau(E_i, (Y_1)_{(2)}) (Y_1)_{(1)} K_{\mu - \alpha_i} \hat{S}^{-1}(X_1) \\
&\qquad + q^{(\alpha_i, \alpha_i)} Y_1 K_\mu \hat{S}^{-1}(E_i X_1). 
\end{align*}
we calculate
\begin{align*}
\kappa(\hat{S}^{-1}(E_i) &
\rightarrow (Y_1 K_\mu \hat{S}^{-1}(X_1)), Y_2 K_\nu \hat{S}(X_2))
 \\
& = - \kappa(E_i \rightarrow (Y_1 K_\mu \hat{S}^{-1}(X_1)), K_i \rightarrow (Y_2 K_\nu \hat{S}(X_2))) \\
&= - q^{(\alpha_i, \gamma_2 + \beta_2)} \kappa(E_i \rightarrow (Y_1 K_\mu \hat{S}^{-1}(X_1)), Y_2 K_\nu \hat{S}(X_2)) \\
&= - q^{(\alpha_i, \gamma_2 + \beta_2)} \biggl(
-q^{(\alpha_i, \beta_1)} q^{(\mu, \nu)/2} \tau(E_i K_i^{-1}, (Y_1)_{(1)}) \tau(X_2, (Y_1)_{(2)}) \tau(X_1, Y_2) \\
&\qquad - q^{(\alpha_i, \beta_1 + \alpha_i)} q^{-(\alpha_i, \mu)} q^{(\mu, \nu)/2} \tau(X_2, Y_1) \tau(X_1 E_i, Y_2) \\
&\qquad + q^{(\alpha_i, \beta_1)} q^{(\nu, -\alpha_i)/2}  q^{(\mu - \alpha_i, \nu)/2} \tau(E_i, (Y_1)_{(2)}) \tau(X_2, (Y_1)_{(1)}) \tau(X_1, Y_2) \\
&\qquad + q^{(\alpha_i, \alpha_i)} q^{(\mu, \nu)/2} \tau(X_2, Y_1) \tau(E_i X_1, Y_2)\biggr) \\
&= - q^{(\alpha_i, \gamma_2 + \beta_2)} \biggl(
-q^{(\alpha_i, \beta_1)} q^{(\mu, \nu)/2} \tau(E_i K_i^{-1} X_2, Y_1) \tau(X_1, Y_2) \\
&\qquad - q^{(\alpha_i, \beta_1 + \alpha_i)} q^{-(\alpha_i, \mu)} q^{(\mu, \nu)/2} \tau(X_2, Y_1) \tau(X_1 E_i, Y_2) \\
&\qquad + q^{(\alpha_i, \beta_1)} q^{-(\nu, \alpha_i)} q^{(\mu, \nu)/2} \tau(X_2 E_i, Y_1) \tau(X_1, Y_2) \\
&\qquad + q^{(\alpha_i, \alpha_i)} q^{(\mu, \nu)/2} \tau(X_2, Y_1) \tau(E_i X_1, Y_2)\biggr) \\
&=  q^{(\mu, \nu)/2} \tau(E_i X_2, Y_1) \tau(X_1, Y_2) \\
&\qquad + q^{(\alpha_i, \beta_2)} q^{-(\alpha_i, \mu)} q^{(\mu, \nu)/2} \tau(X_2, Y_1) \tau(X_1 E_i, Y_2) \\
&\qquad - q^{(\alpha_i, \beta_2)} q^{(\nu, -\alpha_i)} q^{(\mu, \nu)/2} \tau(X_2 E_i, Y_1) \tau(X_1, Y_2) \\
&\qquad - q^{(\alpha_i, \beta_2 - \beta_1)} q^{(\mu, \nu)/2} \tau(X_2, Y_1) \tau(E_i X_1, Y_2),
\end{align*}
 where in the first term of the last equality we are using the fact that $\tau(E_iX_2,Y_1)\tau(X_1,Y_2)$ is zero unless $-\gamma_1 = \alpha_i+\beta_2$ and $-\gamma_2 = \beta_1$, and similarly for each of the other terms.

In a similar way one proceeds for $ Z=\hat{S}^{-1}( F_i )$. This yields the claim. \end{proof} 

Consider also the bilinear form 
\begin{align*}
(X_1 K_\mu \hat{S}(Y_1), Y_2 K_\nu \hat{S}(X_2)) &= q^{-(\mu, \nu)/2} \tau(X_1, Y_2) \tau(X_2, Y_1) \\
&= \kappa(\hat{S}^{-1}(X_1 K_\mu \hat{S}(Y_1)), Y_2 K_\nu \hat{S}(X_2)). 
\end{align*}
\label{nom:Killing_form_variant}%
on $ U_q(\mathfrak{g}) $. Using 
$$
\hat{S}(X) \leftarrow \hat{S}(Z) = \hat{S}^2(Z_{(2)}) \hat{S}(X) \hat{S}(Z_{(1)}) = \hat{S}(Z_{(1)} X \hat{S}(Z_{(2)})) = \hat{S}(Z \rightarrow X) 
$$
and Proposition \ref{killinginvariance} we get 
$$
(X \leftarrow Z, Y) = (X, Z \rightarrow Y) 
$$
for all $ X, Y, Z \in U_q(\mathfrak{g}) $. 
This will be crucial when we use the quantum Killing form to express the locally finite part of $ U_q(\mathfrak{g}) $ in terms of $ \Poly(G_q) $ 
in Theorem \ref{adjointduality}.

\subsubsection{Computation of the Drinfeld pairing}

In this subsection we use the PBW-basis to compute the Drinfeld pairing. For the main part of the arguments we follow the approach 
of Tanisaki \cite{Tanisakidrinfeld}. 

Let us start with some preliminaries. The $ q $-exponential function $ \exp_q $ is defined by 
$$
\exp_q(x) = \sum_{n = 0}^\infty \frac{q^{n (n - 1)/2}}{[n]_q!} x^n.  
$$
\nomenclature[o$exp_q$]{$\exp_q$}{$q$-exponential function}%
Here $ x $ is a formal variable, and the expression on the right hand side can be viewed as an element of $ \mathbb{Q}(q)[[x]] $. 

\begin{lemma} \label{qexpinverse}
The formal series $ \exp_q(x) $ is invertible, with inverse 
$$
\exp_q(x)^{-1} = \sum_{n = 0}^\infty \frac{q^{-n (n - 1)/2}}{[n]_q!} (-x)^n = \exp_{q^{-1}}(-x).  
$$
\end{lemma} 

\begin{proof} We formally multiply and collect terms of common powers in $ x $, explicitly, 
\begin{align*}
\exp_q(x) \exp_{q^{-1}}(-x) &= 
\biggl(\sum_{m = 0}^\infty \frac{q^{m (m - 1)/2}}{[m]_q!} x^m\biggr)\biggl(\sum_{n = 0}^\infty \frac{q^{-n (n - 1)/2}}{[n]_q!} (-x)^n\biggr) \\
&= \sum_{k = 0}^\infty \sum_{l = 0}^k (-1)^{k - l} \frac{q^{l (l - 1)/2} q^{-(k - l)(k - l - 1)/2}}{[l]_q!} \frac{1}{[k - l]_q!} x^k \\
&= \sum_{k = 0}^\infty (-1)^k q^{(-k^2 + k)/2} \biggl(\sum_{l = 0}^k (-1)^l \frac{q^{(l^2 - l + kl + lk - l^2 - l)/2}}{[k]_q!} 
\begin{bmatrix}
k \\
l
\end{bmatrix}_q
\biggr) x^k \\
&= \sum_{k = 0}^\infty (-1)^k \frac{q^{(-k^2 + k)/2}}{[k]_q!} \biggl(\sum_{l = 0}^k (-1)^l q^{kl - l}
\begin{bmatrix}
k \\
l
\end{bmatrix}_q
\biggr) x^k.  
\end{align*}
Hence the claim follows from Lemma \ref{lembinomialvanishing}. \end{proof}

For $ 1 \leq i \leq N $ and integrable $ U_q(\mathfrak{g}) $-modules $ V, W $, 
define a linear operator $ Z_i $ on $ V \otimes W $ by 
$$
Z_i = \exp_{q_i}((q_i - q_i^{-1})(E_i \otimes F_i)).  
$$
According to Lemma \ref{qexpinverse} the operator $ Z_i $ is invertible with inverse 
$$
Z_i^{-1} = \exp_{q_i^{-1}}(-(q_i - q_i^{-1})(E_i \otimes F_i)).  
$$
Recall moreover the definition of the operators $ \T_i $ acting on integrable $ U_q(\mathfrak{g}) $-modules from Section \ref{secbraid}. 
\begin{prop} \label{Tdiagonal}
Let $ V, W $ be finite dimensional integrable $ U_q(\mathfrak{g}) $-modules and $ 1 \leq i \leq N $. Then we have 
$$
\T_i(v \otimes w) = (\T_i \otimes \T_i) Z_i(v \otimes w)
$$
for all $ v \otimes w \in V \otimes W $. 
\end{prop}

\begin{proof} Since the operators $ \T_i $ and $ Z_i $ are defined in terms of $ U_{q_i}(\mathfrak{g}_i) \subset U_q(\mathfrak{g}) $ 
it suffices to consider the case $ \mathfrak{g} = \mathfrak{sl}(2, \mathbb{K}) $. We shall therefore restrict to this case in the sequel, writing 
$ Z = \exp_q(-(q - q^{-1})(E \otimes F)) $. 

Let us first show that both $ \T $ and $ (\T \otimes \T) Z $ commute in the same way with the diagonal action of $ E $ and $ F $. 
By Lemma \ref{tcommutationlemma} we have 
\begin{align*}
\T(E \cdot (v \otimes w)) &= -K^2 F \cdot \T(v \otimes w) \\
&= (-K^2 F \otimes K^2 - 1 \otimes K^2F) \cdot \T(v \otimes w). 
\end{align*}
On the other hand,  
\begin{align*}
(\T \otimes \T)&Z(E \cdot (v \otimes w)) = (\T \otimes \T)Z(E \cdot v \otimes K^2 \cdot w + v \otimes E \cdot w) \\
&= (\T \otimes \T)((E \otimes K^2) \cdot \exp_q((q - q^{-1})q^2(E \otimes F))(v \otimes w)) \\
&\qquad + (\T \otimes \T) Z((1 \otimes E) \cdot (v \otimes w)) \\
&= -(K^2 F \otimes K^{-2}) (\T \otimes \T)(\exp_q((q - q^{-1})q^2(E \otimes F))(v \otimes w)) \\
&\qquad + (\T \otimes \T) Z((1 \otimes E) \cdot (v \otimes w)). 
\end{align*}
In order to compute $ (\T \otimes \T) Z((1 \otimes E) \cdot (v \otimes w)) $ consider $ x = (q - q^{-1}) (E \otimes F) $. 
Let us prove by induction that
\begin{align*}
[x^n, 1 \otimes E] &= \sum_{r = 0}^{n - 1} (E \otimes (K^{-2} q^{-2r} - K^2 q^{2r})) x^{n - 1} \\
&= (q^{-n + 1} [n]_q (E \otimes K^{-2}) - q^{n - 1} [n]_q (E \otimes K^2)) x^{n - 1} 
\end{align*}
for all $ n \in \mathbb{N} $. 
Indeed, for $ n = 1 $ one checks 
\begin{align*}
[x, 1 \otimes E] &= (q - q^{-1}) (E \otimes [F,E]) = E \otimes (K^{-2} - K^2), 
\end{align*}
and for the inductive step we compute 
\begin{align*}
[x^n, 1 \otimes E] &= x[x^{n - 1}, 1 \otimes E] + [x, 1 \otimes E] x^{n - 1} \\
&= x \sum_{r = 0}^{n - 2} (E \otimes (K^{-2} q^{-2r} - K^2 q^{2r})) x^{n - 2} + (E \otimes (K^{-2} - K^2))x^{n - 1} \\
&= \sum_{r = 0}^{n - 2} (E \otimes K^{-2} q^{-2r - 2} - E \otimes K^2 q^{2r + 2}) x^{n - 1} + (E \otimes (K^{-2} - K^2))x^{n - 1} \\
&= \sum_{r = 0}^{n - 1} (E \otimes (K^{-2} q^{-2r} - K^2 q^{2r})) x^{n - 1}. 
\end{align*}
Therefore we formally obtain 
\begin{align*}
Z(1 \otimes E) &= \sum_{n = 0}^\infty \frac{q^{n (n - 1)/2}}{[n]_q!} x^n(1 \otimes E) \\
&= (1 \otimes E) \exp_q(x) + \sum_{n = 1}^\infty \frac{q^{n (n - 1)/2}}{[n]_q!} q^{-n + 1} [n]_q (E \otimes K^{-2}) x^{n - 1} \\
&\qquad - \sum_{n = 1}^\infty \frac{q^{n (n - 1)/2}}{[n]_q!} q^{n - 1} [n]_q (E \otimes K^2) x^{n - 1} \\
&= (1 \otimes E) \exp_q(x) + (E \otimes K^{-2}) \sum_{n = 1}^\infty \frac{q^{(n - 1) (n - 2)/2}}{[n - 1]_q!} x^{n - 1} \\
&\qquad - (E \otimes K^2) \sum_{n = 1}^\infty \frac{q^{(n - 1)(n - 2)/2}}{[n - 1]_q!} q^{2(n - 1)} x^{n - 1} \\
&= (1 \otimes E) \exp_q(x) + (E \otimes K^{-2}) \exp_q(x) - (E \otimes K^2) \exp_q(q^2 x).  
\end{align*}
Summarizing the above calculations we get 
\begin{align*}
(\T \otimes \T)Z(1 \otimes E) &= -(1 \otimes K^2 F) (\T \otimes \T) Z - (K^2 F \otimes K^2) (\T \otimes \T) Z \\
&+ (K^2 F \otimes K^{-2}) (\T \otimes \T)\exp_q(q^2 x), 
\end{align*}
and hence 
\begin{align*}
(\T \otimes \T)&Z(E \cdot (v \otimes w)) \\
&= -(K^2 F \otimes K^{-2}) (\T \otimes \T)(\exp_q((q - q^{-1})q^2(E \otimes F))(v \otimes w)) \\
&\qquad -(1 \otimes K^2 F) (\T \otimes \T) Z(v \otimes w) - (K^2 F \otimes K^2) (\T \otimes \T) Z(v \otimes w) \\
&+ (K^2 F \otimes K^{-2}) (\T \otimes \T)\exp_q((q - q^{-1})q^2(E \otimes F))(v \otimes w)) \\
&= -(1 \otimes K^2 F) (\T \otimes \T) Z(v \otimes w) - (K^2 F \otimes K^2) (\T \otimes \T) Z(v \otimes w). 
\end{align*}
This shows that both operators commute in the same way with the diagonal action of $ E $. 
In a similar way we compute 
\begin{align*}
\T(F \cdot (v \otimes w)) &= -E K^{-2} \cdot \T(v \otimes w) \\
&= (-E K^{-2} \otimes 1 - K^{-2} \otimes EK^{-2}) \cdot \T(v \otimes w). 
\end{align*}
On the other hand, 
\begin{align*}
(\T \otimes \T)&Z(F \cdot (v \otimes w)) = (\T \otimes \T)Z(F \cdot v \otimes w + K^{-2} \cdot v \otimes F \cdot w) \\
&= (\T \otimes \T) Z((F \otimes 1) \cdot (v \otimes w)) \\ 
&\qquad + (\T \otimes \T)((K^{-2} \otimes F) \cdot \exp_q((q - q^{-1})q^2(E \otimes F))(v \otimes w)) \\
&= (\T \otimes \T) Z((F \otimes 1) \cdot (v \otimes w)) \\ 
&\qquad -(K^2 \otimes EK^{-2}) (\T \otimes \T)(\exp_q((q - q^{-1})q^2(E \otimes F))(v \otimes w)). 
\end{align*}
If $ x = (q - q^{-1}) (E \otimes F) $, then we obtain 
\begin{align*}
[x^n, F \otimes 1] &= \sum_{r = 0}^{n - 1} ((K^2 q^{-2r} - K^{-2} q^{2r}) \otimes F) x^{n - 1} \\
&= (q^{-n + 1} [n]_q K^2 \otimes F - q^{n - 1} [n]_q K^{-2} \otimes F) x^{n - 1} 
\end{align*}
for all $ n \in \mathbb{N} $ by induction. 
Indeed, for $ n = 1 $ one checks 
\begin{align*}
[x, F \otimes 1] &= (q - q^{-1}) ([E, F] \otimes F) = (K^2 - K^{-2}) \otimes F, 
\end{align*}
and for the inductive step we compute 
\begin{align*}
[x^n, F \otimes 1] &= x[x^{n - 1}, F \otimes 1] + [x, F \otimes 1] x^{n - 1} \\
&= x \sum_{r = 0}^{n - 2} ((K^2 q^{-2r} - K^{-2} q^{2r}) \otimes F) x^{n - 2} + ((K^2 - K^{-2}) \otimes F)x^{n - 1} \\
&= \sum_{r = 0}^{n - 2} ((K^2 q^{-2r - 2} - K^{-2} q^{2r + 2}) \otimes F) x^{n - 2} + ((K^2 - K^{-2}) \otimes F)x^{n - 1} \\
&= \sum_{r = 0}^{n - 1} ((K^2 q^{-2r} - K^{-2} q^{2r}) \otimes F) x^{n - 1}. 
\end{align*}
We formally obtain 
\begin{align*}
Z(F \otimes 1) &= \sum_{n = 0}^\infty \frac{q^{n (n - 1)/2}}{[n]_q!} x^n(F \otimes 1) \\
&= (F \otimes 1) \exp_q(x) + \sum_{n = 1}^\infty \frac{q^{n (n - 1)/2}}{[n]_q!} q^{-n + 1} [n]_q (K^2 \otimes F) x^{n - 1} \\
&\qquad - \sum_{n = 1}^\infty \frac{q^{n (n - 1)/2}}{[n]_q!} q^{n - 1} [n]_q (K^{-2} \otimes F) x^{n - 1} \\
&= (F \otimes 1) \exp_q(x) + (K^2 \otimes F) \sum_{n = 1}^\infty \frac{q^{(n - 1) (n - 2)/2}}{[n - 1]_q!} x^{n - 1} \\
&\qquad - (K^{-2} \otimes F) \sum_{n = 1}^\infty \frac{q^{(n - 1)(n - 2)/2}}{[n]_q!} q^{2(n - 1)} x^{n - 1} \\
&= (F \otimes 1) \exp_q(x) + (K^2 \otimes F) \exp_q(x) - (K^{-2} \otimes F) \exp_q(q^2 x).  
\end{align*}
Summarising the above calculations we get 
\begin{align*}
(\T \otimes \T)Z(F \otimes 1) &= -(EK^{-2} \otimes 1) (\T \otimes \T) Z - (K^{-2} \otimes EK^{-2}) (\T \otimes \T) Z \\
&+ (K^2 \otimes EK^{-2}) (\T \otimes \T)\exp_q(q^2 x), 
\end{align*}
and hence 
\begin{align*}
(\T \otimes \T)&Z(F \cdot (v \otimes w)) \\
&= -(EK^{-2} \otimes 1) (\T \otimes \T) Z(v \otimes w) - (K^{-2} \otimes EK^{-2}) (\T \otimes \T) Z(v \otimes w) \\
&+ (K^2 \otimes EK^{-2}) (\T \otimes \T)\exp_q(q^2 (q - q^{-1})(E \otimes F)(v \otimes w)) \\
&-(K^2 \otimes EK^{-2}) (\T \otimes \T)(\exp_q((q - q^{-1})q^2(E \otimes F))(v \otimes w)) \\
&= -(EK^{-2} \otimes 1) (\T \otimes \T) Z(v \otimes w) - (K^{-2} \otimes EK^{-2}) (\T \otimes \T) Z(v \otimes w). 
\end{align*}
That is, both operators commute in the same way with the diagonal action of $ F $. Let us point out that $ \T $ and $ (\T \otimes \T) Z $ 
commute also in the same way with the diagonal action of $ K_\mu $ for $ \mu \in \weights $, as one checks easily by weight 
considerations. 

We show next that $ \T $ and $ (\T \otimes \T) Z $ agree on $ V(m) \otimes V(1/2) $ for all $ m \in \frac{1}{2} \mathbb{N}_0 $. 
For $ m = 0 $ there is nothing to prove. Hence assume that $ m \geq \frac{1}{2} $ and consider the tensor product decomposition 
$$ 
V(m) \otimes V(1/2) \cong V(m + 1/2) \oplus V(m - 1/2). 
$$
Since both operators commute in the same way with the diagonal action of $ U_q(\mathfrak{g}) $ it suffices to show that $ \T $ 
and $ (\T \otimes \T)Z $ agree on 
the lowest weight vectors in the irreducible components $V(m + 1/2)$ and $V(m - 1/2)$ of the tensor product.  Up to a scalar, these lowest weight vectors are given by
\begin{align*}
 v_{(-m-1/2)} & = v_{(-m)} \otimes v_{(-1/2)}, \\
 v_{(-m+1/2)} & = [2m]_q v_{(-m)} \otimes v_{(1/2)} - q^{2m} v_{(-m + 1)} \otimes v_{(-1/2)}, 
\end{align*}
as can be checked by computing the diagonal action of $F$ on each.

For the case of $v_{(-m-1/2)}$, note that $Z$ fixes $v_{(-m)} \otimes v_{(-1/2)}$, and so Proposition \ref{texplicit} gives
\[
 (\T \otimes \T)Z(v_{(-m)} \otimes v_{(-1/2)}) = (-1)^{2m}q^{2m}v_{(m)} \otimes -q v_{(1/2)},
\]
while
\[
 \T(v_{(-m-1/2})) = (-1)^{2m+1} q^{2m+1} v_{(m+1/2)}.
\]
For $v_{(-m+1/2)}$, we
have 
\begin{align*}
Z(v_{(-m + 1/2)}) &= v_{(-m + 1/2)} + (q - q^{-1})(E \otimes F) \cdot v_{(-m + 1/2)} \\
&= [2m]_q v_{(-m)} \otimes v_{(1/2)} - q^{2m} v_{(-m + 1)} \otimes v_{(-1/2)} \\
&\qquad + (q - q^{-1})[2m]_q v_{(-m + 1)} \otimes v_{(-1/2)} \\
&= [2m]_q v_{(-m)} \otimes v_{(1/2)} - q^{-2m} v_{(-m + 1)} \otimes v_{(-1/2)}. 
\end{align*} 
Using the formulas from Proposition \ref{texplicit} we get 
\begin{align*}
(\T &\otimes \T)Z(v_{(-m + 1/2)}) \\
&= (\T \otimes \T)([2m]_q v_{(-m)} \otimes v_{(1/2)} - q^{-2m} v_{(-m + 1)} \otimes v_{(-1/2)}) \\
&= (-1)^{2m} q^{2m} [2m]_q v_{(m)} \otimes v_{(-1/2)} - (-1)^{2m - 1} q^{-2m} q^{(2m - 1)2}v_{(m - 1)} \otimes (-q v_{(1/2)}) \\
&= (-1)^{2m} q^{2m} [2m]_q v_{(m)} \otimes v_{(-1/2)} - (-1)^{2m} q^{2m - 1} v_{(m - 1)} \otimes v_{(1/2)}. 
\end{align*} 
On the other hand, from the definition of $\T$ and Lemma \ref{lemEiFicomult} we obtain
\begin{align*}
&\T(v_{(-m + 1/2)}) = (-1)^{2m - 1} q^{2m - 1} E^{(2m - 1)} \cdot v_{(-m + 1/2)} \\
&= (-1)^{2m - 1} q^{2m - 1} E^{(2m - 1)} \cdot ([2m]_q v_{(-m)} \otimes v_{(1/2)} - q^{2m} v_{(-m + 1)} \otimes v_{(-1/2)}) \\
&= (-1)^{2m - 1} q^{2m - 1} ([2m]_q E^{(2m - 1)} \cdot v_{(-m)} \otimes K^{2(2m - 1)} \cdot v_{(1/2)} \\
&\qquad - q^{2m} E^{(2m - 1)} \cdot v_{(-m + 1)} \otimes K^{2(2m - 1)} \cdot v_{(-1/2)} \\
&\qquad - q^{2m} q^{2m - 2} E^{(2m - 2)} \cdot v_{(-m + 1)} \otimes E K^{2(2m - 2)} \cdot v_{(-1/2)}) \\
&= (-1)^{2m - 1} q^{2m - 1} [2m]_q q^{2m - 1} v_{(m - 1)} \otimes v_{(1/2)} \\
&\qquad + (-1)^{2m} q^{2m - 1} q^{2m} [2m]_q v_{(m)} \otimes q^{-(2m - 1)} v_{(-1/2)} \\
&\qquad + (-1)^{2m} q^{2m - 1} q^{2m} q^{2m - 2} q^{-(2m - 2)} [2m - 1]_q v_{(m - 1)} \otimes v_{(1/2)} \\
&= (-1)^{2m - 1} q^{4m - 2} [2m]_q v_{(m - 1)} \otimes v_{(1/2)} \\
&\qquad + (-1)^{2m} q^{2m} [2m]_q v_{(m)} \otimes v_{(-1/2)} \\
&\qquad + (-1)^{2m} q^{4m - 1} [2m - 1]_q v_{(m - 1)} \otimes v_{(1/2)} \\
&= (-1)^{2m - 1} q^{2m - 1} v_{(m - 1)} \otimes v_{(1/2)} + (-1)^{2m} q^{2m} [2m]_q v_{(m)} \otimes v_{(-1/2)}, 
\end{align*} 
using 
$$
q^{4m - 2}[2m]_q - q^{4m - 1}[2m - 1]_q = q^{2m - 1} 
$$
in the last step. 

Now we want to show that $ \T $ and $ (\T \otimes \T)Z $ agree on $ V(m) \otimes V(n) $ for all $ n \in \frac{1}{2} \mathbb{N}_0 $. 
Since both operators commute in the same way with the diagonal action of $ U_q(\mathfrak{g}) $ it suffices to show that 
they agree on $ V(m) \otimes v_{(-n)} $, where $ v_{(-n)} \in V(n) $ is a lowest weight vector. 
The case $ n = 0 $ is trivial, and the case $ n = 1/2 $ follows from our above calculations. Let $n>1/2$ and assume that the assertion is 
proved for all $ m $ and all $ k < n $. Since $ Z(v \otimes v_{(-n)}) = v \otimes v_{(-n)} $ for all $ v \in V(m) $
we have to show $ \T(v \otimes v_{(-n)}) = \T(v) \otimes \T(v_{(-n)}) $.  
To this end consider the inclusion $ V(n) \subset V(n - 1/2) \otimes V(1/2) $. 
Then, up to a scalar, the vector $ v_{(-n)} $ identifies with $ v_{(-n + 1/2)} \otimes v_{(-1/2)} $, and we obtain 
\begin{align*}
\T(v \otimes v_{(-n + 1/2)} \otimes v_{(-1/2)}) &= \T(v \otimes v_{(-n + 1/2)}) \otimes \T(v_{(-1/2)}) \\
&= \T(v) \otimes \T(v_{(-n + 1/2)}) \otimes \T(v_{(-1/2)}) \\
&= \T(v) \otimes \T(v_{(-n + 1/2)} \otimes v_{(-1/2)})
\end{align*}
by our induction hypothesis.  
This finishes the proof. \end{proof} 

Recall that the automorphism $ \T_i: U_q(\mathfrak{g}) \rightarrow U_q(\mathfrak{g}) $ is induced by conjugation 
with the operators $ \T_i $ on integrable $ U_q(\mathfrak{g}) $-modules. As a consequence of Proposition \ref{Tdiagonal} we obtain 
$$
\hat{\Delta}(\T_i^{-1}(X)) = Z_i^{-1} (\T_i^{-1} \otimes \T_i^{-1}) \hat{\Delta}(X) (\T_i \otimes \T_i) Z_i \\
= Z_i^{-1} (\T_i^{-1} \otimes \T_i^{-1})(\hat{\Delta}(X)) Z_i 
$$
for all $ X \in U_q(\mathfrak{g}) $. 

\begin{lemma} \label{drinfeldpairinghelp}
The algebras $\T_i(U_q(\mathfrak{b}_+))$ and $\T_i(U_q(\mathfrak{b}_-))$ are a right coideal subalgebra and a left coideal subalgebra of $U_q(\lie{g})$, respectively, meaning that they satisfy
\begin{align*}
\hat{\Delta}(\T_i(U_q(\mathfrak{b}_+))) &\subset \T_i(U_q(\mathfrak{b}_+)) \otimes U_q(\mathfrak{g}), \\
\hat{\Delta}(\T_i(U_q(\mathfrak{b}_-))) &\subset U_q(\mathfrak{g}) \otimes \T_i(U_q(\mathfrak{b}_-)).
\end{align*}
Moreover,
\begin{align*}
U_q(\mathfrak{n}_+) \cap \T_i(U_q(\mathfrak{b}_+)) = U_q(\mathfrak{n}_+) \cap \T_i(U_q(\mathfrak{n}_+)), \\
U_q(\mathfrak{n}_-) \cap \T_i(U_q(\mathfrak{b}_-)) = U_q(\mathfrak{n}_-) \cap \T_i(U_q(\mathfrak{n}_-)). 
\end{align*}
\end{lemma} 

\begin{proof} For the first formula we compute 
\begin{align*}
\hat{\Delta}(\T_i(X))) &= (\T_i \otimes \T_i)(Z_i \hat{\Delta}(X) Z_i^{-1}) (\T_i^{-1} \otimes \T_i^{-1}) \\
&= (\T_i \otimes \T_i)(Z_i \hat{\Delta}(X) Z_i^{-1}), 
\end{align*}
and if $ X \in U_q(\mathfrak{b}_+) $ the term on the right hand side is contained in $ \T_i(U_q(\mathfrak{b}_+)) \otimes U_q(\mathfrak{g}) $ by the 
explicit formula for $ Z_i $. 
The second formula is obtained in exactly the same way. 

For the third formula it clearly suffices to show that the space on the left hand side is contained in the space on the right hand side. 
Let $ X \in U_q(\mathfrak{n}_+) \cap \T_i(U_q(\mathfrak{b}_+)) $. Moreover let $ V $ be an integrable $ U_q(\mathfrak{g}) $-module and $ v \in V $. 
If $ \lambda \in \weights^+ $ and $ v_\lambda \in V(\lambda) $ is a highest weight vector we have $ Z_i(v_\lambda \otimes v) = v_\lambda \otimes v $, 
hence the formula after Proposition \ref{Tdiagonal} yields 
$$
\T_i^{-1}(X) \cdot (v_\lambda \otimes v) = Z_i^{-1} (\T_i^{-1} \otimes \T_i^{-1})(\hat{\Delta}(X)) \cdot (v_\lambda \otimes v).  
$$
Since $ X \in U_q(\mathfrak{n}_+) $ we have 
$$
\hat{\Delta}(X) - 1 \otimes X \in \bigoplus_{\mu \in \roots^+ \setminus \{0\}} U_q(\mathfrak{n}_+)_\mu \otimes U_q(\mathfrak{b}_+),  
$$
and combining this with $ X \in \T_i(U_q(\mathfrak{b}_+)) $, the first formula of the present Lemma yields 
$$
\hat{\Delta}(X) - 1 \otimes X \in \bigoplus_{\mu \in \roots^+ \setminus \{0\}} (U_q(\mathfrak{n}_+)_\mu \cap \T_i(U_q(\mathfrak{b}_+))) \otimes U_q(\mathfrak{b}_+).
$$
Hence 
$$
(\T_i^{-1} \otimes \T_i^{-1})(\hat{\Delta}(X)) - 1 \otimes \T_i^{-1}(X) 
\in U_q(\mathfrak{h}) \biggl(\bigoplus_{\mu \in \roots^+ \setminus \{0\}} U_q(\mathfrak{n}_+)_\mu\biggr) \otimes U_q(\mathfrak{g}), 
$$
observing that the automorphism $ \T_i^{-1} $ preserves the sum of all nonzero weight spaces in $ U_q(\mathfrak{g}) $. 
We conclude 
$$
\T_i^{-1}(X) \cdot (v_\lambda \otimes v) = Z_i^{-1} (v_\lambda \otimes \T_i^{-1}(X) \cdot v) = v_\lambda \otimes \T_i^{-1}(X) \cdot v. 
$$
Since $ \T_i^{-1}(X) \in U_q(\mathfrak{b}_+) $ we can write 
$$
\T_i^{-1}(X) = \sum_{\nu \in \roots^+} X_\nu K_\nu
$$
for elements $ X_\nu \in U_q(\mathfrak{n}_+) $. Then 
\begin{align*} 
\T_i^{-1}(X) \cdot (v_\lambda \otimes v) = \sum_{\nu \in \roots^+} X_\nu \cdot (K_\nu \cdot v_\lambda \otimes K_\nu \cdot v) 
&= \sum_{\nu \in \roots^+} q^{(\nu, \lambda)} v_\lambda \otimes X_\nu K_\nu \cdot v. 
\end{align*}
On the other hand, 
\begin{align*} 
v_\lambda \otimes (\T_i^{-1}(X) \cdot v) = v_\lambda \otimes \sum_{\nu \in \roots^+} X_\nu K_\nu \cdot v, 
\end{align*}
and thus 
$$
\sum_{\nu \in \roots^+} q^{(\nu, \lambda)} v_\lambda \otimes X_\nu K_\nu \cdot v = v_\lambda \otimes \sum_{\nu \in \roots^+} X_\nu K_\nu \cdot v. 
$$
This holds for all $ v \in V $, where $ V $ is an arbitrary integrable $ U_q(\mathfrak{g}) $-module. Therefore we get
$$
\sum_{\nu \in \roots^+} (1 - q^{(\nu, \lambda)}) X_\nu K_\nu = 0 
$$
for all $ \lambda \in \weights^+ $. This implies $ X_\nu = 0 $ for $ \nu \neq 0 $, and thus $ \T_i^{-1}(X) \in U_q(\mathfrak{n}_+) $. 
In other words, we obtain $ X \in \T_i(U_q(\mathfrak{n}_+)) $ as desired. 

The final formula is obtained in a similar way. \end{proof} 

\begin{lemma} \label{uqgticoproductcontainment}
We have 
\begin{align*}
\hat{\Delta}(U_q(\mathfrak{n}_+) \cap \T_i(U_q(\mathfrak{n}_+))) &\subset (U_q(\mathfrak{n}_+) \cap \T_i(U_q(\mathfrak{n}_+))) \otimes U_q(\mathfrak{b}_+), \\
\hat{\Delta}(U_q(\mathfrak{n}_-) \cap \T_i(U_q(\mathfrak{n}_-))) &\subset U_q(\mathfrak{b}_-) \otimes (U_q(\mathfrak{n}_-) \cap \T_i(U_q(\mathfrak{n}_-))) 
\end{align*}
for $ 1 \leq i \leq N $. 
\end{lemma} 

\begin{proof} According to the first formula of Lemma \ref{drinfeldpairinghelp} we obtain 
\begin{align*}
\hat{\Delta}(U_q(\mathfrak{n}_+) \cap \T_i(U_q(\mathfrak{n}_+))) &\subset (U_q(\mathfrak{n}_+) \cap \T_i(U_q(\mathfrak{b}_+))) \otimes U_q(\mathfrak{b}_+). 
\end{align*}
Hence the third formula of Lemma \ref{drinfeldpairinghelp} yields the first claim. 
The second assertion is proved in the same way using the second and fourth formulas of Lemma \ref{drinfeldpairinghelp}. \end{proof} 

\begin{lemma} \label{uqgdrinfeldorthogonalityrelations}
We have 
\begin{align*}
U_q(\mathfrak{n}_+) \cap \T_i(U_q(\mathfrak{n}_+)) &= \{X \in U_q(\mathfrak{n}_+) \mid \tau(X, U_q(\mathfrak{n}_-) F_i) = 0 \}, \\
U_q(\mathfrak{n}_-) \cap \T_i(U_q(\mathfrak{n}_-)) &= \{Y \in U_q(\mathfrak{n}_-) \mid \tau(U_q(\mathfrak{n}_+) E_i, Y) = 0 \}
\end{align*}
for $ 1 \leq i \leq N $. 
\end{lemma} 

\begin{proof} Assume first $ X \in U_q(\mathfrak{n}_+) \cap \T_i(U_q(\mathfrak{n}_+)) $. 
Using the skew-pairing property we have
$$
\tau(X, Y F_i) = \tau(X_{(1)}, F_i) \tau(X_{(2)}, Y) 
$$ 
for any $ Y \in U_q(\mathfrak{n}_-) $. 
According to the first part of Lemma \ref{uqgticoproductcontainment},
the first tensor factors of $ \hat{\Delta}(X) $ consist of sums of terms $ Z_\gamma \in U_q(\mathfrak{n}_+) \cap \T_i(U_q(\mathfrak{n}_+)) $ 
of weight $ \gamma $ for some $ \gamma \in (\roots^+ \cap s_i \roots^+) \setminus \{0\} $.
Since  $ \gamma \neq \alpha_i $, we must have $ \rho(Z_\gamma, F_i) = 0 $ for any such $Z_\gamma$, and thus $\tau(X,YF_i)=0$.

Assume conversely that $ X \in U_q(\mathfrak{n}_+) $ satisfies $ \tau(X, U_q(\mathfrak{n}_-) F_i) = 0 $, and let us show $ X \in \T_i(U_q(\mathfrak{n}_+)) $. 
According to Lemma \ref{drinfeldpairinghelp} it suffices to prove
$ \T_i^{-1}(X) \in U_q(\mathfrak{b}_+) $, and for this it is enough to verify 
$$
\T_i^{-1}(X)(v_\lambda \otimes V) \subset v_\lambda \otimes V  
$$
for all integrable modules $ V $ and $ v_\lambda \in V(\lambda) $ for $ \lambda \in \weights^+ $. Let us define $ X_r \in U_q(\mathfrak{b}_+) $ 
for $ r \in \mathbb{N} $ by 
$$
\hat{\Delta}(X) - \sum_{r \in \mathbb{N}} E_i^r \otimes X_r K_i^r 
\in \biggl(\bigoplus_{\mu \in \roots^+ \setminus \mathbb{N} \alpha_i} U_q(\mathfrak{n}_+)_\mu\biggr) \otimes U_q(\mathfrak{b}_+). 
$$
For $ Y \in U_q(\mathfrak{n}_-) $ and $ m \in \mathbb{N} $ we get 
\begin{align*}
0 = \tau(X, Y F_i^m) &= \tau(X_{(1)}, F_i^m) \tau(X_{(2)}, Y) \\ 
&= \sum_{r \in \mathbb{N}} \tau(E_i^r, F_i^m) \tau(X_r K_i^r, Y) = \sum_{r \in \mathbb{N}} \tau(E_i^r, F_i^m) \tau(X_r, Y).  
\end{align*}
Since $ \tau(E_i^r, F_i^m) = 0 $ for $ m \neq r $ we deduce $ X_m = 0 $ for all $ m \in \mathbb{N} $, using the nondegeneracy of $ \tau $ 
obtained in Proposition \ref{lemdrinfeldnondegenerate}. Hence 
$$
\hat{\Delta}(X) 
\in \biggl(\bigoplus_{\mu \in \roots^+ \setminus \mathbb{N} \alpha_i} U_q(\mathfrak{n}_+)_\mu\biggr) \otimes U_q(\mathfrak{b}_+). 
$$
On the other hand, we have $ U_q(\mathfrak{n}_+)_\gamma F_i^m \subset \sum_{r = 0}^m U_q(\mathfrak{h}) F_i^r U_q(\mathfrak{n}_+)_{\gamma - r\alpha_i} $. 
If $ \gamma \in \roots^+ \setminus \mathbb{N}_0 \alpha_i $ we thus get $ U_q(\mathfrak{n}_+)_\gamma F_i^m \cdot v_{\lambda} = 0 $. 

We now compute 
\begin{align*}
\T_i^{-1}(X) \cdot (v_\lambda \otimes V) &= Z_i^{-1} \cdot (\T_i^{-1} \otimes \T_i^{-1})(\hat{\Delta}(X) \cdot (\T_i(v_\lambda) \otimes V)) \\
&= Z_i^{-1} \cdot (\T_i^{-1} \otimes \T_i^{-1})(\hat{\Delta}(X) \cdot (F_i^{(\alpha_i^\vee, \lambda)} \cdot v_\lambda \otimes V)) 
\end{align*}
using the formula for the coproduct of $ \T_i^{-1}(X) $ obtained after Proposition \ref{Tdiagonal} 
and the definition of the action of $ \T_i $ on a highest weight vector. 
From our above considerations we get 
\begin{align*}
\hat{\Delta}(X) \cdot (F_i^{(\alpha_i^\vee, \lambda)} \cdot v_\lambda \otimes V) &= (1 \otimes X) \cdot (F_i^{(\alpha_i^\vee, \lambda)} \cdot v_\lambda \otimes V) 
\subset F_i^{(\alpha_i^\vee, \lambda)} \cdot v_\lambda \otimes V. 
\end{align*}
Hence 
$$
\T_i^{-1}(X) \cdot (v_\lambda \otimes V) \subset Z_i^{-1} \cdot (\T_i^{-1}(F_i^{(\alpha_i^\vee, \lambda)} \cdot v_\lambda) \otimes V) 
= Z_i^{-1} \cdot (v_\lambda \otimes V) = v_\lambda \otimes V 
$$
as desired. \\
The second assertion is proved in a similar way. \end{proof} 

\begin{lemma} \label{uqgdecompositionsimpleroots}
The multiplication of $ U_q(\mathfrak{g}) $ induces linear isomorphisms 
\begin{align*}
U_q(\mathfrak{n}_+) &\cong \bigoplus_{r = 0}^\infty \mathbb{K} E_i^r \otimes (U_q(\mathfrak{n}_+) \cap \T_i(U_q(\mathfrak{n}_+))) \\
U_q(\mathfrak{n}_-) &\cong \bigoplus_{r = 0}^\infty \mathbb{K} F_i^r \otimes (U_q(\mathfrak{n}_-) \cap \T_i(U_q(\mathfrak{n}_-))) 
\end{align*}
for each $ 1 \leq i \leq N $. 
\end{lemma}

\begin{proof} Let $ w_0 = s_{i_1} \cdots s_{i_t} $ be a reduced expression of the longest word of $ W $ beginning with $ i_1 = i $. Then 
according to Theorem \ref{PBWgeneralfield} the vectors $ E_i^{a_1} E_{\beta_2}^{a_2} \cdots E_{\beta_n}^{a_n} $ form a basis of 
$ U_q(\mathfrak{n}_+) $. By construction, the linear span of $ E_{\beta_2}^{a_2} \cdots E_{\beta_n}^{a_n} $ is contained 
in $ U_q(\mathfrak{n}_+) \cap \T_i(U_q(\mathfrak{n}_+)) $. Moreover, the multiplication map 
$ \bigoplus_{r = 0}^\infty \mathbb{K} E_i^r \otimes (U_q(\mathfrak{n}_+) \cap \T_i(U_q(\mathfrak{n}_+))) \rightarrow U_q(\mathfrak{n}_+) $ is 
injective because $ \T_i^{-1}(E_i) \in U_q(\mathfrak{b}_-) $ by the definition of $ \T_i $ 
and $ \T_i^{-1}(U_q(\mathfrak{n}_+) \cap \T_i(U_q(\mathfrak{n}_+))) \subset U_q(\mathfrak{n}_+) $. 
This yields the first claim. 

The second isomorphism is obtained in a similar way. \end{proof} 

Using the notation $U_q(\mathfrak{n}_\pm)[w]$ which we introduced before Definition \ref{def:quantum_root_vectors}, we may rephrase the assertion of Lemma \ref{uqgdecompositionsimpleroots} as 
\begin{align*}
U_q(\mathfrak{n}_+) \cap \T_i(U_q(\mathfrak{n}_+)) &= \T_i(U_q(\mathfrak{n}_+)[s_i w_0]), \\
U_q(\mathfrak{n}_-) \cap \T_i(U_q(\mathfrak{n}_-)) &= \T_i(U_q(\mathfrak{n}_-)[s_i w_0]) 
\end{align*} 
for all $ 1 \leq i \leq N $, where $ w_0 \in W $ is the longest element of $ W $.

\begin{definition}
 \label{def:HC_map}
 The Harish-Chandra map is the linear map $ \P: U_q(\mathfrak{g}) \rightarrow U_q(\mathfrak{h}) $ 
 \nomenclature[o$P$2]{$\P$}{Harish-Chandra map}%
 given by $ \hat{\epsilon} \otimes \id \otimes \hat{\epsilon} $ under the triangular isomorphism 
 $ U_q(\mathfrak{g}) \cong U_q(\mathfrak{n}_-) \otimes U_q(\mathfrak{h}) \otimes U_q(\mathfrak{n}_+) $. 
\end{definition}

\begin{prop} \label{DrinfeldpairingHCformula}
Let $ \gamma \in \roots^+ $ and assume $ X \in U_q(\mathfrak{n}_+)_\gamma, Y \in U_q(\mathfrak{n}_-)_{-\gamma} $. Then 
$$
\P(XY) - \tau(X,Y) K_{-\gamma} \in \sum_{\nu \in \roots^+ \setminus \{0\}} \mathbb{K} K_{2 \nu - \gamma}. 
$$
\end{prop} 

\begin{proof} Let us write 
$$
\hat{\Delta}(X) = \sum_r X^{[1]}_r \otimes X^{[2]}_r K_{\gamma_r}, \qquad \hat{\Delta}(Y) = \sum_s K_{-\delta_s} Y^{[1]}_s \otimes Y^{[2]}_s
$$
where $ \gamma_r, \delta_s \in \roots^+ $ and $ X^{[1]}_r \in U_q(\mathfrak{n}_+)_{\gamma_r}, X^{[2]}_r \in U_q(\mathfrak{n}_+)_{\gamma - \gamma_r} $, 
and similarly $ Y^{[1]}_s \in U_q(\mathfrak{n}_-)_{-\gamma + \delta_s}, Y^{[2]}_s \in U_q(\mathfrak{n}_-)_{-\delta_s} $. From 
the structure of the comultiplication we see that 
\begin{align*}
(\id \otimes \hat{\Delta})\hat{\Delta}(X)& - \sum_r X^{[1]}_r \otimes K_{\gamma_r} \otimes X^{[2]}_r K_{\gamma_r} 
\end{align*}
is contained in $ U_q(\mathfrak{n}_+) \otimes U_q(\mathfrak{h}) (U_q(\mathfrak{n}_+) \cap \ker(\hat{\epsilon})) \otimes U_q(\mathfrak{b}_+) $. 
Similarly, 
\begin{align*}
(\hat{\Delta} \otimes \id)\hat{\Delta}(Y)& - \sum_s K_{-\delta_s} Y^{[1]}_s \otimes K_{-\delta_s} \otimes Y^{[2]}_s 
\end{align*}
is contained in $ U_q(\mathfrak{b}_-) \otimes (U_q(\mathfrak{n}_+) \cap \ker(\hat{\epsilon})) U_q(\mathfrak{h}) \otimes U_q(\mathfrak{n}_-) $.  
Hence according to the commutation relations in $ U_q(\mathfrak{g}) $ from Lemma \ref{lemdrinfeldformcommutation} 
and the definition of $ \P $ we obtain 
\begin{align*}
\P(XY) &= \tau(X_{(1)}, \hat{S}^{-1}(Y_{(1)})) \P(Y_{(2)} X_{(2)}) \tau(X_{(3)}, Y_{(3)}) \\
&= \sum_{r,s} \tau(X_r^{[1]}, \hat{S}^{-1}(K_{-\delta_s} Y^{[1]}_s)) K_{\gamma_r - \delta_s} \tau(X^{[2]}_r K_{\gamma_r}, Y_s^{[2]}) \\
&= \sum_{r,s} \tau(X_r^{[1]}, \hat{S}^{-1}(K_{-\delta_s} Y^{[1]}_s)) \tau(X^{[2]}_r K_{\gamma_r}, Y_s^{[2]}) K_{2 \gamma_r - \gamma};
\end{align*}
here we note that only summands $ r, s $ such that $ \gamma_r + \delta_s = \gamma $ contribute. 
The term for $ \gamma_r = 0 $ is 
$$
\sum_s \tau(1, K_{-\delta_s} Y^{[1]}_s) \tau(X, Y_s^{[2]}) = \tau(X, Y). 
$$
This yields the claim. \end{proof} 

We prove now that the Drinfeld pairing is invariant in a suitable sense. 

\begin{theorem} \label{drinfeldpairinginvariance}
Let $ X \in U_q(\mathfrak{n}_+) \cap \T_i(U_q(\mathfrak{n}_+)) $ and $ Y \in U_q(\mathfrak{n}_-) \cap \T_i(U_q(\mathfrak{n}_-)) $. Then 
$$
\tau(\T_i^{-1}(X), \T_i^{-1}(Y)) = \tau(X, Y).  
$$
\end{theorem}

\begin{proof} We may assume $ X \in U_q(\mathfrak{n}_+)_\gamma \cap \T_i(U_q(\mathfrak{n}_+)) $, 
$ Y \in U_q(\mathfrak{n}_-)_{-\gamma} \cap \T_i(U_q(\mathfrak{n}_-)) $ for some $ \gamma \in \roots^+ $. 
Applying Proposition \ref{DrinfeldpairingHCformula} to $ \T_i^{-1}(X), \T_i^{-1}(Y) $, we see that it suffices to show 
$$
\P(\T_i^{-1}(X Y)) - \tau(X, Y) K_{-s_i \gamma} \in \sum_{\nu \in \roots^+ \setminus \{0\}} \mathbb{K} K_{2\nu - s_i \gamma}, 
$$
since we then obtain the claim by comparing the coefficients of $ K_{-s_i \gamma} $. 

As in the proof of Proposition \ref{DrinfeldpairingHCformula} we may write 
$$
\hat{\Delta}(X) = \sum_r X^{[1]}_r \otimes X^{[2]}_r K_{\gamma_r}, \qquad \hat{\Delta}(Y) = \sum_s K_{-\delta_s} Y^{[1]}_s \otimes Y^{[2]}_s
$$
where $ \gamma_r, \delta_s \in \roots^+ \cap s_i \roots^+ $ and $ X^{[1]}_r \in U_q(\mathfrak{n}_+)_{\gamma_r}, Y^{[2]}_s \in U_q(\mathfrak{n}_-)_{-\delta_s} $.
According to Lemma \ref{uqgticoproductcontainment} we have in fact $ X^{[1]}_r \in U_q(\mathfrak{n}_+)_{\gamma_r} \cap \T_i(U_q(\mathfrak{n}_+)), 
X^{[2]}_r \in U_q(\mathfrak{n}_+)_{\gamma - \gamma_r} $, and similarly $ Y^{[1]}_s \in U_q(\mathfrak{n}_-)_{-\gamma + \delta_s}, 
Y^{[2]}_s \in U_q(\mathfrak{n}_-)_{-\delta_s} \cap \T_i(U_q(\mathfrak{n}_-)) $. 
Using again Lemma \ref{uqgticoproductcontainment} and Lemma \ref{uqgdecompositionsimpleroots} 
we find elements $ X^{[1]}_{rm} \in U_q(\mathfrak{n}_+)_{\gamma_r - m \alpha_i} \cap \T_i(U_q(\mathfrak{n}_+)) $ such that 
\begin{align*} 
\hat{\Delta}(X^{[1]}_r) - \sum_{m \in \mathbb{N}_0} X^{[1]}_{rm} \otimes E_i^{(m)} K_{\gamma_r - m \alpha_i} 
\end{align*} 
is contained in $ U_q(\mathfrak{n}_+) \otimes U_q(\mathfrak{b}_+)(U_q(\mathfrak{n}_+) \cap \T_i(U_q(\mathfrak{n}_+)) \cap \ker(\hat{\epsilon})) $. 
Then
\begin{align*} 
&(\hat{\Delta} \otimes \id) \hat{\Delta}(X) 
- \sum_r \sum_{m \in \mathbb{N}_0} X^{[1]}_{rm} \otimes E_i^{(m)} K_{\gamma_r - m \alpha_i} \otimes X^{[2]}_r K_{\gamma_r} 
\end{align*}
is contained in 
$ U_q(\mathfrak{n}_+) \otimes U_q(\mathfrak{b}_+)(U_q(\mathfrak{n}_+) \cap \T_i(U_q(\mathfrak{n}_+)) \cap \ker(\hat{\epsilon})) \otimes U_q(\mathfrak{b}_+) $. 

Similarly, we find $ Y^{[2]}_{sm} \in U_q(\mathfrak{n}_-)_{-\delta_s + m \alpha_i} \cap \T_i(U_q(\mathfrak{n}_-)) $ such that 
\begin{align*}
\hat{\Delta}(Y^{[2]}_s) - \sum_{m \in \mathbb{N}_0} K_{-\delta_s + m \alpha_i} F_i^{(m)} \otimes Y^{[2]}_{sm} 
\end{align*} 
is contained in $ (U_q(\mathfrak{n}_-) \cap \T_i(U_q(\mathfrak{n}_-)) \cap \ker(\hat{\epsilon})) U_q(\mathfrak{b}_-) \otimes U_q(\mathfrak{n}_-) $, 
and then 
\begin{align*} 
(\id \otimes \hat{\Delta}) \hat{\Delta}(Y) 
- \sum_s \sum_{m \in \mathbb{N}_0} K_{-\delta_s} Y^{[1]}_s \otimes K_{-\delta_s + m \alpha_i} F_i^{(m)} \otimes Y^{[2]}_{sm} 
\end{align*} 
is contained in 
$ U_q(\mathfrak{b}_-) \otimes (U_q(\mathfrak{n}_-) \cap \T_i(U_q(\mathfrak{n}_-)) \cap \ker(\hat{\epsilon}))U_q(\mathfrak{b}_-) \otimes U_q(\mathfrak{n}_-) $. 

Due to the commutation relations from Lemma \ref{lemdrinfeldformcommutation} and the invariance properties of $ \tau $ we get 
\begin{align*} 
&XY - \sum_{m \in \mathbb{N}_0} \sum_{\gamma_r + \delta_s = \gamma + m \alpha_i} 
\tau(X^{[1]}_{rm}, \hat{S}^{-1}(Y^{[1]}_s)) K_{-\delta_s + m\alpha_i} F_i^{(m)} E_i^{(m)} K_{\gamma_r - m \alpha_i} 
\tau(X^{[2]}_r, Y^{[2]}_{sm}) \\
&\in (U_q(\mathfrak{n}_-) \cap \T_i(U_q(\mathfrak{n}_-)) \cap \ker(\hat{\epsilon}))U_q(\mathfrak{g}) 
+ U_q(\mathfrak{g})(U_q(\mathfrak{n}_+) \cap \T_i(U_q(\mathfrak{n}_+)) \cap \ker(\hat{\epsilon})).  
\end{align*} 
In the sequel, this relation will be evaluated in two ways. 

Firstly, by the definition of $ \P $ we obtain 
\begin{align*} 
\P(XY) &= \sum_{\gamma_r + \delta_s = \gamma} \tau(X^{[1]}_{r0}, \hat{S}^{-1}(Y^{[1]}_s)) K_{-\gamma + 2\gamma_r} \tau(X^{[2]}_r, Y^{[2]}_{s0}) . 
\end{align*} 
Therefore 
$$ 
\tau(X, Y) = \sum_{\gamma_r = 0, \delta_s = \gamma} \tau(X^{[1]}_{r0}, \hat{S}^{-1}(Y^{[1]}_s)) \tau(X^{[2]}_r, Y^{[2]}_{s0}) 
$$
according to Proposition \ref{DrinfeldpairingHCformula}. 

Secondly, applying $ \T_i^{-1} $ to the above formula gives 
\begin{align*}
\T_i^{-1}(XY) &- \sum_{m \in \mathbb{N}_0} \sum_{\gamma_r + \delta_s = \gamma + m \alpha_i} 
\tau(X^{[1]}_{rm}, \hat{S}^{-1}(Y^{[1]}_s)) \tau(X^{[2]}_r, Y^{[2]}_{sm}) \times \\
&\qquad K_{s_i(-\delta_s + m \alpha_i)} E_i^{(m)} F_i^{(m)} K_{s_i(\gamma_r - m \alpha_i)} \\
&\in (U_q(\mathfrak{n}_-) \cap \ker(\hat{\epsilon}))U_q(\mathfrak{g}) 
+ U_q(\mathfrak{g})(U_q(\mathfrak{n}_+) \cap \ker(\hat{\epsilon})), 
\end{align*}
taking into account
\begin{align*}
\T_i^{-1}(F_i^{(m)} E_i^{(m)}) &= \frac{1}{[m]_q!^2}(-K_i E_i)^m(-F_i K_i^{-1})^m = \frac{1}{[m]_q!^2} E_i^m F_i^m = E_i^{(m)} F_i^{(m)}
\end{align*}
by the definition of $ \T_i $. According to the commutation relations in Proposition \ref{EFpowercommutation} we obtain 
\begin{align*}
E_i^{(m)} F_i^{(m)} - 
\begin{bmatrix}
K_i; 0 \\
m
\end{bmatrix}_{q_i}
\in 
U_q(\mathfrak{g})(U_q(\mathfrak{n}_+) \cap \ker(\hat{\epsilon})) + (U_q(\mathfrak{n}_-) \cap \ker(\hat{\epsilon}))U_q(\mathfrak{g})
\end{align*}
where we recall that
\begin{align*}
\begin{bmatrix}
K_i; 0\\
m
\end{bmatrix}_{q_i}
= \prod_{j = 1}^m \frac{q_i^{1 - j} K_i - q_i^{-(1 - j)} K_i^{-1}}{q_i^j - q_i^{-j}}. 
\end{align*}
Hence, 
\begin{align*}
\T_i^{-1}(XY) &- \sum_{m \in \mathbb{N}_0} \sum_{\gamma_r + \delta_s = \gamma + m \alpha_i} 
\tau(X^{[1]}_{rm}, \hat{S}^{-1}(Y^{[1]}_s)) \tau(X^{[2]}_r, Y^{[2]}_{sm}) 
\begin{bmatrix}
K_i; 0 \\
m
\end{bmatrix}_{q_i}
K_{s_i(\gamma_r - \delta_s)} \\
&\in (U_q(\mathfrak{n}_-) \cap \ker(\hat{\epsilon}))U_q(\mathfrak{g}) 
+ U_q(\mathfrak{g})(U_q(\mathfrak{n}_+) \cap \ker(\hat{\epsilon})), 
\end{align*}
and therefore 
\begin{align*}
\P(\T_i^{-1}(XY)) &= \sum_{m \in \mathbb{N}_0} \sum_{\gamma_r + \delta_s = \gamma + m \alpha_i} 
\tau(X^{[1]}_{rm}, \hat{S}^{-1}(Y^{[1]}_s)) \tau(X^{[2]}_r, Y^{[2]}_{sm}) 
\begin{bmatrix}
K_i; 0 \\
m
\end{bmatrix}_{q_i}
K_{s_i(\gamma_r - \delta_s)}. 
\end{align*}
Observe that 
\begin{align*}
\begin{bmatrix}
K_i; 0 \\
m
\end{bmatrix}_{q_i}
\in K_{-m \alpha_i} \biggl(\mathbb{K}^\times + \sum_{l \in \mathbb{N}} \mathbb{K} K_{2l \alpha_i}\biggr). 
\end{align*}
We may therefore write $ \P(\T_i^{-1}(XY)) $ as a linear combination of elements of the form $ K_{s_i(\gamma_r - \delta_s) - m\alpha_i} K_{2l \alpha_i} $ 
where $ l \in \mathbb{N}_0 $. 

Let us recall that $ \gamma_r \in \roots^+ \cap s_i \roots^+ $ by assumption, which implies $ s_i \gamma_r \in \roots^+ $. 
Now $ \gamma_r + \delta_s = \gamma + m \alpha_i $ yields
$$
\gamma_r - \delta_s + m \alpha_i = -\gamma + 2 \gamma_r, 
$$
and thus
$$
s_i(\gamma_r - \delta_s) - m \alpha_i 
= -s_i \gamma + 2s_i \gamma_r. 
$$
Recall also that $ X^{[1]}_{rm} \in U_q(\mathfrak{n}_+)_{\gamma_r - m \alpha_i} $. Hence if $ \gamma_r = 0 $ then $ X^{[1]}_{rm} = 0 $ 
unless $ m = 0 $. Combining these considerations we obtain 
\begin{align*}
\P(\T_i^{-1}(XY)) &\in K_{-s_i \gamma}\biggl( \sum_{\gamma_r = 0, \delta_s = \gamma} 
\tau(X^{[1]}_{r0}, \hat{S}^{-1}(Y^{[1]}_s)) \tau(X^{[2]}_r, Y^{[2]}_{s0})  + \sum_{\nu \in \roots^+ \setminus \{0\}} \mathbb{K} K_{2\nu} \biggr) \\
&= K_{-s_i \gamma}\biggl(\tau(X, Y) + \sum_{\nu \in \roots^+ \setminus \{0\}} \mathbb{K} K_{2\nu} \biggr), 
\end{align*}
using our previous formula for $ \tau(X, Y) $. This finishes the proof. \end{proof} 

\begin{lemma} \label{drinfeldpairingspecialmult}
If $ X \in U_q(\mathfrak{n}_+)[s_i w_0], Y \in U_q(\mathfrak{n}_-)[s_i w_0] $ and $ r, s \in \mathbb{N}_0 $, then 
\begin{align*}
\tau(\T_i(X) E_i^r, \T_i(Y) F_i^s) &= \delta_{rs} \tau(\T_i(X), \T_i(Y)) \tau(E_i^r, F_i^r) 
\end{align*} 
for any $ 1 \leq i \leq N $. 
\end{lemma}

\begin{proof} Let us first show 
$$
\tau(E_i^r, F_i^r) 
 = q_i^{r(r - 1)/2} [r]_{q_i}!\;  \tau(E_i, F_i)^r
 = \left(\prod_{k = 1}^r \frac{q_i^{2k} - 1}{q_i^2 - 1}\right) \tau(E_i, F_i)^r 
$$
for $ r \in \mathbb{N} $. For $ r = 1 $ this claim is trivially true. Assume that it holds for $ r $, then using Lemma \ref{lemEiFicomult} we obtain 
\begin{align*}
\tau(E_i^{r + 1}, F_i^{r + 1}) &= \tau(E_i^r, (F_i^{r + 1})_{(1)}) \tau(E_i, (F_i^{r + 1})_{(2)}) \\
&= q_i^{r} [r + 1]_{q_i} \tau(E_i^r, F_i^r K_i^{-1}) \tau(E_i, F_i) \\
&= q_i^{(r+1)r/2} [r+1]_{q_i}!\;  \tau(E_i, F_i)^{r+1}.
\end{align*}
This proves the first equality, and the second is elementary.

In order to prove the Lemma we use induction on $ r + s $. For $ r + s = 0 $, that is, $ r = s = 0 $, the 
claim trivially holds. 

Assume now $ r + s = k > 0 $, and that the claim holds for $ k - 1 $. 
Recall from the remark after Lemma \ref{uqgdecompositionsimpleroots} that 
$ \T_i(U_q(\mathfrak{n}_+)[s_i w_0]) = U_q(\mathfrak{n}_+) \cap \T_i(U_q(\mathfrak{n}_+)) $ 
and $ \T_i(U_q(\mathfrak{n}_-)[s_i w_0]) = U_q(\mathfrak{n}_-) \cap \T_i(U_q(\mathfrak{n}_-)) $, 
so that $ \T_i(X) \in U_q(\mathfrak{n}_+) \cap \T_i(U_q(\mathfrak{n}_+)) $ and $ \T_i(Y) \in U_q(\mathfrak{n}_-) \cap \T_i(U_q(\mathfrak{n}_-)) $. 
If $ r = 0 $ then Lemma \ref{uqgdrinfeldorthogonalityrelations} yields $ \tau(\T_i(X), \T_i(Y) F_i^s) = 0 $, 
and similarly we get $ \tau(\T_i(X) E_i^r, \T_i(Y)) = 0 $ if $ s = 0 $. 

Let us therefore assume that both $ r, s $ are positive, and let us consider the case $ s \geq r $. 
Lemma \ref{lemEiFicomult} yields 
\begin{align*}
&\tau(\T_i(X) E_i^r, \T_i(Y) F_i^s) = \tau(\T_i(X)_{(2)} (E_i^r)_{(2)}, \T_i(Y) F_i^{s - 1}) \tau(\T_i(X)_{(1)} (E_i^r)_{(1)}, F_i) \\
&= \sum_{l = 0}^r q_i^{l(r - l)} 
\begin{bmatrix}
r \\
l
\end{bmatrix}_{q_i}
\tau(\T_i(X)_{(2)} E_i^{r - l} K_i^l, \T_i(Y) F_i^{s - 1}) \tau(\T_i(X)_{(1)} E_i^l, F_i) .
\end{align*}
According to Lemma \ref{uqgticoproductcontainment} we know $ \hat{\Delta}(\T_i(X)) \in (\T_i(U_q(\mathfrak{n}_+)) \cap U_q(\mathfrak{n}_+)) \otimes U_q(\mathfrak{b}_+) $, and so the above terms are all zero except when $l=1$.  Therefore, 
\begin{align*}
&\tau(\T_i(X) E_i^r, \T_i(Y) F_i^s) \\
&= q_i^{r - 1} [r]_{q_i} \tau(\T_i(X)_{(2)} E_i^{r - 1} K_i, \T_i(Y) F_i^{s - 1}) \tau(\T_i(X)_{(1)} E_i, F_i) \\
&= q_i^{r - 1} [r]_{q_i} \tau(\T_i(X)_{(2)} E_i^{r - 1} K_i, \T_i(Y) F_i^{s - 1}) \times \\
&\qquad (\tau(\T_i(X)_{(1)}, F_i) \tau(E_i, 1) + \tau(\T_i(X)_{(1)}, K_i^{-1}) \tau(E_i, F_i)) \\
&= q_i^{r - 1} [r]_{q_i} \tau(\T_i(X)_{(2)} E_i^{r - 1} K_i, \T_i(Y) F_i^{s - 1}) 
\tau(\T_i(X)_{(1)}, 1) \tau(E_i, F_i) \\
&= q_i^{r - 1} [r]_{q_i} \tau(\T_i(X) E_i^{r - 1} K_i, \T_i(Y) F_i^{s - 1}) \tau(E_i, F_i) \\
&= q_i^{r - 1} [r]_{q_i} \tau(\T_i(X) E_i^{r - 1}, \T_i(Y) F_i^{s - 1}) \tau(E_i, F_i). 
\end{align*} 
Now the inductive hypothesis and our considerations at the beginning imply
\begin{align*}
q_i^{r - 1} &[r]_{q_i} \tau(\T_i(X) E_i^{r - 1}, \T_i(Y) F_i^{s - 1}) \tau(E_i, F_i) \\
&= \delta_{r,s} q_i^{r - 1} [r]_{q_i} \tau(\T_i(X), \T_i(Y)) \tau(E_i^{r - 1}, F_i^{s - 1}) \tau(E_i, F_i) \\
&= \delta_{r,s} q_i^{r - 1} [r]_{q_i} q_i^{(r - 1)(r - 2)/2} [r - 1]_{q_i}! \tau(\T_i(X), \T_i(Y)) \tau(E_i, F_i)^{r - 1} \tau(E_i, F_i) \\
&= \delta_{r,s} \tau(\T_i(X), \T_i(Y)) \tau(E_i^r, F_i^r). 
\end{align*}
This yields the claim. \end{proof} 

\begin{theorem} \label{PBWdrinfeldorthogonal}
The PBW-basis vectors are orthogonal with respect to the Drinfeld pairing. More precisely, we have 
$$
\tau(E^{a_n}_{\beta_n} \cdots E^{a_1}_{\beta_1}, F^{b_n}_{\beta_n} \cdots F^{b_1}_{\beta_1}) = \prod_{k = 1}^n
\delta_{a_k, b_k} \tau(E_{\beta_k}^{a_k}, F_{\beta_k}^{b_k}). 
$$
\end{theorem} 

\begin{proof} We shall prove more generally that the PBW-vectors associated to any reduced expression 
$ w = s_{i_1} \cdots s_{i_t} $ of $ w \in W $ satisfy the above orthogonality relations. This clearly implies 
the desired statement. 

We use induction on the length of $ w $. For $ l(w) = 0 $ the claim is trivial, and for $ l(w) = 1 $ it follows 
immediately from the construction of the pairing. Assume that $ l(w) > 1 $ and write 
\begin{align*}
X &= (\T_{i_2} \cdots \T_{i_{t - 1}})(E_{i_t}^{a_t}) \cdots \T_{i_2}(E_{i_3}^{a_3}) E_{i_2}^{a_2} \\
Y &= (\T_{i_2} \cdots \T_{i_{n - 1}})(F_{i_t}^{s_t}) \cdots \T_{i_2}(F_{i_3}^{s_3}) F_{i_2}^{s_2}.  
\end{align*} 
Since $ \T_{i_1}(X) \in U_q(\mathfrak{n}_+) $ the remark after Lemma \ref{uqgdecompositionsimpleroots}
shows that $ X \in U_q(\mathfrak{n}_+)[s_{i_1} w_0] $ where $ w_0 $ is the longest element of $ W $. 
Similarly we obtain $ Y \in U_q(\mathfrak{n}_-)[s_{i_1} w_0] $. 

Hence by Lemma \ref{drinfeldpairingspecialmult} and Theorem \ref{drinfeldpairinginvariance} we calculate
\begin{align*}
\tau(\T_{i_1}(X) E_{i_1}^{a_1}, \T_{i_1}(Y) F_{i_1}^{b_1}) &= \tau(\T_{i_1}(X), \T_{i_1}(Y)) \tau(E_{i_1}^{a_1}, F_{i_1}^{b_1}) \\
&= \tau(X, Y) \tau(E_{i_1}^{a_1}, F_{i_1}^{b_1}) = \delta_{a_1, b_1} \tau(E_{i_1}^{a_1}, F_{i_1}^{b_1}) \tau(X, Y). 
\end{align*} 
Now we can apply the induction hypothesis to $ s_{i_1} w = s_{i_2} \cdots s_{i_t} $. This yields the claim. \end{proof} 

To conclude this subsection we complete the description obtained in Theorem \ref{PBWdrinfeldorthogonal} with the following formula.  
We use the notation
\[
 q_\beta = q^{(\beta,\beta)/2}
\]
\label{nom:q_beta}%
for any root $\beta\in\bf\Delta$.

\begin{lemma} \label{drinfeldpairingevaluate}
For any $ 1 \leq k \leq n $ and $ r \geq 0 $ we have 
$$
\tau(E_{\beta_k}^r, F_{\beta_k}^r) = (-1)^r q_{\beta_k}^{r (r - 1)/2} \frac{[r]_{q_{\beta_k}}!}{(q_{\beta_k} - q_{\beta_k}^{-1})^r},
$$
and therefore
\[
\tau(E^{a_n}_{\beta_n} \cdots E^{a_1}_{\beta_1}, F^{b_n}_{\beta_n} \cdots F^{b_1}_{\beta_1}) = \prod_{k = 1}^n
\delta_{a_k, b_k} 
(-1)^{a_k} q_{\beta_k}^{a_k (a_k - 1)/2} 
\frac{[a_k]_{q_{\beta_k}}!}{(q_{\beta_k} - q_{\beta_k}^{-1})^{a_k}}. 
\]
\end{lemma} 

\begin{proof} According to Theorem \ref{drinfeldpairinginvariance} it suffices to prove the first assertion in the case that $ \beta_k = \alpha_i $ is a simple root. 
In this case, the computation in the proof of Lemma \ref{drinfeldpairingspecialmult} combined with Lemma \ref{tauskewpairingproperties} yields
$$
\tau(E_i^r, F_i^r) = q_i^{r(r - 1)/2} [r]_{q_i}! \; \tau(E_i, F_i)^r = (-1)^r q_i^{r(r - 1)/2} [r]_{q_i}! \; \frac{1}{(q_i - q_i^{-1})^r}
$$ 
as desired. \end{proof}

\subsubsection{Nondegeneracy of the Drinfeld pairing} 

The computations of the previous subsection are key to proving nondegeneracy of the Drinfeld pairing and the quantum Killing form.  
We keep the assumptions and notation from above.

\begin{theorem} \label{lemdrinfeldnondegenerate}
The Drinfeld pairing $ \tau: U_q(\mathfrak{b}_+) \times U_q(\mathfrak{b}_-) \rightarrow \mathbb{K} $ and the quantum Killing 
form $ \kappa: U_q(\mathfrak{g}) \times U_q(\mathfrak{g}) \rightarrow \mathbb{K} $ are nondegenerate. 
\end{theorem} 

\begin{proof} For the first claim it suffices to show that the restriction of $ \tau $ 
to $ U_q(\mathfrak{b}_+)_\beta \times U_q(\mathfrak{b}_-)_{-\beta} $ is nondegenerate. 
From Lemma \ref{tauskewpairingproperties} we know that $ \tau $ is diagonal with respect to the tensor 
product decomposition $ U_q(\mathfrak{b}_\pm)_{\pm \beta} \cong U_q(\mathfrak{n}_\pm)_{\pm \beta} \otimes U_q(\mathfrak{h}) $. 

Since $ q $ is not a root of unity the characters $ \chi_\mu $ of $ U_q(\mathfrak{h}) $ 
given by $ \chi_\mu(K_\lambda) = q^{(\mu, \lambda)} $ for $ \mu \in \weights $ are pairwise distinct, and hence linearly independent 
by Artin's Theorem. 
It follows that the restriction of $ \tau $ to $ U_q(\mathfrak{h}) $ is nondegenerate. 
Hence the claim follows from Theorem \ref{PBWdrinfeldorthogonal} combined with the PBW-Theorem \ref{PBWgeneralfield}.  

The second assertion is verified in a similar way using the triangular decomposition of $ U_q(\mathfrak{g}) $ 
and nondegeneracy of $ \tau $. \end{proof} 

We remark that Theorem \ref{lemdrinfeldnondegenerate} shows also that the Rosso form on $ U_q(\mathfrak{b}_-) \bowtie U_q(\mathfrak{b}_+) $ 
is nondegenerate. In particular, this form does not pass to the quotient $ U_q(\mathfrak{g}) $. Therefore the quantum Killing 
cannot be obtained directly from the Rosso form on the double.

\subsection{The quantum Casimir element and simple modules} 
\label{sec:Casimir}

In this section we resume our study of finite dimensional weight modules for $ U_q(\mathfrak{g}) $. 
We define a Casimir element $ \Omega $ in a certain algebraic completion of $ U_q(\mathfrak{g}) $ in order to prove 
irreducibility of the modules $ L(\mu) $ constructed in Theorem \ref{thmfdintegrable}. Throughout we assume 
that $ q =s ^L \in \mathbb{K}^\times $ is not a root of unity.

According to Theorem \ref{lemdrinfeldnondegenerate} the restriction of the Drinfeld pairing $ \tau $ 
to $ U_q(\mathfrak{n}_+)_\gamma \times U_q(\mathfrak{n}_-)_{-\gamma} $ is nondegenerate for every $ \gamma \in \roots^+ $. Therefore, we 
can choose bases consisting of $ b^\gamma_j \in U_q(\mathfrak{n}_+)_\gamma $ and $ a^\gamma_k \in U_q(\mathfrak{n}_-)_{-\gamma} $ such 
that $ \tau(b^\gamma_i, a^\gamma_j) = \delta_{ij} $. For $ \gamma \in \roots^+ $ we define 
$$
C_\gamma = \sum_j a^\gamma_j \otimes b^\gamma_j, \qquad \Omega_\gamma = \sum_j \hat{S}^{-1}(a^\gamma_j) b^\gamma_j, 
$$
in particular we obtain $ C_\gamma = 1 \otimes 1 $ and $ \Omega_\gamma = 1 $ for $ \gamma = 0 $. By definition, the quantum Casimir 
element is the formal sum 
$$
\Omega = \sum_{\gamma \in \roots^+} \Omega_\gamma, 
$$
\label{nom:Casimir}%
which is not an element of $ U_q(\mathfrak{g}) $, but rather of a suitable completion of this algebra. We shall not discuss any completions here, 
instead we will only consider the ``image'' of $ \Omega $ in $ \End(V) $ for representations $ V $ of $ U_q(\mathfrak{g}) $ 
in which for any $ v \in V $ all but finitely many terms $ \Omega_\alpha $ for $ \alpha \in \roots^+ $ act by zero. 
Note that this applies in particular to all Verma modules, and hence also all 
quotients of Verma modules. Indeed, if $ v \in M(\lambda) $ for some $ \lambda \in \mathfrak{h}^*_q $ 
then we have $ U_q(\mathfrak{n}_+)_\gamma \cdot v = 0 $ for all $ \gamma \in \roots^+ $ sufficiently large. 
Hence the Casimir operator $ \Omega $ yields a well-defined endomorphism of $ M(\lambda) $.

\begin{lemma} \label{Casimirhelplemma}
Let $ \gamma \in \roots^+ $. Then we have 
$$
(1 \otimes F_i) C_{\gamma} - C_{\gamma} (1 \otimes F_i) = C_{\gamma - \alpha_i}(F_i \otimes K_i) - (F_i \otimes K_i^{-1}) C_{\gamma - \alpha_i}
$$
for all $ 1 \leq i \leq N $, where we interpret $ C_{\gamma - \alpha_i} = 0 $ if $ \gamma - \alpha_i \notin \roots^+ $. 
\end{lemma} 

\begin{proof} For $ \gamma - \alpha_i \notin \roots^+ $ the element $ F_i $ commutes with all $ b^\gamma_j $ since the latter 
must be sums of monomials in $ E_1, \dots, E_N $ not containing any powers of $ E_i $. Hence the equality holds in this case.

Consider now the case that $ \gamma = \beta + \alpha_i $ for some $ \beta \in \roots^+ $. 
Since the pairing $ \tau $ is nondegenerate it is enough to compare the pairing of both sides against $ X \in U_q(\mathfrak{n}_+)_{\alpha_i + \beta} $ 
with respect to $ \tau $ in the first tensor factor. For the left hand side this yields 
\begin{align*}
\sum_j \tau(X, a^{\beta + \alpha_i}_j) [F_i, b^{\beta + \alpha_i}_j] &= [F_i, X] = F_i X - X F_i.  
\end{align*}
To compute the right hand side let us write $ \hat{\Delta}(X) = \sum_\gamma c_\gamma X^{[1]}_\gamma \otimes X^{[2]}_{\beta + \alpha_i - \gamma} K_\gamma $ with 
monomials $ X^{[1]}_\gamma, X^{[2]}_{\beta + \alpha_i - \gamma} $ in the generators $ F_k $ of weight $ \gamma $ and $ \beta + \alpha_i - \gamma $, respectively,  
and coefficients $ c_\gamma $. We obtain
\begin{align*}
&\sum_j \tau(X, a^{\beta}_j F_i) b^{\beta}_j K_i - \tau(X, F_i a^{\beta}_j) K_i^{-1} b^{\beta}_j \\
&= \sum_{j, \gamma} c_\gamma \tau(X^{[1]}_\gamma, F_i) \tau(X^{[2]}_{\beta + \alpha_i - \gamma} K_\gamma, a^{\beta}_j) b^{\beta}_j K_i 
- c_\gamma \tau(X^{[1]}_\gamma, a^{\beta}_j) \tau(X^{[2]}_{\beta + \alpha_i - \gamma} K_\gamma, F_i) K_i^{-1} b^{\beta}_j \\
&= \sum_j - c_{\alpha_i} \frac{1}{q_i - q_i^{-1}} \tau(X^{[2]}_\beta, a^{\beta}_j) b^{\beta}_j K_i 
+ c_\beta \frac{1}{q_i - q_i^{-1}} \tau(X_\beta^{[1]}, a^{\beta}_j) K_i^{-1} b^{\beta}_j \\
&= -\frac{1}{q_i - q_i^{-1}} (c_{\alpha_i} X^{[2]}_\beta K_i - c_\beta K_i^{-1} X_\beta^{[1]}),
\end{align*}
using Lemma \ref{tauskewpairingproperties} in the second step.  
Moreover, combining
$$ 
(\hat{\Delta} \otimes \id)\hat{\Delta}(F_i) = F_i \otimes 1 \otimes 1 + K_i^{-1} \otimes F_i \otimes 1 + K_i^{-1} \otimes K_i^{-1} \otimes F_i, 
$$ 
with the commutation relations from Lemma \ref{lemdrinfeldformcommutation} gives 
\begin{align*} 
X F_i &= c_{\alpha_i} \tau(\hat{S}(X^{[1]}_{\alpha_i}), F_i) X^{[2]}_{\beta} K_i 
+ F_i X + c_\beta \tau(X^{[2]}_{\alpha_i} K_\beta, F_i) K_i^{-1} X^{[1]}_\beta \\
&= -c_{\alpha_i} \tau(X^{[1]}_{\alpha_i}, F_i) X^{[2]}_{\beta} K_i 
+ F_i X + c_\beta \tau(X^{[2]}_{\alpha_i}, F_i) K_i^{-1} X^{[1]}_\beta \\
&= c_{\alpha_i} \frac{1}{q_i - q_i^{-1}} X^{[2]}_{\beta} K_i 
+ F_i X - c_\beta \frac{1}{q_i - q_i^{-1}} K_i^{-1} X^{[1]}_\beta, 
\end{align*}
taking into account that $ X^{[1]}_{\alpha_i} = E_i = X^{[2]}_{\alpha_i} $ in the last step. 
Comparing the above expressions yields the claim. \end{proof} 

\begin{lemma} \label{Casimirlemma}
Let $ \lambda \in \mathfrak{h}^*_q $. Then for all $ \beta \in \roots^+ $ and $ Y \in U_q(\mathfrak{n}_-)_{-\beta} $ we have 
$$
\Omega Y = q^{2(\rho, \beta) - (\beta, \beta)} Y \Omega K_{2 \beta}
$$
in $ \End(M(\lambda)) $. 
\end{lemma} 

\begin{proof} Assume that $ Y_1 \in U_q(\mathfrak{n}_-)_{-\beta_1}, Y_2 \in U_q(\mathfrak{n}_-)_{-\beta_2} $ for $ \beta_1, \beta_2 \in \roots^+ $ 
satisfy the above relation. Then 
\begin{align*}
\Omega Y_1 Y_2 &= q^{2(\rho, \beta_1) - (\beta_1, \beta_1)} Y_1 \Omega K_{2 \beta_1} Y_2 \\
&= q^{2(\rho, \beta_1) - (\beta_1, \beta_1)} q^{-(2 \beta_1, \beta_2)} Y_1 \Omega Y_2 K_{2 \beta_1} \\
&= q^{2(\rho, \beta_1) - (\beta_1, \beta_1)} q^{-(2 \beta_1, \beta_2)} q^{2(\rho, \beta_2) - (\beta_2, \beta_2)}  Y_1 Y_2 \Omega K_{2 (\beta_1 + \beta_2)} \\
&= q^{2(\rho, \beta_1 + \beta_2) - (\beta_1 + \beta_2, \beta_1 + \beta_2)} Y \Omega K_{2 (\beta_1 + \beta_2)}. 
\end{align*}
Hence using induction on the length of $ \beta $ it suffices to consider the case $ Y = F_i $ for some $ 1 \leq i \leq N $. 
Notice that in this case $ \beta = \alpha_i $ and $ 2 (\rho, \alpha_i) = (\alpha_i, \alpha_i) $. Hence the desired relation reduces to 
$$
\Omega F_i = F_i \Omega K_i^2.
$$
This in turn is a consequence of Lemma \ref{Casimirhelplemma}. Indeed, if we apply $ \hat{S}^{-1} \otimes \id $ to the relation from Lemma \ref{Casimirhelplemma} 
and multiply the tensor factors, we obtain
\begin{align*}
\sum_j &\hat{S}^{-1}(a_j^{\gamma}) F_i b^{\gamma}_j - \hat{S}^{-1}(a^{\gamma}_j) b^{\gamma}_j F_i \\
&= \sum_j \hat{S}^{-1}(F_i) \hat{S}^{-1}(a^{\gamma - \alpha_i}_j) b^{\gamma - \alpha_i}_j K_i 
- \hat{S}^{-1}(a^{\gamma - \alpha_i}_j) \hat{S}^{-1}(F_i) K_i^{-1} b^{\gamma - \alpha_i}_j
\end{align*}
for all $ \gamma \in \roots^+ $. Summing over $ \gamma $ this yields 
$$ 
\sum_{\gamma \in \roots^+} \sum_j \hat{S}^{-1}(a_j^{\gamma}) b^{\gamma}_j F_i 
= \sum_{\gamma \in \roots^+} \sum_j F_i K_i \hat{S}^{-1}(a^{\gamma - \alpha_i}_j) b^{\gamma - \alpha_i}_j K_i, 
$$
which implies $ \Omega F_i = F_i \Omega K_i^2 $ as desired. \end{proof} 

Let us say that $ \gamma \in \mathfrak{h}_q^* $ is a primitive weight of a weight module $ M $ if there exists 
a primitive vector $ v \in M $ of weight $ \gamma $. 

\begin{lemma} \label{Casimirvalue}
Let $ \lambda \in \mathfrak{h}_q^* $. If $ \gamma $ is a primitive weight of $ M(\lambda) $ 
then $ \gamma = \lambda - \beta $ for some $ \beta \in \roots^+ $ satisfying
$$
q^{2(\rho + \lambda, \beta) - (\beta, \beta)} = 1. 
$$
\end{lemma} 

\begin{proof}
If $ v \in M(\lambda) $ is a primitive vector then $ \Omega \cdot v = v $ since all the terms $ \Omega_\alpha $ for $ \alpha \neq 0 $ act on $ v $ 
by zero in this case. At the same time we have $ v = Y \cdot v_\lambda $ for some $ Y \in U_q(\mathfrak{n}_-)_{-\beta} $. Hence 
$$
v = \Omega \cdot v = \Omega Y \cdot v_\lambda = q^{2(\rho, \beta) - (\beta, \beta)} Y \Omega K_{2 \beta} \cdot v_\lambda 
= q^{2(\rho + \lambda, \beta) - (\beta, \beta)} v 
$$ 
by Lemma \ref{Casimirlemma}. This yields the claim. \end{proof}

\begin{lemma} \label{weightsrholemma}
Let $ \lambda \in \weights^+ $ and assume $ \gamma = \lambda - \beta $ for some $ \beta \in \roots^+ $ satisfies 
$$ 
(\gamma + \rho, \gamma + \rho) = (\lambda + \rho, \lambda + \rho). 
$$
If $ \gamma + \rho \in \weights^+ $ then $ \gamma = \lambda $. 
\end{lemma} 

\begin{proof} The given relation can be rewritten as 
$$ 
(\lambda - \gamma, \lambda + \rho) + (\gamma + \rho, \lambda - \gamma) = (\lambda + \rho, \lambda + \rho) 
+ (\gamma + \rho, \lambda - \gamma - \lambda - \rho) = 0. 
$$ 
Note that we have $ \lambda + \rho \in \weights^+ $. If $ \gamma + \rho \in \weights^+ $ then
since $ \lambda - \gamma \in \roots^+ $ it follows that both terms on the left hand side are positive, and therefore vanish separately. 
The equality $ (\lambda - \gamma, \lambda + \rho) = 0 $ implies $ \lambda - \gamma = 0 $, since $ \lambda + \rho $ is a 
linear combination of fundamental weights with all coefficients being strictly positive. 
\end{proof}  

\begin{theorem} \label{thmfdirreducible}
Let $ \mu \in \weights^+ $. Then the largest integrable quotient module $ L(\mu) $ of $ M(\mu) $ is irreducible, and hence $ L(\mu) \cong V(\mu) $. 
\end{theorem} 

\begin{proof} If $ L(\mu) $ is not simple then it must contain a proper submodule $ V $. Inside $ V $ 
we can find a primitive vector $ v $ of highest weight $ \gamma = \mu - \beta$ for some $ \beta \in \roots^+ $.  Note that $\gamma \in\weights^+$ by Theorem \ref{thmfdintegrable}.
  
Then Lemma \ref{Casimirvalue} implies $ 2 (\beta, \mu + \rho) - (\beta, \beta) = 0 $ and thus
\begin{align*}
(\gamma + \rho, \gamma + \rho) &= (\mu + \rho - \beta, \mu + \rho - \beta) \\
&= (\mu + \rho, \mu + \rho) - 2(\beta, \mu + \rho) + (\beta, \beta) = (\mu + \rho, \mu + \rho). 
\end{align*}
According to Lemma \ref{weightsrholemma} we conclude $ \gamma = \mu $. This means $ V = L(\mu) $, contradicting our assumption that $ V $ 
is a proper submodule. Hence $ L(\mu) $ is indeed simple, and therefore the canonical quotient map $ L(\mu) \rightarrow V(\mu) $ is 
an isomorphism. \end{proof}

We write $ \pi_\mu: U_q(\mathfrak{g}) \rightarrow \End(V(\mu)) $ for the algebra homomorphism
corresponding to the module structure on $ V(\mu) $. 
\nomenclature{$\pi_\mu$}{representation of $U_q(\lie{g})$ corresponding to the simple module $V(\mu)$ with $\mu\in\weights^+$}

\begin{lemma} \label{pilambdadense}
For every $ \mu \in \weights^+ $ we have $ \pi_\mu(U_q(\mathfrak{g})) = \End(V(\mu)) $. 
\end{lemma} 

\begin{proof} Assume that $ C \in \End(V(\mu)) $ is contained in the commutant of $ \pi_\mu(U_q(\mathfrak{g})) $, 
or equivalently, $ C $ is an intertwiner of $ V(\mu) $ with itself. Then $ C $ must map the highest weight vector $ v_\mu $ to a multiple 
of itself. It follows that $ C $ acts as a scalar on all of $ V(\mu) $ because $ v_\mu $ is cyclic. 
Now the Jacobson density theorem yields the claim. \end{proof}

\subsubsection{Complete reducibility} 

Having determined the irreducible finite dimensional weight modules over $ U_q(\mathfrak{g}) $, we obtain the following complete reducibility result 
in analogy to the situation for classical universal enveloping algebras, compare also the proof of Proposition \ref{slq2completeirreducibility}.

\begin{theorem} \label{thmuqgfiniterepdecomp}
Every finite dimensional integrable module over $ U_q(\mathfrak{g}) $ decomposes 
into a direct sum of irreducible highest weight modules $ V(\mu) $ for weights $ \mu \in \weights^+ $. 
\end{theorem} 

\begin{proof} We show that every extension $ 0 \rightarrow K \rightarrow E \rightarrow Q \rightarrow 0 $ of finite dimensional integrable modules 
splits. Notice that the extension automatically admits a $ U_q(\mathfrak{h}) $-linear splitting since all the modules 
involved are weight modules. 

In the first step consider the case that $ K = V(\mu) $ and $ Q = V(\lambda) $ are simple. If $ \mu = \lambda $ then 
we can lift the highest weight vector of $ V(\lambda) $ to a primitive vector in $ E $, which induces a splitting of the 
sequence by the characterisation of $ V(\lambda) $ as the largest integrable quotient of $ M(\lambda) $ from Theorem \ref{thmfdirreducible}.  

Assume now $ \mu \neq \lambda $. If $ \lambda \not< \mu $, that is, we do not have  $ \lambda = \mu - \nu $ for a nonzero element $ \nu \in \roots^+ $ then 
$ \lambda $ does not appear as a weight in $ V(\mu) $. 
Again we can lift the highest weight vector of $ V(\lambda) $ to a primitive vector in $ E $ and obtain a splitting. 
Finally, assume $ \mu \not< \lambda $, that is we do not have $ \lambda = \mu + \nu $ for some nonzero element $ \nu \in \roots^+ $. 
We consider the dual sequence $ 0 \rightarrow V(\lambda)^\vee \rightarrow E \rightarrow V(\mu)^\vee \rightarrow 0 $. Since 
$ V(\mu)^\vee \cong V(\mu) $ and $ V(\lambda)^\vee \cong V(\lambda) $ we obtain a splitting by our previous argument. 
Applying duality again shows that the original sequence is split as well. 

In the second step we consider arbitrary $ K $ and $ Q $ and show the existence of a splitting by induction on the dimension of $ E $. 
For $ \dim(E) = 0 $ the assertion is clear. Assume that every extension of dimension less than $ n $ splits and 
let $ \dim(E) = n $. If $ K = 0 $ or $ Q = 0 $ there is nothing to prove. Otherwise both $ K $ and $ Q $ have dimension $ < n $ and hence are direct sums 
of simple modules. Considering each simple block in $ Q $ and its preimage in $ E $ independently we see that we can restrict to the case that 
$ Q $ is simple. If $ K $ is simple then the claim follows from our above considerations. Otherwise there exists a proper simple 
submodule $ L \subset K $. Then $ 0 \rightarrow K/L \rightarrow E/L \rightarrow Q \rightarrow 0 $ 
is an extension, and since $ \dim(E/L) < n $ it is split. That is, $ Q $ is a direct summand in $ E/L $. Considering 
$ 0 \rightarrow L \rightarrow E \rightarrow E/L \rightarrow 0 $, we 
thus obtain a splitting $ Q \rightarrow E $ according to the first part of the proof. This map splits the sequence 
$ 0 \rightarrow K \rightarrow E \rightarrow Q \rightarrow 0 $, which finishes the proof. \end{proof}

As a consequence of Lemma \ref{lem:reduction_to_integrable}, we immediately obtain the analogous result for finite dimensional weight modules.

\begin{cor}
 \label{coruqgfiniterepdecomp}
 Every finite dimensional weight module over $U_q(\lie{g})$ decomposes into a direct sum of irreducible highest weight modules $V(\lambda)$ for weights $\lambda\in\weights_q^+$.
\end{cor}

Using Theorem \ref{thmuqgfiniterepdecomp} we can also strengthen Theorem \ref{uqgfdseparation}.

\begin{theorem} \label{uqgseparation}
The representations $ V(\mu) $ of $ U_q(\mathfrak{g}) $ for $ \mu \in \weights^+ $ separate points. 
More precisely, if $ X \in U_q(\mathfrak{g}) $ satisfies $ \pi_\mu(X) = 0 $ in $ \End(V(\mu)) $ for all $ \mu \in \weights^+ $ then $ X = 0 $. 
\end{theorem}

\begin{proof} According to Theorem \ref{uqgfdseparation} the finite dimensional integrable representations separate the points of $ U_q(\mathfrak{g}) $. 
However, every such representation is a finite direct sum of representations $ V(\mu) $ by Theorem \ref{thmuqgfiniterepdecomp}. This yields the claim. \end{proof}

\subsection{Quantized algebras of functions} \label{secoqg}

We define quantized algebras of functions and discuss their duality with quantized universal enveloping algebras. It is 
assumed throughout that $ q = s^L \in \mathbb{K}^\times $ is not a root of unity. 

\begin{definition} \label{defogq}
Let $ \mathfrak{g} $ be a semisimple Lie algebra. The quantized algebra of functions $ \Poly(G_q) $ consists of all matrix coefficients of 
finite dimensional integrable $ U_q(\mathfrak{g}) $-modules, that is, 
$$
\Poly(G_q) = \bigoplus_{\mu \in \weights^+} \End(V(\mu))^*, 
$$
\nomenclature[o$O(G_q)$]{$\Poly(G_q)$}{quantized algebra of polynomial functions on $G_q$}%
with the direct sum on the right hand side being an algebraic direct sum. 
\end{definition} 

The direct sum decomposition of $ \Poly(G_q) $ in Definition \ref{defogq} is called the Peter-Weyl decomposition. 

By construction, there exists a bilinear pairing between $ U_q(\mathfrak{g}) $ and $ \Poly(G_q) $ given by evaluation of linear functionals. 
We will write $ (X, f) $ for the value of this pairing on $ X \in U_q(\mathfrak{g}) $ and $ f \in \Poly(G_q) $, 
\label{nom:OGq-pairing}%
and we point out that this pairing is nondegenerate. Indeed, $ (X, f) = 0 $ for all $ X \in U_q(\mathfrak{g}) $ implies $ f = 0 $ 
by Lemma \ref{pilambdadense}. On the other hand, $ (X, f) = 0 $ for all $ f \in \Poly(G_q) $ implies $ X = 0 $ by Theorem \ref{uqgseparation}. 

We define the product and coproduct of $ \Poly(G_q) $ such that 
\begin{align*}
(XY, f) &= (X, f_{(1)}) (Y, f_{(2)}), \qquad (X, fg) = (X_{(2)}, f) (X_{(1)}, g), 
\end{align*}
and the antipode $ S $, unit $ 1 $ and counit $ \epsilon $ of $ \Poly(G_q) $ are determined by 
\begin{align*}
(\hat{S}(X), f) &= (X, S^{-1}(f)), \qquad (\hat{S}^{-1}(X), f) = (X, S(f)) 
\end{align*} 
and 
\begin{align*}
\hat{\epsilon}(X) &= (X, 1), \qquad \epsilon(f) = (1, f) 
\end{align*} 
for $ X, Y \in U_q(\mathfrak{g}) $ and $ f, g \in \Poly(G_q) $. 
We shall write $ \Delta $ for the coproduct of $ \Poly(G_q) $, and in the above formulas 
we have used the Sweedler notation $ \Delta(f) = f_{(1)} \otimes f_{(2)} $ for $ f \in \Poly(G_q) $. 
With the above structure maps, the quantized algebra of functions becomes a Hopf algebra. 

Using the transpose of the canonical isomorphism $ \End(V(\mu)) \cong V(\mu) \otimes V(\mu)^* $ we obtain an isomorphism 
$$
\Poly(G_q) \cong \bigoplus_{\mu \in \weights^+} V(\mu)^* \otimes V(\mu). 
$$
More precisely, if $ v \in V(\mu) $ and $ v^* \in V(\mu)^* $ we shall write $ \bra v^*| \bullet |v \ket \in V(\mu)^* \otimes V(\mu) \subset \Poly(G_q) $ 
\nomenclature[o$(\ ,\ )$7]{$\bra v^*\vert\bullet\vert v\ket$}{matrix coefficient}%
for the matrix coefficient determined by 
$$ 
(X, \bra v^*| \bullet |v \ket) = v^*(\pi_\mu(X)(v)). 
$$
The right regular action of $ U_q(\mathfrak{g}) $ on $ \Poly(G_q) $ given by $ X \hit f = f_{(1)} (X, f_{(2)}) $ 
\label{nom:hit_Uqg}%
corresponds to $ X \hit \bra v^*| \bullet |v \ket = \bra v^*| \bullet |X \cdot v \ket $ in this picture.  
The left regular action given by $ X \cdot f = f \hitby \hat{S}(X) = (\hat{S}(X), f_{(1)}) f_{(2)} $
corresponds to $ X \cdot \bra v^*| \bullet |v \ket = \bra X \cdot v^*| \bullet |v \ket $. 
Here we work with the action $ (X \cdot v^*)(v) = v^*(\hat{S}(X) \cdot v) $ on the dual space $ V(\mu)^* $ as usual. 
%
%

If $ e^\mu_1, \dots, e^\mu_n $ is a basis of $ V(\mu) $ with dual basis $ e_\mu^1, \dots, e_\mu^n \in V(\mu)^* $, so that $ e_\mu^i(e^\mu_j) = \delta_{ij} $, 
then the matrix coefficients $ u^\mu_{ij} = \bra e_\mu^i| \bullet | e^\mu_j \ket $
\nomenclature{$u^\mu_{ij}$}{matrix coefficients of the simple integrable module $V(\mu)$ with respect to dual bases, $ u^\mu_{ij} = \bra e_\mu^i\vert \bullet \vert e^\mu_j \ket $}%
form a basis of $ \End(V(\mu))^* $. 
Notice that we have 
$$
\Delta(u^\mu_{ij}) = \sum_{k = 1}^n u^\mu_{ik} \otimes u^\mu_{kj}
$$
for all $ \mu \in \weights^+ $ and $ 1 \leq i,j \leq n = \dim(V(\mu)) $. 
By the remarks preceding Definition \ref{def:Hopf_algebra}, this shows that $ \Poly(G_q) $ is a cosemisimple Hopf algebra. By Proposition \ref{cosemisimpleintegral}, it admits a
left and right invariant integral $ \phi: \Poly(G_q) \rightarrow \mathbb{K} $, specifically 
$$
\phi(u^\mu_{ij}) = 
\begin{cases}
1 & \text{ if } \mu = 0 \\
0 & \text{ else.} 
\end{cases} 
$$
\label{nom:OGq_integral}%
Since $ \Poly(G_q) $ is a cosemisimple Hopf algebra, and hence in particular a regular multiplier Hopf algebra with integrals, it admits a
dual multiplier Hopf algebra $ \D(G_q) $, compare Subsection \ref{sec:dual_multiplier_Hopf_algebra}. Explicitly, we have 
$$
\D(G_q) = \bigoplus_{\mu \in \weights^+} \End(V(\mu)), 
$$
\nomenclature[o$D(G_q)$]{$\D(G_q)$}{dual multiplier Hopf algebra to $\Poly(G_q)$}%
such that the canonical skew-pairing $ (x, f) $ for $ x \in \D(G_q) $ and $ f \in \Poly(G_q) $ corresponds to the obvious pairing between $ \End(V(\mu)) $ 
and $ \End(V(\mu))^* $ in each block. The product and coproduct of $ \D(G_q) $ are fixed such that 
\begin{align*}
(xy, f) &= (x, f_{(1)}) (y, f_{(2)}), \qquad (x, fg) = (x_{(2)}, f) (x_{(1)}, g), 
\end{align*}
and the antipode $ \hat{S} $, counit $ \hat{\epsilon} $ of $ \D(G_q) $ and unit $ 1 \in \M(\D(G_q)) $ are determined by 
\begin{align*}
(\hat{S}(x), f) &= (x, S^{-1}(f)), \qquad (\hat{S}^{-1}(x), f) = (x, S(f)) 
\end{align*} 
and 
\begin{align*}
\hat{\epsilon}(x) &= (x, 1), \qquad \epsilon(f) = (1, f) 
\end{align*} 
for $ x, y \in \D(G_q) $ and $ f, g \in \Poly(G_q) $. We shall write $ \hat{\Delta} $ for the coproduct of $ \D(G_q) $, and use 
the Sweedler notation $ \hat{\Delta}(x) = x_{(1)} \otimes x_{(2)} $ for $ x \in \D(G_q) $. 

It follows immediately from the construction that there is a canonical homomorphism $ U_q(\mathfrak{g}) \rightarrow \M(\D(G_q)) $ of multiplier 
Hopf algebras, namely
\[
 X \mapsto \prod_{\mu\in\weights^+} \pi_\mu(X). 
\]
This homomorphism is injective according to Theorem \ref{uqgseparation}. 
In the sequel it will often be useful to view elements of $ U_q(\mathfrak{g}) $ as multipliers of $ \D(G_q) $. 

Given a basis of matrix coefficients $ u^\mu_{ij} $ as above, the algebra $ \D(G_q) $ admits a linear basis $ \omega^\mu_{ij} $ 
\label{nom:dual_basis_Gq}%
such that 
$$
(\omega^\mu_{ij}, u^\nu_{kl}) = \delta_{\mu \nu} \delta_{ik} \delta_{jl}.  
$$
Note that with respect to the algebra structure of $ \D(G_q) $ we have 
$$
\omega^\mu_{ij} \omega^\nu_{kl} = \delta_{\mu \nu} \delta_{jk} \omega^\mu_{il} 
$$
that is, the elements $ \omega^\mu_{ij} $ are matrix units for $ \D(G_q) $.

We note that the algebra $ \Poly(G_q) $ can be viewed as a deformation of the coordinate ring of the 
affine algebraic variety associated with the group $ G $ corresponding to $ \mathfrak{g} $.

\subsection{The universal $ R $-matrix} \label{subsecuniversalrmatrix}

In this subsection we discuss universal $ R $-matrices and the dual concept of universal $ r $-forms. We derive an explicit formula for the universal $ r $-form 
of $ \Poly(G_q) $ in terms of the PBW-basis of $ U_q(\mathfrak{g}) $. 

\subsubsection{Universal $ R $-matrices} 

Let us first define what it means for a regular multiplier Hopf algebra to be quasitriangular. 
\begin{definition} \label{defquasitriangular}
Let $ H $ be a regular multiplier Hopf algebra. Then $ H $ is called called quasitriangular if there exists an 
invertible element $ \R \in \M(H \otimes H) $, called a universal $R$-matrix,
such that 
$$
\Delta^{\cop}(x) = \R \Delta(x) \R^{-1} 
$$ 
for all $ x \in H $ and 
$$ 
(\Delta \otimes \id)(\R) = \R_{13} \R_{23}, \quad (\id \otimes \Delta)(\R) = \R_{13} \R_{12}. 
$$
\end{definition} 
Here $ \Delta^{\cop} $ denotes the opposite comultiplication on $ H $, given by $ \Delta^{\cop} = \sigma \Delta $, 
where $ \sigma: H \otimes H \rightarrow H \otimes H$, $\sigma(x \otimes y) = y \otimes x $ is the flip map. 
Also, we are using the standard leg-numbering notation, whereby $\R_{12}=\R\otimes1$,  $\R_{23} = 1 \otimes \R$ and $\R_{13} = (\sigma\otimes\id)\R_{23}$ in $\M(H \otimes H \otimes H)$. 
\label{nom:leg-numbering}%

Notice that the final two relations in Definition \ref{defquasitriangular} imply 
$$ 
(\epsilon \otimes \id)(\R) = 1 = (\id \otimes \epsilon)(\R). 
$$ 
Moreover we record the formulas 
$$
(S \otimes \id)(\R) = \R^{-1} = (\id \otimes S^{-1})(\R), 
$$
obtained using the antipode property applied to the final two relations in Definition \ref{defquasitriangular} and the previous formulas. 

From a technical point of view it is sometimes easier to view the $ R $-matrix as a linear functional on the dual side. 
If $ H $ is a Hopf algebra and $ l,k: H \rightarrow \mathbb{K} $ are linear functionals on $ H $ we define the convolution $ l * k  $ to be the 
linear functional on $ H $ defined by 
$$
(l * k)(f) = l(f_{(1)}) k(f_{(2)})
$$
\label{nom:convolution}%
for $ f \in H $. It is straightforward to check that the counit $ \epsilon $ is a neutral element with respect to convolution. Accordingly, 
we say that $ l: H \rightarrow \mathbb{K} $ is convolution invertible if there exists a linear functional $ k: H \rightarrow \mathbb{K} $ such that 
\begin{align*} 
l(f_{(1)}) k(f_{(2)}) = \epsilon(f) = k(f_{(1)}) l(f_{(2)})
\end{align*}
for all $ f \in H $. In this case the functional $ k $ is uniquely determined and we write $ k = l^{-1} $.

\begin{definition} \label{defrform}
Let $ H $ be a Hopf algebra. Then $ H $ is called coquasitriangular if there exists 
a linear functional $ r: H \otimes H \rightarrow \mathbb{K} $, called a universal $r$-form, which is invertible with respect to the convolution multiplication and satisfies 
\begin{bnum}
\item[a)] $ r(f \otimes gh) = r(f_{(1)} \otimes g) r(f_{(2)} \otimes h) $ 
\item[b)] $ r(fg \otimes h) = r(f \otimes h_{(2)}) r(g \otimes h_{(1)}) $
\item[c)] $ r(f_{(1)} \otimes g_{(1)}) f_{(2)} g_{(2)} = g_{(1)} f_{(1)} r(f_{(2)} \otimes g_{(2)}) $
\end{bnum}
for all $ f,g,h \in H $. 
\end{definition} 

If $ H $ is a cosemisimple Hopf algebra then it is in particular a regular multiplier Hopf algebra with integrals. Moreover, 
the linear dual space $ H^* $ identifies with $ \M(\hat{H}) $ in a natural way, and similarly $ (H \otimes H)^* \cong \M(\hat{H} \otimes \hat{H}) $. 
In particular, a linear map $ r: H \otimes H \rightarrow \mathbb{K} $ corresponds to 
a unique element $ \R \in \M(\hat{H} \otimes \hat{H}) $ such that $ (\R, f \otimes g) = r(f \otimes g) $. 
With these observations in mind, the relation between universal $ r $-forms and universal $ R $-matrices is as follows. 

\begin{lemma} \label{skewpairingduality}
Let $ H $ be a cosemisimple Hopf algebra and let $ r: H \otimes H \rightarrow \mathbb{K} $ be a linear form, and let $ \R$ be the corresponding element of $ \M(\hat{H} \otimes \hat{H}) $ as explained above. Then 
\begin{bnum}
\item[a)] $ r $ is convolution invertible iff $ \R $ is invertible. 
\item[b)] $ r(f \otimes gh) = r(f_{(1)} \otimes g) r(f_{(2)} \otimes h) $ for all $ f,g,h \in H $ iff $ (\id \otimes \hat{\Delta})(\R) = \R_{13} \R_{12} $. 
\item[c)] $ r(fg \otimes h) = r(f \otimes h_{(2)}) r(g \otimes h_{(1)}) $ for all $ f,g,h \in H $ iff $ (\hat{\Delta} \otimes \id)(\R) = \R_{13} \R_{23} $. 
\item[d)] $ r(f_{(1)} \otimes g_{(1)}) g_{(2)} f_{(2)} = f_{(1)} g_{(1)} r(f_{(2)} \otimes g_{(2)}) $ for all $ f,g,h \in H $ iff 
$ \hat{\Delta}^{\cop}(x) = \R \hat{\Delta}(x) \R^{-1} $ for all $ x \in \hat{H} $. 
\end{bnum} 
In particular, $ r $ is a universal $ r $-form for $ H $ iff $ \R $ is a universal $ R $-matrix for $ \hat{H} $. 
\end{lemma} 

\begin{proof} $ a) $ This follows immediately from the definition of convolution. The convolution inverse $ r^{-1} $ of $ r $ 
corresponds to the inverse $ \R^{-1} $ of $ \R $. 

$ b) $ Using nondegeneracy of the pairing between $ H \otimes H $ and $ \hat{H} \otimes \hat{H} $, this follows from 
\begin{align*} 
r(f \otimes gh) &= (\R, f \otimes gh) = ((\id \otimes \hat{\Delta})(\R), f \otimes h \otimes g) 
\end{align*}
and 
\begin{align*} 
r(f_{(1)} \otimes g) r(f_{(2)} \otimes h) = (\R, f_{(1)} \otimes g) (\R, f_{(2)} \otimes h) = (\R_{13} \R_{12}, f \otimes h \otimes g) 
\end{align*}
for $ f,g,h \in H $. 

$ c) $ This follows from 
\begin{align*} 
r(fg \otimes h) &= (\R, fg \otimes h) = ((\hat{\Delta} \otimes \id)(\R), g \otimes f \otimes h) 
\end{align*}
and 
\begin{align*}
r(f \otimes h_{(2)}) r(g \otimes h_{(1)}) &= (\R, f \otimes h_{(2)}) (\R, g \otimes h_{(1)}) = (\R_{13} \R_{23}, g \otimes f \otimes h) 
\end{align*}
for all $ f,g,h \in H $. 

$ d) $ For any $ y \in \hat{H} $ and $ f,g,h $ we compute 
\begin{align*} 
r(f_{(1)} \otimes g_{(1)}) (y, g_{(2)} f_{(2)}) &= (\R, f_{(1)} \otimes g_{(1)}) (\hat{\Delta}(y), f_{(2)} \otimes g_{(2)}) 
= (\R \hat{\Delta}(y), f \otimes g) 
\end{align*}
and 
\begin{align*}
(y, f_{(1)} g_{(1)}) r(f_{(2)} \otimes g_{(2)}) = (\hat{\Delta}(y), g_{(1)} \otimes f_{(1)}) (\R, f_{(2)} \otimes g_{(2)}) 
= (\hat{\Delta}^{\cop}(y) \R, f \otimes g).  
\end{align*}
Since the pairing between $ \hat{H} $ and $ H $ is nondegenerate, the relation $ \R \hat{\Delta}(x) = \hat{\Delta}^{\cop}(x) \R $ for all $ x \in \hat{H} $ 
is equivalent to 
$$ 
r(f_{(1)} \otimes g_{(1)}) g_{(2)} f_{(2)} = f_{(1)} g_{(1)} r(f_{(2)} \otimes g_{(2)}) 
$$ 
for all $ f,g,h \in H $. \end{proof}

\subsubsection{The $ R $-matrix for $ U_q(\mathfrak{g}) $} 
 
The quantized universal enveloping algebra $ U_q(\mathfrak{g}) $ is not quite quasitriangular in the sense of 
Definition \ref{defquasitriangular}, but the multiplier Hopf algebra $ \D(G_q) $ is. We assume that $ q = s^L \in \mathbb{K}^\times $ is not a root of unity.

Define $ \Poly(B^\pm_q) $
\nomenclature[o$O(B_q)$]{$\Poly(B^\pm_q)$}{image of $\Poly(G_q)$ under the restriction map $U_q(\lie{g})^* \to U_q(\lie{b}_\pm)^*$}%
to be the image of $ \Poly(G_q) $ under the natural 
projection $ U_q(\mathfrak{g})^* \rightarrow U_q(\mathfrak{b}_\pm)^* $. Then there is a unique Hopf algebra structure on $ \Poly(B^\pm_q) $ such that 
the projection maps $ \Poly(G_q) \rightarrow \Poly(B^\pm_q) $ are Hopf algebra homomorphisms. 

Recall the definition of the Drinfeld pairing from Definition \ref{def:Drinfeld_pairing}. 

\begin{prop} \label{tauboreliso}
The Drinfeld pairing $ \tau $ induces isomorphisms $\iota_\pm: U_q(\mathfrak{b}_\pm) \to \Poly(B^\mp_q) $ of Hopf algebras, such that
\[
\iota_+(X)(Y) = \tau(\hat{S}(X), Y), \qquad \iota_-(Y)(X) = \tau(X, Y) ,
\]
\nomenclature{$\iota_\pm$}{Hopf algebra isomorphisms $\iota_\pm: U_q(\mathfrak{b}_\pm) \to \Poly(B^\mp_q) $}%
for all $X\in U_q(\mathfrak{b}_+)$, $Y\in U_q(\mathfrak{b}_-)$.
\end{prop}

\begin{proof} By Theorem \ref{lemdrinfeldnondegenerate} $ \tau $ is a nondegenerate skew-pairing, so the maps $ \iota_\pm: U_q(\mathfrak{b}_\pm) \rightarrow U_q(\mathfrak{b}_\mp)^* $ are injective algebra homomorphisms.

To show that $\iota_+$ is an isomorphism, we consider the direct sum decompositions
\begin{align*}
  U_q(\mathfrak{b}_+) &= \bigoplus_{\beta\in\weights^+,~\mu\in\weights} A_{\beta,\mu}, \\
  \Poly(B_q^-) &= \bigoplus_{\beta\in\weights^+,~\mu\in\weights} B_{\beta,\mu},
\end{align*}
where
\[
  A_{\beta,\mu} = \hat{S}^{-1}(U_q(\lie{n}_+)_\beta K_\mu) = \{ \hat{S}^{-1}(XK_\mu) \mid X\in U_q(\mathfrak{n}_+)_\beta \} 
\]
and
\[
  B_{\beta,\mu} =  \{ \bra v^+|\bullet|w^+\ket \mid w^+\in V(\lambda)_{-\mu},\ v^+ \in V(\lambda)^*_{\mu+\beta} \text{ for some } \lambda \in \weights^+\}.
\]
We claim that $\iota_+$ restricts to an isomorphism between the finite dimensional subspaces $A_{\beta,\mu}$ and $B_{\beta,\mu}$ for each $\beta\in\weights^+$, $\mu\in\weights$.

Let $X\in U_q(\lie{n}_+)_\beta$.  Using Lemma \ref{tauskewpairingproperties} we see that for every $Y \in U_q(\mathfrak{n}_-)$ and $\nu\in\weights$,
\begin{align*}
 \iota_+(\hat{S}^{-1}(XK_\mu))(YK_\nu) 
  = \tau(XK_\mu, Y K_\nu ) 
  &= q^{-(\mu,\nu)}\tau(X, Y ) \\
  &= q^{-(\mu,\nu)}\tau(XK_\mu, Y ) \\
  &= q^{-(\mu,\nu)} \iota_+(\hat{S}^{-1}(XK_\mu))(Y)
\end{align*}

Therefore, by the orthogonality of the PBW basis (Theorem \ref{PBWdrinfeldorthogonal}), the functional $\varphi = \iota_+(\hat{S}^{-1}(XK_\mu)) \in U_q(\lie{b}_-)^*$
satisfies
\begin{equation}
 \label{eq:iota_property}
 \varphi(YK_\nu) = \delta_{\beta\gamma} q^{-(\mu,\nu)} \varphi(Y)  \qquad \text{for every }Y \in U_q(\mathfrak{n}_-)_{-\gamma},\ \gamma\in\weights^+,\ \nu\in\weights.
\end{equation}
Conversely, any linear functional $\varphi\in U_q(\mathfrak{b}_-)^*$ satisfying the  property  \eqref{eq:iota_property} is completely determined by its restriction to $U_q(\lie{n}_-)_{-\beta}$, so  the orthogonality of the PBW basis implies that it is the image of a unique element $\hat{S}^{-1}(XK_\mu) \in A_{\beta,\mu}$.

Therefore, it suffices to show that the set of $\varphi$ satisfying \eqref{eq:iota_property} agrees with the subspace of matrix coefficients $\bra v^+|\bullet|w^+\ket $ with $w^+\in V(\lambda)_{-\mu}$, $v^+ \in V(\lambda)^*_{\mu+\beta}$, as claimed.  For such a matrix coefficient, we have
\begin{align*}
  (YK_\nu, \bra v^+|\bullet|w^+\ket) 
    &=  v^+ ( YK_\nu \cdot w^+ ) \\
    &= q^{-(\nu,\mu)}  v^+ ( Y \cdot w^+ ) 
    = q^{-(\nu,\mu)} (Y, \bra v^+|\bullet|w^+\ket),
\end{align*}
so that it satisfies \eqref{eq:iota_property}.  Conversely, let $\varphi\in U_q(\mathfrak{b}_-)^*$ satisfy \eqref{eq:iota_property}.  By Lemma \ref{lem:L_looks_like_Uqn} we can choose $\lambda\in \weights^+$ sufficiently large such that,
after fixing an arbitrary nonzero $w^+ \in V(\lambda)_{-\mu}$, the map
\[
  U_q(\mathfrak{n}_-)_{-\beta} \to V(\lambda)_{-\mu-\beta}, \qquad Y \mapsto Y\cdot w^+
\]
is injective.  Therefore we can find $v^+ \in V(\lambda)^*_{\mu+\beta}$ such that
\[
 (Y, \bra v^+| \bullet |w^+\ket) = \varphi(Y)
\]
for every $Y\in U_q(\mathfrak{n}_-)_{-\beta}$.  Since $\bra v^+| \bullet |w^+\ket$ vanishes on $U_q(\mathfrak{n}_-)_{-\gamma}$ for all $\gamma \neq \beta$ by weight considerations, it follows that $\bra v^+| \bullet |w^+\ket = \varphi$ as an element of $U_q(\mathfrak{b}_-)^*$.  This completes the proof for $\iota_+$.

A similar argument shows that $\iota_-$ is an isomorphism.

Finally, using again the skew-pairing property of $ \tau $ 
it is straightforward to check that the maps $ \iota_\pm $ are compatible with the
comultiplications. 
\end{proof} 

Since $ U_q(\mathfrak{b}_\pm) \cong \Poly(B^\mp_q) $, the Drinfeld pairing 
$ \tau: U_q(\mathfrak{b}_+) \otimes U_q(\mathfrak{b}_-) \rightarrow \mathbb{K} $ induces a 
pairing $ \Poly(G_q) \otimes \Poly(G_q) \rightarrow \mathbb{K} $ as follows, compare \cite{Gaitsgoryrmatrix}. 

\begin{prop} \label{rformconstruction}
The linear map $ r: \Poly(G_q) \otimes \Poly(G_q) \rightarrow \mathbb{K} $ given by 
\begin{equation*}
\xymatrix{
\Poly(G_q) \otimes \Poly(G_q) \ar@{->}[r] & \Poly(B^-_q) \otimes \Poly(B^+_q) \cong U_q(\mathfrak{b}_+) \otimes U_q(\mathfrak{b}_-) 
\ar@{->}[r]^{\quad \qquad \quad \qquad \quad \tau} & \mathbb{K}
}
\end{equation*} 
is a universal $ r $-form for $ \Poly(G_q) $. Here the first arrow is the flip map composed with the canonical projection, 
and $ \tau $ is the Drinfeld pairing. 
That is 
\[
 r(f \otimes g) = \tau(\iota_+^{-1}(g), \iota_-^{-1}(f))
\]
where, by slight abuse of notation, we identify $ f $ and $ g $ with their images in $ \Poly(B^+_q) $ and $ \Poly(B^-_q) $, respectively. 
\end{prop} 

\begin{proof} 
Since $ \tau $ is a skew-pairing we see that $ r $ is a convolution invertible linear map satisfying conditions \emph{a)} and \emph{b)}
of Definition \ref{defrform}. 

It remains to verify condition \emph{c)} of Definition \ref{defrform}. That is, we have to check 
\begin{align*} 
r(f_{(1)} \otimes g_{(1)}) g_{(2)} f_{(2)} = f_{(1)} g_{(1)} r(f_{(2)} \otimes g_{(2)}) 
\end{align*} 
for all $ f, g \in \Poly(G_q) $. Since the canonical pairing between $ U_q(\mathfrak{g}) $ and $ \Poly(G_q) $ is nondegenerate, 
it suffices to show 
\begin{align*} 
r(f_{(1)} \otimes g_{(1)}) (X, g_{(2)} f_{(2)}) = (X, f_{(1)} g_{(1)}) r(f_{(2)} \otimes g_{(2)}) 
\end{align*} 
for all $ X \in U_q(\mathfrak{g}) $. Moreover, since $ U_q(\mathfrak{g}) $ is generated as a Hopf algebra by $ U_q(\mathfrak{b}_+) $ and $ U_q(\mathfrak{b}_-) $, 
it suffices to consider $ X \in U_q(\mathfrak{b}_\pm) $. Abbreviating $ Z = \iota_-^{-1}(f) \in U_q(\mathfrak{b}_-) $ and taking $ X \in U_q(\mathfrak{b}_+) $, 
we indeed compute 
\begin{align*} 
r(f_{(1)} \otimes g_{(1)}) (X, g_{(2)} f_{(2)}) &= \tau(\iota_+^{-1}(g_{(1)}), \iota_-^{-1}(f_{(1)})) \tau(X, \iota_-^{-1}(g_{(2)} f_{(2)})) \\
&= \tau(\iota_+^{-1}(g_{(1)}), \iota_-^{-1}(f_{(1)})) \tau(X_{(1)}, \iota_-^{-1}(f_{(2)})) \tau(X_{(2)}, \iota_-^{-1}(g_{(2)})) \\
&= \tau(\iota_+^{-1}(g_{(1)}), Z_{(1)}) \tau(X_{(1)}, Z_{(2)}) \tau(X_{(2)}, \iota_-^{-1}(g_{(2)})) \\
&= (\hat{S}(Z_{(1)}), g_{(1)}) \tau(X_{(1)}, Z_{(2)}) (X_{(2)}, g_{(2)}) \\
&= \rho(\hat{S}(Z_{(2)}), X_{(1)}) (\hat{S}(Z_{(1)}) X_{(2)}, g) \\
&= (X_{(1)} \hat{S}(Z_{(2)}), g) \rho(\hat{S}(Z_{(1)}), X_{(2)}) \\
&= \tau(X_{(2)}, Z_{(1)}) (X_{(1)}, g_{(1)}) (\hat{S}(Z_{(2)}), g_{(2)}) \\
&= \tau(X_{(2)}, Z_{(1)}) \tau(X_{(1)}, \iota_-^{-1}(g_{(1)})) \tau(\iota_+^{-1}(g_{(2)}), Z_{(2)}) \\
&= \tau(X_{(2)}, \iota_-^{-1}(f_{(1)})) \tau(X_{(1)}, \iota_-^{-1}(g_{(1)})) \tau(\iota_+^{-1}(g_{(2)}), \iota_-^{-1}(f_{(2)})) \\
&= \tau(X, \iota_-^{-1}(f_{(1)} g_{(1)})) \tau(\iota_+^{-1}(g_{(2)}), \iota_-^{-1}(f_{(2)})) \\
&= (X, f_{(1)} g_{(1)}) r(f_{(2)} \otimes g_{(2)}), 
\end{align*} 
using Lemma \ref{lemdrinfeldformcommutation}. A similar argument works for $ X \in U_q(\mathfrak{b}_-) $. This finishes the proof. \end{proof} 

Let us now give an explicit formula for the universal $ r $-form of $ \Poly(G_q) $.   
We will use the notation $H_i \in \mathfrak{h}$ ($1\leq i \leq N$)
\nomenclature{$H_i$}{simple coroot in $\mathfrak{g}$, identified with $\alpha_i^\vee$ via the pairing $(\ , \ )$}%
for the simple coroots of $\mathfrak{g}$, identified with $\alpha_i^\vee \in \roots^\vee$ via the pairing:
\[
 \lambda(H_i) = (\lambda,\alpha_i^\vee)
\]
for all $\lambda \in \weights$. 
We also fix quantum root vectors $E_{\beta_1}, \ldots, E_{\beta_n}$, $F_{\beta_1},\ldots, F_{\beta_n}$ associated to a reduced expression of $w_0\in W$, as in Definition \ref{def:quantum_root_vectors}.

\begin{theorem} \label{Rmatrixformula}
The multiplier Hopf algebra $ \D(G_q) $ has a universal $ R $-matrix given by 
$$
\R = q^{\sum_{i,j = 1}^N B_{ij}(H_i \otimes H_j)} 
\prod_{\alpha \in {\bf \Delta}^+} \exp_{q_\alpha}((q_\alpha - q_\alpha^{-1})(E_\alpha \otimes F_\alpha))
$$
where $ H_1, \dots, H_N $ are the simple coroots of $ \mathfrak{g} $, the matrix $ (B_{ij}) $ is the inverse of $ (C_{ij}) = (d_j^{-1} a_{ij}) $,
%
%
and $ q_\alpha = q_i $ if $ \alpha $ and $ \alpha_i $ lie on the same $ W $-orbit. Moreover the factors on the right hand side appear 
in the order $ \beta_1, \dots, \beta_n $. \\
Equivalently, the above formula determines a universal $ r $-form for $ \Poly(G_q) $.
\end{theorem} 

\begin{proof} 
 The formula requires a little explanation.  Note that $C_{ij}=(\alpha_i^\vee,\alpha_j^\vee) =(H_i,H_j)$, so that
 \[
  \sum_{i,j = 1}^N B_{ij}(H_i \otimes H_j) = \sum_{i = 1}^N H_i \otimes \varpi_i = \sum_{i = 1}^N \varpi_i \otimes H_i,
 \]
 where we recall that the fundamental weights $\varpi_i$ form the dual basis to the simple coroots $H_i$.
 The term $q^{\sum_{i,j = 1}^N B_{ij}(H_i \otimes H_j)}$ should thus be interpreted as the element in $\M(\D(G_q)\otimes\D(G_q))$ such that
 \[
  q^{\sum_{i,j = 1}^N B_{ij}(H_i \otimes H_j)}\cdot (v\otimes w) = q^{(\mu,\nu)} v\otimes w,
 \]
  whenever $V$ and $W$ are integrable $U_q(\mathfrak{g})$-modules and $v\in V$, $w\in W$ have weight $\mu$ and $\nu$ respectively,

We shall show that the universal $ r $-form obtained in Proposition \ref{rformconstruction} is given by the stated formula. 
It is technically more convenient however, to compute the inverse $ r^{-1} $, given by 
$$ 
r^{-1}(f \otimes g) = \tau(\hat{S}(\iota_+^{-1}(g)), \iota_-^{-1}(f)) = \rho(\iota_-^{-1}(f)), \iota_+^{-1}(g)). 
$$
Note that according to Lemma \ref{qexpinverse} the inverse of the element $ \R $ above is given by 
$$ 
\R^{-1} = \left( \prod_{\alpha \in {\bf \Delta}^+} \exp_{q_\alpha^{-1}}(-(q_\alpha - q_\alpha^{-1})(E_\alpha \otimes F_\alpha)) \right) 
q^{-\sum_{i,j = 1}^N B_{ij}(H_i \otimes H_j)}, 
$$
where now the ordering of the positive roots in the product is reversed to $ \beta_n, \dots, \beta_ 1 $.

To compute $ r^{-1} $ it is sufficient to consider $r^{-1}(f \otimes g)$ for elements of the form
\[
f = \iota_-(F_{\beta_n}^{t_n} \cdots F_{\beta_1}^{t_1} K_\nu), \qquad 
g = \iota_+(\hat{S}^{-1}(E_{\beta_n}^{s_n} \cdots E_{\beta_1}^{s_1} K_\mu)). 
\]
Using the orthogonality relations for the Drinfeld pairing in Lemma \ref{drinfeldpairingevaluate}, we calculate
\begin{align*}
 r^{-1}(f \otimes g) 
  &= \tau(E_{\beta_n}^{s_n} \cdots E_{\beta_1}^{s_1}, F_{\beta_n}^{t_n} \cdots F_{\beta_1}^{t_1}) q^{-(\mu, \nu)} \\
  &= \sum_{a_1, \dots, a_n = 0}^\infty \tau(E_{\beta_n}^{a_n} \cdots E_{\beta_1}^{a_1}, F_{\beta_n}^{a_n} \cdots F_{\beta_1}^{a_1})^{-1} \\
&\qquad \times \tau(E_{\beta_n}^{s_n} \cdots E_{\beta_1}^{s_1}, F_{\beta_n}^{a_n} \cdots F_{\beta_1}^{a_1}) \;
 \tau(E_{\beta_n}^{a_n} \cdots E_{\beta_1}^{a_1}, F_{\beta_n}^{t_n} \cdots F_{\beta_1}^{t_1}) q^{-(\mu, \nu)} \\
  &= \sum_{a_1, \dots, a_n = 0}^\infty \prod_{k = 1}^n (-1)^{a_k}
\frac{q_{\beta_k}^{-a_k (a_k - 1)/2}}{[a_k]_{q_{\beta_k}}!}(q_{\beta_k} - q_{\beta_k}^{-1})^{a_k} \\
&\qquad \times (F_{\beta_n}^{a_n} \cdots F_{\beta_1}^{a_1},g) \;
  (E_{\beta_n}^{a_n} \cdots E_{\beta_1}^{a_1}, f) q^{-(\mu, \nu)} \\
&= (\prod_{\alpha \in {\bf \Delta}^+} \exp_{q_\alpha^{-1}}(-(q_\alpha - q_\alpha^{-1})(E_\alpha \otimes F_\alpha)), f \otimes g) q^{-(\mu, \nu)}.
\end{align*}
Now, according to the proof of Proposition \ref{tauboreliso}, we have
\[
 f = \bra v^-| \bullet| w^- \ket, \qquad g = \bra v^+| \bullet |w^+ \ket 
\]
where $w^-$ has weight $-\nu$ and $w^+$ has weight $-\mu$.  Therefore
\[
f_{(1)} \otimes g_{(1)} (q^{-\sum_{i,j = 1}^N B_{ij}(H_i \otimes H_j)}, f_{(2)} \otimes g_{(2)}) = q^{-(\mu, \nu)} f \otimes g, 
\]
and we conclude that
\[
  r^{-1}(f \otimes g) = (\R^{-1}, f \otimes g).
\]
This completes the proof.
\end{proof} 

According to Theorem \ref{Rmatrixformula}, we have in particular $ \hat{\Delta}^{\cop}(X) = \R \hat{\Delta}(X) \R^{-1} $ for 
all $ X \in U_q(\mathfrak{g}) \subset \M(\D(G_q)) $. Since $ \M(\D(G_q)) $ can be viewed 
as an algebraic completion of $ U_q(\mathfrak{g}) $ this makes the latter share most properties of Hopf algebras 
which are quasitriangular in the precise sense of definition \ref{defquasitriangular}. In the literature this is often phrased by 
saying that $ U_q(\mathfrak{g}) $ is quasitriangular, however note that the universal $ R $-matrix $ \R $ from Theorem \ref{Rmatrixformula} is 
not contained in $ U_q(\mathfrak{g}) \otimes U_q(\mathfrak{g}) $. 

Let us next introduce $ l $-functionals.

\begin{definition}
	\label{def:l-functionals}
For $ f \in \Poly(G_q) $, we define the $l$-functionals $ l^\pm(f) $ on $ \Poly(G_q) $ by 
$$
(l^+(f), h) = (\R, h \otimes f), \qquad (l^-(f), h) = (\R^{-1}, f \otimes h). 
$$
\nomenclature{$l^\pm$}{$l$-functionals $l^\pm:\Poly(G_q)\to U_q(\lie{g})$}
\end{definition}

Since $ \R $ is contained in $ \M(\D(G_q) \otimes \D(G_q)) $ these functionals 
can be naturally viewed as elements of $ \M(\D(G_q)) $. We can strengthen this observation as follows. 

\begin{lemma} \label{lfunctionallemma}
For any $ f \in \Poly(G_q) $ the $ l $-functionals $ l^\pm(f) $ are contained in $ U_q(\mathfrak{g}) \subset \M(\D(G_q)) $.  
\end{lemma}

\begin{proof}
		Consider the case where $g = \bra v | \bullet | w \ket $, $h = \bra v' | \bullet | w' \ket$ are matrix coefficients with $v\in V(\lambda)^*$, $w\in V(\lambda)$, $v'\in V(\lambda')^*$, $w'\in V(\lambda')$ for some $\lambda, \lambda'\in\weights^+$.  Suppose moreover that $w$ and $w'$ have weight $\mu$ and $\nu$, respectively.  As in the proof of Theorem \ref{Rmatrixformula}, we calculate
	\[
	 (q^{\sum_{i,j} B_{ij}(H_i \otimes H_j)}, h \otimes g) = q^{(\mu,\nu)} \epsilon(h)\epsilon(g) = (\epsilon(g)K_\mu,h).
	\]
	From this we easily deduce $ l^+(f) \in U_q(\mathfrak{g}) $ for any $ f \in \Poly(G_q) $. 
	
	The proof for $ l^-(f) $ is analogous. 
\end{proof}

The conclusion of Lemma \ref{lfunctionallemma} is one of the reasons to work with the \emph{simply connected} 
version of $ U_q(\mathfrak{g}) $ with generators $ K_\mu $ for all $ \mu \in \weights $. The corresponding assertion fails if one considers the version of 
the quantized universal enveloping algebra with Cartan part containing only elements $ K_\mu $ for $ \mu \in \roots $. 

We record the following basic properties of $ l $-functionals. 

\begin{lemma} \label{lfunctionalproperties}
The maps $ l^\pm: \Poly(G_q) \rightarrow U_q(\mathfrak{g}) $ are Hopf algebra homomorphisms. Explicitly, we have 
\begin{align*}
l^\pm(fg) &= l^\pm(f) l^\pm(g) \\
\hat{\Delta}(l^\pm(f)) &= l^\pm(f_{(1)}) \otimes l^\pm(f_{(2)}) \\
\hat{S}(l^\pm(f)) &= l^\pm(S(f)) 
\end{align*}
for all $ f, g \in \Poly(G_q) $. 
\end{lemma} 
\begin{proof}
For $ h \in \Poly(G_q) $ we compute
\begin{align*}
(l^+(fg), h) &= (\R, h \otimes fg) \\
&= ((\id \otimes \hat{\Delta})(\R), h \otimes g \otimes f) \\
&= (\R_{13} \R_{12}, h \otimes g \otimes f) \\
&= (\R, h_{(1)} \otimes f) (\R, h_{(2)} \otimes g) \\
&= (l^+(f) l^+(g), h),
\end{align*}
and similarly 
\begin{align*}
(\hat{\Delta}(l^+(f)), g \otimes h) 
&= ((\hat{\Delta} \otimes \id)(\R), g \otimes h \otimes f) \\
&= (\R_{13} \R_{23}, g \otimes h \otimes f) \\
&= (\R, g \otimes f_{(1)}) (\R, h \otimes f_{(2)}) \\
&=  (l^+(f_{(1)}) \otimes l^+(f_{(2)}), g \otimes h).
\end{align*}
The relation concerning the antipodes follows from these two formulas; explicitly we calculate
$$
(\hat{S}(l^+(f)), h) = ((\hat{S} \otimes \id)(\R), h \otimes f) = ((\id \otimes \hat{S}^{-1})(\R), h \otimes f) = (l^+(S(f)), h).   
$$
The assertions for $ l^- $ are obtained in a similar way. \end{proof}

\subsection{The locally finite part of $ U_q(\mathfrak{g}) $} 
\label{sec:locally_finite_part}

In this subsection we discuss some results on the locally finite part of $ U_q(\mathfrak{g}) $. For more details we refer to 
Section 5.3 in \cite{Josephbook} and \cite{Baumannseparation}, \cite{Calderothesis}. Throughout we assume that $ q = s^{2L} \in \mathbb{K}^\times $ 
is not a root of unity.

Recall that $ U_q(\mathfrak{g}) $ is a left module over itself with respect to the adjoint action given by 
$$
\ad(X)(Y) = X \rightarrow Y = X_{(1)} Y \hat{S}(X_{(2)}) 
$$
\label{nom:left_adjoint_action2}
for $ X, Y \in U_q(\mathfrak{g}) $. The space of invariant elements of $ U_q(\mathfrak{g}) $ with respect to the 
adjoint action is the centre $ ZU_q(\mathfrak{g}) $ of $ U_q(\mathfrak{g}) $. 
\nomenclature{$ZU_q(\lie{g})$}{centre of $U_q(\lie{g})$}%

The locally finite part $ FU_q(\mathfrak{g}) $
\nomenclature{$FU_q(\lie{g})$}{locally finite part of $U_q(\lie{g})$}%
is the sum of all finite dimensional $ U_q(\mathfrak{g}) $-submodules 
of $ U_q(\mathfrak{g}) $ with respect to the adjoint action. In contrast to the classical case it turns out that 
$ FU_q(\mathfrak{g}) $ is a proper subspace of $ U_q(\mathfrak{g}) $.

\begin{lemma}
	\label{lem:FUq_coideal}
	The locally finite part $FU_q(\lie{g})$ is a left coideal subalgebra of $U_q(\lie{g})$, that is, it is a sublagebra and satisfies $\Delta(FU_q(\lie{g})) \subset U_q(\lie{g}) \otimes FU_q(\lie{g})$.
\end{lemma}

\begin{proof}
	The fact that $FU_q(\lie{g})$ is a subalgebra follows from the fact that the adjoint action is compatible with multiplication.  Now let $Y\in FU_q(\lie{g})$ and $X\in U_q(\lie{g})$.  Note that
	\begin{align*}
	  Y_{(1)} \otimes (X\rightarrow Y_{(2)})
	   &= S(X_{(1)})X_{(2)} Y_{(1)} S(X_{(5)})X_{(6)} \otimes X_{(3)}Y_{(2)}S(X_{(4)}) \\
	   &= (S(X_{(1)})\otimes 1)\Delta(X_{(2)}\rightarrow Y) (X_{(3)}\otimes 1).
	\end{align*}  
	Since $U_q(\lie{g})\rightarrow Y$ is finite dimensional it follows that $\Delta(Y) \in U_q(\lie{g}) \otimes FU_q(\lie{g})$.
\end{proof}

Recall that the anti-automorphism $ \tau $ of $ U_q(\mathfrak{g}) $ defined in 
Lemma \ref{deftau} satisfies
\begin{align*}
\tau(X \rightarrow Y) &= \tau(X_{(1)} Y \hat{S}(X_{(2)})) \\
&= \tau(\hat{S}(X_{(2)})) \tau(Y) \tau(X_{(1)}) \\
&= \tau(\hat{S}(X_{(2)})) \tau(Y) \hat{S}(\tau(\hat{S}(X_{(1)}))) \\
&= \tau(\hat{S}(X)) \rightarrow \tau(Y) 
\end{align*}
for all $ X, Y \in U_q(\mathfrak{g}) $. In particular, the involution $ \tau $ preserves the locally finite part $ FU_q(\mathfrak{g}) $ 
of $ U_q(\mathfrak{g}) $.

More generally, for an arbitrary $ U_q(\mathfrak{g}) $-module $ M $ we shall denote by $ FM \subset M $
\nomenclature{$FM$}{locally finite part of a $U_q(\lie{g})$-module $M$}%
its locally finite part, that is, 
the subspace consisting of all elements $ m $ such that $ U_q(\mathfrak{g}) \cdot m $ is finite dimensional. 
We shall say that the action of $ U_q(\mathfrak{g}) $ on $ M $ is locally finite if $ FM = M $. Of course, this is the 
case in particular if $ M = V(\mu) $ for $ \mu \in \weights^+ $. 

Let us also define the coadjoint action of $ U_q(\mathfrak{g}) $ on $U_q(\mathfrak{g})^*$ by 
\[
 (X\rightarrow f)(Y) = f(\hat{S}(X_{(1)}) Y X_{(2)})
\]
\label{nom:coadjoint_action}
for $ X, Y \in U_q(\mathfrak{g}) $ and $ f \in U_q(\mathfrak{g})^* $. 
If $f=\bra v | \bullet | w \ket \in \Poly(G_q)$ is a matrix coefficient, then
\[
 X \rightarrow \bra v | \bullet | w \ket = \bra X_{(1)}\cdot v | \bullet | X_{(2)}\cdot w \ket.
\]
Therefore the coadjoint action preserves $\Poly(G_q)$  and is locally finite on $\Poly(G_q)$. Note that for $f\in\Poly(G_q)$,
\[
  (X\rightarrow f)(Y) = (Y,X\rightarrow f)
\]
for all $X,Y\in U_q(\lie{g})$.

We define a linear map $ J: U_q(\mathfrak{g}) \rightarrow U_q(\mathfrak{g})^* $ by 
$$ 
 J(X)(Y) = \kappa(\hat{S}^{-1}(Y), X) 
$$ 
\nomenclature{$J$}{linear isomorphism $ J: FU_q(\mathfrak{g}) \rightarrow \Poly(G_q)$}%
where $ \kappa $ is the quantum Killing form, see Definition \ref{defquantumkilling}. 
Using $ J $ we can determine the structure of $ FU_q(\mathfrak{g}) $, compare Section 7.1 in \cite{Josephbook} and \cite{Baumannseparation}. 

\begin{theorem} \label{adjointduality}
The linear map $ J $ defines an isomorphism from $ FU_q(\mathfrak{g}) $ onto $ \Poly(G_q) $
compatible with the adjoint and coadjoint actions, respectively. Moreover 
$$
FU_q(\mathfrak{g}) = \bigoplus_{\mu \in \weights^+} U_q(\mathfrak{g}) \rightarrow K_{2 \mu}  
$$
as a subspace of $ U_q(\mathfrak{g}) $. 
\end{theorem}

\begin{proof} We follow the discussion in \cite{Calderothesis}. According to Theorem \ref{lemdrinfeldnondegenerate} the 
map $ J: U_q(\mathfrak{g}) \rightarrow U_q(\mathfrak{g})^* $ is injective. Moreover $ J $ intertwines the 
adjoint and coadjoint actions since 
by the $\ad$-invariance of $\kappa$ from Proposition \ref{killinginvariance} we obtain
\begin{align*}
(X \rightarrow J(Y))(Z) &= J(Y)(\hat{S}(X_{(1)}) Z X_{(2)}) \\
&= \kappa(\hat{S}^{-1}(X_{(2)}) \hat{S}^{-1}(Z) X_{(1)}, Y) \\
&= \kappa(\hat{S}^{-1}(X) \rightarrow \hat{S}^{-1}(Z), Y) \\
&= \kappa(\hat{S}^{-1}(Z), X \rightarrow Y) \\
&= J(X \rightarrow Y)(Z)
\end{align*}
for all $ X, Y, Z \in U_q(\mathfrak{g}) $.  

Let $ \mu \in \weights^+ $ and let $ w_0 \in W $ be the longest element of the Weyl group. Using the definition of the quantum Killing form 
we check that $ J(K_{-2 w_0 \mu}) $ vanishes on monomials $ X K_\nu \hat{S}(Y) $ for which $ X \in U_q(\mathfrak{n}_+) $ or 
$ Y \in U_q(\mathfrak{n}_-) $ is contained in the kernel of $ \hat{\epsilon} $. 
With this in mind, if $ v^\mu $ is a lowest weight vector of $ V(\mu) $ and 
$ v^\mu_*$ is the highest weight vector in $ V(\mu)^* $ such that
$ v^\mu_*(v^\mu) = 1 $, one can confirm that
$$ 
J(K_{-2 w_0 \mu}) = \bra v^\mu_*|\bullet |v^\mu \ket
$$
by comparing both sides on $ U_q(\mathfrak{h}) $.

Since $ V(\mu) $ is irreducible, the vector $ \bra v^\mu_*|\bullet |v^\mu \ket $ 
is cyclic for the coadjoint action. We conclude that $ J $ induces an 
isomorphism $ U_q(\mathfrak{g}) \rightarrow K_{-2 w_0 \mu} \cong \End(V(\mu))^* $. 
Since the spaces $ \End(V(\mu))^* $ form a direct sum in $ \Poly(G_q) \subset U_q(\mathfrak{g})^* $ it follows that the sum of 
the spaces $ U_q(\mathfrak{g}) \rightarrow K_{-2 w_0 \mu} $ 
is a direct sum, and the resulting subspace of $ FU_q(\mathfrak{g}) $ is isomorphic to $ \Poly(G_q) $ via $ J $. 

To finish the proof it suffices to show that $ J(FU_q(\mathfrak{g})) $ is contained in $ \Poly(G_q) $. For 
this let $ f \in J(FU_q(\mathfrak{g})) $ be arbitrary. 
Consider 
the left and right regular
actions of $U_q(\mathfrak{g})$ on  $U_q(\mathfrak{g})^*$ given by
\[
   (X \cdot f)(Y) = f(YX), \qquad (f\cdot X)(Y) = f(XY) 
\]
for $ X, Y \in U_q(\mathfrak{g}) $. 
Inspecting the definition of the quantum Killing form shows that $ U_q(\mathfrak{b}_-) \cdot f $ 
and $ f \cdot U_q(\mathfrak{b}_+) $ are finite dimensional subspaces of $ U_q(\mathfrak{g})^* $. 
Moreover $ f $ is contained in the locally finite part of $ U_q(\mathfrak{g})^* $ with respect to the coadjoint action, 
so the relation $ X \cdot f = X_{(2)} \rightarrow (f \cdot X_{(1)}) $ shows that $ U_q(\mathfrak{b}_+) \cdot f $ is finite 
dimensional as well. 
Hence $ U_q(\mathfrak{g}) \cdot f $ is finite dimensional. 

From the definition of the quantum Killing form $ \kappa $ it follows that the linear 
functional $ J(Y K_\lambda \hat{S}(X)) \in U_q(\mathfrak{g})^* $ for $ X \in U_q(\mathfrak{n}_+)_\nu, Y \in U_q(\mathfrak{n}_-) $ 
has weight $ \nu - \frac{1}{2} \lambda \in \frac{1}{2}\weights \subset \mathfrak{h}_q^* $ with respect to the left regular action.
In particular, $ J(U_q(\mathfrak{g})) $ is a weight module with respect to the left regular action, and 
therefore the same holds for $ U_q(\mathfrak{g}) \cdot f $. According to 
Corollary \ref{coruqgfiniterepdecomp}, 
we conclude
that $ U_q(\mathfrak{g}) \cdot f $ is isomorphic to a finite direct sum of 
modules $ V(\mu_j) $ for $ \mu_1, \dots, \mu_r \in \weights_q^+ $, and since $\weights_q^+\cap \frac{1}{2}\weights = \weights^+$, we have in fact $ \mu_1, \dots, \mu_r \in \weights^+ $.
In particular, we get $ X \cdot f = 0 $ for 
any $ X \in I(\mu_1) \cap \cdots \cap I(\mu_r) $, the intersection of the kernels of the irreducible representations 
corresponding to $ \mu_1, \dots, \mu_r $. 
Therefore $0 = (X\cdot f)(1) = f(X)$ for all such $ X $, which means that $ f $ can be viewed 
as a linear functional on the quotient $ U_q(\mathfrak{g})/(I(\mu_1) \cap \cdots \cap I(\mu_r)) \cong \End(V(\mu_1)) \oplus 
\cdots \oplus \End(V(\mu_r)) $. In other words, we have $ f \in \Poly(G_q) $ as desired. This finishes the proof. \end{proof} 

As a consequence of Theorem \ref{adjointduality} we obtain the following properties of $ FU_q(\mathfrak{g}) $. 

\begin{prop} \label{locallyfinitecartan}
The intersection of $ U_q(\mathfrak{h}) $ and $ FU_q(\mathfrak{g}) $ is linearly spanned by the elements $ K_\lambda $ for $ \lambda \in 2 \weights^+ $. 
\end{prop} 

\begin{proof} Let us first show 
$$ 
E_j^r \rightarrow K_\lambda = \prod_{k = 1}^r (1 - q^{(\lambda - (k - 1) \alpha_j, \alpha_j)}) E_j^r K_{\lambda - r \alpha_j}, 
$$
using induction on $ r \in \mathbb{N} $. 
We have 
\begin{align*}
E_j \rightarrow K_\lambda &= [E_j, K_\lambda] K_j^{-1} = (E_j - K_\lambda E_j K_{-\lambda}) K_{\lambda - \alpha_j} 
= (1 - q^{(\lambda, \alpha_j)}) E_j K_{\lambda - \alpha_j}, 
\end{align*} 
and in the inductive step we compute 
\begin{align*}
&E_j^r \rightarrow K_\lambda = [E_j, E_j^{r - 1} \rightarrow K_\lambda] K_j^{-1} \\ 
&= \prod_{k = 1}^{r - 1} (1 - q^{(\lambda - (k - 1) \alpha_j, \alpha_j)})
(E_j E_j^{r - 1} K_{\lambda - (r - 1) \alpha_j} K_{- \alpha_j} - E_j^{r - 1} K_{\lambda - (r - 1)\alpha_j} E_j K_{-\alpha_j}) \\ 
&= \prod_{k = 1}^{r - 1} (1 - q^{(\lambda - (k - 1) \alpha_j, \alpha_j)}) 
(E_j^r K_{\lambda - r \alpha_j} - q^{(\lambda - (r - 1) \alpha_j, \alpha_j)} E_j^r K_{\lambda - r \alpha_j}) \\ 
&= \prod_{k = 1}^r (1 - q^{(\lambda - (k - 1) \alpha_j, \alpha_j)}) E_j^r K_{\lambda - r \alpha_j}, 
\end{align*} 
so that the formula indeed holds. 

Since $ q $ is not a root of unity it follows that $ E_j^r \rightarrow K_\lambda = 0 $ for some $ r \in \mathbb{N} $ iff 
$ (\alpha_j, \lambda) = n_j(\alpha_j, \alpha_j) $ for some integer $ n_j \geq 0 $. 
Writing $ \lambda = \lambda_1 \varpi_1 + \cdots + \lambda_N \varpi_N $ this amounts to 
$$ 
\lambda_j = (\alpha_j^\vee, \lambda) = \frac{2}{(\alpha_j, \alpha_j)} (\alpha_j, \lambda) = 2 n_j. 
$$
We conclude that $ K_\lambda \in FU_q(\mathfrak{g}) $ implies $ \lambda \in 2 \weights^+ $. 
In fact, the argument shows more generally that an element of $ U_q(\mathfrak{h}) $ which is contained 
in $ FU_q(\mathfrak{g}) $ is necessarily a linear combination of elements $ K_\lambda $ with $ \lambda \in 2 \weights $.

Conversely, according to Theorem \ref{adjointduality} the elements $ K_\lambda $ for $ \lambda \in 2 \weights^+ $ are indeed 
contained in $ FU_q(\mathfrak{g}) $. This finishes the proof. \end{proof} 

Observe that Proposition \ref{locallyfinitecartan} shows in particular that $ FU_q(\mathfrak{g}) $ in not a Hopf subalgebra 
of $ U_q(\mathfrak{g}) $. 

The following result illustrates that the subalgebra $ FU_q(\mathfrak{g}) \subset U_q(\mathfrak{g}) $ is reasonably large. 
For more precise information in this direction see Section 7.1 in \cite{Josephbook}. 

\begin{lemma} \label{uqgequalsfuqguqh}
We have $ U_q(\mathfrak{g}) = FU_q(\mathfrak{g}) U_q(\mathfrak{h}) = U_q(\mathfrak{h}) FU_q(\mathfrak{g}) $.  
\end{lemma} 

\begin{proof} We first note that $ K_\mu FU_q(\mathfrak{g}) K_{-\mu} = \ad(K_\mu)(FU_q(\mathfrak{g})) = FU_q(\mathfrak{g}) $ for all $ \mu \in \weights $,  
and hence $ FU_q(\mathfrak{g}) U_q(\mathfrak{h}) = U_q(\mathfrak{h}) FU_q(\mathfrak{g}) $. It follows that 
the latter is a subalgebra. 

In order to finish the proof it therefore suffices to show that the generators $ E_i, F_i $ are 
contained in $ FU_q(\mathfrak{g}) U_q(\mathfrak{h}) = U_q(\mathfrak{h}) FU_q(\mathfrak{g}) $. 
From the computation in the proof of Proposition \ref{locallyfinitecartan} we conclude that
$$
E_j \rightarrow K_{2 \rho} = (1 - q^{(2\rho, \alpha_j)}) E_j K_{2 \rho - \alpha_j} = (1 - q_{j}^2) E_j K_{2 \rho - \alpha_j} 
$$
is contained in $ FU_q(\mathfrak{g}) $, and hence $ E_j \in FU_q(\mathfrak{g}) U_q(\mathfrak{h}) $. 
In the same way we see that 
$$
F_j \rightarrow K_{2 \rho} = - K_{2 \rho} F_j + F_j K_{2 \rho} = (1 - q^{-(2 \rho, \alpha_j)}) E_j K_{2 \rho} = (1 - q_{j}^{-2}) E_j K_{2 \rho} 
$$
is contained in $ FU_q(\mathfrak{g}) $, and hence $ F_j \in FU_q(\mathfrak{g}) U_q(\mathfrak{h}) $. \end{proof}

Recall next the definition of the $ l $-functionals $ l^\pm(f) $ for $ f \in \Poly(G_q) $ from Definition \ref{def:l-functionals}. 
We define a linear map $ I: \Poly(G_q) \rightarrow U_q(\mathfrak{g}) $ by 
$$
I(f) = l^-(f_{(1)}) l^+(S(f_{(2)})). 
$$
\nomenclature{$I$}{linear isomorphism $I:\Poly(G_q) \rightarrow FU_q(\mathfrak{g})$}%
For the following result compare \cite{Baumannseparation}. 

\begin{prop} \label{adjointdualityplus} 
The map $ I $ induces an isomorphism $ \Poly(G_q) \rightarrow FU_q(\mathfrak{g}) $, and this isomorphism is the inverse 
of the isomorphism $ J: FU_q(\mathfrak{g}) \rightarrow \Poly(G_q) $ constructed above. 
\end{prop} 

\begin{proof} By the definition of the $ l $-functionals we have 
$$ 
(I(f), g) = (\R_{12}^{-1} \R_{21}^{-1}, f \otimes g) 
$$
for $ f, g \in \Poly(G_q) $.  
Here we write $ \R_{12} = \R $ and $ \R_{21} = \sigma(\R) $ where $ \sigma $ is the flip map. Using 
\begin{align*}
\hat{\Delta}(X) \R_{12}^{-1} \R_{21}^{-1}  &= \R_{12}^{-1} \hat{\Delta}^{\cop}(X) \R_{21}^{-1} = \R_{12}^{-1} \R_{21}^{-1} \hat{\Delta}(X)
\end{align*}
for all $ X \in U_q(\mathfrak{g}) $ we compute 
\begin{align*}
(I(X \rightarrow f), g) &= (\hat{S}(X_{(1)}), f_{(1)}) (I(f_{(2)}), g) (X_{(2)}, f_{(3)}) \\ 
&= (\hat{S}(X_{(1)}), f_{(1)}) (\R_{12}^{-1} \R_{21}^{-1}, f_{(2)} \otimes g) (X_{(2)}, f_{(3)}) \\ 
&= (\hat{S}(X_{(1)}), f_{(1)}) (\R_{12}^{-1} \R_{21}^{-1} \hat{\Delta}(X_{(2)}), f_{(2)} \otimes g_{(1)}) (\hat{S}(X_{(3)}), g_{(2)}) \\ 
&= (\hat{S}(X_{(1)}), f_{(1)}) (\hat{\Delta}(X_{(2)}) \R_{12}^{-1} \R_{21}^{-1}, f_{(2)} \otimes g_{(1)}) (\hat{S}(X_{(3)}), g_{(2)}) \\ 
&= (X_{(1)}, g_{(1)}) (\R_{12}^{-1} \R_{21}^{-1}, f \otimes g_{(2)}) (\hat{S}(X_{(2)}), g_{(3)}) \\ 
&= (X_{(1)}, g_{(1)}) (I(f), g_{(2)}) (\hat{S}(X_{(2)}), g_{(3)}) = (X \rightarrow I(f), g) 
\end{align*}
for $ f, g \in \Poly(G_q) $. Hence $ I $ is $ U_q(\mathfrak{g}) $-linear with respect to the coadjoint and adjoint actions, respectively. This 
means in particular that the image of $ I $ is contained in $ FU_q(\mathfrak{g}) $ since the coadjoint action on $ \Poly(G_q) $ is locally finite. 

Let $ v^\mu \in V(\mu) $ be a lowest weight vector and $ v^\mu_* \in V(\mu)^* $ a highest weight vector such that $ v^\mu_*(v^\mu) = 1 $. 
From the explicit description of the $ R $-matrix in Theorem \ref{Rmatrixformula} we see that 
$ l^-(\bra v^\mu_*| \bullet| v \ket) = (v^\mu_*, v) K_{-w_0 \mu} $ for all $ v \in V(\mu) $.
Similarly, we obtain $ l^+(\bra v^*|\bullet| v^\mu \ket) = (v^*, v^\mu) K_{w_0 \mu} $ for all $ v^* \in V(\mu)^* $. 
Taking a basis $ v_1, \dots, v_n $ of $ V(\mu) $ with dual basis $ v_1^*, \dots, v_n^* \in V(\mu)^* $ this shows 
\begin{align*}
I(\bra v^\mu_*|\bullet| v^\mu \ket) &= \sum_{j = 1}^n l^-(\bra v^\mu_*|\bullet| v_j\ket) l^+(S(\bra v_j^*|\bullet| v^\mu \ket)) \\
&= \sum_{j = 1}^n (v^\mu_*, v_j) (v_j^*, v^\mu) K_{-2 w_0 \mu} = K_{-2 w_0 \mu}. 
\end{align*}
Combining this with the results from Theorem \ref{adjointduality} we conclude that $ I $ is indeed the inverse of $ J $. 
In particular $ I $ is an isomorphism. \end{proof} 

Let $ g, f \in \Poly(G_q) $. Using that $ I $ is compatible with the adjoint and coadjoint actions and properties of 
$ l $-functionals we compute 
\begin{align*}
I(g) &I(f) = I(g) l^-(f_{(1)}) l^+(S(f_{(2)})) \\ 
&= l^-(f_{(2)}) (\hat{S}^{-1}(l^-(f_{(1)})) \rightarrow I(g)) l^+(S(f_{(3)})) \\ 
&= l^-(f_{(2)}) I(\hat{S}^{-1}(l^-(f_{(1)})) \rightarrow g) l^+(S(f_{(3)})) \\ 
&= l^-(f_{(2)}) l^-((\hat{S}^{-1}(l^-(f_{(1)})) \rightarrow g)_{(1)}) 
l^+(S((\hat{S}^{-1}(l^-(f_{(1)})) \rightarrow g)_{(2)})) l^+(S(f_{(3)})) \\
&= l^-(f_{(2)} (\hat{S}^{-1}(l^-(f_{(1)})) \rightarrow g)_{(1)}) l^+(S(f_{(3)}(\hat{S}^{-1}(l^-(f_{(1)})) \rightarrow g)_{(2)})) \\
&= I(f_{(2)} \hat{S}^{-1}(l^-(f_{(1)})) \rightarrow g). 
\end{align*}
In a similar way we obtain 
\begin{align*}
I(f) &I(g) = l^-(f_{(1)}) l^+(S(f_{(2)})) I(g) \\ 
&= l^-(f_{(1)}) (l^+(S(f_{(3)})) \rightarrow I(g)) l^+(S(f_{(2)})) \\ 
&= l^-(f_{(1)}) I(l^+(S(f_{(3)})) \rightarrow g) l^+(S(f_{(2)})) \\ 
&= l^-(f_{(1)}) l^-((l^+(S(f_{(3)})) \rightarrow g)_{(1)}) l^+(S(l^+(S(f_{(3)})) \rightarrow g)_{(2)}) l^+(S(f_{(2)})) \\ 
&= l^-(f_{(1)} (l^+(S(f_{(3)})) \rightarrow g)_{(1)}) l^+(S(f_{(2)} (l^+(S(f_{(3)})) \rightarrow g)_{(2)})) \\ 
&= I(f_{(1)} \hat{S}(l^+(f_{(2)})) \rightarrow g). 
\end{align*}
for all $ f, g \in \Poly(G_q) $. 

In particular, the map $ I: \Poly(G_q) \rightarrow FU_q(\mathfrak{g}) $ is not an algebra homomorphism. 
However, assume that $ g $ is an element of the invariant part $ \Poly(G_q)^{G_q} $ of $ \Poly(G_q) $ with respect to 
the coadjoint action. Then the computation above 
shows that $ I $ satisfies the multiplicativity property 
$$ 
  I(fg) = I(f) I(g)  = I(g)I(f)
$$ 
for all $ f \in \Poly(G_q) $. In particular, we obtain the following.

\begin{prop}
	\label{prop:I_centre}
  The map $ I $ induces an algebra isomorphism between $ \Poly(G_q)^{G_q} $ and $ ZU_q(\mathfrak{g}) $. 
\end{prop}

\subsection{The centre of $ U_q(\mathfrak{g}) $ and the Harish-Chandra homomorphism} \label{seczuqg}

In this section we describe the structure of the centre $ ZU_q(\mathfrak{g}) $ of $ U_q(\mathfrak{g}) $. 
For more details we refer to \cite{Baumanncenter}. Throughout we assume that $ q = s^{L} \in \mathbb{K}^\times $ is not a root of unity. 

Recall from Proposition \ref{prop:I_centre} that we have an algebra isomorphism $ \Poly(G_q)^{G_q} \cong ZU_q(\mathfrak{g}) $.  Therefore, we start by studying the algebra $\Poly(G_q)^{G_q}$.

\begin{definition}
	\label{def:quantum_trace}
	Let $V$ be a finite dimensional integrable $U_q(\mathfrak{g})$-module $V$.  The quantum trace $t_V \in \Poly(G_q)$ is defined by
	\[
	 t_V(X) = \tr_V(XK_{-2\rho})
	\]
	\nomenclature{$t_V$}{quantum trace of a finite dimensional integrable module $V$}%
	\nomenclature{$t_\mu$}{quantum trace of the finite dimensional integrable module $V(\mu)$}%
	for all $X \in U_q(\mathfrak{g})$, where $\tr_V$ denotes the unnormalized trace on $\End(V)$.  If $V=V(\mu)$ for some $\mu\in\weights^+$, we will write $t_\mu$ for $t_{V(\mu)}$.
\end{definition}

\begin{lemma}
	\label{lem:quantum_traces}
	The set of quantum traces $t_\mu$ for $\mu\in\weights^+$ is a linear basis of $\Poly(G_q)^{G_q}$.
\end{lemma}

\begin{proof}
	First, one checks that an element $f \in \Poly(G_q)$ is invariant for the coadjoint action if and only if
	\[
	 (\hat{S}^2(X)Y , f) = (YX,f)
	\]
	for all $X,Y \in U_q(\mathfrak{g})$. 
	Recall from Lemma \ref{lem:S2}, that $\hat{S}^2(X) = K_{2\rho}XK_{-2\rho}$.  Using this, and $f\cdot K_{2\rho} = (K_{2\rho}, f_{(1)})f_{(2)}$, we see that the above condition is equivalent to
	\[
	 (XY,f\cdot K_{2\rho}) = (K_{2\rho}XY,f) = (K_{2\rho}YX,f) = (YX, f\cdot K_{2\rho})
	\]
	for all $X,Y \in U_q(\mathfrak{g})$.  This holds if and only if $f\cdot K_{2\rho}$ is a linear combination of traces $\tr_{V(\mu)}$. 
\end{proof}

Note that the quantum traces satisfy
\begin{align*}
t_{V\oplus W} &= t_V + t_W,
&t_{V\otimes W} &= t_Vt_W,
\end{align*}
for all finite dimensional integrable modules $V$, $W$.
It follows that $\Poly(G_q)^{G_q}$ is isomorphic to the polynomial algebra 
$ \mathbb{K}[t_{\varpi_1}, \dots, t_{\varpi_N}] $.

For $ \mu \in \weights^+ $ let us define the Casimir element $ z_\mu \in U_q(\mathfrak{g}) $
\nomenclature{$z_\mu$}{Casimir element in $ZU_q(\lie{g})$ associated to $\mu\in\weights^+$}%
by 
$ z_\mu = I(t_\mu) $. As a direct consequence of Proposition \ref{prop:I_centre} and Lemma \ref{lem:quantum_traces} we obtain the following result. 

\begin{theorem} \label{uqgcenter}
The centre $ ZU_q(\mathfrak{g}) $ of $ U_q(\mathfrak{g}) $ is canonically isomorphic to the polynomial ring $ \mathbb{K}[z_{\varpi_1}, \dots, z_{\varpi_N}] $. 
A linear basis of $ ZU_q(\mathfrak{g}) $ is given by the Casimir elements $ z_\mu $ for $ \mu \in \weights^+ $.  
\end{theorem}

Recall that the Harish-Chandra map is the linear map $ \P: U_q(\mathfrak{g}) \rightarrow U_q(\mathfrak{h}) $ defined by $ \P = \hat{\epsilon} \otimes \id \otimes \hat{\epsilon} $ with respect to the triangular decomposition 	$ U_q(\mathfrak{g}) \cong U_q(\mathfrak{n}_-) \otimes U_q(\mathfrak{h}) \otimes U_q(\mathfrak{n}_+) $.  Our next goal is the quantum analogue of  Harish Chandra's Theorem, which says that $\P$ restricts to an algebra isomorphism of $ZU_q(\lie{g})$ onto a subalgebra of $U_q(\lie{h})$.
	
In order to specify this subalgebra, we will need to introduce a shifted Weyl group action on $U_q(\lie{h})$.  Note that there is a standard action of $W$ on $U_q(\lie{h})$ by algebra automorphisms, defined by
\[
 w (K_\mu) = K_{w\mu},
\]
\label{nom:unshifted_action_on_Uqh}%
for $w\in W$ and $\mu\in\weights$.  The shifted version of this action is defined as follows.

\begin{definition}
	\label{def:shifted_action_on_Ugh}
	Let $\gamma$ be the algebra automorphism of $U_q(\lie{h})$ defined by
	\[
	  \gamma (K_\mu) = q^{(\rho,\mu)} K_\mu \qquad \qquad (\mu\in \weights).
	\]
	\nomenclature{$\gamma$}{automorphism of $U_q(\lie{h})$}%
	The shifted action of $W$ on $U_q(\lie{h})$ is the conjugate of the standard action above by $\gamma$.  Specifically, the shifted action of $w\in W$ is given by
	\[
	  w.K_\mu = q^{(\rho,w\mu-\mu)}K_{w\mu}.
	\]
	\label{nom:shifted_action_on_Uqh}%
  for $w\in W$.
\end{definition}

Note that the automorphism $\gamma$ of $U_q(\lie{h})$ induces a map $\gamma^*$ on the group of characters of $U_q(\lie{h})$ and hence on the parameter space $\lie{h}_q^*$.  Specifically, for the character $\chi_\lambda$ associated to $\lambda \in \lie{h}_q^*$ we have
\[
 \gamma^*(\chi_\lambda) (K_\mu) = \chi_\lambda (\gamma(K_\mu)) = q^{(\lambda+\rho,\mu)}
\]
for all $\mu\in\weights$, so that $\gamma^*$ acts on $\lie{h}_q^*$ by translation by $\rho$.

We also introduce the algebra homomorphism $\tau : \mathbb{K}[\weights] \to U_q(\lie{h})$ defined by
\[
  \tau (\mu) = K_{2\mu}
\]
\nomenclature{$\tau$}{Algebra homomorphism $\tau:\mathbb{K}[\weights] \to U_q(\lie{h})$, $\tau (\mu) = K_{2\mu}$}%
for $\mu\in\weights$.  Note that the subalgebra $\im(\tau) \subset U_q(\lie{h})$ is preserved by the shifted Weyl group action.

The following result appears in various frameworks in \cite{TanisakiHarishChandra, deConciniKacrootof1a, JLlocalfiniteness}, see also \cite{Rossorepgq}.

\begin{theorem} \label{Rossohc}
The Harish-Chandra map $ \P: U_q(\mathfrak{g}) \rightarrow U_q(\mathfrak{h}) $ 
restricts to an algebra isomorphism $ ZU_q(\mathfrak{g}) \cong \im(\tau)^W $, the fixed point subalgebra of $\im(\tau)$ under the shifted Weyl group action.

\end{theorem}

\begin{proof} 
It is equivalent to prove the map $\gamma^{-1}\circ\P$ defines an algebra isomorphism of $ZU_q(\lie{g})$ onto the fixed point subalgebra of $\im(\tau)$ under the standard action of $W$.

Let $\mu\in\weights^+$ and consider the central element  $z_\mu$.  Fix a basis $ e^\mu_1, \dots, e^\mu_n $ of weight vectors for $ V(\mu) $ with dual basis $ e^1_\mu, \dots, e^n_\mu $ and put $u^\mu_{ij} = \bra e^i_\mu| \bullet | e^\mu_j \ket $ so that the quantum trace is
\[
t_\mu = \sum_{j=1}^n (K_{-2\rho},u^\mu_{jj})u^\mu_{jj}.
\]
Then
$$ 
z_\mu = I(t_\mu) = \sum_{j,k} (K_{-2 \rho}, u^\mu_{jj}) l^-(u^\mu_{jk}) l^+(S(u^\mu_{kj})). 
$$
Note that $ l^-(u^\mu_{rs}) $ is contained in $ U_q(\mathfrak{b}_-) $, and similarly $ l^+(S(u^\mu_{rs})) $ is contained in $ U_q(\mathfrak{b}_+) $. 
Inspecting the explicit form of the universal $ R $-matrix in Theorem \ref{Rmatrixformula} we get 
$$ 
\P(z_\mu) = \sum_{\nu \in \weights(V(\mu))} q^{(-2 \rho, \nu)} K_{-\nu} K_{-\nu} = \sum_{\nu \in \weights(V(\mu))} q^{( \rho, -2 \nu)} K_{-2 \nu}, 
$$ 
where $ \nu $ runs over the set of weights $ \weights(V(\mu)) $ of $ V(\mu) $, counted with multiplicities.  We obtain
\[
 \gamma^{-1}\circ\P(z_\mu) = \sum_{\nu \in \weights(V(\mu))} K_{-2 \nu},
\]
which is obviously in $\im(\tau)$ and invariant under the standard action of $W$.  Moreover, for any $\lambda,\mu\in\weights^+$ we have
\begin{align*}
 (\gamma^{-1}\circ\P(z_\mu)) (\gamma^{-1}\circ\P(z_\lambda))
   &= \sum_{\substack{ \nu' \in \weights(V(\mu)), \\ \nu'' \in \weights(V(\lambda))}} K_{-2 \nu'}K_{-2\nu''} \\
   &= \hspace{-2ex} \sum_{\nu \in \weights(V(\mu)\otimes V(\lambda))} \hspace{-2ex} K_{-2 \nu} \\
   & = \gamma^{-1}\circ\P(z_\mu z_\lambda),
\end{align*}
so $\gamma\otimes\P$ is an algebra homomorphism.

Using an induction argument with respect to the partial ordering on $\weights^+$, one checks that 
the elements $ \gamma^{-1}\circ\P(z_\mu) $ for $ \mu \in \weights^+ $ form a basis of the fixed point subalgebra of $ \im(\tau) $ with respect to the standard action. 
Combining this with Theorem \ref{uqgcenter} finishes the proof.
\end{proof} 

Let us illustrate these considerations in the case of $ U_q(\mathfrak{sl}(2, \mathbb{K})) $. In this case 
one can check directly that the Casimir element $ z_{1/2} $ associated with the fundamental representation $ V(1/2) $ is, up to a scalar 
multiple, given by  
$$
C = FE + \frac{qK^2 + q^{-1}K^{-2}}{(q - q^{-1})^2}. 
$$
\nomenclature{$C$}{fundamental Casimir element for $ U_q(\mathfrak{sl}(2, \mathbb{K})) $}%
Let us verify that $ C $ is in the centre of $ U_q(\mathfrak{sl}(2, \mathbb{K})) $. We compute 
\begin{align*}
EC &= [E, F]E + FEE + \frac{q^{-1}K^2E + qK^{-2}E}{(q - q^{-1})^2} \\
&= FEE + \frac{(q - q^{-1}) K^2E + (q^{-1} - q) K^{-2}E}{(q - q^{-1})^2} + \frac{q^{-1}K^2E + qK^{-2}E}{(q - q^{-1})^2} \\
&= FEE + \frac{qK^2E + q^{-1}K^{-2}E}{(q - q^{-1})^2} = CE, 
\end{align*}
similarly we get 
\begin{align*}
CF &= F[E, F] + FFE + \frac{q^{-1}FK^2 + qFK^{-2}}{(q - q^{-1})^2} \\
&= FFE + \frac{(q - q^{-1}) FK^2 + (q^{-1} - q) FK^{-2}}{(q - q^{-1})^2} + \frac{q^{-1}FK^2 + qFK^{-2}}{(q - q^{-1})^2} \\
&= FFE + \frac{qFK^2 + q^{-1}FK^{-2}}{(q - q^{-1})^2} = FC, 
\end{align*}
and the relation $ CK = KC $ is obvious. 

Let us now return to the general theory. Consider a Verma module $ M(\lambda) $ for $ \lambda \in \mathfrak{h}^*_q $. If $ Z \in ZU_q(\mathfrak{g}) $ 
and $ v_\lambda \in M(\lambda) $ is the highest weight vector, then $ Z \cdot v_\lambda $ is again a highest weight vector, and therefore must be of the 
form $ \xi_\lambda(Z) v_\lambda $,
for some linear map $\xi_\lambda: ZU_q(\lie{g}) \to \mathbb{K}$.	
Note that
\[
 Z\cdot v_\lambda = \P(Z)\cdot v_\lambda = \chi_\lambda(\P(Z)) v_\lambda,
\]
where $\chi_\lambda$ is the character given by $\chi_\lambda(K_\nu) = q^{(\lambda, \nu)} $ for all $\nu\in\weights$, and therefore
\[
 \xi_\lambda = \chi_\lambda \circ \P.
\]  
Moreover, 
since $ v_\lambda $ is a cyclic vector for $ M(\lambda) $ it follows that $ Z $ acts by multiplication by $ \xi_\lambda(Z) $ on all of $ M(\lambda) $. 
It is easy to check that $ \xi_\lambda(YZ) = \xi_\lambda(Y) \xi_\lambda(Z) $ for all $ Y,Z \in ZU_q(\mathfrak{g}) $, 
which means that $ \xi_\lambda $ is a character. 

\begin{definition} 
\label{def:central_character}
Let $ \lambda \in \mathfrak{h}^*_q $. 
The character $ \xi_\lambda: ZU_q(\mathfrak{g}) \rightarrow \mathbb{K} $ defined above is called the central character associated with $ \lambda $. 
\nomenclature{$\xi_\lambda$}{central character associated to $\lambda\in\lie{h}_q^*$}%
\end{definition}

From the previous considerations it follows that if $ M(\mu) \subset M(\lambda) $ is a submodule then $ \xi_\mu = \xi_\lambda $. 

We shall now introduce the notion of linkage, which is designed to analyze the structure of Verma modules.  For a classical complex semisimple Lie algebra $\lie{g}$, two Verma modules $M(\lambda)$ and $M(\mu)$ with highest weights $\lambda,\mu\in\lie{h}^*$ have the same central character if and only if $\lambda$ and $\mu$ are in the same orbit of the shifted Weyl group action, in which case they are called $W$-linked,  see \cite{HumphreysO}.  For $U_q(\lie{g})$ the group $W$ needs to be enlarged slightly to take into account the existence of elements in $\lie{h}^*_q$ of order $2$.

\begin{definition}
	\label{def:order_2}
	We define $\mathbf{Y}_q = \{\zeta\in\lie{h}_q^* \mid 2\zeta=0\}$.
\end{definition}

If characteristic of the ground field $\mathbb{K}$ is not equal to $2$, there is a canonical isomorphism
\[
 \mathbf{Y}_q \cong \roots^\vee / 2\roots^\vee
\]
which associates to the class of $\alpha^\vee \in \roots^\vee$ the character $\zeta_{\alpha^\vee} \in \lie{h}_q^*$ defined by
\[
  q^{(\zeta_{\alpha^\vee},\mu)} = (-1)^{(\alpha^\vee,\mu)},
\]
for $\mu\in\roots$.  Here we are using the $q$-exponential notation for $\lie{h}_q^*$ from Definition \ref{def:hq-star}.  This formula is well-defined since the character $\zeta_{\alpha^\vee}$ is trivial if $\alpha^\vee\in 2\roots^\vee$.  Moreover, if $\zeta\in\lie{h}_q^*$ is any element of order two then we have $q^{(\zeta,\mu)}=\pm1$ for all $\mu\in\weights$, which implies that there exists $\alpha^\vee\in\roots^\vee$ such that $\zeta = \zeta_{\alpha^\vee}$, as given by the above formula.

In the case where the ground field is $\mathbb{C}$, where we make the identification $\lie{h}_q^* = \lie{h}^* / i\hbar^{-1} \roots^\vee$, the natural identification is  $\mathbf{Y}_q \cong \half i\hbar^{-1} \roots^\vee / i\hbar^{-1} \roots^\vee$.
\label{nom:order_2_C}%

Note that if $\mathbb{K}$ is an exponential field and if there exists $z\in\mathbb{K}$ such that $q^z = -1$ then we have $\zeta_{\alpha^\vee} = z\alpha^\vee$.  
Here we are using the canonical map $\lie{h}^* \to \lie{h}^*_q$ described in Subsection \ref{sec:weights}.

Now we define the extended Weyl group, see Section 8.3.2 in \cite{Josephbook}.  Note that the action of the Weyl group $W$ on $\lie{h}_q^*$ restricts to $\mathbf{Y}_q$.

\begin{definition} 
	\label{def:extended_Weyl_group}
	The extended Weyl group $ \hat{W} $ is defined as the semidirect product 
	$$ 
	\hat{W} = \mathbf{Y}_q \rtimes W
	$$ 
	\nomenclature[o$W$22]{$\hat{W}$}{extended Weyl group $\hat{W} = \mathbf{Y}_q \rtimes W$}%
	with respect to the action of $ W $ on $ \mathbf{Y}_q $. 
\end{definition}

Observe that the extended Weyl group is a finite group. Explicitly, the group structure of $ \hat{W} $ is 
$$
(\zeta, v)(\eta, w) = (\zeta + v \eta, vw)
$$
for all $ \zeta,\eta \in \mathbf{Y}_q $ and $ v, w \in W $.

We define two actions of $\hat{W}$ on $\mathfrak{h}_q^*$ by
\[
  (\zeta,w)\lambda = w\lambda+\zeta
\]
\label{nom:extended_Weyl_group_action}%
and
\begin{align*}
(\zeta,w).\lambda & = w.\lambda+\zeta = w(\lambda+\rho) - \rho +\zeta,
\end{align*}
\label{nom:shifted_extended_Weyl_group_action}%
for $\lambda\in\mathfrak{h}_q^*$.  The latter is called the shifted action of $\hat{W}$ on $\mathfrak{h}_q^*$.

\begin{definition} \label{defwlinkage}
We say that $ \mu, \lambda \in \mathfrak{h}^*_q $ are $ \hat{W} $-linked if $ \hat{w} . \lambda = \mu $ 
for some $ \hat{w} \in \hat{W} $. 
\end{definition}

Notice that $ \hat{W} $-linkage is an equivalence relation on $ \mathfrak{h}^*_q $. 
This should be compared to the usual notion of $W$-linkage, see \cite{HumphreysO}. 
Indeed, the two notions coincide if we restrict our attention to integral weights.

\begin{lemma}
	Two integral weights $ \mu, \lambda \in \weights \subset \mathfrak{h}^*_q $ are $ \hat{W} $-linked iff they are $ W $-linked.
\end{lemma}

\begin{proof}
	If $\mu,\lambda\in\weights$ with $(\zeta,w).\lambda = w.\lambda +\zeta = \mu$ then $\zeta  \in\weights \subset \lie{h}_q^*$.  But then for all $\nu\in\weights$, we have $q^{(\zeta,\nu)}\in s^\mathbb{Z}$ and also
	$
	 q^{(\zeta,\nu)}  =  \pm1 .
	$
	Since $q$ is not a root of unity, this implies that $\zeta=0$.
\end{proof}

We shall now state and prove the following analogue of Harish-Chandra's Theorem, compare Section 1.10 in \cite{HumphreysO}. 

\begin{theorem} \label{chilinkage}
Two elements $ \mu, \lambda \in \mathfrak{h}^*_q $ are $ \hat{W} $-linked iff $ \xi_\mu = \xi_\lambda $. 
\end{theorem}

\begin{proof} Assume first that $ \lambda $ and $ \mu = \hat{w} . \lambda $ are linked, where $ \hat{w} = (\zeta, w) \in \hat{W} $. 
In order to show $ \xi_{\hat{w} . \lambda} = \xi_{\lambda} $ it suffices due to Theorem \ref{Rossohc} to show $ \chi_{\hat{w} . \lambda} = \chi_\lambda $, where 
both sides are viewed as characters on the fixed point subalgebra $\im(\tau)^W$ under the shifted Weyl group action. If 
$$ 
K = \sum_{\nu \in \weights} c_\nu K_{2 \nu} 
= \sum_{\nu \in \weights} c_\nu q^{(\rho, 2w\nu - 2 \nu)} K_{2 w \nu} 
$$ 
is contained in $ \im(\tau)^W $ then we obtain
\begin{align*}
\chi_\mu(K) = \chi_{w . \lambda + \zeta }(K) &= \sum_{\nu \in \weights} c_\nu q^{(\rho, 2w\nu - 2 \nu)} q^{(w . \lambda + \zeta , 2 w \nu)} \\
&= \sum_{\nu \in \weights} c_\nu q^{(\rho, 2w\nu - 2 \nu)} q^{(w . \lambda, 2 w \nu)} \\
&= \sum_{\nu \in \weights} c_\nu q^{(\rho, 2w\nu - 2\nu)} q^{(w(\lambda + \rho) - \rho, 2 w \nu)} \\
&= \sum_{\nu \in \weights} c_\nu q^{(\rho,- 2\nu)} q^{(w(\lambda + \rho), 2 w \nu)} \\ 
&= \sum_{\nu \in \weights} c_\nu q^{(\rho,- 2\nu)} q^{(\lambda + \rho, 2 \nu)} \\
&= \sum_{\nu \in \weights} c_\nu q^{(\lambda, 2\nu)} = \chi_\lambda(K) 
\end{align*} 
as desired. 
 
Let us now assume that $ \mu $ and $ \lambda $ are not $ \hat{W} $-linked. 
This implies that the shifted Weyl group orbits $W.(2\mu)$ and $W.(2\lambda)$ are disjoint.
Let us write $ 2 W . \mu \cup 2 W . \lambda = \{2 \eta_1, \dots, 2 \eta_n \} \subset \mathfrak{h}^*_q $. 
Each element of this set corresponds to a character $ \weights \rightarrow \mathbb{K}^\times $ 
sending $ \nu \in \weights $ to $ q^{(2 \eta_i, \nu)} $. 
Therefore Artin's Theorem on linear independence of characters allows us to find elements $ \nu_1, \dots, \nu_n \in \weights $ such that the matrix $ (q^{(2 \eta_i, \nu_j)}) $ is invertible, 
and hence scalars $ c_1, \dots, c_n \in \mathbb{K} $ such that 
$$ 
\sum_{j = 1}^n c_j q^{(\eta_i, 2\nu_j)} = 
\begin{cases} 
1 & \text{ if } \eta_i \in W . \mu \\
0 & \text{ if } \eta_i \in W . \lambda. 
\end{cases}
$$

If we write $ K = \sum_j c_j K_{2 \nu_j} $ then $ L = \frac{1}{|W|} \sum_{w \in W} w . K $ is contained in $ \im(\tau)^W $. Moreover, since 
\begin{align*}
\chi_{\eta_i}(w . K) &= \sum_j c_j \chi_{\eta_i}(w . K_{2 \nu_j}) \\
&= \sum_j c_j q^{(\rho, 2 w \nu_j - 2 \nu_j)} q^{(\eta_i, 2 w \nu_j)} \\
&= \sum_j c_j q^{(w^{-1}(\eta_i + \rho) - \rho, 2 \nu_j)} = \sum_j c_j q^{(w^{-1} . \eta_i, 2 \nu_j)}  = \chi_{w^{-1} . \eta_i}(K)
\end{align*}
we have $ \chi_{\eta_i}(L) = 1 $ for $ \eta_i \in W. \mu $ and $ \chi_{\eta_i}(L) = 0 $ for $ \eta_i \in W. \lambda $. 
Taking the preimage $ Z \in ZU_q(\mathfrak{g}) $ of $ L \in \im(\tau)^W $ under the Harish-Chandra isomorphism 
therefore gives $ \xi_\lambda(Z) = 1 $ and $ \xi_\mu(Z) = 0 $. Hence $ \xi_\mu \neq \xi_\lambda $. \end{proof} 

Let us conclude this subsection with the following result on the structure of characters of $ ZU_q(\mathfrak{g}) $. 

\begin{prop} \label{uqgcentralcharacters}
Assume that $ \mathbb{K} $ is an algebraically closed field. Then any character $ \chi: ZU_q(\mathfrak{g}) \rightarrow \mathbb{K} $ is of the form $ \xi_\lambda $ for some $ \lambda \in \mathfrak{h}^*_q $. 
\end{prop} 

\begin{proof} Using Theorem \ref{Rossohc} we can identify $ ZU_q(\mathfrak{g}) $ with the $ W $-invariant part $ A^W $ of the Laurent polynomial 
ring $ A = \im \tau = \mathbb{K}[K_{2 \varpi_1}^{\pm 1}, \dots, K_{2 \varpi_N}^{\pm 1}] $, where $W$ is acting by the shifted Weyl group action. It suffices to show that any 
character $ \chi: A^W \rightarrow \mathbb{K} $ is of the form $ \chi = \chi_\lambda $ for $ \lambda \in \mathfrak{h}^*_q $.

Since $ W $ is a finite group the ring extension $ A^W \subset A $ 
is integral, that is, each element of $ A $ is the root of a monic polynomial with coefficients in $ A^W $. In fact, given $ f \in A $ 
we may consider $ p(x) = \prod_{w \in W} (x - w . f) $. 

To proceed further we need to invoke some commutative algebra for integral ring extensions. 
More precisely, due to Theorem 5.10 in \cite{AtiyahMcdonald} the maximal ideal $ \ker(\chi) $ is of the form $ \ker(\chi) = A^W \cap \mathfrak{p} $ 
for some prime ideal $ \mathfrak{p} $ of $ A $, and Corollary 5.8 in \cite{AtiyahMcdonald} implies that $ \mathfrak{p} $ is a maximal ideal. 
Since $ \mathbb{K} $ is algebraically closed we conclude that $ \mathfrak{p} = \ker(\eta) $ for some character $ \eta: A \rightarrow \mathbb{K} $ 
extending $ \chi $.

Using once again that $\mathbb{K}$ is algebraically closed, we can find $ \lambda \in \mathfrak{h}^*_q $ such that $ \eta(K_\mu) = q^{(\mu, \lambda)} = \chi_\lambda(K_\mu) $ 
for all elements $ K_\mu \in A $. Clearly this yields $ \chi = \chi_\lambda $ as desired. \end{proof}

\subsection{Noetherianity} 

In this section we show that $ U_q(\mathfrak{g}) $ and some related algebras are Noetherian. We assume throughout that $ q = s^L \in \mathbb{K}^\times $ 
is not a root of unity. 

\subsubsection{Noetherian algebras} 

Recall that an algebra $ A $ is called left (right) Noetherian if it satisfies the ascending chain condition (ACC) for 
left (right) ideals. The ACC says that any ascending chain $ I_1 \subset I_2 \subset I_3 \subset \cdots $ of left (right) ideals 
of $ A $ becomes eventually constant, that is, satisfies $ I_n = I_{n + 1} = I_{n + 2} = \cdots $ for some $ n \in \mathbb{N} $. 
The algebra $ A $ is called Noetherian if it is both left and right Noetherian. 

We will have to work with graded and filtered algebras. Let $ A $ be an algebra and $ P $ an additive abelian semigroup, 
and assume that there exists a direct sum decomposition 
$$
A = \bigoplus_{\mu \in P} A_\mu 
$$
into linear subspaces $ A_\mu $ such that $ 1 \in A_0 $ and $ A_\mu A_\nu \subset A_{\mu + \nu} $ for all $ \mu, \nu \in P $. In this case 
we call $ A $ a graded algebra. A left (right) ideal $ I $ of a graded algebra $ A $ is called graded if $ I = \bigoplus_{\mu \in P} I_\mu $ 
where $ I_\mu = A_\mu \cap I $. 

Assume now in addition that the abelian semigroup $ P $ is ordered, that is, equipped with an order relation $ \leq $ such that $ \mu \leq \nu $ 
implies $ \mu + \lambda \leq \nu + \lambda $ for all $ \mu, \nu, \lambda \in P $. 
We shall also assume that $ P $ has cancellation, so that $ \mu + \lambda = \nu + \lambda $ for some $ \lambda \in P $ implies $ \mu = \nu $. 
A $ P $-filtration on an algebra $ A $ is a family of linear subspaces $ \F^\mu(A) \subset A $
\label{nom:filtration}%
indexed by $ \mu \in P $ 
such that $ 1 \in \F^0(A) $, 
$ \F^\mu(A) \subset \F^\nu(A) $ if $ \mu \leq \nu $, 
$$
\bigcup_{\mu \in P} \F^\mu(A) = A,  
$$
and $ \F^\mu(A) \F^\nu(A) \subset \F^{\mu + \nu}(A) $ for all $ \mu, \nu \in P $. In this case we also say that $ A $ is a filtered algebra. 
If $ A $ is a filtered algebra, then the associated graded algebra is 
$$
gr_\F(A) = \bigoplus_{\mu \in P} \F^\mu(A)/\F^{< \mu}(A), 
$$
where 
\begin{align*}
\F^{< \mu}(A) = \sum_{\nu < \mu} \F^{\nu}(A), 
\end{align*}
and we write $ \nu <\mu $ if $ \nu \leq \mu $ and $ \nu \neq \mu $. 
The associated graded algebra is indeed a graded algebra in a natural way with graded subspaces $ gr_\F(A)_\mu = \F^\mu(A)/\F^{< \mu}(A) $ 
labelled by $ P $. 

Recall that $P$ is called well-founded if every subset of $P$ contains a minimal element.  This is equivalent to the descending chain condition: any infinite sequence $\mu_1 \geq \mu_2 \geq \cdots$ in $P$ is eventually constant.

We say that a filtration $ \F^\mu(A) $ of an algebra $ A $ is locally bounded below if for each nonzero element $ a \in A $ the set of all $ \mu \in P $ 
with $ a \in \F^\mu(A) $ has a minimal element.
Recall also that an algebra $ A $ is called a domain if it has no zero-divisors, that is, if $ ab = 0 $ in $ A $ implies $ a = 0 $ or $ b = 0 $. 

Let us record the following basic fact, compare Appendix I.12 in \cite{BrownGoodearllectures} and Chapter V in \cite{Jacobsonliealgebras}. 

\begin{lemma} \label{filteredgradednoether}
Let $ A $ be a filtered algebra with a filtration $ (\F^\mu(A))_{\mu \in P} $ with respect to an abelian semigroup $ P $.
Then the following holds. 
\begin{bnum}
\item[a)] Suppose $ P $ is well-founded. If the associated graded algebra $ gr_\F(A) $ is left (right) Noetherian then $ A $ is left (right) Noetherian. 
\item[b)] Suppose $P$ is totally ordered and the filtration is locally bounded below. If $ gr_\F(A) $ is a domain then $ A $ is a domain. 
\end{bnum}
\end{lemma} 

\begin{proof} $ a) $ We shall consider left ideals only, the proof for right ideals is completely analogous. 

Note that if $ I \subset A $ is a left ideal then setting $ \F^\mu(I) = I \cap \F^\mu(A) $ for $ \mu \in P $ defines a filtration of $ I $, with associated graded
$$
gr_\F(I) = \bigoplus_{\mu \in P} gr_\F(I)_\mu = \bigoplus_{\mu \in P} \F^\mu(I)/\F^{< \mu}(I).
$$
There is a canonical inclusion of $\F^\mu(I)/\F^{<\mu}(I)$ into $\F^\mu(A)/\F^{<\mu}(A)$, and in this way $gr_\F(I)$ becomes a graded left ideal of $gr_\F(A)$.

Assume now that $ I_1 \subset I_2 \subset \cdots $ is an increasing chain of left ideals in $ A $. Then the associated graded left ideals  fit into an ascending chain $gr_\F(I_1) \subset gr_\F(I_2) \subset \cdots$. By assumption, 
this chain stabilizes eventually, so that $gr_\F(I_n) = gr_\F(I_{n+1}) = \cdots$ for some $ n \in \mathbb{N} $. 

We shall show that in general if $ I \subset J $ are left ideals of $ A $ such that $ gr_\F(I) = gr_\F(J) $ then $ I = J $. This will clearly 
imply that our original chain $ I_1 \subset I_2 \subset \cdots $ is eventually constant.

Suppose $I \neq J$.  Since $P$ is well-founded, we can find $\mu\in P$ minimal such that $\F^\mu(I) \neq \F^\mu(J)$.  Choose $a\in \F^\mu(J)\setminus\F^\mu(I)$.  Since $gr_\F(I) = gr_\F(J)$, there exists $b\in \F^\mu(I)$ such that $b-a\in\F^{<\mu}(J)=\F^{<\mu}(I)$.  But then $a = b + (a-b) \in \F^\mu(I)$, contradicting the choice of $a$.  Therefore $I=J$ as claimed.

$ b) $ Assume $ a, b \in A $ are nonzero elements such that $ a b = 0 $. Since $ P $ is locally bounded below we may pick $ \mu, \nu \in P $ minimal 
such that $ a \in \F^\mu(A) $ and $ b \in \F^\nu(A) $. The corresponding cosets in $ gr_\F(A) $ multiply to zero, so that $ a + \F^{< \mu}(A) = 0 $ 
or $ b + \F^{< \nu}(A) = 0 $ in $ gr_\F(A) $. This means $ a \in \F^{< \mu}(A) $ or $ b \in \F^{< \nu}(A) $. Since $ P $ is totally ordered 
this implies $ a \in \F^\lambda(A) $ for some $ \lambda < \mu $ or $ b \in \F^\eta(A) $ for some $ \eta < \nu $, 
contradicting the choice of $ \mu $ and $ \nu $. Hence $ A $ is a domain. \end{proof}

We remark that in the proof of part $a)$ it suffices to assume that $gr_\F(A)$ is Noetherian for graded ideals.  It follows that a graded algebra $A$ is Noetherian if and only if it is Noetherian for graded ideals.  For we may equip $A$ with the filtration $ \F $ given by 
$$ 
\F^\mu(A) = \bigoplus_{\nu \leq \mu} A_\nu, 
$$ 
so that $ A \cong gr_\F(A) $ as graded algebras. Then the proof of Lemma \ref{filteredgradednoether} $a)$ implies that $ A $ is left Noetherian if $ gr_\F(A) = A $ 
is Noetherian for graded left ideals. Conversely, if $ A $ is left Noetherian then it is clearly also Noetherian for graded left ideals.

Let $ A $ be an algebra and let $ \theta \in \Aut(A) $ be an algebra automorphism. 
Given $ \theta $ we can form the semigroup crossed product $ A \rtimes_\theta \mathbb{N}_0 $,
which has $ A \otimes \mathbb{K}[t] $ as underlying vector 
space, with elements written as $ a \otimes t^m = a \rtimes t^m $, and multiplication given by 
$$
(a \rtimes t^m)(b \rtimes t^n) = a \theta^m(b) \rtimes t^{m + n}
$$
for $ a,b \in A $ and $ m,n \in \mathbb{N}_0 $. 

\begin{lemma} \label{noetherdomainlemma}
Let $ A $ be an algebra, let $ \theta \in \Aut(A) $ and let $ A \rtimes_\theta \mathbb{N}_0 $ be the associated semigroup crossed product. 
\begin{bnum} 
\item[a)] If $ A $ is left (right) Noetherian then $ A \rtimes_\theta \mathbb{N}_0 $ is left (right) Noetherian. 
\item[b)] If $ A $ is a domain then $ A \rtimes_\theta \mathbb{N}_0 $ is a domain. 
\end{bnum}
\end{lemma} 

\begin{proof} $ a) $ The argument is a variant of the Hilbert Basis Theorem, for a proof see Theorem 1.2.9 in \cite{MCRnoetherian}. 

$ b) $ Let $ a = \sum_{i = 0}^m a_i \rtimes t^i, b = \sum_{j = 0}^n b_j \rtimes t^j $ be nonzero elements of $ A \rtimes_\theta \mathbb{N}_0 $, 
and assume without loss of generality that $ a_m $ and $ b_n $ are nonzero. Observe that 
$$ 
ab = \sum_{i = 0}^m \sum_{j = 0}^n a_i \theta^i(b_j) \rtimes t^{i + j}.   
$$ 
Since $ \theta $ is an automorphism we have $ \theta^k(b_n) \neq 0 $ for all $ k \in \mathbb{N}_0 $. 
However, $ ab = 0 $ implies $ a_m \theta^m(b_n) = 0 $, contradicting our assumption that $ A $ is a domain. 
It follows that $ A \rtimes_\theta \mathbb{N}_0 $ is a domain. \end{proof} 

Let us call an algebra $ A $ a skew-polynomial algebra if it is generated by finitely many elements $ y_1, \dots, y_m $ and relations
$$
y_i y_j = q_{ij} y_j y_i,  
$$
where $ q_{ij} \in \mathbb{K}^\times $ are invertible scalars for all $ 1 \leq  i,j \leq m $ with $q_{ij}=q_{ji}^{-1}$  and $q_{ii}=1$ for all $i,j$.  

\begin{lemma} \label{skewpolynomialnoether}
Any skew-polynomial algebra is a Noetherian domain. 
\end{lemma}

\begin{proof} We use induction on the number of generators. For $ m = 0 $ we have that $ A = \mathbb{K} $, so that the assertion is obvious. 
Assume now that the claim has been proved if $ k < m $ for some $ m \in \mathbb{N} $, and let $ y_1, \dots, y_m $ be generators of $ A $ as above. 
Then the subalgebra $ B \subset A $ generated by $ y_1, \dots, y_{m - 1} $ is a skew-polynomial algebra and hence a Noetherian domain by the 
inductive hypothesis. 
Moreover, it is easy to see that $ \theta(y_j) = q_{mj} y_j $ defines an algebra automorphism of $ B $. 
one checks that the corresponding semigroup crossed $ B \rtimes_\theta \mathbb{N}_0 $ is isomorphic to $ A $.   
Hence the assertion follows from Lemma \ref{noetherdomainlemma}. \end{proof} 

The following result is Proposition I.8.17 in \cite{BrownGoodearllectures}, which we reproduce here for the convenience of the reader. 

\begin{prop} \label{genrelnoether}
Let $ A $ be an algebra generated by elements $ u_1, \dots, u_m $, and assume that there are scalars $ q_{ij} \in \mathbb{K}^\times $ and 
$ \alpha_{ij}^{st}, \beta_{ij}^{st} \in \mathbb{K} $ such that 
$$
u_i u_j - q_{ij} u_j u_i = \sum_{s = 1}^{j - 1} \sum_{t = 1}^m \alpha^{st}_{ij} u_s u_t + \beta^{st}_{ij} u_t u_s 
$$
for all $ 1 \leq j < i \leq m $. Then $ A $ is a Noetherian. 
\end{prop} 

\begin{proof} The idea is to construct a filtration $ \F^n(A) $ for $ n \in \mathbb{N}_0 $ of $ A $ such that $ gr_\F(A) $ is generated by 
elements $ y_1, \dots, y_m $ such that $ y_i y_j = q_{ij} y_j y_i $ for all $ 1 \leq i,j \leq m $. 

To this end let 
$$ 
d_i = 2^m - 2^{m - i}
$$
for $ i = 1, \dots, m $, so that $ d_1 < d_2 < \cdots < d_m $. For $ 1 < j < m $ we have 
$$
2^{m - j - 1} + 2^{m - j} = (1 + 2) \cdot 2^{m - j - 1} < 4 \cdot 2^{m - j - 1} < 2^{m - j + 1} + 1,  
$$ 
and hence 
\begin{align*}
d_{j - 1} + d_m &= 2^m - 2^{m - j + 1} + 2^m - 1 \\
&= 2^{m + 1} - (2^{m - j + 1} + 1) < 2^{m + 1} - (2^{m - j - 1} + 2^{m - j}) = d_{j + 1} + d_j. 
\end{align*}
We thus get for $ i > j > s $ and any $ t $ the relation 
$$
d_s + d_t \leq d_{j - 1} + d_m < d_{j + 1} + d_j \leq d_i + d_j,  
$$
using the fact that $ d_s \leq d_{j - 1} $ and $ d_t \leq d_m $ in the first step. 

Define $ \F^0(A) = \mathbb{K} 1 $ and let $ \F^d(A) $ for $ d > 0 $ be the linear subspace spanned by all 
monomials $ u_{i_1} \cdots u_{i_r} $ such that $ d_{i_1} + \cdots + d_{i_r} \leq d $. 
It is straightforward to check that this defines an $\mathbb{N}_0$-filtration on $ A $. Let $ y_i $ be the coset of $ u_i $ in $ gr_\F(A) $. 
Each nonzero homogeneous component of $ gr_\F(A) $ 
is spanned by the cosets of those monomials $ u_{i_1} \cdots u_{i_r} $ such that $ d_{i_1} + \cdots + d_{i_r}  = d $. In particular, 
the elements $ y_1, \dots, y_m $ generate $ gr_\F(A) $ as an algebra. Since $ u_i u_j - q_{ij} u_j u_i $ for $ i > j $ 
is a linear combination of monomials $ u_s u_t $ and $ u_t u_s $ with $ s < j $, and all such monomials have filtration degree strictly 
smaller that $ d_i + d_j $, we see that 
$$ 
y_i y_j - q_{ij} y_j y_i = 0 
$$ 
in $ gr_\F(A) $. That is, $ gr_\F(A) $ is a quotient of a skew-polynomial algebra, and according to Lemma \ref{skewpolynomialnoether} this means that $ gr_\F(A) $ 
is Noetherian. Due to Lemma \ref{filteredgradednoether} this implies that $ A $ itself is a Noetherian. \end{proof}

\subsubsection{Noetherianity of $ U_q(\mathfrak{g}) $} 

In this subsection we use the PBW-basis to show that $ U_q(\mathfrak{g}) $ and a few other algebras we have encountered in our study of quantized 
universal enveloping algebras so far are Noetherian. These results are originally due to de Concini-Kac, see Section 1 in \cite{deConciniKacrootof1}. 

\begin{prop} \label{uqgnoetheriandomain}
The algebras $ U_q(\mathfrak{n}_\pm), U_q(\mathfrak{b}_\pm) $ and $ U_q(\mathfrak{g}) $ are Noetherian domains. 
\end{prop} 

\begin{proof} We will prove the claim for $ A = U_q(\mathfrak{g}) $, reproducing the argument in I.6 of \cite{BrownGoodearllectures}. The other cases are treated 
in a similar manner. Let us remark that the characteristic zero assumption made in \cite{BrownGoodearllectures} is not needed. 

We have the PBW-basis vectors 
$$
X({\bf t}, \mu, {\bf s}) = F_{\beta_n}^{t_n} \cdots F_{\beta_1}^{t_1} K_\mu E_{\beta_n}^{s_n} \cdots E_{\beta_1}^{s_1}
$$
where $ {\bf t} = (t_1, \dots, t_n), {\bf s} = (s_1, \dots, s_n) \in \mathbb{N}_0^n $. If $ \nu = \sum \nu_i \alpha_i \in \roots $ then we define 
the height of $ \nu $ by $ ht(\nu) = \nu_1 + \cdots + \nu_N \in \mathbb{Z} $. Let us define moreover the height of $ X({\bf t}, \mu, {\bf s}) $ by 
$$
ht(X({\bf t}, \mu, {\bf s})) = (t_1 + s_1) ht(\beta_1) + \cdots + (t_N + s_N) ht(\beta_N), 
$$
\label{nom:ht-PBW}%
so that the height of this vector is the difference of the heights of its graded pieces in $ U_q(\mathfrak{n}_\pm) $. 
We define the total degree of a PBW-basis vector by 
$$
d(X({\bf t}, \mu, {\bf s})) = (s_n, \dots, s_1, t_n, \dots, t_1, ht(X({\bf t}, \mu, {\bf s}))) \in \mathbb{N}_0^{2n + 1}. 
$$
Let us give the additive semigroup $ P = \mathbb{N}_0^{2n + 1} $ the lexicographical ordering from right to left, so that 
$ e_1 < e_2 < \cdots < e_{2n + 1} $, where $ e_j \in \mathbb{N}_0^{2n + 1} $ denotes the $ j $-th standard basis vector. 
This turns $ P $ into a well-ordered abelian semigroup, in particular $ P $ is totally ordered and well-founded. 
Moreover let $ \F^{\bf m}(A) \subset A $ for $ {\bf m} \in \mathbb{N}_0^{2n + 1} $ be the linear span of all PBW-basis vectors such 
that $ d(X({\bf t}, \mu, {\bf s})) \leq (m_1, \dots, m_{2n + 1}) = {\bf m} $ with respect to the lexicographical ordering. 

We claim that the spaces $ \F^{\bf m}(A) $ for $ {\bf m} \in P $ define a locally bounded below filtration of $ A $. Most of the required properties 
are straightforward, except that $ Y \in \F^{\bf m}(A), Z \in \F^{\bf n}(A) $ implies $ YZ \in \F^{\bf m + n}(A) $. 
The main point to check here is what happens when monomials in the $ E_{\beta_j} $ or $ F_{\beta_j} $ need to be reordered. Note that the formula 
obtained in Proposition \ref{uqnplusqcommutationrelations} shows that $ E_{\beta_i} E_{\beta_j} - q^{(\beta_i, \beta_j)} E_{\beta_j} E_{\beta_i} $ 
for $ i < j $ is a sum of monomials of the same height as $ E_{\beta_i} E_{\beta_j} $. According to the definition of the order in $ P $, all these monomials 
have strictly smaller total degree than $ E_{\beta_j} E_{\beta_i} $. 

Using this observation we see that the associated graded algebra $ gr_\F(A) $ is the algebra with generators 
$ E_{\beta_1}, \dots, E_{\beta_n}, K_\mu, F_{\beta_1}, \dots, F_{\beta_n} $ and relations 
\begin{align*}
K_\mu K_\nu &= K_\nu K_\mu, \\
K_\mu E_{\beta_i} &= q^{(\mu, \beta_i)} E_{\beta_i} K_\mu, \\
K_\mu F_{\beta_i} &= q^{-(\mu, \beta_i)} F_{\beta_i} K_\mu, \\
E_{\beta_i} F_{\beta_j} &= F_{\beta_j} E_{\beta_i}, \\
E_{\beta_i} E_{\beta_j} &= q^{(\beta_i, \beta_j)} E_{\beta_j} E_{\beta_i}, \\
F_{\beta_i} F_{\beta_j} &= q^{(\beta_i, \beta_j)} F_{\beta_j} F_{\beta_i}
\end{align*}
for all $ \mu, \nu \in \weights $ and $ 1 \leq i \leq n $ or $ 1 \leq i < j \leq n $, respectively. Therefore $ gr_\F(A) $ is a skew-polynomial ring in the 
generators $ E_{\beta_1}, \dots, E_{\beta_n}, F_{\beta_1}, \dots, F_{\beta_n} $ and $ K_{\varpi_1}, \dots, K_{\varpi_N} $. As such 
it is a Noetherian domain by Lemma \ref{skewpolynomialnoether}, and hence the same holds for $ A $ by Lemma \ref{filteredgradednoether}. \end{proof}

\subsubsection{Noetherianity of $ \Poly(G_q) $} 

We show here that the coordinate algebra $ \Poly(G_q) $ is Noetherian, again following the proof in Section I.8 of \cite{BrownGoodearllectures}. 
The techniques developed in this context will be used again in the next subsection to obtain Noetherianity of $ FU_q(\mathfrak{g}) $. 

We need some preparations. For any $ \lambda \in \weights^+ $ we fix a linear basis of weight vectors $ e^\lambda_1, \dots, e^\lambda_m $ for $ V(\lambda) $ 
with dual basis $ e_\lambda^1, \dots, e_\lambda^m \in V(\lambda)^* $. We write $ \epsilon_j $ for the weight of $ e^\lambda_j $, and we may assume 
without loss of generality that the vectors are ordered in a non-ascending order, so that $ \epsilon_i > \epsilon_j $ implies $ i < j $. 
Let us denote by $ u^\lambda_{ij} = \bra e_\lambda^i| \bullet| e^\lambda_j \ket $ the corresponding matrix coefficients. 

\begin{prop} \label{ogqcommutation}
Let $ \mu, \nu \in \weights^+ $ and $ 1 \leq i,j \leq m, 1 \leq k,l \leq n $. 
With the notations as above, there are scalars $ \alpha^{ijkl}_{rs}, \beta^{ijkl}_{uv} \in \mathbb{K} $ such that 
\begin{align*}
q^{-(\epsilon_j, \epsilon_l)} &u^\mu_{ij} u^\nu_{kl} + \sum_{r = 1}^{j - 1} \sum_{s = l + 1}^m \alpha^{ijkl}_{rs} u^\mu_{ir} u^\nu_{ks} 
= q^{-(\epsilon_i, \epsilon_k)} u^\nu_{kl} u^\mu_{ij} + \sum_{u = i + 1}^m \sum_{v = 1}^{k - 1} \beta^{ijkl}_{uv} u^\nu_{vl} u^\mu_{uj}
\end{align*}
in $ \Poly(G_q) $. Moreover $ \alpha^{ijkl}_{rs} = 0 $ unless $ \epsilon_r > \epsilon_j $ and $ \epsilon_s < \epsilon_l $. 
Similarly $ \beta^{ijkl}_{uv} = 0 $ unless $ \epsilon_u < \epsilon_i $ and $ \epsilon_v > \epsilon_k $. 
\end{prop}

\begin{proof} From general properties of the universal $ R $-matrix we obtain 
$$
f_{(1)} g_{(1)} (\R, f_{(2)} \otimes g_{(2)}) = (\R, f_{(1)} \otimes g_{(1)}) g_{(2)} f_{(2)} 
$$
for all $ f, g \in \Poly(G_q) $. Inspecting the description of the universal $ R $-matrix in Theorem \ref{Rmatrixformula} 
for $ f = u^\mu_{ij} $ and $ g = u^\nu_{kl} $ yields the desired formula. \end{proof} 

\begin{theorem} \label{ogqnoetherian}
The algebra $ \Poly(G_q) $ is Noetherian. 
\end{theorem}

\begin{proof} Note that the modules $ V(\varpi_1), \dots, V(\varpi_N) $ generate all irreducible finite dimensional weight modules of $ U_q(\mathfrak{g}) $ 
in the sense that any $ V(\mu) $ for $ \mu \in \weights^+ $ is a submodule of a suitable tensor product of the $ V(\varpi_j) $. Let us pick bases of weight 
vectors $ e^k_1, \dots, e^k_{n_k} $ for each $ V(\varpi_k) $ with dual bases $ e_k^1, \dots, e_k^{n_k} $, and write $ \epsilon_j $ for 
the weight of $ e^k_j $. The matrix elements $ u^k_{ij} = \bra e_k^i| \bullet | e^k_j \ket $ for $ k = 1, \dots, N $ and $ 1 \leq i,j \leq n_k $ 
generate $ \Poly(G_q) $ as an algebra. 

Let $ X $ be the collection of all these matrix elements $ u^k_{ij} $. We shall order the elements of $ X $ in a list $ u_1, \dots, u_m $ 
such that the following condition holds: For $ u_a = u^r_{ij}, u_b = u^s_{kl} $ we have $ b < a $ if either $ \epsilon_k < \epsilon_i $, 
or $ \epsilon_k = \epsilon_i $ and $ \epsilon_l > \epsilon_j $. 
According to Proposition \ref{ogqcommutation} we then obtain scalars $ q_{ij} \in \mathbb{K}^\times $ and $ \alpha^{st}_{ij}, \beta^{st}_{ij} \in \mathbb{K} $ 
such that 
$$
u_i u_j = q_{ij} u_j u_i + \sum_{s = 1}^{j - 1} \sum_{t = 1}^m \alpha_{ij}^{st} u_s u_t + \beta_{ij}^{st} u_t u_s 
$$
for all $ 1 \leq j < i \leq m $. Therefore Lemma \ref{genrelnoether} yields the claim. \end{proof} 

\subsubsection{Noetherianity of $ FU_q(\mathfrak{g}) $} 

In this subsection we discuss Noetherianity of the locally finite part of $ U_q(\mathfrak{g}) $. 
Joseph \cite{Josephbook} relies on tricky filtration arguments for this, we shall instead 
give a proof based on the link between $ FU_q(\mathfrak{g}) $ and $ \Poly(G_q) $ obtained in Theorem \ref{adjointduality}. 

\begin{theorem} \label{fuqgnoetherian}
The algebra $ FU_q(\mathfrak{g}) $ is Noetherian. 
\end{theorem}

\begin{proof}
From the calculations after Proposition \ref{adjointdualityplus} it follows that we can identify the opposite algebra of $ FU_q(\mathfrak{g}) $ with $ \Poly(G_q) = A $, 
if the latter is equipped with the multiplication 
\begin{align*}
f \bullet g &= f_{(2)} \hat{S}^{-1}(l^-(f_{(1)})) \rightarrow g 
\end{align*}
where 
$$ 
X \rightarrow g = (\hat{S}(X_{(1)}), g_{(1)}) g_{(2)} (X_{(2)}, g_{(3)}) 
$$ 
is the coadjoint action of $ U_q(\mathfrak{g}) $ on $ \Poly(G_q) $. 
Equivalently, according to Lemma \ref{lfunctionalproperties} and the definition of the $ l $-functionals we can write 
\begin{align*}
f \bullet g &= (l^-(f_{(2)}), g_{(1)}) f_{(3)} g_{(2)} (\hat{S}^{-1}(l^-(f_{(1)})), g_{(3)}) \\
&= (\R^{-1}, f_{(2)} \otimes g_{(1)}) f_{(3)} g_{(2)} (\R^{-1}, f_{(1)} \otimes S(g_{(3)})) \\ 
&= g_{(1)} f_{(2)} (\R^{-1}, f_{(3)} \otimes g_{(2)}) (\R^{-1}, f_{(1)} \otimes S(g_{(3)})), 
\end{align*}
using properties of the universal $ R $-matrix in the final step.

We shall construct a certain filtration of $ A $ as follows. 
The ordered semigroup will be
\[
 P = \{(\mu,\lambda) \in \weights^+ \times \weights \mid
   \lambda \text{ is a weight of } V(\mu) \}
\]
with the product ordering, namely
\[
 (\mu,\lambda) \leq (\mu',\lambda') \quad \text{iff} \quad 
   \mu \leq \mu' \text { and } \lambda \leq \lambda'.
\]
Note that this partial order is well-founded.

To define the filtration, we begin with a $P$-grading on the underlying vector space of $A$.  
Let us write $C(\mu) = \End(V(\mu))^* \subset \Poly(G_q)$ for the space of matrix coefficients of $V(\mu)$.  Thus, as a linear space we have 
\[
  A = \Poly(G_q) = \bigoplus_{\mu\in\weights^+} C(\mu)
\]
and this determines a $\weights^+$-grading on $A$.  At the same time we have
a $ \weights $-grading on the underlying vector space of $A$ given by the right regular action 
of $ U_q(\mathfrak{h}) $, namely, $f\in \Poly(G_q)$ has weight $\lambda$ with respect to this grading if
\[
  K_\nu \hit f = f_{(1)} (K_\nu,f_{(2)}) = q^{(\nu,\lambda)}f
\]
for all $\nu\in\weights$.  Let us write $\Poly(G_q)_\lambda$ for the subspace of all vectors of weight $\lambda$
and $C(\mu)_\lambda = C(\mu)\cap \Poly(G_q)_\lambda$.   We therefore obtain a direct sum decomposition
\[
 A \cong \bigoplus_{(\mu,\lambda)\in P} C(\mu)_\lambda,
\]
and the $P$-filtration of $A$ is defined by
\[
  \F^{(\mu,\lambda)}(A) = 
   \hspace{-2ex}
   \bigoplus_{(\mu',\lambda') \leq (\mu,\lambda)} 
   \hspace{-2ex} C(\mu')_{\lambda'}.   
\]

We claim that this is an algebra filtration of $A$.  Suppose $f\in C(\mu_1)_{\lambda_1}$ and $g\in C(\mu_2)_{\lambda_2}$.  Note that $f\bullet g \in C(\mu_1)C(\mu_2)$.  Since the irreducible components of $V(\mu_1) \otimes V(\mu_2)$ have highest weights less than or equal to $\mu_1+\mu_2$, we see that 
\[
 f\bullet g \in \bigoplus_{\mu\leq\mu_1+\mu_2} C(\mu).
\]
Moreover, the description of $ \R^{-1} $ in the proof of Theorem \ref{Rmatrixformula} shows that
\[
f \bullet g = (\R^{-1}, f_{(2)} \otimes g_{(1)}) f_{(3)} g_{(2)} (\R^{-1}, f_{(1)} \otimes S(g_{(3)}))
\]
is a sum of terms whose weights with respect to the right regular action of $U_q(\lie{h})$ are of the form $\lambda_1+\lambda_2-\gamma$ with $\gamma\in\roots^+$.  

Therefore 
\[
 f\bullet g \in 
 \hspace{-3ex} \bigoplus_{(\mu,\lambda)\leq(\mu_1+\mu_2,\lambda_1+\lambda_2)}
 \hspace{-3ex}
 C(\mu)_\lambda = \F^{(\mu_1+\mu_2,\lambda_1+\lambda_2)}(A)
\]
as required.

We write $B=gr_\F(A)$ for the associated graded algebra.  We again have
$B \cong \bigoplus_{(\mu,\lambda)\in P} C(\mu)_\lambda = \Poly(G_q)$ as a vector space, and we write $f\circ g$ for the product in $B$ of $f,g\in\Poly(G_q)$.   Then for $f\in C(\mu_1)_{\lambda_1}$ and $g\in C(\mu_2)_{\lambda_2}$, we have that $f \circ g$ is the projection onto  $C(\mu_1+\mu_2)_{\lambda_1+\lambda_2}$ of $f\bullet g$. 
In order to give an explicit formula for this multiplication, we introduce the notation $f\cdot g \in C(\mu+\nu)$ 
for the projection of the usual $\Poly(G_q)$-product $fg$ onto the component $C(\mu+\nu)$, where $f\in C(\mu)$ and $g\in C(\nu)$.  
Using the above formula for $f\bullet g$ and again considering the formula for $\R^{-1}$ from Theorem \ref{Rmatrixformula} we obtain
\begin{align*}
 f \circ g 
  &=   (\R^{-1}, f_{(2)} \otimes g_{(1)}) f_{(3)} \cdot g_{(2)} (q^{\sum_{i,j = 1}^N B_{ij} (H_i \otimes H_j)}, f_{(1)} \otimes g_{(3)})  \\
  &= g_{(1)} \cdot f_{(2)} (\R^{-1}, f_{(3)} \otimes g_{(2)}) (q^{\sum_{i,j = 1}^N B_{ij} (H_i \otimes H_j)}, f_{(1)} \otimes g_{(3)}).
\end{align*}

Let us calculate this product in terms of matrix coefficients.  Let $ e^\mu_1, \dots, e^\mu_m $ be a basis of weight 
vectors for $ V(\mu) $ with dual basis $ e_\mu^1, \dots, e_\mu^m \in V(\mu)^* $. We write $ \epsilon_j $ for the weight of $ e^\mu_j $, and assume 
that the vectors are ordered in a non-ascending order, so that $ \epsilon_i > \epsilon_j $ implies $ i < j $. 
Write $ u^\mu_{ij} = \bra e_\mu^i| \bullet| e^\mu_j \ket $ for the corresponding matrix coefficients. 
Similarly, for $ \nu \in \weights^+ $ we fix a basis $ e^\nu_1, \dots, e^\nu_n $ of $ V(\nu) $ with dual basis $ e_\nu^1, \dots, e_\nu^n \in V(\nu)^* $ 
with the same properties. 
Using this notation one obtains 
\begin{align*}
u^\mu_{ij} \circ u^\nu_{kl} &= q^{(\epsilon_i, \epsilon_l - \epsilon_k)} u^\mu_{ij} \cdot u^\nu_{kl} 
+ \sum_{r = i + 1}^m \sum_{s = 1}^{k - 1} \alpha^{ijkl}_{rs} u^\mu_{rj} \cdot u^\nu_{sl} \\
&= q^{(\epsilon_l, \epsilon_i - \epsilon_j)} u^\nu_{kl} \cdot u^\mu_{ij} + \sum_{u = 1}^{j - 1} \sum_{v = l + 1}^{n} \beta^{ijkl}_{uv} u^\nu_{kv} \cdot u^\mu_{iu} 
\end{align*}
with certain scalars $ \alpha^{ijkl}_{rs}, \beta^{ijkl}_{uv} \in \mathbb{K} $. Moreover $ \alpha^{ijkl}_{rs} = 0 $ unless $ \epsilon_r < \epsilon_i $ and $ \epsilon_s > \epsilon_k $. Similarly $ \beta^{ijkl}_{uv} = 0 $ unless $ \epsilon_u > \epsilon_j $ and $ \epsilon_v < \epsilon_l $.

The above formulas show in particular that the space $ C(\mu) \circ C(\nu) $ is 
contained in $C(\mu)\cdot C(\nu)$. Using an induction argument on the first of these formulas, we also obtain 
\begin{align*}
u^\mu_{ij} \circ u^\nu_{kl} &= q^{(\epsilon_i, \epsilon_l - \epsilon_k)} u^\mu_{ij} \cdot u^\nu_{kl} 
+ \sum_{r = i + 1}^m \sum_{s = 1}^{k - 1} \gamma^{ijkl}_{rs} u^\mu_{rj} \circ u^\nu_{sl} 
\end{align*}
for certain coefficients $ \gamma^{ijkl}_{rs} \in \mathbb{K} $ 
with $ \gamma^{ijkl}_{rs} = 0 $ unless $ \epsilon_r < \epsilon_i $ and $ \epsilon_s > \epsilon_{k} $.
This implies $ C(\mu) \cdot C(\nu) \subset C(\mu) \circ C(\nu) $ and thus $ C(\mu) \circ C(\nu) = C(\mu) \cdot C(\nu)  = C(\mu+\nu) $. 

In a similar way, the second formula for the product $ \circ $ from above, with the roles of $ u^\mu_{ij} $ and $ u^\nu_{kl} $ swapped, yields 
\begin{align*}
u^\nu_{kl} \circ u^\mu_{ij} &= 
q^{(\epsilon_j, \epsilon_k - \epsilon_l)} u^\mu_{ij} \cdot u^\nu_{kl} 
+ \sum_{r = i}^m \sum_{s = 1}^k \sum_{u = 1}^{l - 1} \sum_{v = j + 1}^m \delta^{ijkl}_{rsuv} u^\mu_{rv} \circ u^\nu_{su} 
\end{align*}
for certain coefficients $ \delta^{ijkl}_{rsuv} \in \mathbb{K} $ with $\delta^{ijkl}_{rsuv} =0$ unless $\epsilon_u > \epsilon_l$, $\epsilon_v < \epsilon_j$, $\epsilon_r\leq\epsilon_i$ and $\epsilon_s\geq\epsilon_k$. 
Setting $ q_{ijkl} = q^{(\epsilon_j + \epsilon_i, \epsilon_k -\epsilon_l)} $ this yields 
\begin{align*}
u^\nu_{kl} \circ u^\mu_{ij} - q_{ijkl} u^\mu_{ij} \circ u^\nu_{kl}  
&= \sum_{r = i}^m \sum_{s = 1}^k \sum_{u = 1}^{l - 1} \sum_{v = j + 1}^m \delta^{ijkl}_{rsuv} u^\mu_{rv} \circ u^\nu_{su} \\ 
&\qquad - \sum_{r = i + 1}^m \sum_{s = 1}^{k - 1} q_{ijkl} \gamma^{ijkl}_{rs} u^\mu_{rj} \circ u^\nu_{sl}. 
\end{align*}

We now wish to show that $B$ satisfies the hypotheses of Proposition \ref{genrelnoether}. 
Consider the matrix coefficients $ u^k_{ij} = \bra e_{\varpi_k}^i| \bullet | e^{\varpi_k}_j \ket $ for $ k = 1, \dots, N $ and $ 1 \leq i,j \leq n_k $.  
From the relation $ C(\mu) \circ C(\nu) = C(\mu+\nu) $ for all $ \mu, \nu \in \weights^+ $ one sees immediately that these elements generate $ B $ as an algebra.

Let $ X $ be the collection of all the elements $ u^k_{ij} $. 
We list the elements of $ X $ in an ordered sequence $ u_1, \dots, u_m $ such that the following condition holds: 
For $ u_a = u^r_{ij}, u_b = u^s_{kl} $ we have $ b < a $ if either $ \epsilon_k < \epsilon_i $, or $ \epsilon_k = \epsilon_i $ and $ \epsilon_l < \epsilon_j $. 
According to our above considerations we obtain elements $ q_{ij} \in \mathbb{K}^\times $ and $ \alpha^{st}_{ij}  \in \mathbb{K} $ such that 
$$
u_i \circ u_j = q_{ij} u_i \circ u_j + \sum_{s = 1}^{j - 1} \sum_{t = 1}^m \alpha_{ij}^{st} u_s \circ u_t 
$$
for all $ 1 \leq j < i \leq m $. Therefore Proposition \ref{genrelnoether} shows that $ B $ is Noetherian. 

According to Lemma \ref{filteredgradednoether} it follows that $ A $ is Noetherian. 
\end{proof}

Since $ FU_q(\mathfrak{g}) \subset U_q(\mathfrak{g}) $ is a subalgebra it is immediate from Proposition \ref{uqgnoetheriandomain} that $ FU_q(\mathfrak{g}) $ 
is a domain as well. 

Let us remark that Noetherianity of $ U_q(\mathfrak{g}) $ can be deduced from Theorem \ref{fuqgnoetherian} as follows, independently of 
Proposition \ref{uqgnoetheriandomain}. Firstly, extend $ FU_q(\mathfrak{g}) $ with abstract Cartan generators $ L_\mu $ for $ \mu \in \weights $, commuting with 
the elements of $ FU_q(\mathfrak{g}) $ like the Cartan generators $ K_\mu $ of $ U_q(\mathfrak{g}) $. The resulting algebra $ B $ is Noetherian by 
Theorem \ref{fuqgnoetherian} and Lemma \ref{noetherdomainlemma}, and hence Noetherianity of $ U_q(\mathfrak{g}) $ follows by observing that $ U_q(\mathfrak{g}) $ 
is naturally a quotient of $ B $.

\subsection{Canonical bases} 

In this section we give a brief summary of the theory of canonical bases. This theory, due to Lusztig \cite{Lusztigcanonicalbases} 
and Kashiwara \cite{Kashiwaracrystal1}, \cite{Kashiwaracrystal2}, is devoted to studying the $ q \rightarrow 0 $ limit of the quantized 
universal enveloping algebra $ U_q(\mathfrak{g}) $. For our purposes we only need relatively basic aspects of the theory. A thorough exposition 
can be found in \cite{Lusztigbook}, see also \cite{HKqgcb}. 

In the literature on canonical bases the deformation parameter $ q $ is usually taken to be a transcendental variable over $ \mathbb{Q} $. We will also 
be interested in the specialization to non-root of unity invertible elements in an arbitrary ground field $ \mathbb{K} $.

\subsubsection{Crystal bases} 

In this subsection we discuss without proofs the existence and uniqueness of crystal bases for integrable $ U_q(\mathfrak{g}) $-modules. 

We shall work with the following definition of an abstract crystal. 

\begin{definition} 
	\label{def:crystal}
	Let $ I = \{1, \dots, N \} $ be a finite set. A \emph{crystal} is a set $ B $
	\nomenclature{$B$}{underlying set of a crystal}%
	together with maps $ \tilde{e}_i, \tilde{f}_i: B \sqcup \{0\} \rightarrow B \sqcup \{0\} $ for all $ i \in I $  such that the following conditions hold. 
	\begin{bnum}
		\item[a)] $ \tilde{e}_i(0) = 0 = \tilde{f}_i(0) $ for all $ i \in I $. 
		\item[b)] For any $ i \in I $ and $ b \in B $ there exists $ n \in \mathbb{N} $ such that $ \tilde{e}^n_i(b) = 0 = \tilde{f}^n_i(b) $. 
		\item[c)] For any $ i \in I $ and $ b, c \in B $ we have $ c = \tilde{f}_i(b) $ iff $ \tilde{e}_i(c) = b $. 
	\end{bnum} 
	Given crystals $ B_1, B_2 $, a (strict) morphism from $ B_1 $ to $ B_2 $ is a map $ g: B_1 \sqcup \{0 \} \rightarrow B_2 \sqcup \{0 \} $ 
	such that $ g(0) = 0 $ and $ g $ commutes with all operators $ \tilde{e}_i, \tilde{f}_i $. 
\end{definition}  

In the above definition, the symbol $ \sqcup $ stands for disjoint union. For an element $ b \in B $ one sets 
$$
\varepsilon_i(b) = \text{max} \{n \geq 0 \mid \tilde{e}_i^n(b) \neq 0 \}, \qquad \varphi_i(b) = \text{max} \{n \geq 0 \mid \tilde{f}^n_i(b) \neq 0 \}. 
$$
Let $ \weights $ be the free abelian group abstractly generated by elements $ \varpi_1, \dots, \varpi_N $. 
If $ B $ is a crystal we define a map $ wt: B \rightarrow \weights $ by 
$$
wt(b) = \sum_{i \in I} (\varphi_i(b) - \varepsilon_i(b)) \varpi_i, 
$$
and refer to $ wt(b) $ as the weight of $ b $. 

If $ B_1, B_2 $ are crystals then the direct sum $ B_1 \oplus B_2 $ is the crystal with underlying set $ B_1 \sqcup B_2 $, with the 
operators $ \tilde{e}_i, \tilde{f}_i: B_1 \sqcup B_2 \sqcup \{0\} \rightarrow B_1 \sqcup B_2 \sqcup \{0\} $ induced by the corresponding 
operators for $ B_1 $ and $ B_2 $. 

If $ B_1, B_2 $ are crystals then the tensor product $ B_1 \otimes B_2 $ is the crystal with underlying set $ B_1 \times B_2 $, and elements 
written as $ (b_1, b_2) = b_1 \otimes b_2 $, together with the action 
\begin{align*}
\tilde{e}_i(b_1 \otimes b_2) &= 
\begin{cases} 
\tilde{e}_i(b_1) \otimes b_2 & \text{ if } \varphi_i(b_1) \geq \varepsilon_i(b_2) \\ 
b_1 \otimes \tilde{e}_i(b_2) & \text{ if } \varphi_i(b_1) < \varepsilon_i(b_2) 
\end{cases} 
\\
\tilde{f}_i(b_1 \otimes b_2) &= 
\begin{cases} 
\tilde{f}_i(b_1) \otimes b_2 & \text{ if } \varphi_i(b_1) > \varepsilon_i(b_2) \\ 
b_1 \otimes \tilde{f}_i(b_2) & \text{ if } \varphi_i(b_1) \leq \varepsilon_i(b_2) 
\end{cases}
\end{align*}
where we interpret $ b_1 \otimes 0 = 0 = 0 \otimes b_2 $ for all $ b_1 \in B_1, b_2 \in B_2 $.  
The tensor product $ B_1 \otimes B_2 $ is again a crystal, and the operation of taking tensor products is associative. 

Given a crystal $ B $, an element $ b \in B $ is called a highest weight vector if $ \tilde{e}_i(b) = 0 $ for all $ i $. Similarly, 
$ b \in B $ is called a lowest weight vector if $ \tilde{f}_i(b) = 0 $ for all $ i $. 

For the rest of this subsection we work over the field $ \mathbb{K} = \mathbb{Q}(s) $ where $ s $ is an indeterminate and $ q = s^L $ as before, 
deviating slightly from most of the literature in which $ \mathbb{Q}(q) $ is taken as base field. 
We fix the semisimple Lie algebra $ \mathfrak{g} $ 
associated with a finite Cartan matrix $ A = (a_{ij}) $. Our choice of $ \mathbb{Q}(s) $ is necessary because we work with the \emph{simply connected}
version of $ U_q(\mathfrak{g}) $; however this does not affect the constructions and arguments in a serious way. 

Given an integrable $ U_q(\mathfrak{g}) $-module $ M $ and $ 1 \leq i \leq N $, one can write every element $ m \in M $ 
of weight $ \lambda \in \weights $ uniquely in the form 
$$
m = \sum_{n \geq 0} F_i^{(n)} \cdot m_n 
$$
where $ m_n \in M_{\lambda + n \alpha_i} $ satisfies $ E_i \cdot m_n = 0 $, and we recall that $ F_i^{(n)} = F_i^n/[n]_{q_i}! $. 
On a vector $ m \in M_\lambda $ written in the above form, the Kashiwara operators are defined by 
\begin{align*}
\tilde{e}_i(m) &= \sum_{n \geq 1} F_i^{(n - 1)} \cdot m_n \\
\tilde{f}_i(m) &= \sum_{n \geq 0} F_i^{(n + 1)} \cdot m_n,   
\end{align*}
\label{nom:Kashiwara_operators}
and this is extended linearly to all of $ M $. Hence we obtain linear operators $ \tilde{e}_i, \tilde{f}_i: M \rightarrow M $ for $ i = 1, \dots, N $ 
in this way. 

Consider the algebra $ \A_0 $
\nomenclature[o$A_0$]{$\A_0$}{localization of $\mathbb{Q}[s]$ at $s=0$}%
obtained by localizing the polynomial ring $ \mathbb{Q}[s] $ at the maximal ideal 
generated by $ s $, corresponding to the point $ s = 0 $ on the affine line. Explicitly, the elements of $ \A_0 $ 
can be written in the form $ f(s)/g(s) $ where $ f(s), g(s) \in \mathbb{Q}[s] $ and $ g(0) \neq 0 $. 

Let us also recall that $ \beta $ denotes the field automorphism of $ \mathbb{Q}(s) $ which maps $ s $ to $ s^{-1} $,
and define $ \A_\infty = \beta(\A_0) \subset \mathbb{Q}(s) $. 
This can be viewed as a localization at $ \infty $, which will be more convenient for 
us than localization at $ 0 $. 

The following definition of crystal bases corresponds to the notion of basis at $ \infty $ in the sense of Chapter 20 in \cite{Lusztigbook}. 

\begin{definition} 
	Let $ M $ be an integrable $ U_q(\mathfrak{g}) $-module. A crystal basis $ (\L, \B) $ for $ M $ is a free $ \A_\infty $-module 
	$ \L \subset M $
	such that $ \mathbb{Q}(s) \otimes_{\A_\infty} \L = M $, together with a basis $ \B $
	of the vector space $ \L/s^{-1} \L $ 
	over $ \mathbb{Q} $ such that the following conditions hold. 
	\begin{bnum} 
		\item[a)] For any $ \mu \in \weights $, the space $ \L_\mu = M_\mu \cap \L $
		satisfies 
		$ \mathbb{Q}(s) \otimes_{\A_\infty} \L_\mu = M_\mu $, and $ \B_\mu = \B \cap \L_\mu/s^{-1} \L_\mu $
		is a basis of $ \L_\mu/s^{-1} \L_\mu $ over $ \mathbb{Q} $. 
		\item[b)] The Kashiwara operators $ \tilde{e}_i, \tilde{f}_i $ on $ M $ leave $ \L $ invariant and induce on $ \L/s^{-1}\L $ 
		operators which leave $ \B \sqcup \{0 \} $ invariant, such that $ \B $ becomes a crystal. 
	\end{bnum} 
\end{definition} 

Consider the simple module $ M = V(\lambda) $ for $ \lambda \in \weights^+ $ and let $\L(\lambda)$
be the $ \A_\infty $-submodule of $ V(\lambda) $ 
spanned by all vectors of the form $ \tilde{f}_{i_1} \cdots \tilde{f}_{i_l}(v_\lambda) $ where $ i_1, \dots, i_l \in I $ 
and $ \tilde{f}_1, \dots, \tilde{f}_N $ denote the Kashiwara operators. Moreover let $ \B(\lambda) $ 
be the collection of all nonzero cosets in $ \L(\lambda)/s^{-1} \L(\lambda) $ of the form $ \tilde{f}_{i_1} \cdots \tilde{f}_{i_l}(v_\lambda) $. 
The fact that $ (\L(\lambda), \B(\lambda)) $ is a crystal basis for $ V(\lambda) $ is crucial in the proof of the following foundational result
due to Lusztig and Kashiwara, see \cite{Kashiwaracrystal2}. 

\begin{theorem} \label{crystalmain}
	For every integrable $ U_q(\mathfrak{g}) $-module $ M $ there exists a crystal basis $ (\L, \B) $. 
	Moreover, if $ (\L_1, \B_1), (\L_2, \B_2) $ are crystal bases of $ M $ then there exists an automorphism $ f: M \rightarrow M $ of $ U_q(\mathfrak{g}) $-modules 
	which restricts to an isomorphism of crystals $ f: \B_1 \rightarrow \B_2 $. 
\end{theorem} 

The proof of Theorem \ref{crystalmain} relies on the \emph{grand loop} argument. We refer to \cite{Kashiwaracrystal2} for the details, see also 
\cite{HKqgcb} and Chapter 5 in \cite{Josephbook}. At the same time, one can construct suitable crystal bases for $ U_q(\mathfrak{n}_-) $ 
and for all Verma modules.

\subsubsection{Global bases} 

In this subsection we describe how to use crystal bases to obtain global bases for $ U_q(\mathfrak{g}) $-modules in the case of an arbitrary base 
field $ \mathbb{K} $ provided that $ q = s^L \in \mathbb{K}^\times $ is not a root of unity. We follow the exposition in Chapter 6 of \cite{HKqgcb}. 

Initially, we shall work over $ \mathbb{K} = \mathbb{Q}(s) $ and consider several subrings of this field. 
Note that we may view $ \mathbb{Q}(s) $ as the ring of all fractions $ f(s)/g(s) $ with $ f(s), g(s) \in \mathbb{Z}[s] $ and $ g(s) \neq 0 $. 
Note also that the localization $\A_0$ of $ \mathbb{Q}[s] $ at $ s = 0 $ can be identified with the localization $ \mathbb{Z}[s]_0 $ of the ring $ \mathbb{Z}[s] $ of polynomial functions with integer coefficients at $ s = 0 $. 
That is, elements of $ \A_0 $ can be written as fractions $ f(s)/g(s) $ where $ f(s), g(s) \in \mathbb{Z}[s] $ such that $ g(0) \neq 0 $. Let us also 
define $ \A_\infty = \beta(\A_0) $, where we recall that $ \beta: \mathbb{Q}(s) \rightarrow \mathbb{Q}(s) $ denotes the field automorphism determined by $ \beta(s) = s^{-1} $.  The ring $\A_\infty$ can be thought of as localization of $ \mathbb{Q}[s^{-1}] $ at $ s = \infty $. We shall 
write $ \A = \mathbb{Z}[s, s^{-1}] $ as before. Note that $ \A_0, \A_\infty, \A \subset \mathbb{Q}(s) $ are naturally subrings. 

In the sequel we shall consider various lattices in the sense of the following definition. 

\begin{definition}
	If $ R \subset S $ is a subring and $ V $ is a free $ S $-module, then a free $ R $-submodule $ \L \subset V $ is called a free $ R $-lattice of $ V $ if 
	the canonical map $ S \otimes_R \L \rightarrow V $ is an isomorphism. 
\end{definition} 

Assume that $ V $ is a finite dimensional $ \mathbb{Q}(s) $-vector space, 
and let $ \L_0, \L_\infty $ and $ V_\A $ be free $ \A_0 $-, $\A_\infty $- and $ \A $-lattices of $ V $, respectively. 
Even if no compatibility between these lattices is assumed a priori, we automatically have the following properties, compare Section 6.1 in \cite{HKqgcb}. 

\begin{lemma} \label{balancedtriplehelp}
	In the above situation, the canonical map $ \A_0 \otimes_{\mathbb{Z}[s]} (V_\A \cap \L_0) \rightarrow \L_0 $ is an isomorphism of $ \A_0 $-modules. 
	Similarly, the canonical map $ \A_\infty \otimes_{\mathbb{Z}[s^{-1}]} (V_\A \cap \L_\infty) \rightarrow \L_\infty $ is an isomorphism of $ \A_\infty $-modules. 
\end{lemma} 

\begin{proof} 
	We shall only prove the assertion for $ \L_0 $, the proof for $ \L_\infty $ is analogous. 
	Since $ \A_0 $ and $ \A $ are localizations of $ \mathbb{Z}[s] $, we see that 
	the inclusion map $ V_\A \cap \L_0 \rightarrow \L_0 $ induces injective 
	maps $ \A \otimes_{\mathbb{Z}[s]} (V_\A \cap \L_0) \rightarrow \A \otimes_{\mathbb{Z}[s]} \L_0 $ and 
	$ \A_0 \otimes_{\mathbb{Z}[s]} (V_\A \cap \L_0) \rightarrow \A_0 \otimes_{\mathbb{Z}[s]} \L_0 $. 
	The latter identifies with the canonical map $ \A_0 \otimes_{\mathbb{Z}[s]} (V_\A \cap \L_0) \rightarrow \L_0 $ 
	since $ \A_0 \otimes_{\mathbb{Z}[s]} \L_0 \cong \A_0 \otimes_{\A_0} \L_0 \cong \L_0 $. 
	
	Hence it suffices to show that the map $ \A_0 \otimes_{\mathbb{Z}[s]} (V_\A \cap \L_0) \rightarrow \L_0 $ is surjective. 
	Since $ V_\A \subset V $ is an $ \A $-lattice any element of $ V $, in particular any element $ v \in \L_0 $, can be written in the form $ \frac{f(s)}{g(s)} u $ 
	where $ f(s), g(s) \in \mathbb{Z}[s], g(s) \neq 0 $ and $ u \in V_\A $. 
	That is, there exists a nonzero element $ g(s) \in \mathbb{Z}[s] $ such that $ g(s) v \in V_\A $. 
	Since $ V_\A $ is an $ \A $-module, upon dividing by a suitable power of $ s $ we may assume without loss of generality that $ g(0) \neq 0 $. Hence we have 
	$$
	\frac{1}{g(s)} \otimes g(s) v \in \A_0 \otimes_{\mathbb{Z}[s]} (V_\A \cap \L_0),  
	$$
	and this element maps to $ v $ under the canonical map as desired. 
\end{proof} 

Let us introduce the notion of a balanced triple. 
\begin{definition}
	Let $ V $ be a finite dimensional $ \mathbb{Q}(s) $-vector space. Moreover let $ \L_0 \subset V $ be a free $ \A_0 $-lattice, $ V_\A \subset V $ a free 
	$ \A $-lattice, and $ \L_\infty \subset V $ a free $ \A_\infty $-lattice. 
	\nomenclature[o$L_0$]{$\L_0$}{free $\A_0$-lattice in a balanced triple}%
  \nomenclature[o$L_\infty$]{$\L_\infty$}{free $\A_\infty$-lattice in a balanced triple}%
	\nomenclature[o$V_A$]{$V_\A$}{free $\A$-lattice in a balanced triple}%
	If we define 
	$$
	\L = \L_0 \cap V_\A \cap \L_\infty
	$$
	then $ (\L_0, V_\A, \L_\infty) $ is called a balanced triple for $ V $ provided the following conditions hold. 
	\begin{bnum}
		\item[a)] $ \L $ is a free $ \mathbb{Z} $-lattice for the $ \A_0 $-module $ \L_0 $. 
		\item[b)] $ \L $ is a free $ \mathbb{Z} $-lattice for the $ \A $-module $ V_\A $. 
		\item[c)] $ \L $ is a free $ \mathbb{Z} $-lattice for the $ \A_\infty $-module $ \L_\infty $. 
	\end{bnum}
\end{definition} 

The conditions for $ (\L_0, V_\A, \L_\infty) $ to be a balanced triple mean that $ \L $ has finite rank, and that the canonical multiplication maps 
induce isomorphisms
$$
\A_0 \otimes_\mathbb{Z} \L \cong \L_0, \quad \A \otimes_\mathbb{Z} \L \cong V_\A, \quad \A_\infty \otimes_\mathbb{Z} \L \cong \L_\infty. 
$$

\begin{prop} \label{balancedtriplechar}
	Let $ V $ be a finite dimensional $ \mathbb{Q}(s) $-vector space. Moreover let $ \L_0 \subset V $ be a free $ \A_0 $-lattice, $ V_\A \subset V $ a free 
	$ \A $-lattice, and $ \L_\infty \subset V $ a free $ \A_\infty $-lattice. Then the following conditions are equivalent. 
	\begin{bnum}
		\item[a)] $ (\L_0, V_\A, \L_\infty) $ is a balanced triple for $ V $. 
		\item[b)] The canonical map $ \L \rightarrow \L_0/s \L_0 $ is an isomorphism. 
		\item[c)] The canonical map $ \L \rightarrow \L_\infty/s^{-1} \L_\infty $ is an isomorphism. 
	\end{bnum}
\end{prop}

\begin{proof} 
	$ a) \Rightarrow b) $ We have canonical isomorphisms 
	\begin{align*}
	\L \cong \mathbb{Z} \otimes_\mathbb{Z} \L &\cong \mathbb{Z} \otimes_{\A_0} \A_0 \otimes_\mathbb{Z} \L 
	\cong \mathbb{Z} \otimes_{\A_0} \L_0 \cong \A_0/s \A_0 \otimes_{\A_0} \L_0 \cong \L_0/s \L_0. 
	\end{align*}
	This yields the claim. 
	
	$ b) \Rightarrow a) $ We shall first prove by induction that the canonical map $ m_k: \L \rightarrow \L_0 \cap V_\A \cap s^k \L_\infty $ given by $ m_k(v) = s^k v $ 
	induces an isomorphism 
	$$
	\biggl(\bigoplus_{k = 0}^n \mathbb{Z} s^k \biggr) \otimes_\mathbb{Z} \L \cong \L_0 \cap V_\A \cap s^n \L_\infty
	$$
	for any $ n \in \mathbb{N}_0 $. For $ n = 0 $ this is obvious, so assume that the assertion holds for $ n - 1 $ for some $ n > 0 $. Then 
	we have a canonical isomorphism 
	$$
	\biggl(\bigoplus_{k = 1}^n \mathbb{Z} s^k \biggr) \otimes_\mathbb{Z} \L \cong s \L_0 \cap V_\A \cap s^n \L_\infty
	$$
	and a commutative diagram 
	$$
	\xymatrix{
		0 \ar@{->}[r] & \biggl(\bigoplus_{k = 1}^n \mathbb{Z} s^k \biggr) \otimes_\mathbb{Z} \L 
		\ar@{->}[r] \ar@{->}[d]^\cong & \biggl(\bigoplus_{k = 0}^n \mathbb{Z} s^k \biggr) \otimes_\mathbb{Z} \L 
		\ar@{->}[d] \ar@{->}[r] & \L \ar@{->}[d]^\cong \ar@{->}[r] & 0 \\
		0 \ar@{->}[r] & s \L_0 \cap V_\A \cap s^n \L_\infty \ar@{->}[r] & \L_0 \cap V_\A \cap s^n \L_\infty \ar@{->}[r] & \L_0/s \L_0 
	}
	$$
	with exact rows. Hence the $ 5 $-Lemma shows that the middle vertical arrow is an isomorphism, which yields the inductive step. 
	
	As a consequence, we have 
	$$
	\biggl(\bigoplus_{k = 0}^{b - a} \mathbb{Z} s^k \biggr) \otimes_\mathbb{Z} \L \cong \L_0 \cap V_\A \cap s^{b - a} \L_\infty
	$$
	for all $ b \geq a $, which implies 
	$$
	\biggl(\bigoplus_{k = a}^{b} \mathbb{Z} s^k \biggr) \otimes_\mathbb{Z} \L \cong s^a \L_0 \cap V_\A \cap s^b \L_\infty
	$$
	under the canonical map. Since $ \L_\infty $ is an $ \A_\infty $-lattice in $ V $ we have $ \bigcup_{n = 0}^\infty s^n \L_\infty = V $, 
	and hence the above yields 
	\begin{align*}
	\mathbb{Z}[s] \otimes_\mathbb{Z} \L &\cong V_\A \cap \L_0, \\
	\mathbb{Z}[s, s^{-1}] \otimes_\mathbb{Z} \L &\cong V_\A, \\
	\mathbb{Z}[s^{-1}] \otimes_\mathbb{Z} \L &\cong V_\A \cap \L_\infty. 
	\end{align*}
	The first of these isomorphisms implies 
	$$
	\A_0 \otimes_\mathbb{Z} \L \cong \A_0 \otimes_{\mathbb{Z}[s]} \mathbb{Z}[s] \otimes_\mathbb{Z} \L \cong \A_0 \otimes_{\mathbb{Z}[s]} (V_\A \cap \L_0) \cong \L_0, 
	$$
	using Lemma \ref{balancedtriplehelp} in the last step. In the same way one obtains 
	$$
	\A_\infty \otimes_\mathbb{Z} \L \cong \A_\infty \otimes_{\mathbb{Z}[s^{-1}]} \mathbb{Z}[s^{-1}] \otimes_\mathbb{Z} \L
	\cong \A_\infty \otimes_{\mathbb{Z}[s^{-1}]} (V_\A \cap \L_\infty) \cong \L_\infty, 
	$$
	and we conclude that $ (\L_0, V_\A, \L_\infty) $ is a balanced triple. 
	
	$ a) \Leftrightarrow c) $ is proved in the same way. 
\end{proof} 

Assume that $ (\L_0, V_\A, \L_\infty) $ is a balanced triple for the finite dimensional $ \mathbb{Q}(s) $-vector space $ V $, and let 
$ G: \L_\infty/s^{-1} \L_\infty \rightarrow \L $
be the inverse of the isomorphism $ \L \cong \L_\infty/s^{-1} \L_\infty $ 
obtained in Proposition \ref{balancedtriplechar}. If $ \B $ is a $ \mathbb{Z} $-basis of $ \L_\infty/s^{-1} \L_\infty $ then the 
vectors $ G(b) $
for $ b \in \B $ form an $ \A $-basis $ G(\B) $ of $ V_\A $ and a $ \mathbb{Q}(s) $-basis of $ V $. 
Indeed, writing $ \L = \bigoplus_{b \in \B} \mathbb{Z} G(b) $ we obtain 
$$
V_\A = \A \otimes_\mathbb{Z} \L = \A \otimes_\mathbb{Z} \bigoplus_{b \in \B} \mathbb{Z} G(b) = \bigoplus_{b \in \B} \A G(b), 
$$
and the claim for $ V $ follows similarly from $ V \cong \mathbb{Q}(s) \otimes_\A V_\A $. One calls $ G(\B) $ the global basis associated 
to the local basis $ \B $. 

We continue to work over $ \mathbb{K} = \mathbb{Q}(s) $. Recall the definition of the integral form $ U_q^\A(\mathfrak{g}) $ of 
$ U_q(\mathfrak{g}) $ from Definition \ref{def:integral_form} and that $ U_q^\A(\mathfrak{n}_-) $ is the $ \A $-subalgebra of $ U_q(\mathfrak{n}_-) $ generated by 
the divided powers $ F_i^{(m)} $ for $ i = 1, \dots, N $ and $ m \in \mathbb{N}_0 $. 
The automorphism $ \beta $ of $ U_q(\mathfrak{g}) $ defined in Lemma \ref{defbarinvolution} restricts to an automorphism of $ U_q^\A(\mathfrak{g}) $ 
preserving $ U_q^\A(\mathfrak{n}_-) $. Moreover, for $ \lambda \in \weights^+ $ we obtain a well-defined $ \mathbb{Q} $-linear automorphism 
$ \beta_{V(\lambda)}: V(\lambda) \rightarrow V(\lambda) $ by setting $ \beta_{V(\lambda)}(Y \cdot v_\lambda) = \beta(Y) \cdot v_\lambda $, for any $Y\in U_q(\lie{g})$. 
\label{nom:bar_on_V}%
Indeed, view $ V(\lambda) $ as a module over $ U_q(\mathfrak{g}) $ with action $ X \cdot_\beta v = \beta(X) \cdot v $, and write $ V(\lambda)^\beta $ for 
this module. Then $ V(\lambda)^\beta $ is irreducible and of highest weight $ \lambda $, and hence must be isomorphic to $ V(\lambda) $ 
as $ U_q(\mathfrak{g}) $-module. The corresponding intertwiner is precisely the desired map $ \beta_{V(\lambda)} $. 

Let $ (\L(\lambda), \B(\lambda)) $ be the crystal basis of $ V(\lambda) $ as discussed before Theorem \ref{crystalmain}, and 
let us also write $ \L(\lambda) = \L_\infty(\lambda) $. 
If we define $ \L_0(\lambda) = \beta_{V(\lambda)}(\L(\lambda)) $, then $ \L_0(\lambda) $ is a free $ \A_0 $-lattice of $ V(\lambda) $. Moreover set 
$$
V(\lambda)_\A = U_q^\A(\mathfrak{g}) \cdot v_\lambda = U_q^\A(\mathfrak{n}_-) \cdot v_\lambda. 
$$
Then we have $ \beta_{V(\lambda)}(V(\lambda)_\A) = V(\lambda)_\A $. 

\begin{theorem} \label{globalmain}
	Let $ \lambda \in \weights^+ $. With the notation as above, $ (\L_0(\lambda), V(\lambda)_\A, \L_\infty(\lambda)) $ is a balanced triple for $ V(\lambda) $. 
\end{theorem} 

The proof of Theorem \ref{globalmain} can be found in \cite{Kashiwaracrystal2}, see also Chapter 6 of \cite{HKqgcb} and Section 6.2 in \cite{Josephbook}.  
According to the discussion after Proposition \ref{balancedtriplechar} we obtain global basis elements $ G(b) $ for $ b \in \B(\lambda) $ such that 
$$ 
V(\lambda)_\A = \bigoplus_{b \in \B(\lambda)} \A G(b),  
$$
and the elements $ G(b) $ also form a basis of $ V(\lambda) $ as a $ \mathbb{Q}(s) $-vector space. 

Let us now consider the case where $ \mathbb{K} $ is arbitrary and $ q = s^L \in \mathbb{K}^\times $ is not a root of unity. We write $ U_q(\mathfrak{g}) $ 
for the quantized universal enveloping algebra over $ \mathbb{K} $. 
Then the canonical ring homomorphism $ \A = \mathbb{Z}[s, s^{-1}] \rightarrow \mathbb{K} $ 
induces an $\A$-linear map from the $ \A $-module $ V(\lambda)_\A $ as defined above into the $ \mathbb{K} $-vector space $ V(\lambda) $, the irreducible highest weight module 
of $ U_q(\mathfrak{g}) $ corresponding to $ \lambda \in \weights^+ $. The image of the resulting $ \mathbb{K} $-linear 
map $ \iota_\lambda: \mathbb{K} \otimes_\A V(\lambda)_\A \rightarrow V(\lambda) $ is a nonzero submodule of $ V(\lambda) $, and hence $ \iota_\lambda $ 
is surjective by the irreducibility of $ V(\lambda) $. Conversely, $ \mathbb{K} \otimes_\A V(\lambda)_\A $ is clearly an integrable $ U_q(\mathfrak{g}) $-module, 
and hence a quotient of $ V(\lambda) $ by Theorem \ref{thmfdintegrable} and Theorem \ref{thmfdirreducible}. It follows that $ \iota_\lambda $ is in fact an isomorphism. 

We thus obtain the following result. 
\begin{theorem} \label{globalbasis}
	Let $ \mathbb{K} $ be a field and assume $ q = s^L \in \mathbb{K}^\times $ is not a root of unity. For any $ \lambda \in \weights^+ $ the 
	elements $ G(b) $ for $ b \in \B(\lambda) $ form a basis of $ V(\lambda) $ as a $ \mathbb{K} $-vector space. 
\end{theorem}

\subsection{Separation of Variables} 

In this section we prove a key result on the structure of the locally finite part of $ U_q(\mathfrak{g}) $, originally due to Joseph and 
Letzter \cite{JLseparation}, see also Section 7.3 in \cite{Josephbook}. We follow the approach developed by 
Baumann \cite{Baumannseparation}, \cite{Baumanncanonicalconjugating}. Throughout this section we assume that $ q = s^L \in \mathbb{K}^\times $ is 
not a root of unity.

\subsubsection{Based modules} 

The key part of the proof of separation of variables further below relies on some results involving canonical bases. We shall collect 
these facts here, and refer to \cite{Lusztigbook} for the proofs. In this subsection we work over $ \mathbb{K} = \mathbb{Q}(s) $ with $ q = s^L $ as 
before. 

Let us first introduce the notion of a based module, see Section 27.1 in \cite{Lusztigbook}. 
By definition, an involution on a $ U_q(\mathfrak{g}) $-module $ M $ is a $ \mathbb{Q} $-linear automorphism $ \beta_M: M \rightarrow M $ such that 
$$ 
\beta_M(X \cdot m) = \beta(X) \cdot \beta_M(m) 
$$ 
\label{nom:bar_on_M}%
for all $ X \in U_q(\mathfrak{g}) $ and $ m \in M $, where $ \beta $ is the bar involution of $ U_q(\mathfrak{g}) $ as in Lemma \ref{defbarinvolution}. 
For instance, if $ M = V(\lambda) $ for $ \lambda \in \weights^+ $ and $ v_\lambda \in V(\lambda) $ is a highest weight vector 
then the map $ \beta_M $ defined by $ \beta_M(X \cdot v_\lambda) = \beta(X) \cdot v_\lambda $ for $X\in U_q(\lie{g})$, as discussed before Theorem \ref{globalmain}, 
defines an involution of $ M $. 

\begin{definition} 
	An integrable $ U_q(\mathfrak{g}) $-module $ M $ with involution $ \beta_M $ together with a $ \mathbb{Q}(s) $-basis $ B $ is called a based module 
	if the following conditions hold. 
	\begin{bnum} 
		\item[a)] $ B \cap M_\mu $ is a basis of $ M_\mu $ for any $ \mu \in \weights $. 
		\item[b)] The $ \A $-submodule $ M_\A $ generated by $ B $ is stable under $ U_q^\A(\mathfrak{g}) $. 
		\item[c)] We have $ \beta_M(b) = b $ for all $ b \in B $. 
		\item[d)] The $ \A_\infty $-submodule $ \L_M $ generated by $ B $ together with the image $ \B $ of $ B $ in $ \L_M/s^{-1}\L_M $ 
		forms a crystal basis 
		for $ M $. 
	\end{bnum} 
\end{definition}  

A morphism of based modules from $ (M, B_M) $ to $ (N, B_N) $ is a $ U_q(\mathfrak{g}) $-linear map $ f: M \rightarrow N $ such that 
$ f(b) \in B_N \cup \{0 \} $ for all $ b \in B_M $ and $ B \cap \ker(f) $ is a basis of the kernel $ \ker(f) $. 
In this case the kernel of $ f $ naturally inherits the structure of a based module. 

For any $ \lambda \in \weights^+ $, the simple module $ V(\lambda) $, equipped with the involution fixing $ v_\lambda $, 
together with the corresponding global basis as in Theorem \ref{globalbasis} is a based module. 
In this way simple integrable $ U_q(\mathfrak{g}) $-modules can be viewed as based modules. 

Let us now sketch the construction of a suitable based module structure on the tensor product $ M \otimes N $ of two based modules $ M, N $, 
see Chapter 27.3 in \cite{Lusztigbook}. 
One first defines an involution $ \beta_{M \otimes N} $ on $ M \otimes N $ out of the involutions 
$ \beta_M $ and $ \beta_N $ and the quasi-$ R $-matrix for $ U_q(\mathfrak{g}) $. 
Note here that the ordinary tensor product of the involutions $ \beta_M $ and $ \beta_N $ is not compatible 
with the $ U_q(\mathfrak{g}) $-module structure on $ M \otimes N $ in general. 
The construction of an appropriate basis for $ M \otimes N $ is then characterized by the following result, see Theorem 24.3.3 in \cite{Lusztigbook}. 

\begin{theorem} \label{basedtensorproduct}
	Let $ (M, B) $ and $ (N, C) $ be based modules. Moreover let $ \L $ be the $ \mathbb{Z}[s^{-1}] $-submodule of $ M \otimes N $ 
	generated by all elements $ b \otimes c $ for $ (b,c) \in B \times C $. Then for any $ (b,c) \in B \times C $ there exists a unique 
	element $ b \diamond c \in \L $
	such that $ \beta_{M \otimes N}(b \diamond c) = b \diamond c $ and $ b \diamond c - b \otimes c \in s^{-1} \L $. 
	Moreover, the elements $ b \diamond c $ form a basis $ B \diamond C $ of $ M \otimes N $ which turn the latter into a based module with 
	involution $ \beta_{M \otimes N} $. 
\end{theorem}

The basis of the tensor product $ M \otimes N $ as in Theorem \ref{basedtensorproduct} will also be referred to as the canonical basis. 

Recall that if $ M $ is a finite dimensional $ U_q(\mathfrak{g}) $-module then the dual module is the dual space $ M^* $ with the $ U_q(\mathfrak{g}) $-module 
structure defined by $ (X \cdot f)(m) = f(\hat{S}(X) \cdot m) $ for all $ m \in M $. We have $ (M \otimes N)^* \cong N^* \otimes M^* $ naturally 
if $ M, N $ are finite dimensional. 

Assume that $ N = V(\mu) $ for some $ \mu \in \weights^+ $ with highest weight vector $ v_\mu $. 
Then $ M = V(\mu)^* $ is an irreducible $ U_q(\mathfrak{g}) $-module 
with lowest weight $ -\mu $, and we denote by $ v^\mu \in V(\mu)^* $ the lowest weight vector satisfying $ v^{\mu}(v_\mu) = 1 $. Since $ V(\mu) $ is a 
simple $ U_q(\mathfrak{g}) $-module there exists a unique $ U_q(\mathfrak{g}) $-linear map $ ev_\mu: V(\mu)^* \otimes V(\mu) \rightarrow \mathbb{K} $ 
such that $ ev_\mu(v^\mu \otimes v_\mu) = 1 $. 

For $ \mu, \lambda \in \weights^+ $ let $ p_{\mu \lambda}: V(\mu) \otimes V(\lambda) \rightarrow V(\mu + \lambda) $ 
be the 
unique $ U_q(\mathfrak{g}) $-linear map satisfying $ p_{\mu \lambda}(v_\mu \otimes v_\lambda) = v_{\mu + \lambda} $. 
Similarly, let $ i_{\mu \lambda}: V(\mu + \lambda) \rightarrow V(\mu) \otimes V(\lambda) $
be the unique $ U_q(\mathfrak{g}) $-linear 
map satisfying $ i_{\mu \lambda}(v_{\mu + \lambda}) = v_\mu \otimes v_\lambda $. 
By construction we then have $ p_{\mu \lambda} i_{\mu \lambda} = \id $. 
The transpose of $ p_{\mu \lambda} $ determines a $ U_q(\mathfrak{g}) $-linear map 
$ p_{\mu \lambda}^*: V(\mu + \lambda)^* \rightarrow (V(\mu) \otimes V(\lambda))^* \cong V(\lambda)^* \otimes V(\mu)^* $ 
satisfying $ p_{\mu \lambda}^*(v^{\mu + \lambda}) = v^\lambda \otimes v^\mu $. 
Similarly, the transpose of $ i_{\mu \lambda} $ determines a $ U_q(\mathfrak{g}) $-linear map 
$ i_{\mu \lambda}^*: V(\lambda)^* \otimes V(\mu)^* \cong (V(\mu) \otimes V(\lambda))^* \rightarrow V(\mu + \lambda)^* $ 
satisfying $ i_{\mu \lambda}^*(v^\lambda \otimes v^\mu) = v^{\mu + \lambda} $.  

Let us define $ T_\lambda: V(\mu + \lambda)^* \otimes V(\mu + \lambda) \rightarrow V(\mu)^* \otimes V(\mu) $
as the composition 
$$
\xymatrix{
	V(\mu + \lambda)^* \otimes V(\mu + \lambda) \ar@{->}[r]^{\!\!\!\!\!\!\!\!\!\!\!\!\!\! p_{\lambda \mu}^* \otimes i_{\lambda \mu}} 
	& V(\mu)^* \otimes V(\lambda)^* \otimes V(\lambda) \otimes V(\mu) \ar@{->}[r]^{\qquad \quad ev_\lambda} & V(\mu)^* \otimes V(\mu). 
}
$$
Notice that $ T_\lambda $ maps $ v^{\mu + \lambda} \otimes v_{\mu + \lambda} $ to $ v^{\mu} \otimes v_{\mu} $. 
Since the vector $ v^{\mu} \otimes v_{\mu} $ generates $ V(\mu)^* \otimes V(\mu) $ as a $ U_q(\mathfrak{g}) $-module we see that $ T_\lambda $ is surjective.  
We also observe that $T_\lambda\circ T_{\lambda'} =  T_{\lambda+\lambda'}$ for any $\lambda,\lambda'\in\weights^+$, as a map from $V(\mu+\lambda+\lambda')^*\otimes V(\mu+\lambda+\lambda')$ to $V(\mu)^*\otimes V(\mu)$.

The following result is a translation of Proposition 27.3.5 in \cite{Lusztigbook}. 

\begin{prop} \label{tbasedmodulemorphism}
	The map $ T_\lambda $ is a morphism of based modules. 
\end{prop} 

For any $ \nu \in \weights^+ $ we shall identify $ V(\nu) \cong V(\nu)^{**} $ such that $ v_\nu $ is mapped to the linear form which evaluates $ v^\nu $ to $ 1 $, and consider the corresponding $ U_q(\mathfrak{g}) $-linear isomorphism 
$$ 
(V(\nu)^* \otimes V(\nu))^* \cong V(\nu)^* \otimes V(\nu)^{**} \cong V(\nu)^* \otimes V(\nu). 
$$
Then the transpose of $ T_\lambda $ identifies with a $ U_q(\mathfrak{g}) $-linear map 
$ \delta_\lambda: V(\mu)^* \otimes V(\mu) \rightarrow V(\mu + \lambda)^* \otimes V(\mu + \lambda) $ satisfying 
$$ 
\delta_\lambda(v^\mu \otimes v_\mu)(v^{\mu + \lambda} \otimes v_{\mu + \lambda}) = 1. 
$$
More precisely, we have 
\begin{align*} 
\delta_\lambda(v^\mu \otimes v_\mu) &= v^{\mu + \lambda} \otimes v_{\mu + \lambda} \\
&= i_{\lambda \mu}^*(v^\mu \otimes v^\lambda) \otimes p_{\lambda \mu}(v_\lambda \otimes v_\mu) \\
&= (i_{\lambda \mu}^* \otimes p_{\lambda \mu})(\id \otimes ev_\lambda^* \otimes \id)(v^\mu \otimes v_\mu) 
\end{align*}
with these conventions. 

Consider the partial order on $ \weights^+ $ given by declaring $ \mu \preceq \nu $ if $ \nu - \mu \in \weights^+ $.
\nomenclature[o$<$]{$\preceq$}{partial order on $\weights^+$ given by $ \mu \preceq \nu $ if $ \nu - \mu \in \weights^+ $}%
Note that $ \weights^+ $ 
with this partial order is a lattice, in particular, every pair $ \mu, \nu \in \weights^+ $ has a least upper bound. 

The following Lemma is a reformulation of Lemma 5 in \cite{Baumanncanonicalconjugating}, which is based in turn on Proposition 25.1.10 in \cite{Lusztigbook}. 

\begin{lemma} \label{semigroupactionhelp} 
	Let $ \mu \in \weights^+ $. For each vector $ y $ in the canonical basis of the based module $ V(\mu)^* \otimes V(\mu) $ 
	there exists a greatest element $ \lambda(y) $
	in the lattice $ \weights^+ $ with respect to the partial order $ \preceq $ 
	such that $ T_{\lambda(y)}(y) \neq 0 $. 
\end{lemma} 

Note that the weight $ \lambda(y) $ in Lemma \ref{semigroupactionhelp} is uniquely determined by $ y $, due to the fact that greatest elements 
in partially ordered sets are necessarily unique.

\subsubsection{Further preliminaries} 

In this subsection we collect some additional preparations for the proof of the main theorem presented in the next subsection. Throughout 
this subsection the ground field $ \mathbb{K} $ is arbitrary and $ q = s^L \in \mathbb{K}^\times $ is not a root of unity. 

Let us first discuss a certain filtration of $ \Poly(G_q) $. 
We define 
$$ 
ht(\mu) = \mu_1 + \cdots + \mu_N 
$$ 
\label{nom:ht-weights}%
if $ \mu = \mu_1 \alpha_1 + \cdots + \mu_N \alpha_N \in \weights $. For $ \mu \in \roots $ this agrees 
with the height as in the proof of Proposition \ref{uqgnoetheriandomain}, however, we allow here arbitrary $ \mu \in \weights $. 
It can be shown that $ ht(\mu) \in \frac{1}{2}  \mathbb{N}_0$ for all $ \mu \in \weights^+ $, 
compare Table 1 in Section 13.2 of \cite{HumphreysLie}. 

We obtain an algebra filtration of $ A = \Poly(G_q) $ indexed by $ \frac{1}{2} \mathbb{N}_0 $ by setting 
$$
\F^m(A) = \bigoplus_{ht(\mu) \leq m} V(\mu)^* \otimes V(\mu)
$$
\label{nom:OGq_filtration}%
for $ m \in \frac{1}{2} \mathbb{N}_0 $. 
We call this filtration the height filtration of $ \Poly(G_q) $. 

The associated graded algebra $ E = gr_\F(\Poly(G_q)) $
is canonically isomorphic to $ \Poly(G_q) $ as a left and right $ U_q(\mathfrak{g}) $-module, 
and thus also as a $ U_q(\mathfrak{g}) $-module with the coadjoint action from Section \ref{sec:locally_finite_part}.
We will always view $ \Poly(G_q) $ and $ E $ as a $ U_q(\mathfrak{g}) $-module with this action in the sequel. 

Let us write 
$$
E = \bigoplus_{\mu \in \weights^+} E(\mu) 
$$
with $ E(\mu) = V(\mu)^* \otimes V(\mu) $.
Using the maps $ p_{\mu \lambda}, i_{\mu \lambda} $ introduced before Proposition \ref{tbasedmodulemorphism} we can write 
the product $ x \cdot y \in E $ of elements $ x, y \in E $ as follows.

\begin{lemma} \label{Emultiplication}
	With the notation as above, the multiplication of $ E $ restricts to linear maps  
	$ E(\mu) \otimes E(\lambda) \rightarrow E(\mu + \lambda) $ for $ \mu, \lambda \in \weights^+ $, explicitly given by the formula
	$$
	\bra f| \bullet| v \ket \cdot \bra g| \bullet | w \ket = \bra i_{\lambda \mu}^*(f \otimes g)| \bullet| p_{\lambda \mu}(w \otimes v) \ket 
	$$
	for $ f \in V(\mu)^*, v \in V(\mu), g \in V(\lambda)^*, w \in V(\lambda) $. \\
	Moreover, right multiplication with the quantum trace $ t_\lambda \in E(\lambda) $ identifies with the map $ \delta_\lambda $ defined above,
	up to multiplication by $ q^{(2 \rho, \lambda)} $. 
\end{lemma} 

\begin{proof} 
	Using the identification of $ E $ with $ \Poly(G_q) $ as vector spaces and the pairing between $ U_q(\mathfrak{g}) $ and $ \Poly(G_q) $ we have 
	\begin{align*} 
	(X, \bra f| \bullet| v \ket \cdot \bra g| \bullet | w \ket) &= 
	(g \otimes f)(i_{\lambda \mu}(X \cdot p_{\lambda \mu}(w \otimes v))) \\
	&= (X, \bra i_{\lambda \mu}^*(f \otimes g)| \bullet| p_{\lambda \mu}(w \otimes v) \ket)
	\end{align*}
	for $ X \in U_q(\mathfrak{g}) $, since multiplication in $ E $ only takes into account the highest weight component of the product in $ \Poly(G_q) $.

Next, recall that in our present conventions we identify $V(\lambda)^{**}$ with $V(\lambda)$ such that highest weight vector $v_\lambda\in V(\lambda)$ pairs with $v^\lambda\in V(\lambda)^*$ to give $1$.  This means, more generally, that a vector $v\in V(\lambda)$ corresponds to the form on $V(\lambda)^*$ given by evaluation on $q^{(-2\rho,\lambda)} K_{2\rho}\cdot v$.  It follows that the
 quantum trace $ t_\lambda \in V(\lambda)^* \otimes V(\lambda) $ 
	from Definition \ref{def:quantum_trace} evaluates on $ f \otimes v \in V(\lambda)^* \otimes V(\lambda) $ 
	as $ t_\lambda(f \otimes v) = q^{(-2 \rho, \lambda)} ev_\lambda(f \otimes v) $. 
	That is, we have $ ev_\lambda^*(1) = q^{(2\rho, \lambda)} t_\lambda $, and combining this with the computations at the end of the previous 
	subsection yields the claim. 
\end{proof} 

Let us next review some facts about graded vector spaces. 
Recall that if $ A, B $ are vector spaces filtered by $ \mathbb{N}_0 $ then a linear map $ f: A \rightarrow B $ is called a morphism of filtered vector spaces 
if $ f(\F^n(A)) \subset \F^n(B) $ for all $ n $. In this case $ f $ induces a map $ gr(f): gr(A) \rightarrow gr(B) $
of the associated graded 
vector spaces. Moreover, $ f $ is an isomorphism provided $ gr(f) $ is an isomorphism. 

If $ A, B $ are vector spaces filtered by $ \mathbb{N}_0 $ we obtain a filtration of $ A \otimes B $ by setting 
$$
\F^n(A \otimes B) = \sum_{k + l = n} \F^k(A) \otimes \F^l(B).  
$$
Moreover, the canonical map $ gr(A) \otimes gr(B) \rightarrow gr(A \otimes B) $ is an isomorphism in this case.  

Combining these observations we obtain the following basic fact. 

\begin{lemma} \label{separationgradinglemma}
	Let $ H, Z $ and $ A $ be $ \mathbb{N}_0 $-filtered vector spaces and let $ m: H \otimes Z \rightarrow A $ be a morphism 
	of filtered vector spaces. If the induced map $ gr(m): gr(H) \otimes gr(Z) \rightarrow gr(A) $ is an isomorphism, then 
	$ m $ is an isomorphism as well. 
\end{lemma} 

Finally, we need a result on the structure of certain mapping spaces. 

\begin{lemma} \label{separationmultiplicitylemma} 
	Let $ \mu, \nu, \lambda \in \weights^+ $. Then the linear map 
	$$ 
	\phi: \Hom_{U_q(\mathfrak{g})}(V(\lambda) \otimes V(\nu)^*, V(\mu)) \rightarrow V(\mu)_{\lambda - \nu} 
	$$ 
	given by $ \phi(f) = f(v_\lambda \otimes v^\nu) $ is injective and 
	$$ 
	\im(\phi) = \{v \in V(\mu)_{\lambda - \nu} \mid E_i^{(\nu, \alpha_i^\vee) + 1} \cdot v = 0 \text{ for all } i = 1, \dots, N \}. 
	$$
\end{lemma}

\begin{proof} 
	From weight considerations we see that $ \phi $ is well-defined, that is, $ f(v_\lambda \otimes v^\nu) $ is indeed contained 
	in $ V(\mu)_{\lambda - \nu} $ for $ f \in \Hom_{U_q(\mathfrak{g})}(V(\lambda) \otimes V(\nu)^*, V(\mu)) $. Since $ v_\lambda \otimes v^\nu $ is 
	a cyclic vector for $ V(\lambda) \otimes V(\nu)^* $ it follows that $ \phi $ is injective. 
	
	In order to determine the image of $ \phi $ let us abbreviate 
	$$ 
	U = \{v \in V(\mu)_{\lambda - \nu} \mid E_i^{(\nu, \alpha_i^\vee) + 1} \cdot v = 0 \text{ for all } i = 1, \dots, N \}. 
	$$
	We observe that $ E_i^k \cdot (v_\lambda \otimes v^\nu) = v_\lambda \otimes E_i^k \cdot v^\nu $ for all $ k \in \mathbb{N}_0 $, 
	hence the smallest power of $ E_i $ killing $ v_\lambda \otimes v^\nu $ agrees with 
	the smallest power of $ E_i $ killing $ v^\nu $. This is determined by the weight of $ \nu $ with respect to 
	$ U_{q_i}(\mathfrak{g}_i) \subset U_q(\mathfrak{g}) $, and equals $ E_i^{(\nu, \alpha_i^\vee) + 1} $. 
	It follows that $ \im(\phi) $ is indeed contained in $ U $. 
	
	To prove the reverse inclusion $ U \subset \im(\phi) $ we construct a linear map 
	$$ 
	\psi: U \rightarrow \Hom_{U_q(\mathfrak{g})}(V(\lambda) \otimes V(\nu)^*, V(\mu)) \cong 
	\Hom_{U_q(\mathfrak{g})}(V(\lambda), \Hom(V(\nu)^*, V(\mu))) 
	$$ 
	as follows. Let $ u \in U $ be given and define a $ U_q(\mathfrak{n}_+) $-linear map $ T(u): N(-\nu) \rightarrow V(\mu) $ 
	by $ T(u)(X \cdot v_{-\nu}) = X \cdot u $ for $ X \in U_q(\mathfrak{n}_+) $. 
	Here $ N(-\nu) $ is the universal lowest weight module with lowest weight $ -\nu $ generated by the lowest weight vector $ v_{-\nu} $.
	Since $ E_i^{(\nu, \alpha_i^\vee) + 1} \cdot u = 0 $ 
	for all $ i = 1, \dots, N $, it follows from Theorem \ref{thmfdirreducible} that $ T(u) $ factorizes through $  V(\nu)^* $. Moreover, for 
	any $ v \in V(\nu)^* $ we have 
	\begin{align*}
	(K_\eta \rightarrow T(u))(v) &= q^{(\eta, \nu) + (\eta, \lambda - \nu)} T(u)(v) = q^{(\eta, \lambda)} T(u)(v), \\
	(E_i \rightarrow T(u))(v) &= E_i \cdot T(u)(K_i^{-1} \cdot v) - T(u)(E_i K_i^{-1} \cdot v) = 0
	\end{align*}
	for $ \eta \in \weights $ and $ i = 1, \dots, N $. 
	We conclude that $ T(u) \in \Hom(V(\nu)^*, V(\mu)) $ is a highest weight vector of weight $ \lambda $. 
	Since $ \Hom(V(\nu)^*, V(\mu)) $ is integrable, Theorem \ref{thmfdirreducible} shows 
	that setting $ \psi(u)(v_\lambda) = T(u) $ determines a well-defined element $ \psi(u) $ of $ \Hom_{U_q(\mathfrak{g})}(V(\lambda), \Hom(V(\nu)^*, V(\mu))) $.
	Using the canonical identification with $ \Hom_{U_q(\mathfrak{g})}(V(\lambda) \otimes V(\nu)^*, V(\mu)) $
	it is straightforward to check that $ \phi(\psi(u)) = u $,  
	and we thus conclude $ \im(\phi) = U $ as desired. 
\end{proof} 

As a consequence of Lemma \ref{separationmultiplicitylemma} we see that 
$ \Hom_{U_q(\mathfrak{g})}(\End(V(\lambda)), V(\mu)) $ for $ \mu, \lambda \in \weights^+ $ can be identified with the subspace of $ V(\mu)_0 $ 
consisting of all vectors $ v $ satisfying $ E_i^{(\lambda, \alpha_i^\vee) + 1} \cdot v = 0 $ for $ i = 1, \dots, N $. 
If the coefficients $ m_1, \dots, m_N $ in $ \lambda = \sum_{i = 1}^N m_i \varpi_i $ are sufficiently large the latter 
condition becomes vacuous. In particular, the multiplicity of $V(\mu)$ in $\End(V(\lambda))$ is given by $ [\End(V(\lambda)): V(\mu)] = \dim(V(\mu))_0 $ for all $ \lambda $ of the 
form $ \lambda = \mu + \eta $ for $ \eta \in \weights^+ $.

\subsubsection{Separation of Variables} 

We shall now prove the main result of this section, originally due to Joseph and Letzter \cite{JLseparation}. 
As mentioned above, we follow closely the approach by Baumann \cite{Baumannseparation}, \cite{Baumanncanonicalconjugating}. 

\begin{theorem}[Separation of Variables] \label{separation}
	Assume $ q \in \mathbb{K}^\times $ is not a root of unity. There exists a linear subspace $ \mathbb{H} \subset FU_q(\mathfrak{g}) $,
	\nomenclature[o$H$]{$\mathbb{H}$}{linear space with $ \mathbb{H} \otimes ZU_q(\mathfrak{g}) \cong FU_q(\mathfrak{g}) $ in the separation of variables}%
	invariant under 
	the adjoint action, such that the multiplication map $ \mathbb{H} \otimes ZU_q(\mathfrak{g}) \rightarrow FU_q(\mathfrak{g}) $ is an isomorphism. 
	Moreover, for any $ \mu \in \weights^+ $ we have $ [\mathbb{H}: V(\mu)] = \dim(V(\mu)_0) $ for the multiplicity of the isotypical component 
	of type $ \mu $ of $ \mathbb{H} $. 
\end{theorem} 

\begin{proof} 
	Due to Theorem \ref{adjointdualityplus} and the remarks following it, it suffices to find a linear subspace $ H \subset \Poly(G_q) $, 
	invariant under the coadjoint action, such that the multiplication map $ H \otimes \Poly(G_q)^{G_q} \rightarrow \Poly(G_q) $ is an isomorphism. 
	
	As explained in the proof of Theorem \ref{uqgcenter}, the algebra $ Z = \Poly(G_q)^{G_q} $ has a linear basis consisting of the quantum 
	traces $ t_\lambda $ for $ \lambda \in \weights^+ $, and can be identified with the polynomial 
	algebra $ \mathbb{K}[t_{\varpi_1}, \dots, t_{\varpi_N}] $. In view of Lemma \ref{Emultiplication} it will be convenient to 
	rescale $ t_\lambda $ and work with $ \theta_\lambda = q^{(2\rho, \lambda)} t_\lambda $ instead, and to 
	identify $ Z \cong \mathbb{K}[\theta_{\varpi_1}, \dots, \theta_{\varpi_N}] $. 
	If $ \F $ denotes the height 
	filtration of $ \Poly(G_q) $ discussed in the previous subsection, then each $ \theta_\lambda $ can also be viewed as an element 
	of the associated graded algebra $ E = gr_\F(\Poly(G_q)) $ in a natural way. We have 
	$ \theta_\lambda \theta_\eta = \theta_{\lambda + \eta} $ in $ gr_\F(\Poly(G_q)) $ for all $ \lambda, \eta \in \weights^+ $, 
	and $ \Z = gr_\F(\Poly(G_q)^{G_q}) $ can be identified with the polynomial algebra $ \mathbb{K}[\theta_{\varpi_1}, \dots, \theta_{\varpi_N}] $ 
	as well. 
	
	Let $ \mu \in \weights^+ $ and recall the notation $ E(\mu) = V(\mu)^* \otimes V(\mu) \subset E $ introduced before Lemma \ref{Emultiplication}. 
	According to Lemma \ref{Emultiplication}, the map $ \delta_\lambda: V(\mu)^* \otimes V(\mu) \rightarrow 
	V(\lambda + \mu)^* \otimes V(\lambda + \mu) $ obtained by the transposition of $ T_\lambda $ 
	identifies with the map $ m_\lambda: E(\mu) \rightarrow E(\lambda + \mu) $ given by $ m_\lambda(x) = x \cdot \theta_\lambda $ for any $ \lambda \in \weights^+ $. 
	We note that $ m_\lambda $ is injective since $ T_\lambda $ is surjective. 
	
	Let us momentarily work over $ \mathbb{K} = \mathbb{Q}(s) $. 
	With the identification $ E(\mu) \cong V(\mu)^* \otimes V(\mu)^{**} \cong (V(\mu)^* \otimes V(\mu))^* $ used in the discussion after 
	Proposition \ref{tbasedmodulemorphism}, we obtain a linear basis $ B(\mu) $ of $ E(\mu) $ dual to the 
	canonical basis of $ V(\mu)^* \otimes V(\mu) $. If we set 
	$$ 
	B = \bigcup_{\mu \in \weights^+} B(\mu),  
	$$ 
	then Proposition \ref{tbasedmodulemorphism} shows that $ m_\lambda $ restricts to an injective map $ m_\lambda: B \rightarrow B $ for all $ \lambda \in \weights^+ $. 
	
	If $ x $ is an element of the canonical basis of $ V(\mu)^* \otimes V(\mu) $ let us write $ x^\vee $ for the dual element of $ B $, 
	determined by the relation $ x^\vee(y) = \delta_{xy} $ for all canonical basis elements $ y $. 
	We define a map $ \epsilon: B \rightarrow B $ by stipulating $ \epsilon(x^\vee) = T_{\lambda(x)}(x)^\vee $, where $ \lambda(x) \in \weights^+ $ is 
	the weight determined in Lemma \ref{semigroupactionhelp}. Then we have 
	$$
	x^\vee(y) = \delta_{xy} = \delta_{T_{\lambda(x)}(x), T_{\lambda(x)}(y)} 
	= T_{\lambda(x)}(x)^\vee(T_{\lambda(x)}(y)) 
	= m_{\lambda(x)}(T_{\lambda(x)}(x)^\vee)(y)
	$$
	for any $ y $, so that $ x^\vee = m_{\lambda(x)}(\epsilon(x^\vee)) = \epsilon(x^\vee) \cdot \theta_{\lambda(x)} $. 
	Let us define
	\begin{align*}
	B_\H &= \{\epsilon(b) \mid b \in B\} \\
	B_\Z &= \{\theta_\lambda \mid \lambda \in \weights^+ \}. 
	\end{align*}
	Note that since $ \epsilon(\epsilon(b)) = \epsilon(b) $ for all $ b \in B $ 
	we can write $ B_\H = \{b \in B \mid \epsilon(b) = b \} $. 
	Now assume $ x^\vee = y^\vee \cdot \theta_{\lambda} $ for some $ y^\vee \in B_\H $ and $ \theta_\lambda \in B_\Z $. 
	Then $ 1 = x^\vee(x) = m_{\lambda}(y^\vee)(x) = y^\vee(T_\lambda(x)) $ implies that $ T_{\lambda}(x) = y $ is nonzero. 
	Due to Lemma \ref{semigroupactionhelp} we thus have $ \lambda \preceq \lambda(x) $. Since $ y^\vee \in B_\H $ we also get $T_{\lambda+\mu}(x) = T_\mu(y) = 0 $ 
	for any nonzero element $ \mu \in \weights^+ $. This implies $ \lambda = \lambda(x) $ and hence $ y = T_{\lambda(x)}(x) $, 
	or equivalently $ y^\vee = \epsilon(x^\vee) $. 
	Summarizing this discussion, we conclude that multiplication in $ E $ induces a bijection $ m: B_\H \times B_\Z \rightarrow B $. 
	
	Let $ E_\A \subset E $ be the $ \A $-linear span of the basis $ B $. 
	Similarly, let $ \H_\A \subset E_\A $ the $ \A $-linear span of $ B_\H $, and let $ \Z_\A \subset E $ be the $ \A $-linear span of $ B_\Z $. 
	From our above considerations we conclude that multiplication in $ E $ induces an isomorphism 
	$$ 
	\H_\A \otimes_\A \Z_\A \rightarrow E_\A 
	$$ 
	of $ \A $-modules. 
	
	Now let $ \mathbb{K} $ be again an arbitrary field and $ q = s^L \in \mathbb{K}^\times $ not a root of unity. We write $ \Poly(G_q) $ 
	for the $ \mathbb{K} $-algebra of matrix coefficients, and 
	let $ Z = \Poly(G_q)^{G_q} \subset \Poly(G_q) $ be the invariant subspace with respect to the coadjoint action. 
	As before we denote by $ E = gr_\F(\Poly(G_q)) $ the associated graded algebra for the height filtration of $ \Poly(G_q) $. 
	The canonical bases from Theorem \ref{basedtensorproduct} determine $ \mathbb{K} $-bases for the $ U_q(\mathfrak{g}) $-modules $ V(\mu)^* \otimes V(\mu) $,  
	and we write $ B_\mathbb{K} $ for the dual $ \mathbb{K} $-basis of $ \Poly(G_q) $, and hence of  $ E $, which is obtained in this way. 
	Mapping the elements $ b \in B \subset E_\A $ to the corresponding basis vectors 
	of $ B_\mathbb{K} $ in $ E $ yields a $ \mathbb{K} $-algebra isomorphism $ E_\A \otimes_\A \mathbb{K} \cong E $ 
	by Lemma \ref{Emultiplication}. 
	
	Let $ H \subset \Poly(G_q) $ be the $ \mathbb{K} $-linear span of the image of $ B_\H $ in $ \Poly(G_q) $. 
	If $ \H = gr_\F(H) \subset E $ denotes the associated graded of $ H $ and $ \Z = gr_\F(Z) \subset E $ the associated graded of $ Z $, 
	then we have $ \Z_\A \otimes_\A \mathbb{K} \cong \Z $ and $ \H_\A \otimes_\A \mathbb{K} \cong \H $ under the identification $ E_\A \otimes_\A \mathbb{K} \cong E $. 
	Hence the isomorphism $ \H_\A \otimes_\A \Z_\A \rightarrow E_\A $ of $ \A $-modules induces an isomorphism
	$$
	\H \otimes \Z \rightarrow E 
	$$
	of $ \mathbb{K} $-vector spaces, given by multiplication in $ E $.  
	Therefore, we may apply Lemma \ref{separationgradinglemma}
	to conclude that multiplication in $ \Poly(G_q) $ induces an isomorphism 
	$$
	H \otimes \Poly(G_q)^{G_q} = H \otimes Z \rightarrow \Poly(G_q)
	$$
	as required. 
	
	In order to determine the multiplicity of $ V(\mu) $ inside $ H $ observe that $ [H: V(\mu)] = [gr_\F(H): V(\mu)] = [\H: V(\mu)] $. 
	Thanks to the decomposition $ B_\H \times B_\Z \cong B $ described above, for any $\lambda\in\weights^+$ the canonical basis vectors $b\in B(\lambda)$ which belong to submodules of $E(\lambda)$ of highest weight $\mu\in\weights^+$ can all be written in the form $b=\epsilon(b)\cdot\theta_\nu$ for some $\nu\in\weights^+$, where $\epsilon(b) \in B_\H$ belongs to a submodule of $\H$ of highest weight $\mu$.  It follows that the multiplicity $ [\H: V(\mu)] $ is given by the maximal 
	value of $ [\End(V(\lambda)): V(\mu)] $ as $ \lambda $ runs through $ \weights^+ $. 
	According to the remarks after Lemma \ref{separationmultiplicitylemma}, we thus obtain $ [\H: V(\mu)] = \dim(V(\mu))_0 $. This finishes the proof. 
\end{proof} 

We remark that each basis element of the space $ \mathbb{H} \subset FU_q(\mathfrak{g}) $ 
corresponding to $ H \subset \Poly(G_q) $ 
constructed in the proof of Theorem \ref{separation} is contained in a subspace of the form $ U_q(\mathfrak{g}) \rightarrow K_{2\lambda} $ 
for some $ \lambda \in \weights^+ $. Note finally that since $ ZU_q(\mathfrak{g}) \subset FU_q(\mathfrak{g}) $ commutes pointwise 
with all elements of $ FU_q(\mathfrak{g}) $, Theorem \ref{separation} shows that multiplication also induces a $ U_q(\mathfrak{g}) $-linear isomorphism 
$ ZU_q(\mathfrak{g}) \otimes \mathbb{H} \cong FU_q(\mathfrak{g}) $.

\newpage 
\section{Complex semisimple quantum groups} \label{chcomplexqg}

In this chapter we introduce our main object of study, namely complex semisimple quantum groups. We complement the 
discussion with some background material on locally compact quantum groups in general, and on compact quantum groups arising from 
$ q $-deformations in particular. 

Throughout this chapter we work over the complex numbers $ \mathbb{C} $. If $ \H $ is a Hilbert space we write $ \LH(\H) $
\nomenclature[o$L(H)$]{$\LH(\H)$}{bounded operators on a Hilbert space $\H$}%
for the algebra of 
bounded operators on $ \H $, and denote by $ \KH(\H) $
\nomenclature[o$K(H)$]{$\KH(\H)$}{compact operators on a Hilbert space $\H$}%
the algebra of compact operators. If $ A $ is a $ C^* $-algebra we write $ M(A) $
\nomenclature{$M(A)$}{$C^*$-algebraic multiplier algebra of a $C^*$-algebra $A$}%
for
the multiplier algebra of $ A $ in the sense of $ C^* $-algebras; this needs to be distinguished from the algebraic multiplier algebra $ \M(A) $ 
in chapter \ref{chhopf}. By slight abuse of notation, we use the symbol $ \otimes $ to denote algebraic tensor products, tensor products of Hilbert spaces, 
minimal tensor products of $ C^* $-algebras, or spatial tensor products of von Neumann algebras. It should always be clear from 
the context which tensor product is used. 
If $ X $ is a subset of a Banach space $ B $ we write $ [X] \subset B $ for the closed linear span 
of $ X $. For general background on $ C^* $-algebras and von Neumann algebras we refer to \cite{Murphybook}.

\subsection{Locally compact quantum groups} \label{seclcqg}

In this section we review some basic definitions and facts from the theory of locally compact quantum groups \cite{KVLCQG}. 

\subsubsection{Hopf $ C^* $-algebras} 

Let us start with basic definitions and constructions related to Hopf $ C^* $-algebras.   

\begin{definition} \label{defhopfcstar}
A Hopf $ C^* $-algebra is a $ C^* $-algebra $ H $ together with an injective nondegenerate $ * $-homomorphism
$ \Delta: H \rightarrow M(H \otimes H) $ such that the diagram
$$
\xymatrix{
H \ar@{->}[r]^{\Delta} \ar@{->}[d] & M(H \otimes H) \ar@{->}[d]^{\id \otimes \Delta} \\
M(H \otimes H) \ar@{->}[r]^{\!\!\!\!\!\!\!\!\! \Delta \otimes \id} & M(H \otimes H \otimes H)
     }
$$
is commutative and $ [\Delta(H)(1 \otimes H)] = H \otimes H = [(H \otimes 1)\Delta(H)] $. 
\end{definition}

Comparing Definition \ref{defhopfcstar} with the algebraic definition of a multiplier Hopf algebra in Definition \ref{defmhopf}, 
we note that the density conditions in the former can be thought of as a replacement of the requirement that 
the Galois maps are isomorphisms in the latter. 

If $ H $ is a Hopf $ C^* $-algebra we write $ H^\cop $ for the Hopf-$ C^* $-algebra obtained by equipping
$ H $ with the opposite comultiplication $ \Delta^\cop = \sigma \Delta $,  where as previously, we use $\sigma$ to denote the flip map.
\label{nom:flip2}%

A unitary corepresentation of a Hopf-$ C^* $-algebra $ H $ on a Hilbert space $ \E $ is a unitary $ X \in 
	 M(H \otimes \KH(\E)) $ satisfying
$$
(\Delta \otimes \id)(X) = X_{13} X_{23}
$$
in $ M(H\otimes H\otimes \mathbb{K}(\E)) $, where we are using leg-numbering notation as described after Definition \ref{defquasitriangular}. 
More generally, we may define a corepresentation of $H$ with values in a $C^*$-algebra $A$ as a unitary $X \in M(H \otimes A)$ satisfying $(\Delta \otimes \id)(X) = X_{13} X_{23}$.

A universal dual of $ H $ is a Hopf-$ C^* $-algebra $ \hat{H} $ together with a unitary
corepresentation $ \X \in M(H \otimes \hat{H}) $ of $H$ with values in $\hat{H}$ satisfying the following universal property:
for every Hilbert space $ \E $ and every unitary corepresentation $ X \in M(H\otimes\KH(\E)) $ there exists a unique nondegenerate $ * $-homomorphism
$ \pi_X: \hat{H} \rightarrow \LH(\E) $ such that $ (\id \otimes \pi_X)(\X) = X $. 

We shall be exclusively interested in Hopf $ C^* $-algebras arising from locally compact quantum groups.

\subsubsection{The definition of locally compact quantum groups} 

The theory of locally compact quantum groups has been axiomatized by Kustermans and Vaes \cite{KVLCQG}. 

Let $ \phi $ be a normal, semifinite and faithful weight on a von Neumann algebra $ M $.  Let $M_+$ denote the set of positive elements of $M$. 
\nomenclature[o$M_+$]{$M_+$}{set of positive elements of a von Neumann algebra $\M$}%
We use the standard notation
$$
\M^+_\phi = \{ x \in M_+| \phi(x) < \infty \}, \qquad \N_\phi = \{ x \in M| \phi(x^* x) < \infty \}
$$
\nomenclature[o$M_\phi$]{$\M^+_\phi$}{set of positive integrable elements with respect to $\phi$}%
\nomenclature[o$N_\phi$]{$\N_\phi$}{set of square-integrable elements with respect to $\phi$}%
and write $ M_*^+ $ for the space of positive normal linear functionals on $ M $.
\nomenclature{$M_+^*$}{set of positive normal linear functionals on a von Neumann algebra $M$}%
Assume that
$ \Delta: M \rightarrow M \otimes M $ is a normal unital $ * $-homomorphism. The weight $ \phi $ is called left invariant
with respect to $ \Delta $ if
$$
\phi((\omega \otimes \id)\Delta(x)) = \phi(x) \omega(1)
$$
for all $ x \in \M_\phi^+ $ and $ \omega \in M_*^+ $. Similarly, a normal semifinite faithful weight $\psi$ is right-invariant if 
$$
\psi((\id \otimes \omega)\Delta(x)) = \psi(x) \omega(1)
$$
for all $ x \in \M_\phi^+ $ and $ \omega \in M_*^+ $.

\begin{definition}\label{defqg}
A locally compact quantum group $ G $ is given by a von Neumann algebra $ L^\infty(G) $
\nomenclature{$L^\infty(G)$}{von Neumann algebra of functions of a locally compact quantum group}%
together with a normal unital $ * $-homomorphism
$ \Delta: L^\infty(G) \rightarrow L^\infty(G) \otimes L^\infty(G) $ satisfying the coassociativity relation
$$
(\Delta \otimes \id)\Delta = (\id \otimes \Delta)\Delta,
$$
and normal semifinite faithful weights $ \phi $ and $ \psi $ on $ L^\infty(G) $ which are left and right invariant, respectively.
The weights $\phi$ and $\psi$ are called the left and right Haar weights of $G$.
\label{nom:Haar_weights}
\end{definition}

Our notation for locally compact quantum groups is intended to make clear how ordinary locally compact groups can be viewed as quantum groups.
Indeed, if $ G $ is a locally compact group, then the algebra $ L^\infty(G) $ of essentially bounded measurable functions on $ G $ together with the comultiplication
$ \Delta: L^\infty(G) \rightarrow L^\infty(G) \otimes L^\infty(G) $ given by
$$
\Delta(f)(s,t) = f(st)
$$
defines a locally compact quantum group. The weights $ \phi $ and $ \psi $ are given by integration with respect to left and right Haar measures, respectively. 

For a general locally compact quantum group $ G $ the notation $ L^\infty(G) $ is purely formal.
Similar remarks apply to the $ C^* $-algebras $ C^*_\mx(G), C^*_\red(G) $ and $ C^\mx_0(G), C^\red_0(G) $ associated to $ G $ that we discuss below.
It is convenient to view all of them as different incarnations of the quantum group $ G $. 

Let $ G $ be a locally compact quantum group and let $ \Lambda: \N_\phi \rightarrow L^2(G) $
\label{nom:lc_L2G}
be the GNS-construction for the weight $ \phi $.
Throughout we shall only consider quantum groups for which $ L^2(G) $ is a separable Hilbert space. 
One obtains a unitary $ W_G = W $ on $ L^2(G) \otimes L^2(G) $ such that
$$
W^*(\Lambda(x) \otimes \Lambda(y)) = (\Lambda \otimes \Lambda)(\Delta(y)(x \otimes 1))
$$
\label{nom:multiplicative_unitary}%
for all $ x, y \in \N_\phi $. This unitary is multiplicative, which means that $ W $ satisfies the pentagonal equation
$$
W_{12} W_{13} W_{23} = W_{23} W_{12}.
$$
From $ W $ one can recover the von Neumann algebra $ L^\infty(G) $ as the strong closure of the algebra
$ (\id \otimes \LH(L^2(G))_*)(W) $ where $ \LH(L^2(G))_* $ denotes the space of normal linear functionals on $ \LH(L^2(G)) $. Moreover
one has
$$
\Delta(x) = W^*(1 \otimes x) W
$$
for all $ x \in M $. The algebra $ L^\infty(G) $ has an antipode which
is an unbounded, $ \sigma $-strong* closed linear map $ S $ given by $ S(\id \otimes \omega)(W) = (\id \otimes \omega)(W^*) $
\label{nom:vN-antipode}
for $ \omega \in \LH(L^2(G))_* $. 
\nomenclature[o$L(L^2(G))_*$]{$\LH(L^2(G))_*$}{predual of $\LH(L^2(G))$}%
Moreover there is a polar decomposition $ S = R \tau_{-i/2} $ where $ R $ is an
antiautomorphism of $ L^\infty(G) $ called the unitary antipode and $ (\tau_t) $ is a strongly continuous one-parameter group
of automorphisms of $ L^\infty(G) $ called the scaling group.
The unitary antipode satisfies $ \sigma(R \otimes R) \Delta = \Delta R $. 

The group-von Neumann algebra $ \mathcal{L}(G) $
\nomenclature[o$L(G)$]{$\mathcal{L}(G)$}{group von Neumann algebra of a locally compact quantum group}
of the quantum group $ G $ is the strong closure of the algebra $ (\LH(L^2(G))_* \otimes \id)(W) $ 
with the comultiplication $ \hat{\Delta}: \mathcal{L}(G) \rightarrow \mathcal{L}(G) \otimes \mathcal{L}(G) $ given by
$$
\hat{\Delta}(y) = \hat{W}^* (1 \otimes y) \hat{W},
$$
\nomenclature[o$W$12]{$\hat{W}$}{multiplicative unitary of the Pontrjagin dual}
\label{nom:dual_vN-coproduct}%
where $ \hat{W} = \Sigma W^* \Sigma $ and $ \Sigma \in \LH(L^2(G) \otimes L^2(G)) $ is the flip map. 
\nomenclature{$\Sigma$}{flip map on $L^2(G) \otimes L^2(G)$}%
It defines
a locally compact quantum group $ \hat{G} $
\label{nom:lc-dual}%
which is called the dual of $ G $. The GNS construction for the left invariant weight
$ \hat{\phi} $ of the dual quantum group
\label{nom:dual_vN_Haar_weight}%
can be identified with a map $ \hat{\Lambda}: \N_{\hat{\phi}} \rightarrow L^2(G) $
such that we have $ \mathcal{L}(G) = L^\infty(\hat{G}) $. 
\nomenclature[o$\Lambda$2]{$\hat{\Lambda}$}{GNS map $\N_{\hat{\phi}} \rightarrow L^2(G)$ of the dual of a locally compact quantum group}

We will mainly work with the $ C^* $-algebras associated to the locally compact quantum group $ G $. The
algebra $ [(\id \otimes \LH(L^2(G))_*)(W)] $ is a strongly dense $ C^* $-subalgebra of $ L^\infty(G) $
which we denote by $ C^\red_0(G) $. 
\nomenclature[o$CG$11]{$C^\red_0(G)$}{reduced $C^*$-algebra of continuous functions vanishing at infinity on a locally compact quantum group $G$}%
Dually, the algebra $ [(\LH(L^2(G))_* \otimes \id)(W)] $ is a strongly dense $ C^* $-subalgebra of $ \mathcal{L}(G) $ 
which we denote by $ C^*_\red(G) $.
\nomenclature[o$CG$21]{$C^*_\red(G)$}{reduced group $C^*$-algebra of a locally compact quantum group $G$}%
These algebras are called the reduced algebra of continuous functions vanishing at infinity
on $ G $ and the reduced group $ C^* $-algebra of $ G $, respectively. One has
$ W \in M(C^\red_0(G) \otimes C^*_\red(G)) $. 
Restriction of the comultiplications on $ L^\infty(G) $ and $ \mathcal{L}(G) $
turns $ C^\red_0(G) $ and $ C^*_\red(G) $ into Hopf $ C^* $-algebras.

For every locally compact quantum group $ G $ there exists a universal dual $ C^*_\mx(G) $ of $ C_0^\red(G) $ and
a universal dual $ C^\mx_0(G) $ of $ C^*_\red(G) $, respectively \cite{Kustermansuniversal}.
\nomenclature[o$CG$22]{$C^*_\mx(G)$}{maximal group $C^*$-algebra of a locally compact quantum group $G$}%
\nomenclature[o$CG$12]{$C^\mx_0(G)$}{maximal $C^*$-algebra of continuous functions vanishing at infinity on a locally compact quantum group $G$}%
We call $ C^*_\mx(G) $ the maximal group $ C^* $-algebra of $ G $ and
$ C_0^\mx(G) $ the maximal algebra of continuous functions on $ G $ vanishing at infinity.
Since $ L^2(G) $ is assumed to be separable the $ C^* $-algebras $ C^\mx_0(G), C^\red_0(G) $ and
$ C^*_\mx(G), C^*_\red(G) $ are separable. The quantum group $ G $ is called compact if $ C^\mx_0(G) $ is unital, and it is called
discrete if $ C^*_\mx(G) $ is unital. In the compact case we also write $ C^\mx(G) $ and $ C^\red(G) $ instead of $ C^\mx_0(G) $ and 
$ C^\red_0(G) $, respectively. 

In general, we have a surjective morphism $ \hat{\pi}: C^*_\mx(G) \rightarrow C^*_\red(G) $ of Hopf-$ C^* $-algebras
associated to the left regular corepresentation $ W \in M(C_0(G) \otimes C^*_\red(G)) $. Similarly, there is
a surjective morphism $ \pi: C^\mx_0(G) \rightarrow C^\red_0(G) $.
We will call the quantum group $ G $ amenable if $ \hat{\pi}: C^*_\mx(G) \rightarrow C^*_\red(G) $
is an isomorphism and coamenable if $ \pi: C^\mx_0(G) \rightarrow C^\red_0(G) $
is an isomorphism. If $ G $ is amenable or coamenable, respectively, we also write $ C^*(G) $ or $ C_0(G) $
for the corresponding $ C^* $-algebras. For more information on amenability for locally compact quantum groups see \cite{BedosTuset}. 
\nomenclature[o$CG$23]{$C^*(G)$}{group $C^*$-algebra of an amenable locally compact quantum group $G$}%
\nomenclature[o$CG$13]{$C_0(G)$}{$C^*$-algebra of continuous functions vanishing at infinity on a coamenable locally compact quantum group $G$}%

\subsection{Algebraic quantum groups} \label{secaqg}

The analytical theory of locally compact quantum groups simplifies considerably if one restricts attention to examples that are essentially determined 
algebraically.  This is the case for the class of \emph{algebraic quantum group} in the sense of van Daele \cite{vDadvances}. 
The concept of an algebraic quantum group is a variant of the notion of a regular multiplier Hopf algebra with integrals in which one adds $ * $-structures.

\subsubsection{The definition of algebraic quantum groups}  

Recall the notion of a regular multiplier Hopf algebra with integrals from Definitions \ref{defmhopf} and \ref{defmhopfintegrals}.
Let us now introduce the notion of a multiplier Hopf $ * $-algebra \cite{vDadvances}.  

\begin{definition}
A multiplier Hopf $ * $-algebra is a regular multiplier Hopf algebra $ H $ which is equipped with a $ * $-structure such that 
$ \Delta: H \rightarrow \M(H \otimes H) $ is a $ * $-homomorphism. 
\end{definition} 

Let $ H $ be a multiplier Hopf $ * $-algebra. Then the counit $ \epsilon: H \rightarrow \mathbb{C} $ is a $ * $-homomorphism and the 
antipode $ S: H \rightarrow H $ satisfies 
\[
 S(S(f^*)^*) = f 
\] 
for all $ f \in H $, see Section 5 in \cite{vDmult}. Let us remark that 
the regularity condition in the definition of a regular multiplier Hopf algebra is in fact automatic in the $ * $-algebraic situation. 

\begin{definition} \label{defaqg}
An algebraic quantum group is a multiplier Hopf $ * $-algebra $ H = \CF^\infty_c(G) $
\nomenclature[o$C^\infty(G)$]{$\CF^\infty_c(G)$}{algebra of functions on an algebraic quantum group}%
such that there exists a positive left invariant 
integral $ \phi: H \rightarrow \mathbb{C} $ and a positive right invariant integral $ \psi: H \rightarrow \mathbb{C} $. 
We will also refer to the virtual object $G$ as an algebraic quantum group.
\end{definition} 

Here a linear functional $ \omega: H \rightarrow \mathbb{C} $ is called positive if $ \omega(f^* f) \geq 0 $ for all $ f \in H $. 
We note that (positive) left/right invariant integrals are always unique up to a (positive) scalar, see Section 3 in \cite{vDadvances}.

In a similar way as in the definition of a locally compact quantum group, our notation $ H = \CF^\infty_c(G) $ is meant to suggest 
that $ H $ should be thought of as an algebra of compactly supported smooth functions on an underlying object $ G $, 
and by slight abuse of language which we will sometimes also refer to the latter as an algebraic quantum group. 
In contrast to the situation for locally compact quantum groups the situation is not quite as clean here; for instance, if $ G $ is a Lie group 
then the algebra $ C^\infty_c(G) $ is typically not a multiplier Hopf algebra.

The duality theory for a regular multiplier Hopf algebra with integrals $ H $ discussed in Section \ref{secintegrals} 
is compatible with the positivity requirement for Haar functionals, see \cite{vDadvances}. In particular, if $ H $ is an algebraic 
quantum group and we consider the $ * $-structure on $ \hat{H} $ defined by 
$$
(x^*, f) = \overline{(x, S(f)^*)}
$$
for $ f \in H $ and $ x \in \hat{H} $, then the dual $ \hat{H} $ is an algebraic quantum group as well. 
When using the notation $H=\CF^\infty_c(G)$ we will write either $\hat{H} = \DF(G)$,
\nomenclature[o$D(G)$]{$\DF(G)$}{group algebra of an algebraic quantum group}%
which we refer to as the group algebra of $G$, or $\hat{H}=\CF^\infty_c(\hat{G})$,
\nomenclature[o$C^\infty(G-hat)$]{$\CF^\infty_c(\hat{G})$}{algebra of functions on the Pontrjagin dual of an algebraic quantum group $\hat{G}$}%
where $\hat{G}$ is called the Pontrjagin dual of $G$.
\label{nom:alg_dual}%

One has
the following version of Theorem \ref{bidualitymh}. 

\begin{theorem}[Biduality Theorem for algebraic quantum groups] \label{bidualityaqg}
Let $ H $ be an algebraic quantum group. Then the dual of $ \hat{H} $ is isomorphic to $ H $ as an algebraic quantum group. 
\end{theorem}

Let us note that for an algebraic quantum group $H$, the modular element $ \delta \in M(H) $ of Subsection \ref{sec:defintegrals}, is a positive element, see \cite{Kustermansanalyticalgebraic}. 
We write $ \hat{\delta} \in M(\hat{H}) $ for the modular element of the dual $ \hat{H} $.
\nomenclature[o$\delta$2]{$\hat{\delta}$}{modular element of the Pontrjagin dual of an algebraic quantum group}%

\subsubsection{Algebraic quantum groups on the Hilbert space level} 

In this subsection we explain how to associate a locally compact quantum group to any algebraic quantum group. 
A detailed exposition of this is material can be found in the work of Kustermans and van Daele \cite{KvD}, \cite{Kustermansanalyticalgebraic}. 

Let $ G $ be an algebraic quantum group and let $ \phi: \CF^\infty_c(G) \rightarrow \mathbb{C} $ be a left invariant integral. 
We write $ L^2(G) $ 
\label{nom:L2_alg}%
for the Hilbert space completion of $ \CF^\infty_c(G) $ with respect to the scalar product 
$$
\bra f, g \ket = \phi(f^* g), 
$$
\label{nom:GNS-ip}%
and we let $ \Lambda: \CF^\infty_c(G) \rightarrow L^2(G) $ be the GNS map. 
\label{nom:alg_GNS}%
Then one can define a unitary operator $ W $ on $ L^2(G) \otimes L^2(G) $ by 
$$
W(\Lambda(f) \otimes \Lambda(g)) = \Lambda(S^{-1}(g_{(1)}) f) \otimes \Lambda(g_{(2)}),  
$$
\label{nom:W_alg}%
using the inverse of the antipode $ S $ of $ \CF^\infty_c(G) $. The inverse of $ W $ is given by 
$$
W^*(\Lambda(f) \otimes \Lambda(g)) = \Lambda(g_{(1)} f) \otimes \Lambda(g_{(2)}), 
$$
which formally agrees with the definition given in the case of locally compact quantum groups. 
It is straightforward to check that $ W $ is multiplicative, that is, we have the pentagon relation $ W_{12} W_{13} W_{23} = W_{23} W_{12} $ 
as in Section \ref{seclcqg}. 

Note that the action of $ \CF^\infty_c(G) $ on itself by left multiplication induces 
a $ * $-homomorphism $ \lambda: \CF^\infty_c(G) \rightarrow \LH(L^2(G)) $, explicitly given by 
$$
\lambda(f)(\Lambda(g)) = \Lambda(fg) 
$$
\nomenclature{$\lambda$}{GNS representation $\lambda:\CF^\infty_c(G) \to \mathbb{L}(L^2(G))$ of an algebraic quantum group}%
for $ f, g \in \CF^\infty_c(G) $. 
To see this we first claim that the left regular action of $ \CF^\infty_c(G) $ on $ \CF^\infty_c(G) \subset L^2(G) $ can be written in the form
$$ 
\lambda(f) = (\id \otimes \omega_{f_{(1)}^*, S(f_{(2)}) \chi})(W), 
$$ 
where $ \chi \in \CF^\infty_c(G) $ is any element satisfying $ \phi(\chi) = 1 $, and $ \omega_{h,k}(T) = \bra \Lambda(h), T \Lambda(k) \ket $
for all $ h, k \in \CF^\infty_c(G) $. Indeed, we have 
\begin{align*}
(\id \otimes \omega_{f_{(1)}^*, S(f_{(2)}) \chi})(W)\Lambda(h) &= \Lambda(S^{-1}(S(f_{(3)}) \chi_{(1)}) h) \bra f_{(1)}^*, S(f_{(2)}) \chi_{(2)} \ket \\
&= \Lambda(S^{-1}(S(f_{(3)}) \chi_{(1)}) h) \phi(f_{(1)} S(f_{(2)}) \chi_{(2)}) \\
&= \Lambda(S^{-1}(\chi_{(1)}) f h) \phi(\chi_{(2)}) \\
&= \Lambda(f h) \phi(\chi) = \lambda(f) \Lambda(h) 
\end{align*}
for all $ h \in \CF^\infty_c(G) $. In particular, $ \lambda(f) $ extends naturally to a bounded operator on $ L^2(G) $. 
It is then straightforward to check that $ \lambda $ yields in fact a faithful $ * $-representation of $ \CF^\infty_c(G) $ on $ L^2(G) $. 

Moreover we have 
$$
\Delta(f) = W^*(1 \otimes f) W
$$
for all $ f \in \CF^\infty_c(G) $, where we identify $ f $ with $ \lambda(f) \in \LH(L^2(G)) $. Indeed, one computes 
\begin{align*}
(W^* (1 \otimes f) W)(\Lambda(g) \otimes \Lambda(h)) &= (W^* (1 \otimes f))(\Lambda(S^{-1}(h_{(1)}g) \otimes \Lambda(h_{(2)})) \\
&= W^*(\Lambda(S^{-1}(h_{(1)}g)) \otimes \Lambda(fh_{(2)})) \\
&= \Lambda(f_{(1)} g) \otimes \Lambda(f_{(2)} h) \\
&= \Delta(f)(\Lambda(g) \otimes \Lambda(h)) 
\end{align*}
for all $ g, h \in \CF^\infty_c(G) $. 

We shall next identify the dual multiplier Hopf algebra $ \DF(G) $ inside $ \LH(L^2(G)) $.
Recall from Section \ref{secintegrals} that $ \CF^\infty_c(G) $ is linked with its dual $ \DF(G) $ by Fourier transform. 
We will also write $ \DF(G) = \CF^\infty_c(\hat{G})  $ for the dual. We shall sometimes refer to $\DF(G)$ as the group algebra of the algebraic quantum group $G$.

\begin{lemma} \label{Fouriertransform}
Let $ G $ be an algebraic quantum group. The Fourier transform $ \F: \CF^\infty_c(G) \rightarrow \DF(G) $ given by $ \F(f)(h) = \phi(hf) $ 
induces an isometric linear isomorphism $ L^2(G) \rightarrow L^2(\hat{G}) $. 
\end{lemma} 

\begin{proof} For $ f, g, h \in \CF^\infty_c(G) $ we have 
\begin{align*}
(\F(f)^* \F(g), h) 
&= (\F(f)^*, h_{(1)}) \phi(h_{(2)} g) \\ 
&= (\F(f)^*, h_{(1)} g_{(2)} S^{-1}(g_{(1)})) \phi(h_{(2)} g_{(3)}) \\ 
&= (\F(f)^*, S^{-1}(g_{(1)})) \phi(h g_{(2)}) \\ 
&= \overline{\phi(g_{(1)}^* f)} \phi(h g_{(2)}) \\ 
&= \phi(f^*g_{(1)}) (\F(g_{(2)}), h). 
\end{align*}
Hence we obtain 
$$
\hat{\phi}(\F(f)^* \F(g)) = \phi(f^*g_{(1)}) \hat{\phi}(\F(g_{(2)})) = \phi(f^* g) 
$$
for all $ f,g \in \CF^\infty_c(G) $. This shows that $ \F $ extends to an isometric isomorphism with respect to the canonical scalar products.
\end{proof}

Using the Fourier transform from Lemma \ref{Fouriertransform} we can transport the left regular representation of $ \CF^\infty_c(\hat{G}) = \DF(G) $ 
on $ L^2(\hat{G}) $ to $ L^2(G) $ as follows. 
Define a linear map $ \hat{\lambda} $ from $ \DF(G) $ to the space of linear endomorphisms of $ \CF^\infty_c(G) $, viewed 
as a subspace of $ L^2(G) $, by the formula 
$$
\hat{\lambda}(x) \Lambda(f) = (\hat{S}(x), f_{(1)}) \Lambda(f_{(2)}) = (x, S^{-1}(f_{(1)})) \Lambda(f_{(2)}),
$$
\nomenclature[o$\lambda$2]{$\hat{\lambda}$}{representation $ \hat{\lambda}: \DF(G) \rightarrow \LH(L^2(G)) $ of the group algebra of an algebraic quantum group on $L^2(G)$}%
for $x\in\DF(G)$, $f\in\CF^\infty_c(G)$.
Then for all $h\in\CF^\infty_c(G)$ we have 
\begin{align*}
(\F(\hat{\lambda}(x) \Lambda(f)), h) &= (x, S^{-1}(f_{(1)})) \phi(h f_{(2)}) \\ 
&= (x, S^{-1}(S(h_{(1)}) h_{(2)}f_{(1)})) \phi(h_{(3)} f_{(2)}) \\ 
&= (x, h_{(1)}) \phi(h_{(2)} f) \\ 
&= (x, h_{(1)}) (\F(f), h_{(2)}) \\ 
&= (x \F(f), h), 
\end{align*}
which means that $ \F \hat{\lambda}(x) \F^{-1} $ corresponds to the GNS-representation of $ \CF^\infty_c(\hat{G}) $ on 
$ L^2(\hat{G}) $. In particular, we obtain a faithful $ * $-representation $ \hat{\lambda}: \DF(G) \rightarrow \LH(L^2(G)) $ using the above construction. 

In terms of the multiplicative unitary $ W $, the comultiplication $ \hat{\Delta} $ for $ \DF(G) $ is determined by the formula
$$
\hat{\Delta}(x) = \hat{W}^*(1 \otimes x) \hat{W}
$$
where $ \hat{W} = \Sigma W^* \Sigma $, and we identify $ x $ with $ \hat{\lambda}(x) \in \LH(L^2(G)) $.

Inspecting the above formulas we see that we obtain Hopf $ C^* $-algebras $ C^\red_0(G) $ and $ C^*_\red(G) = C^\red_0(\hat{G}) $ inside $ \LH(L^2(G)) $ 
by taking the closures of $ \lambda(\CF^\infty_c(G)) $ and $ \hat{\lambda}(\DF(G)) $. These algebras identify with 
the legs of the multiplicative unitary $ W $. Moreover, these constructions are compatible with the multiplier Hopf algebra structures 
of $ \CF^\infty_c(G) $ and $ \DF(G) $, respectively.  
In a similar way one obtains von Neumann algebras $ L^\infty(G) $ and $ \L(G) $ with comultiplications by taking the weak closures 
of $ \lambda(\CF^\infty_c(G)) $ and $ \hat{\lambda}(\DF(G)) $, compare the constructions in Section \ref{seclcqg}. 

The key result due to Kustermans and van Daele is that these operator algebras define a locally compact quantum group with multiplicative unitary $ W $, 
see Section 6 in \cite{KvD}. Let us phrase this as follows. 

\begin{theorem} 
Let $ G $ be an algebraic quantum group. 
With the notation as above, the left/right invariant integrals on $ \CF^\infty_c(G) $ extend to left/right invariant weights on $ L^\infty(G) $. 
In particular, $ G $ canonically defines a locally compact quantum group.  
\end{theorem} 

To conclude this subsection, let us verify that the  unitary $ \hat{W} = \Sigma W^* \Sigma $ indeed corresponds to the fundamental multiplicative 
unitary for $ \hat{G} $ under Fourier transform. 
For any $f,g \in \CF^\infty_c(G)$, we have 
\begin{align*}
(\F \otimes \F)(\Sigma W^* \Sigma)(\Lambda(f) \otimes \Lambda(g)) 
&= (\F \otimes \F)(\Lambda(f_{(2)}) \otimes \Lambda(f_{(1)}g)) \\ 
&= \hat{\Lambda}(\F(f_{(2)})) \otimes \hat{\Lambda}(\F(f_{(1)}g)) . 
\end{align*}
Now, for any $h,k\in\CF^\infty_c(G)$ we can compute
\begin{align*}
(\F(f_{(2)}) \otimes \F(f_{(1)}g), h \otimes k) &= \phi(h f_{(2)}) \phi(kf_{(1)} g) \\ 
&= \phi(h_{(3)} f_{(2)}) \phi(kS(h_{(1)}) h_{(2)} f_{(1)} g) \\ 
&= \phi(k S(h_{(1)}) g) \phi(h_{(2)} f) \\
&= (\F(g), k S(h_{(1)})) (\F(f), h_{(2)}) \\
&= (\F(g)_{(1)}, S(h_{(1)})) (\F(f), h_{(2)}) (\F(g)_{(2)}, k) \\
&= (\hat{S}^{-1}(\F(g)_{(1)}) \F(f) \otimes \F(g)_{(2)}, h \otimes k) ,
\end{align*}
and therefore the previous calculation gives
\begin{align*}
(\F \otimes \F)(\Sigma W^* \Sigma)(\Lambda(f) \otimes \Lambda(g)) 
&= \hat{\Lambda}(\F(f_{(2)})) \otimes \hat{\Lambda}(\F(f_{(1)}g)) \\ 
& = (\hat{\Lambda} \otimes \hat{\Lambda})(\hat{S}^{-1}(\F(g)_{(1)}) \F(f) \otimes \F(g)_{(2)})  \\ 
&  = \hat{W}(\hat{\Lambda}(\F(f)) \otimes \hat{\Lambda}(\F(g))).
\end{align*}

This yields the claim.

\subsubsection{Compact quantum groups} 
\label{sec:compact_quantum_groups}

In this subsection we briefly sketch the theory of compact quantum groups. For more information we refer to \cite{Woronowiczleshouches}, \cite{MvDnotes} 
and \cite{KS}. 

There are various ways in which the concept of a compact quantum group can be defined. In our set-up, it is convenient to consider 
compact quantum groups as a special case of algebraic quantum groups as in Definition \ref{defaqg}. Historically, the development took place in the 
opposite order, in fact, the invention of algebraic quantum groups was strongly motivated by the theory of compact quantum groups 
and attempts to generalize it, see \cite{vDadvances}. 

\begin{definition} 
A compact quantum group is an algebraic quantum group $ H $ such that the underlying algebra of $ H $ is unital. 
\end{definition} 

We shall also write $ H = \CF^\infty(K) $ in this case and refer to $ K $ as a compact quantum group. 
\nomenclature[o$C^\infty(K)$]{$\CF^\infty(K)$}{algebra of polynomial functions on a compact quantum group}%
Note that the comultiplication is a $ * $-homomorphism $ \Delta: H \rightarrow H \otimes H $, so that $ H $ 
is in particular a Hopf $ * $-algebra. 

Moreover, by definition there exists a positive left invariant integral $ \phi: \CF^\infty(K) \rightarrow \mathbb{C} $ 
and a positive right invariant integral $ \psi: \CF^\infty(K) \rightarrow \mathbb{C} $ such that $ \phi(1) = 1 = \psi(1) $.
\label{nom:Haar-compact}%
We have 
$$
\phi(f) = \phi(f) \psi(1) = \psi((\id \otimes \phi)\Delta(f)) = (\psi \otimes \phi)\Delta(f) 
= \psi(f) \phi(1) = \psi(f) 
$$
for all $ f \in H $, so that in fact $ \phi = \psi $. We refer to this left and right invariant functional as the Haar state of $ \CF^\infty(K) $. 

In particular, due to Proposition \ref{cosemisimpleintegral} the Hopf algebra $ \CF^\infty(K) $ is cosemisimple. 
That is, we can write 
$$
\CF^\infty(K) \cong \bigoplus_{\lambda \in \Lambda} M_{n_\lambda}(\mathbb{C})^* 
$$
as a direct sum of simple matrix coalgebras. 
Moreover, this isomorphism can be chosen so that the standard matrix coefficient functionals $u^\lambda_{ij} \in M_{n_\lambda}(\mathbb{C})^* $
\label{nom:matrix_coeff_general}%
satisfy
$ (u^\lambda_{ij})^* = S(u^\lambda_{ji}) $ 
for all $ \lambda \in \Lambda $. Equivalently, the matrix $ u^\lambda = (u^\lambda_{ij}) \in M_{n_\lambda}(\CF^\infty(K)) $ is unitary for all $\lambda\in\Lambda$.

For each $ \lambda \in \Lambda $ there exists a unique positive invertible matrix $ F_{\lambda} \in M_{n_\lambda}(\mathbb{C}) $ 
\label{nom:F_general}%
such that $S^2(u^\lambda) = F_\lambda u^\lambda F_\lambda^{-1}$
and $ \tr(F_\lambda) = \tr(F_\lambda^{-1}) $. Here we write $ S^2(u^\lambda) $ for the matrix obtained by applying $ S^2 $ entrywise to $ u^\lambda $, 
and we consider the unnormalized standard trace $ \tr $ on $ M_{n_\lambda}(\mathbb{C}) $. 
\label{nom:tr}%

If we fix matrix coefficients $ u^\lambda_{ij} $ as above then we have the Schur orthogonality relations 
$$
\phi(u^\beta_{ij} (u^\gamma_{kl})^*) = \delta_{\beta \gamma} \delta_{ik} \, \frac{(F_\beta)_{lj}}{\tr(F_\beta)}, 
\qquad 
\phi((u^\beta_{ij})^* u^\gamma_{kl}) = \delta_{\beta \gamma} \delta_{jl} \, \frac{(F_\beta^{-1})_{ki}}{\tr(F_\beta)}, 
$$
compare for instance Chapter 11 in \cite{KS}. 

The Schur orthogonality relations imply that the dual algebraic quantum group $ \DF(K) $ can be written 
as a direct sum 
$$
\DF(K) = \bigoplus_{\lambda \in \Lambda} M_{n_\lambda}(\mathbb{C}) 
$$
of matrix algebras, such that the pairing $ (x, f) $ for $ x \in \DF(K), f \in \CF^\infty(K) $ is given by evaluation in each component. 
As usual, we write $ x \hit h = (x, h_{(2)}) h_{(1)} $ and $ h \hitby x = (x, h_{(1)}) h_{(2)} $
\label{nom:hit_DK}%
for the natural left and right actions of $\DF(K)$ on $\CF^\infty_c(K)$.  If we define the element
\[
  F = \bigoplus_{\lambda\in\Lambda} F_\lambda \in\M(\DF(K)),
\] 
\label{nom:F_general}%
where the matrices $F_\lambda$ are those from the Schur orthogonality formulas above, then the Haar state satisfies
\[
  \phi(fg) 
 	= \phi(g(F\hit f \hitby F))
 	= \phi((F^{-1}\hit g \hitby F^{-1}) f)
\]
for all $ f, g \in \CF^\infty(K) $.

We obtain a basis of $ \DF(K) $ consisting of the functionals $ \omega^{\mu}_{ij} $ defined by 
$$
(\omega^\mu_{ij}, u^\eta_{kl}) = \delta_{\mu \eta} \delta_{ik} \delta_{jl}.  
$$
\label{nom:dual_basis_K}%
With this notation, the multiplicative unitary $ W \in \M(\CF^\infty(K) \otimes \DF(K)) $ can be written as 
$$
W = \sum_{\mu, i,j} u^\mu_{ij} \otimes \omega^\mu_{ij}.  
$$
Indeed, we compute 
\begin{align*}
\sum_{\mu, i,j} (\lambda(u^\mu_{ij}) \otimes \hat{\lambda}(\omega^\mu_{ij}))(\Lambda(f) \otimes \Lambda(g)) 
&= \sum_{\mu, i,j} (\omega^\mu_{ij}, S^{-1}(g_{(1)})) \Lambda(u^\mu_{ij} f) \otimes \Lambda(g_{(2)}) \\
&= \Lambda(S^{-1}(g_{(1)}) f) \otimes \Lambda(g_{(2)}) 
\end{align*}
for all $ f, g \in \CF^\infty(K) $. In a similar way one obtains 
$$
W^{-1} = \sum_{\mu, i,j} S(u^\mu_{ij}) \otimes \omega^\mu_{ij} = \sum_{\mu, i,j} u^\mu_{ij} \otimes \hat{S}^{-1}(\omega^\mu_{ij})  
$$
for the inverse of $ W $. 

Positive left and right invariant Haar functionals for $ \DF(K) $ are given by 
$$
\hat{\phi}(x) = \sum_{\mu \in \Lambda} \tr(F_\mu) \tr(F_\mu^{-1} x), 
\qquad \hat{\psi}(x) = \sum_{\mu \in \Lambda} \tr(F_\mu) \tr(F_\mu x), 
$$
\label{nom:dual_Haar_weights}%
respectively. Note here that the positive matrices $ F_\lambda \in M_{n_\lambda}(\mathbb{C}) $  are naturally elements of $ \DF(K) $, and that $ F_\lambda^{\pm 1} x $ is contained 
in $ M_{n_\lambda}(\mathbb{C}) \subset \DF(K) $ for any $ x \in \DF(K) $. 
We note the formulas
\begin{align*}
 \hat{\phi}(xy) &= \hat{\phi}(F y F^{-1}x),
 & \hat{\psi}(xy) &= \hat{\psi}(F^{-1} y F x),
\end{align*}
for all $x,y\in\DF(K)$.
We remark finally that $ \hat{\delta} = F^2 $ is the modular element of the dual quantum group $ \DF(K) $.
In particular, $ \DF(K) $ is unimodular iff $ F = 1 $.

\subsubsection{The Drinfeld double of algebraic quantum groups} 
\label{sec:AQC-double}

In this subsection we discuss the Drinfeld double construction in the framework of algebraic quantum groups. The Drinfeld double 
of regular multiplier Hopf algebras was already treated in Section \ref{secddalgebraic}. Here we shall explain how 
to incorporate $ * $-structures in the construction, and also approach it from the dual point of view. 

Let $ K $ and $ L $ be algebraic quantum groups with group algebras $ \DF(K) $ and $ \DF(L) $, respectively. 
Recall from Section \ref{secddalgebraic} that in order to form the Drinfeld double $ \DF(K) \bowtie \DF(L) $ 
one needs a skew-pairing $ \tau: \DF(K) \times \DF(L) \rightarrow \mathbb{C} $. We shall say that the skew-pairing $ \tau $ is unitary if 
$$
\tau(x^*, y^*) = \overline{\tau^{-1}(x, y)} 
$$
for all $ x \in \DF(K), y \in \DF(L) $, where $\tau^{-1}$ denotes the convolution inverse of $\tau$ as in Section \ref{secddalgebraic}. 

Let us also introduce the notion of a unitary bicharacter. 
\begin{definition} \label{defbicharacter}
Let $ K $ and $ L $ be algebraic quantum groups. A (unitary) bicharacter for $ K,L $ is a (unitary) invertible 
element $ U \in \M(\CF^\infty_c(K) \otimes \CF^\infty_c(L)) $ such that 
$$
(\Delta_K \otimes \id)(U) = U_{13} U_{23}, \qquad (\id \otimes \Delta_L)(U) = U_{13} U_{12}
$$
and 
$$
(\epsilon_K \otimes \id)(U) = 1, \qquad (\id \otimes \epsilon_L)(U) = 1. 
$$
\end{definition} 

In the same way as in the discussion of universal $ R $-matrices in Section \ref{subsecuniversalrmatrix} one checks that a bicharacter $ U $ satisfies 
$$
(S_K \otimes \id)(U) = U^{-1} = (\id \otimes S_L^{-1})(U).  
$$

The notion of a bicharacter is dual to the concept of a skew-pairing in the following sense. 
\begin{prop} \label{bicharacterchar}
Assume that $ K $ and $ L $ are algebraic quantum groups. If $ U \in \M(\CF^\infty_c(K) \otimes \CF^\infty_c(L)) $ is a bicharacter
then $ \tau_U: \DF(K) \times \DF(L) \rightarrow \mathbb{C} $ given by 
$$ 
\tau_U(x, y) = (x \otimes y, U^{-1}) 
$$
\nomenclature{$\tau_U$}{skew-pairing associated to a bicharacter $U$}%
is a skew-pairing. Every skew-pairing 
$ \DF(K) \times \DF(L) \rightarrow \mathbb{C} $ arises in this way from a bicharacter. 
Moreover $ U $ is unitary iff the skew-pairing $ \tau_U $ is unitary. 
\end{prop}   

\begin{proof} Let us sketch the argument. If $ U $ is a bicharacter then one can check in the same way as in the proof of Lemma \ref{skewpairingduality} that 
$ \tau_U $ yields a skew-pairing. 

To see that every skew-pairing arises from a bicharacter let $ \tau: \DF(K) \times \DF(L) \rightarrow \mathbb{C} $ be given. 
Using regularity one checks that the formulas
\begin{align*} 
(x \otimes y, U^{-1}_\tau(f \otimes g)) &= \tau(x_{(2)}, y_{(2)}) (x_{(1)}, f) (y_{(1)}, g) \\
(x \otimes y, (f \otimes g) U^{-1}_\tau) &= (x_{(2)}, f) (y_{(2)}, g) \tau(x_{(1)}, y_{(1)}) 
\end{align*}
determine an invertible multiplier $ U_\tau \in \M(\CF^\infty_c(K) \otimes \CF^\infty_c(L)) $. 
Moreover, one checks that this multiplier 
satisfies the conditions in Definition \ref{defbicharacter}. The skew-pairing $ \tau $ is reobtained by applying the above construction to $ U_\tau $. 

For the last claim use the relation $ (S_K \otimes S_L)(U) = U $ 
to compute 
$$
\tau_U(x^*, y^*) = (x^* \otimes y^*, U^{-1}) = \overline{(x \otimes y, (S_K \otimes S_L)(U^{-1})^*)}
= \overline{(x \otimes y, (U^{-1})^*)}
$$
and 
$$
\overline{\tau_U^{-1}(x, y)} = \overline{(x \otimes y, U)}.  
$$
Comparing these expressions yields the assertion. 
\end{proof} 

Let $ K $ and $ L $ be algebraic quantum groups and assume $ U \in \M(\CF^\infty_c(K) \otimes \CF^\infty_c(L)) $ is a bicharacter. 
As in Section \ref{secddalgebraic} we can then construct the Drinfeld double $ \DF(K) \bowtie \DF(L) $ using the skew-pairing $ \tau_U $ 
from Proposition \ref{bicharacterchar}. 
Recall that this is the regular multiplier Hopf algebra 
$$ 
\DF(K \bowtie L) = \DF(K) \bowtie \DF(L), 
$$ 
\label{nom:alg_double}%
equipped with the tensor product comultiplication and the multiplication given by 
\begin{align*}
(x \bowtie f)(y \bowtie g) &= x\, \tau_U(y_{(1)}, f_{(1)}) y_{(2)} \bowtie f_{(2)} \tau_U (\hat{S}_K(y_{(3)}), f_{(3)}) g \\
&= x \,\tau_U(y_{(1)}, f_{(1)}) y_{(2)} \bowtie f_{(2)} \tau_U(y_{(3)}, \hat{S}_L^{-1}(f_{(3)})) g. 
\end{align*}
The counit of $ \DF(K \bowtie L) $ is given by 
$$
\hat{\epsilon}_{K \bowtie L}(x \bowtie f) = \hat{\epsilon}_K(x) \hat{\epsilon}_L(f) 
$$
for $ x \in \DF(K), f \in \DF(L) $. The antipode of $ \DF(K \bowtie L) $ is defined by 
\begin{align*}
\hat{S}_{K \bowtie L}(x \bowtie f) &= (1 \bowtie \hat{S}_L(f))(\hat{S}_K(x) \bowtie 1) \\
&= \tau_U(\hat{S}(x_{(3)}), \hat{S}_L(f_{(3)})) \hat{S}_K(x_{(2)}) \bowtie \hat{S}_L(f_{(2)}) \tau_U(\hat{S}_K(x_{(1)}), f_{(1)}). 
\end{align*}
If $ U $ is unitary we can define a $ * $-structure on $ \DF(K \bowtie L) $ by 
$$
(x \bowtie f)^* 
 = \tau_U(x^*_{(1)}, f^*_{(1)}) x^*_{(2)} \bowtie f^*_{(2)} \tau_U^{-1}(x^*_{(3)}, f^*_{(3)})
  = (1 \bowtie f^*)(x^* \bowtie 1) .  
$$
To check antimultiplicativity of this $ * $-structure one computes 
\begin{align*}
((x \bowtie f)&(y \bowtie g))^* = (x \tau_U(y_{(1)}, f_{(1)}) y_{(2)} \bowtie f_{(2)} \tau_U^{-1}(y_{(3)}, f_{(3)}) g)^* \\ 
&= \overline{\tau_U(y_{(1)}, f_{(1)})} (xy_{(2)} \bowtie f_{(2)} g)^* \overline{\tau_U^{-1}(y_{(3)}, f_{(3)})} \\ 
&= \tau_U^{-1}(y_{(1)}^*, f_{(1)}^*) (xy_{(2)} \bowtie f_{(2)} g)^* \tau_U(y_{(3)}^*, f_{(3)}^*) \\ 
&= \tau_U^{-1}(y_{(1)}^*, f_{(1)}^*) (1 \bowtie g^*) \tau_U(y^*_{(2)}, f^*_{(2)}) (y^*_{(3)} \bowtie f^*_{(3)}) \\
&\qquad \times \tau_U^{-1}(y^*_{(4)}, f^*_{(4)}) (x^* \bowtie 1) \tau_U(y_{(5)}^*, f_{(5)}^*) \\ 
&= (1 \bowtie g^*)(y^* \bowtie f^*)(x^* \bowtie 1) \\
&= (1 \bowtie g^*)(y^* \bowtie 1) (1 \bowtie f^*)(x^* \bowtie 1) \\
&= (y \bowtie g)^* (x \bowtie f)^*, 
\end{align*}
where we use
$$
\overline{\tau_U(x, f)} = \tau_U^{-1}(x^* \otimes f^*). 
$$
For involutivity note that 
$$ 
(x \bowtie f)^{**} = ((1 \bowtie f^*)(x^* \bowtie 1))^* = x \bowtie f 
$$ 
using that $ * $ is antimultiplicative. 
Similarly, to check that $ \hat{\Delta}_{K \bowtie L} $ is a $ * $-homomorphism one calculates 
\begin{align*}
\hat{\Delta}_{K \bowtie L}((x \bowtie f)^*) &= ((1 \bowtie f^*_{(1)}) \otimes (1 \bowtie f^*_{(2)}))((x_{(1)}^* \bowtie 1) \otimes (x_{(2)}^* \bowtie 1)) \\
&= (x^*_{(1)} \bowtie f^*_{(1)}) \otimes (x_{(2)}^* \bowtie f_{(2)}^*) \\
&= \hat{\Delta}_{K \bowtie L}(x \bowtie f)^*. 
\end{align*}

Finally, a left Haar integral for $ \DF(K \bowtie L) $ is given by 
$$ 
\hat{\phi}_{K \bowtie L}(x \bowtie f) = \hat{\phi}_K(x) \hat{\phi}_L(f), 
$$ 
where $ \hat{\phi}_K $ and $ \hat{\phi}_L $ are left Haar integrals for $ \DF(K) $ and $ \DF(L) $, respectively. 
Similarly, a right Haar integral for $ \DF(K \bowtie L) $ is given by the tensor product of right Haar integrals for $ \DF(K) $ and $ \DF(L) $. 

If $ L $ is a discrete quantum group then it is not hard to check that the resulting functionals are positive provided one starts 
from positive integrals for $ K $ and $ L $, respectively. For the question of positivity in general see \cite{DvDdouble}. 

From general theory, we obtain the dual multiplier Hopf algebra $ \CF^\infty_c(K \bowtie L) $ of $ \DF(K \bowtie L) $. 
\label{nom:codouble}%
Sometimes this is referred to as the \emph{Drinfeld codouble}, but we shall call $ \CF^\infty_c(K \bowtie L) $ the algebra of 
functions on the Drinfeld double $ K \bowtie L $. Explicitly, the structure of $ \CF^\infty_c(K \bowtie L) $ looks as follows. 

\begin{prop} \label{doublestructure}
Let $ K $ and $ L $ be algebraic quantum groups and assume that $ U \in \M(\CF^\infty_c(K) \otimes \CF^\infty_c(L)) $ is a unitary bicharacter. 
Then the algebra 
$$
\CF^\infty_c(K \bowtie L) = \CF^\infty_c(K) \otimes \CF^\infty_c(L),  
$$
equipped with the comultiplication 
$$
\Delta_{K \bowtie L} = (\id \otimes \sigma \otimes \id)(\id \otimes \ad(U) \otimes \id)(\Delta_K \otimes \Delta_L)
$$
is an algebraic quantum group.  
Moreover, the counit of $ \CF^\infty(K \bowtie L) $ is the tensor product counit $ \epsilon_{K \bowtie L} = \epsilon_K \otimes \epsilon_L $, and the 
antipode is given by 
$$
S_{K \bowtie L}(f \otimes x) = U^{-1} (S_K(f) \otimes S_L(x)) U = (S_K \otimes S_L)(U (f \otimes x) U^{-1})
$$
The algebraic quantum group $ \CF^\infty_c(K \bowtie L) $ is dual to the double $ \DF(K) \bowtie \DF(L) $, the latter being constructed with respect to the 
skew-pairing $ \tau_U $ corresponding to $ U $. 
\end{prop} 

\begin{proof}
Note first that $ \ad(U) $ is conjugation with the bicharacter $ U $ and $ \sigma $ denotes the flip map 
in the above formula for the comultiplication.  

Let us abbreviate $ G = K \bowtie L $. In order to prove the Proposition we shall show that the dual of $ \DF(G) = \DF(K) \bowtie \DF(L) $ 
can be identified as stated. 

Firstly, using $ \hat{\phi}_G(x \bowtie f) = \hat{\phi}_K(x) \hat{\phi}_L(f) $ we compute 
\begin{align*} 
\hat{\F}_G&(y \bowtie g)(x \bowtie f) = \hat{\phi}_G(x \tau_U(y_{(1)}, f_{(1)}) y_{(2)} \bowtie f_{(2)} \tau_U^{-1}(y_{(3)}, f_{(3)}) g) \\
&= \tau_U(y_{(1)}, f_{(1)}) \hat{\phi}_K(xy_{(2)}) \hat{\phi}_L(f_{(2)} g) \tau_U(\hat{S}_K(y_{(3)}), f_{(3)}) \\
&= \tau_U(y_{(1)}, \hat{S}_L^{-1}(g_{(1)})) \hat{\phi}_K(xy_{(2)}) \hat{\phi}_L(f_{(1)} g_{(2)}) \tau_U(\hat{S}_K(y_{(3)}), f_{(2)}) \\
&= \tau_U(y_{(1)}, \hat{S}_L^{-1}(g_{(1)})) \hat{\phi}_K(xy_{(2)}) \hat{\phi}_L(f_{(1)} g_{(2)}) \tau_U(y_{(3)}, S_L^{-1}(f_{(2)} g_{(3)} S_L(g_{(4)}))) \\
&= \tau_U(y_{(1)}, \hat{S}_L^{-1}(g_{(1)})) \hat{\phi}_K(xy_{(2)}) \hat{\phi}_L(f g_{(2)}) \tau_U(y_{(3)}, g_{(3)} \hat{\delta}_L^{-1}) \\
&= (x \otimes f, \tau_U^{-1}(y_{(1)}, g_{(1)}) \hat{\F}_K(y_{(2)}) \otimes \hat{\F}_L(g_{(2)}) \tau_U(y_{(3)}, g_{(3)} \hat{\delta}_L^{-1})) 
\end{align*} 
for all $ x,y \in \DF(K), f, g, \in \DF(L) $. Here $ \hat{\delta}_L $ denotes the modular element of $ \DF(L) $. 
This implies 
$$
\hat{\F}_G(\tau_U(y_{(1)}, g_{(1)}) y_{(2)} \bowtie g_{(2)} \tau_U^{-1}(y_{(3)}, g_{(3)} \hat{\delta}_L^{-1})) = \hat{\F}_K(y) \otimes \hat{\F}_L(g). 
$$
In particular, we can identify the underlying vector space of the dual quantum group of $ \DF(G) $ 
with the space $ \CF^\infty_c(G) = \CF^\infty_c(K) \otimes \CF^\infty_c(L) $ 
such that the canonical pairing between $ \DF(G) $ and $ \CF^\infty_c(G) $ becomes 
$$
(y \bowtie g, f \otimes x) = (y,f) (g,x) 
$$
for all $ y \in \DF(K), f \in \CF^\infty_c(K), g \in \DF(L), x \in \CF^\infty_c(L) $. 
Moreover, from the definition of the comultiplication in $ \DF(G) $ it is clear that the algebra structure 
of the dual multiplier Hopf algebra $ \CF^\infty_c(G) = \CF^\infty_c(K) \otimes \CF^\infty_c(L)  $ 
is the tensor product algebra structure. 

Let us next identify the comultiplication of $ \CF^\infty_c(G) $. By general theory, this is determined by 
\begin{align*}
((y \bowtie g) & \otimes (z \bowtie h), \Delta_G(f \otimes x)) = 
(y \tau_U(z_{(1)}, g_{(1)}) z_{(2)} \bowtie g_{(2)} \tau_U^{-1}(z_{(3)}, g_{(3)}) h, f \otimes x) \\
&= (y \otimes \tau_U(z_{(1)}, g_{(1)}) z_{(2)} \bowtie g_{(2)} \otimes \tau_U^{-1}(z_{(3)}, g_{(3)}) h, (\Delta_K \otimes \Delta_L)(f \otimes x)) \\
&= (y \otimes (z_{(1)} \otimes g_{(1)}, U^{-1}) z_{(2)} \bowtie g_{(2)} \otimes (z_{(3)} \otimes g_{(3)}, U) h, (\Delta_K \otimes \Delta_L)(f \otimes x)) \\
&= (y \otimes z \otimes g \otimes h, (\id \otimes \ad(U) \otimes \id)(\Delta_K \otimes \Delta_L)(f \otimes x)) \\
&= ((y \bowtie g) \otimes (z \bowtie h), (\id \otimes \sigma \otimes \id)(\id \otimes \ad(U) \otimes \id)(\Delta_K \otimes \Delta_L)(f \otimes x)),  
\end{align*}
so that we obtain 
$$
\Delta_G = (\id \otimes \sigma \otimes \id)(\id \otimes \ad(U) \otimes \id)(\Delta_K \otimes \Delta_L) 
$$
as claimed. 
Remark that it is immediate from unitarity of $ U $ that $ \Delta_G: \CF^\infty_c(G) \rightarrow \M(\CF^\infty_c(G) \otimes \CF^\infty_c(G)) $ 
is indeed an essential $ * $-algebra homomorphism. 

The formula $ \epsilon_G = \epsilon_K \otimes \epsilon_L $ for the counit of $ \CF^\infty_c(G) $ can also be deduced from duality. 
Alternatively, we may check the counit property of $ \epsilon_G $ directly and compute
\begin{align*}
(\epsilon_G \otimes \id)\Delta_G
&= (\epsilon_K \otimes \epsilon_L \otimes \id \otimes \id)(\id \otimes \sigma \otimes \id)(\id \otimes \ad(U) \otimes \id)(\Delta_K \otimes \Delta_L) \\
&= (\epsilon_K \otimes \id \otimes \epsilon_L \otimes \id)(\id \otimes \ad((\id \otimes \epsilon_L)(U)) \otimes \id)(\Delta_K \otimes \Delta_L) \\
&= (\epsilon_K \otimes \id \otimes \epsilon_L \otimes \id)(\Delta_K \otimes \Delta_L) \\
&= (\id \otimes \epsilon_K \otimes \id \otimes \epsilon_L)(\Delta_K \otimes \Delta_L) \\
&= (\id \otimes \epsilon_K \otimes \id \otimes \epsilon_L)(\id \otimes \ad((\epsilon_K \otimes \id)(U)) \otimes \id)(\Delta_K \otimes \Delta_L) \\
&= (\id \otimes \id \otimes \epsilon_K \otimes \epsilon_L)(\id \otimes \sigma \otimes \id)(\id \otimes \ad(U) \otimes \id)(\Delta_K \otimes \Delta_L) \\
&= (\id \otimes \epsilon_G)\Delta_G, 
\end{align*}
using the relations $ (\epsilon_K \otimes \id)(U) = \id $ and $ (\id \otimes \epsilon_L)(U) = \id $. 

To verify the formula for the antipode we check 
\begin{align*}
&m_G(S_G \otimes \id) \Delta_G = m_G(S_G \otimes \id)(\id \otimes \sigma \otimes \id)(\id \otimes \ad(U) \otimes \id)(\Delta_K \otimes \Delta_L) \\
&=(m_K \otimes m_L)(S_K \otimes \id \otimes S_L \otimes \id) \ad(U_{13} U_{23})(\Delta_K \otimes \Delta_L) \\
&=(m_K \otimes m_L)(S_K \otimes \id \otimes S_L \otimes \id) \ad((\Delta_K \otimes \id)(U)_{123})(\Delta_K \otimes \Delta_L) \\
&= m_L(S_L \otimes \id) (\epsilon_K \otimes \Delta_L) \\
&= \epsilon_K \otimes \epsilon_L,  
\end{align*} 
using the antipode axioms for $ K $ and $ L $, respectively. A similar computation shows 
\begin{align*}
&m_G(\id \otimes S_G) \Delta_G = m_G(\id \otimes S_G)(\id \otimes \sigma \otimes \id)(\id \otimes \ad(U) \otimes \id)(\Delta_K \otimes \Delta_L) \\
&=(m_K \otimes m_L)(\id \otimes S_K \otimes \id \otimes S_L) \ad(U_{24} U_{23})(\Delta_K \otimes \Delta_L) \\
&=(m_K \otimes m_L)(\id \otimes S_K \otimes \id \otimes S_L) \ad((\id \otimes \Delta_L)(U)_{234})(\Delta_K \otimes \Delta_L) \\
&= m_K(\id \otimes S_K) (\Delta_K \otimes \epsilon_L) \\
&= \epsilon_K \otimes \epsilon_L. 
\end{align*} 

According to \cite{DvDdouble} and the duality theory for algebraic quantum groups, there exists a positive left invariant functional 
and a positive right invariant functional on $ \CF^\infty_c(G) $. Hence $ \CF^\infty_c(G) $ is again an algebraic quantum group. 
\end{proof} 

Explicitly, the left Haar functional for $ G = K \bowtie L $ are obtained by combining the left Haar functional for $ K $ and a twisted version of the 
left Haar functional for $ L $, depending on the modular properties of the quantum groups and the pairing involved. We refer to \cite{BV} for a detailed 
analysis. 

In the special case of the Drinfeld double of a compact quantum group we will write down an explicit formula for a (left and right) invariant 
integral further below.

\subsection{Compact semisimple quantum groups} 
\label{sec:compact_ssqgs}

Let $ \mathfrak{g} $ be a semisimple Lie algebra. In order to pass from the algebraic theory of quantized universal enveloping algebras developed 
in Chapter \ref{chuqg} to the analytical setting one has to introduce a $ * $-structure on $ U_q(\mathfrak{g}) $ and $ \Poly(G_q) $. 

In the case of $ U_q(\mathfrak{g}) $ we shall work with the $ * $-structure defined as follows. 

\begin{lemma} \label{defuqgstar}
Assume $ q \in \mathbb{C}^\times $ is real and $ q \neq \pm 1 $. Then there is a unique $ * $-structure on $ U_q(\mathfrak{g}) $ satisfying
\begin{align*}
E_i^* = K_i F_i, \qquad F_i^* = E_i K_i^{-1}, \qquad K_\lambda^* = K_\lambda
\end{align*}
for $ i = 1, \dots, N $ and $ \lambda \in \weights $.   The comultiplication $\Delta$ is a $*$-homomorphism for this $*$-structure.
\end{lemma} 

\begin{proof} This is done in the same way as for the definition of the algebra antiautomorphism $ \tau $ in Lemma \ref{deftau}. Note 
that $ * $ acts in the same way on generators. The only difference 
is that $ * $ is extended anti-linearly, whereas $ \tau $ is extended linearly to general elements.
\end{proof}

\begin{definition}
	We will write $U_q^\mathbb{R}(\lie{k})$ to signify $U_q(\lie{g})$ with the Hopf $*$-structure of the lemma above.
	\nomenclature{$U_q^\mathbb{R}(\lie{k})$}{quantized enveloping algebra of the compact form $K$ of $G$}%
\end{definition}

Here,
$ U_q^\mathbb{R}(\mathfrak{k}) $ should be viewed as the universal enveloping algebra of the \emph{complexification} of the Lie 
algebra $ \mathfrak{k} $ of the compact real form $ K $ of the simply connected group $ G $ corresponding to $ \mathfrak{g} $. 

We have the following compatibility of the $ \R $-matrix with the $ * $-structure. 

\begin{lemma} \label{rmatrixreal}
The universal $ R $-matrix of $ U_q(\mathfrak{g}) = U_q^\mathbb{R}(\mathfrak{k}) $ satisfies 
$$
(\hat{S} \otimes \hat{S})(\R^*) = \R_{21}, 
$$
where $ \R_{21} $ is obtained from $ \R $ by flipping the tensor factors. 
\end{lemma}

\begin{proof} Due to Theorem \ref{Rmatrixformula} the universal $ R $-matrix is the product of the 
Cartan part $ q^{\sum_{i,j = 1}^N B_{ij}(H_i \otimes H_j)} $ 
and the nilpotent part $ \prod_{\alpha \in {\bf \Delta}^+} \exp_{q_\alpha}((q_\alpha - q_\alpha^{-1})(E_\alpha \otimes F_\alpha)) $. Since 
the antipode and the $ * $-structure are both antimultiplicative, it suffices to show that applying the $ * $-structure followed by 
$ \hat{S} \otimes \hat{S} $ to each individual factor of these elements switches legs in the tensor product. For the Cartan part this is 
obvious since $ \hat{S} \otimes \hat{S} $ introduces two minus signs which cancel out, and the $ * $-structure leaves the Cartan generators $ H_k $ 
fixed. 

For the nilpotent factors note that 
\begin{align*}
(\hat{S} \otimes \hat{S})((E_i \otimes F_i)^*) &= (\hat{S} \otimes \hat{S})(K_i F_i \otimes E_i K_i^{-1}) \\
&= (-K_i F_i K_i^{-1}) \otimes (-K_i E_i K_i^{-1}) \\
&= F_i \otimes E_i 
\end{align*}
for all $ i = 1, \dots, N $. This yields the claim for all factors corresponding to simple roots. 

For the factors corresponding to arbitrary positive roots 
we use Theorem \ref{deftiautomorphism} and the relations 
\begin{align*} 
\hat{S}(E_i^*) &= \hat{S}(K_i F_i) = (-K_i F_i) K_i^{-1} = - q^{-(\alpha_i, \alpha_i)} F_i = -q_i^{-2} F_i \\
\hat{S}(F_i^*) &= \hat{S}(E_i K_i^{-1}) = K_i(- E_i K_i^{-1}) = -q^{(\alpha_i, \alpha_i)} E_i = - q_i^2 E_i.  
\end{align*}
More precisely, we compute 
\begin{align*} 
q_i^2 \T_i(\hat{S}(E_i^*)) = - \T_i(F_i) &= E_i K_i^{-1} 
= q_i^{2} K_i^{-1} E_i 
= -\hat{S}((K_i F_i)^*) = \hat{S}(\T_i(E_i)^*) 
\end{align*}
and 
\begin{align*} 
q^{a_{ij}} &\T_i(\hat{S}(E_j^*)) = - q_i^{a_{ij}} q_j^{-2} \T_i(F_j) \\
&= - q_i^{a_{ij}} q_j^{-2} \sum_{k = 0}^{-a_{ij}} (-1)^k q_i^{-k} F_i^{(-a_{ij} - k)} F_j F_i^{(k)} \\
&= \sum_{k = 0}^{-a_{ij}} (-1)^{-a_{ij} - k} q_i^{-a_{ij} - k} (-1)^{-a_{ij} + 1} q_i^{2 a_{ij}} q_j^{-2} F_i^{(-a_{ij} - k)} F_j F_i^{(k)} \\
&= \sum_{k = 0}^{-a_{ij}} (-1)^{-a_{ij} - k} q_i^{-a_{ij} - k} \hat{S}((E_i^{(-a_{ij} - k)} E_j E_i^{(k)})^*) \\
&= \sum_{k = 0}^{-a_{ij}} (-1)^k q_i^k \hat{S}((E_i^{(k)} E_j E_i^{(-a_{ij} - k)})^*) = \hat{S}(\T_i(E_j)^*) 
\end{align*}
for $ i \neq j $. Similarly, one obtains 
\begin{align*} 
q^{-a_{ij}} &\T_i(\hat{S}(F_j^*)) = \hat{S}(\T_i(F_j)^*) 
\end{align*}
for all $ i,j $. An induction argument then yields 
\begin{align*}
(\hat{S} \otimes \hat{S})((E_\alpha \otimes F_\alpha)^*) &= F_\alpha \otimes E_\alpha
\end{align*}
for all positive roots $ \alpha $. 
\end{proof}

Recall the definition of $ \Poly(G_q) $ from Section \ref{secoqg}. By construction, there exists a canonical 
bilinear pairing between $ U_q^\mathbb{R}(\mathfrak{k}) $ and $ \Poly(G_q) $. We shall introduce a $ * $-structure on $ \Poly(G_q) $ 
by stipulating 
$$
(x, f^*) = \overline{(\hat{S}^{-1}(x)^*, f)}
$$
for all $ x \in U_q^\mathbb{R}(\mathfrak{k}) $. In this way $ \Poly(G_q) $ becomes a Hopf $ * $-algebra.

\begin{definition}
	\label{def:CKq}
	We write $\CF^\infty(K_q)$ to denote $\Poly(G_q)$ with the above $*$-structure.
	\nomenclature[o$C^\infty(K_q)$]{$\CF^\infty(K_q)$}{algebra of polynomial functions on the compact semisimple quantum group 
	$K_q$}%
  \nomenclature{$K_q$}{quantization of a compact semisimple Lie group $K$}%
\end{definition}

The canonical pairing between $ U_q^\mathbb{R}(\mathfrak{k}) $ and $ \CF^\infty(K_q) $ satisfies 
\label{nom:K_q-pairing}%
\begin{align*}
(xy, f) &= (x, f_{(1)}) (y, f_{(2)}), \qquad (x, fg) = (x_{(2)}, f) (x_{(1)}, g) 
\end{align*}
and 
\begin{align*}
(\hat{S}(x), f) &= (x, S^{-1}(f)), \qquad (\hat{S}^{-1}(x), f) = (x, S(f)) 
\end{align*} 
for $ f, g \in \CF^\infty(K_q), x, y \in U_q^\mathbb{R}(\mathfrak{k}) $, and the compatibility with the $ * $-structures is given by 
\begin{align*}
(x, f^*) &= \overline{(\hat{S}^{-1}(x)^*, f)}, \qquad (x^*, f) = \overline{(x, S(f)^*)}. 
\end{align*} 
In addition, 
\begin{align*}
(x, 1) &= \hat{\epsilon}(x), \qquad (1, f) = \epsilon(f), 
\end{align*} 
for $ f \in \CF^\infty(K_q) $ and $ x \in U_q^\mathbb{R}(\mathfrak{k}) $. We may summarize this by saying that 
the canonical pairing is a skew-pairing of the Hopf $ * $-algebras $ U_q^\mathbb{R}(\mathfrak{k}) $ and $ \CF^\infty(K_q) $. 

We will also need to make use of the reverse skew-pairing between $ \CF^\infty(K_q) $ and $ U_q^\mathbb{R}(\mathfrak{k}) $ defined by 
$$
(f, x) = (\hat{S}(x), f) = (x, S^{-1}(f))
$$
\label{nom:K_q-reverse-pairing}%
for $ f \in \CF^\infty(K_q) $ and $ x \in U_q^\mathbb{R}(\mathfrak{k}) $. Then
\begin{align*}
(f^*, x) &= (\hat{S}(x), f^*) = \overline{(x^*, f)} = \overline{(\hat{S}(\hat{S}(x)^*), f)} = 
\overline{(f, \hat{S}(x)^*)}, 
\end{align*}
so that this pairing is again compatible with the $ * $-structures. As in Subsection \ref{sec:dual_multiplier_Hopf_algebra}, we stress that our skew-pairings $U_q^\mathbb{R}(\lie{k}) \times \CF^\infty(K_q) \to \mathbb{C}$ and $ \CF^\infty(K_q) \times U_q^\mathbb{R}(\lie{k}) \to \mathbb{C}$ need to be distinguished, and in particular $(f,x) \neq (x,f)$ in general.

Recall from Definition \ref{defogq} that the Hopf algebra $ \CF^\infty(K_q) $ is defined as the span of the matrix coefficients of the irreducible integrable representations $ V(\mu) $ of $ U_q^\mathbb{R}(\mathfrak{k}) $ for $ \mu \in \weights^+ $.  In particular, it is cosemisimple and so admits a left and right invariant integral $\phi$ given by projection onto the coefficient of the trivial corepresentation, see Proposition \ref{cosemisimpleintegral}.
In order to see that $ \CF^\infty(K_q) $ defines a compact quantum 
group $ K_q $ it remains to check that  $ \phi $ is positive. 

Following the discussion in Section 11.3 of \cite{KS}, this is equivalent to proving that 
the irreducible integrable representations $V(\mu)$ of $U_q^\mathbb{R}(\lie{k})$ are all unitarizable.
Here, a representation $V$ is called unitarizable if it admits a positive definite sesquilinear form which is invariant in the sense that
$$
 \bra X \cdot v, w \ket = \bra v, X^* \cdot w \ket
$$
for all $X\in U_q^\mathbb{R}(\lie{k})$ and all $v,w\in V$.

Let us explain how this can be done. 
First, we claim that each module $ V(\mu) $ for $ \mu \in \weights^+ $ 
can be equipped with an essentially unique invariant sesquilinear form $ \bra \;, \; \ket $. 
\label{nom:V-ip}%
Indeed, let $ \overline{V(\mu)^*} $ be the $ U_q(\mathfrak{g}) $-module defined on the conjugate vector space of the dual $ V(\mu)^* $ 
by setting $ (X \cdot \overline{f})(v) = \overline{f}(X^* \cdot v) $ for all $ v \in V(\mu) $. 
Then $ \overline{V(\mu)^*} $ is an irreducible highest weight module of highest weight $ \mu $, and hence $ \overline{V(\mu)^*} \cong V(\mu) $. 
In particular, there exists a unique hermitian sesquilinear form $ \bra \;, \; \ket $ on $ V(\mu) $ such that $ \bra v_\mu, v_\mu \ket = 1 $. 

\begin{prop} 
Let $q\neq 1$ be a strictly positive real number.
For each $ \mu \in \weights^+ $ the hermitian form $ \bra \;, \; \ket $ on $ V(\mu) $ constructed above is positive definite. 
\end{prop} 

\begin{proof}
 We shall use a continuity argument in the parameter $ q $. 
In order to emphasize the dependence on $ q $ we write $ V(\mu)_q $ for the irreducible highest weight module of highest weight $ \mu $ 
associated with $ q $. 
Then $ V(\mu)_q = V(\mu)_\A \otimes_\A \mathbb{C} $ where $ V(\mu)_\A $ is the integral form of $ V(\mu) $ and $ s \in \A $ acts on $ \mathbb{C} $ 
such that $ s^L = q $. Let us write $ \bra \;, \; \ket = \bra \;, \; \ket_q $ for the hermitian sesquilinear form on $ V(\mu)_q $ 
satisfying $ \bra v_\mu, v_\mu \ket_q = 1 $.

Choosing a free $ \A $-basis of $ V(\mu)_\A $ we see that the forms $ \bra \;, \; \ket_q $ can be viewed as a continuous family 
of hermitian forms depending on $ q \in (0, \infty) $ on a fixed vector space. 

Recall from Proposition \ref{quantumtoclassical} that the specialisation $ U_1(\mathfrak{g}) $ of $ U_q^\A(\mathfrak{g}) $ at $ 1 $ 
maps onto the classical universal enveloping algebra $ U(\mathfrak{g}) $ of $ \mathfrak{g} $ over $ \mathbb{C} $. 
This map is compatible with $ * $-structures if we consider the $ * $-structure on $ U(\mathfrak{g}) $ given by 
\begin{align*}
E_i^* = F_i, \qquad F_i^* = E_i, \qquad H_i^* = H_i
\end{align*}
for $ i = 1, \dots, N $, corresponding to the compact real form or $ \mathfrak{g} $. 

The representation $ V(\mu)_1 $ of $ U_q^\A(\mathfrak{g}) $ at $ q = 1 $ 
correspond to the irreducible highest weight representation of weight $ \mu $ of the classical Lie algebra $ \mathfrak{g} $. 
In particular, the sesquilinear form $ \bra \;, \; \ket_1 $ is positive definite. 
By continuity, 
we conclude that $ \bra \;, \; \ket_q $ is positive definite for all $ q \in (0,\infty) $. 
\end{proof}

Let us specify the element $F \in \DF(K_q)$ which appears in the Schur orthogonality relations, see Subsection \ref{sec:compact_quantum_groups}.

\begin{lemma}
	\label{lem:F_is_K2rho}
 We have $F=K_{-2\rho}$. 
\end{lemma}

\begin{proof}
  Let us fix an orthonormal weight basis $(e_i)$ for $V(\mu)$ and let $M = (M_{ij})$ denote the matrix of $K_{-2\rho}$ acting on this basis, so that $M_{ij} = \bra e_i , K_{-2\rho}\cdot e_j \ket = (K_{-2\rho}, u^\mu_{ij})$.  We will show that $M$ satisfies the defining properties of $F_\mu$ from Subsection \ref{sec:compact_quantum_groups}.  Clearly, $M$ is positive.  If $w_0\in W$ is the longest element of the Weyl group then $w_0\rho=-\rho$ and so by the Weyl group invariance of the set of weights of $V(\mu)$ we have 
  \begin{align*}
  \tr(M) &= \sum_{\nu \in \weights(V(\mu))} q^{(-2\rho,\nu)} 
  = \sum_{\nu \in \weights(V(\mu))} q^{(2\rho,w_0\nu)} 
  = \tr(M^{-1}),
  \end{align*}
  where the sum is over all weights of $V(\mu)$ counted with multiplicities.  Finally, using Lemma \ref{lem:S2} we have for any $X\in U_q^\mathbb{R}(\lie{k})$,
  \begin{align*}
   (X, S^{2}(u^\mu_{ij}))
     & = (\hat{S}^{-2}(X), u^\mu_{ij}) \\
     & = (K_{-2\rho} X K_{2\rho}, u^\mu_{ij}) \\ 
     & = \left(X, \sum_{k,l} (K_{-2\rho}, u^\mu_{ik}) u^\mu_{kl} (K_{2\rho}, u^\mu_{lj}) \right).
  \end{align*}
  It follows that $S^2(u^\mu) = M u^\mu M^{-1}$.
  This completes the proof.
\end{proof}

We shall write $ C(K_q) $ for the Hopf $ C^* $-algebra of functions on $ K_q $ obtained as the completion of $ \CF^\infty(K_q) $ 
\nomenclature{$C(K_q)$}{$C^*$-algebra of continuous functions on $K_q$}%
in the GNS-representation of the left Haar weight. Let us remark that $ K_q $ is coamenable, so that there is no need to 
distinguish between maximal and reduced algebras of functions in this case, see Corollary 5.1 in \cite{Banicafusion}.

To conclude this section, let us explain how the classical maximal torus $T$
\nomenclature{$T$}{maximal torus of $K$ or $K_q$}%
of $K$ appears also as a quantum subgroup of $K_q$.  Let $ U_q^\mathbb{R}(\mathfrak{t})$
\nomenclature{$U_q^\mathbb{R}(\lie{t})$}{Cartan subalgebra $U_q(\lie{h})$ with $*$-structure inherited from the real form $U_q^\mathrm{R}(\lie{k})$}%
be the Hopf algebra $U_q(\lie{h})$ equipped with the $*$-structure induced from $U_q^\mathbb{R}(\mathfrak{k}) $.  Define $\CF^\infty(T) \subset U_q^\mathbb{R}(t)^*$
\nomenclature[o$C^\infty(T)$]{$\CF^\infty(T)$}{algebra of polynomial functions on the maximal torus $T$ of $K$ or $K_q$}%
to be the restriction of elements of $\CF^\infty(K_q) \subset U_q^\RR(\lie{k})^*$ to $U_q^\mathbb{R}(\lie{t})$.  In light of definition \ref{defogq}, the image of this restriction map is the set of linear combination of characters of the form
\[
  \chi_\lambda ( K_\mu) = q^{(\mu,\lambda)}
\]
with $\lambda \in \weights$.  We therefore have a canonical isomorphism of $\CF^\infty(T)$ with the group algebra $\CC[\weights]$ which is isomorphic to the polynomials on the maximal torus $T$ of $K$. 

This homomorphism induces a surjection $ C(K_q) \rightarrow C(T) $ on the $ C^* $-level.

\subsection{Complex semisimple quantum groups} \label{seccqg}

In this section we define our main objects of study, namely $ q $-deformations of complex semisimple Lie groups. The construction of 
these quantum groups relies on the Drinfeld double construction. We also discuss some variants of the quantum groups we will be 
interested in, related to connected components and certain central extensions. Throughout this 
section we fix a positive deformation parameter $ q = e^h $ different from $ 1 $.

\subsubsection{The definition of complex quantum groups} 
\label{sec:CQG_definition}

The definition of complex semisimple quantum groups is a special case of the Drinfeld double construction explained in Section \ref{secaqg}. 

\begin{definition} \label{defcqg}
Let $ G $ be a simply connected semisimple complex Lie group and let $ K $ be a maximal compact subgroup. 
The complex semisimple quantum group $ G_q $ is the Drinfeld double 
$$
G_q = K_q \bowtie \hat{K}_q 
$$
\nomenclature{$G_q$}{$q$-deformation of a complex semisimple Lie group $G$}%
of $ K_q $ and its Pontrjagin dual $ \hat{K}_q $ with respect to the canonical bicharacter $ W \in \M(\CF^\infty(K_q) \otimes \DF(K_q)) $ 
given by the multiplicative unitary of $ K_q $. 
\end{definition}

Let us explicitly write down some of the structure maps underlying the quantum group $ G_q $, based on the general discussion in Section \ref{secaqg}.  
By definition, $ G_q $ is the algebraic quantum group with underlying algebra 
$$ 
\CF^\infty_c(G_q) = \CF^\infty_c(K_q \bowtie \hat{K}_q) = \CF^\infty(K_q) \otimes \DF(K_q),  
$$ 
\nomenclature[o$C^\infty(G_q)$]{$\CF^\infty_c(G_q)$}{algebra of functions $\CF^\infty_c(G_q) = \CF^\infty(K_q \bowtie \hat{K}_q)$ on the algebraic quantum group $G_q$}%
equipped with the comultiplication
$$
\Delta_{G_q} = (\id \otimes \sigma \otimes \id)(\id \otimes \ad(W) \otimes \id)(\Delta \otimes \hat{\Delta}). 
$$
\nomenclature{$\Delta_{G_q}$}{coproduct on $\CF_c^\infty(G_q)$}%
Note that we can identify
$$ 
\CF^\infty_c(G_q) = \bigoplus_{\mu \in \weights^+} \CF^\infty(K_q) \otimes \LH(V(\mu)),  
$$ 
so that the algebraic multiplier algebra of $ \CF^\infty_c(G_q) $ is 
$$
\CF^\infty(G_q) = \M(\CF^\infty_c(G_q)) = \prod_{\mu \in \weights^+} \CF^\infty(K_q) \otimes \LH(V(\mu)). 
$$
\nomenclature[o$C^\infty(G_q)$2]{$\CF^\infty(G_q)$}{multiplier algebra of $\CF^\infty_c(G_q)$}%
The counit of $ \CF_{c}^\infty(G_q) $ is the tensor product counit $ \epsilon_{G_q} = \epsilon \otimes \hat{\epsilon} $,
\nomenclature{$\epsilon_{G_q}$}{counit of $\CF^\infty_c(G_q)$}%
and the antipode is given by 
$$
S_{G_q}(f \otimes x) = W^{-1} (S(f) \otimes \hat{S}(x)) W = (S \otimes \hat{S})(W (f \otimes x) W^{-1}).
$$
\nomenclature[o$S$3]{$S_{G_q}$}{antipode of $\CF^\infty_c(G_q)$}%

Like its classical counterpart, the quantum group $ G_q $ is unimodular \cite{PWlorentz}. 
\begin{prop} \label{doubleunimodular}
A positive left and right invariant integral on $ \CF^\infty_c(G_q) $ is given by 
$$
\phi_{G_q}(f \otimes x) = \phi(f) \hat{\psi}(x), 
$$
for $f\in\CF^\infty(K_q)$, $x\in\DF(K_q)$.
\nomenclature{$\phi_{G_q}$}{left and right invariant integral $\phi_{G_q} = \phi\otimes\hat{\psi}$ on $\CF^\infty_c(G_q)$}%
\end{prop} 

\begin{proof}
It follows from the end of Subsection \ref{sec:compact_quantum_groups} that the right invariant Haar integral $\hat\psi$ on $\DF(K_q)$ is given by
$$ 
 \hat{\psi}(x) = \hat{\phi}(F^2 x)
$$ 
for all $ x \in \DF(K_q) $, where $ F = K_{-2 \rho} $ is the invertible element obtained from the Schur orthogonality relations, see Lemma \ref{lem:F_is_K2rho}. 
Moreover, observe that
\begin{align*}
\sum_{\mu, i,j} (F^{-1} \hit u^\mu_{ij} \hitby F^{-1}) \otimes F \omega^\mu_{ij} F = \sum_{\mu, i,j} u^\mu_{ij} \otimes \omega^\mu_{ij} = W 
\end{align*}
due to the dual basis property of $ W $. 
This implies
\begin{align*}
(\phi \otimes \hat{\phi})(W(f \otimes x)) 
&= \sum_{\mu, i,j} \phi((F^{-1} \hit u^\mu_{ij} \hitby F^{-1}) f) \hat{\phi}(F \omega^\mu_{ij} F x) \\
&= \sum_{\mu, i,j} \phi(f u^\mu_{ij}) \hat{\phi}(F^2 x \omega^\mu_{ij}) \\
&= \sum_{\mu, i,j} \phi(f u^\mu_{ij}) \hat{\psi}(x \omega^\mu_{ij}) \\
&= (\phi \otimes \hat{\psi})((f \otimes x)W) 
\end{align*}
for all $ f \in \CF^\infty(K_q), x \in \DF(K_q) $. 

Since $ (\id \otimes \hat{\Delta})(W) = W_{13} W_{12} $ we thus obtain 
\begin{align*}
&(\id \otimes \phi_G)\Delta_G(f \otimes x)
= (\id \otimes \phi_G)(\id \otimes \sigma \otimes \id)(\id \otimes \ad(W) \otimes \id)(\Delta \otimes \hat{\Delta})(f \otimes x)  \\ 
&= (\id \otimes \phi \otimes \id \otimes \hat{\psi}) \ad(W_{23}) (\Delta \otimes \hat{\Delta})(f \otimes x) \\ 
&= (\id \otimes \phi \otimes \id \otimes \hat{\psi}) \ad(W_{24}^{-1} (\id \otimes \hat{\Delta})(W)_{234}) (\Delta \otimes \hat{\Delta})(f \otimes x) \\ 
&= (\id \otimes \phi \otimes \id \otimes \hat{\phi}) \ad((\id \otimes \hat{\Delta})(W)_{234}) (\Delta \otimes \hat{\Delta})(f \otimes x) \\ 
&= (\id \otimes \phi \otimes \hat{\phi}) \ad(W_{23}) (\Delta \otimes \id)(f \otimes x) \\ 
&= (\id \otimes \phi \otimes \hat{\psi})(\Delta \otimes \id)(f \otimes x) \\ 
&= (\phi \otimes \hat{\psi})(f \otimes x) = \phi_G(f \otimes x)
\end{align*}
for all $ f \in \CF^\infty_c(K_q), x \in \DF(K_q) $ as desired. Hence $ \phi_G $ is left invariant. 
Similarly, since $ (\Delta \otimes \id)(W) = W_{13} W_{23} $ we obtain 
\begin{align*}
&(\phi_G \otimes \id)\Delta_G(f \otimes x)
= (\phi_G \otimes \id)(\id \otimes \sigma \otimes \id)(\id \otimes \ad(W) \otimes \id)(\Delta \otimes \hat{\Delta})(f \otimes x)  \\ 
&= (\phi \otimes \id \otimes \hat{\psi} \otimes \id) \ad(W_{23}) (\Delta \otimes \hat{\Delta})(f \otimes x) \\ 
&= (\phi \otimes \id \otimes \hat{\psi} \otimes \id) \ad(W_{13}^{-1} (\Delta \otimes \id)(W)_{123}) (\Delta \otimes \hat{\Delta})(f \otimes x) \\ 
&= (\phi \otimes \id \otimes \hat{\phi} \otimes \id) \ad((\Delta \otimes \id)(W)_{123}) (\Delta \otimes \hat{\Delta})(f \otimes x) \\ 
&= (\phi \otimes \hat{\phi} \otimes \id) \ad(W_{12}) (\id \otimes \hat{\Delta})(f \otimes x) \\ 
&= (\phi \otimes \hat{\psi} \otimes \id)(\id \otimes \hat{\Delta})(f \otimes x) \\ 
&= (\phi \otimes \hat{\psi})(f \otimes x) = \phi_G(f \otimes x), 
\end{align*}
which means that $ \phi_G $ is also right invariant. 
\end{proof}

The dual $ \DF(G_q) $ of $ \CF^\infty_c(G_q) $ in the sense of algebraic quantum groups is given by 
$$ 
\DF(G_q) = \DF(K_q) \bowtie \CF^\infty(K_q), 
$$ 
\nomenclature[o$D(G_q)$]{$\DF(G_q)$}{group algebra $\DF(G_q) = \DF(K_q) \bowtie \CF^\infty(K_q)$ of $G_q$}%
equipped with the tensor product comultiplication and the multiplication given by 
\begin{align*}
(x \bowtie f)(y \bowtie g) &= x (y_{(1)}, f_{(1)}) y_{(2)} \bowtie f_{(2)} (\hat{S}(y_{(3)}), f_{(3)}) g \\
&= x (y_{(1)}, f_{(1)}) y_{(2)} \bowtie f_{(2)} (y_{(3)}, S^{-1}(f_{(3)})) g, 
\end{align*}
using the natural skew-pairing between $ \DF(K_q) $ and $ \CF^\infty(K_q) $. 

It is important to note that when relating $ \CF^\infty_c(G_q) $ with $ \DF(G_q) $ one should work with the pairing $\DF(G_q) \times \CF^\infty_c(G_q) \to \mathbb{C}$ given by
$$
 (y\bowtie g, f\otimes x) = (y,f)(g,x),
$$
\label{nom:G_q-pairing}%
for $ f,g \in \CF^\infty(K_q), x, y \in \DF(K_q) $, where we recall that we write $(g,x) = (\hat{S}(x),g) = (x,S^{-1}(g))$ for the reverse pairing of $ g \in \CF^\infty(K_q) $ and $ x \in \DF(K_q) $. 
The appearance of the antipode in this context is due to the fact that 
the comultiplication of $ \DF(G_q) $ is supposed to be the transpose of the opposite multiplication on $ \CF^\infty_c(G_q) $, so that 
one has to flip the roles of $ \CF^\infty(K_q) $ and $ \DF(K_q) $ in the second tensor factor of the pairing. 

The antipode of $ \DF(G_q) $ is defined by 
\begin{align*}
 \hat{S}_{G_q}(x \bowtie f) &= (1 \bowtie S(f))(\hat{S}(x) \bowtie 1) \\
&= (\hat{S}(x_{(3)}), S(f_{(3)})) \hat{S}(x_{(2)}) \bowtie S(f_{(2)}) (\hat{S}(x_{(1)}), f_{(1)}), \\ 
&= (x_{(3)}, f_{(3)}) \hat{S}(x_{(2)}) \bowtie S(f_{(2)}) (\hat{S}(x_{(1)}), f_{(1)}), 
\end{align*}
\nomenclature[o$S$4]{$\hat{S}_{G_q}$}{antipode of $\DF(G_q)$}%
and the $ * $-structure on $ \DF(G_q) $ is given by 
$$
(x \bowtie f)^* = (1 \bowtie f^*)(x^* \bowtie 1) 
= (x^*_{(1)}, f^*_{(1)}) x^*_{(2)} \bowtie f^*_{(2)} (\hat{S}(x^*_{(3)}), f^*_{(3)}).  
$$
A left Haar integral for $ \DF(G_q) $ is given by 
$$ 
\hat{\phi}_{G_q}(x \bowtie f) = \hat{\phi}(x) \phi(f)
$$ 
\nomenclature[o$\phi_{G_q}$2]{$\hat{\phi}_{G_q}$}{left Haar integral $\hat{\phi}_{G_q} = \hat{\phi}\otimes\phi$ on $\DF(G_q)$}%
\nomenclature[o$\psi_{G_q}$2]{$\hat{\psi}_{G_q}$}{right Haar integral $\hat{\psi}_{G_q} = \hat{\psi}\otimes\psi$ on $\DF(G_q)$}%
where $ \hat{\phi} $ and $ \phi $ are left Haar integrals for $ \DF(K_q) $ and $ \CF^\infty(K_q) $, respectively. 
Similarly, a right Haar integral for $ \DF(G_q) $ is given by the tensor product of right Haar integrals for $ \DF(K_q) $ and $ \CF^\infty(K_q) $. 

Since $ G_q $ is an algebraic quantum group, the general theory outlined in Section \ref{secaqg} implies that it defines 
a locally compact quantum group. Explicitly, the Hopf $ C^* $-algebra of functions on $ G_q $ is given by
$$ 
C_0(G_q) = C(K_q) \otimes C^*(K_q).  
$$ 
\nomenclature[o$CG$16]{$C_0(G_q)$}{$C^*$-algebra of continuous functions vanishing at infinity on $G_q$}%
Note that since the compact quantum group $ K_q $ is both amenable and coamenable there is no need to distinguish between maximal and 
reduced versions of $ C(K_q) $ and $ C^*(K_q) $ here. 
The quantum group $ G_q $ is coamenable as well because $ \epsilon \otimes \hat{\epsilon} $ is a bounded counit for $ C_0(G_q) = C_0^{\red}(G_q) $, 
so that $ C_0^{\mx}(G_q) \cong C_0^{\red}(G_q) $, compare Theorem 3.1 in \cite{BedosTuset}. 

The full and reduced Hopf $ C^* $-algebras $ C^*_\mx(G_q) $ and $ C^*_\red(G_q) $
may be constructed using the general theory of locally compact 
quantum groups, or more directly as completions of $ \DF(G_q) $. 
We remark that the maximal and reduced group $ C^* $-algebras of $ G_q $ are not isomorphic, that is, the quantum group $ G_q $ is not amenable. 

It is not apparent from Definition \ref{defcqg} why $ G_q $ should be understood as a quantum deformation of the complex Lie group $ G $. Nevertheless, this is indeed 
the case, and it is a basic instance of the quantum duality principle, see \cite{Drinfeldicm}. 
From this point of view, the discrete part $ \hat{K}_q $ of the Drinfeld double corresponds to the subgroup $ AN $ in the Iwasawa decomposition 
$ G = KAN $ of $ G $. In fact, the group $ G $ can be viewed as  the classical double of the Poisson-Lie group $ K $ with 
its standard Poisson structure, see \cite{KorogodskiSoibelmanbook}. 
As we shall discuss further below, the deformation aspect of the theory of complex quantum groups is most visible 
if one works with algebras of polynomial functions instead. 
This is in fact the starting point taken by Podle\'s and Woronowicz in \cite{PWlorentz}.

\subsubsection{The connected component of the identity} 
\label{sec:connected_component}

	Although the classical group $G$ is connected, the quantization $G_q$ defined above behaves more like an almost connected group.  Specifically, we shall see here that each of the quantum groups $G_q$ admits a quotient to a finite classical group.
	This has implications in the representation theory of $G_q$.  In particular, it leads to a finite group of one-dimensional representations.

Recall that we are using $G$ to denote the connected, simply connected Lie group associated to the complex semisimple Lie algebra $\lie{g}$.  The centre $Z$ of $G$
\nomenclature{$Z$}{centre of $G$}%
is a finite subgroup of the maximal compact torus $T$ in $K$.  As such, we can identify $Z$ also as a quantum subgroup of the compact quantum group $K_q$, via the surjective morphism of Hopf $*$-algebras $ \pi_Z: \CF^\infty(K_q) \rightarrow \CF^\infty(T) \rightarrow \CF^\infty(Z) $.
\nomenclature{$\pi_Z$}{projection $\CF^\infty(K_q) \to \CF^\infty(Z)$ corresponding to the central subgroup $Z$ of $K_q$}%

There is a canonical isomorphism between $Z$ and the quotient group $\weights^\vee/\roots^\vee$, where
\[
 \weights^\vee = \{ \mu \in \lie{h}^* \mid (\mu,\alpha) \in \mathbb{Z} \text{ for all } \alpha \in \roots \}
\]
\nomenclature[o$P^\vee$]{$\weights^\vee$}{coweight lattice}%
denotes the coweight lattice associated to the root system of $\lie{g}$.  We can explicitly realize $\weights^\vee/\roots^\vee$ as a central quantum subgroup of $K_q$ as follows.  Consider the finite subgroup $(i\hbar^{-1}\weights^\vee)/(i\hbar^{-1}\roots^\vee)$ of $\lie{h}_q^*$.  
The corresponding algebra characters $K_{i\hbar^{-1}\gamma}$
\label{nom:K-central}%
with $\gamma\in\weights^\vee/\roots^\vee$
have the property that if $f = \bra v'|\bullet|v\ket \in \End(V(\mu))^*$ is a matrix coefficient of any simple module $V(\mu)$, with $v$ a vector of weight $\nu$ then
\[
 (K_{i\hbar^{-1}\gamma}, f) = q^{i\hbar^{-1}(\gamma,\nu)} \bra v' , v \ket
   = e^{2\pi i(\gamma,\nu)} \epsilon(f) = e^{2\pi i(\gamma,\mu)} \epsilon(f) ,
\]
since $\mu-\nu \in \roots$.  That is to say, $K_{i\hbar^{-1}\gamma}$ evaluates on any matrix coefficient in $\End(V(\mu))^*$ as a constant multiple $e^{2\pi i(\gamma,\mu)}$ of the counit.
It follows that
\[
 K_{i\hbar^{-1}\gamma}\hit f  = f_{(1)} (K_{i\hbar^{-1}\gamma}, f_{(2)}) = (K_{i\hbar^{-1}\gamma}, f_{(1)}) f_{(2)}= f \hitby K_{i\hbar^{-1}\gamma}
\]
for all $f\in \CF^\infty(K_q)$, and thus that $K_{i\hbar^{-1}\gamma}$ is a central group-like element in $\M(\DF(K_q))$.  The subalgebra of $\M(\DF(K_q))$ spanned by these elements is a Hopf $*$-algebra isomorphic to the group algebra $\DF(\weights^\vee/\roots^\vee) \cong \DF(Z)$.
The above calculation also implies that the projection $\pi_Z:\CF^\infty(K_q) \to \CF^\infty(Z)$ satisfies
\[
 \pi_Z(f_{(1)}) \otimes f_{(2)} = \pi_Z(f_{(2)}) \otimes f_{(1)}
\]
for all $f\in \CF^\infty(K_q)$.

\begin{prop}
	\label{prop:Gq_to_Zhat}
	The linear map $\DF(G_q) \to \DF(\hat{Z})$ defined by
	\[
	 \hat\epsilon \otimes \pi_Z : \DF(G_q) = \DF(K_q) \bowtie \CF^\infty(K_q) \to \CF^\infty(Z) \cong \DF(\hat{Z})
	\]
	is a surjective morphism of algebraic quantum groups.
\end{prop}

\begin{proof}
	Let $x\bowtie f, y\bowtie g \in \DF(K_q) \bowtie \CF^\infty(G_q)$. 	Using the above property of $\pi_Z$, we have
	\begin{align*}
	 (\hat{\epsilon}\otimes\pi_Z) ((x\bowtie f)(y\bowtie g))
	  & = (y_{(1)},f_{(1)}) \, (\hat{S}(y_{(3)}),f_{(3)}) \, 
	  \hat{\epsilon}(xy_{(2)}) \pi_Z(f_{(2)}g) \\
	  &= (y_{(1)},f_{(1)}) \, (\hat{S}(y_{(3)}),f_{(2)}) \, 
	  \hat{\epsilon}(x)\hat{\epsilon}(y_{(2)}) \pi_Z(f_{(3)})\pi_Z(g) \\
	  &= (y_{(1)}\hat{S}(y_{(2)}),f_{(1)}) \hat{\epsilon}(x) \pi_Z(f_{(2)}) \pi_Z(g) \\
	  &= \hat{\epsilon}(y)\hat{\epsilon}(x) \pi_Z(f) \pi_Z(g)\\
	  &= (\hat{\epsilon}\otimes\pi_Z) (x\bowtie f) \, (\hat{\epsilon}\otimes\pi_Z) (y\bowtie g).	  
	\end{align*}
	The map $\hat{\epsilon}\otimes\pi_Z$ clearly respects the coproduct, and is therefore a morphism of multiplier Hopf-algebras.  It is surjective since the quotient map $\CF^\infty(K_q) \to \CF(Z)$ is surjective, and it respects the involution since 
	\[
	(\hat{\epsilon}\otimes\pi_Z)((x\bowtie f)^*) = ((\hat{\epsilon}\otimes\pi_Z)(1\bowtie f^*) )
	((\hat{\epsilon}\otimes\pi_Z)(x^*\bowtie 1) )= \overline{\hat{\epsilon}(x)} \pi_Z(f)^*
	\]
	for all $x\bowtie f \in \DF(G_q)$.
\end{proof}

Since the elements $K_{i\hbar^{-1}\gamma}$ for $\gamma\in\weights^\vee/\roots^\vee$ are group-like in $\DF(Z) = \CF^\infty(\hat{Z})$, we can use the morphism of Proposition \ref{prop:Gq_to_Zhat} to define a finite group of algebra characters of $\DF(G_q)$, indexed by $Z \cong \weights^\vee/\roots^\vee$, as follows.

\begin{definition}
	\label{def:Gq-characters}
	For each $\gamma\in \weights^\vee/\roots^\vee$, we define a character $\hat{\epsilon}_\gamma:\DF(G_q) \to \mathbb{C}$ by
	\[
	  \hat{\epsilon}_\gamma (x\bowtie f) = (K_{i\hbar^{-1}\gamma}, (\hat{\epsilon}\otimes \pi_Z) (x\bowtie f))
	   = \hat{\epsilon}(x) (K_{i\hbar^{-1}\gamma},f).
	\]
	\nomenclature[o$\epsilon_\gamma$]{$\hat{\epsilon}_\gamma$}{algebra character of $\DF(G_q)$}%
\end{definition}

These are non-degenerate $*$-characters, so they correspond to one-dimensional unitary representations of the quantum group $G_q$.  Note that we have $(\hat{\epsilon}_\gamma \otimes \hat{\epsilon}_{\gamma'})\circ\hat{\Delta} = \hat{\epsilon}_{\gamma\gamma'}$ for all $\gamma, \gamma'\in  \weights^\vee/\roots^\vee$.

We shall denote by $ PK_q $ the quotient quantum group $ PK_q = K_q/Z $.  Explicitly, $PK_q$ is defined via its algebra of functions
\begin{align*}
\CF^\infty(PK_q) & = \{f \in \CF^\infty(K_q) \mid (\id\otimes\pi_Z)\Delta(f) = f \otimes 1 \} \\
&= \bigoplus_{\mu\in\weights^+\cap\roots} \End(V(\mu))^*.
\end{align*}
\nomenclature{$PK_q$}{quantum group $K_q/Z$, corresponding to compact group $PK$ of adjoint type}%
Geometrically, this corresponds to the 
projective version of $ K_q $, which is a quantization of the group of adjoint type associated with $ \mathfrak{k} $.  
The dual multiplier Hopf algebra is
\[
\DF(PK_q) = \bigoplus_{\mu\in\weights^+\cap\roots} \End(V(\mu)),
\]
and the obvious projection
$\hat{\pi} : \DF(K_q) \to \DF(PK_q)$
defines a morphism of algebraic quantum groups.

Let us define the ``connected component'' $ G_q^0 $ of $ G_q $ to be the quantum double $ G_q^0 = K_q \bowtie \widehat{PK}_q $. That is, 
$$
\CF^\infty_c(G_q^0) = \CF^\infty(K_q) \otimes \DF(PK_q)
$$
\nomenclature{$G_q^0$}{``connected component'' of $G_q$, defined by $ G_q^0 = K_q \bowtie \widehat{PK}_q $}%
with the coproduct twisted by the bicharacter $ U = (\id \otimes \hat{\pi})(W_{K_q})$, where $ W_{K_q} \in \M(\CF^\infty(K_q) \otimes \DF(K_q)) $ 
is the fundamental multiplicative unitary for $ K_q $ and $ \hat{\pi}: \DF(K_q) \rightarrow \M(\DF(PK_q)) $ is the projection homomorphism.  The dual algebraic quantum group $\DF(G_q^0) = \DF(K_q) \bowtie \CF^\infty(PK_q)$ satisfies
\[
 \DF(G_q^0) = \{ u\in\DF(G_q) \mid \hat\epsilon_\gamma(u) = \hat\epsilon_{G_q}(u) \text{ for all } \gamma \in \weights^\vee/\roots^\vee \},
\]
which is to say that the one-dimensional representations of $G_q$ from Definition \ref{def:Gq-characters} all become trivial upon restriction to the quantum subgroup $G_q^0$.

\subsection{Polynomial functions} 

Our next aim is to describe the Hopf $ * $-algebra of holomorphic and antiholomorphic polynomial functions on the complex quantum group $ G_q $. 
We keep our general assumptions from the previous section. 

Using the universal $ R $-matrix of $ U_q(\mathfrak{g}) $ we obtain a skew-pairing $ r^{-1} $ between $ \Poly(G_q) $ and $ \Poly(G_q) $ given by 
$$
r^{-1}(f,g) = (\R^{-1}, f \otimes g), 
$$
see Section \ref{subsecuniversalrmatrix}. We also set $ r(f, g) = (\R, f \otimes g) = r^{-1}(S(f), g) $. The functional $ r $ 
will be referred to as the universal $ r $-form on $ \Poly(G_q) $. 

Starting from the universal $ r $-form, we may form the Hopf algebraic quantum double of $ \Poly(G_q) $. More precisely, we shall consider the Hopf algebra 
$$ 
\Poly^\mathbb{R}(G_q) = \Poly(G_q) \bowtie \Poly(G_q),  
$$ 
\nomenclature[o$O^R(G_q)$]{$\Poly^\mathbb{R}(G_q)$}{polynomial algebra of $G_q$, $\Poly^\mathbb{R}(G_q)= \Poly(G_q)\bowtie\Poly(G_q)$}%
with underlying vector space $ \Poly(G_q) \otimes \Poly(G_q) $ together with the tensor product coalgebra structure and the multiplication 
$$
(f \bowtie g)(h \bowtie k) = (\R^{-1}, h_{(1)} \otimes g_{(1)}) fh_{(2)} \bowtie g_{(2)} k (\R, h_{(3)} \otimes g_{(3)}). 
$$

\begin{lemma} 
The algebra $ \Poly^\mathbb{R}(G_q) = \Poly(G_q) \bowtie \Poly(G_q) $ becomes a Hopf $ * $-algebra with the $ * $-structure 
$$
(f \bowtie g)^* = g^* \bowtie f^*. 
$$
\end{lemma}

\begin{proof}
Due to the relation $ (\hat{S} \otimes \hat{S})(\R^*) = \R_{21} $ obtained in Lemma \ref{rmatrixreal} 
the universal $ r $-form on $ \Poly(G_q) = \Poly^\mathbb{R}(K_q)$ is real, that is, 
$$ 
(\R, f^* \otimes g^*) = \overline{(\R, g \otimes f)}. 
$$ 
We compute 
\begin{align*}
((f \bowtie g)(h \bowtie k))^* &= \overline{(\R^{-1}, h_{(1)} \otimes g_{(1)})} (k^* g_{(2)}^* \bowtie h_{(2)}^* f^*) \overline{(\R, h_{(3)} \otimes g_{(3)})} \\ 
&= (\R^{-1}, g^*_{(1)} \otimes h^*_{(1)}) (k^* g_{(2)}^* \bowtie h^*_{(2)} f^*) (\R, g^*_{(3)} \otimes h^*_{(3)}) \\ 
&= (k^* \bowtie h^*)(g^* \bowtie f^*) \\ 
&= (h \bowtie k)^*(f \bowtie g)^* 
\end{align*} 
using 
\begin{align*}
(\R^{-1}, g^* \otimes h^*) &= (\R, S^{-1}(g^*) \otimes h^*) \\
&= (\R, S(g)^* \otimes h^*) \\
&= \overline{(\R, h \otimes S(g))} = \overline{(\R^{-1}, h \otimes g)}. 
\end{align*}
Hence $ \Poly^\mathbb{R}(G_q) $ is indeed a $ * $-algebra with the above $ * $-structure. 
It is clear that the comultiplication on $ \Poly^\mathbb{R}(G_q) $ is a $ * $-homomorphism. 
\end{proof} 

The Hopf $ * $-algebra $ \Poly^\mathbb{R}(G_q) $ should be thought of as a quantum version of the algebra of polynomial functions on the real Lie group 
underlying $ G $. The two copies of $ \Poly(G_q) $ inside $ \Poly^\mathbb{R}(G_q) $ may be interpreted as holomorphic 
and antiholomorphic polynomials. 

We recall from Lemma \ref{lfunctionalproperties} that the $ l $-functionals on $ \Poly(G_q) $, given by 
$$
l^+(f)(h) = (\R, h \otimes f), \qquad l^-(f)(h) = (\R^{-1}, f \otimes h)
$$
for $ f \in \Poly(G_q) $, 
satisfy 
$
\hat{\Delta}(l^\pm(f)) = l^\pm(f_{(1)}) \otimes l^\pm(f_{(2)}) 
$
and 
$
l^\pm(fg) = l^\pm(f) l^\pm(g) 
$
for $ f, g \in \Poly(K_q) $. 

\begin{lemma} \label{lfunctionalstar}
For $ f \in \CF^\infty(K_q) $ we have 
$$
l^\pm(f)^* = l^\mp(f^*)
$$
in $ U_q^\mathbb{R}(\mathfrak{k}) $. 
\end{lemma} 

\begin{proof} We compute 
\begin{align*}
(l^+(f)^*, h) &= \overline{(l^+(f), S(h)^*)} \\ 
&= \overline{(\R, S(h)^* \otimes f)} \\ 
&= (\R, f^* \otimes S(h)) \\ 
&= (\R^{-1}, f^* \otimes h) \\ 
&= (l^-(f^*), h)   .
\end{align*}
Similarly, 
\begin{align*}
(l^-(f)^*, h) &= \overline{(l^-(f), S(h)^*)} \\ 
&= \overline{(\R^{-1}, f \otimes S(h)^*)} \\ 
&= \overline{(\R, f \otimes h^*)} \\ 
&= (\R, h \otimes f^*) \\ 
&= (l^+(f^*), h)   
\end{align*}
as claimed. 
\end{proof}

The following result should be viewed as a quantum analogue of the fact that complex-valued polynomial functions 
on the group $ G $ embed into the algebra of smooth functions on $ G $, compare \cite{Hodgesdouble}.

\begin{prop} \label{polynomialtosmooth}
The linear map $ i: \Poly^\mathbb{R}(G_q) \rightarrow 
\CF^\infty(G_q) = \M(\CF^\infty_c(G_q)) $ given by 
$$
i(f \bowtie g) = f_{(1)} g_{(1)} \otimes l^-(f_{(2)}) l^+(g_{(2)}) 
$$
\nomenclature{$i$}{homomorphism of multiplier Hopf $*$-algebras $i:\Poly^\mathbb{R}(G_q) \to \CF^\infty(G_q)$}%
is a nondegenerate homomorphism of multiplier Hopf $ * $-algebras.
\end{prop}

\begin{proof} 
Let $f,g,h,k \in \Poly(G_q)$.
Recall from Lemma \ref{skewpairingduality} that we have the relation $(\R,f_{(1)}\otimes g_{(1)}) g_{(2)}f_{(2)}
  = f_{(1)}g_{(1)} (\R, f_{(2)}\otimes g_{(2)})$ and hence
\[
  (\R^{-1},f_{(1)}\otimes g_{(1)}) f_{(2)}g_{(2)}
   = g_{(1)} f_{(1)} (\R^{-1}, f_{(2)}\otimes g_{(2)}).
\]
Moreover, from the calculation
\begin{align*}
  (\R, f_{(1)}\otimes g_{(1)}) \, & (l^+(g_{(2)}) l^-(f_{(2)}) , k) \\
    &= (\R, f_{(1)}\otimes g_{(1)}) \,(\R, k_{(1)}\otimes g_{(2)})
    \, (\R^{-1}, f_{(2)} \otimes k_{(2)} ) \\
    &= (\R, k_{(1)}f_{(1)}\otimes g) \, (\R^{-1}, f_{(2)} \otimes k_{(2)} ) \\
    &= (\R^{-1},f_{(1)}\otimes k_{(1)}) \, (\R, f_{(2)}k_{(2)} \otimes g) \\
    &= (\R^{-1},f_{(1)}\otimes k_{(1)}) \, (\R, f_{(2)} \otimes g_{(2)}) \, (\R, k_{(2)}\otimes g_{(1)}) \\
    &= ( l^-(f_{(1)}) l^+(g_{(1)}), k) \, (\R, f_{(2)} \otimes g_{(2)})
\end{align*}
we obtain the exchange relation 
\[
(\R, f_{(1)}\otimes g_{(1)}) l^+(g_{(2)}) l^-(f_{(2)})
= l^-(f_{(1)})l^+(g_{(1)}) (\R, f_{(2)} \otimes g_{(2)}).
\]
With these formulas, we compute
\begin{align*} 
i&((f \bowtie g)(h \bowtie k)) \\
 &= (\R^{-1}, h_{(1)} \otimes g_{(1)})\, i( fh_{(2)} \bowtie g_{(2)} k ) \, (\R, h_{(3)} \otimes g_{(3)}) \\
 &= (\R^{-1}, h_{(1)} \otimes g_{(1)})\, f_{(1)} h_{(2)} g_{(2)} k_{(1)} \otimes l^-(f_{(2)}) l^-(h_{(3)}) l^+(g_{(3)}) l^+( k_{(2)}) \, (\R, h_{(4)} \otimes g_{(4)}) \\
 &= f_{(1)} g_{(1)} h_{(1)} k_{(1)}
 \, (\R^{-1}, h_{(2)} \otimes g_{(2)}) \otimes 
  (\R, h_{(3)} \otimes g_{(3)}) \, l^-(f_{(2)}) l^+(g_{(4)})  l^-(h_{(4)}) l^+(k_{(2)})  \\
 &= f_{(1)} g_{(1)} h_{(1)} k_{(1)} \otimes l^-(f_{(2)}) l^+(g_{(2)}) l^-(h_{(2)}) l^+(k_{(2)}) \\
 &= i(f \bowtie g) i(h \bowtie k).
\end{align*}
We conclude that $ i $ is an algebra homomorphism. 

To check the coalgebra homomorphism property it is enough to consider 
elements of the form $ f \bowtie 1 $ and $ 1 \bowtie g $. Using Lemma \ref{lfunctionalproperties} we compute 
\begin{align*} 
\Delta_{G_q} i(f \bowtie 1) &= \Delta_{G_q}(f_{(1)} \otimes l^-(f_{(2)})) \\
&  = (\id \otimes \sigma \otimes \id) \,
  (f_{(1)} \otimes W(f_{(2)} \otimes l^-(f_{(3)}))W^{-1} \otimes l^-(f_{(4)}) \\
&= (f_{(1)} \otimes l^-(f_{(2)})) \otimes (f_{(3)} \otimes l^-(f_{(4)})) \\
&= (i \otimes i)(f_{(1)} \bowtie 1 \otimes f_{(2)} \bowtie 1) \\
&= (i \otimes i)\Delta(f \bowtie 1) 
\end{align*}
for $ f \in \Poly(G_q) $, using 
\begin{align*}
(W (f_{(1)} \otimes l^-(f_{(2)})), x \otimes h) &= (W, x_{(1)} \otimes h_{(1)}) (f_{(1)} \otimes l^-(f_{(2)}), x_{(2)} \otimes h_{(2)}) \\
&= (h_{(1)}, x_{(1)}) (f_{(1)}, x_{(2)}) (l^-(f_{(2)}), h_{(2)}) \\
&= (h_{(1)} f_{(1)}, x) (\R^{-1}, f_{(2)} \otimes h_{(2)}) \\
&= (\R^{-1}, f_{(1)} \otimes h_{(1)}) (f_{(2)} h_{(2)}, x) \\
&= (l^-(f_{(1)}), h_{(1)}) (f_{(2)}, x_{(1)}) (h_{(2)}, x_{(2)}) \\
&= (f_{(2)} \otimes l^-(f_{(1)}), x_{(1)} \otimes h_{(1)}) (W, x_{(2)} \otimes h_{(2)}) \\
&= ((f_{(2)} \otimes l^-(f_{(1)})) W, x \otimes h). 
\end{align*}
In the same way one checks $ \Delta_{G_q} i(1 \bowtie g) = (i \otimes i)\Delta(1 \bowtie g) $ for $ g \in \Poly(G_q) $. 

For the $ * $-compatibility notice that we have 
\begin{align*}
i((f \otimes g)^*) &= i(g^* \otimes f^*) \\ 
&= g_{(1)}^* f_{(1)}^* \otimes l^-(g^*_{(2)}) l^+(f^*_{(2)}) \\ 
&= (f_{(1)} g_{(1)})^* \otimes l^+(g_{(2)})^* l^-(f_{(2)})^* \\ 
&= i(f \otimes g)^*, 
\end{align*}
due to Lemma \ref{lfunctionalstar}. 
\end{proof}

\subsection{The quantized universal enveloping algebra of a complex group} \label{subsecuqrg}

In this section we introduce the quantized universal enveloping algebra of the complex quantum group $ G_q $ and discuss some related constructions and 
results. We keep our general assumptions from the previous section. 

Recall that the group algebra $ \DF(G_q) $ can be identified with the Drinfeld double 
$$
\DF(G_q) = \DF(K_q) \bowtie \CF^\infty(K_q) 
$$
with respect to the canonical pairing between $ \DF(K_q) $ and $ \CF^\infty(K_q) $. 

\begin{prop} \label{finiteseparating}
The linear map $ \iota: \DF(G_q) \rightarrow \M(\DF(K_q) \otimes \DF(K_q)) $ given by 
$$ 
\iota(x \bowtie f) = \hat{\Delta}(x) (l^-(f_{(1)}) \otimes l^+(f_{(2)})) 
$$ 
\label{nom:iota-DGq}%
is an injective essential algebra homomorphism. 
\end{prop}

\begin{proof} 
	We compute 
\begin{align*}
\iota&((x \bowtie f)(y \bowtie g)) = \iota(x (y_{(1)}, f_{(1)}) y_{(2)} \bowtie f_{(2)} (\hat{S}(y_{(3)}), f_{(3)}) g) \\ 
&= (y_{(1)}, f_{(1)}) \hat{\Delta}(xy_{(2)}) (\hat{S}(y_{(3)}), f_{(4)}) (l^-(f_{(2)} g_{(1)}) \otimes l^+(f_{(3)} g_{(2)})) \\ 
&= \hat{\Delta}(x) (y_{(1)}, f_{(1)}) \hat{\Delta}(y_{(2)}) (\hat{S}(y_{(3)}), f_{(4)}) (l^-(f_{(2)}) l^-(g_{(1)}) \otimes l^+(f_{(3)}) l^+(g_{(2)})) \\ 
&= \hat{\Delta}(x) (l^-(f_{(1)}) \otimes l^+(f_{(2)})) \hat{\Delta}(y) (l^-(g_{(1)}) \otimes l^+(g_{(2)})) \\
&= \iota(x \bowtie f) \iota(y \bowtie g)
\end{align*} 
for $ x \bowtie f, y \bowtie g \in \DF(G_q) $, using 
\begin{align*}
((y_{(1)}, f_{(1)}) &\hat{\Delta}(y_{(2)}) (l^-(f_{(2)}) \otimes l^+(f_{(3)})), h \otimes k) \\
&= (y_{(1)}, f_{(1)}) (y_{(2)}, k_{(1)} h_{(1)}) (\R^{-1}, f_{(2)} \otimes h_{(2)})(\R, k_{(2)} \otimes f_{(3)}) \\
&= (y, k_{(1)} h_{(1)} f_{(1)}) (\R^{-1}, f_{(2)} \otimes h_{(2)})(\R, k_{(2)} \otimes f_{(3)}) \\
&= (\R^{-1}, f_{(1)} \otimes h_{(1)})(\R, k_{(1)} \otimes f_{(2)}) (y, f_{(3)} k_{(2)} h_{(2)}) \\
&= (\R^{-1}, f_{(1)} \otimes h_{(1)})(\R, k_{(1)} \otimes f_{(2)}) (y_{(1)}, k_{(2)} h_{(2)}) (y_{(2)}, f_{(3)}) \\
&= ((l^-(f_{(1)}) \otimes l^+(f_{(2)})) \hat{\Delta}(y_{(1)}) (y_{(2)}, f_{(3)}), h \otimes k). 
\end{align*} 
Hence $ \iota $ is an algebra homomorphism, and it is evident that $ \iota $ is essential.

We next show that $\iota$ is injective.  For this, note that the image of $\iota$ belongs to the space
\begin{align*}
M = \{ m \in \M(\DF(K_q) \otimes \DF(K_q)) \mid (z \otimes 1)m \in & \DF(K_q) \otimes \DF(K_q) \text{ for all } z \in \DF(K_q) \}   
\end{align*}
Note that the Galois map on $\DF(K_q) \otimes \DF(K_q)$ defined by
$$
 x\otimes y  \mapsto x\hat{S}(y_{(1)})\otimes y_{(2)}
$$
extends to a linear isomorphism $\gamma:M\to M$.  It suffices therefore to prove the injectivity of the composition
\[
 i:\DF(G_q) \to \DF(G_q) \stackrel{\iota}{\longrightarrow} M \stackrel{\gamma}{\longrightarrow} M,
\]
where the first map is the linear isomorphism sending $x\bowtie f$ to $(1\bowtie f)(x\bowtie 1)$.
We calculate
\[
  i(x \bowtie f) = l^-(f_{(1)}) l^+(S(f_{(2)})) \otimes l^+(f_{(3)}) x = I(f_{(1)}) \otimes l^+(f_{(2)}) x.
\]
where $I:\Poly(G_q) \to FU_q(\lie{g})$ is the isomorphism of Proposition \ref{adjointdualityplus}.  

Consider then an element in the kernel of $i$.  We can write it as a finite linear combination
\[
  \sum_{\mu\in\weights^+} \sum_{j,k}  x^\mu_{jk} \bowtie u^\mu_{jk},
\]
where $u^\mu_{jk} = \bra e^\mu_j | \bullet | e^\mu_k \ket$ are the matrix coefficients with respect to an orthonormal weight basis $(e^\mu_j)$ of $V(\mu)$, and $x^\mu_{jk}$ are some elements of $\DF(K_q)$.  Then by assumption we have
\[
  i\big(\sum_{\mu\in\weights^+} \sum_{j,k} x^\mu_{jk} \bowtie u^\mu_{jk}\, \big) = \sum_{\mu\in\weights^+} \sum_{j,r,k} I(u^\mu_{jr}) \otimes l^+(u^\mu_{rk})x_{jk}^\mu = 0.
\]
Since $I$ is an isomorphism and the $u^\mu_{jr}$ are linearly independent, this implies  that for all fixed $\mu,j,r$,
\[
  \sum_{k} l^+(u^\mu_{rk})x^\mu_{jk} = 0.
\]

Let us write $\epsilon_k$ for the weight of $e^\mu_k$. 
From the form of the universal $R$-matrix, see Theorem \ref{Rmatrixformula}, we see that $l^+(u^\mu_{k'k}) = 0$ unless $\epsilon_k \geq \epsilon_{k'}$.  Moreover, $l^+(u^\mu_{kk})$ is invertible in $U_q(\lie{g})$, while $l^+(u^\mu_{k'k})=0$ if $\epsilon_{k'} = \epsilon_k$ but $k'\neq k$.

Fix $\mu$ and $j$ and suppose $x^\mu_{jk} \neq0$ for some $k$.  Choose $k'$ with $\epsilon_{k'}$ maximal amongst all $\epsilon_k$ where $x^\mu_{jk}\neq0$.  In this case we have $l^+(u^\mu_{k'k})x^\mu_{jk} =0$ for all $k$ except $k=k'$ and hence $x^\mu_{jk'}=0$, a contradiction.  It follows that $x^\mu_{jk}=0$ for all $\mu,j,k$, as desired.
\end{proof}

Let us now give the definition of the quantized universal enveloping algebra of a complex group. 

\begin{definition} 
The quantized universal enveloping algebra of the real Lie algebra underlying $ \mathfrak{g} $ is 
$$ 
U_q^\mathbb{R}(\mathfrak{g}) = U_q^\mathbb{R}(\mathfrak{k}) \bowtie \CF^\infty(K_q),
$$
\nomenclature{$U_q^\mathbb{R}(\lie{g})$}{quantized enveloping algebra of $\lie{g}$ as a real Lie algebra}%
equipped with the standard Hopf $ * $-algebra structure. 
\end{definition} 

	Explicitly, the multiplication in $ U_q^\mathbb{R}(\mathfrak{g}) $ is given by 
	$$
	(X \bowtie f)(Y \bowtie g) = X (Y_{(1)}, f_{(1)}) Y_{(2)} \bowtie f_{(2)} (\hat{S}(Y_{(3)}), f_{(3)}) g, 
	$$
	for $ X, Y \in U_q^\mathbb{R}(\mathfrak{k}) $ and $ f, g \in \CF^\infty(K_q) $ and the $*$-structure is given by
	\[
	(X\bowtie f)^* =  (1\bowtie f^*) (X^*\bowtie 1)
	= (X_{(1)}^*, f_{(1)}^*) \, X_{(2)}^* \bowtie f_{(2)}^* \,
	(\hat{S}(X_{(3)}^*),f_{(3)}^*).
	\]

We view $ U_q^\mathbb{R}(\mathfrak{g}) $ as a substitute of the universal enveloping algebra of the real Lie algebra $ \mathfrak{g} $. 
Note that $ U_q^\mathbb{R}(\mathfrak{g}) \subset \M(\DF(G_q)) $ naturally. 

The following result is originally due to Kr\"ahmer \cite{KraehmerFRT} and Arano \cite{Aranocomparison}. 

\begin{lemma} \label{uqrgdiagonalembedding}
The linear map $ \iota: U_q^\mathbb{R}(\mathfrak{g}) \rightarrow U_q(\mathfrak{g}) \otimes U_q(\mathfrak{g}) $ given by 
$$
\iota(X \bowtie f) = \hat{\Delta}(X) (l^-(f_{(1)}) \otimes l^+(f_{(2)})) 
$$
is an injective algebra homomorphism, and its image is 
\begin{align*}
\iota(U_q^\mathbb{R}(\mathfrak{g})) &= (FU_q(\mathfrak{g}) \otimes 1) \hat{\Delta}(U_q(\mathfrak{g})) \\
&= (1 \otimes \hat{S}^{-1}(FU_q(\mathfrak{g}))) \hat{\Delta}(U_q(\mathfrak{g})) \\
&= \hat{\Delta}(U_q(\mathfrak{g}))(FU_q(\mathfrak{g}) \otimes 1)  \\
&= \hat{\Delta}(U_q(\mathfrak{g})) (1 \otimes \hat{S}^{-1}(FU_q(\mathfrak{g}))). 
\end{align*}
In particular, we have an algebra isomorphism 
$$
U_q^\mathbb{R}(\mathfrak{g}) \cong FU_q(\mathfrak{g}) \rtimes U_q(\mathfrak{g})
$$
where $ U_q(\mathfrak{g}) $ acts on $ FU_q(\mathfrak{g}) $ via the adjoint action. 
\end{lemma} 

\begin{proof}
Since $ U_q^\mathbb{R}(\mathfrak{g}) \subset \M(\DF(G_q)) $ it follows from Proposition \ref{finiteseparating} that $ \iota $ defines an 
injective algebra homomorphism $ U_q^\mathbb{R}(\mathfrak{g}) \rightarrow U_q(\mathfrak{g}) \otimes U_q(\mathfrak{g}) $. This fact 
was proved in a slightly different way in \cite{KraehmerFRT}. 

For the second claim compare Section 4 in \cite{Aranocomparison}. We have 
\begin{align*}
\iota((1 \bowtie f)(X \bowtie 1)) &= (l^-(f_{(1)}) \otimes l^+(f_{(2)})) \hat{\Delta}(X) \\
&= (l^-(f_{(1)}) l^+(S(f_{(2)})) l^+(f_{(3)}) \otimes l^+(f_{(4)})) \hat{\Delta}(X) \\
&= (I(f_{(1)}) \otimes 1) \hat{\Delta}(l^+(f_{(2)})) \hat{\Delta}(X) \subset (FU_q(\mathfrak{g}) \otimes 1) \hat{\Delta}(U_q(\mathfrak{g})),  
\end{align*} 
where $ I $ is the isomorphism from Proposition \ref{adjointdualityplus}. 
This proves the inclusion $ \iota(U_q^\mathbb{R}(\mathfrak{g})) \subset (FU_q(\mathfrak{g}) \otimes 1) \hat{\Delta}(U_q(\mathfrak{g})) $. 
Conversely, for $ X = I(f) \in FU_q(\mathfrak{g}) $ and $ Y \in U_q(\mathfrak{g}) $ we have 
\begin{align*}
(I(f) \otimes 1) \hat{\Delta}(Y) &= (l^-(f_{(1)}) l^+(S(f_{(2)})) \otimes 1) \hat{\Delta}(Y) \\
&= (l^-(f_{(1)}) \otimes l^+(f_{(2)}))(l^+(S(f_{(4)})) \otimes l^+(S(f_{(3)}))) \hat{\Delta}(X) 
\subset \iota(U_q^\mathbb{R}(\mathfrak{g})). 
\end{align*} 
Using the relations $ X \otimes 1 = (1 \otimes \hat{S}^{-1}(X_{(3)}))(X_{(1)} \otimes X_{(2)}) $ 
and $ 1 \otimes \hat{S}^{-1}(X) =  (X_{(3)} \otimes 1)(\hat{S}^{-1}(X_{(2)}) \otimes \hat{S}^{-1}(X_{(1)})) $,  
and the fact that $ FU_q(\mathfrak{g}) $ satisfies $ \hat{\Delta}(FU_q(\mathfrak{g})) \subset U_q(\mathfrak{g}) \otimes FU_q(\mathfrak{g}) $ we obtain 
\begin{align*}
(FU_q(\mathfrak{g}) \otimes 1) \hat{\Delta}(U_q(\mathfrak{g})) = (1 \otimes \hat{S}^{-1}(FU_q(\mathfrak{g}))) \hat{\Delta}(U_q(\mathfrak{g})). 
\end{align*} 
The remaining equalities follow from the fact that $ \iota(U_q^\mathbb{R}(\mathfrak{g})) $ is an algebra, combined with 
\begin{align*}
(Y \otimes 1) \hat{\Delta}(X) &= X_{(3)} \hat{S}^{-1}(X_{(2)}) Y X_{(1)} \otimes X_{(4)} \\
&= X_{(2)} \hat{S}^{-1}(X_{(1)}) \rightarrow Y \otimes X_{(3)} \in \hat{\Delta}(U_q(\mathfrak{g})) (FU_q(\mathfrak{g}) \otimes 1)
\end{align*}
and 
\begin{align*}
(1 \otimes \hat{S}^{-1}(Y)) \hat{\Delta}(\hat{S}^{-1}(X)) &= \hat{S}^{-1}(X_{(4)}) \otimes \hat{S}^{-1}(X_{(3)}) X_{(2)} \hat{S}^{-1}(Y) \hat{S}^{-1}(X_{(1)}) \\
&= \hat{S}^{-1}(X_{(3)}) \otimes \hat{S}^{-1}(X_{(2)}) \hat{S}^{-1}(X_{(1)} \rightarrow Y) \\
&\qquad \in \hat{\Delta}(U_q(\mathfrak{g})) (1 \otimes \hat{S}^{-1}(FU_q(\mathfrak{g})))
\end{align*}
for $ Y \in FU_q(\mathfrak{g}) $ and $ X \in U_q(\mathfrak{g}) $.

Finally, we have a linear 
isomorphism $ FU_q(\mathfrak{g}) \rtimes U_q(\mathfrak{g}) \rightarrow (FU_q(\mathfrak{g}) \otimes 1) \hat{\Delta}(U_q(\mathfrak{g})) $ 
given by $ \gamma(X \rtimes Y) = (X \otimes 1) \hat{\Delta}(Y) $, and since 
\begin{align*}
\gamma((X \rtimes Y)(X' \rtimes Y')) &= \gamma(X Y_{(1)} X' \hat{S}(Y_{(2)}) \rtimes Y_{(3)} Y') \\
&= X Y_{(1)} X' Y'_{(1)} \otimes Y_{(2)} Y'_{(2)} \\
&= (X \otimes 1) \hat{\Delta}(Y) (X' \otimes 1) \hat{\Delta}(Y') = \gamma(X \rtimes Y) \gamma(X' \rtimes Y')
\end{align*}
the map $ \gamma $ is compatible with multiplication. 
\end{proof} 

We remark that $ \iota $ is not a homomorphism of coalgebras; in fact there is no bialgebra structure on $ U_q(\mathfrak{g}) \otimes U_q(\mathfrak{g}) $ 
for which $ \iota $ becomes a homomorphism of coalgebras. Let us also point out that $ \iota $ is not a $ * $-homomorphism in a natural way.

Using Lemma \ref{uqrgdiagonalembedding} we can determine the centre of $ U_q^\mathbb{R}(\mathfrak{g}) $. 

\begin{lemma} \label{uqrgcenter}
The centre $ ZU_q^\mathbb{R}(\mathfrak{g}) $ of $ U_q^\mathbb{R}(\mathfrak{g}) $ is isomorphic to $ ZU_q(\mathfrak{g}) \otimes ZU_q(\mathfrak{g}) $.
\nomenclature{$ZU_q^\mathbb{R}(\lie{g})$}{centre of $U_q^\mathbb{R}(\lie{g})$}%
\end{lemma} 

\begin{proof} 
Lemma \ref{uqrgdiagonalembedding} shows that $ ZU_q(\mathfrak{g}) \otimes ZU_q(\mathfrak{g}) \subset Z(\im(\iota)) $, where 
$ \iota: U_q^\mathbb{R}(\mathfrak{g}) \rightarrow U_q(\mathfrak{g}) \otimes U_q(\mathfrak{g}) $ is as above. 

Conversely, any element $ X \in \im(\iota) $ which commutes with $ \im(\iota) $ must in particular 
commute with the elements $ K_{2\mu} \otimes 1 $ and $ 1 \otimes K_{-2\mu} $ for $ \mu \in \weights^+ $. It follows that $ X $ 
has weight $ 0 $ with respect to the diagonal action of $ U_q(\mathfrak{h}) $ on both left and right tensor factors. 
Hence $ X $ commutes with $ U_q(\mathfrak{h}) FU_q(\mathfrak{g}) = U_q(\mathfrak{g}) $ in both tensor factors, and therefore is 
contained in $ ZU_q(\mathfrak{g}) \otimes ZU_q(\mathfrak{g}) $. \end{proof} 

We have a bilinear pairing between $ U_q(\mathfrak{g}) \otimes U_q(\mathfrak{g}) $ and $ \Poly(G_q) \bowtie \Poly(G_q) $ given by 
$$
(X \otimes Y, f \bowtie g) = (X, g) (Y, f).  
$$
\label{nom:UU-OO-pairing}%
Using this pairing, elements of $ U_q(\mathfrak{g}) \otimes U_q(\mathfrak{g}) $ can be viewed as linear functionals on $ \Poly(G_q) \bowtie \Poly(G_q) $. 
For $ X \bowtie f \in U_q^\mathbb{R}(\mathfrak{g}) $ and $ g \bowtie h \in \Poly^\mathbb{R}(G_q) $ we compute 
\begin{align*}
(\iota(X \bowtie f), &g \bowtie h) = (X_{(1)} l^-(f_{(1)}) \otimes X_{(2)} l^+(f_{(2)}), g \bowtie h) \\
&= (X_{(1)}, h_{(1)}) (\R^{-1}, f_{(1)} \otimes h_{(2)}) (X_{(2)}, g_{(1)}) (\R, g_{(2)} \otimes f_{(2)}) \\
&= (X_{(2)}, g_{(1)})(X_{(1)}, h_{(1)}) (\R^{-1}, g_{(2)} \otimes S^{-1}(f_{(2)})) (\R, S^{-1}(f_{(1)}) \otimes h_{(2)}) \\
&= (X, g_{(1)} h_{(1)}) (l^-(g_{(2)}) l^+(h_{(2)}), S^{-1}(f)) \\
&= (X \bowtie f, g_{(1)} h_{(1)} \otimes l^-(g_{(2)}) l^+(h_{(2)})) \\
&= (X \bowtie f, i(g \bowtie h)) 
\end{align*}
where $i:\Poly^\mathbb{R}(G_q) \to \CF^\infty(G_q)$ is the homomorphism from Proposition \ref{polynomialtosmooth}.

\subsection{Parabolic quantum subgroups} 
\label{sec:parabolic_subgroups}

In this section we describe briefly how to obtain quantum analogues of parabolic subgroups in complex semisimple Lie groups. We continue to use 
the notation introduced in previously. 

Recall that $ \simpleroots = \{\alpha_1, \dots, \alpha_N \} $ denotes the set of simple roots, and let $ S \subset \simpleroots $ be a subset. 
We obtain a corresponding Hopf $ * $-subalgebra 
$ U_q^\mathbb{R}(\mathfrak{k}_S) \subset U_q^\mathbb{R}(\mathfrak{k}) $ generated by all $ K_\lambda $ for $ \lambda \in \weights $ 
\nomenclature{$U_q(\lie{k}_s)$}{Hopf $*$-subalgebra of $U_q(\lie{k})$ corresponding to a subset of simple roots $S$}%
together with the generators $ E_i, F_i $ for $ \alpha_i \in S $. 
Notice that $ U_q^\mathbb{R}(\mathfrak{k}_\Sigma) = U_q^\mathbb{R}(\mathfrak{k}) $ 
and $ U_q^\mathbb{R}(\mathfrak{k}_\emptyset) = U_q^\mathbb{R}(\mathfrak{t}) $. 
The inclusion $ U_q^\mathbb{R}(\mathfrak{k}_S) \rightarrow U_q^\mathbb{R}(\mathfrak{k}) $ induces a map 
$ U_q^\mathbb{R}(\mathfrak{k})^* \rightarrow U_q^\mathbb{R}(\mathfrak{k}_S)^* $, and we denote by $ \CF^\infty(K_{q, S}) $
\nomenclature[o$C^\infty(K_q,S)$]{$\CF^\infty(K_{q,S})$}{Hopf $*$-subalgebra of $\CF^\infty(K_q)$ corresponding to a subset $S$ of simple roots}%
\nomenclature{$K_{q,S}$}{compact quantum subgroup of $K_q$ corresponding to a subset $S$ of simple roots}%
the Hopf $ * $-algebra 
obtained as the image of $ \CF^\infty(K_q) $ under this map. 
In this way, $ K_{q,S} $ is a closed quantum subgroup of $ K_q $. 
Notice that $ K_{q, \Sigma} = K_q $, and that $ K_{q, \emptyset} = T $ is the classical maximal torus inside $ K_q $.

We define the parabolic quantum subgroup $ P_q \subset G_q $ associated to the set $ S $ as the Drinfeld double 
$$ 
\CF^\infty_c(P_q) = \CF^\infty_c(K_{q, S} \bowtie \hat{K}_q) = \CF^\infty(K_{q,S}) \otimes \DF(K_q) 
$$ 
\nomenclature{$P_q$}{parabolic quantum subgroup of $G_q$}%
\nomenclature[o$C^\infty(P_q)$]{$\CF_c^\infty(P_q)$}{algebra of functions on the parabolic quantum subgroup $P_q$}
\nomenclature{$\pi_S$}{quotient map $\pi_S:\CF^\infty(K_q) \to \CF^\infty(K_{q,S})$}%
with the coproduct induced from $ \CF^\infty_c(G_q) $. That is, in the formula for the comultiplication of $ \CF^\infty_c(G_q) $ one has to 
replace $ W $ by $ (\pi_S \otimes \id)(W) \in \M(\CF^\infty(K_{q,S}) \otimes \DF(K_q)) $, where $ \pi_S: \CF^\infty(K_q) \rightarrow \CF^\infty(K_{q, S}) $ 
is the canonical quotient map. 

In a similar way we construct the quantized universal enveloping algebras of the corresponding real Lie groups. More precisely, 
we define the Hopf $ * $-algebra $ U_q^\mathbb{R}(\mathfrak{p}) $ associated to the parabolic quantum subgroup $ P_q $ by 
$$
U_q^\mathbb{R}(\mathfrak{p}) = U_q^\mathbb{R}(\mathfrak{k}_S) \bowtie \CF^\infty(K_q) 
$$
\nomenclature{$U_q^\mathbb{R}(\lie{p})$}{quantized enveloping algebra associated to the parabolic quantum subgroup $P_q$}%
where the pairing between $ U_q^\mathbb{R}(\mathfrak{k}_S) $ and $ \CF^\infty(K_q) $ 
used in the definition of the double is induced from the canonical pairing 
between $ U_q^\mathbb{R}(\mathfrak{k}) $ and $ \CF^\infty(K_q) $. By construction, we have a canonical inclusion homomorphism 
$ U_q^\mathbb{R}(\mathfrak{p}) \rightarrow U_q^\mathbb{R}(\mathfrak{g}) $. 

The quantum Borel subgroup $B_q \subset G_q$ is defined to be the parabolic subgroup corresponding to $S=\emptyset$.
 Explicitly, 
$$ 
\CF^\infty_c(B_q) = \CF^\infty_c(T \bowtie \hat{K}_q) = \CF^\infty(T) \otimes \DF(K_q)
$$
\nomenclature{$B_q$}{quantum Borel subgroup $B_q = T\bowtie \hat{K}_q$ of $G_q$}%
with the comultiplication twisted by the bicharacter as explained above. The corresponding quantized universal enveloping algebra is 
$$
U_q^\mathbb{R}(\mathfrak{b}) = U_q^\mathbb{R}(\mathfrak{t}) \bowtie \CF^\infty(K_q). 
$$
\nomenclature{$U_q^\mathbb{R}(\lie{b})$}{quantized enveloping algebra associated to the quantum Borel subgroup $B_q$}%
In the study of principal series representations we will only work with the Borel quantum subgroup and its quantized universal enveloping algebra.

\newpage  

\section{Category $ \O $} \label{chcato}

In this chapter we study some aspects of the representation theory of quantized universal enveloping algebras with applications to the theory 
of Yetter-Drinfeld modules and complex quantum groups. Our main goal is a proof of the Verma module annihilator Theorem, following the 
work of Joseph and Letzter. We refer to \cite{FarkasLetztersurvey} for a survey of the ideas involved in the proof and background. 

In this chapter, unless explicitly stated otherwise, we shall work over $ \mathbb{K} = \mathbb{C} $ and assume that $ 1 \neq q = e^h $ is positive. 
We shall also use the notation $ \hbar = \tfrac{h}{2 \pi} $. 
Recall that in this situation, we have the identification $\lie{h}_q^* = \lie{h}^*/i\hbar^{-1}\roots^\vee$, see Subsection \ref{sec:weights}.

\subsection{The definition of category $ \O $} \label{secdefo}

In this section we introduce category $ \O $ for $ U_q(\mathfrak{g}) $, compare Section 4.1.4 in \cite{Josephbook} and \cite{AMcategoryo}, \cite{HumphreysO}.

We start with the following definition, compare \cite{AMcategoryo}. 
\begin{definition} \label{defcategoryO}
\nomenclature[o$O$]{$\O$}{category $\O$ for $U_q(\lie{g})$}%
A left module $ M $ over $ U_q(\mathfrak{g}) $ is said to belong to category $ \O $ if 
\begin{bnum}
\item[a)] $ M $ is finitely generated as a $ U_q(\mathfrak{g}) $-module. 
\item[b)] $ M $ is a weight module, that is, a direct sum of its weight spaces $ M_\lambda $ for $ \lambda \in \mathfrak{h}^*_q $. 
\item[c)] The action of $ U_q(\mathfrak{n}_+) $ on $ M $ is locally nilpotent. 
\end{bnum}
Morphisms in category $ \O $ are all $ U_q(\mathfrak{g}) $-linear maps. 
\end{definition} 

Note that category $\O$ is closed under taking submodules and quotient modules.  Specifically, finite generation passes to submodules by Noetherianity of $U_q(\lie{g})$, see Proposition \ref{uqgnoetheriandomain}, and submodules and quotients of weight modules are again weight modules.
Local nilpotency is obvious in either case.

Due to finite generation, any module $ M $ in category $ \O $ satisfies $ \dim M_\lambda < \infty $ for all $ \lambda \in \mathfrak{h}^*_q $. 
We define the formal character of $ M $ by setting 
$$
\ch(M) = \sum_{\lambda \in \mathfrak{h}^*_q} \dim(M_\lambda) e^\lambda. 
$$
\nomenclature[o$ch$]{$\ch$}{formal character of a module in category $\O$}%
Here the expression on the right hand side is interpreted as a formal sum.

More generally, we will consider formal sums of the form $\sum_{\lambda\in\lie{h}_q^*} f(\lambda) e^\lambda$ where $f:\lie{h}^*_q \to \ZZ$ is any integer valued function whose support lies in a finite union of sets of the form $\nu - \roots^+$ with $\nu\in\lie{h}_q^*$.  For such formal sums, there is a well-defined convolution product given by
\[
  \left(\sum_{\lambda\in\lie{h}_q^*} f(\lambda)e^\lambda \right)
  \left(\sum_{\mu\in\lie{h}_q^*} g(\mu)e^\mu \right)   = \sum_{\lambda,\mu\in\lie{h}_q^*} f(\lambda) g(\mu) e^{\lambda+\mu}.
\]

In particular, let us define
\[
   p = \prod_{\beta\in{\bf \Delta}^+} \left(\sum_{m=0}^\infty e^{-m\beta}\right) = \sum_{\nu\in\roots^+} P(\nu)e^{-\nu},
\]
\nomenclature{$p$}{formal character of $M(0)$}%
where $P$ is Kostant's partition function,
 \[
  P(\nu) = \left| \{ (r_1,\ldots,r_n)\in\NN_0^n \mid r_1\beta_1 + \cdots + r_n\beta_n = \nu \}\right|.
\]
\label{nom:Kostant_partition_function2}%
Here $\beta_1,\ldots,\beta_n$ are the positive roots.
If $ \mu \in \mathfrak{h}^*_q $ and $ M = M(\mu) $ is the Verma module with highest weight $ \mu $ then Proposition \ref{prop:Uqn_weight_multiplicities} immediately shows that the character of $ M $ is given by 
\[
  \ch(M(\mu)) = e^\mu p.
\]

\subsubsection{Category $ \O $ is Artinian} 
\label{sec:O_Artinian}

In this subsection we discuss finiteness properties of category $ \O $. 

Let us first show that every module $ M $ in $ \O $ admits a finite filtration by highest weight modules, compare Section 1.2 in \cite{HumphreysO}. 
Recall that a $ U_q(\mathfrak{g}) $-module $ M $ is called a highest weight module of highest weight $ \lambda $ if there exists a primitive vector $ v_\lambda $ 
in $ M $ which generates $ M $. 
Recall also that we have a partial order on $\lie{h}_q^*$ defined by $\lambda\leq\mu$ if and only if $\mu-\lambda\in\roots^+$, see Section \ref{secvermauqg}.

\begin{lemma} \label{lemascendinghighestweightmodules}
Let $ M $ be a nonzero module in $ \O $. Then $ M $ has a finite filtration 
$$
0 = M_0 \subset M_1 \subset \cdots \subset M_n = M 
$$
such that each subquotient $ M_j/M_{j - 1} $ for $ 1 \leq j \leq n $ is a highest weight module. 
\end{lemma} 

\begin{proof}  According to condition $ a) $ in Definition \ref{defcategoryO} we see that $ M $ is generated by finitely many weight vectors. 
Consider the $ U_q(\mathfrak{n}_+) $-module $ V $ generated by a finite generating set of weight vectors, which by condition $ c) $ 
is finite dimensional. 

We prove the claim by induction on $ \dim(V) $, the case $ \dim(V) = 1 $ being trivial. 
For the inductive step pick a weight in $ \mathfrak{h}^*_q $ which is maximal among the weights appearing in $ V $. 
Then we can find a corresponding primitive vector $ v \in M $ and obtain an associated submodule $ M_0 = U_q(\mathfrak{g}) \cdot v \subset M $. 
The quotient $ M/M_0 $ 
is again finitely generated by the $ U_q(\mathfrak{n}_+) $-module $ V/\mathbb{C}v $. 
We can therefore apply our inductive hypothesis to obtain a filtration of the desired type 
for $ M/M_0 $, and pulling this filtration back to $ M $ yields the assertion. \end{proof}

Recall that a module $ M $ over a ring $ R $ is called Noetherian (Artinian) if every ascending chain $ M_1 \subset M_2 \subset M_3 \subset \cdots $ of submodules 
(every descending chain $ M_1 \supset M_2 \supset \cdots $ of submodules) becomes stationary, that is, if there exists $ n \in \mathbb{N} $ 
such that $ M_n = M_{n + 1} = M_{n + 2} = \cdots $. 

We show that every module in $ \O $ is Artinian and Noetherian, compare Section 1.11 in \cite{HumphreysO}. 

\begin{theorem} \label{catoartinian}
Every module $ M $ in category $ \O $ is Artinian and Noetherian. Moreover $ \dim \Hom_{U_q(\mathfrak{g})}(M,N) < \infty $ for all $ M, N \in \O $. 
\end{theorem} 

\begin{proof}  According to Lemma \ref{lemascendinghighestweightmodules}, any nonzero module $ M $ in $ \O $ has a finite filtration with subquotients 
given by highest weight modules. Hence it suffices to treat the case where $ M = M(\lambda) $ for $ \lambda \in \mathfrak{h}^*_q $ is a Verma module. 

Consider the finite dimensional 
subspace $ V = \bigoplus_{\hat{w} \in \hat{W}} M(\lambda)_{\hat{w} . \lambda} $ of $ M = M(\lambda) $, where we recall that $ \hat{W} $ is the 
extended Weyl group defined in Section \ref{seczuqg}. 
Assume that $ N_1 \subset N_2 $ is a proper inclusion of submodules of $ M $. 
Then $ ZU_q(\mathfrak{g}) $ acts on $ N_1, N_2 $ and $ N_2/N_1 $ by the central character $ \xi_\lambda $. 
Since $ N_2 \subset M $ the module $ N_2/N_1 $ contains a primitive vector $ v_\mu $ of some weight $ \mu \leq \lambda $. 
Therefore $ \xi_\mu = \xi_\lambda $. According to Theorem \ref{chilinkage} this implies $ \mu = w . \lambda $ for some $ w \in \hat{W} $. 
We conclude $ N_2 \cap V \neq 0 $ and $ \dim(N_2 \cap V) > \dim(N_1 \cap V) $. Since $ V $ is finite dimensional this means that any 
strictly ascending or descending chain of submodules will have length at most $ \dim(V) $. 
In particular, every module $ M $ in $ \O $ is both Artinian and Noetherian. 

Any $ U_q(\mathfrak{g}) $-linear map $ M \rightarrow N $ is determined by its values on a finite generating set of weight vectors. Since the weight 
spaces of $ N $ are finite dimensional we conclude $ \dim \Hom_{U_q(\mathfrak{g})}(M,N) < \infty $ for all $ M, N \in \O $. \end{proof}  

Of course, Noetherianity of modules in $ \O $ follows also from the fact that $ U_q(\mathfrak{g}) $ is Noetherian, see Theorem \ref{uqgnoetheriandomain}. 

Due to Theorem \ref{catoartinian} one can apply Jordan-H\"older theory to category $ \O $. More precisely, every module $ M \in \O $ 
has a decomposition series $ 0 = M_0 \subset M_1 \subset \cdots \subset M_n = M $ such that all subquotients $ M_{j + 1}/M_j $ are simple 
highest weight modules. Moreover, the number of subquotients isomorphic to $ V(\lambda) $ for $ \lambda \in \mathfrak{h}^*_q $ is independent 
of the decomposition series
and will be denoted by $ [M: V(\lambda)] $. Note that this gives the formula 
$$
\ch(M) = \sum_{\lambda \in \mathfrak{h}^*_q} [M: V(\lambda)] \ch(V(\lambda)) 
$$
for the character of $ M $, with only finitely many nonzero terms on the right hand side. 

In particular, for any $ \mu \in \mathfrak{h}^*_q $ the character of the Verma module $ M(\mu) $ can be written as 
$$
\ch(M(\mu)) = \sum_{\lambda \in \mathfrak{h}^*_q} [M(\mu): V(\lambda)] \ch(V(\lambda)). 
$$
Note that $ [M(\mu): V(\mu)] = 1 $, and moreover $ [M(\mu): V(\lambda)] = 0 $ unless $\lambda\leq \mu$ and $\lambda$ is $\hat{W}$-linked to $\mu$, due to Theorem \ref{chilinkage}.  It follows that
the character of the simple module $ V(\mu) $ can be expressed in the form 
$$
\ch(V(\mu)) = \sum_{\lambda \in W.\mu} m_\lambda \ch(M(\lambda)) 
$$
for certain integers $ m_\lambda $. For $ \mu \in \weights^+ $ these coefficients are given as follows. 

\begin{prop} \label{charactersimple}
Let $ \mu \in \weights^+ $. Then we have 
$$
\ch(V(\mu)) = \sum_{w \in W} (-1)^{l(w)} \ch(M(w . \mu)).  
$$
\end{prop} 

\begin{proof} 
The above discussion shows that
\[
 \ch(V(\mu)) = \sum_{w \in W} m_{w.\mu} \ch(M(w.\mu))
\]
for some integers $m_{w.\mu}$.  We need to show that $m_{w.\mu} = (-1)^{l(w)}$.

Recall that $\ch(M(\mu)) = e^\mu p$  
where $p = \prod_{\beta\in{\bf\Delta}^+} (1 + e^{-\beta} + e^{-2\beta}+\cdots)$.  We introduce the formal sum
\[
  q = \prod_{\beta\in{\bf\Delta}^+} (e^{\beta/2} - e^{-\beta/2})
    = e^\rho \prod_{\beta\in{\bf\Delta}^+} (1 - e^{-\beta}),
\]
and note that the right hand expression gives $pq=e^\rho$.  Therefore, we obtain
\[
  \ch(V(\mu)) q = \sum_{w \in W} m_{w.\mu} e^{w.\mu} pq = \sum_{w \in W} m_{w.\mu} e^{w\mu}.
\]

Consider the Weyl group action on finite formal sums defined by $ w e^\lambda = e^{w \lambda} $ for $ w \in W $ 
and $ \lambda \in \mathfrak{h}^*_q $.  Note that $q$ is alternating with respect to this action, meaning that
$
 w q = (-1)^{l(w)}q
$
for all $w\in W$.  
Moreover, $\ch(V(\mu))$ is $W$-invariant by Proposition \ref{prop:weights_Weyl_invariant}.  Therefore, $\ch(V(\mu)) q$ is alternating, and hence $m_{w.\mu} = (-1)^{l(w)} m_\mu$ for all $w\in W$.  Since $m_\mu=1$, the result follows.
\end{proof}

As a consequence, the weight spaces of the irreducible module $ V(\mu) $ for $ \mu \in \weights^+ $ have the same dimensions 
as in the classical case, and the dimension of $ V(\mu) $ is given by the Weyl dimension formula.

\subsubsection{Duality} 

The appropriate duality operation in category $\O$ is different from the usual duality for modules over Hopf algebras defined 
in terms of the antipode.

We define duals in category $ \O $ in the same way as we did for general weight modules in Section \ref{secvermauqg}.  Namely, if
$ M $ is a module in category $ \O $ we let the dual of $ M $ be the $ U_q(\mathfrak{g}) $-module 
$$ 
M^\vee = \bigoplus_{\lambda \in \mathfrak{h}_q^*} \Hom(M_\lambda, \mathbb{C}), 
$$ 
\label{nom:restricted_dual2}%
with the left $ U_q(\mathfrak{g}) $-module structure given by 
$$
(X \cdot f)(v) = f(\tau(X) \cdot v), 
$$
where $ \tau $ is the involution from Definition \ref{deftau}. 
Notice that we have a canonical isomorphism $ M^{\vee \vee} \cong M $ since all weight spaces are finite dimensional and $ \tau $ is involutive. 
If $ 0 \rightarrow K \rightarrow M \rightarrow Q \rightarrow 0 $ is an exact sequence of modules in category $ \O $ then 
the dual sequence $ 0 \rightarrow Q^\vee \rightarrow M^\vee \rightarrow K^\vee \rightarrow 0 $ is again exact. 

Let us show that the dual $ M^\vee $ of a module $ M $ in category $ \O $ is again in category $ \O $. It is clear that $ M^\vee $ is 
a weight module such that $ (M^\vee)_\lambda = \Hom(M_\lambda, \mathbb{C}) $. Local nilpotency of the action of $ U_q(\mathfrak{n}_+) $ 
follows from the fact that the elements in $ M^\vee $ are supported on only finitely many weight spaces. 

To see that $ M^\vee $ is finitely generated we argue as follows, see \cite{HumphreysO}. 
Assume $ M^\vee $ is not finitely generated. Then we can find a strictly increasing infinite sequence 
$ 0 = U_0 \subset U_1 \subset \cdots \subset M^\vee $ of finitely generated submodules $ U_k \subset M^\vee $. Setting $ Q_j = M^\vee/U_j $ 
we obtain a corresponding infinite sequence $ M^\vee = M^\vee/0 \rightarrow Q_1 \rightarrow Q_2 \rightarrow \cdots $ of quotient 
modules of $ M^\vee $ and surjective module maps. By exactness of the duality functor this in turn leads to a strictly decreasing 
infinite sequence of submodules of $ M $. As we have seen in the proof of Theorem \ref{catoartinian} this is impossible. 

In summary, we conclude that sending $ M $ to $ M^\vee $ defines a contravariant involutive self-equivalence of category $ \O $.

\subsubsection{Dominant and antidominant weights} 
\label{sec:dominant_and_antidominant}

In this subsection we discuss the notion of dominant and antidominant weights in analogy to the classical theory, see Chapter 10 in \cite{HumphreysO}. 
There are some new features in the quantum setting due to the exponentiation in the Cartan part of the quantized universal 
enveloping algebra. 

Recall from Definition \ref{def:extended_Weyl_group} that the extended Weyl group is defined as $\hat{W} = \mathbf{Y}_q \rtimes W$, where $\mathbf{Y}_q$ 
denotes the subgroup of elements of order at most $2$ in $\lie{h}_q^*$.  Under the identification $\lie{h}_q^* = \lie{h}^*/i\hbar^{-1}\roots^\vee$ 
we have $\mathbf{Y}_q = \tfrac{1}{2}i\hbar^{-1}\roots^\vee / i\hbar^{-1}\roots^\vee$.  The extended Weyl group acts on $\lie{h}_q^*$ by
\[
(\zeta,w) \lambda = w\lambda+\zeta,
\]
for $\zeta\in\tfrac{1}{2}i\hbar^{-1}\roots^\vee / i\hbar^{-1}\roots^\vee$, $ w \in W $ and $ \lambda \in \lie{h}_q^* $, and also by the shifted action
\[
(\zeta,w) \,.\, \lambda = w.\lambda+\zeta = w(\lambda+\rho)-\rho+\zeta.
\]
Two elements of $ \lie{h}_q^* $ are $ \hat{W} $-linked if they lie in the same orbit of the shifted $ \hat{W} $-action.

Following the notation in \cite{HumphreysO}, we define 
$$
{\bf \Delta}_{[\lambda]} = \{\alpha \in {\bf \Delta} \mid q_\alpha^{(\lambda, \alpha^\vee)} \in \pm q_\alpha^{\mathbb{Z}} \}.
$$
\nomenclature[o$\Delta_\lambda$]{${\bf \Delta}_{[\lambda]}$}{set of roots which are integral with respect to $\lambda\in\lie{h}_q^*$}%
Here we are writing $q_\alpha^{(\lambda,\mu)}$ to denote $q^{(\lambda,d_\alpha\mu)}$ for $\mu\in\weights$.
Note that $ {\bf \Delta}_{[\mu]} = {\bf \Delta}_{[\lambda]} $ if $ \mu \in \lambda + \weights $,
so we may equivalently characterize $ {\bf \Delta}_{[\lambda]} $ by 
$$ 
{\bf \Delta}_{[\lambda]} = \{\alpha \in {\bf \Delta} \mid q_\alpha^{(\lambda + \rho, \alpha^\vee)} \in \pm q_\alpha^{\mathbb{Z}} \}. 
$$
Next we define 
$$
\hat{W}_{[\lambda]} = \{\hat{w} \in \hat{W} \mid \hat{w} \lambda - \lambda \in \roots \}.
$$
\nomenclature[o$W$23]{$\hat{W}_{[\lambda]}$}{subgroup of $\hat{W}$ preserving $\lambda$ modulo $\roots$}%
Note that if $\mu=\lambda+\nu$ with $\nu\in\weights$ then for any $\hat{w} = (\zeta,w) \in \hat{W}_{[\lambda]}$ we have
\[
(\zeta,w) \mu- \mu = (w \lambda + \zeta - \lambda) + (w\nu - \nu) = ((\zeta,w)\lambda - \lambda) + (w\nu - \nu) \in \roots.
\]
It follows that $\hat{W}_{[\mu]} = \hat{W}_{[\lambda]}$ whenever $\mu \in \lambda+\weights$, and therefore
$$ 
\hat{W}_{[\lambda]} = \{\hat{w} \in W \mid \hat{w} \,.\, \lambda - \lambda \in \roots\} .
$$ 
We define $W_{[\lambda]}$ as the image of $\hat{W}_{[\lambda]}$ under the canonical projection $\hat{W} \to W$. 
\nomenclature[o$W$24]{$W_{[\lambda]}$}{projection of $\hat{W}_{[\lambda]}$ to $W$}%
Note that $ \hat{W}_{[\lambda]} \cap (\tfrac{1}{2}i\hbar^{-1}\roots^\vee / i\hbar^{-1}\roots^\vee\rtimes \{e\})$ is trivial, so the projection 
map $\hat{W}_{[\lambda]} \to W_{[\lambda]}$ is in fact an isomorphism.
Explicitly, we have
\begin{align*}
W_{[\lambda]} &= \{w \in W \mid w \lambda - \lambda \in \roots +\tfrac{1}{2}i\hbar^{-1}\roots^\vee / i\hbar^{-1}\roots^\vee \} \\
&= \{w \in W \mid w \,.\, \lambda - \lambda \in \roots +\tfrac{1}{2}i\hbar^{-1}\roots^\vee / i\hbar^{-1}\roots^\vee \}.
\end{align*}

In the next proposition, we write $ E $ for the $ \mathbb{R} $-span of the root system $ {\bf \Delta} \subset \mathfrak{h}^* $. 
For $\lambda\in\lie{h}^*$ we write $\lambda = \Re(\lambda)+i\Im(\lambda)$ where $\Re(\lambda)$, $\Im(\lambda)$ belong to $E$.  
Similarly, we decompose $\lambda\in\lie{h}_q^*$, under the understanding that $\Im(\lambda)$ is only defined modulo $\hbar^{-1}\roots^\vee$.

\begin{prop} 
\label{prop:root_subsystem}
Let $ \lambda \in \mathfrak{h}^*_q $. Then $ {\bf \Delta}_{[\lambda]} $ is a root system in 
its $ \mathbb{R} $-span $ E(\lambda) \subset E $, where $ E(\lambda) $ is equipped with the inner product induced from $ E $. 
Moreover $ W_{[\lambda]} $ is the Weyl group of $ {\bf \Delta}_{[\lambda]} $. 
\end{prop}

\begin{proof} In order to verify that $ {\bf \Delta}_{[\lambda]} $ is a root system in $ E(\lambda) $ it is enough to 
check $ s_\beta {\bf \Delta}_{[\lambda]} \subset {\bf \Delta}_{[\lambda]} $ for all $ \beta \in {\bf \Delta}_{[\lambda]} $. 
So let $ \alpha \in {\bf \Delta}_{[\lambda]} $ and observe that $ (s_\beta \alpha)^\vee = s_\beta \alpha^\vee $ since the action of $ s_\beta $ 
on $ E $ is isometric. Therefore 
$$
q_\alpha^{(\lambda, (s_\beta \alpha)^\vee)} = q_\alpha^{(\lambda, s_\beta \alpha^\vee)} 
= q_\alpha^{(\lambda, \alpha^\vee) - (\lambda, \beta^\vee) (\beta, \alpha^\vee)} \in \pm q_\alpha^{\mathbb{Z}}
$$
as desired, using that both $ \alpha $ and $ \beta $ are contained in $ {\bf \Delta}_{[\lambda]} $. Since $ q_\alpha = q_{s_\beta \alpha} $ 
this shows $ s_\beta \alpha \in {\bf \Delta}_{[\lambda]} $. 
	
In order to prove that $ W_{[\lambda]} $ is the Weyl group of $ {\bf \Delta}_{[\lambda]} $ we proceed in a similar way as in Section 3.4 in \cite{HumphreysO}.  
We show first that $ \alpha \in {\bf \Delta} $ satisfies $ \alpha \in {\bf \Delta}_{[\lambda]} $ iff $ s_\alpha \in W_{[\lambda]} $. 
Indeed, we have 
$ s_\alpha . \lambda - \lambda = -(\lambda + \rho, \alpha^\vee) \alpha $ in $ \mathfrak{h}^*_q $, and the condition 
$ q_\alpha^{(\lambda + \rho, \alpha^\vee)} \in \pm q_\alpha^{\mathbb{Z}} $ is equivalent to $ -(\lambda + \rho, \alpha^\vee) \alpha \in \roots $ 
in $ \mathfrak{h}^*_q/\frac{1}{2} i\hbar^{-1} \roots^\vee $. 
Thus, the Weyl group of $\Delta_{[\lambda]}$ is contained in $W_{[\lambda]}$, and to complete the proof it suffices to show that $W_{[\lambda]}$ is generated by the reflections it contains.
	
We first consider the case $ \lambda \in E $. We introduce the affine Weyl group $ W_a = \roots \rtimes W $, where $\roots$ acts upon $E$ by translations and $W$ 
by the standard Weyl group action.
Note that an element $w\in W$ belongs to $W_{[\lambda]}$ if and only if $(\beta , w)\lambda = \lambda $ for some $\beta\in\roots$. 
By Theorem 4.8 of \cite{HumphreysReflectiongroups} the  subgroup of $W_a$ which fixes any $\lambda\in E$ is generated by its reflections, and the result follows.
	
Now let $\lambda\in\lie{h}_q^* = \lie{h}^* / i\hbar^{-1}\roots^\vee$ be arbitrary, and lift it to an element $\lambda_0+i\lambda_1 \in \lie{h}^*$ 
with $\lambda_0,\lambda_1\in E$. An element $w\in W$ belongs to $W_{[\lambda]}$ if and only if
\[
w\lambda_0 - \lambda_0 \in \roots \qquad \text{ and } \qquad
w\lambda_1 - \lambda_1 \in \tfrac{1}{2}\hbar^{-1}\roots^\vee .
\]
By the previous case, the first condition is equivalent to $w\in W_{[\lambda_0]}$, which is the Weyl group of the root system ${\bf\Delta}_{[\lambda_0]}$ in the subspace $E(\lambda_0)$.  
	
Consider the rescaled coroot system
\[
{\bf\Delta}' = \{ \tfrac{1}{2}\hbar^{-1}\alpha^\vee \mid \alpha\in{\bf\Delta}_{[\lambda_0]} \} \subset E(\lambda_0).
\]
It is a root system with the same Weyl group $W_{[\lambda_0]}$. 
We write $\mathbf{Q}'$ for the lattice it generates and $W'_a = \mathbf{Q}'\rtimes W_{[\lambda_0]}$ for the associated affine Weyl group.  
Let us decompose $\lambda_1$ as 
\[
\lambda_1 = \lambda_1' + \lambda_1'' \in E(\lambda_0) \oplus E(\lambda_0)^\perp.
\]
Since $W_{[\lambda_0]}$  fixes $E(\lambda_0)^\perp$, we see that if $w\in W_{[\lambda_0]}$ then the above condition on $\lambda_1$ is equivalent to
\[
w \lambda_1 - \lambda_1 = w\lambda_1' - \lambda_1' \in \tfrac{1}{2}\hbar^{-1}\roots^\vee \cap E(\lambda_0) = \roots',
\]
and therefore $(\gamma , w)\lambda_1' = \lambda_1'$
for some $\gamma \in \roots'$.  As before, the subgroup of $W_a'$ fixing $\lambda_1'$ is generated by its reflections.  Thus $W_{[\lambda]}$ is generated by its reflections, and the result follows.
\end{proof} 

Let us now introduce the concept of dominant and antidominant weights. We recall that we use the 
convention $\mathbb{N} = \{1,2,\dots \} $ and $\mathbb{N}_0 = \mathbb{N} \cup \{0\} $.
	
\begin{definition} \label{defdominantantidominant}
Let $ \lambda \in \mathfrak{h}^*_q $. Then 
\begin{bnum} 
\item[a)] $ \lambda $ is dominant if 
$ q_\alpha^{(\lambda + \rho, \alpha^\vee)} \notin \pm q_\alpha^{- \mathbb{N}} $ 
for all $ \alpha \in {\bf \Delta}^+ $. 
\item[b)] $ \lambda $ is antidominant if 
$ q_\alpha^{(\lambda + \rho, \alpha^\vee)} \notin \pm q_\alpha^{\mathbb{N}} $
for all $ \alpha \in {\bf \Delta}^+ $. 
\end{bnum}
\end{definition} 

We remark that $ -\rho $ is both dominant and antidominant. Note that $ \lambda \in \mathfrak{h}^*_q $ is not dominant iff there exists 
some $ \alpha \in {\bf \Delta}^+ $ such that $ q_\alpha^{(\lambda + \rho, \alpha^\vee)} \in \pm q_\alpha^{-\mathbb{N}} $, 
or equivalently, such that $ (\lambda + \rho, \alpha^\vee) = -m + \frac{1}{2} i \hbar_\alpha^{-1} \mathbb{Z}$ for some $ m \in \mathbb{N} $, where we are using the notation $\hbar_\alpha = d_\alpha\hbar$. 
Similarly, $ \lambda \in \mathfrak{h}^*_q $ is not antidominant iff there exists 
some $ \alpha \in {\bf \Delta}^+ $ such that $ q_\alpha^{(\lambda + \rho, \alpha^\vee)} \in \pm q_\alpha^{-\mathbb{N}} $, 
or equivalently, such that $ (\lambda + \rho, \alpha^\vee) = m + \frac{1}{2} i \hbar_\alpha^{-1} \mathbb{Z} $ for some $ m \in \mathbb{N} $. 

We remark that the terminology introduced in Definition \ref{defdominantantidominant} is in tension with standard terminology used for dominant integral weights. 
More precisely, $ \lambda \in \weights $ is dominant in the sense of Definition \ref{defdominantantidominant} iff 
$ (\lambda + \rho, \alpha^\vee) \geq 0 $ for all $ \alpha \in {\bf \Delta}^+ $, whereas $ \lambda $ is dominant integral if $ (\lambda, \alpha^\vee) \geq 0 $ 
for all $ \alpha \in {\bf \Delta}^+ $. 
We shall refer to weights in $ \weights^+ $ as dominant integral weights, which means that being dominant and integral is not the same thing 
as being dominant integral. 

Let us introduce some further notation for reflections in the extended Weyl group $\hat{W} = \mathbf{Y}_q \rtimes W$, where we recall that $\mathbf{Y}_q = \tfrac{1}{2}i\hbar^{-1}\roots^\vee / i\hbar^{-1}\roots^\vee$. 
If $ k \in \mathbb{Z} $ and $ \alpha \in {\bf \Delta}_{[\lambda]} $, we write $ s_{k, \alpha} = (\frac{k}{2} i \hbar^{-1} \alpha^\vee, s_\alpha) \in \hat{W} $.
\nomenclature{$s_{k,\alpha}$}{affine reflection in the extended Weyl group $\hat{W}$}%
Explicitly, this element acts on $\nu\in\lie{h}_q^*$ by 
$$
s_{k, \alpha} . \nu = s_\alpha . \nu + \tfrac{k}{2} i \hbar^{-1} \alpha^\vee.
$$
Note that $ s_{k, \alpha} = s_{k + 2m, \alpha} $ for all $ m \in \mathbb{Z} $, so we shall 
frequently parametrize $ s_{k, \alpha} $ by $ k \in \mathbb{Z}_2 $.

In the sequel we will use the notation $ \lambda \leq \mu $ in $ \mathfrak{h}^*_q/\tfrac{1}{2} i \hbar^{-1} \roots^\vee $
\label{nom:leq_in_hqstar}%
to say that $ \mu - \lambda $ is congruent to an element of $ \roots^+ $ modulo $ \tfrac{1}{2} i \hbar^{-1} \roots^\vee $. In other words, this corresponds
to the natural order relation in the quotient $ \mathfrak{h}^*_q/\tfrac{1}{2} i \hbar^{-1} \roots^\vee $ of $ \mathfrak{h}^*_q $. 

\begin{prop} \label{antidominantchar}
Let $ \lambda \in \mathfrak{h}^*_q $. Then $ \lambda $ is antidominant if and only if one of the following equivalent conditions hold. 
\begin{bnum} 
\item[a)] $ (\lambda + \rho, \alpha) \leq 0 $ in $ \mathbb{C}/\frac{1}{2} i \hbar^{-1} \mathbb{Z} $ for all $ \alpha \in {\bf \Delta}^+_{[\lambda]} $, 
that is, the real part of $ (\lambda + \rho, \alpha) $ is non-positive, and the imaginary part of $ (\lambda + \rho, \alpha)  $ 
is contained in $ \frac{1}{2} \hbar^{-1} \mathbb{Z} $. 
\item[b)] $ \lambda \leq s_\alpha . \lambda $ in $ \mathfrak{h}^*_q/\tfrac{1}{2} i\hbar^{-1} \roots^\vee $ for all $ \alpha \in {\bf \Delta}^+_{[\lambda]} $. 
\item[c)] $ \lambda \leq w . \lambda $ in $ \mathfrak{h}^*_q/\tfrac{1}{2} i\hbar^{-1} \roots^\vee $ for all $ w \in W_{[\lambda]} $. 
\item[d)] If $ \hat{w} \in \hat{W}_{[\lambda]} $ satisfies $ \hat{w} . \lambda \leq \lambda $ in $ \mathfrak{h}^*_q $ then $ \hat{w} = e $. 
\end{bnum} 
\end{prop}

\begin{proof}  
Recall that $\alpha \in \Delta_{[\lambda]} $ if and only if $ q_\alpha^{(\lambda + \rho, \alpha^\vee)} \in \pm q_\alpha^{\mathbb{Z}}$.  Therefore $\lambda$ is antidominant if and only if for every $\alpha\in \Delta_{[\lambda]}^+ $ we have in fact $q_\alpha^{(\lambda+\rho,\alpha^\vee)} \in \pm q_\alpha^{-\mathbb{N}_0}$.  This is equivalent to condition $a)$.

Now we prove the equivalence of the four listed conditions.

$ a) \Leftrightarrow b) $ For any root $ \alpha \in {\bf \Delta}^+ $ we have 
$$
s_\alpha . \lambda = s_\alpha(\lambda + \rho) - \rho = \lambda - (\lambda + \rho, \alpha^\vee) \alpha = \lambda - (\lambda + \rho, \alpha) \alpha^\vee. 
$$
This shows the desired equivalence.

$ c) \Rightarrow b) $ is trivial. 

$ b) \Rightarrow c) $ We use induction on the length of $ w \in W_{[\lambda]} $. For $ w = e $ the claim obviously holds. 
Hence assume $ w = v s_i $ for some $ v \in W_{[\lambda]} $ with $ l(v) < l(w) $ and $ s_i $ the reflection associated to some simple root $\alpha_i$ in $\Delta_{[\lambda]}$. Note that, in this case, $ w \alpha < 0 $. We have 
$$
\lambda - w . \lambda = (\lambda - v . \lambda) + v  (\lambda - s_i . \lambda),
$$
where in the second term $v$ is acting by the unshifted action.
The first term satisfies $ \lambda - v . \lambda \leq 0 $ by our inductive hypothesis. For the second term we obtain
$$
v  (\lambda - s_i . \lambda) = w(s_i . \lambda - \lambda) = -(\lambda + \rho, \alpha_i^\vee) w \alpha_i \leq 0 
$$
as well.

$ c) \Rightarrow d) $ Assume $ \hat{w} = (\beta^\vee, w) $ satisfies $ \hat{w} . \lambda \leq \lambda $. We have 
$$
\hat{w} . \lambda - \lambda = w . \lambda - \lambda + \tfrac{1}{2} i \hbar^{-1} \beta^\vee, 
$$
so that $ w . \lambda \leq \lambda $ in $ \mathfrak{h}^*_q/\tfrac{1}{2} i \hbar^{-1} \roots^\vee $. Hence condition $ c) $ 
implies $ w = e $. Now $ (\beta^\vee, e) . \lambda = \lambda + \frac{1}{2} i \hbar^{-1} \beta^\vee \leq \lambda $ yields $ \beta^\vee = 0 $. 

$ d) \Rightarrow a) $ Fix $ \alpha \in {\bf \Delta}^+_{[\lambda]} $ and consider the affine reflection $ s_{k, \alpha} = (k \alpha^\vee, s_\alpha) $. 
Since $ s_{k, \alpha} \in \hat{W}_{[\lambda]} $ we automatically have 
$$ 
-(\lambda + \rho, \alpha^\vee) \alpha + \tfrac{k}{2} i \hbar^{-1} \alpha^\vee =  s_{k, \alpha} . \lambda - \lambda \in \roots 
$$ 
inside $ \mathfrak{h}^*_q/\tfrac{1}{2} i \hbar^{-1} \roots^\vee $, that is, 
$ (\lambda + \rho, \alpha^\vee) \alpha \in \roots \subset \mathfrak{h}^*_q/\tfrac{1}{2} i \hbar^{-1} \roots^\vee $. 
In particular, $ (\Im(\lambda) + \rho, \alpha) \in \frac{1}{2} i \hbar^{-1} \mathbb{Z} $. 
If $ (\Re(\lambda) + \rho, \alpha^\vee) > 0 $ then we get $ s_{k, \alpha} . \lambda < \lambda $ for suitable choice of $ k $, 
which is impossible by $ d) $. Hence $ (\Re(\lambda) + \rho, \alpha^\vee) \leq 0 $. 

This finishes the proof. \end{proof}  

An analogous result holds for dominant weights. 

\begin{prop} \label{dominantchar}
Let $ \lambda \in \mathfrak{h}^*_q $. Then $ \lambda $ is dominant if and only if one of the following equivalent conditions hold. 
\begin{bnum} 
\item[a)] $ (\lambda + \rho, \alpha) \geq 0 $ in $ \mathbb{C}/\frac{1}{2} i \hbar^{-1} \mathbb{Z} $ for all $ \alpha \in {\bf \Delta}^+_{[\lambda]} $, 
that is, the real part of $ (\lambda + \rho, \alpha) $ is non-negative, and the imaginary part of $ (\lambda + \rho, \alpha)  $ 
is contained in $ \frac{1}{2} \hbar^{-1} \mathbb{Z} $.
\item[b)] $ \lambda \geq s_\alpha . \lambda $ in $ \mathfrak{h}^*_q/\tfrac{1}{2} i\hbar^{-1} \roots^\vee $ for all $ \alpha \in {\bf \Delta}^+_{[\lambda]} $. 
\item[c)] $ \lambda \geq w . \lambda $ in $ \mathfrak{h}^*_q/\tfrac{1}{2} i\hbar^{-1} \roots^\vee $ for all $ w \in W_{[\lambda]} $. 
\item[d)] If $ \hat{w} \in \hat{W}_{[\lambda]} $ satisfies $ \hat{w} . \lambda \geq \lambda $ in $ \mathfrak{h}^*_q $ then $ \hat{w} = e $. 
\end{bnum} 
\end{prop}

\begin{proof}  This is a direct translation of the proof of Proposition \ref{antidominantchar}. \end{proof}

\subsection{Submodules of Verma modules} 

In this section we examine the existence of submodules of Verma modules. We shall discuss in particular an analogue of Verma's Theorem 
and a necessary and sufficient condition for the irreducibility of Verma modules. 

The socle $ \Soc(M) $
of a module $ M $ is the sum of all simple submodules of $ M $. Recall from Theorem \ref{uqgnoetheriandomain} 
that the algebra $ U_q(\mathfrak{n}_-) $ is Noetherian without zero-divisors. Since $ M(\lambda) $ for $ \lambda \in \mathfrak{h}^*_q $ 
is a free $ U_q(\mathfrak{n}_-) $-module, the following result follows from general facts, see Section 4.1 in \cite{HumphreysO}. 

\begin{lemma} \label{Vermasubmoduleintersection}
Let $ \lambda \in \mathfrak{h}^*_q $. Any two nonzero submodules of $ M(\lambda) $ have a nonzero intersection. 
\end{lemma} 

\begin{proof} 
Assume $ M, N \subset M(\lambda) $ are nonzero submodules. 
Then $ M \supset U_q(\mathfrak{n}_-) m \cdot v_\lambda $ and $ N \supset U_q(\mathfrak{n}_-) n \cdot v_\lambda $  for some nonzero 
elements $ m,n \in U_q(\mathfrak{n}_-) $. Since $ U_q(\mathfrak{n}_-) $ is Noetherian and has no zero divisors, 
the left ideals $ U_q(\mathfrak{n}_-) n $ and $ U_q(\mathfrak{n}_-) m $ must have nontrivial intersection, see Section 4.1 in \cite{HumphreysO}. 
\end{proof}  

We recall that $V(\mu)$ denotes the unique simple quotient of the Verma module $M(\mu)$.

\begin{lemma} \label{Vermasocle}
Let $ \lambda \in \mathfrak{h}^*_q $. Then the socle of $ M(\lambda) $ is of the form 
$$ 
\Soc(M(\lambda)) \cong V(\mu) \cong M(\mu) 
$$ 
for a unique weight $ \mu \leq \lambda $. 
\end{lemma} 

\begin{proof} 
According to Lemma \ref{Vermasubmoduleintersection}, the socle of $ M(\lambda) $ must be simple, therefore of the form $ V(\mu) $ for 
some $ \mu \in \mathfrak{h}^*_q $. From the weight space structure of $ M(\lambda) $ it is immediate that $ \mu \leq \lambda $. 
Moreover $ V(\mu) \subset M(\lambda) $ is necessarily a Verma module since $ M(\lambda) $ is free as a $ U_q(\mathfrak{n}_-) $-module. 
Hence $ V(\mu) \cong M(\mu) $. 
\end{proof}  

Let us remark that since $ M(\lambda) $ has finite length it is not hard to check that $ \Soc(M(\lambda)) = M(\mu) \subset V $ for any 
submodule $ V \subset M(\lambda) $. 

We shall now study more general Verma submodules of $ M(\lambda) $. Recall the affine reflections $s_{k, \alpha_i}$ defined before Proposition \ref{antidominantchar}.

\begin{lemma} \label{Vermahelpembedding}
Let $ \lambda \in \mathfrak{h}^*_q $ and assume $ s_{k, \alpha_i} . \lambda \leq \lambda $ for some $ 1 \leq i \leq N $ and $ k \in \mathbb{Z}_2 $. 
Then there exists an embedding of $ M(s_{k, \alpha_i} . \lambda) $ into $ M(\lambda) $. 
\end{lemma} 

\begin{proof}  
Note that 
\[
\lambda - s_{k,\alpha_i}.\lambda = (\lambda+\rho,\alpha_i^\vee)\alpha_i - \tfrac{k}{2} i \hbar^{-1} \alpha_i^\vee, 
\]
so the assumption $ s_{k, \alpha_i} . \lambda \leq \lambda $ 
implies
\begin{align*}
&\Re(\lambda+\rho,\alpha_i^\vee) \in \mathbb{N}_0, \\
&\Im(\lambda,\alpha_i^\vee) = \Im(\lambda+\rho,\alpha_i^\vee) \equiv \tfrac{k}{2}\hbar^{-1}d_i^{-1}  \pmod{\hbar^{-1}d_i^{-1}\mathbb{Z}}.
\end{align*}
Put $m = \Re(\lambda+\rho,\alpha_i^\vee) = \Re(\lambda,\alpha_i^\vee)+1$.  We claim that the vector $F_i^m v_\lambda$ is primitive.  Indeed, for $j\neq i$ we have 
\[
E_j F_i^m v_\lambda = F_i^m E_j v_\lambda = 0,
\]
while for $j=i$, Lemma \ref{EFbasic} gives
\begin{align*}
E_i F_i^m \cdot v_\lambda 
&= F_i^m E_i v_\lambda + [m]_q F_i^{m-1} \frac{q_i^{-m+1}K_i - q_i^{m-1}K_i^{-1}}{q_i-q_i^{-1}}\cdot v_\lambda \\
&= [m]_q \frac{q^{-d_i\Re(\lambda,\alpha_i^\vee)}q^{(\lambda,\alpha_i)} 
- q^{d_i\Re(\lambda,\alpha_i^\vee)}q^{-(\lambda,\alpha_i)}}{q_i-q_i^{-1}} F_i^{m-1}\cdot v_\lambda \\
&= [m]_q \frac{q^{-i\Im(\lambda,\alpha_i)} - q^{i\Im(\lambda,\alpha_i)}}{q_i-q_i^{-1}} F_i^{m-1} \cdot v_\lambda \\
&= 0.
\end{align*}
The weight of $ F_i^m \cdot v_\lambda $ is 
\begin{align*}
\lambda - \Re(\lambda+\rho,\alpha_i^\vee)\alpha_i &= s_i.\lambda + \Im(\lambda+\rho,\alpha_i^\vee)\alpha_i = s_{k, \alpha_i} . \lambda  \in \mathfrak{h}_q^*.
\end{align*}
Hence we obtain an embedding as desired. 
\end{proof}

Let us recall some definitions related to the Weyl group $ W $, see Appendix A.1 in \cite{Josephbook}. For $ w \in W $ set 
$$
S(w) = \{\alpha \in {\bf \Delta}^+ \mid w \alpha \in {\bf \Delta}^- \} = {\bf \Delta}^+ \cap w^{-1} {\bf \Delta}^-. 
$$
\nomenclature{$S(w)$}{$= {\bf \Delta}^+ \cap w^{-1} {\bf \Delta}^-$ for a Weyl group element $w$}%
If $ u, v \in W $ and $ \alpha \in {\bf \Delta}^+ $ we write $ u \overset{\alpha}{\leftarrow} v $ or simply $ u \leftarrow v $
\nomenclature[o$<$]{$\overset{\alpha}{\leftarrow}$, $\leftarrow$}{adjacency in the Bruhat order}%
if $ u = s_\alpha v $ and $ l(u) = l(v) - 1 $. In this case one has $ \alpha \in S(v^{-1}) $. 
The Bruhat order on $ W $ is defined by saying that $ u \leq v $
\nomenclature[o$<$]{$\leq$, $<$}{Bruhat order on the Weyl group}%
if $ u = v $ or if there 
is a chain $ u = w_1 \leftarrow \cdots \leftarrow w_n = v $. As usual we write $ u < v $ if $ u \leq v $ and $ u \neq v $.

\begin{lemma} \label{Vermahelpintegral}
Let $ \mu \in \weights^+ $ and assume $ u, v \in W $ are such that $ u \leq v $ for the Bruhat order. Then there exists an embedding 
of $ M(v . \mu) $ into $ M(u . \mu) $. 
\end{lemma} 

\begin{proof}  We use induction on the length of $ v $. The case $ l(v) = 0 $, corresponding to $ v = e $ is trivial. For the inductive step assume 
$ s_i v < v $ for some $ i $. Then 
\[
s_i.(v.\mu) - v.\mu
= s_iv(\mu+\rho) - v(\mu+\rho)
= -(v(\mu+\rho),\alpha_i^\vee)
= -(\mu+\rho,v^{-1}\alpha_i^\vee) > 0,
\]
which means $ v . \mu \leq s_i . (v . \mu) $. Hence $ M(v. \mu) $ is a submodule of $ M(s_i v . \mu) $ 
according to Lemma \ref{Vermahelpembedding}. 

Properties of the Bruhat order imply that we have either $ u \leq s_i v $ or $ s_i u \leq s_i v $, see for instance Proposition A.1.7 in \cite{Josephbook}. 

In the first case we obtain an embedding $ M(s_i v . \mu) \subset M(u. \mu) $ 
by our inductive hypothesis since $ s_i v $ has length $ l(v) - 1 $. Combining this with the embedding $ M(v . \mu) \subset M(s_i v . \mu) $ 
obtained above yields the claim in this situation. 

In the second case we may assume $ s_i u < u $, since otherwise we are again in the previous case.  
Consider the embedding $ \pi: M(v . \mu) \rightarrow M(s_i v . \mu) \rightarrow M(s_i u. \mu) $ 
obtained from the inductive hypothesis.
Let $ x $ be a highest weight vector for $M(v . \mu)$. 
Note that $ K_i \cdot x = q^{(\alpha_i, v . \mu)} x $ with 
\[ 
(\alpha_i, v . \mu) = (v^{-1}(\alpha_i),\mu+\rho) - (\alpha_i,\rho) < 0 .
\]
Therefore, repeated application of Lemma \ref{EFbasic} shows
that $ E_i^r F_i^r \cdot x$ is a nonzero multiple of $ x $ 
for all $ r \in \mathbb{N}_0 $.  On the other hand, let $y$ be a highest weight vector for $ M(s_i u .\mu) $, and note from the proof of 
Lemma \ref{Vermahelpembedding} that a highest weight vector of the submodule $M(u.\mu) \subset M(s_iu.\mu)$ is given by $F_i^my$ for some $m\in\NN$. 
Since $\pi(x)\in U_q(\mathfrak{n}_-)\cdot y$, repeated application of the Serre relations shows that for $r\gg 0$ we 
have $F_i^r\pi(x) \in U_q(\mathfrak{n}_-)F_i^m\cdot y = M(u.\mu)$.  Therefore, $\pi(E_i^rF_i^r\cdot x) \in M(u.\mu)$ and hence $\pi$ maps $M(v.\mu)$ 
into $M(u.\mu)$, as desired.
\end{proof}  

\begin{definition}
	\label{def:strongly_linked}
Let $ \mu, \lambda \in \mathfrak{h}^*_q $. We say that $ \mu \uparrow \lambda $
\nomenclature[o$<$]{$\uparrow$}{strong linkage in $\lie{h}_q^*$}%
if $ \mu = \lambda $ or if there exists a chain of positive roots 
$ \alpha_1, \dots, \alpha_r \in {\bf \Delta}^+ $ and $ k_1, \dots, k_r \in \mathbb{Z} $ such that 
$$
\mu = s_{k_1, \alpha_1} \cdots s_{k_r, \alpha_r} . \lambda < s_{k_2, \alpha_2} \cdots s_{k_r, \alpha_r} . \lambda < \cdots < s_{k_r, \alpha_r} . \lambda < \lambda.  
$$ 
We say that $ \mu $ is \emph{strongly linked} to $ \lambda $ if $ \mu \uparrow \lambda $.  
\end{definition}

Note that $ s_{k, \alpha} . \nu < \nu $ for some $k\in\ZZ$ if and only if
\[
q_\alpha^{(\nu + \rho, \alpha^\vee)} \in \pm q_\alpha^{\mathbb{N}}.
\]

We are now ready to prove the following analogue of Verma's Theorem, compare section 4.4.9 in \cite{Josephbook} and Section 4.7 in \cite{HumphreysO}. 

\begin{theorem} \label{Vermatheorem}
Let $ \mu, \lambda \in \mathfrak{h}^*_q $. 
If $ \mu $ is strongly linked to $ \lambda $ then $ M(\mu) \subset M(\lambda) $, in particular $ [M(\lambda): V(\mu)] \neq 0 $. 
\end{theorem}

\begin{proof} 
In order to prove the claim it is enough to consider the case $ \mu = s_{k, \alpha} . \lambda $ 
for some $ \alpha \in {\bf \Delta}^+ $ with $ q_\alpha^{(\lambda + \rho, \alpha^\vee)} \in \pm q_\alpha^{\mathbb{N}} $. 

For $ \alpha \in {\bf \Delta}^+ $ and $ n \in \mathbb{N} $ consider the sets 
$$
\Lambda_{n, \alpha} = \{\lambda \in \mathfrak{h}_q^* \mid q_\alpha^{(\lambda + \rho, \alpha^\vee)} = \pm q_\alpha^n \}
$$
and 
$$
X_{n, \alpha} = \{\lambda \in \mathfrak{h}^*_q \mid \Hom_{U_q(\mathfrak{g})}(M(\lambda - n \alpha), M(\lambda)) \neq 0 \}. 
$$
Note that $ X_{n, \alpha} \subset \Lambda_{n, \alpha} $. 
Our aim is to show $ X_{n, \alpha} = \Lambda_{n, \alpha} $.

We begin by showing that $ \weights \cap \Lambda_{n, \alpha} \subset X_{n, \alpha} $. 
Let $ \lambda \in \weights \cap \Lambda_{n, \alpha} $, and choose $ w \in W $ such that $ \mu = w^{-1} s_\alpha . \lambda \in \weights^+ $. 
Then also $ w^{-1} s_\alpha . \lambda + \rho \in \weights^+ $ and 
$$ 
(w^{-1} s_\alpha . \lambda + \rho, w^{-1} \alpha^\vee) = (w^{-1} s_\alpha(\lambda + \rho), w^{-1} \alpha^\vee) 
= -(\lambda + \rho, \alpha^\vee) = - n, 
$$
which means $ w^{-1} \alpha \in {\bf \Delta}^- $. This implies $ w > s_\alpha w $, 
so we obtain $ M(s_\alpha . \lambda) = M(w . \mu) \subset M(s_\alpha w . \mu) = M(\lambda) $ by Lemma \ref{Vermahelpintegral}. Hence 
$ \lambda \in X_{n, \alpha} $, and we conclude $ \weights \cap \Lambda_{n, \alpha} \subset X_{n, \alpha} $. 

To complete the proof, we use a Zariski density argument.
If $ B = \mathbb{C}[x_1, \dots, x_d]/I $ is a commutative algebra then the algebraic variety $ \Spm(B) $ 
is the set of all points $ x = (x_1, \dots, x_d) \in \mathbb{C}^d $ such that $ f(x_1, \dots, x_d) = 0 $ for all $ f \in I $. 
Observe that $ \mathfrak{h}^*_q $ can be identified with the variety $ \Spm(B) $ of $ B = U_q(\mathfrak{h}) $ 
such that $ \lambda \in \mathfrak{h}^*_q $ corresponds to the point of $ \Spm(B) $ given by the character $ \chi_\lambda $ on $ B $. 
Clearly, $ \Lambda_{n,\alpha} $ is a Zariski closed subspace.

We claim that $ X_{n, \alpha} \subset \mathfrak{h}^*_q $ is also Zariski closed.  
Indeed, for each $ \gamma \in \roots^+ $ fix a basis $ Y_1^\gamma, \dots, Y_r^\gamma $ 
of $ U_q(\mathfrak{n}_-)_{-\gamma} $. From the defining relations of $ U_q(\mathfrak{g}) $ it is clear that we can find Laurent polynomials $ p^k_{ij} $ 
in the generators $ K_1, \dots, K_N $ 
such that 
$$ 
[E_k, Y^\gamma_i] = \sum_j Y^{\gamma - \alpha_k}_j p^k_{ij}.  
$$ 
Now consider $ Y \in U_q(\mathfrak{n}_-)_{-\gamma} $ and write $ Y = \sum c_i Y^\gamma_i $. 
Then we have 
$$
E_k Y \cdot v_\lambda = [E_k, Y] \cdot v_\lambda = \sum_{i,j} c_i Y^{\gamma - \alpha_k}_j p^k_{ij} \cdot v_\lambda 
= \sum_{i,j} c_i \chi_\lambda(p^k_{ij}) Y^{\gamma - \alpha_k}_j \cdot v_\lambda
$$
Since $ M(\lambda) $ is a free $ U_q(\mathfrak{n}_-) $-module it follows that $ Y \cdot v_\lambda $ 
is the highest weight vector for a submodule $ M(\lambda - \gamma) \subset M(\lambda) $ iff $ \sum_i c_i \chi_\lambda(p^k_{ij}) = 0 $ for all $ k,j $. 
Thus $ \lambda\in X_{n, \alpha} $ if and only if there exists a nontrivial kernel of the linear map
\[
(c_i)_i \mapsto \left( \sum_i c_i \chi_\lambda(p^k_{ij}) \right)_{j,k},
\]
whose coefficients are polynomial in the character $\chi_\lambda$. The existence of such a kernel is determined by the vanishing of the determinants of minors of 
the corresponding matrix $ (\chi_\lambda(p^k_{ij}))_{(j,k), i} $. 
Applying this to $ \gamma = n \alpha $ we conclude that $ X_{n, \alpha} $ is an algebraic subset of $ \Spm(B) $. 

According to our above considerations $ X_{n, \alpha} \cap \weights = \Lambda_{n, \alpha} \cap \weights $. 
Since $ \weights \cap \Lambda_{n, \alpha} \subset \Lambda_{n, \alpha} $ is Zariski dense 
we conclude $ X_{n, \alpha} = \Lambda_{n, \alpha} $ as desired. \end{proof}  

As a consequence of Theorem \ref{Vermatheorem} we obtain the following characterization of simple Verma modules, 
compare Section 17.4 in \cite{deConciniProcesiqg}.

\begin{theorem} \label{Vermacomplexirreducible}
Let $ \lambda \in \mathfrak{h}^*_q $. The Verma module $ M(\lambda) $ is simple iff $ \lambda $ is antidominant. 
\end{theorem} 

\begin{proof}  If $ \lambda $ is not antidominant there exists an $ \alpha \in {\bf \Delta}^+ $ 
such that $ (\lambda + \rho, \alpha^\vee) \in m + \frac{1}{2}  i \hbar_\alpha^{-1} \mathbb{Z} $ for some $ m \in \mathbb{N} $, where as before we use the notation $\hbar_\alpha = d_\alpha\hbar$. In this case 
$ s_{k, \alpha} . \lambda < \lambda $ for suitable $ k $, and according to Theorem \ref{Vermatheorem} there exists an 
embedding $ M(s_{k, \alpha} . \lambda) \hookrightarrow M(\lambda) $. In particular, $ M(\lambda) $ is not irreducible. 

Conversely, assume that $ \lambda $ is antidominant. Then due to Proposition \ref{antidominantchar} the weight $ \lambda $ is minimal in 
its $ \hat{W}_{[\lambda]} $-orbit. By Harish-Chandra's Theorem \ref{chilinkage}, all irreducible subquotients of $ M(\lambda) $ must be of the 
form $ V(\mu) $ for some $ \mu \leq \lambda $ in the $ \hat{W}_{[\lambda]} $-orbit of $ \lambda $. 
We conclude $ \mu = \lambda $, and hence $ M(\lambda) $ is irreducible. \end{proof}  

Let us also characterize projective Verma modules, compare 
Section 8.2 in \cite{Josephbook}. 

\begin{prop} \label{Vermaprojective}
Let $ \lambda \in \mathfrak{h}^*_q $. Then the following conditions are equivalent. 
\begin{bnum} 
\item[a)] $ M(\lambda) $ is a projective object in $ \O $.
\item[b)] If $ \Hom_{U_q(\mathfrak{g})}(M(\lambda), M(\mu)) \neq 0 $ for $ \mu \geq \lambda $ then $ \lambda = \mu $. 
\item[c)] $ \lambda $ is dominant. 
\end{bnum}
\end{prop}

\begin{proof} 
$ a) \Rightarrow b) $ Assume $ f: M(\lambda) \rightarrow M(\mu) $ is a nonzero homomorphism and $ \mu > \lambda $.  
Note that $ f $ is injective and set $ M = M(\mu)/f(I(\lambda)) $, where $I(\lambda)$ is the maximal proper submodule in $M(\lambda)$. Then $ f $ induces 
an embedding $ V(\lambda) \subset M $, and hence a surjection $ M^\vee \rightarrow V(\lambda) $. 
Since $ M(\lambda) $ is projective we get a nonzero map $ M(\lambda) \rightarrow M^\vee $ 
lifting the canonical projection $ M(\lambda) \rightarrow V(\lambda) $, 
or equivalently, a nonzero map $ M \rightarrow M(\lambda)^\vee $. This induces a nonzero map $ M(\mu) \rightarrow M(\lambda)^\vee $.  This is impossible, 
since the $ \mu $-weight space of $ M(\lambda)^\vee $ is trivial.

$ b) \Rightarrow c) $ Assume $ \lambda $ is not dominant. Then according to Proposition \ref{dominantchar} 
there exists some $ \mu > \lambda $ such that $ \lambda \uparrow \mu $, 
and hence an embedding $ M(\lambda) \rightarrow M(\mu) $ by Theorem \ref{Vermatheorem}. This contradicts condition $ b) $. 

$ c) \Rightarrow a) $ Assume that $ \pi: M \rightarrow N $ is a surjective morphism in $ \O $ and let $ f: M(\lambda) \rightarrow N $ be given. 
We may assume that $ f $ is nonzero. 
Moreover, since the centre acts by $ \xi_\lambda $ on $ M(\lambda) $ the same is true for the action of $ ZU_q(\mathfrak{g}) $ on $ f(M(\lambda)) $. 
Using primary decomposition we may assume that $ Z - \xi_\lambda(Z) $ acts nilpotently on $ M $ and $ N $ for all $ Z \in ZU_q(\mathfrak{g}) $. 

Let $ v \in M $ be a vector of weight $ \lambda $ such that $ \pi(v) = f(v_\lambda) $ and let $ V = U_q(\mathfrak{n}_+) \cdot v \subset M $ be 
the $ U_q(\mathfrak{n}_+) $-submodule of $ M $ generated by $ v $. Since $ M $ is in $ \O $ the space $ V $ is finite dimensional. Choose a primitive 
vector $ v_\mu \in V $, for some weight $ \mu \geq \lambda $. Then $ \xi_\lambda = \xi_\mu $, from which we 
conclude that $ \mu $ and $ \lambda $ are $ \hat{W} $-linked. Due to dominance of $ \lambda $ this means $ \mu = \lambda $. 
Hence $ v $ itself is a primitive vector, which means that there exists a unique homomorphism $ F: M(\lambda) \rightarrow M $ 
such that $ F(v_\lambda) = v $. Clearly $ F $ is a lift for $ f $ as desired. 
\end{proof}

\subsection{The Shapovalov determinant} 

In this section we study the Shapovalov form for $ U_q(\mathfrak{g}) $. This form and its determinant are important tools in 
the study of Verma modules. We refer to \cite{HeckenbergerYamaneshapovalov} for the analysis of Shapovalov forms in greater generality. 

Let us recall from Lemma \ref{deftau} that the involutive algebra anti-automorphism $ \tau $ of $ U_q(\mathfrak{g}) $ keeps $ U_q(\mathfrak{h}) $ 
pointwise fixed and interchanges $ U_q(\mathfrak{b}_+) $ and $ U_q(\mathfrak{b}_-) $. 
Recall also that the Harish-Chandra map is the linear map $ \P: U_q(\mathfrak{g}) \rightarrow U_q(\mathfrak{h}) $ 
given by $ \hat{\epsilon} \otimes \id \otimes \hat{\epsilon} $ with respect to the triangular decomposition 
$ U_q(\mathfrak{g}) \cong U_q(\mathfrak{n}_-) \otimes U_q(\mathfrak{h}) \otimes U_q(\mathfrak{n}_+) $, see Proposition \ref{uqgtriangular}. 

We shall call the bilinear form $ \Sh: U_q(\mathfrak{n}_-) \times U_q(\mathfrak{n}_-) \rightarrow U_q(\mathfrak{h}) $ defined by 
$$
\Sh(X,Y) = \P(\tau(X) Y)
$$
\nomenclature[o$Sh$]{$\Sh$}{Shapovalov form for $U_q(\lie{g})$}%
the Shapovalov form for $ U_q(\mathfrak{g}) $. 
One checks that $ \P \tau = \P $, and using $ \tau^2 = \id $ this implies immediately that $ \Sh $ is symmetric. 
Moreover we have $ \Sh(X, Y) = 0 $ if $ X, Y $ are of different weight. 

We denote the restriction of $ \Sh $ to the weight space $ U_q(\mathfrak{n}_-)_{-\mu} $ for  $\mu \in \roots^+$ 
by $ \Sh_\mu $. 
In the following, we will use the determinant of the form $\Sh_\mu$ as follows.  Fix a basis $ x_1, \dots, x_{n_\mu} $ of $ U_q(\mathfrak{n}_-)_{-\mu} $, and denote by $\det(\Sh_\mu)$
\nomenclature[o$det(Sh_\mu)$]{$\det(\Sh_\mu)$}{$U_q(\lie{h})$-valued determinant of the Shapovalov form $\Sh_\mu$}%
the determinant of the matrix $ (\Sh_\mu(x_i, x_j)) \in M_{n_\mu}(U_q(\mathfrak{h})) $.  This only depends on the choice of basis up to multiplication by a nonzero scalar in $\mathbb{C}$, which we shall disregard in the following discussion.
Moreover,by inspecting the definition of $\tau$ and the commutation relations for $E_i$ and $F_j$, one sees that the Shapovalov determinant $ \det(\Sh_\mu) $ is in fact contained in the subalgebra of $ U_q(\mathfrak{h}) $ generated 
by the elements $ K_j^{\pm 1} $ for $ j = 1, \dots, N $. 

Let $ \lambda \in \mathfrak{h}^*_q $. 
The contravariant bilinear form on $ M(\lambda) $ is defined as the unique bilinear form $ \Sh^\lambda: M(\lambda) \times M(\lambda) \rightarrow \mathbb{C} $ satisfying $ \Sh^\lambda(v_\lambda, v_\lambda) = 1 $ and 
$$ 
\Sh^\lambda(\tau(X) \cdot v, w) = \Sh^\lambda(v, X \cdot w) 
$$ 
\nomenclature[o$Sh^\lambda$]{$\Sh^\lambda$}{contravariant bilinear form on the Verma module $M(\lambda)$}%
for all $ v,w \in M(\lambda) $ and $ X \in U_q(\mathfrak{g}) $. 
Different weight spaces of $ M(\lambda) $ are orthogonal with respect to the form, that is, we have $ \Sh^\lambda(M(\lambda)_\mu, M(\lambda)_\nu) = 0 $ 
for $ \mu \neq \nu $. Notice that 
\[
\Sh^\lambda(X\cdot v_\lambda,Y\cdot v_\lambda) = \chi_\lambda(\Sh(X,Y))
\]
for all $ X, Y \in U_q(\mathfrak{n}_-) $, 
where we recall that $ \chi_\lambda(K_\mu) = q^{(\lambda, \mu)} $. 

The radical of $ \Sh^\lambda $ is the subspace $ R^\lambda \subset M(\lambda) $
consisting of all $ u $ satisfying 
$ \Sh^\lambda(u, v) = 0 $ for all $ v \in M(\lambda) $. 

\begin{lemma} \label{radicalmaximalsubmodule}
Let $ \lambda \in \mathfrak{h}^*_q $. The radical of $ \Sh^\lambda $ agrees with the maximal proper submodule $ I(\lambda) $ 
of $ M(\lambda) $. 
\end{lemma} 

\begin{proof}  From the invariance property of $ \Sh^\lambda $ it is clear that the radical $ R^\lambda $ is a submodule of $ M(\lambda) $. 
It is also clear that the highest weight vector $ v_\lambda $ is not contained in $ R^\lambda $, and hence $ R^\lambda \subset I(\lambda) $. 
Assume that the quotient $ M(\lambda)/R^\lambda $ is not simple. Then $ M(\lambda)/R^\lambda $ must contain a primitive vector, 
that is, a nonzero vector $ v $ of weight $ \lambda - \mu $ 
for some $ \mu \in \roots^+ \setminus \{0\} $ which is annihilated by $ U_q(\mathfrak{n}_+) \cap \ker(\hat{\epsilon}) $. This implies 
$$
\Sh^\lambda(v, U_q(\mathfrak{n}_-) \cdot v_\lambda) = \Sh^\lambda(U_q(\mathfrak{n}_+) \cdot v, v_\lambda) = 0 
$$
and hence $ v = 0 $ by nondegeneracy of the induced form on $ M(\lambda)/R^\lambda $. This is a contradiction, so that $ M(\lambda)/R^\lambda $ must 
be simple. Hence we obtain $ R^\lambda = I(\lambda) $ as claimed. \end{proof}  

Instead of working with the form $ \Sh(X,Y) $ it is sometimes technically more convenient to consider the modified form 
$ \S: U_q(\mathfrak{n}_-) \times U_q(\mathfrak{n}_-) \rightarrow U_q(\mathfrak{h}) $ given by 
$$ 
\S(X,Y) = \P(\Omega(X) Y),  
$$ 
where $ \Omega $ is the algebra anti-homomorphism from Lemma \ref{defOmega}. As before, if $ \nu \in \roots^+ $ we write $ \S_\nu $ 
for the 
restriction of $ \S $ to $ U_q(\mathfrak{n}_-)_{-\nu} $. The determinant of $ \Sh_\nu $ differs from the one of $ \S_\nu $ 
only by an invertible element of $ U_q(\mathfrak{h}) $.

\begin{lemma} \label{shapovalovorderestimate}
Let $ \nu \in \roots^+ $ and $ X, Y \in U_q(\mathfrak{n}_-)_{-\nu} $. Then 
$$ 
\S_\nu(X, Y) \in \sum_{\substack{\beta, \gamma \in \roots^+\\ \beta + \gamma = \nu}} \mathbb{C} K_{\beta - \gamma}.  
$$ 
\end{lemma} 

\begin{proof}  We may assume without loss of generality that both $ X $ and $ Y $ are monomials in the generators $ F_1, \dots, F_N $. 
Let us proceed by induction on the height $ht(\nu) =  \nu_1 + \cdots + \nu_N $ of $ \nu = \nu_1 \alpha_1 + \cdots + \nu_N \alpha_N $. 
Note that for $ ht(\nu) = 0 $ the claim is obvious. 

Fix $m\geq0$ and assume that the claim holds for all terms of height $ m  $. Let $X,Y \in U_q(\lie{n}_-)_{-\nu}$ be monomials in  $F_1,\ldots,F_N$ with $ht(\nu)=m+1$, and write $ X = F_i Z $ for suitable $ i $ and $ Z \in U_q(\mathfrak{n}_-)_{-\nu + \alpha_i} $. 
Let $ r $ be the total number of factors $ F_i $ appearing in the monomial $ Y $. If $ 1 \leq k \leq r $ then for the $ k $-th occurrence of $ F_i $ in $ Y $ 
we decompose $ Y = Y^{(a)}_k F_i Y^{(b)}_k $, where $ Y^{(a)}_k Y^{(b)}_k $ are monomials in the generators $ F_1, \dots, F_N $ 
containing a total number of factors $ F_i $ one less than for $ Y $. We obtain 
\begin{align*}
\S_\nu(X, Y) = \sum_{k = 1}^r \P(\Omega(Z) Y^{(a)}_k [E_i, F_i] Y^{(b)}_k) &\in 
\sum_{k = 1}^r (\mathbb{C} K_i + \mathbb{C} K_i^{-1}) \P(\Omega(Z) Y^{(a)}_k Y^{(b)}_k). 
\end{align*} 
The latter is contained in 
\begin{align*}
(\mathbb{C} K_i + \mathbb{C} K_i^{-1}) \sum_{\beta, \gamma \in \roots^+, \beta + \gamma = \nu - \alpha_i} \mathbb{C} K_{\beta - \gamma} 
= \sum_{\beta, \gamma \in \roots^+, \beta + \gamma = \nu} \mathbb{C} K_{\beta - \gamma}
\end{align*}
according to our inductive hypothesis. This yields the claim. \end{proof}  

Let us enumerate the positive roots of $\mathfrak{g}$ again as $\beta_1,\ldots,\beta_n$. Recall that the Kostant partition 
function $ P: \roots \to \NN_0 $ is given by
\[
 P(\nu) = \left| \{ (r_1,\ldots,r_n)\in\NN_0^n \mid r_1\beta_1 + \cdots + r_n\beta_n = \nu \}\right|.
\]

\begin{lemma} \label{shapmultiplicity} 
We have an equality of weights
$$
\dim(U_q(\mathfrak{n}_-)_{-\nu}) \nu = \sum_{j = 1}^n \sum_{m = 1}^\infty P(\nu - m\beta_j) \beta_j
$$
for any $ \nu \in \roots^+ $.  
\end{lemma}

\begin{proof} 
According to the PBW-Theorem \ref{PBWgeneralfield},  
the space $ U_q(\mathfrak{n}_-)_{-\nu} $ is spanned by all PBW-vectors $ F_{\beta_1}^{r_1} \cdots F_{\beta_n}^{r_n} $ 
such that $ r_1 \beta_1 + \cdots + r_n \beta_n = \nu $. Therefore,
\[
\dim(U_q(\mathfrak{n}_-)_{-\nu}) \nu 
= \hspace{-1cm} \sum_{{(r_1,\ldots,r_n)\in\NN_0^n \atop \text{with } r_1 \beta_1 +\cdots + r_n\beta_n =\nu}}
\hspace{-1cm} r_1 \beta_1 +\cdots + r_n\beta_n.
\]
Separating the coefficients of $\beta_j$ in this sum gives
\begin{align*}
\sum_{j=1}^n \sum_{m=1}^\infty & \left| \{ (r_1,\ldots,r_n)\in\NN_0^n \mid r_1 \beta_1 +\cdots + r_n\beta_n =\nu  \text{ with } r_j=m\}\right| \, m\beta_j \\
&= \sum_{j=1}^n \sum_{m=1}^\infty \left| \{ (r_1,\ldots,r_n)\in\NN_0^n \mid r_1 \beta_1 +\cdots + r_n\beta_n =\nu \text{ with } r_j\geq m\}\right| \, \beta_j, \\
&= \sum_{j=1}^n \sum_{m=1}^\infty \left| \{ (r_1,\ldots,r_n)\in\NN_0^n \mid r_1 \beta_1 +\cdots + r_n\beta_n =\nu - m\beta_j \}\right| \, \beta_j, \\
&= \sum_{j=1}^n \sum_{m=1}^\infty P(\nu-m\beta_j) \beta_j,
\end{align*}
as claimed.
\end{proof}

We will also need a result from commutative algebra. Let $ B \cong \mathbb{C}[x_1, \dots, x_d]/I $ 
be a finitely generated commutative algebra over $ \mathbb{C} $. 
Recall that the algebraic variety $ \Spm(B) $
of $ B $ is the set of all points $ p = (p_1, \dots, p_d) $ in $ \mathbb{C}^d $ such that $ f(p_1, \dots, p_d) = 0 $ 
for all $ f \in I $. The elements of $ B $ can be viewed as polynomial functions on $ \Spm(B) $, and we write $ b_p $ for the value 
of $ b \in B $ at $ p \in \Spm(B) $.

If $ y(x) \in B[x] $ is a polynomial we write $ y(0) \in B $ for its image under evaluation at $ 0 $, and 
we use the same notation for matrices with entries in $ B[x] $. 

The following two Lemmas are taken from \cite{HeckenbergerYamaneshapovalov}. 

\begin{lemma} \label{shapaghelp} 
Let $ B $ be an integral domain and $ Y \in M_n(B[x]) $ for some $ n \in \mathbb{N} $. Then 
there exist $ 0 \leq k \leq n $, matrices $ M_1, M_2 \in M_n(B) $ with nonzero determinants, a matrix $ M \in M_n(B[x]) $ and a nonzero element $ b \in B $ 
such that 
$$
M_1 Y M_2 = x M + b D_k
$$
where $D_k = \mathrm{diag}(1,\ldots,1,0,\ldots,0)$ is the diagonal matrix projecting onto the first $k$ coordinates. 
\end{lemma} 

\begin{proof} 
Let $ \Quot(B) $ be the field of fractions of $ B $. By linear algebra there exist invertible matrices $ C_1, C_2 \in M_n(\Quot(B)) $ such that 
$$
C_1 Y(0) C_2 = D_k 
$$
for some $ 0 \leq k \leq n $. Let $ b_1, b_2 $ be nonzero elements such that $ M_1 = b_1 C_1 $ and $ M_2 = b_2 C_2 $ are contained in $ M_n(B) $. 
Write $ Y = Y(0) + x Y' $ for suitable $ Y' \in M_n(B(x)) $. If we set $ b = b_1 b_2 $ and $ M = M_1 Y' M_2 $ then the claim follows. 
\end{proof}  

The rank of a matrix $ Y $ with entries in an integral domain $ B $ is denoted by $ rk(Y) $;
by definition this is the usual rank of $ Y $ where $ Y $ 
is viewed as a matrix over the field of fractions $ \Quot(B) $. 

\begin{lemma} \label{shapag} 
Let $ B $ be an integral domain which is at the same time a finitely generated commutative algebra over $ \mathbb{C} $, 
and let $ n \in \mathbb{N} $. If there exists $ 0 \leq r \leq n $ such that 
$ Y(x) = (y_{ij}(x)) \in M_n(B[x]) $ satisfies $ rk(Y(0)_p) \leq r $ for all points $ p $ in a 
nonempty Zariski-open subset of the affine variety $ \Spm(B) $ of $ B $, then $ \det(Y(x)) = x^{n - r} b $ for some $ b \in B[x] $. 
\end{lemma}

\begin{proof} 
According to Lemma \ref{shapaghelp} there exists $ 0 \leq k \leq n $ such that $ M_1 Y M_2 = x M + b D_k $ for matrices $ M_1, M_2 $ and $ b \in B $ 
such that $ \det(M_1), \det(M_2) $ are nonzero. Let $ V $ be a non-empty Zariski-open set of $ \Spm(B) $ such 
that $ \det(M_1)_p \neq 0, \det(M_2)_p \neq 0, b_p \neq 0 $ and $ rk (Y(0)_p) \leq r $ for all $ p \in V $. This exists by assumption since $ \Spm(B) $ 
is irreducible. 
Since $rk(Y(0)_p) \leq r$ for $p \in V$ we get $k \leq r $. Hence 
$$
\det(M_1) \det(Y) \det(M_2) = x^{n - r} c 
$$
for some $ c \in B[x] $. Since $ \det(M_1), \det(M_2) \in B $ and $ B $ is an integral domain we conclude $ \det(Y) \in x^{n - r} B[x] $. 
\end{proof}  

For the following result compare Theorem 1.9 in \cite{deConciniKacrootof1} and Theorem 8.1 in \cite{HeckenbergerYamaneshapovalov}.  

\begin{theorem} \label{shapovalov}
For $ \nu \in \roots^+ $, the Shapovalov determinant $ \det(\Sh_\nu) $ 
of $ U_q(\mathfrak{g}) $ is given, up to multiplication by an invertible element of $U_q(\mathfrak{h})$, by 
$$
\det(\Sh_\nu) 
= \prod_{\beta \in {\bf \Delta}^+} \prod_{m = 1}^\infty (q_\beta^{2 (\rho, \beta^\vee)} K_\beta - q_\beta^{2m} K_\beta^{-1})^{P(\nu-m \beta)},
$$
where $ P $ is the Kostant partition function.
\end{theorem}

\begin{proof} 
Note first that if $ \nu \in \roots^+ $ is fixed, then for any $ \beta \in {\bf \Delta}^+ $ the expression 
$ P(\nu - m\beta) $ is zero for large values of $ m $. Therefore the right hand side of the asserted formula is in fact a finite product. 

As explained before Lemma \ref{shapovalovorderestimate}, it suffices to prove the claim for the form $ \S_\nu $ instead of $ \Sh_\nu $. 

Let $ D $ denote the dimension of $U_q(\mathfrak{n}_-)_{-\nu}$ and fix a basis $ x_1, \dots, x_D$ for this space. 
From Lemma \ref{shapovalovorderestimate} we have 
$$ 
\S_\nu(x_i, x_j) \in \sum_{\beta, \gamma \in \roots^+, \beta + \gamma = \nu} \mathbb{C} K_{\beta - \gamma}, 
$$
and so the determinant $ \det(\S_\nu) $, which is an order $ D $ polynomial in these terms, is contained in the linear subspace 
$$
L_\nu = \sum_{\beta, \gamma \in \roots^+, \beta + \gamma = D \nu} \mathbb{C} K_{\beta - \gamma} \subset U_q(\mathfrak{h}).
$$
For all $ 1 \leq j \leq n $ and $ m \in \mathbb{N} $ the polynomials 
$$
q_{\beta_j}^{2 (\rho, \beta_j^\vee)} K_{\beta_j} - q_{\beta_j}^{2m} K_{\beta_j}^{-1}
$$
have mutually distinct irreducible factors.
According to Lemma \ref{shapmultiplicity}, the double product in the theorem is contained in $ L_\nu $ as well. Therefore it suffices to show $ \det(\S_\nu) \neq 0 $ and that each term
$ (q_{\beta_j}^{2 (\rho, \beta_j^\vee)} K_{\beta_j} - q_{\beta_j}^{2m} K_{\beta_j}^{-1})^{P(\nu-m \beta_j)} $ divides $ \det(S_\nu) $. 

Let $ V \subset \mathfrak{h}^*_q $ be the affine subvariety of $ \mathfrak{h}^*_q $ defined by $ \det(\S_\nu) = 0 $. Then $ \lambda \in V $ if and only if 
\[
\chi_\lambda(\det(\Sh(x_i,x_j)_{ij}))=\det(\Sh^\lambda(x_i\cdot v_\lambda,x_j\cdot v_\lambda)_{ij})=0.
\]
By Lemma \ref{radicalmaximalsubmodule} this implies that $M(\lambda)$ is not simple, so by 
Theorem \ref{Vermacomplexirreducible} $ \lambda $ is not antidominant.  Thus
\begin{align*}
\chi_\lambda(q_{\beta_j}^{2 (\rho, \beta_j^\vee)} K_{\beta_j} - q_{\beta_j}^{2m} K_{\beta_j}^{-1})
&= q_{\beta_j}^{2(\rho,\beta_j^\vee)} q^{(\lambda,\beta_j)} - q_{\beta_j}^{2m} q^{-(\lambda,\beta_j)} \\
&= q^{-(\lambda,\beta_j)} \left( q_{\beta_j}^{2(\rho+\lambda,\beta_j^\vee)} - q_{\beta_j}^{2m}\right) = 0
\end{align*}
for some $\beta_j\in\mathbf{\Delta}^+$ and some $ m \in \NN $. It follows that $ V $ is a subset of the union of varieties defined by the zero sets of the polynomials 
$$
q_{\beta_j}^{2 (\rho, \beta_j^\vee)} K_{\beta_j} - q_{\beta_j}^{2m} K_{\beta_j}^{-1}
$$
for $ j = 1, \dots, n $ and $ m \in \mathbb{N} $. 

Using that $ \det(\S_\nu) $ is contained in $ L_\nu $, and hence a linear combination of terms of the form $ K_{2 \beta - D \nu} $ for $ \beta \in \roots^+ $,
we conclude that there exists an element $ f \in \mathbb{C}[K_1^{\pm 1}, \dots, K_N^{\pm 1}] $ which is invertible on $ V $ such that 
$$
\det(\S_\nu) = f \prod_{j = 1}^n 
\prod_{m = 1}^\infty (q_{\beta_j}^{2 (\rho, \beta_j^\vee)} K_{\beta_j} - q_{\beta_j}^{2m} K_{\beta_j}^{-1})^{N(\nu, \beta_j, m)}, 
$$
for certain exponents $ N(\nu, \beta_j, m) \in \mathbb{N}_0 $. It follows 
in particular that $ \det(\S_\nu) $ is nonzero. In order to finish the proof it suffices to show $ N(\nu, \beta_j, m) = P(\nu -m\beta_j) $ 
for all $ j = 1, \dots, n $ and $ m \in \mathbb{N} $. 

Fix $ 1 \leq j \leq n $ and $ m \in \mathbb{N} $, and choose $ w \in W $ and $ 1 \leq i \leq N $ such that $ \beta_j = w \alpha_i $. Consider the 
algebra $ A = \mathbb{C}[K_{w \alpha_1}^{\pm 1}, \dots, K_{w \alpha_N}^{\pm 1}] $. 
Moreover set 
$$
B = \mathbb{C}[K_{w \alpha_1}^{\pm 1}, \dots, K_{w \alpha_{i - 1}}^{\pm 1}, K_{w \alpha_i}, K_{w \alpha_{i + 1}}^{\pm 1}, \dots, K_{w \alpha_N}^{\pm 1}]. 
$$
If we put $ x = q_{\beta_j}^{2 (\rho, \beta_j^\vee)} K_{\beta_j} - q_{\beta_j}^{2m} K_{\beta_j}^{-1} $ then 
$
A = B[x]. 
$
Let $ Y = (\S_\nu(x_i, x_j)) \in M_D(A) = M_D(B[x]) $, and 
write $ Y(0) \in M_D(B) $ for the image of $ Y $ under the map induced by evaluation at $ 0 $. 

If $ \lambda \in \mathfrak{h}^*_q $ is an element in the variety of $ A/(x) \cong B $ then we have 
\[
\chi_\lambda ( q_{\beta_j}^{2 (\rho, \beta_j^\vee)} K_{\beta_j} - q_{\beta_j}^{2m} K_{\beta_j}^{-1}) = 0.
\]
By a calculation similar to that above, this implies
\[
(\lambda + \rho, \beta_j^\vee) = m + \frac{k}{2} i \hbar_{\beta_j}^{-1}
\]
for some $ k \in \ZZ $, and so $ s_{k, \beta_j} . \lambda = \lambda - m \beta_j < \lambda $.
According to Theorem \ref{Vermatheorem} 
there exists a submodule $ M(\lambda - m \beta_j) \subset M(\lambda) $. This submodule is contained in $ I(\lambda) $, 
and hence the radical of $ \S_\nu $ evaluated at $ \lambda $ contains a subspace of dimension 
\[
\dim(  M(\lambda-m\beta_j)_{\lambda-\nu}) = \dim (U_q(\mathfrak{n}_-)_{-\nu+m\beta_j}) = P(\nu-m\beta_j) .
\]
In other words, the rank of $ Y(0) $ evaluated at a point $ p $ of the variety of $ B = A/(x) $ satisfies 
$ rk(Y(0)_p) \leq D - P(\nu - m \beta_j) $. 
Hence $ \det(Y) =  x^{P(\nu - m \beta_j)} b $ for some $ b \in B $ according to Lemma \ref{shapag}. 
In particular, $  x^{P(\nu - m \beta_j)} $ is a factor of $ \det(S_\nu) $. This finishes the proof. 
\end{proof}

\subsection{Jantzen filtration and the BGG Theorem} 

In this section we discuss the Jantzen filtration and the BGG Theorem for $ U_q(\mathfrak{g}) $, giving us information about the composition factors 
of Verma modules. 

We shall first formulate the Jantzen filtration. Consider $ B = \mathbb{C}[T, T^{-1}] $, that is, the algebra of Laurent polynomials with 
coefficients in $ \mathbb{C} $. Informally, we will think of $ T $ as $ s^t $ where $ s \in \mathbb{C} $ such that $s^L = q$ and $t$ is an indeterminate. 
The ring $ B $ is a principal ideal domain, and we write $ \mathbb{L} $ for its field of quotients. 
\nomenclature[o$L$2]{$\mathbb{L}$}{field of quotients of $\CH[T,T^{-1}]$}%
Then $ \mathbb{L} $ and $ q \in \mathbb{L}^\times $ satisfy our general assumptions 
in the discussion of the quantized universal enveloping algebra 
in Chapter \ref{chuqg}. Let us write $ U_q^\mathbb{L}(\mathfrak{g}) $
for the quantized universal enveloping algebra defined over $ \mathbb{L} $. 
If $ \lambda \in \mathfrak{h}^*_q = \lie{h}^*/i\hbar^{-1}\mathbf{Q}^\vee $
is a weight with respect to the ground field $\mathbb{C}$,
we define a character $ \chi_{\lambda_T}: U_q^\mathbb{L}(\mathfrak{h}) \rightarrow \mathbb{L} $ 
by the formula 
$$ 
\chi(K_\mu) = q^{(\mu, \lambda)} T^{L (\mu, \rho)} = q^{(\mu, \lambda + \rho t)}, 
$$
where the second equality is obtained by formally writing $ T = s^t $. 

Let $ \M $ be a free module of finite rank over $ B $. We write $ M = \M \otimes_B \mathbb{C} $ 
for the vector space induced from the evaluation homomorphism $ \epsilon: B \rightarrow \mathbb{C}, \epsilon(f) = f(1) $. 
If $ S: \M \times \M \rightarrow B $ is a symmetric $ B $-bilinear form let us set
$$
\M^i = \{x \in \M \mid S(x, \M) \subset (T - 1)^i B \},
$$
and similarly $ M^i = \M^i \otimes_B \mathbb{C} $ for all $ i $. For any $ x \in B $ let $ v(x) \in \mathbb{N}_0 $ be the largest number $ n $ 
such that $ x $ can be written in the form $ x = y(T - 1)^n $ for some $ y \in B $. 

With this notation in place we have the following result, which is a special case of the Lemma in Section 5.6 of \cite{HumphreysO}. 

\begin{lemma} \label{Jantzenfiltrationhelp}
Let $ \M $ be a free module of finite rank over $ B = \mathbb{C}[T, T^{-1}] $. Assume that $ S: \M \times \M \rightarrow B $ is a 
nondegenerate symmetric $ B $-bilinear form on $ \M $, and let $ \det(S) $ be its determinant with respect to some basis of $ \M $. 
Then the following holds. 
\begin{bnum}
\item[a)] We have $ v(\det(S)) = \sum_{i > 0} \dim_\mathbb{C}(M^i) $.
\item[b)] For each $ i \geq 0 $ the map $ S_i: \M^i \times \M^i \rightarrow B $ given by 
$$ 
S_i(x,y) = (T - 1)^{-i} S(x,y) 
$$ 
is well-defined and induces a nondegenerate $ B $-bilinear form on $ M^i/M^{i + 1} $ with values in $ \mathbb{C} $. 
\end{bnum}
\end{lemma} 

\begin{proof}  We follow the treatment in \cite{HumphreysO}, adapted to the special case at hand for the convenience of the reader. 

$ a) $ Assume that $ \M $ is free of rank $ r $. The dual module $ \M^* = \Hom_B(\M, B) $ is again free of rank $ r $, as is the submodule 
$ \M^\vee \subset \M^* $ consisting of all linear functionals of the form $ x^\vee $ where $ x^\vee(y) = S(x,y) $ for $ x \in \M $. 
Using the structure theory for modules over principal ideal domains we find a basis $ e_1, \dots, e_r $ of $ \M $ and elements $ d_1, \dots, d_r \in B $ 
such that $ d_i e^i $ is a basis of $ \M^\vee $, where $ e^1, \dots, e^r $ is the dual basis of $ \M^* $ defined by $ e^i(e_j) = \delta_{ij} $. 
We then obtain $ \det(S) = d_1 \cdots d_r $ up to an invertible scalar. Moreover, there is another basis $ f_1, \dots f_r $ of $ \M $ such that 
$ f_i^\vee = d_i e^i $. 
If $ f = \sum a_j f_j \in \M $ and $ n_j = v(d_j) $ we get $ f \in \M^i $ iff $ v(S(e_j, f)) \geq i $ for all $ j $ 
iff $ v(a_j) \geq i - n_j $ for all $ j $. Hence $ \M^i $ is spanned by the elements $ f_j $ for which $ i \leq n_j $ and the $ (T - 1)^{i - n_j} f_j $ 
for which $ i > n_j $. In particular $ M^i = \M^i \otimes_B \mathbb{C} = 0 $ for sufficiently large $ i $. Moreover we obtain
$$
v(\det(S)) = \sum_{j = 1}^r n_j = \sum_{i > 0} |\{j| i \leq n_j\}| = \sum_{i > 0} \dim(M^i)  
$$
as desired. 

$ b) $ By construction, the form $ S_i $ takes values in $ B \subset \mathbb{L} $. We need to check that it descends to a form on $ M^i $. 
For this consider $ x = (T - 1) y \in \M^i $ and compute 
$$
S_i(x, \M^i) \subset (T - 1)^{-i} S(x, \M^i) \subset (T - 1)^{-i + 1} (\M, \M^i) \subset 
(T - 1) B,  
$$
so that $ S_i(x, \M^i) = 0 $ in $ B \otimes_B \mathbb{C} = \mathbb{C} $. Hence $ S^i $ induces a bilinear form $ M^i \times M^i \rightarrow \mathbb{C} $. 
In order to compute the radical $ R^i \subset M^i $ of the latter notice that $ M^{i + 1} \subset R^i $. 
Moreover, the coset of $ f_j $ as above is contained in $ R^i $ iff $ i < n_j $. Since the $ f_j $ with $ i \leq n_j $ form a basis of $ M^i $ we 
see that the form is nondegenerate on $ M^i/M^{i + 1} $, so that $ R^i = M^{i + 1} $. \end{proof}  

If $ M $ is a $ U_q(\mathfrak{g}) $-module then a bilinear form on $ M $ will be called contravariant if 
$$
(\tau(X) \cdot m, n) = (m, X \cdot n) 
$$
for all $ x \in U_q(\mathfrak{g}) $ and $ m,n \in M $. 

\begin{theorem}[Jantzen filtration] \label{Jantzenfiltration}
Let $ \lambda \in \mathfrak{h}^*_q $. Then there exist a filtration 
$$
M(\lambda) = M(\lambda)^0 \supset M(\lambda)^1 \supset M(\lambda)^2 \supset \cdots 
$$
of $ M(\lambda) $ with $ M(\lambda)^i = 0 $ for $ i $ large, such that the following conditions hold. 
\begin{bnum} 
\item[a)] Each nonzero quotient $ M(\lambda)^i/M(\lambda)^{i + 1} $ has a nondegenerate contravariant bilinear form. 
\item[b)] $ M(\lambda)^1 = I(\lambda) $ is the maximal submodule of $ M(\lambda) $. 
\item[c)] (Jantzen Sum Formula) The formal characters of the modules $ M(\lambda)^i $ satisfy 
$$
\sum_{i > 0} \ch(M(\lambda)^i) = \sum_{\substack{\alpha \in {\Delta}^+, k \in \mathbb{Z}_2 \\ s_{k, \alpha} . \lambda < \lambda}} \ch(M(s_{k, \alpha} . \lambda)). 
$$
\end{bnum} 
\end{theorem}

\begin{proof} 
We work over $ B = \mathbb{C}[T, T^{-1}] $ and its field of quotients $ \mathbb{L} $ as above. 
Given $ \lambda \in \mathfrak{h}^*_q  = \mathfrak{h}^*/i\hbar^{-1} \roots^\vee $ consider 
the character $ \chi_{\lambda_T} $ where $ \lambda_T = \lambda + t \rho $. 
We claim that $ M(\lambda_T) $ is a simple $ U_q^{\mathbb{L}}(\mathfrak{g}) $-module.  Indeed, by Theorem \ref{chilinkage} any submodule of $M(\lambda_T)$ 
must have highest weight in the same $\hat{W}$-linkage class as $\lambda_T$.  But if $\hat{w} = (\zeta,w) \in \hat{W}$ is nontrivial then
\[
\hat{w}.(\lambda + t\rho) - (\lambda + t\rho) = (w\lambda + w\rho - \lambda - \rho + \zeta) + t(w\rho-\rho) \notin \roots,
\]
so the $\hat{W}$-linkage class of $\lambda_T$ is trivial.
Hence the unique contravariant bilinear form on $ M(\lambda_T) $ with values in $ \mathbb{L} $ 
such that $ (v_\lambda, v_\lambda) = 1 $ is nondegenerate. 

Using the PBW Theorem \ref{PBWgeneralfield} we obtain from 
the $ B $-form $ U_q^B(\mathfrak{g}) = B \otimes U_q(\mathfrak{g}) $ a $ B $-form $ \M(\lambda_T) $ 
of the Verma module $ M(\lambda_T) $ in which all weight subspaces are free of finite rank. 
For $ \mu \in \roots^+ $ define 
$$
\M(\lambda_T)^i_{\lambda_T - \mu}= \{x \in \M(\lambda_T)_{\lambda_T - \mu} \mid (x, \M(\lambda_T)_{\lambda_T - \mu}) \subset (T - 1)^i B \} 
$$
and 
$$
\M(\lambda_T)^i = \sum_{\mu \in \roots^+} \M(\lambda_T)^i_{\lambda_T - \mu}.  
$$
It is straightforward to check that $ \M(\lambda_T)^i $ is a $ U_q^B(\mathfrak{g}) $-submodule of $ \M(\lambda_T) $, and that 
these submodules form a decreasing filtration. 

Setting $ T = 1 $ yields a decreasing filtration $ M(\lambda)^i = \M(\lambda_T)^i/(T - 1)\M(\lambda_T)^i $ of $ M(\lambda) $ such 
that $ M^0(\lambda) = M(\lambda) $. Thanks to Lemma \ref{Jantzenfiltrationhelp} the quotients of this filtration acquire nondegenerate contravariant bilinear forms. 
Since $ M(\lambda) $ has finite length 
we see that $ M(\lambda)^i = 0 $ for $ i \gg 0 $. This proves part $ a) $. 

For part $ b) $ we only need to observe that the contravariant form on $ M(\lambda)/M(\lambda)^1 $ is nondegenerate, and therefore $ M(\lambda)^1 = I(\lambda) $ 
by Lemma \ref{radicalmaximalsubmodule}. 

It remains to prove the Jantzen Sum Formula in part $ c) $. Consider $ \nu \in \roots^+ $ and denote the determinant of the Shapovalov 
form of $ \M(\lambda_T) $ on the $ \lambda_T - \nu $ weight space by $ \det(\Sh_\nu) $. 
Since $ B $ contains $ \mathbb{C} $, the formula for the Shapovalov determinant in Theorem \ref{shapovalov} continues to 
hold over $ U_q^B(\mathfrak{g}) $. More precisely, we obtain
$$
\chi_{\lambda_T}(\det(\Sh_\nu)) 
= \prod_{\alpha \in {\bf \Delta}^+} \prod_{m = 1}^\infty (q^{2 (\lambda + \rho + t \rho, \alpha)} - q^{m(\alpha, \alpha)})^{P(\nu-m\alpha)}, 
$$
up to an invertible element of $ B $. Applying $ v $ to
the term
\[ x = q^{2 (\lambda + \rho + t \rho, \alpha)} - q^{m(\alpha, \alpha)} = q^{2(\lambda + \rho, \alpha)} T^{2L(\rho, \alpha)} - q^{m (\alpha, \alpha)} 
\]
gives $ 0 $ unless $ q^{m(\alpha, \alpha)} = q^{2(\lambda + \rho, \alpha)} $, or equivalently, unless $ (\lambda + \rho, \alpha^\vee) = m $ 
modulo $ \frac{1}{2} i \hbar_\alpha^{-1} \mathbb{Z}$, where as usual $\hbar_\alpha=d_\alpha\hbar$. In this case 
we get $ v(x) = 1 $. The contribution of the $ (\lambda - \nu) $-weight space to the valuation is 
therefore $ P(\nu - (\lambda + \rho, \alpha^\vee) \alpha) $, where $ P $ denotes Kostant's partition function. 
Using Lemma \ref{Jantzenfiltrationhelp} we obtain 
\begin{align*}
\sum_{i > 0} \ch(M(\lambda)^i) &= \sum_{\nu \in \roots^+} \sum_{\substack{\alpha \in {\Delta}^+, k \in \mathbb{Z}_2 \\ s_{k, \alpha} . \lambda < \lambda}} 
P(\nu - (\lambda + \rho, \alpha^\vee) \alpha) e^{\lambda - \nu} \\
&= \sum_{\substack{\alpha \in {\Delta}^+, k \in \mathbb{Z}_2 \\ s_{k, \alpha} . \lambda < \lambda}} \sum_{\nu \in \roots^+} 
P(\nu) e^{\lambda - (\lambda + \rho, \alpha^\vee) \alpha - \nu} \\ 
&= \sum_{\substack{\alpha \in {\Delta}^+, k \in \mathbb{Z}_2 \\ s_{k, \alpha} . \lambda < \lambda}} \ch(M(s_{k, \alpha} . \lambda)) 
\end{align*}
as desired. 
\end{proof}  

Notice that the Jantzen filtration reduces to the trivial filtration $ M(\lambda) = M(\lambda)^0 \supset M(\lambda)^1 = 0 $ 
in the case that $ M(\lambda) $ is simple. The Jantzen sum formula does not provide any information in this case. 

Let $ \mu, \lambda \in \mathfrak{h}^*_q $. Recall that $\mu$ is strongly linked to $\lambda$, written $ \mu \uparrow \lambda $, if $ \mu = \lambda $ or if there exists a chain of positive roots 
$ \alpha_1, \dots, \alpha_r \in {\bf \Delta}^+ $ and $ k_1, \dots k_r \in \mathbb{Z}_2 $ such that 
$$
\mu = s_{k_1, \alpha_1} \cdots s_{k_r, \alpha_r} . \lambda < s_{k_2, \alpha_2} \cdots s_{k_r, \alpha_r} . \lambda < \cdots < s_{k_r, \alpha_r} . \lambda < \lambda.  
$$ 
Using the Jantzen filtration we shall now prove an analogue of the BGG Theorem, compare Section 5.1 in \cite{HumphreysO}.  

\begin{theorem}[The BGG Theorem] \label{BGGtheorem}
Let $ \mu, \lambda \in \mathfrak{h}^*_q $. Then $ [M(\lambda): V(\mu)] \neq 0 $ iff $ \mu $ is strongly linked to $ \lambda $. 
\end{theorem} 

\begin{proof} 
From Theorem \ref{Vermatheorem} we know that $ \mu \uparrow \lambda $ implies $ [M(\lambda): V(\mu)] \neq 0 $. 

For the converse let us use the Jantzen filtration and induction on the number $ k $ of weights $ \mu $ linked to $ \lambda $ satisfying $ \mu \leq \lambda $. 
If $ k = 1 $ then $ \lambda $ is minimal in its linkage class, so there is nothing to prove. Assume now that the claim is proved for all weights 
and some $ k $. Suppose $ [M(\lambda): V(\mu)] > 0 $ for $ \mu < \lambda $. This means $ [M(\lambda)^1: V(\mu)] > 0 $ for the first step $ M(\lambda)^1 $ 
in the Jantzen filtration of $ M(\lambda) $. The sum formula in Theorem \ref{Jantzenfiltration} forces $ [M(s_{k, \alpha} . \lambda): V(\mu)] > 0 $ 
for some $ \alpha \in {\bf \Delta}_{[\lambda]}^+, k \in \mathbb{Z}_2 $. By our inductive hypothesis, there 
exist $ \alpha_1, \dots, \alpha_r \in {\bf \Delta}^+, k_1, \dots, k_r \in \mathbb{Z}_2 $ such that 
$$
\mu = s_{k_1, \alpha_1} \cdots s_{k_r, \alpha_r} s_{k, \alpha} . \lambda < s_{k_2, \alpha_2} \cdots s_{k_r, \alpha_r} s_{k, \alpha} . \lambda 
< \cdots < s_{k_r, \alpha_r} s_{k, \alpha} . \lambda < s_{k, \alpha} . \lambda.  
$$ 
Appending this chain with $ s_{k, \alpha} . \lambda < \lambda $ yields $ \mu \uparrow \lambda $ as desired. 
\end{proof}

\subsection{The PRV determinant} \label{secprv}

In this section we discuss the Parthasarathy-Ranga Rao-Varadarajan determinants of $ U_q(\mathfrak{g}) $.  
In the classical setting they were introduced in \cite{PRVannals}. For more information on the quantum case we refer to \cite{Josephbook}.

\subsubsection{The spaces $ F \Hom(M, N) $} \label{sec:FHomMN}

Before discussing the PRV-determinants we will need a couple of preparations regarding the spaces $ F \Hom(M, N) $, which we now describe.

Let $ M, N $ be $ U_q(\mathfrak{g}) $-modules. The space of $\CC$-linear maps $ \Hom(M, N) $ becomes a $U_q(\mathfrak{g})$-module by setting 
$$
(X \rightarrow T)(m) = X_{(1)} \cdot T(\hat{S}(X_{(2)}) \cdot m) 
$$
\label{nom:adjoint_action_on_Hom}%
for $ T \in \Hom(M, N) $, $X\in U_q(\mathfrak{g})$ and $ m \in M $.  We will refer to this as the adjoint action of $U_q(\mathfrak{g})$.
We denote by $ F\Hom(M, N) $ 
\nomenclature{$F\Hom(M,N)$}{locally finite part of $\Hom(M,N)$ with respect to the adjoint action}%
the locally finite part of $ \Hom(M, N) $ with respect to the adjoint action. 

We can also define a $U_q(\mathfrak{g})$-bimodule structure on $\Hom(M,N)$ by
$$
(Y \cdot T \cdot Z)(m) = Y \cdot T(Z \cdot m).
$$
\label{nom:Hom_bimodule}%
for $ Y,Z\in U_q(\mathfrak{g})$, $T\in\Hom(M,N)$, $m\in M$.
By restriction, $F\Hom(M,N)$ becomes an $F U_q(\mathfrak{g})$-bimodule.  
Indeed, notice that 
\begin{align*}
X \rightarrow (Y \cdot T \cdot Z)(m) &= X_{(1)} Y \cdot T(Z \hat{S}(X_{(2)}) \cdot m) \\
&= (X_{(1)} \rightarrow Y) \cdot (X_{(2)} \rightarrow T)((X_{(3)} \rightarrow Z) \cdot m)
\end{align*}
for all $ m \in M $. Hence if $ Y,Z \in FU_q(\mathfrak{g}) $ and $ T \in F\Hom(M,N) $, then $ Y \cdot T \cdot Z $ is again contained in $ F\Hom(M, N) $. 

The $ FU_q(\mathfrak{g}) $-bimodule structure on $ F\Hom(M, N) $ is compatible with the adjoint $ U_q(\mathfrak{g}) $-action
in the following sense.

\begin{definition} \label{defadcompatibility}
	Let $ V $ be an $ FU_q(\mathfrak{g}) $-bimodule, with left and right actions denoted by $ X \cdot v $ and $ v \cdot X $ 
	for $ X \in FU_q(\mathfrak{g}) $ and $ v \in V $, respectively. We say that a left module structure of $ U_q(\mathfrak{g}) $ on $ V $, written $ X \rightarrow v $ 
	for $ X \in U_q(\mathfrak{g}) $ and $ v \in V $, is compatible if 
	$$
	(X_{(1)} \rightarrow v) \cdot X_{(2)} = X \cdot v 
	$$
	for all $ X \in FU_q(\mathfrak{g}) $ and $ v \in V $.  
	Equipped with such a structure, $V$ will be called a $U_q(\lie{g})$-compatible $FU_q(\lie{g})$-bimodule.
\end{definition}

Note that the formula in Definition \ref{defadcompatibility} makes sense since $ FU_q(\mathfrak{g}) $ is a left coideal, 
see Lemma \ref{lem:FUq_coideal}. Morally, the compatibility condition can be written as 
$$
X \rightarrow v = X_{(1)} \cdot v \cdot \hat{S}(X_{(2)}); 
$$
however, the right hand side of this formula might not be well-defined, even if we consider $ X \in FU_q(\mathfrak{g}) $. 
Nevertheless, inspired by this formula, we shall always refer to a compatible action of $U_q(\lie{g})$ on an $FU_q(\lie{g})$-bimodule as the adjoint action.

Next, we introduce similar structures associated to the linear dual $(M\otimes N)^*$, where again $M$ and $N$ are in category $\O$. 
Recall the algebra antiautomorphism and coalgebra automorphism $\tau$ of $U_q(\mathfrak{g})$ from Lemma \ref{deftau}.
We introduce a $U_q(\mathfrak{g})$-module structure on $(M\otimes N)^*$ by
$$
(X \rightarrow \varphi)(m \otimes n) = \varphi(\hat{S}(X_{(2)}) \cdot m \otimes \tau(X_{(1)}) \cdot n) 
$$
\label{nom:adjoint_action_on_MxN}%
for $X\in U_q(\mathfrak{g})$, $\varphi\in(M\otimes N)^*$ and $m\otimes n \in M\otimes N$.
Again, we refer to this as the adjoint action of $U_q(\mathfrak{g})$ on  $(M\otimes N)^*$.
We also have a $U_q(\lie{g})$-bimodule structure  on $(M\otimes N)^*$ given by
\[
 (Y\cdot\varphi\cdot Z)(m \otimes n) = \varphi(Z\cdot m \otimes \tau(Y)\cdot n),
\]
\label{nom:MxN_bimodule}%
for $Y,Z\in U_q(\lie{g})$.  

We denote by $F((M\otimes N)^*)$
\nomenclature{$F((M\otimes N)^*)$}{locally finite part of $(M\otimes N)^*$ with respect to the adjoint action}%
the locally finite part of $(M\otimes N)^*$ with respect to the adjoint action of $U_q(\lie{g})$ defined above.  The following lemma shows that by restricting the above actions, $F((M\otimes N)^*)$ becomes a $U_q(\lie{g})$-compatible $FU_q(\lie{g})$-bimodule in the sense of Definition \ref{defadcompatibility}.

\begin{lemma} \label{locallyfiniteduallemma}
Let $ M, N $ be in category $ \O $. Then there exists an isomorphism of $FU_q(\lie{g})$-bimodules
\begin{align*}
F \Hom(M, N^\vee) \stackrel{ \cong }{\to} F((M \otimes N)^*),
\end{align*}
which is $U_q(\mathfrak{g})$-linear for the adjoint actions denoted by $\rightarrow$ above. Explicitly, the isomorphism is the restriction of the natural map $ \gamma: \Hom(M, N^*) \rightarrow (M \otimes N)^* $ defined by 
	\[
	\gamma(T)(m \otimes n) = T(m)(n)
	\]
for $T\in\Hom(M, N^*)$ and $m\in M$, $n\in N$.
\end{lemma} 

\begin{proof}  
To begin with, let us show that $ F\Hom(M,N^\vee) = F\Hom(M,N^*)$.  Here, we are equipping $N^*$ with the obvious extension of the $U_q(\mathfrak{g})$-action 
on $N^\vee$, namely
\[
(X \cdot f)(n) = f(\tau(X)\cdot n), \qquad f \in N^*, ~X\in U_q(\mathfrak{g}), ~n\in N.
\]
Let $T\in F\Hom(M,N^*)$. For $K \in U_q(\mathfrak{h})$, we have
\[
K \cdot (T(m)) = (K_{(1)}\rightarrow T)(K_{(2)}\cdot m).
\]
Since $T$ is locally finite and $m\in M$ is contained in a finite sum of weight spaces, we see that $T(m)$ generates a finite dimensional $U_q(\mathfrak{h})$-submodule of $N^*$ and so must belong to $N^\vee$.

Now consider the linear map $ \gamma: \Hom(M, N^*) \rightarrow (M \otimes N)^* $ defined above.
General facts about tensor products imply that $ \gamma $ is a linear isomorphism. Since 
\begin{align*} 
(X \rightarrow \gamma(T))(m \otimes n) 
  &= \gamma(T)(\hat{S}(X_{(2)}) \cdot m \otimes \tau(X_{(1)}) \cdot n) \\
  &= T(\hat{S}(X_{(2)}) \cdot m)(\tau(X_{(1)}) \cdot n) \\
   &= (X_{(1)} \cdot T(\hat{S}(X_{(2)} \cdot m))(n) \\
   &= \gamma(X \rightarrow T)(m \otimes n) 
\end{align*}
for $X\in U_q(\lie{g})$, $ m \in M, n \in N $ we see that $ \gamma $ is $ U_q(\mathfrak{g}) $-linear and so restricts to an isomorphism on the locally finite parts. 

Finally, for any $Y,Z \in U_q(\mathfrak{g})$, $T\in\Hom(M,N^*)$ and $m\otimes n\in M\otimes N$ we have
 \begin{align*}
  (Y \cdot \gamma(T) \cdot Z)(m \otimes n)
  &= \gamma(T)(Z \cdot m \otimes \tau(Y) \cdot n) \\ 
   &= T(Z \cdot m)(\tau(Y) \cdot n) \\ 
   &= (Y\cdot T(Z \cdot m))(n) \\
   &=  \gamma(Y \cdot T \cdot Z)(m \otimes n).
 \end{align*}
Since $F\Hom(M,N^\vee)$ is invariant for the left and right actions of $FU_q(\lie{g})$, we deduce that the same is true for $F((M\otimes N)^*)$, and $\gamma: F \Hom(M, N^\vee) \rightarrow F((M \otimes N)^*) $ is an isomorphism of $ FU_q(\mathfrak{g})$-bimodules. 
\end{proof}

Consider the automorphism $ \theta = \tau \hat{S} = \hat{S}^{-1} \tau $.
\label{nom:theta1}%
It is an algebra homomorphism and a coalgebra antihomomorphism. Note also that $ \theta $ 
is involutive, since $ \theta^2 = \tau \hat{S} \hat{S}^{-1} \tau = \tau^2 = \id $. 

\begin{lemma} \label{locallyfinitedualtransposelemma}
Let $ M, N $ be in category $ \O $. Then the flip map $ M \otimes N \rightarrow N \otimes M $ induces a linear isomorphism 
\begin{align*}
\alpha: F((M \otimes N)^*) \rightarrow F((N \otimes M)^*)
\end{align*}
\nomenclature[o$\alpha$]{$\alpha$}{linear isomorphism $\alpha: F((M \otimes N)^*) \rightarrow F((N \otimes M)^*)$}%
such that
$$
\alpha(Y \cdot \varphi \cdot Z) = \tau(Z) \cdot \alpha(\varphi) \cdot \tau(Y), \qquad \alpha(X \rightarrow \varphi) = \theta(X) \rightarrow \alpha(\varphi)
$$
for all $ Y,Z \in FU_q(\mathfrak{g}) $ and $ X \in U_q(\mathfrak{g}) $. \\
Similarly, there exists a linear isomorphism 
\begin{align*}
\beta: F \Hom(M, N) \rightarrow F \Hom(N^\vee, M^\vee) 
\end{align*}
\nomenclature[o$\beta$3]{$\beta$}{linear isomorphism $\beta: F \Hom(M, N) \rightarrow F \Hom(N^\vee, M^\vee) $}%
such that 
$$
\beta(Y \cdot T \cdot Z) = \tau(Z) \cdot \beta(T) \cdot \tau(Y), \qquad \beta(X \rightarrow T) = \theta(X) \rightarrow \beta(T)
$$
for all $ Y,Z \in FU_q(\mathfrak{g}) $ and $ X \in U_q(\mathfrak{g}) $. 
\end{lemma} 

\begin{proof}  
We define $ \alpha: (M \otimes N)^* \rightarrow (N \otimes M)^* $ by $ \alpha(\varphi)(n \otimes m) = \varphi(m \otimes n) $. Clearly $ \alpha $ 
is a linear isomorphism. To check that this isomorphism is compatible with the locally finite parts 
we compute 
\begin{align*}
\theta(X) \rightarrow \alpha(\varphi)(n \otimes m) &= \alpha(\varphi)(\hat{S}(\theta(X_{(1)})) \cdot n \otimes \tau(\theta(X_{(2)})) \cdot m) \\
&= \alpha(\varphi)(\tau(X_{(1)}) \cdot n \otimes \hat{S}(X_{(2)}) \cdot m) \\
&= \varphi(\hat{S}(X_{(2)}) \cdot m \otimes \tau(X_{(1)}) \cdot n) \\
&= \alpha(X \rightarrow \varphi)(n \otimes m). 
\end{align*}
Hence $ \alpha $ restricts to an isomorphism $ F((M \otimes N)^*) \rightarrow F((N \otimes M)^*) $. 
Moreover we have 
\begin{align*}
\alpha(Y \cdot \varphi \cdot Z)(n \otimes m) 
  &= \varphi(Z\cdot m \otimes \tau(Y)\cdot n) \\
  &= \alpha(\varphi)(\tau(Y)\cdot n \otimes Z\cdot m ) \\
  &= (\tau(Z)\cdot \alpha(\varphi) \cdot \tau(Y))(n\otimes m)
\end{align*}
for all $ Y, Z \in FU_q(\mathfrak{g}) $. 

In order to prove the second assertion it suffices to construct an isomorphism $ \beta: F \Hom(M, N^\vee) \rightarrow F\Hom(N, M^\vee) $ with the 
claimed properties. This can be done by combining the map $ \alpha $ from the first part of the proof with the isomorphism $\gamma$
from Lemma \ref{locallyfiniteduallemma}. 
\end{proof}

\subsubsection{Conditions for $ F\Hom(M,N) = 0 $}

In this subsection we use Gelfand-Kirillov dimension as a tool to derive 
sufficient conditions for $ F \Hom(M, N) = 0 $ for certain modules $ M, N \in \O $. For background on the 
Gelfand-Kirillov dimension we refer to \cite{BKgkdimension}. 

Let $ A $ be an algebra over $ \mathbb{C} $ and let $ V \subset A $ be a linear subspace. We write $ V^n $ for the linear subspace generated by all $ n $-fold 
products $ a_1 \cdots a_n $ with $ a_1, \dots, a_n \in V $, in addition we set $ V^0 = \mathbb{C} $. Assume that $ A $ is finitely generated and that $ V $ 
is a finite dimensional generating subspace, so that $ \bigcup_{n = 0}^\infty V^n = A $. 
Moreover let $ M $ be a finitely generated $ A $-module with generating set $ M^0 $, so that $ \bigcup_{n = 0}^\infty M^n = M $ where $ M^n = V^n \cdot M^0 $. 
Then the Gelfand-Kirillov dimension of $ M $ over $ A $ is defined by
$$
d_A(M) = \limsup_{n \rightarrow \infty} \frac{\log(\dim(M^n))}{\log(n)}. 
$$
\nomenclature[o$d(M)$1]{$d_A(M)$}{Gelfand-Kirillov dimension of a finitely generated module $M$ over $A$}%
It is not hard to check that this does not depend on the choice of the generating subspaces $ V $ and $ M^0 $. 

Let us collect some general facts regarding the GK dimension.
\begin{lemma} \label{gkgeneral}
Let $ A $ be a finitely generated algebra. 
\begin{bnum}
\item[a)] If $ 0 \rightarrow K \rightarrow E \rightarrow Q \rightarrow 0 $ is an extension of finitely generated $ A $-modules then 
$ d_A(K) \leq d_A(E) $ and $ d_A(Q) \leq d_A(E) $. 
\item[b)] If $ I \subset A $ is a left ideal containing an element $ a \in I $ which is not a right zero divisor then $ d_A(A/I) \leq d_A(A) - 1 $. 
\end{bnum}
\end{lemma} 
\begin{proof}  $ a) $ Pick a generating subspace $ K^0 $ for $ K $ and extend it to a generating subspace $ E^0 $ for $ E $. Then $ K^n \subset E^n $ 
for all $ n $ and hence $ \dim(K^n) \leq \dim(E^n) $. Similarly, if $ \pi: E \rightarrow Q $ denotes the quotient map then $ Q^0 = \pi(E^0) $ 
is a generating subspace for $ Q $, and we have $ \dim(Q^n) \leq \dim(E^n) $ for all $ n $. This yields the claim. 

$ b) $ We follow the proof of Theorem 3.4 in \cite{BKgkdimension}. Without loss of generality we may assume that $ d_A(A) $ is finite. 
Let $ V \subset A $ be a finite dimensional generating subspace containing $ 1 $. 
For each $ n $ let $ D_n \subset V^n $ be a complement of $ V^n \cap I $. If $ \pi: A \rightarrow A/I $ is the quotient map then we get 
$ \pi(D_n) = \pi(V^n) $. 

We claim that $ D_n + D_n a + \cdots D_n a^m $ is a direct sum for all $ m \in \mathbb{N} $. Indeed, from $ x_0 + x_1 a + \cdots + x_m a^m = 0 $ 
with $ x_j \in D_n $ for all $ j $ we get $ \pi(x_0) = 0 $ because $ a \in I $. Hence $ x_0 = 0 $, and since $ a $ is not a right zero divisor 
we deduce $ x_1 + x_2 a + \cdots + x_m a^m = 0 $. Therefore the assertion follows by induction. 

Next observe $ D_n + D_n a + \cdots D_n a^n \subset V^{2n} $ 
so that $ \dim(V^{2n}) \geq n \dim(D_n) = n \dim(\pi(V^n)) $. 
Since $ W = V^2 $ is again a generating subspace for $ A $ and $ W^n = V^{2n} $ we conclude
$$
d_A(A) = \limsup_{n \rightarrow \infty} \frac{\log(\dim(V^{2n}))}{\log(n)} \geq \limsup_{n \rightarrow \infty} \frac{\log(n \dim(\pi(V^n)))}{\log(n)} 
= 1 + d_A(A/I) 
$$
as desired. \end{proof}  

We will be interested in the case that $ A = U_q(\mathfrak{n}_-) $ or $ A = U_q(\mathfrak{b}) $ and $ M $ is a $ U_q(\mathfrak{g}) $-module contained in $ \O $, 
viewed as an $ A $-module. Note that $ d_A(M) $ will be the same for both choices of $ A $, and we will write $ d(M) = d_A(M) $
\nomenclature[o$d(M)$2]{$d(M)$}{Gelfand-Kirillov dimension of a module $M\in\O$ over $U_q(\lie{n}_-)$ or $U_q(\lie{b})$}%
in either case. 

In the next lemma, recall that $ w_0 $ denotes the longest element of the Weyl group.

\begin{lemma} \label{gkspecial} 
Let $ A = U_q(\mathfrak{n}_-) $. 
\begin{bnum}
\item[a)] If $ \lambda \in \mathfrak{h}^*_q $ then $ d(M(\lambda)) = l(w_0) $. 
\item[b)] If $ Q $ is a proper quotient of a Verma module then $ d(Q) < l(w_0) $. 
\item[c)] If $ M $ is in $ \O $ and $\mu \in \weights^+$ we have $ d(M \otimes V(\mu)) = d(M) $. 
\end{bnum}
\end{lemma} 

\begin{proof}
$ a) $ Every Verma module is a free $ A $-module of rank $ 1 $. Hence the claim can be obtained by invoking a filtration of $ A $ as in the 
proof of Proposition \ref{uqgnoetheriandomain}, taking into account Satz 5.5 in \cite{BKgkdimension}.  

$ b) $ This follows again from the fact that every Verma module is a free $ A $-module of rank $ 1 $ combined with $ A $ having no zero-divisors 
and part $ b) $ of Lemma \ref{gkgeneral}. 

$ c) $ We use that $ M $, being in $ \O $, is finitely generated as an $ A $-module. If $ V \subset A $ is the linear span of $ 1 $ and the 
generators $ F_1, \dots, F_N $ 
then $ V $ is a generating subspace of $ A $. Let $ M^0 $ be a finite dimensional generating subspace of $ M $ over $ A $. 
Then $ M^0 \otimes V(\mu) $ is a generating subspace for $ M \otimes V(\mu) $, and we have 
$$ 
V^n \cdot (M^0 \otimes V(\mu)) \subset V^n \cdot M^0 \otimes V(\mu) = M^n \otimes V(\mu), 
$$
which implies $ d(M \otimes V(\mu)) \leq d(M) $.  
If $ v^\mu \in V(\mu) $ is a lowest weight vector then $ M \otimes v^\mu $ is isomorphic to $ M $ 
as an $ A $-module, so that $ d(M \otimes V(\mu)) \geq d(M) $. Hence we get $ d(M \otimes V(\mu)) = d(M) $ as desired. 
\end{proof}  

\begin{lemma} \label{lem:FHom_0}
For $M,N$ in $\O$, the following are equivalent:
\begin{bnum}
\item[a)] $F\Hom(M,N) = 0$.
\item[b)] $\Hom_{U_q(\mathfrak{g})}(V(\mu)\otimes M,N) = 0$ for every $\mu\in\weights^+$.
\item[c)] $\Hom_{U_q(\mathfrak{g})}(M,V(\mu)\otimes N) = 0$ for every $\mu\in\weights^+$.
\end{bnum}
\end{lemma}

\begin{proof}
The map $\kappa : \Hom(M,N) \otimes V(\mu)^* \to \Hom(V(\mu) \otimes M , N)$ defined by
\[
\kappa(T \otimes f)(v \otimes m) = f(v) T(m) 
\]
is a $U_q(\mathfrak{g})$-linear isomorphism:
\begin{align*}
(X \rightarrow \kappa(T \otimes f))(v\otimes m) 
&= f(\hat{S}(X_{(3)})\cdot v) \, X_{(1)}\cdot T(\hat{S}(X_{(2)})\cdot m) \\
&= (X_{(1)} \rightarrow T)(m) \, (X_{(2)}\rightarrow f)(v).
\end{align*}
Therefore, if $\Hom(M,N)$ contains a finite dimensional $U_q(\mathfrak{g})$-submodule of highest weight $\mu$ then $\Hom(V(\mu)\otimes M,N)$ contains a trivial submodule, and conversely.  An element of $\Hom(V(\mu)\otimes M,N)$ spans a trivial submodule if and only if it is $U_q(\mathfrak{g})$-linear.  This proves that $a)$ and $b)$ are equivalent.

A similar argument 
proves the equivalence of $a)$ and $c)$.
\end{proof}

Note that  $V(\mu)\otimes M$ can be replaced by $M\otimes V(\mu)$ in $b)$ and $c)$ of the above lemma, since they are isomorphic as $U_q(\mathfrak{g})$-modules.

\begin{prop} \label{GKlemma} 
Let $ M, N $ be in category $ \O $. 
\begin{bnum} 
\item[a)] If $ d(M) < d(K) $ for every simple submodule $ K $ of $ N $ we have $ F\Hom(M, N) = 0 $.
\item[b)] If $ d(N) < d(Q) $ for every simple quotient $ Q $ of $ M $ we have $ F\Hom(M, N) = 0 $.
\end{bnum}
\end{prop} 

\begin{proof} 
$ a) $ By Lemma \ref{lem:FHom_0}, it suffices to show $ \Hom_{U_q(\mathfrak{g})}(M \otimes V(\mu), N) = 0 $ for all $ \mu \in \weights^+ $. 
Let $ f: M \otimes V(\mu) \rightarrow N $ be a nonzero $ U_q(\mathfrak{g}) $-linear map. 
According to part $ c) $ of Lemma \ref{gkspecial} we have $ d(M \otimes V(\mu)) = d(M) $, and we have $ d(f(M \otimes V(\mu)) \leq d(M \otimes V(\mu)) = d(M) < d(K) $
for all simple submodules $ K \subset N $ by assumption. On the other hand, if $ K \subset f(M \otimes V(\mu)) $ is a simple submodule 
then $ d(K) \leq d(f(M \otimes V(\mu))) $ by Lemma \ref{gkgeneral}. This is a contradiction, thus giving the claim.

$ b) $ In a similar fashion it is enough to show $ \Hom_{U_q(\mathfrak{g})}(M, V(\mu) \otimes N) = 0 $ for all $ \mu \in \weights^+ $. 
Let $ f: M \rightarrow V(\mu) \otimes N $ be a nonzero $ U_q(\mathfrak{g}) $-linear map. Using $ V(\mu) \otimes N \cong N \otimes V(\mu) $ 
we have $ d(N) = d(V(\mu) \otimes N) $ according to part $ c) $ of Lemma \ref{gkspecial}, and hence $ d(V(\mu) \otimes N) < d(Q) \leq d(f(M))$ 
for every simple quotient $ Q $ of $ f(M) $ by assumption. 
In addition, $ d(f(M)) \leq d(V(\mu) \otimes N) $ by Lemma \ref{gkgeneral} because $ f(M) \subset V(\mu) \otimes N $ is a submodule. Again we 
obtain a contradiction. 
\end{proof}  

As a consequence we arrive at the following result. 

\begin{prop} \label{GKhelp} 
Let $ \lambda \in \mathfrak{h}^*_q $. If $ Q $ is a proper quotient of a Verma module then $ F\Hom(Q, M(\lambda)) = 0 $. 
Similarly, if $ M(\lambda) $ is simple then $ F\Hom(M(\lambda), Q) = 0 $. 
\end{prop} 

\begin{proof} 
If $ Q $ is a proper quotient of a Verma module then $ d(Q) < l(w_0) $ according to Lemma \ref{gkspecial}. 
In particular, we have $ d(Q) < d(K) $ for the unique simple submodule $ K = \Soc(M(\lambda)) \subset M(\lambda) $, 
see Lemma \ref{Vermasocle}. Hence the first assertion follows from Proposition \ref{GKlemma} $ a) $. 
Similarly, the second claim follows from Proposition \ref{GKlemma} $ b) $. 
\end{proof}

\subsubsection{Multiplicities in $ F \Hom(M, N) $}

If $ V $ is any $ U_q(\mathfrak{g}) $-module and $ \lambda \in \mathfrak{h}^*_q $, then $ V \otimes M(\lambda) $ is a free $ U_q(\mathfrak{n}_-) $-module. 
To see this, let $V_\tau$ denote $V$ equipped with the trivial action of $ U_q(\mathfrak{n}_-) $ and observe that the 
isomorphism $ V_\tau \otimes M(\lambda) \cong V \otimes M(\lambda) $ defined by
\[
v \otimes Y \cdot v_\lambda \mapsto Y_{(1)} \cdot v \otimes Y_{(2)} \cdot v_\lambda, 
\]
for $ v \in V $, $ Y \in U_q(\mathfrak{n}_-) $ is $ U_q(\mathfrak{n}_-) $-linear.

\begin{lemma} \label{Vermatensorfiltration} 
Let $ \lambda \in \mathfrak{h}^*_q $ and $ \nu \in \weights^+ $. Then there exists a finite decreasing filtration 
$$
V(\nu) \otimes M(\lambda) = M_0 \supset M_1 \supset \cdots \supset M_k \supset 0 
$$
with quotients isomorphic to $ M(\lambda + \gamma) $, where $ \gamma $ runs over all weights of $ V(\nu) $ counted with 
multiplicity. 
\end{lemma} 

\begin{proof} 
We can filter $ V(\nu) = V_0 \supset V_1 \supset \cdots V_k \supset 0 $ as a $ U_q(\mathfrak{b}_+) $-module by one-dimensional 
quotients $ \mathbb{C}_\gamma $ where $ U_q(\mathfrak{b}_+) $ acts on $ \mathbb{C}_\gamma $ by the character $ \chi_\gamma $. 
Using the above observations, tensoring this 
filtration with $ M(\lambda) $ yields the claim. 
\end{proof}   

Recall that if $ M $ is an integrable $ U_q(\mathfrak{g}) $-module and $ \nu \in \weights^+ $ we write $ [M: V(\nu)] $ for the multiplicity of $ V(\nu) $ in $ M $. 

\begin{prop} \label{principalseriesmultiplicity}
Let $ \lambda, \eta \in \mathfrak{h}^*_q $ and $ \nu \in \weights^+ $. The the following properties hold. 
\begin{bnum} 
\item[a)] We always have 
$$
[F((M(\lambda) \otimes M(\eta))^*): V(\nu)] = \dim(V(\nu)_{\eta - \lambda}).
$$
\item[b)] If $ M(\eta) $ is simple then 
$$
[F \Hom(M(\lambda), M(\eta)): V(\nu)] = \dim(V(\nu)_{\eta - \lambda}). 
$$
\item[c)] If $ M(\lambda) $ is projective then 
$$
[F \Hom(M(\lambda), M(\eta)): V(\nu)] = \dim(V(\nu)_{\eta - \lambda}). 
$$
\end{bnum} 
\end{prop} 

\begin{proof}  
$ a) $ Recall that the appropriate $U_q(\mathfrak{g})$-module structure on $ (M(\lambda) \otimes M(\eta))^* $ is given by
\[
(X \rightarrow \varphi)(m\otimes n) = \varphi(\hat{S}(X_{(2)})\cdot m \otimes \tau (X_{(1)})\cdot n)
\]
for $X\in U_q(\mathfrak{g})$, $\varphi\in(M(\lambda) \otimes M(\eta))^*$ and $m\otimes n \in M(\lambda)\otimes M(\eta)$. This can be rewritten as
\[
(X \rightarrow \varphi)(m\otimes n) = \varphi(\hat{S}(X)\cdot(m\otimes n)),
\]
if we define the $U_q(\mathfrak{g})$-action on $M(\lambda) \otimes M(\eta)$ by 
\[
X\cdot(m\otimes n) = X_{(1)}\cdot m \otimes \hat{S}\tau(X_{(2)})\cdot n.
\]
Similarly, let us define a $U_q(\mathfrak{g})$-action on $U_q(\mathfrak{g}) \otimes U_q(\mathfrak{g})$ by
\[
X\cdot(Y\otimes Z) = X_{(1)}Y \otimes \hat{S}\tau(X_{(2)})Z
\]
for $X,Y,Z\in U_q(\mathfrak{g})$. Using the triangular decomposition, see Theorem \ref{uqgtriangular}, one can check that this action induces an isomorphism 
\begin{align*}
U_q(\mathfrak{g}) \otimes_{U_q(\mathfrak{h})} (U_q(\mathfrak{b}_+) \otimes U_q(\mathfrak{b}_+))
&\cong U_q(\mathfrak{g}) \otimes U_q(\mathfrak{g}) \\
X \otimes(Y\otimes Z) &\mapsto X \cdot (Y\otimes Z)
\end{align*}
of left $U_q(\mathfrak{g})$-modules.  

We obtain isomorphisms of left $U_q(\mathfrak{g})$-modules as follows
\begin{align*}
(M(\lambda) \otimes M(\eta))^*
&\cong \Hom_{U_q(\mathfrak{b}_+) \otimes U_q(\mathfrak{b}_+)} (U_q(\mathfrak{g}) \otimes U_q(\mathfrak{g}), (\CC_\lambda \otimes \CC_\eta)^* ) \\
&\cong \Hom_{U_q(\mathfrak{b}_+) \otimes U_q(\mathfrak{b}_+)} 
(U_q(\mathfrak{g}) \otimes_{U_q(\mathfrak{h})} (U_q(\mathfrak{b}_+) \otimes U_q(\mathfrak{b}_+))), (\CC_\lambda \otimes \CC_\eta)^* ) \\
&\cong \Hom_{U_q(\mathfrak{h})} ( U_q(\mathfrak{g}) , (\CC_\lambda \otimes \CC_\eta)^*),
\end{align*}
where the $\Hom$ spaces are spaces of morphisms of right modules.  In the last line, the right $U_q(\mathfrak{h})$-action on $(\CC_\lambda \otimes \CC_\eta)^*$ is given by
\[
((\phi\otimes\psi)\cdot H)(v\otimes w) = (\phi\otimes\psi) (\chi_\lambda(H_{(1)})v \otimes \chi_\eta(\hat{S}\tau(H_{(2)}))w) 
= \chi_{\lambda-\eta}(H) (\phi\otimes\psi)(v\otimes w),
\]
so that the second factor is isomorphic to $\CC_{\lambda-\mu}^*$ as a right $U_q(\mathfrak{h})$-module.  Thus
\[
(M(\lambda) \otimes M(\eta))^* \cong ( U_q(\mathfrak{g}) \otimes_{U_q(\mathfrak{h})} \CC_{\lambda-\mu})^*.
\]

We now have 
\begin{align*}
\Hom_{U_q(\mathfrak{g})}(V(\nu), & (U_q(\mathfrak{g}) \otimes_{U_q(\mathfrak{h})} \mathbb{C}_{\lambda - \eta})^*) \\
&\cong \Hom_{U_q(\mathfrak{g})}(V(\nu) \otimes (U_q(\mathfrak{g}) \otimes_{U_q(\mathfrak{h})} \mathbb{C}_{\lambda - \eta}), \mathbb{C}) \\
&\cong \Hom_{U_q(\mathfrak{g})}((U_q(\mathfrak{g}) \otimes_{U_q(\mathfrak{h})} \mathbb{C}_{\lambda - \eta}) \otimes V(\nu), \mathbb{C}) \\
&\cong \Hom_{U_q(\mathfrak{g})}(U_q(\mathfrak{g}) \otimes_{U_q(\mathfrak{h})} \mathbb{C}_{\lambda - \eta}, V(\nu)^*) \\
&\cong \Hom_{U_q(\mathfrak{h})}(\mathbb{C}_{\lambda - \eta}, V(\nu)^*) \\
&\cong  \Hom_{U_q(\mathfrak{h})}(\mathbb{C}_{\eta - \lambda}, V(\nu)) .
\end{align*}
This yields the first claim.

$ b) $ If $ M(\eta) $ is simple we have $ M(\eta)^\vee \cong V(\eta) \cong M(\eta) $. Hence Lemma \ref{locallyfiniteduallemma} yields
$$ 
F((M(\lambda) \otimes M(\eta))^*) \cong F((M(\lambda) \otimes M(\eta)^\vee)^*) \cong 
F \Hom(M(\lambda), M(\eta)),  
$$ 
and part $ a) $ implies $ [F \Hom(M(\lambda), M(\eta)): V(\nu)] = \dim(V(\nu))_{\eta - \lambda} $ for all $ \nu \in \weights^+ $. 

$ c) $ We have
\begin{align*}
[F \Hom(M(\lambda), & M(\eta)): V(\nu)] \\
&= \dim \Hom_{U_q(\mathfrak{g})}(V(\nu), \Hom(M(\lambda), M(\eta))) \\
&= \dim \Hom_{U_q(\mathfrak{g})}(V(\nu) \otimes M(\lambda), M(\eta)) \\
&= \dim \Hom_{U_q(\mathfrak{g})}(M(\lambda), V(\nu)^* \otimes M(\eta)) .
\end{align*}

By Lemma \ref{Vermatensorfiltration}, the module $V(\nu)^*\otimes M(\eta)$ admits a finite length filtration with quotients isomorphic to $M(\eta+\gamma)$, 
where $\gamma$ runs over all weights of $V(\nu)^*$ counted with multiplicity.  If $M(\lambda)$ is projective, this yields a corresponding filtration 
of $\Hom_{U_q(\mathfrak{g})}(M(\lambda), V(\nu)^* \otimes M(\eta))$, and therefore
\begin{align*}
[F \Hom(M(\lambda), & M(\eta)): V(\nu)] 
&= \sum_{\gamma \in \weights(V(\nu)^*)} \!\! \dim \Hom_{U_q(\mathfrak{g})}(M(\lambda), M(\eta + \gamma)),
\end{align*}
where $\weights(V(\nu)^*)$ denotes the weights of $V(\nu)^*$ counted with multiplicity.
Using Proposition \ref{Vermaprojective} we have 
\[
\dim \Hom_{U_q(\mathfrak{g})}(M(\lambda), M(\eta + \gamma)) = \delta_{\lambda, \eta+\gamma} .
\]
We therefore conclude that
\[
[F \Hom(M(\lambda), M(\eta)): V(\nu)] = \dim(V(\nu)^*_{\lambda-\eta}) = \dim(V(\nu)_{\eta - \lambda}). 
\]
\end{proof}

Note that Proposition \ref{principalseriesmultiplicity} implies in particular 
$$
[F((M(\lambda) \otimes M(\eta))^*): V(\nu)] = 0
$$
for all $ \nu \in \weights^+ $ if $ \eta - \lambda $ is not contained in $ \weights $. That is, we have $ F((M(\lambda) \otimes M(\eta))^*) = 0 $ 
in this case.

\subsubsection{Hilbert-Poincar\'e series for the locally finite part of $ U_q(\mathfrak{g}) $}

In this subsection we collect further auxiliary considerations needed for the description of the PRV determinant. We follow closely the 
treatment in \cite{Josephbook}.  

Let us first recall that the height of a weight $ \mu \in \weights $ is defined by 
$ ht(\mu) = c_1 + \cdots + c_N \in \frac{1}{2} \mathbb{Z} $ if $ \mu = c_1 \alpha_1 + \cdots + c_N \alpha_N $,
\label{nom:ht-weights2}%
and that $ ht(\mu) \geq 0 $ 
for $ \mu \in \weights^+ $.  
For our purposes below it will be notationally convenient to work with the degree of $ \mu $ defined by $ \deg(\mu) = 2 ht(\mu) = ht(2 \mu) $ 
instead. 
The advantage of $ \deg $ over $ ht $ is that the former takes only integral values for $ \mu \in \weights $. 

We define an $ \mathbb{N}_0 $-filtration of the vector space $ FU_q(\mathfrak{g}) $ by 
$$
\F^k(FU_q(\mathfrak{g})) = \bigoplus_{\substack{\mu \in \weights^+ \\ \deg(\mu) \leq k}} U_q(\mathfrak{g}) \rightarrow K_{2\mu} 
$$
for $ k \in \mathbb{N}_0 $, and refer to this as the degree filtration of $ FU_q(\mathfrak{g}) $. In fact, this is essentially the 
filtration of $ \Poly(G_q) $ used in the proof of Theorem \ref{separation}, transported to $ FU_q(\mathfrak{g}) $ via the isomorphism from 
Proposition \ref{adjointdualityplus}.
We note that the degree filtration is not compatible with the algebra structure of $ FU_q(\mathfrak{g}) $ in 
an obvious way. 

The associated graded components of the degree filtration are 
$$
gr_k(FU_q(\mathfrak{g})) = \F^k(FU_q(\mathfrak{g}))/\F^{k - 1}(FU_q(\mathfrak{g})), 
$$
observing that $ \F^{-1}(FU_q(\mathfrak{g})) = 0 $. We note that the graded components $ gr_k(FU_q(\mathfrak{g})) $ are finite dimensional 
for all $ k \in \mathbb{N}_0 $. Using the isomorphism $ FU_q(\mathfrak{g}) \cong \mathbb{H} \otimes ZU_q(\mathfrak{g}) $ 
in Theorem \ref{separation} we obtain an induced $ \mathbb{N}_0 $-filtration of the vector 
space $ \mathbb{H} \cong \mathbb{H} \otimes 1 \subset FU_q(\mathfrak{g}) $. We shall write 
$$ 
\mathbb{H}_k = (\F^k(FU_q(\mathfrak{g})) \cap \mathbb{H})/(\F^{k - 1}(FU_q(\mathfrak{g})) \cap \mathbb{H})
$$ 
for the corresponding graded components. 

As for any filtration with finite dimensional graded components, we may define associated Hilbert-Poincar\'e series as the 
generating function of the corresponding graded dimensions as follows. 

\begin{definition} \label{defhilbertpoincare}
Let $ \mu \in \weights^+ $. With the notation as above, the $ \mu $-Hilbert-Poincar\'e series of $ FU_q(\mathfrak{g}) $ and $ \mathbb{H} $ are 
the formal power series defined by 
$$
h^{FU_q(\mathfrak{g})}_\mu(z) = \sum_{k = 0}^\infty [gr_k(FU_q(\mathfrak{g})): V(\mu)] z^k, 
\qquad h^\mathbb{H}_\mu(z) = \sum_{k = 0}^\infty [\mathbb{H}_k: V(\mu)] z^k, 
$$
respectively. 
\end{definition}

Recall from the proof of Theorem \ref{adjointduality} that $U_q(\lie{g})\rightarrow K_{2\lambda}$ is isomorphic to $\End(V(-w_0\lambda))^*
\cong (V(\lambda)^* \otimes V(\lambda))^*
\cong V(\lambda)\otimes V(\lambda)^* 
\cong \End(V(\lambda))$ as a $U_q(\lie{g})$-module,
where we are using the fact that $V\otimes W \cong W\otimes V$ for any finite dimensional $U_q(\lie{g})$-modules $V,W$.
Therefore we can also write
\[
h^{FU_q(\mathfrak{g})}_\mu(z) = \sum_{\lambda\in\weights^+} [\End(V(\lambda)):V(\mu)] z^{\deg(\lambda)}.
\]

In the remainder of this subsection we shall derive a formula for the formal derivative of $ h^\mathbb{H}_\mu(z) $ at $ z = 1 $, presented 
in Proposition \ref{PRVdegreeformula} further below. 

Recall that if $ M $ is in category $ \O $ and $ M = M_0 \supset M_1 \supset \cdots \supset M_n = 0 $ a filtration with 
simple subquotients as in Lemma \ref{lemascendinghighestweightmodules}, then $ [M: V(\mu)] $ is the number of subquotients $ M_i/M_{i + 1} $ 
isomorphic to $ V(\mu) $. 

\begin{lemma} \label{PRVdegreeformulahelphelp}
Let $ \mu, \lambda \in \weights^+ $. Then 
$$
[\End(V(\lambda)): V(\mu)] = \sum_{v, w \in W} (-1)^{l(v) + l(w)} P(v.\lambda - w.\lambda + \mu), 
$$
where $ P $ denotes Kostant's partition function. 
\end{lemma}

\begin{proof}
If $ V(-\lambda) \cong V(\lambda)^* $ denotes the integrable simple $ U_q(\mathfrak{g}) $-module with lowest weight $ -\lambda $, we 
have $ \End(V(\lambda)) \cong V(\lambda) \otimes V(-\lambda) \cong V(-\lambda) \otimes V(\lambda) $. 
Using Lemma \ref{Vermatensorfiltration} and Proposition \ref{charactersimple} we therefore obtain 
\begin{align*}
[\End(V(\lambda)): V(\mu)] &= [V(-\lambda) \otimes V(\lambda): V(\mu)] \\
&= \sum_{w \in W} (-1)^{l(w)} [V(-\lambda) \otimes M(w . \lambda): V(\mu)] \\
&= \sum_{w \in W} \sum_{\nu \in \weights(V(-\lambda))} (-1)^{l(w)} [M(w . \lambda + \nu): V(\mu)],
\end{align*}
where in the last line the sum is over the weights $\nu$ of $V(-\lambda)$ counted with multiplicities.
Since $ \mu \in \weights^+ $ we have $ [M(w . \lambda + \nu): V(\mu)] = 1 $ if $ \nu = \mu - w.\lambda $ and $ [M(w . \lambda + \nu): V(\mu)] = 0 $ 
otherwise by Theorem \ref{thmfdintegrable}, so that we arrive at 
\begin{align*}
[\End(V(\lambda)): V(\mu)] &= \sum_{w \in W} (-1)^{l(w)} \dim(V(-\lambda)_{\mu - w . \lambda}) \\
&= \sum_{w \in W} (-1)^{l(w)} \dim(V(\lambda)_{w . \lambda - \mu}) \\
&= \sum_{v, w \in W} (-1)^{l(v) + l(w)} P(v.\lambda - w.\lambda + \mu) 
\end{align*}
as claimed, using that $ \dim(V(-\lambda)_\nu) = \dim(V(\lambda)_{-\nu}) $ for all weights $ \nu $ and Proposition \ref{charactersimple} again. 
\end{proof}

Let us write the group ring $ \mathbb{C}[\weights] $ of the abelian group $ \weights $ additively, with basis $ e^\lambda $
\label{nom:e_lambda}%
for $ \lambda \in \weights $. 

\begin{definition} \label{defQ}
	We define 
	\begin{align*}
	Q^{FU_q(\mathfrak{g})}(z) &= \sum_{w \in W} (-1)^{l(w)} e^{w . 0} \prod_{j = 1}^N \frac{1}{1 - z^{\deg(\varpi_j)} e^{w \varpi_j - \varpi_j}}, \\
	Q^\mathbb{H}(z) &= \sum_{w \in W} (-1)^{l(w)} e^{w . 0} \prod_{j = 1}^N \frac{1 - z^{\deg(\varpi_j)}}{1 - z^{\deg(\varpi_j)} e^{w \varpi_j - \varpi_j}}, 
	\end{align*}
	both viewed as elements of the function field $ \mathbb{F}(z) $, where $ \mathbb{F} $ is the field of quotients of the integral domain $ \mathbb{C}[\weights] $. 
\end{definition} 

We note that expanding the geometric series in Definition \ref{defQ} we can write 
\begin{align*}
Q^{FU_q(\mathfrak{g})}(z) &= \sum_{\lambda \in \weights^+} \sum_{w \in W} (-1)^{l(w)} e^{w . \lambda - \lambda} z^{\deg(\lambda)}, \\
Q^\mathbb{H}(z) &= \sum_{\lambda \in \weights^+} \sum_{w \in W} (-1)^{l(w)} e^{w . \lambda - \lambda} z^{\deg(\lambda)} \prod_{j = 1}^N (1 - z^{\deg(\varpi_j)}),  
\end{align*}
as formal power series with coefficients in $ \mathbb{C}[\weights] $.

In the next lemma, we will need to consider power series with coefficients belonging to infinite formal series in the $e^\lambda$ of the type introduced in Section \ref{secdefo}.  That is, we will be considering formal sums of the form
\[
  f(z) = \sum_{n=0}^\infty \sum_{\lambda\in\lie{h}_q^*} f_n(\lambda) e^\lambda z^n,
\]
where for each fixed $n$ the support of $f_n:\lie{h}_q^*\to\mathbb{Z}$ is contained in a finite union of sets of the form $\nu-\roots^+$ for $\nu\in\lie{h}_q^*$.  As in Section \ref{secdefo}, such series admit a well-defined convolution product.

\begin{lemma} \label{PRVdegreeformulahelp1}
Let $ \mu \in \weights^+ $ and $ h^{FU_q(\mathfrak{g})}_\mu(z), h^\mathbb{H}_\mu(z) $ be the Poincar\'e-Hilbert series for $ FU_q(\mathfrak{g}) $ and $ \mathbb{H} $, respectively. 
Then 
\begin{bnum}
\item[a)] $ h^{FU_q(\mathfrak{g})}_\mu(z) $ is the coefficient of $ e^0 $ in $ \ch(V(\mu)) Q^{FU_q(\mathfrak{g})}(z) $. 
\item[b)] $ h^\mathbb{H}_\mu(z) $ is the coefficient of $ e^0 $ in $ \ch(V(\mu)) Q^\mathbb{H}(z) $. 
\end{bnum}
Here $ Q^{FU_q(\mathfrak{g})}(z) $ and $ Q^\mathbb{H}(z) $ are viewed as formal power series. 
\end{lemma}

\begin{proof}
	$ a) $ 
	Recall from Section \ref{secdefo} that we defined the formal series
  \[
	 p = \ch(M(0)) = \sum_{\nu\in\weights} P(\nu)e^{-\nu}.
	\]
  Using Proposition \ref{charactersimple} we can write
  \begin{align*}
   \ch(V(\mu))Q^{FU_q(\lie{g})}(z) 
    &= \sum_{w\in W} (-1)^{l(w)} p e^{w.\mu} Q^{FU_q(\lie{g})}(z) \\
    &= \sum_{\lambda\in\weights^+} \sum_{v,w\in W} (-1)^{l(v)+l(w)} p e^{w.\mu + v.\lambda - \lambda} z^{\deg(\lambda)}\\
    &= \sum_{\lambda\in\weights^+} \sum_{w \in W} (-1)^{l(w)} \ch(V(\lambda)) e^{w(\mu+\rho)-\rho-\lambda} z^{\deg(\lambda)}.
  \end{align*}
  For each $n\in\mathbb{N}_0$, the coefficient of $z^n$ is the above series is a finite formal sum, so we may consider the action of the Weyl group  upon it, given by $w(e^\lambda) = e^{w\lambda}$.  Moreover, the $e^0$ term in each summand is invariant under the Weyl group action, so the coefficient of $e^0$ in the above sum is the same as the coefficient of $e^0$ in
  \begin{align*}
  \sum_{\lambda\in\weights^+} & \sum_{w \in W} (-1)^{l(w)} \ch(V(\lambda)) e^{\mu+\rho-w^{-1}(\rho+\lambda)} z^{\deg(\lambda)} \\
  &=  \sum_{\lambda\in\weights^+} \sum_{w \in W} (-1)^{l(w)} \ch(V(\lambda)) e^{\mu-w^{-1}.\lambda} z^{\deg(\lambda)} \\
  &= \sum_{\lambda\in\weights^+} \sum_{v,w\in W} (-1)^{l(v)+l(w)} p e^{v.\lambda} e^{\mu- w^{-1}.\lambda} z^{\deg(\lambda)} \\
  &= \sum_{\nu\in\weights} \sum_{\lambda\in\weights^+} \sum_{v,w\in W} (-1)^{l(v)+l(w)} P(\nu) e^{v.\lambda - w^{-1}.\lambda + \mu-\nu} z^{\deg(\lambda)},
  \end{align*}
  where in the second equality we have again used Proposition \ref{charactersimple}.
  Replacing $w$ by $w^{-1}$ in the sum and using Lemma \ref{PRVdegreeformulahelphelp}, we calculate the coefficient of $e^0$ to be
	\begin{align*}
	 \sum_{\lambda\in\weights^+}  \sum_{v,w\in W}  (-1)^{l(v)+l(w)}  &P(v.\lambda-w.\lambda+\mu) z^{\deg(\lambda)} \\
	   & = \sum_{\lambda\in\weights^+} [\End(V(\lambda)):V(\mu)] z^{\deg(\lambda)}.
	\end{align*}
	From the formula immediately after Definition \ref{defhilbertpoincare}, this is $h^{FU_q(\mathfrak{g})}_\mu(z)$.
	
	$ b) $ Using Theorem \ref{uqgcenter} it is easy to see that the centre $ ZU_q(\mathfrak{g}) $ has Hilbert-Poincar\'e 
	series 
	\[
	 h^{ZU_q(\lie{g})}(z) = \prod_{j = 1}^N \left(\sum_{k=0}^\infty z^{k\deg(\varpi_j)}\right) 
	\]
	with respect to the degree filtration of $ FU_q(\mathfrak{g}) $, and therefore
	\[
	 h^{ZU_q(\lie{g})}(z) \prod_{j = 1}^N (1 - z^{\deg(\varpi_j)}) = e^0.
	\]
	Due to Theorem \ref{separation} we therefore get 
	$$
	h^{\mathbb{H}}_\mu(z) = h^{\mathbb{H}}_\mu(z) h^{ZU_q(\lie{g})}(z) \prod_{j = 1}^N (1 - z^{\deg(\varpi_j)}) = h^{FU_q(\mathfrak{g})}_\mu(z) \prod_{j = 1}^N (1 - z^{\deg(\varpi_j)}).  
	$$
  Hence the assertion follows from $ a) $. 
\end{proof}

\begin{lemma} \label{PRVdegreeformulahelp2}
The derivative of the rational function $ Q^\mathbb{H}(z) $ defined above 
satisfies
$$
(\partial_z Q^\mathbb{H})(1) = \sum_{j = 1}^N \frac{e^{-\alpha_j}}{1 - e^{-\alpha_j}} \deg(\varpi_j)
= \sum_{j = 1}^N \sum_{n = 1}^\infty e^{-n\alpha_j} \deg(\varpi_j). 
$$
\end{lemma} 

\begin{proof}
For every $ \lambda \in \weights^+ $ the stabilizer over $ \lambda $ with respect to the action of $ W $ is generated by 
the simple reflections it contains, see Theorem 1.12 (a) of \cite{HumphreysReflectiongroups}. 
In particular, the stabilizer of $ \varpi_j $ is the subgroup generated by all $ s_i $ for $ i \neq j $. 
Hence for each $ w \in W $ with $ l(w) > 1 $ 
there are at least two fundamental weights $\varpi_j$ for which $w\varpi_j \neq \varpi_j$.  As a consequence, the product 
	\[
	 \prod_{j=1}^N \frac{1-z^{\deg(\varpi_j)}}{1-z^{\deg{\varpi_j}}e^{w\varpi_j-\varpi_j}}
  \]
vanishes at $z=1$ to at least order $2$ when $l(w)>1$.

It follows that only the summands 
with $ l(w) = 1 $ in the definition of $ Q^\mathbb{H}(z) $ contribute to $ (\partial_z Q^\mathbb{H})(1) $. 
Using $ s_i . 0 = - \alpha_i $ and $ s_i \varpi_j = \varpi_j -\delta_{ij} \alpha_i $ we therefore obtain 
$$ 
(\partial_z Q^\mathbb{H})(1) = \sum_{j = 1}^N \frac{e^{-\alpha_j}}{1 - e^{-\alpha_j}} \deg(\varpi_j),
$$
as claimed.
Expanding the geometric series on the right hand side of this formula yields the second equality. 
\end{proof}

\begin{lemma} \label{PRVdegreeformulahelp3}
Let $ \mu \in \weights^+ $ and $ n \in \mathbb{N} $. Then we have 
$$
\sum_{j = 1}^N \dim(V(\mu)_{n \alpha_j}) \deg(\varpi_j) = \frac{1}{2} \sum_{\alpha \in {\bf \Delta}^+} \dim(V(\mu)_{n \alpha}) \deg(\alpha). 
$$
\end{lemma}

\begin{proof}
We begin with the case where $\lie{g}$ is a simple Lie algebra, so that the root system ${\bf\Delta}$ is irreducible.  
In this case, any two roots of the same length are conjugate under the Weyl group action, see Lemma C in Section 10.4 of \cite{HumphreysLie}, so the decomposition of ${\bf\Delta}$ into $W$-orbits is given by ${\bf\Delta} = {\bf \Delta}_s \cup {\bf \Delta}_l $, where ${\bf\Delta}_s$ is the set of short roots, meaning roots of minimal length, and ${\bf\Delta}_l$ is its complement, the possibly empty set of long roots.
We shall show that the contributions of all roots of a given length on both sides of the claimed formula agree. 

More precisely, denote by $ {\bf \Delta}^+_s $ and $ {\bf \Delta}^+_l $ the sets of short and long roots in $ {\bf \Delta}^+ $, respectively, 
and let $ \Sigma_s = {\bf \Delta}^+_s \cap \Sigma $ and $ \Sigma_l = {\bf \Delta}^+_l \cap \Sigma $ be the sets of short simple roots and long simple roots. 
Since the dimension of weight spaces is constant on Weyl group orbits we have $ \dim(V(\mu)_{n \alpha_j}) = \dim(V(\mu)_{n \alpha}) $ 
for any $ \alpha_j \in \Sigma_s $ and $ \alpha \in {\bf \Delta}^+_s $, or $ \alpha_j \in \Sigma_l $ and $ \alpha \in {\bf \Delta}^+_l $. 
Hence it suffices to show 
\begin{align*}
\sum_{\alpha_j \in \Sigma_s} \varpi_j &= \frac{1}{2} \sum_{\alpha \in {\bf \Delta}_s^+} \alpha, 
\qquad \sum_{\alpha_j \in \Sigma_l} \varpi_j = \frac{1}{2} \sum_{\alpha \in {\bf \Delta}_l^+} \alpha.  
\end{align*}
For the first equality, write $ \rho_s $ for the right hand side and notice that a simple reflection $ s_j $ satisfies $ s_j \rho_s = \rho_s - \alpha_j $ 
if $ \alpha_j \in \Sigma_s $, and $ s_j \rho_s = \rho_s $ otherwise. Hence $ (\alpha_j^\vee, \rho_s) = 1 $ for $ \alpha_j \in \Sigma_s $ and 
$ (\alpha_j^\vee, \rho_s) = 0 $ otherwise. However, the same relations characterize the sum on the left hand side. 

The second equality can be proved in the same way, or by using the first formula and the fact that the sum of all fundamental weights equals the half-sum
of all positive roots. 

Finally, if ${\bf\Delta}$ is not irreducible, then the two sums in the Lemma can be decomposed into sums corresponding to each of the irreducible components, and the result follows.
\end{proof}

Assembling the above results we obtain the following formula. 

\begin{prop} \label{PRVdegreeformula}
Let $ \mu \in \weights^+ $. Then we have 
$$
(\partial_z h^\mathbb{H}_\mu)(1) = \frac{1}{2} \sum_{n = 1}^\infty \sum_{\alpha \in {\bf \Delta}^+} \dim(V(\mu)_{n \alpha}) \deg(\alpha). 
$$
\end{prop}

\begin{proof}
According to Lemma \ref{PRVdegreeformulahelp1} we know that $ (\partial_z h^\mathbb{H}_\mu)(1) $ equals the coefficient of $ e^0 $ 
in $ \ch(V(\mu)) (\partial_z Q^\mathbb{H})(1) $. Due to Lemma \ref{PRVdegreeformulahelp2} this coefficient is given by 
$$
\sum_{j = 1}^N \sum_{n = 1}^\infty \dim(V(\mu)_{n \alpha_j}) \deg(\varpi_j). 
$$
Hence the claim follows from Lemma \ref{PRVdegreeformulahelp3}. 
\end{proof}

\subsubsection{The quantum PRV determinant} 

According to Theorem \ref{separation} there exists a $ U_q(\mathfrak{g}) $-invariant linear subspace $ \mathbb{H} \subset FU_q(\mathfrak{g}) $ 
such that the multiplication map $ \mathbb{H} \otimes ZU_q(\mathfrak{g}) \rightarrow FU_q(\mathfrak{g}) $ is a $ U_q(\mathfrak{g}) $-linear isomorphism. 
For $ \mu \in \weights^+ $ let 
$$ 
\mathbb{H}^\mu \cong \bigoplus_{j = 1}^m V(\mu) 
$$ 
be the $ \mu $-isotypical component of $ \mathbb{H} $. It is also shown in Theorem \ref{separation} that the multiplicity $ m $ is given by 
$ m = [\mathbb{H}: V(\mu)] = \dim(V(\mu)_0) $. Choose a basis $ v_1, \dots, v_m $ of the zero weight space of $ V(\mu) $, 
and write $ v_{ij} \in \mathbb{H}^\mu $ for the vector $ v_j $ in the $ i $-th copy of $ V(\mu) $ with respect to the above identification. 
According to the remarks after Theorem \ref{separation}, each of the copies of $ V(\mu) $ in $ \mathbb{H}^\mu $ can be assumed to 
be contained in a subspace of the form $ U_q(\mathfrak{g}) \rightarrow K_{2 \nu_i} $ for some $ \nu_1, \dots, \nu_m \in \weights^+ $.

If we identify $ \mathbb{H}^\mu $ with $ \mathbb{H}^\mu \otimes 1 \subset FU_q(\mathfrak{g}) $ via the above isomorphism, 
then the PRV determinant of $ U_q(\mathfrak{g}) $ associated with $ \mu $ is defined as the determinant of the matrix $ \P_\mu = (\P(v_{ij})) $, 
where $ \P: U_q(\mathfrak{g}) \rightarrow U_q(\mathfrak{h}) $ denotes the Harish-Chandra map. 
Note that this determinant, which we denote by $\det(\P_\mu)\in U_q(\lie{h})$,
\nomenclature[o$det(\P_\mu)$]{$\det(\P_\mu)$}{PRV determinant}%
is independent of the choice of basis $ v_1, \dots, v_m $ up to an invertible scalar in $ \mathbb{C}$.
Our aim in this subsection is to compute $ \det(\P_\mu) $, at least up to multiplication by an invertible element of $U_q(\lie{h})$. We start with a result on annihilators, see Lemma 7.1.10 in \cite{Josephbook}. 

\begin{lemma} \label{annihilatorlemma}
Let $ H \subset FU_q(\mathfrak{g}) $ be a submodule with respect to the adjoint action, and let $ \lambda \in \mathfrak{h}_q^* $. 
Then $ \chi_\lambda(\P(H)) = 0 $ iff $ H \subset \ann_{FU_q(\mathfrak{g})}(V(\lambda)) $. 
\end{lemma} 

\begin{proof}  Recall that $ \chi_\lambda $ denotes the character of $ U_q(\mathfrak{h}) $ given by $ \chi_\lambda(K_\beta) = q^{(\lambda, \beta)} $. 
Let us write $ N \subset V(\lambda) $ for the linear span of all weight spaces in $ V(\lambda) $ 
apart from the highest weight space $ V(\lambda)_\lambda = \mathbb{C} v_\lambda $.

We claim that $ \chi_\lambda(\P(H)) = 0 $ iff $ H \cdot v_\lambda \subset N $. Indeed, if $ H \cdot v_\lambda \subset N $ then the projection $ \P(H) $ 
of $ H $ onto $ U_q(\mathfrak{h}) = 1 \otimes U_q(\mathfrak{h}) \otimes 1 \subset U_q(\mathfrak{n}_-) \otimes U_q(\mathfrak{h}) \otimes U_q(\mathfrak{n}_+) $ 
must satisfy $ \chi_\lambda(\P(H)) = 0 $ since it acts on $ v_\lambda $ by the character $ \chi_\lambda $. 
Conversely, if $ \chi_\lambda(\P(H)) = 0 $ then only terms in $ H \cdot v_\lambda $ which lower weights survive, 
which means $ H \cdot v_\lambda \subset N $.

Assume $ \chi_\lambda(\P(H)) = 0 $ and let $ X \in H $ be a weight vector of weight $ \mu $ with respect to the adjoint action. Then 
$$
X F_j = q^{(\mu,\alpha_j)} K_j^{-1} X K_j F_j = q^{(\mu,\alpha_j)} (F_j X - (F_j \rightarrow X))
$$
by definition of the adjoint action. We conclude $ H U_q(\mathfrak{n}_-) = U_q(\mathfrak{n}_-) H $. 
Since $ \chi_\lambda(\P(H)) = 0 $ implies $ H \cdot v_\lambda \subset N $ this implies
$$ 
H \cdot V(\lambda) = H U_q(\mathfrak{n}_-) \cdot v_\lambda = U_q(\mathfrak{n}_-) H \cdot v_\lambda \subset N. 
$$
In addition, because $ H $ is $ \ad $-stable we have $ U_q(\mathfrak{g}) H = H U_q(\mathfrak{g}) $ and 
$ U_q(\mathfrak{g}) H \cdot V(\lambda) = H U_q(\mathfrak{g}) \cdot V(\lambda) = H \cdot V(\lambda) $. 
That is, the vector space $ H \cdot V(\lambda) \subset V(\lambda) $ is a $ U_q(\mathfrak{g}) $-submodule. 
Since $ H \cdot V(\lambda) $ is contained in $ N $ and $ V(\lambda) $ is simple we conclude $ H \cdot V(\lambda) = 0 $. 
That is, we have $ H \subset \ann_{U_q(\mathfrak{g})}(V(\lambda)) $, and therefore in particular 
$ H \subset \ann_{FU_q(\mathfrak{g})}(V(\lambda)) $. 

Conversely, if $ H \subset \ann_{FU_q(\mathfrak{g})}(V(\lambda)) $ we obtain $ H \cdot v_\lambda = 0 \subset N $, 
and our initial argument implies $ \chi_\lambda(\P(H)) = 0 $. 
\end{proof}  

We note that Lemma \ref{annihilatorlemma} implies in particular that $ \det(\P_\mu) $ is nonzero for any $ \mu \in \weights^+ $. In fact, 
otherwise we would find a nonzero $ \ad $-stable subspace $ H \subset \mathbb{H}^\mu \subset \mathbb{H} $ such that $ \chi_\lambda(\P(H)) = 0 $ 
for all $ \lambda \in \mathfrak{h}^*_q $, and hence $ H $ would act by zero on all irreducible modules $ V(\lambda) $. 
This is impossible due to Theorem \ref{uqgseparation}.

For $ \alpha \in {\bf \Delta}^+ $ and $ n \in \mathbb{N} $ let us define 
$$
\Gamma_{n, \alpha} = \{\lambda \in \mathfrak{h}^*_q \mid q_\alpha^{(\lambda + \rho, \alpha^\vee)} = q_\alpha^n 
\text{ and } q_\beta^{(\lambda + \rho, \beta^\vee)} \notin \pm q_\beta^{\mathbb{Z}} \text{ for all } \beta \in {\bf \Delta}^+, \beta \neq \alpha \}. 
$$
If $ \lambda \in \Gamma_{n, \alpha} $ then this description implies in particular that $ s_\alpha . \lambda = s_{0, \alpha} . \lambda $ is antidominant, see Definition \ref{defdominantantidominant}.  
According to Theorem \ref{Vermacomplexirreducible} it follows that $ M(s_{0, \alpha} . \lambda) $ is simple. 
Moreover the Jantzen sum formula from Theorem \ref{Jantzenfiltration} shows $ M(\lambda)^i = 0 $ for $ i > 1 $ in this case and 
hence $ M(\lambda)^1 = M(s_{0, \alpha} . \lambda) $. In particular $ M(\lambda)/M(s_{0, \alpha} . \lambda) $ is simple. 

\begin{lemma} \label{PRVmultiplicity}
Let $ \lambda \in \Gamma_{n,\alpha} $. Then 
$$ 
[F\End(V(\lambda)): V(\mu)] = \dim(V(\mu)_0) - \dim(V(\mu)_{n \alpha})
$$ 
for all $ \mu \in \weights^+ $. 
\end{lemma}

\begin{proof}  
Since $ q_\beta^{(\lambda + \rho, \beta^\vee)} \notin \pm q_\beta^{-\mathbb{N}} $ for all $ \beta \in {\bf \Delta}^+ $ by the definition of $ \Gamma_{n,\alpha} $,  
it follows from Proposition \ref{Vermaprojective} that $ M(\lambda) $ is projective. 
Moreover $ V(\lambda) = M(\lambda)/M(s_{0, \alpha} . \lambda) $ is a proper quotient of $ M(\lambda) $, so by Lemma \ref{gkspecial} we have
$ d(V(\lambda)) < d(M(\lambda)) = d(M(s_{0, \alpha} . \lambda)) $. 
Therefore Proposition \ref{GKlemma} b) shows that $F\Hom(M(s_{0,\alpha}.\lambda), V(\lambda)) = 0 $, 
and hence the map $ F \End(V(\lambda), V(\lambda)) \rightarrow F \End(M(\lambda), V(\lambda)) $ is an isomorphism. 
Since $ M(\lambda) $ is projective the functor $ F \Hom(M(\lambda), -) $ is exact, so applying this functor to the short exact sequence 
$$
\xymatrix{
0 \ar@{->}[r] & M(s_{0, \alpha} . \lambda) \ar@{->}[r] & M(\lambda) \ar@{->}[r] & V(\lambda)
\ar@{->}[r] & 0 
     }
$$
gives
\begin{align*}
[F \End(V(\lambda), V(\lambda)): V(\mu)] &= [F \Hom(M(\lambda), M(\lambda)): V(\mu)] \\
&\qquad - [F \Hom(M(\lambda), M(s_{0, \alpha} . \lambda)): V(\mu)].  
\end{align*}
Now Lemma \ref{principalseriesmultiplicity} yields the claim, keeping in mind that $\dim(V(\mu)_{-n\alpha}) = \dim(V(\mu)_{n\alpha})$ since $\pm n\alpha$ are in the same orbit of the Weyl group action.
\end{proof}

\begin{lemma} \label{commutationHC}
Let $ \nu \in \weights^+ $. Then 
$$ 
\P(U_q(\mathfrak{g}) \rightarrow K_{2 \nu}) K_{-2 \nu} 
\quad \subset \hspace{-2ex} \sum_{\gamma \in \nu + \weights(V(-w_0 \nu))}
\hspace{-4ex} \mathbb{C} K_{-2 \gamma} 
\quad \subset \quad \mathbb{C}[K_1^{-2}, \dots, K_N^{-2}], 
$$
where $ \weights(V(-w_0 \nu)) $ denotes the set of weights of $ V(-w_0 \nu) \cong V(\nu)^* $. 
\label{nom:weights_of_V}%
\end{lemma} 

\begin{proof} 
Let us first show $ \P(X \rightarrow K_{2 \nu}) K_{-2 \nu} \in \mathbb{C}[K_1^{-2}, \dots, K_N^{-2}] $ 
for any $ X \in U_q(\mathfrak{g}) $. 
For this it suffices to consider a monomial $ X = E_{i_1} \cdots E_{i_r} K_\lambda F_{j_1} \cdots F_{j_s} $ 
with $ 1 \leq i_1, \dots, i_r \leq N, \lambda \in \weights $ and $ 1 \leq j_1, \dots, j_s \leq N $. 
Using the commutation relations in $ U_q(\mathfrak{g}) $ one checks that $ F_{j_1} \cdots F_{j_s} \rightarrow K_{2 \nu} $ 
is a linear combination of monomials of the form $ K_{2 \lambda} F_{j_\sigma(1)} \cdots F_{j_\sigma(s)} $ where $ \sigma $ is a permutation 
of $ 1, \dots, s $. The same is true after applying the adjoint action of $ K_\lambda $. 
In order for $ \P(X \rightarrow K_{2 \nu}) $ to be nonzero it is therefore necessary that $ E_{i_1} \cdots E_{i_r} $ contains exactly as many 
generators $ E_k $ as $ F_{j_1} \cdots F_{j_s} $ contains generators $ F_k $ for all $ 1 \leq k \leq N $. 
Our assertion follows now from 
$$
E_i \rightarrow F_j = [E_i, F_j] K_i^{-1} = \frac{\delta_{ij}}{q_i - q_i^{-1}} (1-K_i^{-2})
$$
for $ 1 \leq i, j  \leq N $ and the fact that the adjoint action is compatible with multiplication. 

In order to prove the Lemma we recall that the isomorphism $ J: FU_q(\mathfrak{g}) \cong \Poly(G_q) $ from Theorem \ref{adjointduality} maps 
$ U_q(\mathfrak{g}) \rightarrow K_{2 \nu} $ onto $ \End(V(-w_0 \nu))^* $. Since $ J $ is compatible with the adjoint action on $ FU_q(\mathfrak{g}) $ 
and the coadjoint action on $ \Poly(G_q) $ it follows that for $ \P(X \rightarrow K_{2\nu}) $ to be nonzero it is necessary that $ E_{i_1} \cdots E_{i_r} $ 
is contained in $ U_q(\mathfrak{n}_+)_\gamma $ for some $ \gamma \in \roots^+ $ such that $ -\nu + \gamma \in \weights(V(-w_0 \nu)) $. 
Combined with the above considerations this yields the claim. 
\end{proof} 

As pointed out at the start of this subsection, the $ j $-th copy of $ V(\mu) $ in $ \mathbb{H}^\mu $ can be assumed to 
be contained in a subspace of the form $ U_q(\mathfrak{g}) \rightarrow K_{2 \nu_j} $
for some $ \nu_1, \dots, \nu_m \in \weights^+ $. 
Hence if we write $ \nu = \nu_1 + \cdots + \nu_m $ then according to Lemma \ref{commutationHC} it follows that 
$ \det(\P_\mu) $ is contained in $ \mathbb{C}[K_1^{-2}, \dots, K_N^{-2}] K_{2 \nu} \subset U_q(\mathfrak{h}) $. 
For technical reasons it will be convenient to remove the factor $ K_{2 \nu} $ from $ \det(\P_\mu) $ in our considerations below, 
and work with $ \Det(\P_\mu) = \det(\P_\mu) K_{-2 \nu} $
instead. Thus 
we have $ \Det(\P_\mu) \in \mathbb{C}[K_1^{-2}, \dots, K_N^{-2}] \subset \mathbb{C}[ \roots] $ by construction. 

Recall from Subsection \ref{sec:weights} that we identify elements $\lambda\in\lie{h}_q^*$ with the associated algebra characters $\chi_\lambda: U_q(\lie{h}) \cong \mathbb{C}[\weights] \to \mathbb{C}$.  We may also restrict these to characters of $\mathbb{C}[\roots]$.  
For $ \lambda \in \mathfrak{h}^*_q $ we denote 
by $ \mathfrak{m}_\lambda \subset \mathbb{C}[K_1^{\pm 1}, \dots, K_N^{\pm 1}] = \mathbb{C}[\roots] $
the kernel 
of $ \chi_\lambda: \mathbb{C}[\roots] \rightarrow \mathbb{C} $. 
Given $ f \in \mathbb{C}[\roots] $ we shall say that $ \lambda $ is a zero of $ f $ of order $ \geq d $ if 
$ f \in \mathfrak{m}_\lambda^d $ for some $ d \in \mathbb{N} $. We also say that $ \lambda $ is a zero of order $ d $ if $ f \in \mathfrak{m}_\lambda^d $ 
and $ f \notin \mathfrak{m}_\lambda^{d + 1} $. 

\begin{lemma} \label{orderzeroPRVhelp}
Let $ \lambda \in \mathfrak{h}^*_q $ and $ \mu \in \weights^+ $. 
Then $ \lambda $ is a zero of $ \Det(\P_\mu) $ of order $ \geq [\ann_{\mathbb{H}}(V(\lambda)): V(\mu)] $. 
\end{lemma}

\begin{proof}  
For each copy $ V $ of $ V(\mu) $ in $ \ann_{\mathbb{H}}(V(\lambda)) $ we have $ \chi_\lambda(\P(V)) = 0 $ by Lemma \ref{annihilatorlemma}. 
In terms of the elements $ v_{ij} \in \mathbb{H}^\mu $ defined at the start of this subsection this means that the rank of the 
matrix $ \chi_\lambda(\P(v_{ij})) $ satisfies $ rk(\chi_\lambda(\P(v_{ij})) \leq m - [\ann_{\mathbb{H}}(V(\lambda)): V(\mu)] $. 
Hence the same is true for the matrix $ \chi_\lambda(\P(u_{ij})) $ where $ u_{ij} = v_{ij} K_{-2 \nu_j} \in \mathbb{C}[K_1^{-2}, \dots K_N^{-2}] $. 
Since $ \Det(\P_\mu) = \det(\P(u_{ij})) $ this yields the assertion. 
\end{proof}  

\begin{lemma} \label{orderzeroPRV}
Let $ \lambda \in \Gamma_{n,\alpha} $ and $ \mu \in \weights^+ $. 
Then $ \lambda $ is a zero of $ \Det(\P_\mu) $ of order $ \geq \dim(V(\mu)_{n \alpha}) $. 
\end{lemma}

\begin{proof} 
According to our above observations we have $ V(\lambda) = M(\lambda)/M(s_{0, \alpha} . \lambda) $. 
Moreover the action of $ FU_q(\mathfrak{g}) $ on $ V(\lambda) = M(\lambda)/M(s_{0, \alpha} . \lambda) $ clearly defines an 
injection of $ FU_q(\mathfrak{g})/\ann_{FU_q(\mathfrak{g})}(V(\lambda)) $ into $ F \Hom(V(\lambda), V(\lambda)) $. Hence 
$$ 
[\mathbb{H}/\ann_{\mathbb{H}}(V(\lambda)): V(\mu)] \leq \dim(V(\mu)_0) - \dim(V(\mu)_{n \alpha}) = [\mathbb{H}: V(\mu)] - \dim(V(\mu)_{n \alpha})
$$ 
by Lemma \ref{PRVmultiplicity}.  We thus obtain 
\begin{align*}
[\ann_\mathbb{H}(V(\lambda)): V(\mu)] &= [\mathbb{H}: V(\mu)] - [\mathbb{H}/\ann_{\mathbb{H}}(V(\lambda)): V(\mu)] \geq \dim(V(\mu)_{n \alpha}). 
\end{align*}
Therefore the assertion follows from Lemma \ref{orderzeroPRVhelp}. 
\end{proof}  

\begin{lemma} \label{PRVdivisor}
Let $ \mu \in \weights^+ $. For each $ \alpha \in {\bf \Delta}^+ $ and $ n \in \mathbb{N} $, the element
$$
(K_\alpha - q_\alpha^{2n - 2 (\rho, \alpha^\vee)} K_{-\alpha})^{\dim(V(\mu)_{n \alpha})} 
$$
divides $ \Det(\P_\mu) $. 
\end{lemma} 

\begin{proof}
We shall show that $ p^{\dim(V(\mu)_{n \alpha})} $ divides $ \Det(\P_\mu) $ where 
$$ 
p = 1 - q_\alpha^{2n - (2\rho, \alpha^\vee)}K_\alpha^{-2} \in \mathbb{C}[2 \roots]. 
$$ 
Note first that in $ \mathbb{C}[\roots] $ we have the factorization $ p = p_+ p_- $ where 
$$
p_\pm = 1 \pm q_\alpha^{n - (\rho, \alpha^\vee)}K_\alpha^{-1}, 
$$
and for $ \lambda \in \Gamma_{n,\alpha} $ the formula 
$$
1 \pm q_\alpha^{n - (\rho, \alpha^\vee)}K_\alpha^{-1} = 1 \pm q_\alpha^{(\lambda, \alpha^\vee)}K_\alpha^{-1} = 1 \pm  q^{(\lambda, \alpha)}K_{-\alpha}
$$
shows that $ \chi_\lambda(p_-) = 0 $ and $ \chi_\lambda(p_+) \neq 0 $. 

We claim that $ p_\pm $ are irreducible elements in $ \mathbb{C}[\roots] $ 
and that $ \lambda $ is a zero of $ p_- $ of order $ 1 $. This is clear if $ \alpha = \alpha_i $ is a simple root. 
For general $ \alpha $ we can use the Weyl group action on $ \mathbb{C}[\roots] $ to transform $ p_\pm $ into a polynomial of the 
form   $ 1 \pm q^{(\lambda, \alpha)} K_{\alpha_i}^{-1}$ for some simple root $ \alpha_i $, which easily yields the claim in this case as well. 

Consider the linear subspace $ N_\alpha = \{\gamma \in \mathfrak{h}^* \mid (\gamma, \alpha) = 0 \} $ 
of $ \mathfrak{h}^* $. The bilinear form $ (\;,\;) $ on $ \mathfrak{h}^* $ induces an orthogonal direct sum decomposition 
$$
\mathfrak{h}^* = N_\alpha \oplus N_\alpha^\bot, 
$$
such that $ \eta \in \mathfrak{h}^* $ corresponds to $ (\eta - \frac{1}{2} (\eta, \alpha^\vee) \alpha, \frac{1}{2} (\eta, \alpha^\vee) \alpha) $. 
We shall write $ \pi: \mathfrak{h}^* \rightarrow N_\alpha $ for the canonical projection, given by 
$$
\pi(\eta) = \eta - \frac{1}{2} (\eta, \alpha^\vee) \alpha. 
$$

Note that when $\eta\in2\roots$ we have $\half(\eta,\alpha^\vee)\in\mathbb{Z}$, so that $\eta-\pi(\eta) \in \mathbb{Z}\alpha$ and $\pi(\eta)\in N_\alpha\cap\roots$.
It follows that if we define $ L_\alpha = \pi(2 \roots) \subset N_\alpha \cap \roots $ then every element of $ \mathbb{C}[2 \roots] $
can be written in the form $ \sum_{j = -l}^l g_{j} K_\alpha^j $ for some $l\geq 0$ and elements $g_j \in \mathbb{C}[L_\alpha]$.  

In particular, we obtain such an expression for $\Det(\P_\mu)$. 
It follows that we can write 
$$
\Det(\P_\mu) =  K_\alpha^{l} \sum_{j=d}^k  f_j p_-^j, 
$$
where each $ f_j $ is contained in $ \mathbb{C}[L_\alpha] \subset \mathbb{C}[\roots] $, and we have chosen $ d\in\mathbb{\NN}_0 $ to be the smallest integer with $ f_d \neq 0 $.  

Consider the canonical embedding $ \mathbb{C}[L_\alpha] \subset \mathbb{C}[\roots] $. 
Since the restriction of the bilinear form $ (\;, \;) $ to $ N_\alpha $ is nondegenerate, the assignment $ \eta \mapsto \chi_\eta $ induces a 
surjection $ \chi: N_\alpha \rightarrow \Hom(L_\alpha, \mathbb{C}^\times) $ 
such that the diagram 
$$
\xymatrix{
\mathfrak{h}^* \ar@{->}[r]^{\!\!\!\!\!\!\!\!\!\!\! \chi} \ar@{->}[d]^\pi & \Hom(\roots, \mathbb{C}^\times) \ar@{->}[d] \\
N_\alpha \ar@{->}[r]^{\!\!\!\!\!\!\!\!\!\!\! \chi} & \Hom(L_\alpha, \mathbb{C}^\times)
     }
$$
is commutative. Here the right hand vertical map is induced by restriction.

For $ n \in \mathbb{N} $ let us define 
\begin{align*}
Z_{n, \alpha} = \{\lambda \in \mathfrak{h}^* \mid (\lambda & +\rho,  \alpha^\vee) = n \\
& \text{ and } (\lambda + \rho, \beta^\vee) \notin \mathbb{Z} +  {\textstyle\half} i\hbar_\beta^{-1}\mathbb{Z} \text{ for all } \beta \in {\bf \Delta}^+, \beta \neq \alpha \}, 
\end{align*}
where we are using the notation $\hbar_\beta = d_\beta\hbar$.
Note that under the surjection $ \mathfrak{h}^* \rightarrow \mathfrak{h}^*_q \cong \Hom(\weights, \mathbb{C}^\times) $ the set $ Z_{n, \alpha} $ gets mapped to $\Gamma_{n, \alpha}$.
The set $Z_{n,\alpha}$ is open and dense in the affine subspace
\[
  \{\lambda \in \mathfrak{h}^* \mid (\lambda +\rho,  \alpha^\vee) = n \},
\]
and it follows that $ \pi(Z_{n, \alpha}) $ is open and dense in $ N_\alpha $. 
Interpreting $ \Hom(L_\alpha, \mathbb{C}^\times) $ as the set of characters of $ \mathbb{C}[L_\alpha] $ it follows easily that we 
find $ \lambda \in Z_{n, \alpha} $ such that $ \chi_\lambda(f_d) $ is nonzero. 
Since the image of $\lambda$ in $\Gamma_{n,\alpha}$ is a zero of $ \Det(\P_\mu) $ of order $ \geq \dim(V(\mu)_{n \alpha}) $ by 
Lemma \ref{orderzeroPRV}, we conclude that $ d \geq \dim(V(\mu)_{n \alpha}) $.

Finally, recall that $ \Det(\P_\mu) $ is contained in $ \mathbb{C}[2 \roots] $, 
and therefore invariant under the automorphism of $ \mathbb{C}[\roots] $ which maps $ K_j $ to $ -K_j $. It follows 
that $ p_+^d $ divides $ \Det(\P_\mu) $ as well. We have thus shown that the polynomial $ p^d $ divides $ \Det(\P_\mu) $, and this finishes the proof. 
\end{proof}  

The following result is essentially the content of Theorem 8.2.10 in \cite{Josephbook}. 

\begin{theorem} \label{PRVformula}
For $ \mu \in \weights^+ $ we have 
$$
\det(\P_\mu) = \prod_{n \in \mathbb{N}} \prod_{\alpha \in {\bf \Delta}^+} 
(K_\alpha - q_\alpha^{2n - 2 (\rho, \alpha^\vee)} K_{-\alpha})^{\dim(V(\mu)_{n \alpha})}, 
$$
up to multiplication by an invertible element of $ U_q(\mathfrak{h}) $. 
\end{theorem} 

\begin{proof} 
Recall that we may choose the decomposition of $ \mathbb{H}^\mu $ in such a way 
that the $ i $-th copy of $ V(\mu) $ is contained in $ U_q(\mathfrak{g}) \rightarrow K_{2 \nu_i} $ for 
some $ \nu_i \in \weights^+ $ for $ i = 1, \dots, m $, and we define $ \nu = \nu_1 + \cdots + \nu_m $. 
Clearly it suffices to prove the assertion for $ \Det(\P_\mu) = \det(\P_\mu) K_{-2\nu} $ instead of $ \det(\P_\mu) $.  

According to Lemma \ref{PRVdivisor} we have 
$$
\Det(\P_\mu) 
= f \prod_{n \in \mathbb{N}} \prod_{\alpha \in {\bf \Delta}^+} (1 - q_\alpha^{2(n - (\rho, \alpha^\vee))}K_\alpha^{-2})^{\dim(V(\mu)_{n \alpha})}
$$
for some element $ f \in \mathbb{C}[K_1^{-2}, \dots, K_N^{-2}] $.  Our strategy is to use degree considerations to show that $ f $ is actually an invertible scalar. 

Let us first use the definition of $ \deg $ to define a filtration on $ \mathbb{C}[K_1^{-2}, \dots, K_N^{-2}]$ whereby, for a linear combination $X = \sum_j c_jK_{-\eta_j}$ with $ \eta_j \in 2 \roots^+ $ and nonzero coefficients  $c_j\in\mathbb{C}$, we set $\deg(X) = \max_j\deg(\eta_j)$. 
By Lemma \ref{commutationHC} each element of $ \P(U_q(\mathfrak{g}) \rightarrow K_{2 \nu_j}) K_{-2 \nu_j} $ is a linear combination of 
monomials $ K_{-2 \gamma} $ with $ \gamma \in \nu_j + \weights(V(-w_0 \nu_j)) $. 
Since the action of $ -w_0 $ on $ \weights $ preserves the set of positive roots it is easy to see that
$ \deg(-2 w_0 \nu_j) = \deg(2 \nu_j) $, 
so that each element $ X \in \P(U_q(\mathfrak{g}) \rightarrow K_{2 \nu_j}) K_{-2 \nu_j} $ satisfies 
$$ 
 \deg(X) \leq \deg(2\nu_j) + \deg(-2 w_0 \nu_j) = 4 \deg(\nu_j). 
$$ 
It follows that $ \deg(\Det(\P_\mu)) \leq 4 \deg(\nu) $.  

If we consider the Hilbert-Poincar\'e series 
$$ 
h^\mathbb{H}_\mu(z) = \sum_{k = 0}^\infty [\mathbb{H}_k: V(\mu)] z^k
$$ 
for $ \mathbb{H} $ from Definition \ref{defhilbertpoincare}, then by the definition of the $ \nu_i $ we have 
$$ 
(\partial_z h^\mathbb{H}_\mu)(1) = \sum_{k = 0}^\infty k [\mathbb{H}_k: V(\mu)] = \deg(\nu_1) + \cdots + \deg(\nu_m) = \deg(\nu). 
$$
According to Proposition \ref{PRVdegreeformula} 
we therefore obtain 
\begin{align*}
  \deg(\nu) &= \frac{1}{2} \sum_{n \in \mathbb{N}} \sum_{\alpha \in {\bf \Delta}^+} \dim(V(\mu)_{n \alpha}) \deg(\alpha).
\end{align*}
This implies
$$
\deg\bigg(\prod_{n \in \mathbb{N}} \prod_{\alpha \in {\bf \Delta}^+}
 (1 - q_\alpha^{2 (n - (\rho, \alpha^\vee))}K_\alpha^{-2})^{\dim(V(\mu)_{n \alpha})}
\bigg)
= 4 \deg(\nu). 
$$

Combining our above considerations yields 
$$ 
4 \deg(\nu) \deg(f) = \deg(\Det(\P_\mu)) \leq 4 \deg(\nu), 
$$
which implies $ \deg(f) = 0 $. We conclude that $ f $ is a scalar, which must be invertible since $ \Det(\P_\mu) $ is nonzero. This finishes the proof. 
\end{proof}  

We note that Joseph establishes a slightly stronger version of the formula for the PRV determinant given in Theorem \ref{PRVformula}. 
Namely, it is shown in Theorem 8.2.10 of \cite{Josephbook} that the formula holds even up to multiplication by a nonzero scalar of the ground field, and not only 
up to an invertible element of $ U_q(\mathfrak{h}) $. We shall not need this stronger assertion in the sequel. 

In fact, our main application of Theorem \ref{PRVformula} is the following key result on annihilators. 

\begin{prop} \label{annHprop}
Let $ \lambda \in \mathfrak{h}^*_q $. Then the Verma module $ M(\lambda) $ is simple iff 
$$ 
\ann_{FU_q(\mathfrak{g})}(V(\lambda)) \cap \mathbb{H} = 0. 
$$ 
\end{prop} 

\begin{proof} 
From our considerations in Lemma \ref{radicalmaximalsubmodule} we know that
the canonical projection $ M(\lambda) \rightarrow V(\lambda) $ has a nonzero kernel if and only if
$ \chi_\lambda(\det(\Sh_\nu)) = 0 $ for some Shapovalov determinant $ \det(\Sh_\nu) $ with $ \nu \in \roots^+ $. 

We claim that $ \ann_{FU_q(\mathfrak{g})}(V(\lambda)) \cap \mathbb{H} \neq 0 $ 
iff $ \chi_\lambda(\det(\P_\mu)) = 0 $ for some PRV determinant $ \det(\P_\mu) $ with $ \mu \in \weights^+ $. 
Indeed, consider the space $ H = \ann_{FU_q(\mathfrak{g})}(V(\lambda)) \cap \mathbb{H} $. Then $ H $ is invariant under the adjoint action and 
therefore $ \chi_\lambda(\P(H)) = 0 $ by Lemma \ref{annihilatorlemma}. 
Hence if $ H $ is nonzero we have $ \chi_\lambda(\det(\P_\mu)) = 0 $ for some $ \mu $. Conversely, if $ \chi_\lambda(\det(\P_\mu)) = 0 $ there exist 
scalars $ c_1, \dots, c_m $, not all zero, such that $ \sum_i c_i \chi_\lambda(\P(v_{ij})) = 0 $ for all $ j $. Let us write $ V $ for the 
linear span of the vectors $ \sum c_i v_{ij} $ for $ j = 1, \dots, m $. 
Then $ V $ is the weight zero subspace of a uniquely determined $ U_q(\mathfrak{g}) $-invariant subspace $ L \subset \mathbb{H}^\mu $. 
Applying again Lemma \ref{annihilatorlemma} we have $ L \subset H $, so $ H $ is nonzero.

According to Theorem \ref{PRVformula}, we find that $ \chi_\lambda(\det(\P_\mu)) $ is zero for some $ \mu \in \weights^+ $ iff 
\[
q_\alpha^{2(\alpha^\vee, \lambda + \rho)} - q_\alpha^{2n} = 0 
\]
for some $ \alpha \in {\bf \Delta}^+ $ and $ n \in \mathbb{N} $. By Theorem \ref{shapovalov}, this corresponds precisely to the condition for the vanishing 
of $ \chi_\lambda(\det(\Sh_\nu)) $ for some $ \nu \in \roots^+ $. This yields the claim. 
\end{proof}

\subsection{Annihilators of Verma modules} 

This section is devoted to the Verma module annihilator Theorem of Joseph and Letzter \cite{JLVermaannihilator}, see also Theorem 8.3.9 in \cite{Josephbook}. 
We shall rely crucially on results obtained in Chapter \ref{chuqg} and the considerations in previous sections.

Let $ \lambda \in \mathfrak{h}_q^* $ and consider the associated Verma module $ M(\lambda) = U_q(\mathfrak{g}) \otimes_{U_q(\mathfrak{b})} \mathbb{C}_\lambda $, 
where we recall that $ \mathbb{C}_\lambda $
\label{nom:C_lambda2}%
denotes the $ U_q(\mathfrak{b}) $-module associated with the character $ \chi_\lambda $. 
The canonical left action of elements in $ U_q(\mathfrak{g}) $ on $ M(\lambda) $ defines an algebra homomorphism 
$ \phi_\lambda: U_q(\mathfrak{g}) \rightarrow \End(M(\lambda)) $. 

The map $ \phi_\lambda $ is $ U_q(\mathfrak{g}) $-linear with respect to the adjoint action 
$$
 \ad(X)(Y) = X \rightarrow Y = X_{(1)} Y \hat{S}(X_{(2)})
$$
on $ U_q(\mathfrak{g}) $ and the action given by 
$$ 
 (X \cdot T)(m) = X_{(1)} \cdot T(\hat{S}(X_{(2)}) \cdot m) 
$$
\label{nom:End-action}%
on $ \End(M(\lambda)) $. 
In particular, $ \phi_\lambda $ 
induces an algebra homomorphism $ FU_q(\mathfrak{g}) \rightarrow F \End(M(\lambda)) $,
which we again denote by $ \phi_\lambda $.

We write $ \ann_{FU_q(\mathfrak{g})}(M(\lambda)) \subset FU_q(\mathfrak{g}) $
\label{nom:ann2}
for the annihilator of $ M(\lambda) $ viewed 
as a $ FU_q(\mathfrak{g}) $-module, that is, for the kernel of $ \phi_\lambda: FU_q(\mathfrak{g}) \rightarrow F \End(M(\lambda)) $. 

Recall that $ ZU_q(\mathfrak{g}) \subset FU_q(\mathfrak{g}) $ denotes the centre of $ U_q(\mathfrak{g}) $. 
In a similar way we write $ \ann_{ZU_q(\mathfrak{g})}(M(\lambda)) \subset ZU_q(\mathfrak{g}) $ for the annihilator of 
$ M(\lambda) $ viewed as a $ ZU_q(\mathfrak{g}) $-module. 

\begin{theorem}[Verma module Annihilator Theorem] \label{Vermaannihilator}
Let $ \lambda \in \mathfrak{h}_q^* $. Then 
$$ 
\ann_{FU_q(\mathfrak{g})}(M(\lambda)) = FU_q(\mathfrak{g}) \ann_{ZU_q(\mathfrak{g})}(M(\lambda)),
$$ 
and the linear map $ \phi_\lambda: FU_q(\mathfrak{g})/\ann_{FU_q(\mathfrak{g})}(M(\lambda)) \rightarrow F \End(M(\lambda)) $ 
is an isomorphism. 
\end{theorem} 

\begin{proof}  We first claim that
$$ 
FU_q(\mathfrak{g}) = \mathbb{H} + FU_q(\mathfrak{g}) \ann_{ZU_q(\mathfrak{g})}(M(\lambda)), 
$$ 
where $\mathbb{H}$ is the subspace of $FU_q(\mathfrak{g})$ obtained in the Separation of Variables Theorem \ref{separation}.
Indeed, according to Theorem \ref{separation}, we have $ FU_q(\mathfrak{g}) = \mathbb{H} \otimes ZU_q(\mathfrak{g}) $. Moreover, recall that 
the centre $ ZU_q(\mathfrak{g}) $ acts on $ M(\lambda) $ by the central character $ \xi_\lambda $. In particular, 
we have $ Z - \xi_\lambda(Z) 1 \in \ann_{ZU_q(\mathfrak{g})}(M(\lambda)) $ for all $ Z \in ZU_q(\mathfrak{g}) $. 
Hence we can write any element $ h \otimes Z \in \mathbb{H} \otimes ZU_q(\mathfrak{g}) $ in the form 
$$
\xi_\lambda(Z) h \otimes 1 + h \otimes (Z - \xi_\lambda(Z) 1) \in \mathbb{H} + FU_q(\mathfrak{g}) \ann_{ZU_q(\mathfrak{g})}(M(\lambda))
$$ 
as desired. 

Assume first that $ M(\lambda) $ is simple. Then Proposition \ref{annHprop} implies 
$$ 
\ann_{FU_q(\mathfrak{g})}(M(\lambda)) \cap \mathbb{H} = 0. 
$$
It follows that the restriction of the map $ FU_q(\mathfrak{g}) \rightarrow FU_q(\mathfrak{g})/\ann_{FU_q(\mathfrak{g})}(M(\lambda)) $ 
to $ \mathbb{H} $ is an injective map $ \mathbb{H} \rightarrow FU_q(\mathfrak{g})/\ann_{FU_q(\mathfrak{g})}(M(\lambda)) $.  It is also surjective since $FU_q(\mathfrak{g}) = \mathbb{H} + FU_q(\mathfrak{g}) \ann_{ZU_q(\mathfrak{g})}(M(\lambda))$ and $FU_q(\mathfrak{g}) \ann_{ZU_q(\mathfrak{g})}(M(\lambda)) \subset \ann_{FU_q(\mathfrak{g})}(M(\lambda))$.  Therefore, we have an isomorphism
\[
 FU_q(\mathfrak{g})/\ann_{FU_q(\mathfrak{g})}(M(\lambda)) \cong \mathbb{H}. 
\]

Still assuming that $ M(\lambda) $ is simple, Proposition \ref{principalseriesmultiplicity} implies 
$$
[F \End(M(\lambda)): V(\nu)] = \dim(V(\nu)_0)
$$
for any $ \nu \in \weights^+ $, where $ [N: V(\nu)] $ denotes the multiplicity of $ V(\nu) $ inside $ N $ and $ V(\nu)_0 $ is the zero weight 
space of $ V(\nu) $. From Separation of Variables in Theorem \ref{separation} we know that $ [\mathbb{H}: V(\nu)] = \dim(V(\nu)_0) $ as well. 
By comparing the multiplicities of all isotypical components it follows that the inclusion map 
$$
\mathbb{H} \cong FU_q(\mathfrak{g})/\ann_{FU_q(\mathfrak{g})}(M(\lambda)) \rightarrow F \End(M(\lambda))
$$
is an isomorphism. This finishes the proof in the case that $ M(\lambda) $ is simple. 

Now consider an arbitrary $ \lambda \in \mathfrak{h}^*_q $. Then according to Lemma \ref{Vermasocle} there exists $ \lambda' \leq \lambda $ such 
that $ M(\lambda') \subset M(\lambda) $ and $ M(\lambda') $ is simple. Due to Proposition \ref{annHprop} 
we have $ \ann_{FU_q(\mathfrak{g})}(M(\lambda)) \cap \mathbb{H} \subset \ann_{FU_q(\mathfrak{g})}(M(\lambda')) \cap \mathbb{H} = 0 $. 
Hence we obtain 
$$
FU_q(\mathfrak{g})/\ann_{FU_q(\mathfrak{g})}(M(\lambda)) \cong \mathbb{H}. 
$$
in the same way as above. Consider the commutative diagram 
$$
\xymatrix{
F\Hom(M(\lambda), M(\lambda')) \ar@{->}[r] \ar@{->}[d] & F\Hom(M(\lambda), M(\lambda)) \ar@{->}[d] \\
F\Hom(M(\lambda'), M(\lambda')) \ar@{->}[r] & F\Hom(M(\lambda'), M(\lambda))
     }
$$
induced by the inclusion $ M(\lambda') \subset M(\lambda) $. Using Proposition \ref{GKhelp} and exactness of the Hom-functor 
we see that the two vertical maps are injections and the bottom horizontal 
map is an isomorphism. 

This gives us an injection $ F \End(M(\lambda)) \cong F \End(M(\lambda')) $ and we obtain a commutative diagram 
$$
\xymatrix{
\mathbb{H} \cong FU_q(\mathfrak{g})/\ann_{FU_q(\mathfrak{g})}(M(\lambda)) \ar@{->}[r] \ar@{->}[d] & F \End(M(\lambda)) \ar@{->}[d] \\
\mathbb{H} \cong FU_q(\mathfrak{g})/\ann_{FU_q(\mathfrak{g})}(M(\lambda')) \ar@{->}[r] & F \End(M(\lambda')).
     }
$$
The bottom horizontal map is an isomorphism since $ M(\lambda') $ is simple.  The left vertical map is also an isomorphism.  The right vertical map is therefore surjective. Since we have shown it is injective, this map is in fact an isomorphism, and so is the top horizontal map.
This finishes the proof. 
\end{proof}

\newpage 

\section{Representation theory of complex semisimple quantum groups} \label{chreptheory}

In this chapter, we discuss the representation theory of complex semisimple quantum groups.  The appropriate notion of a $G_q$-representation here is that of a Harish-Chandra module for $G_q$, which means an essential $\DF(G_q)$-module with $K_q$-types of finite multiplicity, see Section \ref{sec:G_q-representations}.
In particular, the irreducible unitary representations of $G_q$ belong to this class.
	
Our main focus will be the classification of irreducible Harish-Chandra modules.
This is achieved by studying the (non-unitary) principal series representations of $G_q$ and the intertwining operators between them.
We shall begin, however, with some results on the Verma modules for the quantized universal enveloping algebra $ U_q^\mathbb{R}(\mathfrak{g}) $ associated to a complex semisimple quantum group $G_q$, see Section \ref{sec:UqRg-Verma-modules}. 

\medskip

Although some of the constructions and results presented here work more generally, we shall assume throughout that $ \mathbb{K} = \mathbb{C} $ 
and $ 1 \neq q = e^h $ is positive. We identify $\lie{h}_q^* = \lie{h}^*/i\hbar^{-1}\roots^\vee$ where $ \hbar = \tfrac{h}{2 \pi} $. Moreover we write 
$$
[z]_q = \frac{q^z - q^{-z}}{q - q^{-1}} 
$$
for the $ q $-number associated with $ z \in \mathbb{C} $, and similarly use the notation for $ q $-binomial coefficients as in Section \ref{secqnumbers}.

\subsection{Verma modules for $ U_q^\mathbb{R}(\mathfrak{g}) $} 
\label{sec:UqRg-Verma-modules}

In this section we discuss the theory of Verma modules for the quantized universal enveloping algebra 
\[
	U_q^\mathbb{R}(\mathfrak{g}) = U_q^\RR(\lie{k}) \bowtie \CF^\infty(K_q)
\]  
of the complex quantum group $ G_q $, which was introduced in Section \ref{subsecuqrg}.

\subsubsection{Characters of $U_q^\RR(\lie{b})$}
\label{sec:characters_of_UqRb}

Characters of the algebra $U_q^\RR(\lie{b}) = U_q^\RR(\lie{t})\bowtie \CF^\infty(K_q)$ are parametrized by pairs of weights $ (\mu, \lambda) $ where $ \mu, \lambda \in \mathfrak{h}^*_q $, as we now describe. 

To $ \mu \in \mathfrak{h}^*_q $ we associate the character $ \chi_\mu $ of $ U_q^\mathbb{R}(\mathfrak{t}) = U_q(\mathfrak{h}) $ given by 
$ \chi_\mu(K_\nu) = q^{(\mu, \nu)} $ for all $\nu\in\weights$, as in Subsection \ref{sec:weights}.

To $\lambda\in \mathfrak{h}^*_q$ we associate a character $K_\lambda$ of $\CF^\infty(K_q)$ as follows.
\label{nom:K_lambda}%
Briefly, if $\lambda = \sum_i a_i \alpha_i \in \lie{h}^*$ with $a_i\in\mathbb{C}$ then we may define $K_\lambda = \prod_i K_i^{a_i} \in \M(\DF(K_q))$ by functional calculus.  More explicitly, recalling that the multiplier algebra of $\DF(K_q)$ is 
\[
 \M(\DF(K_q)) = \CF^\infty(\hat{K}_q) 
  = \prod_{\gamma \in \weights^+} \End(V(\gamma)),
\]
we define $K_\lambda \in \M(\DF(K_q))$ to be the element which acts on any vector $v\in V(\gamma)$ of weight $\nu$ by 
\[
 K_\lambda \cdot v = q^{(\lambda,\nu)}v.
\]
This is compatible with the existing notation for elements $K_\lambda$ when $\lambda\in\weights$, and depends only on the class of $\lambda$ in $\lie{h}_q^* = \lie{h}^*/i\hbar^{-1}\roots^\vee$.  
The element $K_\lambda$ determines a character of $\CF^\infty(K_q)$ by
\[
  f \mapsto (f, K_\lambda).
\]
Note that we are using the reverse pairing $(f, K_\lambda) = (\hat{S}(K_\lambda),f) = (K_{-\lambda},f)$ in the definition of this character, 
since we will want to interpret $K_\lambda$ as an element of the function algebra $\CF^\infty(\hat{K}_q)$, see the discussion after Definition \ref{def:CKq},.

We remark that $K_\lambda$ belongs to the subalgebra $\M(\DF(T)) \subset \M(\DF(K_q))$ where $T$ is the maximal torus subgroup.  It is in fact the pullback of a character $K_\lambda$ of $\CF^\infty(T)$ by the quotient map $\CF^\infty(K_q) \to \CF^\infty(T)$.  Moreover, every character of $ \CF^\infty(K_q) $ is of this form, see \cite{Yakimovspectra}.  In this way, the parameter space $\lie{h}^*_q$ is identified with the Chevalley complexification $T_\CC$ of the torus $T$.

Combining the characters $ \chi_\mu $ and $ K_\lambda $ for any $ \mu, \lambda \in \mathfrak{h}^*_q $ we obtain 
a character $ \chi_{\mu, \lambda} $
\nomenclature{$\chi_{\mu,\lambda}$}{character of $U_q^\mathbb{R}(\lie{b})$ associated to $(\mu,\lambda)\in\weights\times\lie{h}_q^*$}%
of $ U_q^\mathbb{R}(\mathfrak{b}) = U_q^\mathbb{R}(\mathfrak{t}) \bowtie \CF^\infty(K_q) $ by setting 
$$ 
\chi_{\mu, \lambda}(X \bowtie f) = \chi_\mu(X) (f, K_\lambda).
$$  
Notice that $ \chi_{\mu, \lambda} $ is typically not $ * $-preserving for the standard $ * $-structure on $ U_q^\mathbb{R}(\mathfrak{b}) $.

Let us conclude with some remarks about the restriction of these characters to the subalgebra $ U_q^\mathbb{R}(\lie{t}) \bowtie \CF^\infty(PK_q)$, which corresponds to the Borel subgroup of the ``connected component'' $G_q^0$ of $G_q$ that was discussed in Subsection \ref{sec:connected_component}.  One sees that the character
$\chi_{\mu,\lambda}$ restricts to the trivial character on $ U_q^\mathbb{R}(\lie{t}) \bowtie \CF^\infty(PK_q)$ if and only if $(\mu,\lambda) = (0,i\hbar^{-1}\gamma)$ for some $\gamma\in\weights^\vee/\roots^\vee$.  It follows that the characters of $ U_q^\mathbb{R}(\lie{t}) \bowtie \CF^\infty(PK_q)$ are indexed by pairs $(\mu,\lambda) \in \lie{h}_q^* \times \lie{h}^*/i\hbar^{-1}\weights^\vee$.

\subsubsection{Verma modules for $U_q^\RR(\lie{g})$}
\label{sec:UqRg-Verma}

The Verma module associated to $ \chi_{\mu, \lambda} $ is defined by 
$$
M(\mu, \lambda) = U_q^\mathbb{R}(\mathfrak{g}) \otimes_{U_q^\mathbb{R}(\mathfrak{b})} \mathbb{C}_{\mu, \lambda}, 
$$
\nomenclature{$M(\mu,\lambda)$}{Verma module for $U_q^\mathbb{R}(\lie{g})$}%
where $ \mathbb{C}_{\mu, \lambda} = \mathbb{C} $
\nomenclature[o$C_{\mu,\lambda}$]{$\CH_{\mu,\lambda}$}{one-dimensional $U_q^\mathbb{R}(\mathfrak{b})$-module}
is the representation space corresponding to the character $ \chi_{\mu, \lambda} $ defined above.  Following previous notation, we will denote by $v_{\mu,\lambda}$ the cyclic vector $1\otimes1 \in M(\mu, \lambda) $.
\nomenclature{$v_{\mu,\lambda}$}{cyclic vector of the $U_q^\mathbb{R}(\lie{b})$-Verma module $M(\mu,\lambda)$}%

Recall from Lemma \ref{uqrgdiagonalembedding} 
that we have an embedding $ U_q^\mathbb{R}(\mathfrak{g}) \subset U_q(\mathfrak{g}) \otimes U_q(\mathfrak{g}) $ of algebras given explicitly by 
$$
\iota(X \bowtie f) = \hat{\Delta}(X) (l^-(f_{(1)}) \otimes l^+(f_{(2)})). 
$$
In this way we may view modules over $ U_q(\mathfrak{g}) \otimes U_q(\mathfrak{g}) $ as modules over $ U_q^\mathbb{R}(\mathfrak{g}) $. 
However, it will be convenient in later sections to use a different embedding of $U_q^\mathbb{R}(\lie{g})$ in $U_q(\lie{g})\otimes U_q(\lie{g})$ as follows.
 
\begin{definition}
	\label{lem:iota_prime}
 We define the embedding $\iota':U_q^\mathbb{R}(\lie{g}) \to U_q(\lie{g})\otimes U_q(\lie{g})$ by
 \[
   \iota'(X\bowtie f) = \hat{\Delta}(X) (l^+(f_{(1)}) \otimes l^-(f_{(2)}))
 \]
 \nomenclature{$\iota'$}{alternative diagonal embedding $\iota': U_q^\mathbb{R}(\lie{g}) \to U_q(\lie{g})\otimes U_q(\lie{g})$}%
 for $X\bowtie f \in U_q^\mathbb{R}(\lie{g})$.
\end{definition} 

Note that $\iota$ and $\iota'$ differ only by the order of the $l$-functionals in the formula.
Again, $\iota'$ is an algebra homomorphism.  

We observe that for any $X\bowtie f \in U_q^\mathbb{R}(\lie{g})$ we have
\begin{align*}
\iota'(X\bowtie f) 
&= \sigma((\hat{S}\otimes\hat{S})(\iota((1\bowtie S^{-1}(f))(\hat{S}^{-1}(X)\bowtie 1)))),
\end{align*}
using Lemma \ref{lfunctionalproperties}.  Thus we can also write
\[
\iota' = \sigma\circ(\hat{S}\otimes\hat{S})\circ \iota \circ \hat{S}_{G_q}^{-1},
\]
where $\sigma$ denotes the flip map.

In the following proposition, $N(\lambda)$ denotes the universal lowest weight module of $U_q(\lie{g})$ with lowest weight $\lambda\in\lie{h}_q^*$, namely
\[
 N(\lambda) = U_q(\lie{g}) \otimes_{U_q(\lie{b}_-)} \mathbb{C}_\lambda
\]
\nomenclature{$N(\lambda)$}{universal lowest weight module of $U_q(\lie{g})$}%
where $\mathbb{C}_\lambda$ denotes the one-dimensional $U_q(\lie{b}_-)$-module upon which $U_q(\lie{n}_-)$ acts trivially and $U_q(\lie{h})$ acts by the character $\chi_\lambda$.  We denote the cyclic vector $1\otimes 1$ by $v^\lambda$.
\nomenclature{$v^\lambda$}{lowest weight vector of the universal lowest weight module $N(\lambda)$}%

\begin{prop}
 \label{prop:Verma_module_iso}
 Let $ (\mu, \lambda) \in \mathfrak{h}^*_q \times \mathfrak{h}^*_q $ and let $l,r\in\mathfrak{h}^*_q$ be such that
\[
  \mu = l-r, \qquad  \lambda= -l-r.
\]
\label{nom:lr}%
We have an isomorphism 
 $$ 
  M(\mu, \lambda) \cong M(l) \otimes N(-r)
 $$ 
 of $ U_q^\mathbb{R}(\mathfrak{g}) $-modules which sends the cyclic vector $v_{\mu,\lambda}$ to $v_l \otimes v^{-r}$.  Here the action of $U_q^\mathbb{R}(\lie{g})$ on $M(l) \otimes N(-r)$ is that induced by the embedding $\iota'$ above. 
\end{prop}

\begin{proof}  
For all $\nu\in\weights$, the action of $K_\nu\in U_q^\mathbb{R}(\mathfrak{t}) $ on $v_l \otimes v^{-r}$ 
is given by 
\begin{align*}
K_\nu \cdot (v_{l} \otimes v^{-r}) &= K_{\nu} \cdot v_{l} \otimes K_{\nu} \cdot v^{-r} \\
&= q^{(l,\nu)} v_{l} \otimes q^{(-r,\nu)} v^{-r} 
= \chi_\mu(K_\nu) v_{l} \otimes v^{-r}.
\end{align*}

We make a similar calculation for the action of a matrix coefficient $ u^\eta_{ij} = \bra e^i_\eta| \bullet| e^\eta_j \ket $
 where $ \eta \in \weights^+ $ and the 
 $e^\eta_j$ form a basis of weight vectors for $ V(\eta) $, with dual basis vectors $ e^i_\eta \in V(\eta)^* $. 
Writing $ \epsilon_j $ for the weight of  $e^\eta_j$ we obtain 
\begin{align*}
u^\eta_{ij} \cdot (v_{l} \otimes v^{-r}) 
&= \sum_{k}  l^+(u^\eta_{ik}) \cdot v_{l} \otimes l^-(u^\eta_{kj}) \cdot v^{-r} \\ 
&= \delta_{ij} l^+(u^\eta_{ii}) \cdot v_{l} \otimes l^-(u^\eta_{ii}) \cdot v^{-r} \\ 
&= \delta_{ij} K_{\epsilon_i} \cdot v_{l} \otimes K_{-\epsilon_i} \cdot v^{-r} \\ 
&= \delta_{ij} q^{(l,\epsilon_i)} v_{l} \otimes q^{(r,\epsilon_i)} v^{-r} \\
&= \delta_{ij} q^{(-\lambda, \epsilon_i)} v_{l} \otimes v^{-r} \\
&= (u^\eta_{ij},K_{\lambda})\, v_{l} \otimes v^{-r}. 
\end{align*}
Here, in the second step, we use that $ l^\pm(u^\eta_{ij}) $ is contained in $ U_q(\mathfrak{b}_\pm)_{\epsilon_j - \epsilon_i} $, and that the vectors $ v_{-l}$ and $v_{r} $ are annihilated by $ U_q(\mathfrak{n}_+)$ and $U_q(\mathfrak{n}_-)$, respectively.  Note also that when $\epsilon_i=\epsilon_j$ we have $ l^\pm(u^\nu_{ij}) = \delta_{ij} K_{\pm \epsilon_i} $. 

By the definition of $ M(\mu, \lambda) $, we thus obtain a $ U_q^\mathbb{R}(\mathfrak{g}) $-linear map 
$ \gamma: M(\mu, \lambda) \rightarrow M(l) \otimes N(-r) $ such that $ \gamma(v_{\mu, \lambda}) = v_l\otimes v^{-r}$.

For surjectivity, we will prove inductively that the subspaces $M(l)_{l-\nu} \otimes N(-r)$ are contained in $\im(\gamma)$ for each $\nu\in\roots^+$.  For the case $\nu=0$, note that the action of $E_i \in U_q^\mathbb{R}(\lie{k})$ on $v_l \otimes v^{-r}$ is given by
\[
 E_i\cdot(v_{l}\otimes v^{-r}) =   v_{l} \otimes E_i \cdot v^{-r}.
\]
It follows that $\mathbb{C}v_l \otimes N(-r) \subset \im(\gamma)$.  Now fix $\nu\in\roots^+$ and suppose that $M(l)_{l-\nu'} \otimes N(-r) \subset \im(\gamma)$ for all $\nu' < \nu$.  
Note that $M(l)_{l-\nu} \otimes N(-r)$ is spanned by elements of the form $F_i\cdot m \otimes n$ with $m\in M(l)_{l-\nu+\alpha_i}$ and $n \in N(-r)$.  Using the action of $F_i$ on $n\otimes m$ we can write 
\[
  (F_i\cdot m) \otimes n =  F_i\cdot(m\otimes n)  - K_i^{-1}\cdot m \otimes F_i\cdot n   \in \im(\gamma),
\]
and surjectivity follows.

For injectivity, we note to begin with that $ U_q^\mathbb{R}(\mathfrak{g}) $ is a free right $ U_q^\mathbb{R}(\mathfrak{b}) $-module
generated by $ U_q(\mathfrak{n}_-) \otimes U_q(\mathfrak{n}_+) \subset U_q(\mathfrak{g}) = U_q^\mathbb{R}(\mathfrak{k}) \subset U_q^\mathbb{R}(\mathfrak{g}) $.  That is, the map
\[
  U_q(\mathfrak{n}_-) \otimes U_q(\mathfrak{n}_+) \to M(\mu,\lambda); \qquad Y\otimes X \mapsto YX\cdot v_{\mu,\lambda}
\]
is an isomorphism.  Consider an element 
\begin{equation*}
 \label{eq:verma_idenitification_injectivity}
  \sum_j Y_j \otimes X_j \in U_q(\mathfrak{n}_-) \otimes U_q(\mathfrak{n}_+)
\end{equation*}
where $ Y_j \in U_q(\mathfrak{n}_-)_{-\nu'_{j}}$, $X_j \in U_q(\mathfrak{n}_+)_{\nu''_{j}} $ for some weights $\nu'_j, \nu''_j \in \roots^+$.
We may arrange the sum so that $\nu''_0$ is maximal among those $\nu''_j$ appearing nontrivially, and $\nu'_0$ is maximal among the weights $\nu'_j$ such that $\nu''_j = \nu''_0$ and $(\nu'_j,\nu''_j)$ appears nontrivially.
By considering the actions of $E_i$ and $F_i$, one sees that 
\[
 \gamma ( \textstyle \sum_j Y_{\nu'_{j}} X_{\nu''_{j}} \cdot v_{\mu,\lambda})
   = \sum_j Y_{\nu'_{j}} X_{\nu''_{j}} \cdot (v_l \otimes v^{-r})
\]
contains a nonzero term in $ M(l)_{l-\nu'_0} \otimes N(-r)_{-r+\nu''_0} $.  
This completes the proof.
\end{proof}  

In order to convert the lowest weight module $N(-r)$ into a highest weight module, we may twist the action of $U_q(\lie{g})$ by the algebra automorphism 
\[
 \theta = \tau\hat{S} = \hat{S}^{-1}\tau,
\]
\label{nom:theta2}%
where $\tau$ is the algebra anti-automorphism defined in Lemma \ref{deftau}.
We record the formulas
\[
 \theta(E_i) = -F_i, \qquad \theta(F_i) = -E_i, \qquad \theta(K_\nu) = K_{-\nu},
\]
from which it follows that $N(-r)$ is isomorphic to $M(r)^\theta$, where the latter denotes the Verma module $M(r)$ endowed with the action obtained by composition with $\theta$.

Let us further note that for any $U_q(\lie{g})$-module $M$, we have $M^{\hat{S}^2} \cong M$, where $M^{\hat{S}^2}$ denotes the same space with action twisted by the automorphism $\hat{S}^2$.  This is because the two modules are intertwined by the action of $K_{2\rho}$, see Lemma \ref{lem:S2}.  Therefore, if we put $\theta' = \hat{S}\tau = \hat{S}^2\theta$ then we also have $N(-r) \cong M(r)^{\theta'}$.  
\nomenclature{$\theta'$}{algebra automorphism, coalgebra anti-automorphism $\theta'=\hat{S}\tau$ of $U_q(\lie{g})$}%

We immediately obtain the following.

\begin{cor}
 \label{cor:Verma_module_iso}
 Let $ (\mu, \lambda) \in \mathfrak{h}^*_q \times \mathfrak{h}^*_q $ and let $l,r\in\mathfrak{h}^*_q$ be such that
\[
  \mu = l-r, \qquad  \lambda= -l-r.
\]
 We have an isomorphism of $U_q^\mathbb{R}(\lie{g})$-modules 
 $$ 
  M(\mu, \lambda) \cong M(l) \otimes M(r)
 $$ 
 which sends the cyclic vector $v_{\mu,\lambda}$ to $v_l \otimes v_{r}$, and where the action of $X \bowtie f \in U_q^\mathbb{R}(\lie{g})$ on $m\otimes n \in M(l) \otimes M(r)$ is given by
 \begin{align*}
  (X\bowtie f)\cdot (m\otimes n) &= (\id\otimes\theta)(\iota'(X\bowtie f))\cdot (m\otimes n).
 \end{align*}
 The same is true if $\theta$ is replaced by $\theta'$.
\end{cor}

We point out that for a given character $\chi_{\mu,\lambda}$ of $U_q^\RR(\mathfrak{b})$ with $\mu,\lambda \in \mathfrak{h}^*_q$, there are $2^N$ choices for the associated parameters $(l,r) \in \mathfrak{h}^*_q \times \mathfrak{h}^*_q$, where $N$ is the rank of $G$. Specifically, 
if we lift $\mu, \lambda$ to $\mathfrak{h}^*$ then we obtain the solutions for $(l,r)$ given by
\[
 \textstyle
 l= -\half (\lambda-\mu) +\half i \hbar^{-1} \alpha^\vee, 
  \qquad
  r = -\half (\lambda+\mu) + \half i \hbar^{-1} \alpha^\vee,
\]
\label{nom:extended_diagonal_action}%
for $\alpha^\vee \in \roots^\vee / 2 \roots^\vee$.
In other words, if we let the extended Weyl group from Definition \ref{def:extended_Weyl_group} 
\[
 \hat{W} = \mathbf{Y}_q \rtimes W = (\textstyle\half i \hbar^{-1} \roots^\vee / i \hbar^{-1} \roots^\vee) \rtimes W 
\]
act diagonally upon the parameter space $\lie{h}_q^* \times \lie{h}_q^*$ by $\hat{w}(l,r) = (\hat{w}l,\hat{w}r)$ then the action of the translation subgroup $\mathbf{Y}_q  = \half i \hbar^{-1} \roots^\vee / i \hbar^{-1} \roots^\vee$ does not change the isomorphism class of the associated $U_q^\mathbb{R}(\lie{g})$-modules $M(l)\otimes M(r)$.

The following criterion allows us to characterize irreducibility of the Verma modules $ M(\mu, \lambda) $ introduced above.   We recall that in our conventions we use $\mathbb{N} = \{1,2,\ldots\}$.

\begin{theorem} \label{Vermadoubleirreducible}
The $U_q^\mathbb{R}(\lie{g})$-module $ M(\mu, \lambda) $ is irreducible if and only if $ (\mu, \lambda) \in \mathfrak{h}^*_q \times \mathfrak{h}^*_q $ satisfies 
$$ 
q_\alpha^{(-\lambda+2\rho, \alpha^\vee)} \neq q_\alpha^{(\pm \mu, \alpha^\vee) + 2m}
$$ 
for all $ \alpha \in \bf{\Delta}^+ $ and all $ m \in \mathbb{N} $. 
\end{theorem}

\begin{proof}

As in Proposition \ref{cor:Verma_module_iso} we choose $l,r \in \mathfrak{h}^*_q$ such that $\mu=l-r$, $\lambda=-l-r$.  
Then $ M(l) \otimes M(r) $ is irreducible as 
a $ U_q(\mathfrak{g}) \otimes U_q(\mathfrak{g}) $-module if and only if $ M(l) $ and $ M(r) $ are both irreducible as $ U_q(\mathfrak{g}) $-modules.
According to Theorem \ref{Vermacomplexirreducible} and Definition \ref{defdominantantidominant}, this is the case iff
\[
   q_\alpha^{2(l+\rho,\alpha^\vee)} \notin  q_\alpha^{2\NN} \quad \text{and}\quad
   q_\alpha^{2(r+\rho,\alpha^\vee)} \notin  q_\alpha^{2\NN}
\]
for all $\alpha\in\bf{\Delta}^+$.  Using $2l=-\lambda+\mu$ and $2r = -\lambda-\mu$, this is equivalent to the condition
\[
 q_\alpha^{(-\lambda+2\rho,\alpha^\vee)} \neq q_\alpha^{(\pm\mu,\alpha^\vee) + 2m}
\]
for all $\alpha\in\mathbf{\Delta}^+$, $m\in\NN$.

It remains to check that $ M(l) \otimes M(r) $ is irreducible as a 
a $ U_q(\mathfrak{g}) \otimes U_q(\mathfrak{g}) $-module if and only if it is irreducible as a $ U_q^\mathbb{R}(\mathfrak{g}) $-module with the action induced by $(\theta\otimes\id)\iota'$,  
compare the argument in the proof of Theorem 3.4 in \cite{Hodgesdouble}.  
For the nontrivial implication assume that $ M(l) \otimes M(r) $ is irreducible as a $ U_q(\mathfrak{g}) \otimes U_q(\mathfrak{g}) $-module. 
According to Lemma \ref{uqrgdiagonalembedding}, the image of $\iota'$ is $(1\otimes\hat{S}(FU_q(\lie{g}))) \hat{\Delta}(U_q(\lie{g}))$.
Thus, Theorem \ref{adjointduality} shows that the image of $(\id\otimes\theta)\iota'$ contains 
$ K_{\nu'} \otimes K_{2\nu''-\nu'}$ for all $ \nu' \in \weights $ and $ \nu'' \in \weights^+ $. It follows that any $ U_q^\mathbb{R}(\mathfrak{g}) $-submodule 
$ V \subset M(l) \otimes M(r) $ is the sum of its $ U_q(\mathfrak{g}) \otimes U_q(\mathfrak{g}) $-weight spaces.

Let $v = \sum_j x_j \otimes y_j \in V$ be a nonzero vector of weight  $(\epsilon_l,\epsilon_r)$ with $\epsilon=\epsilon_l+\epsilon_r$ maximal in the submodule $V$.  We may assume that the $y_j$ are linearly independent  in $M({r})_{\epsilon_{r}}$.  Since
\[
  E_i\cdot(x_j\otimes y_j) = E_i \cdot x_j \otimes K_{i}^{-1}\cdot y_j - x_j \otimes F_i\cdot y_j,
\]
the maximality of $\epsilon$ forces $E_i\cdot x_j=0$ for all $i$, and so each $x_j$ is primitive.  
Since $ M(l) $ is irreducible we conclude that each $ x_j $ is a scalar multiple of $ v_{l} $, and thus $v=v_l \otimes y$ for some $y$. Likewise, using
\[
  F_i \cdot(v_l \otimes y) = F_i\cdot v_l \otimes y - K_i^{-1}\cdot v_l \otimes E_i\cdot y,
\]
the maximality of $\epsilon$ shows that $y$ is primitive, and hence is a nonzero multiple of $v_{r}$.
By Corollary \ref{cor:Verma_module_iso}, $v_l \otimes v_r$ is a cyclic vector for the $U_q^\RR(\mathfrak{g})$-action, so $V=M(l)\otimes M(r)$.  This yields the claim.
\end{proof}

\subsection{Representations of $G_q$} 
\label{sec:G_q-representations}

\subsubsection{Harish-Chandra modules}

In this subsection we define the main notion of $G_q$-representation which we will be studying in the remainder of the chapter.  

We recall that the convolution algebra of $G_q$ is $\DF(G_q) = \DF(K_q) \bowtie \CF^\infty(K_q)$.  Recall also that a $\DF(G_q)$-module $V$ is called essential if the multiplication map $\DF(G_q) \otimes_{\DF(G_q)} V \to V$ is an isomorphism.  
Using the central idempotents in $\DF(K_q)$ one sees that a $\DF(G_q)$-module is essential if and only if it is an essential $\DF(K_q)$-module under restriction.

Recall from Section \ref{subsecuqrg} that the quantized universal enveloping algebra $U_q^\RR(\lie{g}) = U_q^\RR(\lie{k}) \bowtie \CF^\infty(K_q)$ of $G_q$ sits inside the multiplier algebra of $\DF(G_q)$. 
Let us say that a $ U_q^\mathbb{R}(\mathfrak{g}) $-module $ V $ is integrable if it is an integrable module for the action 
of $ U_q^\mathbb{R}(\mathfrak{k}) \subset U_q^\mathbb{R}(\mathfrak{g}) $. The following result is then essentially immediate from the definitions. 

\begin{lemma}
	\label{lem:DG-UqRg-modules}
	There is a canonical isomorphism between the category of essential $ \DF(G_q) $-modules and 
	the category of integrable $ U_q^\mathbb{R}(\mathfrak{g}) $-modules. 
\end{lemma}

\begin{proof}  Every essential $ \DF(G_q) $-module becomes an integrable $ U_q^\mathbb{R}(\mathfrak{g}) $-module via the 
	inclusion $ U_q^\mathbb{R}(\mathfrak{g}) \subset \M(\DF(G_q)) $. 
	
	Conversely, if $ V $ is an integrable $ U_q^\mathbb{R}(\mathfrak{g}) $-module then the action 
	of $ U_q^\mathbb{R}(\mathfrak{k}) \subset U_q^\mathbb{R}(\mathfrak{g}) $ corresponds uniquely to an essential $ \DF(K_q) $-module structure on $ V $, 
	and the latter combines with the action of $ \CF^\infty(K_q) \subset U_q^\mathbb{R}(\mathfrak{g}) $ 
	to turn $ V $ into an essential $ \DF(G_q) $-module. 
	
	These procedures are inverse to each other and compatible with morphisms. \end{proof}

A third structure which is equivalent to an essential $\DF(G_q)$-module is that of a $\CF^\infty(K_q)$-Yetter Drinfeld module.  We recall the general definition.

\begin{definition} \label{defyd}
	Let $ H $ be a Hopf algebra. A Yetter-Drinfeld module over $ H $ is a left $ H $-module $ V $ which is at the same 
	time a left $ H $-comodule such that 
	$$
	(f \cdot v)_{(-1)} \otimes (f \cdot v)_{(0)} = f_{(1)} v_{(-1)} S(f_{(3)}) \otimes f_{(2)} \cdot v_{(0)}
	$$
	for all $ v \in V $ and $ f \in H $. 
	Here we write $ \gamma(v) = v_{(-1)} \otimes v_{(0)} $ for the left coaction of $ H $ on $ V $. 
\end{definition}

If $V$ is a Yetter-Drinfeld module over $\CF^\infty(K_q)$ then we can convert the left $\CF^\infty(K_q)$-coaction into an essential
left action of $\DF(K_q)$ by the formula
\[
 x\cdot v =  (\hat{S}(x),v_{(-1)}) v_{(0)}, \qquad  x\in\DF(K_q).
\]
\label{nom:YD-Kq-action}%
A standard calculation shows that the resulting $\DF(K_q)$-action and the given $\CF^\infty(K_q)$-action satisfy the commutation relations in the Drinfeld double $\DF(G_q)$.
It follows that the category of $\CF^\infty(K_q)$-Yetter-Drinfeld modules is isomorphic to the category of essential $\DF(G_q)$-modules.

Since any essential $\DF(G_q)$-module is an essential $\DF(K_q)$-module $V$ by restriction, it decomposes as a direct sum 
$$
V = \bigoplus_{\gamma \in \weights^+} V^\gamma, 
$$
where $ V^\gamma $ is a direct sum of copies of the simple module with highest weight $ \gamma $. The subspace $ V^\gamma \subset V $ 
is called the isotypical component corresponding to $ \gamma $. 
\nomenclature[s$V^\gamma$]{$V^\gamma$}{isotypical component of a $\DF(G_q)$- or $\DF(K_q)$-module $V$ associated with highest weight $\gamma\in\weights^+$}%

\begin{definition} 
 \label{defadmissible}
 An essential $\DF(K_q)$-module $V $ is called admissible if each isotypical component $ V $ is finite dimensional.  An essential $\DF(G_q)$-module which is admissible as a $\DF(K_q)$-module will be called a Harish-Chandra module over $G_q$.
\end{definition}

We remark that the term ``Harish-Chandra module'' has several different meanings in the classical literature, in particular with or without the imposition of admissibility.  The imposition of admissibility in this definition avoids many technical annoyances, and as in the classical case, all irreducible unitary $G_q$-representations are admissible, see the following section.

\subsubsection{Unitary $G_q$-representations}
\label{sec:unitary-Gq-reps}

By definition, a unitary representation of $G_q$ on a Hilbert space $\H$ is a non-degenerate $ * $-homomorphism $ C^*_{\mx}(G_q) \rightarrow \LH(\H)$.

In particular any unitary representation $ \H $ of $G_q$ becomes a unitary representation of $ K_q $ by restriction. 
A vector $ \xi \in \H $ is called $ K_q $-finite if it is contained in a finite dimensional $K_q$-subrepresentation. 
The Harish-Chandra module $ \HC(\H) $ associated with $ \H $ is the space of all $ K_q $-finite vectors in $ \H $. 
\nomenclature[o$HC(H)$]{$\HC(\H)$}{Harish-Chandra module associated with a unitary $G_q$ representation $\H$}%
Explicitly, this is given by 
$$ 
\HC(\H) = \DF(K_q) \cdot \H = \DF(G_q) \cdot \H \subset \H. 
$$ 
From this description it is clear that $ \HC(\H) $ is dense in $ \H $, and that $\HC(\H)$ naturally becomes an essential module over $ \DF(G_q) $.

The goal of this subsection is to show that the $\DF(G_q)$-module $\HC(\H) $ associated with an irreducible unitary representation $ \H $ of $ G_q $ is admissible, so that it is a Harish-Chandra module in the sense of Definition \ref{defadmissible}. The argument is the same as in the classical case, based on results of Godement.

Firstly, if $ A $ is any algebra then the $ n $-commutator of $ n $ elements $ a_1, \dots, a_n \in A $ is defined by 
$$
[a_1, \dots, a_n] = \sum_{\sigma \in S_n} \sign(\sigma) a_{\sigma(1)} \cdots a_{\sigma(n)}. 
$$
Note that for $ n = 2 $ this reduces to the usual commutator $ [a_1, a_2] = a_1 a_2 - a_2 a_1 $. 
Let us say that $ A $ is $ n $-abelian if all $ n $-commutators $ [a_1, \dots, a_n] $ for $ a_1, \dots, a_n \in A $ vanish. 
Clearly, subalgebras and quotients of $ n $-abelian algebras are again $ n $-abelian. 

By basic linear algebra, every finite dimensional algebra $ A $ is $ n $-abelian for any $ n > \dim(A) $. 
In particular, for every $ k \in \mathbb{N} $ there exists a smallest number $ r(k) \in \mathbb{N} $ such that $ M_k(\mathbb{C}) $ is $ r(k) $-abelian. 

\begin{lemma} 
	With the notation as above, we have $ r(k + 1) > r(k) $ for all $ k \in \mathbb{N} $. 
\end{lemma} 
\begin{proof}  Put $ n = r(k) - 1 $.  Then there exists $ x_1, \dots, x_n \in M_k(\mathbb{C}) $ 
	such that $ X = [x_1, \dots, x_n] \neq 0 $. In particular, if we write $ X = (X_{ij}) $ then there exist indices $ 1 \leq i,j \leq k $ such that $ X_{ij} \neq 0 $. 
	Set 
	$$
	y_l = 
	\begin{pmatrix}
	x_l & 0 \\
	0 & 0
	\end{pmatrix}
	\in M_{k + 1}(\mathbb{C})
	\qquad \text{for } l=1,\ldots n,
	$$
	and $ y_{n + 1} = e_{j, k + 1} $ the standard matrix unit with $1$ in the $(j,k+1)$-position.
	Then $ y_{n + 1} y_l = 0 $ for all $ 1 \leq l \leq n $ and thus 
	$$
	[y_1, \dots, y_{n + 1}] = [y_1, \dots, y_n] y_{n + 1} = [y_1, \dots, y_n] e_{j, k + 1} \neq 0.  
	$$
	We conclude $ r(k + 1) \geq r(k) + 1 $ as claimed. \end{proof}  

We are now ready to prove admissibility of irreducible unitary representations of complex quantum groups, compare \cite{Aranospherical}. 

\begin{theorem} \label{unitaryadmissibility}
	Let $ \H $ be an irreducible unitary representation of $ G_q $. Then the associated $\DF(G_q)$-module $ \HC(\H) $ is admissible, and so a Harish-Chandra module. More precisely, for any $ \mu \in \weights^+ $ the multiplicity of the $ K_q $-type $ \HC(\H)^\mu $ in $ \HC(\H) $ 
	is at most $ \dim(V(\mu)) $. 
\end{theorem} 

\begin{proof}  According to Proposition \ref{finiteseparating}, the embedding 
	$ \DF(G_q) \rightarrow \M(\DF(K_q) \otimes \DF(K_q)) $ restricts to an embedding 
	$$ 
	p_\mu \DF(G_q) p_\mu \rightarrow \prod_{\eta, \nu \in \weights^+} \End(p_\mu \cdot (V(\eta) \otimes V(\nu))). 
	$$
	Let $ v_\eta \in V(\eta) $ and $ v^{\nu} \in V(\nu) $ be the highest weight vector and lowest weight vector, respectively. 
	Then $ v_\eta \otimes v^{\nu} $ is a cyclic vector for $ V(\eta) \otimes V(\nu) $ as a $U_q(\mathfrak{g})$-module, and therefore the 
	linear map $ \Hom_{U_q(\mathfrak{g})}(V(\eta) \otimes V(\nu), V(\mu)) \rightarrow V(\mu) $ sending $ f $ to $ f(v_\eta \otimes v^{\nu}) $ is injective. 
	We conclude that the dimension of $ p_\mu \cdot (V(\eta) \otimes V(\nu)) $ is at most $ d = \dim(V(\mu))^2 $. 
	Therefore the algebra $ A = p_\mu \DF(G_q) p_\mu $ is $ r(d) $-abelian. 
	
	Now let $ \pi: C^*_\mx(G_q) \rightarrow \LH(\H) $ be an irreducible unitary representation. Then the von Neumann algebra  
	$ \pi(C^*_\mx(G_q))'' $ equals $ \LH(\H) $ by irreducibility, and hence $ \pi(p_\mu) \pi(C^*_\mx(G_q))'' \pi(p_\mu) = \LH(\pi(p_\mu) \H) $. 
	If we set $ A = p_\mu \DF(G_q) p_\mu $ 
	this means that the strong closure $ \pi(A)'' $ of $ \pi(A) \subset \LH(p_\mu \H) $ is 
	equal to $ \LH(\pi(p_\mu) \H) $. Since $ \pi(A) $ is $ r(d) $-abelian, the same holds for its strong closure $ \LH(\pi(p_\mu) \H) $. 
	We conclude that $ \pi(p_\mu) \in \LH(\H) $ is a finite rank projection of rank at most $ d $. Hence the multiplicity of $ V(\mu) $ 
	in $ \H $ is at most $ \dim(V(\mu)) $ as claimed. \end{proof}  

Let $ \pi $ be an irreducible unitary representation of $ G_q $ on the Hilbert space $ \H $. 
According to Theorem \ref{unitaryadmissibility}, the image of any element of $ \DF(K_q) \subset \DF(G_q) $ is a finite-rank operator on $ \H $. 
This implies that $ \DF(G_q) $ acts by finite-rank operators. 
In particular, the associated representation $ \pi: C^*_\mx(G_q) \rightarrow \LH(\H) $ takes values in the algebra of compact operators on $ \H $. 

Let us say that a locally compact quantum group $ G $ is type $ I $ if the full group $ C^* $-algebra $ C^*_{\mx}(G) $ is 
a $ C^* $-algebra of type $ I $. As an immediate consequence of the above observations we obtain the following result.  

\begin{cor} 
	Complex semisimple quantum groups are type $ I $. 
\end{cor} 

We remark that the corresponding result for classical semisimple groups is due to Harish-Chandra, see \cite{HCrepssliegroupbanach}.

\subsection{Action of $ U_q^\mathbb{R}(\mathfrak{g}) $ on $ K_q $-types}
\label{sec:action_on_K_q_types}

In this section we collect some facts on the structure of integrable $ U_q^\mathbb{R}(\mathfrak{g}) $-modules, following Section 9.1 
in \cite{Dixmierenveloping}. 

Let $ V $ be an essential $ \DF(G_q) $-module. We may also view $ V $ as an integrable $ U_q^\mathbb{R}(\mathfrak{g}) $-module, see Lemma \ref{lem:DG-UqRg-modules}. 
By definition, the module $ V $ decomposes as a direct sum 
$$
V = \bigoplus_{\gamma \in \weights^+} V^\gamma 
$$
of its $ K_q $-isotypical components.

For $ \sigma \in \weights^+ $ we denote 
\[
 I^\sigma = \ker(\pi_\sigma) = \ann(V(\sigma)) \subset U_q^\mathbb{R}(\mathfrak{k}) ,
\] 
\nomenclature{$I^\sigma$}{annihilator ideal of the finite dimensional $U_q^\mathbb{R}(\lie{k})$-module $V(\sigma)$}%
where $ \pi_\sigma: U_q^\mathbb{R}(\mathfrak{k}) \rightarrow \End(V(\sigma)) $ 
is the homomorphism corresponding to the action on $ V(\sigma) $. By construction we have $ U_q^\mathbb{R}(\mathfrak{k})/I^\sigma \cong \End(V(\sigma)) $. 
Notice that 
$$ 
V^\tau = \{v \in V \mid X \cdot v = 0 \text{ for all } X \in I^\tau \}. 
$$

For $ \sigma, \tau \in \weights^+ $ we define 
$$
U_q^\mathbb{R}(\mathfrak{g})^{\tau, \sigma} 
= \{X \in U_q^\mathbb{R}(\mathfrak{g}) \mid I^\tau X \subset U_q^\mathbb{R}(\mathfrak{g}) I^{\sigma} \}.
$$
\nomenclature{$U_q^\mathbb{R}(\mathfrak{g})^{\tau, \sigma} $}{$= \{X \in U_q^\mathbb{R}(\mathfrak{g}) \mid I^\tau X \subset U_q^\mathbb{R}(\mathfrak{g}) I^{\sigma} \}$}%
By construction, we have $ U_q^\mathbb{R}(\mathfrak{g})^{\tau, \sigma} \cdot V^\sigma \subset V^\tau $ for any $ U_q^\mathbb{R}(\mathfrak{g}) $-module $ V $.

\begin{lemma} \label{Ug-K-types}
Let $\sigma\in\weights^+$. With the notation as above, the following holds. 
\begin{bnum} 
\item[a)] $ U_q^\mathbb{R}(\mathfrak{g}) I^\sigma \subset U_q^\mathbb{R}(\mathfrak{g})^{\tau,\sigma} $ for all $\tau\in\weights^+$.
\item[b)] $ U_q^\mathbb{R}(\mathfrak{g})/U_q^\mathbb{R}(\mathfrak{g}) I^\sigma $ is an integrable
 $ \DF(G_q) $-module with respect to left multiplication. 
Its decomposition into $ K_q $-types is 
\begin{equation*} 
U_q^\mathbb{R}(\mathfrak{g})/U_q^\mathbb{R}(\mathfrak{g}) I^\sigma 
= \bigoplus_{\tau \in \weights^+} U_q^\mathbb{R}(\mathfrak{g})^{\tau,\sigma}/U_q^\mathbb{R}(\mathfrak{g}) I^\sigma.
\end{equation*}
\item[c)] For all $\tau\in\weights^+$, the space
$ U_q^\mathbb{R}(\mathfrak{g})^{\tau,\sigma}$ consists precisely of all elements of $ U_q^\mathbb{R}(\mathfrak{g}) $ 
that map $ V^\sigma $ into $ V^\tau $ for any $ U_q^\mathbb{R}(\mathfrak{g}) $-module $ V $. 
\end{bnum}
\end{lemma} 

\begin{proof}  $ a) $ is obvious. 

$ b) $
Using $U_q^\mathbb{R}(\mathfrak{k})/ I^\sigma \cong \End(V(\sigma))$, one sees that the 
generator $1 \in U_q^\mathbb{R}(\mathfrak{g}) / U_q^\mathbb{R}(\mathfrak{g})I^\sigma$ is $K_q$-finite. 
It follows that $U_q^\mathbb{R}(\mathfrak{g}) / U_q^\mathbb{R}(\mathfrak{g})I^\sigma$ is integrable.
The $ \tau $-isotypical subspace is the subspace annihilated by $ I^\tau $, that is, 
$$
(U_q^\mathbb{R}(\mathfrak{g})/U_q^\mathbb{R}(\mathfrak{g}) I^\sigma)^\tau 
= \{X + U_q^\mathbb{R}(\mathfrak{g}) I^\sigma \mid I^\tau X \subset U_q^\mathbb{R}(\mathfrak{g}) I^\sigma \} 
= U_q^\mathbb{R}(\mathfrak{g})^{\tau,\sigma}/U_q^\mathbb{R}(\mathfrak{g}) I^\sigma.
$$
This yields the claim. 

$ c) $ follows from the decomposition into $K_q$-types in part $ b) $.
\end{proof}

We will mostly be interested in $ U_q^\mathbb{R}(\mathfrak{g})^{\sigma, \sigma} $ for any given $ \sigma \in \weights^+ $. It is easy to check 
that this is an algebra. It contains $ U_q^\mathbb{R}(\mathfrak{g}) I^\sigma $ as an obvious two-sided ideal, 
and $ U_q^\mathbb{R}(\mathfrak{g}) I^\sigma $ acts as zero on every $ V^\sigma $. We will typically factor this ideal out. 
Accordingly, for any $ U_q^\mathbb{R}(\mathfrak{g}) $-module $ V $ we get a map 
$ U_q^\mathbb{R}(\mathfrak{g})^{\sigma, \sigma}/U_q^\mathbb{R}(\mathfrak{g})I^\sigma \rightarrow \End(V^\sigma) $.

\begin{prop} \label{subquotienthelp}
If $ V $ is a simple integrable $ U_q^\mathbb{R}(\mathfrak{g})$-module then $ V^\sigma $ is a simple $ U_q^\mathbb{R}(\mathfrak{g})^{\sigma, \sigma} $-module 
or zero. 
\end{prop}

\begin{proof}  Let $ v \in V^\sigma $ be nonzero. By simplicity of $ V $ we have $ U_q^\mathbb{R}(\mathfrak{g}) \cdot v = V $. Moreover, 
using the above observations and part $ b) $ of Lemma \ref{Ug-K-types} we get 
$ U_q^\mathbb{R}(\mathfrak{g}) \cdot v = \sum_{\tau \in \weights^+} U_q^\mathbb{R}(\mathfrak{g})^{\tau,\sigma} \cdot v $. 
Since $ U_q^\mathbb{R}(\mathfrak{g})^{\tau,\sigma} \cdot v \subset V^\tau $, we 
have $ V^\sigma = (U_q^\mathbb{R}(\mathfrak{g}) \cdot v)^\sigma = U_q^\mathbb{R}(\mathfrak{g})^{\sigma,\sigma} \cdot v $. \end{proof}  

The converse is not true, that is, simplicity of $ V^\sigma $ as a $ U_q^\mathbb{R}(\mathfrak{g})^{\sigma, \sigma} $-module  
does not imply that $ V $ is simple. Nonetheless, we will show the existence of a unique simple $ U_q^\mathbb{R}(\mathfrak{g}) $-module associated 
to each simple $ U_q^\mathbb{R}(\mathfrak{g})^{\sigma, \sigma}/ U_q^\mathbb{R}(\mathfrak{g}) I^\sigma $-module. 

We start with the following construction.
Let $ V $ be a Harish-Chandra module and $ \sigma \in \weights^+ $ such that $ V^\sigma \neq 0 $. Moreover assume that 
$ W \subset V^\sigma $ is a $ U_q^\mathbb{R}(\mathfrak{g})^{\sigma, \sigma} $-submodule. 
Then we define 
$$
W^\mathrm{min} = U_q^\mathbb{R}(\mathfrak{g}) \cdot W, \qquad 
W^\mathrm{max} = \{v \in V \mid U_q^\mathbb{R}(\mathfrak{g}) \cdot v \cap V^\sigma \subseteq W \}. 
$$
The following Proposition shows that $ W^\mathrm{min} $ and $ W^\mathrm{max} $ are minimal and maximal, respectively, among the  
$ U_q^\mathbb{R}(\mathfrak{g}) $-submodules $ V' \subset V $ satisfying $(V')^\sigma = W $. 

\begin{prop} \label{Wminmax}
We have $ (W^\mathrm{min})^\sigma = (W^\mathrm{max})^\sigma = W $, and for any $ U_q^\mathbb{R}(\mathfrak{g}) $-submodule $ V' $ of $ V $ such 
that $ (V')^\sigma = W $, we have $ W^\mathrm{min} \subset V' \subset W^\mathrm{max} $.
\end{prop}

\begin{proof}  
From Lemma \ref{Ug-K-types} $ c) $ we have that $(U_q^\mathbb{R}(\lie{g})\cdot W) \cap V^\sigma = U_q^\mathbb{R}(\lie{g})^{\sigma,\sigma}\cdot W = W$, and it follows that $(W^\mathrm{min})^\sigma = (W^\mathrm{max})^\sigma = W$.  If $V'$ is as given, then $ W^\mathrm{min} = U_q^\mathrm{R}(\lie{g}) \cdot W \subset U_q^\mathrm{R}(\lie{g}) \cdot V' \subset W^\mathrm{max}$, as claimed.
\end{proof}  

\begin{prop} \label{prop:ideal_correspondence}
For $ \sigma \in \weights^+ $ define 
\begin{align*}
\mathcal{M}_\sigma &= \{\text{maximal left ideals of $U_q^\mathbb{R}(\mathfrak{g})$ containing $U_q^\mathbb{R}(\mathfrak{g})I^\sigma$} \}, \\
\mathcal{L}_\sigma &= \{\text{maximal left ideals of $U_q^\mathbb{R}(\mathfrak{g})^{\sigma,\sigma}$ containing $U_q^\mathbb{R}(\mathfrak{g})I^\sigma$} \}.
\end{align*}
Then there is a bijective correspondence $ \phi: \mathcal{M}_\sigma \rightarrow \mathcal{L}_\sigma $ given by
$$
\phi(M) = M \cap U_q^\mathbb{R}(\mathfrak{g})^{\sigma,\sigma}. 
$$
The inverse of $ \phi $ is given by 
$$
\phi^{-1}(L) = \{x \in U_q^\mathbb{R}(\mathfrak{g}) \mid U_q^\mathbb{R}(\mathfrak{g})x \cap U_q^\mathbb{R}(\mathfrak{g})^{\sigma, \sigma} \subset L \}.
$$
\end{prop}

\begin{proof}  
Notice that we have a natural bijection between 
left ideals of $ U_q^\mathbb{R}(\mathfrak{g})$ containing $U_q^\mathbb{R}(\mathfrak{g})I^\sigma $ and 
$ U_q^\mathbb{R}(\mathfrak{g})$-submodules of $ U_q^\mathbb{R}(\mathfrak{g})/U_q^\mathbb{R}(\mathfrak{g})I^\sigma $. 
If $ M $ is a maximal left ideal of $ U_q^\mathbb{R}(\mathfrak{g}) $ then $ 1 \notin M $, and 
therefore $ M \cap U_q^\mathbb{R}(\mathfrak{g})^{\sigma, \sigma} \not = U_q^\mathbb{R}(\mathfrak{g})^{\sigma, \sigma} $.

Similarly, we have a natural bijection between left ideals of $ U_q^\mathbb{R}(\mathfrak{g})^{\sigma, \sigma} $ 
containing $ U_q^\mathbb{R}(\mathfrak{g})I^\sigma $ and $ U_q^\mathbb{R}(\mathfrak{g})^{\sigma, \sigma} $-submodules 
of $ (U_q^\mathbb{R}(\mathfrak{g})/U_q^\mathbb{R}(\mathfrak{g})I^\sigma)^\sigma $. 

Now it suffices to apply Proposition \ref{Wminmax} to $ V = U_q^\mathbb{R}(\mathfrak{g})/U_q^\mathbb{R}(\mathfrak{g})I^\sigma $. 
\end{proof}

\begin{prop}
\label{prop:simple_modules_vs_simple_K-types}
Let $ \sigma \in \weights^+ $.  Consider
\begin{align*}
  \S_\sigma & = \{  \text{simple integrable $ U_q^\mathbb{R}(\mathfrak{g})$-modules $ V $ with $ V^\sigma \neq 0 $, modulo isomorphism} \}, \\
  \T_\sigma &= \{ \text{simple $ U_q^\mathbb{R}(\mathfrak{g})^{\sigma,\sigma}/U_q^\mathbb{R}(\mathfrak{g}) I^\sigma $-modules, modulo isomorphism} \}.
\end{align*} 
There is a natural bijective correspondence $\psi:\S_\sigma \to \T_\sigma$ given by
\[
  \psi(V) = V^\sigma.
\]
\end{prop}

\begin{proof}  
The map $\psi$ is well-defined by  Proposition \ref{subquotienthelp}.

For injectivity, suppose that $V,V'$ are simple integrable $U_q^\mathbb{R}(\lie{g})$-modules with nontrivial $\sigma$-isotypical component, and 
that $ f: V^\sigma \rightarrow (V')^\sigma $ is an isomorphism of $ U_q^\mathbb{R}(\mathfrak{g})^{\sigma, \sigma}/U_q^\mathbb{R}(\mathfrak{g}) I^\sigma $-modules.
Let $ v \in V^\sigma $ be nonzero and set
\begin{itemize}
\item $ L = $ annihilator of $ v $ in $ U_q^\mathbb{R}(\mathfrak{g})^{\sigma,\sigma} $,
\item $ M = $ annihilator of $ v $ in $ U_q^\mathbb{R}(\mathfrak{g}) $,
\item $ L' = $ annihilator of $ f(v) $ in $ U_q^\mathbb{R}(\mathfrak{g})^{\sigma,\sigma} $,
\item $ M'= $ annihilator of $ f(v) $ in $ U_q^\mathbb{R}(\mathfrak{g}) $.
\end{itemize}
These are all maximal left ideals of $U_q^\mathbb{R}(\lie{g})$ or $U_q^\mathbb{R}(\lie{g})^{\sigma,\sigma}$ respectively, containing $U_q^\mathbb{R}(\lie{g})I^\sigma$. 
But 
$$ 
M \cap U_q^\mathbb{R}(\mathfrak{g})^{\sigma, \sigma} = L = L' = M' \cap U_q^\mathbb{R}(\mathfrak{g})^{\sigma, \sigma}, 
$$
so by the correspondence of Proposition \ref{prop:ideal_correspondence} we deduce $ M = M' $. 

For surjectivity, let $ W $ be a simple $ U_q^\mathbb{R}(\mathfrak{g})^{\sigma,\sigma} $-module whose annihilator contains $ U_q^\mathbb{R}(\mathfrak{g})I^\sigma $. 
Moreover let $ w \in W $ be nonzero. Put
\begin{itemize}
\item $ L = $ annihilator of $ w $ in $ U_q^\mathbb{R}(\mathfrak{g})^{\sigma,\sigma} $,
\item $ M = \phi^{-1}(L) $, the maximal left ideal of $ U_q^\mathbb{R}(\mathfrak{g}) $ such that $ M \cap U_q^\mathbb{R}(\mathfrak{g})^{\sigma,\sigma} = L $, 
which is obtained from Proposition \ref{prop:ideal_correspondence},
\item $ V = U_q^\mathbb{R}(\mathfrak{g})/M $.
\end{itemize}
Then $ V $ is simple and by Lemma \ref{Ug-K-types} $ b) $, 
\begin{align*}
V^\sigma &= (U_q^\mathbb{R}(\mathfrak{g})/U_q^\mathbb{R}(\mathfrak{g})I^\sigma)^\sigma/(M/U_q^\mathbb{R}(\mathfrak{g})I^\sigma)^\sigma \\ 
&= U_q^\mathbb{R}(\mathfrak{g})^{\sigma, \sigma}/(M \cap U_q^\mathbb{R}(\mathfrak{g})^{\sigma,\sigma}) \\ 
&= U_q^\mathbb{R}(\mathfrak{g})^{\sigma,\sigma}/L = W.
\end{align*}
This finishes the proof. 
\end{proof}

\subsection{Principal series representations} 

In this section we define principal series representations for the complex quantum group $ G_q $. As in the classical case, the principal series is key to the analysis of the representation theory 
of complex quantum groups. 

\subsubsection{The definition of principal series representations}
\label{sec:principal_series_def} 

Recall from Subsection \ref{sec:characters_of_UqRb} that associated to any $\mu\in\lie{h}_q^*$ we have a character of $U_q^\mathbb{R}(\lie{t}) = U_q(\lie{h})$ given by $\chi_\mu(K_\nu) = q^{(\mu,\nu)}$.  In the case where $\mu\in\weights$ is an integral weight, this character can be obtained as the pairing of $U_q^\mathbb{R}(\lie{t})$ with an appropriately chosen matrix coefficient in $\CF^\infty(K_q)$.
In other words, the character $\chi_\mu$ can be viewed as an element of the function algebra $\CF^\infty(T)$, compare Section \ref{sec:parabolic_subgroups}.  We will denote this function by $e^\mu \in \CF^\infty(T)$,
\label{nom:e_mu}%
so that
\[
(K_\nu, e^\mu) = \chi_\mu(K_\nu) = q^{(\mu,\nu)},
\]
for all $\nu\in\weights$.

In analogy with the classical case, there is 
an ``associated line bundle'' $ \E_\mu $ over the quantum flag variety $ \mathcal{X}_q = K_q/T $. 
This bundle is defined via its space of 
sections
\begin{equation*} \label{eq:line_bundle}
\Gamma(\E_\mu) = \{\xi \in \CF^\infty(K_q) \mid (\id \otimes \pi_T) \Delta(\xi) =  \xi \otimes e^\mu \}, 
\end{equation*}
\nomenclature{$\Gamma(\E_\mu)$}{space of sections of the line bundle $\E_\mu$ over the quantum flag variety}%
where $ \pi_T: \CF^\infty(K_q) \rightarrow \CF^\infty(T) $ denotes the canonical projection map. 
\nomenclature{$\pi_T$}{projection map $\pi_T: \CF^\infty(K_q) \rightarrow \CF^\infty(T) $}%
Note that $ \Gamma(\E_\mu) \subset \CF^\infty(K_q) $ is equal to the subspace of weight $ \mu $ with respect to the 
left $ U_q^\mathbb{R}(\mathfrak{k}) $-action 
$$
X \hit \xi =\xi_{(1)} (X, \xi_{(2)}).  
$$

The space $ \Gamma(\E_\mu) $ is a left $\CF^\infty(K_q)$-comodule with coaction given by $\Delta$.
For $ \lambda \in \mathfrak{h}_q^* $ we define the twisted left adjoint representation of $ \CF^\infty(K_q) $ on 
$ \Gamma(\E_\mu) $ by 
\begin{equation*} 
f \cdot \xi =  f_{(1)} \, \xi \, S(f_{(3)})  \; (K_{2 \rho + \lambda}, f_{(2)}) .
\label{nom:twisted_adjoint}
\end{equation*}
This combination of action and coaction make $\Gamma(\E_\mu)$ into a Yetter-Drinfeld module over $\CF^\infty(K_q)$, see Definition \ref{defyd}, since we have
\begin{align*}
 \Delta(f\cdot \xi)
  &=  f_{(1)} \xi_{(1)} S(f_{(5)}) \otimes f_{(2)} \xi_{(2)} S(f_{(4)}) \; (K_{2\rho+\lambda} , f_{(3)}) \\
  &=  f_{(1)} \xi_{(1)} S(f_{(3)}) \otimes f_{(2)} \cdot \xi_{(2)}
\end{align*}
for all $f\in\CF^\infty(K_q)$, $\xi\in\Gamma(\E_{\mu})$.

We point out that the factor of $2\rho$ in the above formula is included merely to force a shift in the parameter $\lambda$.  This shift is chosen such that for purely imaginary $\lambda$ the resulting Yetter-Drinfeld module is always unitary, see Section \ref{sec:unitary} below.

\begin{definition} \label{defprincipalseries}
We write $\Gamma(\E_{\mu, \lambda})$ 
\nomenclature{$\Gamma(\E_{\mu,\lambda})$}{principal series representation of $\DF(G_q)$ or $U_q^\mathbb{R}(\lie{g})$}%
for the space $\Gamma(\E_\mu)$ equipped with the action and coaction of $\CF^\infty(K_q)$ as above, and call it
the \emph{principal series Yetter-Drinfeld module}, or principal series representation, with parameter $ (\mu, \lambda) \in \weights \times \mathfrak{h}^*_q $.  
\end{definition}

Equivalently, $\Gamma(\E_{\mu, \lambda}) $ can be viewed as an essential module over $ \DF(G_q) $  
by transforming the coaction into the associated left $ \DF(K_q) $-action
$$
x \cdot \xi = (\hat{S}(x), \xi_{(1)}) \xi_{(2)} 
$$
for $x \in \DF(K_q) $ and $ \xi \in \Gamma(\E_\mu) $.  Likewise, replacing $x\in\DF(K_q)$ by $X\in U_q^\mathbb{R}(\lie{k})$ in this formula we obtain the realization of $\Gamma(\E_{\mu, \lambda})$ as
an integrable $U_q^\mathbb{R}(\mathfrak{g})$-module, see Lemma \ref{lem:DG-UqRg-modules}.

  We record the following formula for the multiplicities of $K_q$-types in $\Gamma(\E_{\mu,\lambda})$.

\begin{lemma}
\label{lem:principal_series_multiplicities}
 The multiplicities of the $K_q$-types in $\Gamma(\E_{\mu,\lambda})$ are given by
 \[
  [\Gamma(\E_{\mu,\lambda}) : V(\nu)] = \dim((V(\nu)^*)_\mu) = \dim (V(\nu)_{-\mu}),
 \]
 for all $\nu\in\weights^+$.  
\end{lemma}

\begin{proof}
Note that we have
\[
\CF^\infty(K_q) \cong \bigoplus_{\nu\in\weights^+} (\End V(\nu))^* \cong \bigoplus_{\nu\in\weights^+} V(\nu)^* \otimes V(\nu)
\cong \bigoplus_{\nu\in\weights^+} V(-w_0\nu) \otimes V(-w_0\nu)^*,
\] 
where the left and right regular $U_q^\mathbb{R}(\lie{k})$-actions correspond to the actions on the left and right tensor factors, respectively.  Here, the dual modules are equipped with the contragredient action. 
Therefore, we have an isomorphism of $U_q^\mathbb{R}(\lie{k})$-modules
\[
 \Gamma(\E_{\mu,\lambda}) 
 \cong  \bigoplus_{\nu\in\weights^+} V(\nu) \otimes (V(\nu)^*)_\mu.
\]
This proves the first equality.  The second equality follows from $\dim((V(\nu)^*)_\mu) = \dim((V(\nu)_{-\mu})^*) = \dim (V(\nu)_{-\mu})$.
\end{proof}

Thus, in our conventions, the principal series representation  $\Gamma(\E_{\mu,\lambda})$ has a nontrivial $K_q$-isotypical component of highest weight $\nu$ if and only if $\nu \in (-\mu)^+ + \roots^+$, where $(-\mu)^+$ denotes the unique dominant integral weight which is conjugate to $-\mu$ under the Weyl group action.  
In particular, $\Gamma(\E_{\mu,\lambda})$ has a minimal $K_q$-type of highest weight $(-\mu)^+$ and this minimal $K_q$-type occurs with multiplicity one.

\subsubsection{Compact versus noncompact pictures} 

Definition \ref{defprincipalseries} is referred to as the ``compact picture'' of the principal series representation. 
One may simply accept this definition without motivation and check that it satisfies the Yetter-Drinfeld condition in 
Definition \ref{defyd}. But for more insight, we note that 
principal series representations also admit an interpretation as representations of $ G_q $ induced from characters of the parabolic quantum 
subgroup $ B_q = T \bowtie \hat{K}_q $.

Let $ \mathbb{C}_{\mu, \lambda} $ denote the one dimensional representation of $ B_q $ with character $ \chi_{\mu, \lambda} $ as in the 
definition of the Verma module $ M(\mu, \lambda) $ in Section \ref{sec:characters_of_UqRb}. By definition, the (algebraic) unitarily induced representation of $ G_q $ is
$$
\ind_{B_q}^{G_q}(\mathbb{C}_{\mu,\lambda}) = \{\xi \in \CF^\infty(G_q) \mid
(\id \otimes \pi_{B_q}) \Delta_{G_q}(\xi) = \xi \otimes (e^\mu \otimes K_{2 \rho + \lambda}) \}.
$$
\nomenclature[o$ind$]{$\ind_{B_q}^{G_q}(\mathbb{C}_{\mu,\lambda})$}{principal series representation of $G_q$ as an induced representation}%
Here $ K_{2 \rho + \lambda} $ is viewed as multiplier of $ \DF(K_q) $ inside $ \CF^\infty(G_q) $. Again, the shift by $ 2 \rho $ is included to 
ensure a suitable parametrization with inner products later on. Observe that the coproduct of $ \CF^\infty_c(G_q) $ induces a left $ \DF(G_q) $-module 
structure on $ \ind_{B_q}^{G_q}(\mathbb{C}_{\mu,\lambda}) $ according to the formula 
$$
x \cdot \xi = (\hat{S}_{G_q}(x), \xi_{(1)}) \xi_{(2)} 
$$
for $ x \in \DF(G_q) $ and $ \xi \in \ind_{B_q}^{G_q}(\mathbb{C}_{\mu,\lambda}) $.

On the other hand, recall that
$$ 
\Gamma(\E_{\mu,\lambda}) = \{\xi \in \CF^\infty(K_q) \mid (\id \otimes \pi_T) \Delta(\xi) = \xi \otimes e^\mu \}
$$ 
admits a $\DF(G_q)$-action corresponding to the Yetter-Drinfeld structure from Definition \ref{defprincipalseries}.

\begin{lemma} 
	\label{lem:ext}
Let $ (\mu, \lambda) \in \weights \times \mathfrak{h}^*_q $. The linear 
maps 
\[
 \ext: \Gamma(\E_{\mu, \lambda}) \rightarrow \ind_{B_q}^{G_q}(\mathbb{C}_{\mu, \lambda}) ; \qquad
     \ext(\xi)  = \xi \otimes K_{2\rho+\lambda}
\]
and
\[
 \res: \ind_{B_q}^{G_q}(\mathbb{C}_{\mu, \lambda}) \rightarrow \Gamma(\E_{\mu, \lambda}) ; \qquad
   \res(\sigma) = (\id \otimes \hat{\epsilon})(\sigma) 
\]
are well-defined and are inverse isomorphisms of $\DF(K_q)$-modules.
Moreover, the Yetter-Drinfeld action of $ \CF^\infty(K_q) $ on $ \Gamma(\E_{\mu, \lambda}) $ corresponds to the natural action 
of $ \CF^\infty(K_q) \subset \M(\DF(G_q)) $ on $ \ind_{B_q}^{G_q}(\mathbb{C}_{\mu, \lambda}) $. 
\end{lemma} 

\begin{proof}  Let us check that $ \ext(\xi) $ for $ \xi \in \Gamma(\E_\mu) $ satisfies the correct invariance properties. 
Recall that we can write the multiplicative unitary $ W $ in the form 
$$
W = \sum_{\nu, i,j} u^\nu_{ij} \otimes \omega^\nu_{ij}
$$
where $(u^\nu_{ij})$ is a basis of matrix units for $\CF^\infty(K_q)$ and and $(\omega^\nu_{ij})$ is its dual basis.
Using the definition of $\Delta_{G_q}$ from Subsection \ref{sec:CQG_definition}, we compute 
\begin{align*}
(\id \otimes \pi_{B_q})\Delta_{G_q}(\ext(\xi)) 
&= \sum_{\nu, \eta, i,j,r,s} \xi_{(1)} \otimes \omega^\nu_{ij} K_{2 \rho + \lambda} \omega^\eta_{rs} \otimes 
\pi_{B_q}(u^\nu_{ij} \xi_{(2)} S(u^\eta_{rs}) \otimes K_{2 \rho + \lambda}) \\
&= \sum_{\nu, i,j,s} \xi \otimes \omega^\nu_{ij} K_{2 \rho + \lambda} \omega^\nu_{js} 
\otimes \pi_T(u^\nu_{ij}) \,e^\mu\,  \pi_T( S(u^\nu_{js})) \otimes K_{2 \rho + \lambda} \\ 
&= \sum_{\nu, j} \xi \otimes \omega^\nu_{jj} K_{2 \rho + \lambda} \omega^\nu_{jj} 
\otimes e^\mu  \pi_T(u^\nu_{jj}) S(\pi_T(u^\nu_{jj})) \otimes K_{2 \rho + \lambda} \\ 
&= \sum_{\nu, j} 
\xi \otimes \omega^\nu_{jj} K_{2 \rho + \lambda} \omega^\nu_{jj} \otimes e^\mu \otimes K_{2 \rho + \lambda} \\ 
&= \xi \otimes K_{2 \rho + \lambda} \otimes e^\mu \otimes K_{2 \rho + \lambda} 
= \ext(\xi) \otimes e^\mu \otimes K_{2 \rho + \lambda} 
\end{align*}
in $ \M(\CF^\infty_c(G_q) \otimes \CF^\infty_c(B_q)) $. Hence $ \ext(\xi) $ satisfies the invariance condition 
in the definition of $ \ind_{B_q}^{G_q}(\mathbb{C}_{\mu, \lambda}) $. 

Similarly, for $ \sigma \in \ind_{B_q}^{G_q}(\mathbb{C}_{\mu, \lambda}) $ the element $ \res(\sigma) = (\id \otimes \hat{\epsilon})(\sigma) $ satisfies 
\begin{align*}
(\id \otimes \pi_{T})\Delta(\res(\sigma)) &= 
(\id \otimes \hat{\epsilon} \otimes \pi_T \otimes \hat{\epsilon})\Delta_{G_q}(\sigma) \\ 
&= (\id \otimes \hat{\epsilon} \otimes \id \otimes \hat{\epsilon})(\id \otimes \pi_{B_q}) \Delta_{G_q}(\sigma) \\ 
&= (\id \otimes \hat{\epsilon})(\sigma) \otimes e^\mu \hat{\epsilon}(K_{2 \rho + \lambda}) \\ 
&= \res(\sigma) \otimes e^\mu 
\end{align*}
inside $ \CF^\infty(K_q) \otimes \CF^\infty(T) $. 

It is clear that $\res\circ\ext$ is the identity on $\Gamma(\E_{\mu,\lambda})$.  For the reverse composition, we begin by observing that
\[
 (\id_{\CF^\infty(K_q)} \otimes \hat{\epsilon}\otimes \epsilon \otimes \id_{\DF(K_q)})\Delta_{G_q} = \id_{\CF_c^\infty(G_q)}.
\]
Then for $\sigma \in  \ind_{B_q}^{G_q}(\mathbb{C}_{\mu, \lambda}) $ we obtain
\begin{align*}
 \sigma &=  (\id \otimes \hat{\epsilon}\otimes \epsilon \otimes \id)\Delta_{G_q} (\sigma) \\
        &=  (\id \otimes \hat{\epsilon}\otimes \epsilon \otimes \id) (\id \otimes \pi_{B_q}) \Delta_{G_q} (\sigma) \\
        &= (\id \otimes \hat{\epsilon}\otimes \epsilon \otimes \id) (\sigma \otimes (e^\mu\otimes K_{2\rho+\lambda})) \\
        &= (\id\otimes \hat{\epsilon})(\sigma)\otimes K_{2\rho+\lambda}\\
        &= \ext\circ\res(\sigma).
\end{align*}
Therefore $\ext$ and $\res$ are isomorphisms.

Finally, we consider the actions.
The action of $x = t\bowtie a \in \DF(G_q)$ on $\sigma \in \ind_{B_q}^{G_q}(\mathbb{C}_{\mu, \lambda})$ is given by
\begin{align*}
 (t\bowtie a) \cdot \sigma &= (\hat{S}_{G_q}(t\bowtie a), \sigma_{(1)}) \sigma_{(2)} \\
 &= (t\bowtie a, S_{G_q}^{-1}(\sigma_{(1)})) \sigma_{(2)} \\
 &=  (t\bowtie a, (S^{-1} \otimes \hat{S}^{-1})(W \sigma_{(1)} W^{-1})) \sigma_{(2)} \\
   &= (\hat{S}(t) \bowtie S(a) , W \sigma_{(1)} W^{-1})  \sigma_{(2)},
\end{align*}
according to the formula before Proposition \ref{doubleunimodular}. 
To transfer this to the compact picture, we consider $\xi \in \Gamma(\E_{\mu,\lambda})$ and calculate
\begin{align*}
\res &  ((t \bowtie a) \cdot \ext(\xi )) \\
& 
 = \sum_{\nu, \eta, i,j,r,s} (\hat{S}(t) \bowtie S(a), W(\xi_{(1)} \otimes \omega^\nu_{ij} K_{2 \rho + \lambda} \omega^\eta_{rs})W^{-1}) u^\nu_{ij} \xi_{(2)} S(u^\eta_{rs})  \\
&= \sum_{\substack{\nu, \eta, \alpha, \beta, \\
i,j,r,s,p,q,m,n}} 
(\hat{S}(t) \bowtie S(a), u^\alpha_{pq} \xi_{(1)} u^\beta_{mn} \otimes \omega^\alpha_{pq} 
\omega^\nu_{ij} K_{2 \rho + \lambda} \omega^\eta_{rs} \hat{S}^{-1}(\omega^\beta_{mn})) u^\nu_{ij} \xi_{(2)} S(u^\eta_{rs}) \\
&= \sum_{\nu, \eta, i,j,r,s} (\hat{S}(t_{(2)}), \xi_{(1)})
(\hat{S}(t_{(1)}) \omega^\nu_{ij} K_{2 \rho + \lambda} \omega^\eta_{rs} \hat{S}^{-1}(\hat{S}(t_{(3)})), a) u^\nu_{ij} \xi_{(2)} S(u^\eta_{rs}) \\
&= (\hat{S}(t_{(1)}), a_{(1)}) (t_{(3)}, a_{(5)}) a_{(2)} (t_{(2)} \cdot \xi) S(a_{(4)}) (K_{2 \rho + \lambda}, a_{(3)}) \\
&= (\hat{S}(t_{(1)}), a_{(1)}) (t_{(3)}, a_{(3)}) a_{(2)} \cdot (t_{(2)} \cdot \xi)\\
& = t \cdot a \cdot \xi,  
\end{align*}
where the final lines use the actions of  $ \DF(K_q) $ and $ \CF^\infty(K_q) $ on $\Gamma(\E_{\mu,\lambda})$ given by
$$
t \cdot \xi = (\hat{S}(t), \xi_{(1)}) \xi_{(2)} 
$$
for $ t \in \DF(K_q) $ and 
$$
a \cdot \xi = a_{(1)} \xi S(a_{(3)}) (K_{2 \rho + \lambda}, a_{(2)}), 
$$
for $ a \in \CF^\infty(K_q) $, respectively. \end{proof}  

In our conventions, the principal series module $ \Gamma(\E_{0,-2\rho}) $ corresponds to the 
representation of $ G_q $ induced from the trivial representation of $ B_q $ when disregarding the $ \rho $-shift. Geometrically, this is 
the algebra of functions on the flag variety $ G_q/B_q $ equipped 
with the regular representation. In particular, this algebra contains the constant 
function $ 1 = 1 \otimes 1 \in \CF^\infty(G_q) $, which is invariant under the $ \DF(G_q) $-action, 
\begin{align*}
(x \bowtie f) \cdot 1 &= (\hat{S}_{G_q}(x \bowtie f), 1) 1 = \hat{\epsilon}_{G_q}(x \bowtie f) 1.
\end{align*}
That is, $ \Gamma(\E_{0,-2\rho}) $ contains the one-dimensional trivial representation of $ G_q $ as a submodule.

\subsubsection{The action of the centre on the principal series}
\label{sec:action_of_centre}

Recall from Subsection \ref{sec:connected_component} that the quantum group $G_q$ admits a finite group of one-dimensional unitary representations indexed by the centre $Z \cong \weights^\vee/\roots^\vee$ of $G$.  Explicitly, for each $\gamma \in \weights^\vee/\roots^\vee$ we have the unitary character $\hat{\epsilon}_\gamma : \DF(G_q) \to \mathbb{C}$ defined by
  \[
   \hat{\epsilon}_\gamma(x \bowtie f) = \hat\epsilon(x) (K_{i\hbar^{-1}\gamma}, f),
  \]
see Definition \ref{def:Gq-characters}. 
These characters induce an action of $Z$ on the class of Harish-Chandra modules for $G_q$ whereby $\gamma\in Z$ sends a Harish-Chandra module $H$ to the module $H\otimes\mathbb{C}_{i\hbar^{-1}\gamma}$, with $\mathbb{C}_{i\hbar^{-1}\gamma}$ being the one-dimensional module with action $\hat{\epsilon}_\gamma$.  On the principal series, this action is given as follows.

\begin{lemma}
	\label{lem:Z-action_on_principal_series}
 For any $(\mu,\lambda) \in \weights \times \lie{h}_q^*$ and $\gamma\in\weights^\vee/\roots^\vee$, we have $\Gamma(\E_{\mu,\lambda}) \otimes \mathbb{C}_{i\hbar^{-1}\gamma} = \Gamma(\E_{\mu,\lambda+i\hbar^{-1}\gamma})$.
\end{lemma}

\begin{proof}
 Using the compact picture, $\Gamma(\E_{\mu,\lambda}) \otimes \mathbb{C}_{i\hbar^{-1}\gamma}$ and $\Gamma(\E_{\mu,\lambda+i\hbar^{-1}\gamma})$ are naturally identified as $\DF(K_q)$-modules.  It remains to compute the Yetter-Drinfeld action of $f\in\CF^\infty(K_q)$ on $\xi\otimes1 \in \Gamma(\E_{\mu,\lambda}) \otimes \mathbb{C}_{i\hbar^{-1}\gamma}$.  Using the fact that $K_{i\hbar^{-1}\gamma}\hit f = f \hitby K_{i\hbar^{-1}\gamma}$ we obtain
 \begin{align*}
  f\cdot(\xi \otimes 1)
   &= f_{(1)}\xi S(f_{(3)}) \, (K_{2\rho+\lambda},f_{(2)}) \otimes (K_{i\hbar^{-1}\gamma}, f_{(4)}) \\
   &= f_{(1)}\xi S(f_{(4)}) \, (K_{2\rho+\lambda},f_{(2)})\,(K_{i\hbar^{-1}\gamma}, f_{(3)}) \otimes 1\\
   &= f_{(1)}\xi S(f_{(3)}) \, (K_{2\rho+\lambda+i\hbar^{-1}\gamma}, f_{(2)}) \otimes 1.
 \end{align*}
This completes the proof.
\end{proof}

	This $Z$-action is responsible for the appearance of certain representations of $G_q$ which do not have classical analogues. 
	For instance, in the case of $SL_q(2,\mathbb{C})$, Podle\'s and Woronowicz observed the appearance of a nontrivial one-dimensional representation and of a pair of complementary series \cite{PWlorentz}, as compared to the classical group $SL(2,\mathbb{C})$, which has only the trivial one-dimensional representation and only one complementary series.  This is related to the action of $Z \cong \mathbb{Z}_2$.

In general, the various representations in a single $Z$-orbit all become isomorphic upon restriction to the ``connected component'' $G_q^0$ described in Subsection \ref{sec:connected_component}, since the characters $\hat{\epsilon}_\gamma$ are all trivial on $\DF(G_q^0)$.  As a consequence, the natural parameter space for the principal series representations of $G_q^0$ is $\weights \times \lie{h}^*/i\hbar^{-1}\weights^\vee$, compare the remarks at the end of Subsection \ref{sec:characters_of_UqRb}.

\subsubsection{Duality for principal series modules}
\label{sec:principal_series_duality}

One can introduce a general notion of the dual of a Harish-Chandra module.

	\begin{definition}
		\label{def:dual_HC-module}
		Let $V$ be a Harish-Chandra module over $G_q$, see Definition \ref{defadmissible}.  We equip the dual space $V^*$ with the $\DF(G_q)$-action given by
		\[
		(x\cdot \varphi)(v) = \varphi(\hat{S}_{G_q}(x)\cdot v).
		\] 
		for $x\in\DF(G_q)$, $\varphi\in V^*$, $v\in V$.  The locally finite part $F(V^*) \subset V^*$ for the $\DF(K_q)$-action is again a Harish-Chandra module over $G_q$, called the dual Harish-Chandra module.
	\end{definition}
	
	Note that this defines an exact contravariant functor on the category of Harish-Chandra modules over $G_q$ with $\DF(G_q)$-linear morphisms.  Note also that the double dual of $V$ is isomorphic to $V$.  Explicitly, the double dual of $V$ is naturally identified with $V$ as a vector space, equipped with the $\DF(G_q)$-action given by precomposing the original action by $\hat{S}_{G_q}^2$.  Using the relation
	\[
	 (K_{2\rho}\bowtie 1)(x\bowtie f)(K_{-2\rho}\bowtie 1) =\hat{S}^2(x) \bowtie S^2(f) = \hat{S}_{G_q}^2(x \bowtie f)
	\]
	for all $x \bowtie f \in \DF(G_q)$, we see that the two actions are intertwined by the action of $K_{2\rho}$.

\begin{lemma}
		\label{lem:principal_series_bilinear_pairing}
		Let $(\mu,\lambda) \in \weights \times \lie{h}_q^*$.  The Haar integral $\phi$ on $L^2(K_q)$ induces a non-degenerate bilinear pairing
		\[
		\Gamma(\E_{\mu,\lambda}) \times \Gamma(\E_{-\mu,-\lambda}) \to \mathbb{C};  
		\qquad (\xi,\eta) = \phi(\xi\eta)
		\]
		which is $U_q^\mathbb{R}(\lie{g})$-invariant in the sense that for all $X\bowtie f \in U_q^\mathbb{R}(\lie{g})$,
		\[
		((X\bowtie f)\cdot\xi,\eta) = (\xi, \hat{S}_{G_q} (X\bowtie f) \cdot \eta).
		\]
\end{lemma}
	
\begin{proof}
		Non-degeneracy and the fact that $(\hat{S}(X)\cdot \xi,\eta) = ( \xi, X\cdot \eta) $ for $X\in U_q^\mathbb{R}(\lie{k})$ follow from the basic properties of the Haar integral.  For $f,g\in\CF^\infty(K_q)$ we have $S^2(f) = K_{2\rho} \hit f \hitby K_{-2\rho}$, where we are using the notation
		\[
		  X \hit f = (X,f_{(2)})f_{(1)} \qquad f \hitby X = (X,f_{(1)})f_{(2)}.
		\]
		 We also have the modular property $\phi(fg) = \phi((K_{2\rho} \hit g \hitby K_{2\rho})f)$, see Subsection \ref{sec:compact_quantum_groups} and Lemma \ref{lem:F_is_K2rho}.  Therefore, 
		\begin{align*}
		(f \cdot \xi, \eta) &= (K_\lambda, f_{(2)}) \, \phi(K_{2 \rho} \hit f_{(1)} \xi S(f_{(3)}) \eta) \\
		&= (K_\lambda, f_{(2)}) \, \phi(\xi S(f_{(3)}) \eta f_{(1)} \hitby K_{-2 \rho}) \\
		&= (K_{-\lambda}, S(f_{(2)})) \, \phi(\xi S(f_{(3)}) \eta f_{(1)} \hitby K_{-2 \rho}) \\
		&= (K_{-\lambda}, S(f_{(2)})) \, \phi(\xi S(f_{(3)}) \eta K_{-2\rho} \hit S^2(f_{(1)})) \\
		&= (\xi, S(f) \cdot \eta).
		\end{align*}
		This completes the proof.
\end{proof}

It follows that we have an isomorphism between $\Gamma(\E_{\mu,\lambda})$ and the dual Harish-Chandra module of $\Gamma(\E_{-\mu,-\lambda})$, in the sense of Definition \ref{def:dual_HC-module}.

\subsubsection{Relation between principal series modules and Verma modules}

Recall from Lemma \ref{uqrgdiagonalembedding} and Proposition \ref{polynomialtosmooth} that we have an 
embedding $\iota: U_q^\mathbb{R}(\mathfrak{g}) \rightarrow U_q(\mathfrak{g}) \otimes U_q(\mathfrak{g}) $ of algebras 
and an embedding $ \Poly(G_q) \bowtie \Poly(G_q) \rightarrow \CF^\infty(G_q) $ of multiplier Hopf $ * $-algebras. 
These embeddings are compatible with the pairings 
\begin{align*}
 U_q^\mathbb{R}(\mathfrak{g}) \times \CF^\infty(G_q) & \rightarrow \mathbb{C} \\
 (X \bowtie c, a \otimes t) &= (X, a) (c, t) 
\end{align*}
for $ a \otimes t \in \CF^\infty(K_q) \otimes \D(K_q) \subset \CF^\infty(G_q) $  
and $ X \bowtie c \in U_q^\mathbb{R}(\mathfrak{k}) \bowtie \CF^\infty(K_q) = U_q^\mathbb{R}(\mathfrak{g}) $, 
and 
\begin{align*}
 ( U_q(\mathfrak{g}) \otimes U_q(\mathfrak{g}) ) \times ( \Poly(G_q) \bowtie \Poly(G_q) )
    & \rightarrow \mathbb{C}, \\ 
 (X \otimes Y, f \bowtie g) &= (X, g) (Y, f)
\end{align*} 
for $ X \otimes Y \in U_q(\mathfrak{g}) \otimes U_q(\mathfrak{g}) $ and $ f \bowtie g \in \Poly(G_q) \bowtie \Poly(G_q) $, respectively, 
see the remark at the end of Section \ref{subsecuqrg}.

 \begin{lemma}
  \label{lem:Verma-principal-pairing}
   Let $ (\mu, \lambda) \in \weights \times \mathfrak{h}^*_q $ and consider the Yetter-Drinfeld module $ \Gamma(\E_{\mu, \lambda}) \cong \ind_{B_q}^{G_q}(\mathbb{C}_{\mu, \lambda}) \subset \CF^\infty(G_q) $.  
   The pairing $ U_q^\mathbb{R}(\mathfrak{g}) \times \ind_{B_q}^{G_q}(\mathbb{C}_{\mu, \lambda}) \rightarrow \mathbb{C} $ descends to a well-defined bilinear pairing 
   \[ 
    M(\mu, 2 \rho + \lambda) \times \ind_{B_q}^{G_q}(\mathbb{C}_{\mu, \lambda}) \rightarrow \mathbb{C} .
   \]
   This induces an injective $U_q^\mathbb{R}(\mathfrak{g})$-linear map
   \[
    \Gamma(\E_{\mu,\lambda}) \to M(\mu,2 \rho + \lambda)^*
   \]
   where $M(\mu,2 \rho + \lambda)^*$ is equipped with the action given by
 \[
  (X\cdot T)(m) = T(\hat{S}_{G_q}(X)\cdot m)
 \]
 for $X\in U_q^\mathbb{R}(\mathfrak{g})$, $T\in M(\mu,2 \rho + \lambda)^*$, $m\in M(\mu,2 \rho + \lambda)$.
 \end{lemma}

\begin{proof}
Let 
\begin{align*}
 X \bowtie c &\in U_q^\mathbb{R}(\mathfrak{k}) \bowtie \CF^\infty(K_q) = U_q^\mathbb{R}(\mathfrak{g}),\\
 Y \bowtie d &\in U_q^\mathbb{R}(\mathfrak{t}) \bowtie \CF^\infty(K_q) =U_q^\mathbb{R}(\mathfrak{b}), 
\end{align*}
and $a\in \Gamma(\E_\mu)$, so that
\[
 \ext(a) = a \otimes K_{2\rho + \lambda} \in \ind_{B_q}^{G_q}(\mathbb{C}_{\mu, \lambda}).
\]
Note that since $Y\in U_q^\mathbb{R}(\mathfrak{t})$ we have
$$
 (XY, a) = (X \otimes Y, (\id \otimes \pi_T)\Delta(a)) = (X, a) \chi_\mu(Y) .
$$
Using the fact that $U_q^\mathbb{R}(\mathfrak{t})$ is commutative and cocommutative, we obtain
\begin{align*}
((X \bowtie c)&(Y \bowtie d), a \otimes K_{2\rho + \lambda}) \\
&= \left( (Y_{(1)}, c_{(1)})(\hat{S}(Y_{(3)}), c_{(3)}) X Y_{(2)} \bowtie c_{(2)} d, a \otimes K_{2\rho + \lambda} \right) \\
&= (Y_{(1)}, c_{(1)})(\hat{S}(Y_{(3)}), c_{(3)}) (X Y_{(2)}, a)(c_{(2)} d, K_{2\rho + \lambda}) \\
&= (Y_{(1)}, c_{(1)})(\hat{S}(Y_{(3)}), c_{(3)}) (X, a) \chi_\mu(Y_{(2)}) (K_{-2\rho - \lambda},c_{(2)}) (K_{-2\rho - \lambda}, d) \\
&= (Y_{(1)} K_{-2\rho - \lambda} \hat{S}(Y_{(2)}), c) (X, a) \chi_\mu(Y_{(3)}) (K_{-2\rho - \lambda}, d)   \\
&= (X, a) \chi_\mu(Y) (c, K_{2\rho + \lambda}) (d, K_{2\rho + \lambda}) \\
&= \chi_{\mu, 2 \rho + \lambda}(Y \bowtie d) (X \bowtie c, a \otimes K_{2\rho + \lambda}).
\end{align*}
This proves that the pairing is well-defined.  

For the final statement we note that, with the action of $U_q^\mathbb{R}(\mathfrak{g})$ on $\CF^\infty(G_q)$ given by
\[
 X\cdot f = (\hat{S}_{G_q}(X),f_{(1)}) f_{(2)},
\]
the pairing $ U_q^\mathbb{R}(\mathfrak{g}) \times \CF^\infty(G_q) \to \CC $ satisfies
\[
 (XY,f) = (Y, \hat{S}_{G_q}^{-1}(X)\cdot f).
\]
Therefore, the map 
\[
 \CF^\infty(G_q)  \to U_q^\mathbb{R}(\mathfrak{g})^*; \qquad f \mapsto (\bullet, f)
\]
is $U_q^\mathbb{R}(\mathfrak{g})$-linear, and this restricts to the map stated in the lemma.  Moreover, this map is injective by the nondegeneracy of the pairing.
\end{proof}

The above lemma will allow us to identify the principal series module $\Gamma(\E_{\mu,\lambda}) $ with the locally finite part of $(M(l) \otimes M(r))^*$ for well-chosen values of $l$ and $r$ in $\lie{h}_q^*$.   To do so, we need to clarify the relationship between $U_q^\mathbb{R}(\lie{g})$-modules and $FU_q(\lie{g})$-bimodules with compatible adjoint action of $U_q[\lie{g})$ as in Definition \ref{defadcompatibility}.

To begin with, suppose $V$ is any $U_q^\mathbb{R}(\lie{g})$-module.  The diagonal embedding $\iota: U_q^\mathbb{R}(\lie{g}) \to U_q(\lie{g})\otimes U_q(\lie{g})$ of Lemma \ref{uqrgdiagonalembedding} gives an isomorphism
of $U_q^\mathbb{R}(\lie{g})$ onto the subalgebra $(FU_q(\lie{g})\otimes 1) \hat{\Delta} (U_q(\lie{g}))  = (1\otimes \hat{S}^{-1}(FU_q(\lie{g}))) \hat{\Delta}(U_q(\lie{g}))$.
Therefore, we can define an $FU_q(\lie{g})$-bimodule structure on $V$ by
\[
Y\cdot v = \iota^{-1}(Y\otimes 1)v, 
\qquad
v\cdot Z = \iota^{-1}(1\otimes \hat{S}^{-1}(Z))v.
\]
\label{nom:FUq-actions}%
Moreover, we have a $U_q(\lie{g})$-action on $V$ given by the action of $U_q(\lie{g})= U_q^\mathbb{R}(\lie{k}) \subset U_q^\mathbb{R}(\lie{g})$.  As usual, we refer to this $U_q(\lie{g})$-action as the adjoint action and denote it by $ X \rightarrow v $ for $X\in U_q(\lie{g})$ and $v\in V$.  Note that for $X\in FU_q(\lie{g})$ we have
\begin{align*}
  (X_{(1)}\rightarrow v)\cdot X_{(2)} 
    &= \iota^{-1}((1 \otimes \hat{S}^{-1}(X_{(2)})) \hat{\Delta} (X_{(1)})) v \\
    & = \iota^{-1}(X \otimes 1)  v \\
    &= X\cdot v,
\end{align*}
so the actions indeed satisfy the compatibility condition of Definition \ref{defadcompatibility}.

Next, recall from Subsection \ref{sec:FHomMN} that we defined left, right and adjoint actions of $U_q(\lie{g})$ on $(M\otimes N)^*$ by
\begin{align*}
(Y\cdot \varphi\cdot Z)(m\otimes n) &= \varphi(Z\cdot m \otimes \tau(Y)\cdot n), \\
(X \rightarrow \varphi)(m\otimes n) &= \varphi(\hat{S}(X_{(2)})\cdot m \otimes \tau(X_{(1)}) \cdot n),
\end{align*}
\label{nom:adjoint_action_on_MxN2}%
for $X,Y,Z\in U_q(\lie{g})$, where $\tau$ is the algebra antiautomorphism from Lemma \ref{deftau}.  The locally finite part of $(M\otimes N)^*$ with respect to the adjoint action is denoted $F((M\otimes N)^*)$, and the above actions restrict to make $F((M\otimes N)^*)$ into an $FU_q(\lie{g})$-bimodule with compatible adjoint action of $U_q(\lie{g})$, see Lemma \ref{locallyfiniteduallemma} and the discussion preceding it.

These actions correspond to a $U_q^\mathbb{R}(\lie{g})$-module structure on $F((M\otimes N)^*)$.  Explicitly, we can define an action of $U_q^\mathbb{R}(\lie{g})$ on $(M\otimes N)^*$ by first embedding $U_q^\mathbb{R}(\lie{g})$ into $U_q(\lie{g})\otimes U_q(\lie{g})$ via the map $\iota$ and then letting $U_q(\lie{g})\otimes U_q(\lie{g})$ act on $(M\otimes N)^*$ by
\begin{align*}
  ((Y\otimes Z)\varphi)(m\otimes n) &= (Y\cdot \varphi \cdot \hat{S}(Z))(m\otimes n) \\
   &= \varphi(\hat{S}(Z)\cdot m \otimes \tau(Y)\cdot n)
\end{align*}
for $Y,Z\in U_q(\lie{g})$.   By the above discussion, this $U_q^\mathbb{R}(\lie{g})$-action restricts to the locally finite part $F((M\otimes N)^*)$.

With this in place, we have the following result.

\begin{prop} \label{principalseriesversusverma}
Let $ (\mu, \lambda) \in \weights \times \mathfrak{h}^*_q $ and let $l,r \in \mathfrak{h}^*_q$ such that
$$ 
 \mu=l-r, \qquad \lambda+2\rho = -l-r.
$$
\label{nom:lr_shifted}%
Then there exists an isomorphism of $U_q^\mathbb{R}(\lie{g})$-modules
$$
\Gamma(\E_{\mu, \lambda}) \cong F((M(l) \otimes M(r))^*),
$$
where $F((M(l) \otimes M(r))^*)$ is equipped with the $U_q^\mathbb{R}(\lie{g})$-module structure defined above.
\end{prop}

\begin{proof}  
By Corollary \ref{cor:Verma_module_iso} we can identify the Verma module $ M(\mu, \lambda+2\rho) $ with $ M(l) \otimes M(r) $ as a 
$ U_q^\mathbb{R}(\mathfrak{g}) $-module if we equip $M(l)\otimes M(r)$ with the $U_q^\mathbb{R}(\lie{g})$-action given by
\[
  X\cdot (m\otimes n) = ((\id\otimes\theta')\circ\iota'(X)) (m\otimes n),
\]
for $X\in U_q^\mathbb{R}(\lie{g})$. 
Thus, from Lemma \ref{lem:Verma-principal-pairing} we obtain 
an injective linear map $ \gamma: \Gamma(\E_{\mu, \lambda}) \rightarrow (M(l) \otimes M(r))^* $.   

By the $U_q^\mathbb{R}(\lie{g})$-linearity of the map in Lemma \ref{lem:Verma-principal-pairing}, we see that for any $X\in U_q^\mathbb{R}(\lie{g})$ and $f\in \Gamma(\E_{\mu, \lambda})$ we have
\[
  \gamma(X\cdot f)(m\otimes n) = \gamma(f)(\hat{S}_{G_q}(X)\cdot(m\otimes n))
\]
for all $m\otimes n \in M(l)\otimes M(r)$, where the action of $\hat{S}_{G_q}(X)$ on the right-hand side is given by
\begin{align*}
  \hat{S}_{G_q}(X)\cdot(m\otimes n)
    &= (\id\otimes \theta')\circ\sigma\circ(\hat{S}\otimes\hat{S})\circ\iota(X)(m\otimes n) \\
    &= ((\hat{S}\otimes\tau)(\sigma\circ\iota(X)))(m\otimes n),
\end{align*}
see Corollary \ref{cor:Verma_module_iso}.
Thus, if we write $\iota(X) = \sum_j Y_j \otimes Z_j$ we get
\[
 \gamma(X\cdot f)(m\otimes n) = \sum_j \gamma(f)(\hat{S}(Z_j) \cdot m \otimes \tau(Y_j)\cdot n).
\]
Comparing this with the formula for the $U_q^\mathbb{R}(\lie{g})$-action on $(M\otimes N)^*$ given just before the proposition, we see that $\gamma$ is $U_q^\mathbb{R}(\lie{g})$-linear.

Since the action of $ U_q^\mathbb{R}(\mathfrak{k}) = U_q(\mathfrak{g}) $ on $ \Gamma(\E_{\mu, \lambda}) $ is locally finite and $ \gamma $ 
is $ U_q^\mathbb{R}(\mathfrak{g}) $-linear, the image $ \im(\gamma) $ is contained in the locally finite part of $ (M(l) \otimes M(r))^* $.    Moreover, $\gamma$ is surjective onto each $K_q$-type, since Proposition \ref{principalseriesmultiplicity} and Lemma \ref{lem:principal_series_multiplicities} show that for every $\nu\in\weights^+$,
\[
  [ F((M(l)\otimes M(r))^*) : V(\nu)] = \dim (V(\nu)_{-\mu} )
 = [\Gamma(\E_{\mu,\lambda}) : V(\nu)].
\]
This completes the proof.
\end{proof}

Using Lemma \ref{locallyfiniteduallemma}, we can also reformulate
Proposition \ref{principalseriesversusverma} in terms of the spaces $F\Hom(M(l),M(r)^\vee)$.  Here we will state the result in terms of $FU_q(\lie{g})$-bimodules with compatible $U_q(\lie{g})$-action rather than $U_q^\mathbb{R}(\lie{g})$-modules, using the discussion preceding Proposition \ref{principalseriesversusverma}.

\begin{cor}
	\label{principalseriesversusFHom}
 Let $ (\mu, \lambda) \in \weights \times \mathfrak{h}^*_q $ and let $l,r \in \mathfrak{h}^*_q$ such that
 $$ 
  \mu=l-r, \qquad \lambda+2\rho = -l-r.
 $$
 There is an isomorphism of $FU_q(\lie{g})$-bimodules
 \[
  \Gamma(\E_{\mu,\lambda}) \cong F\Hom(M(l),M(r)^\vee)
 \]
 which respects the adjoint actions of $U_q(\lie{g})$.  Here $F\Hom(M(l),M(r)^\vee)$ is equipped with the actions
 \begin{align*}
 (Y\cdot \varphi \cdot Z) (m)(n)
  & = \varphi (Z\cdot m) (\tau(Y)\cdot n), \\
 (X\rightarrow \varphi)(m)(n) 
  & = \varphi (\hat{S}(X_{(2)})\cdot m) (\tau(X_{(1)})\cdot n)
 \end{align*}
 \label{nom:FHom_bimodule}%
  for $Y,Z\in FU_q(\lie{g})$, $X\in U_q(\lie{g})$, $m\in M(l)$ and $n\in M(r)$.  
\end{cor}

In the sequel we will often identify $ \Gamma(\E_{\mu, \lambda}) $ and $F\Hom(M(l),M(r)^\vee)$, with the parameters matching as in Corollary \ref{principalseriesversusFHom}.

 As we noted in the remarks after Corollary \ref{cor:Verma_module_iso}, the parameters $(l,r)$ and $(l',r')$ correspond to the same pair $(\mu,\lambda) \in \weights \times \lie{h}_q^*$ iff they differ by the diagonal action of the translation subgroup $\mathbf{Y}_q= \half i \hbar^{-1} \roots^\vee / i \hbar^{-1} \roots^\vee$ of the extended Weyl group $\hat{W} = \mathbf{Y}_q \rtimes W$, that is, iff
 \[
  \textstyle (l',r') = (l+\half i \hbar^{-1} \alpha^\vee , r+\half i \hbar^{-1} \alpha^\vee)
 \]
 for some $\alpha^\vee \in \roots^\vee$.

We also note that by applying Lemma \ref{locallyfinitedualtransposelemma},
we have an isomorphism of $U_q(\lie{g})$-compatible $FU_q(\lie{g})$-bimodules
\[
\Gamma(\E_{\mu,\lambda}) \cong F\Hom(M(r),M(l)^\vee)
\]
where now $F\Hom(M(r),M(l)^\vee)$ should be equipped with the actions
\begin{align*}
 (Y\cdot \varphi \cdot Z) (m)(n)
   & = \varphi (\tau(Y)\cdot m) (Z\cdot n), \\
(X\rightarrow \varphi)(m)(n) 
   & = \varphi (\tau(X_{(1)})\cdot m) (\hat{S}(X_{(2)})\cdot n)
\end{align*}
\label{nom:FHom_bimodule2}%
for  $Y,Z\in FU_q(\lie{g})$ and $X\in U_q(\lie{g})$.

\subsection{An equivalence of categories} \label{seccategoryequiv}

In this section we discuss the relation between certain categories of $FU_q(\lie{g})$-bimodules, which are related to the Harish-Chandra modules of Section \ref{sec:G_q-representations}, and subcategories of 
category $ \O $, following Chapter 8.4 in \cite{Josephbook}. The results are due to Joseph and Letzter \cite{JLVermaannihilator}. 
The corresponding theory in the classical setting was developed independently by Bernstein-Gelfand in \cite{BGtensorproducts}, Enright \cite{Enrightlectures},  
and Joseph \cite{Josephdixmierproblem}.

We first introduce certain subcategories of category $ \O $.

\begin{definition}
For $ l \in \mathfrak{h}^*_q $ let $ \O_l $ be the full subcategory of $ \O $ consisting of all modules with weight spaces associated to 
weights in $ l + \weights \subset \mathfrak{h}^*_q $. 
\nomenclature[o$O_l$]{$\O_l$}{subcategory of $\O$ consisting of modules with weights in $ l + \weights$}%
\end{definition} 

Comparing with the set-up in \cite{Josephbook}, note that for all modules $ N $ in $ \O $ the annihilator $ \ann_{ZU_q(\mathfrak{g})}(N) $ 
has finite codimension in $ ZU_q(\mathfrak{g}) $. This is due to the fact that we require modules in $ \O $ to be finitely generated.

Next we introduce a category of $U_q(\lie{g})$-compatible $FU_q(\lie{g})$-bimodules as in Definition \ref{defadcompatibility}, but with some additional structural hypotheses.

\begin{definition} \label{defhcbimodule}
Let $ H $ be a $ FU_q(\mathfrak{g}) $-bimodule with a compatible locally finite action of $ U_q(\mathfrak{g}) $, meaning that
\[
 (X_{(1)} \rightarrow v)\cdot X_{(2)} = X\cdot v
\]
for all $X\in FU_q(\lie{g})$ and $v\in V$. We say 
that $ H $ is a Harish-Chandra bimodule if the following conditions are satisfied. 
\begin{bnum} 
\item[a)] All isotypical components for the adjoint action of $ U_q(\mathfrak{g}) $ are finite dimensional and integrable.
\item[b)] The annihilators of both the left and right actions of $ ZU_q(\mathfrak{g}) $ have finite codimension. 
\item[c)] $ H $ is finitely generated as a right $ FU_q(\mathfrak{g}) $-module. 
\end{bnum}

We write $\HC$ for the category of Harish-Chandra bimodules with morphisms being $ FU_q(\mathfrak{g}) $-bimodule maps which are also $ U_q(\mathfrak{g}) $-linear.
\nomenclature[o$HC$1]{$\HC$}{category of Harish-Chandra bimodules}%
\end{definition}

Let us observe that every irreducible Harish-Chandra module $ V $ naturally defines a Harish-Chandra bimodule. Indeed, according to Definition \ref{defadmissible} Harish-Chandra modules are admissible as $\DF(K_q)$-modules, so that
the isotypical components of $ V $ are finite dimensional and integrable by assumption in this case. Irreducibility implies that the annihilators for the 
action of left and right action of $ ZU_q(\mathfrak{g}) $ have finite codimension because the centre must act by fixed left/right central characters. 
If $ V $ is simple then every finite dimensional $ \DF(K_q) $-submodule $ W $ of $ V $ generates $ V $ as a $U_q^\RR(\lie{g})$-module, 
and the compatibility condition in Definition \ref{defadcompatibility} implies that $ W $ generates $ V $ even as a right $ FU_q(\mathfrak{g}) $-module.

Recall that $ \xi_l: ZU_q(\mathfrak{g}) \rightarrow \mathbb{C} $ is the central character corresponding to $ l \in \lie{h}_q^* $  as in Definition \ref{def:central_character}, and note 
that $ \ker(\xi_l) $ annihilates $ M(l) $.

\begin{definition} 
For $ l \in \mathfrak{h}^*_q $ let $ \HC_l $
\nomenclature[o$HC$2]{$\HC_l$}{subcategory of $\HC$ associated with the central character $\xi_l$}%
be the full subcategory of $ \HC $ consisting of all objects for which the 
annihilator of the right action of $ ZU_q(\mathfrak{g}) $ contains $ \ker(\xi_l) $. 
\end{definition}

We shall now relate Verma modules and compatible bimodules in the Harish-Chandra category, at least with suitable extra conditions 
on both sides. 

\begin{prop} 
	\label{prop:F_l}
For $ l \in \mathfrak{h}^*_q $ we obtain a covariant functor $ \F_l: \O_l \rightarrow \HC_l $ by setting 
$$
\F_l(M) = F \Hom(M(l), M)
$$
\nomenclature[o$F_l$]{$\F_l$}{functor $ \F_l: \O_l \rightarrow \HC_l $}%
for $ M \in \O_l $, equipped with the compatible bimodule structure
	\begin{align*}
	(Y\cdot \varphi \cdot Z)(m)
	& = Y\cdot\varphi (Z\cdot m) \\
	(X\rightarrow \varphi)(m) 
	& = X_{(1)}\cdot\varphi (\hat{S}(X_{(2)})\cdot m)
	\end{align*}
for $Y,Z \in FU_q(\lie{g})$, $X\in U_q(\lie{g})$, $\varphi\in F\Hom(M(l),M)$ and $m\in M(l)$. 
\end{prop}

\begin{proof}  Let us check that $ \F_l(M) $ is indeed contained in $ \HC_l $. 
Firstly, the adjoint action of $ U_q(\mathfrak{g}) $ is clearly locally finite on the right hand side. 
Moreover, since the weights of $M$ belong to $l+\weights$, one easily checks that the adjoint action has weights in $\weights$, so is integrable.
The annihilator of the left action of $ ZU_q(\mathfrak{g}) $ has finite codimension because $ M \in \O $, 
and the annihilator of the right action of $ ZU_q(\mathfrak{g}) $ has finite codimension because it contains $ \ker(\xi_l) $. 

For the finite multiplicity requirement, note that modules in category $\O$ have finite length by Theorem \ref{catoartinian}, and since every simple module embeds in a dual Verma module it is therefore enough to consider the case that $ M = M(r)^\vee $ is a dual Verma module. 
In this case it follows from Lemma  \ref{locallyfiniteduallemma} and Proposition \ref{principalseriesmultiplicity} $a)$. 

Hence $ \F_l(M) $ is indeed contained in $ \HC_l $. Moreover $ \F_l $ clearly maps $ U_q(\mathfrak{g}) $-linear maps 
to $ U_q(\mathfrak{g}) $-linear $ FU_q(\mathfrak{g}) $-bimodule maps. \end{proof}

Note that when $M=M(r)^\vee$ for $r\in\lie{h}^*_q$, the bimodule structure in Proposition \ref{prop:F_l} is exactly the same as that from Corollary \ref{principalseriesversusFHom}.  In other words, with $\mu=l-r$ and $\lambda+2\rho=-l-r$ we have
\[
 F_l(M(r)^\vee) \cong \Gamma(\E_{\mu,\lambda}).
\]

Let $ F $ be the functor on $ U_q(\mathfrak{g}) $-modules that assigns to $ M $ its locally finite part $ FM $. 
This is \emph{not} an exact functor. If $ 0 \rightarrow K \rightarrow E \rightarrow Q \rightarrow 0 $ is exact then clearly 
$ 0 \rightarrow FK \rightarrow FE \rightarrow FQ \rightarrow 0 $ is exact at $ FK $ and $ FE $. Surjectivity 
may fail however. Consider for instance $ E = M(\mu), Q = V(\mu), K = I(\mu) $ for $ \mu \in \weights^+ $. 
Then $ FE = 0 $ and $ FQ = Q $. 

Despite this fact we have the following result.  

\begin{lemma} \label{Flexact}
Let $ l \in \mathfrak{h}^*_q $ be dominant. Then $ \F_l: \O_l \rightarrow \HC_l $ is an exact functor. 
\end{lemma} 

\begin{proof}  
Recall from Definition \ref{defdominantantidominant} that $l\in\lie{h}_q^*$ is called dominant if 
$q_\alpha^{(l+\rho,\alpha^\vee)} \notin \pm q_\alpha^{-\mathbb{N}} $
for all
$	\alpha\in{\bf\Delta}^+$,
where $\mathbb{N} = \{1,2,\ldots\}$.
According to Proposition \ref{Vermaprojective}, dominance of $ l $ implies that $ M(l) $ is projective. 
Hence the invariant part 
of the induced sequence $ 0 \rightarrow \Hom(M(l), K) \rightarrow \Hom(M(l), E) \rightarrow \Hom(M(l), Q) \rightarrow 0 $
is exact. Tensoring with a simple module $ V(\mu) $ for $ \mu \in \weights^+ $ we see that this applies to all $ \mu $-isotypical 
components. Indeed, the $ \mu $-isotypical component of $ \Hom(M(l), N) $ is the invariant part of 
\begin{align*}
\Hom(M(l), N) \otimes V(\mu)^* 
\cong \Hom(M(l), N \otimes V(\mu)^*), 
\end{align*}
so exactness follows again from projectivity of $ M(l) $. \end{proof}

Our next aim is to define a functor in the reverse direction. For this we need some auxiliary considerations. 

Firstly, assume that $ H $ is any $ FU_q(\mathfrak{g}) $-bimodule with a compatible action of $ U_q(\mathfrak{g}) $, 
and let $ V $ be a finite dimensional integrable $ U_q(\mathfrak{g}) $-module. 
Consider the right $ FU_q(\mathfrak{g}) $-module structure 
on the second factor of $ V \otimes H $ and the left $ U_q(\mathfrak{g}) $-module structure given by the diagonal action
\[
  X \rightarrow (v \otimes  h) = (X_{(1)}\cdot v) \otimes (X_{(2)} \rightarrow h)
\]
for $X \in U_q(\lie{g})$.
Let us define a left $ FU_q(\mathfrak{g}) $-module structure on $ V \otimes H $ by setting 
$$
Y \cdot (v \otimes h) = (Y_{(1)} \cdot v) \otimes Y_{(2)} \cdot h,
$$
for $Y \in FU_q(\lie{g})$.  This makes sense since $FU_q(\lie{g})$ is a left coideal, see Lemma \ref{lem:FUq_coideal}.
Then we compute 
\begin{align*}
(Y_{(1)} \rightarrow (v \otimes h)) \cdot Y_{(2)} &= (Y_{(1)} \cdot v) \otimes (Y_{(2)} \rightarrow h) \cdot Y_{(3)} \\
&= (Y_{(1)} \cdot v) \otimes (Y_{(2)} \cdot h) \\
&= Y \cdot (v \otimes h), 
\end{align*}
so that the resulting $ FU_q(\mathfrak{g}) $-bimodule $ V \otimes H $ is compatible with the $ U_q(\mathfrak{g}) $-action. 

\begin{prop} 
Let $ V $ be a finite dimensional integrable $ U_q(\mathfrak{g}) $-module and assume that $ H \in \HC $. 
Then the compatible $ FU_q(\mathfrak{g}) $-bimodule $ V \otimes H $ constructed above is again contained in $ \HC $. 
For $ H \in \HC_l $ we have $ V \otimes H \in \HC_l $.  
\end{prop} 

\begin{proof}  We have to verify that the conditions in Definition \ref{defhcbimodule} are satisfied for $ V \otimes H $. For condition $ a) $ notice that $ V $ is 
locally finite and integrable for the left $ U_q(\mathfrak{g}) $-module structure, and that local finiteness and integrability are preserved under tensor products. 

The key point to check is condition $ b) $. For this we may assume that $ V = V(\mu) $ is simple. 
The annihilator of $ ZU_q(\mathfrak{g}) $ for the left action on $ H $ 
being of finite codimension implies that we have 
\[
 \ann_{ZU_q(\mathfrak{g})}(H) \supseteq (\ker(\xi_{l_1}) \cdots \ker(\xi_{l_m}))^r 
\]
for some weights $ l_1, \dots, l_m \in \mathfrak{h}^*_q $ and $ r \in \mathbb{N} $. 
Note that, thanks to the decomposition series of $V\otimes M(l_j)$ from Lemma \ref{Vermatensorfiltration}, we have
\begin{align*}
 \ann_{ZU_q(\mathfrak{g})}(V(\mu)\otimes M(l_j))
  &\supseteq \prod_{\nu \in \weights(V(\mu))} (\ann_{ZU_q(\mathfrak{g})} M(l_j+\nu))^s \\
   &= \prod_{\nu \in \weights(V(\mu))} (\ker(\xi_{l_j+\nu}))^s
\end{align*}
for sufficiently large $s$.  It follows that 
\[
\ann_{ZU_q(\mathfrak{g})}(V(\mu)\otimes H)
 \supseteq \prod_j \prod_{\nu \in \weights(V(\mu))} \ker(\xi_{l_j+\nu})^{s'} \]
for some $s'$ and hence $\ann_{ZU_q(\mathfrak{g})}(V(\mu)\otimes H)$ has finite codimension.

Finite codimension for the annihilator of the right action of $ ZU_q(\mathfrak{g}) $ is obvious since $ H \in \HC $.

Condition $ c) $ is evidently satisfied because the right action of $ FU_q(\mathfrak{g}) $ on $ V \otimes H $ 
only sees the factor $ H $. We conclude that $ V \otimes H $ is contained in $ HC $. 

The final assertion regarding $ \HC_l $ is again immediate because the right action of $ ZU_q(\mathfrak{g}) $ 
depends only on the second tensor factor of $ V \otimes H $.
\end{proof}

Let $ H \in \HC $ be given. Since $ H $ is finitely generated as a right $ FU_q(\mathfrak{g}) $-module we find a finite dimensional integrable $ U_q(\mathfrak{g}) $-submodule $ V \subset H $ such that the canonical map $ V \otimes FU_q(\mathfrak{g}) \rightarrow H $ is surjective.  

If $ H \in \HC_l $ then 
the right annihilator of $ H $ contains $ \ker(\xi_l) \subset ZU_q(\mathfrak{g}) $, and hence 
Theorem \ref{Vermaannihilator} implies that the the projection $ V \otimes FU_q(\mathfrak{g}) \rightarrow H $ 
factorizes through $ V \otimes F \End(M(l)) $.

\begin{lemma} \label{lemma:Tl_functor}
Let $ l \in \mathfrak{h}^*_q $ and let $ H \in \HC_l $. 
The left action of $ FU_q(\mathfrak{g}) $ on the first leg of $ H \otimes_{FU_q(\mathfrak{g})} M(l) $ extends to a left action of $ U_q(\mathfrak{g}) $ 
via the formula
\[
  X\cdot (h\otimes m) = X_{(1)}\rightarrow h \otimes X_{(2)}\cdot m,
\]
and we obtain a functor $ \T_l: \HC_l \rightarrow \O_l $ by setting 
$$
\T_l(H) = H \otimes_{FU_q(\mathfrak{g})} M(l)
$$
\nomenclature[o$T_l$]{$\T_l$}{functor $ \T_l: \HC_l \rightarrow \O_l $}%
for $ H \in \HC_l $. 
\end{lemma}

\begin{proof} 
Note that when $X\in FU_q(\lie{g})$, we have
\[
 (X\cdot h)\otimes m  =  (X_{(1)}\rightarrow h)\cdot X_{(2)} \otimes m = X_{(1)}\rightarrow h \otimes X_{(2)}\cdot m,
\]
for all $h\otimes m \in H \otimes_{FU_q(\lie{g})} M(l)$.  We need to show that the expression on the right-hand side remains well-defined on the balanced tensor product when $X\in U_q(\lie{g})$.

Since $H$ is a Harish-Chandra bimodule, we have a direct sum decomposition $H = \bigoplus_{\mu\in \weights} H_\mu$ into weight spaces for the adjoint action.  If $h\in H_\mu$, $m\in M(l)_{\lambda}$
then the above calculation shows $(K_\nu\cdot h)\otimes m = q^{(\mu+\lambda,\nu)} h\otimes m$
for every $\nu\in2\weights^+$. It follows that these elements $K_{\nu}\in FU_q(\lie{g})$ are simultaneously diagonalizable on $H\otimes_{FU_q(\lie{g})} M(l)$ with eigenspaces
\[
  (H\otimes_{FU_q(\lie{g})} M(l))_{\gamma} =
  \bigoplus_{\mu+\lambda = \gamma} H_\mu \otimes M(l)_{\lambda}.
\]
Now it makes sense to define the action of any $K_\nu\in U_q(\lie{h})$ by
$K_\nu \cdot v = q^{(\gamma,\nu)} v$ whenever $v\in (H\otimes_{FU_q(\lie{g})} M(l))_{\gamma}$, and this agrees with the formula stated in the lemma.  By Lemma \ref{uqgequalsfuqguqh}, this action is well-defined for all $X\in U_q(\lie{g})$.

The resulting $ U_q(\mathfrak{g}) $-module $ \T_l(H) $ is clearly a weight module with weights contained in 
$ l + \weights $. According to our considerations above there exists some finite dimensional integrable $ U_q(\mathfrak{g}) $-module $ V $ 
such that $ H $ can be written as a quotient of $ V \otimes FU_q(\mathfrak{g}) $. 
Hence $ \T_l(H) = H \otimes_{FU_q(\mathfrak{g})} M(l) $ is a quotient of $ V \otimes M(l) $, 
and the latter is in category $ \O_l $. It follows that $ \T_l(H) $ is contained in $ \O_l $. 

It is clear that $ \T_l $ maps morphisms in $ \HC_l $ to morphisms in $ \O_l $. \end{proof}  

We obtain a Frobenius reciprocity relation for the functors $ \F_l $ and $ \T_l $. 

\begin{prop} \label{fltladjointduality}
For $ H \in \HC_l $ and $ N \in \O_l $ there is a natural isomorphism
$$
\Hom_{U_q(\mathfrak{g})}(\T_l(H), N) \cong \Hom_{\HC}(H, \F_l(N)), 
$$
that is, the functor $ \T_l $ is left adjoint to $ \F_l $. 
\end{prop} 

\begin{proof}  The standard $ \Hom $-tensor adjunction yields an isomorphism 
$$ 
\phi: \Hom_{FU_q(\mathfrak{g})\text{-} FU_q(\mathfrak{g})}(H, \Hom(M(l), N)) \rightarrow
\Hom_{FU_q(\mathfrak{g})}(H \otimes_{FU_q(\mathfrak{g})} M(l), N), 
$$ 
sending $ f \in \Hom_{FU_q(\mathfrak{g})\text{-} FU_q(\mathfrak{g})}(H, \Hom(M(l), N)) $ to $\phi(f)$ given by $ \phi(f)(h \otimes m) = f(h)(m) $. 
The inverse isomorphism is given by $ \psi(g)(h)(m) = g(h \otimes m) $.  

Let us show that $ \phi $ maps morphisms in $ \HC $ to $ U_q(\mathfrak{g}) $-linear maps and $ \psi $ maps 
$ U_q(\mathfrak{g}) $-linear maps to morphisms in $ \HC $. 
We do this by verifying that the unit and counit of the above adjunction have these properties. For $ H \in \HC_l $ the unit 
$ u: H \rightarrow \F_l \T_l(H) = F \Hom(M(l), H \otimes_{FU_q(\mathfrak{g})} M(l)) $ is given by $ u(h)(m) = h \otimes m $. 
We compute 
\begin{align*}
u(X \rightarrow h)(m) &= (X \rightarrow h) \otimes m \\
&  = (X_{(1)}\rightarrow h) \otimes (X_{(2)} \hat{S}(X_{(3)})) m \\ 
&= X_{(1)} \cdot (h \otimes \hat{S}(X_{(2)}) \cdot m) \\
&= X_{(1)} \cdot u(h)(\hat{S}(X_{(2)}) \cdot m) = (X \rightarrow u(h))(m).
\end{align*}
The counit $ e: \T_l \F_l(M) = F\Hom(M(l), M) \otimes_{FU_q(\mathfrak{g})} M(l) \rightarrow M $ is given by $ e(T \otimes m) = T(m) $.  For any $\nu\in\weights$ we have 
\begin{align*}
e(K_\nu \cdot (T \otimes m)) &= q^{(\nu, \mu + \lambda)} e(T \otimes m) \\
&= q^{(\nu, \mu + \lambda)} T(m) \\
&= (K_\nu \rightarrow T)(K_\nu \cdot m) \\
&= K_\nu \cdot T(K_{\nu}^{-1} \cdot (K_\nu \cdot m)) \\
&= K_\nu \cdot T(m) \\
&= K_\nu \cdot e(T \otimes m) 
\end{align*}
if $ T $ has weight $ \mu $ and $ m $ has weight $ \lambda $.  In light of Lemma \ref{uqgequalsfuqguqh}, this suffices to conclude that $e$ is $U_q(\lie{g})$-linear.

We conclude that the map $ \phi $ constructed above restricts to an isomorphism
\begin{align*}
\Hom_{\HC}(H, F\Hom(M(l), N)) \rightarrow
\Hom_{U_q(\mathfrak{g})}(H \otimes_{FU_q(\mathfrak{g})} M(l), N)
\end{align*}
as desired. \end{proof}  

\begin{prop} \label{unitisomorphism}
Assume $ l \in \mathfrak{h}^*_q $ is dominant. 
Then the canonical map $ H \rightarrow \F_l \T_l(H) $ induced from the adjunction of $ \F_l $ and $ \T_l $ is an isomorphism 
for $ H \in \HC_l $. 
\end{prop} 

\begin{proof}  First, take $ H = F \End(M(l)) $.  
By the Verma module annihilator theorem \ref{Vermaannihilator} we have $ F \End(M(l)) = FU_q(\mathfrak{g})/\ann_{FU_q(\mathfrak{g})}(M(l)) $  and so the canonical multiplication map $ F \End(M(l)) \otimes_{FU_q(\mathfrak{g})} M(l) \rightarrow M(l) $ is an isomorphism.   It follows that 
$$ 
\T_l(H) = F \End(M(l)) \otimes_{FU_q(\mathfrak{g})} M(l) \cong M(l).
$$  
Hence we get $ \F_l \T_l(H) = \F_l(M(l)) = F \End(M(l)) = H $ in this case. 

Next assume that $ H = V \otimes F \End(M(l)) $ for some finite dimensional integrable $ U_q(\mathfrak{g}) $-module $ V $. 
Then by our previous considerations
$$ 
\T_l(H) \cong (V \otimes F \End(M(l))) \otimes_{FU_q(\mathfrak{g})} M(l) = V \otimes \T_l(F \End(M(l))) = V \otimes M(l), 
$$ 
and hence 
\begin{align*}
\F_l \T_l(H) &\cong \F_l(V \otimes M(l)) \\
&= F\Hom(M(l), V \otimes M(l)) \\
&\cong V \otimes F\End(M(l)) .
\end{align*}
This means $ \F_l \T_l(H) \cong H $ in this case as well. 

Finally, suppose that $ H \in \HC_l $ is arbitrary. According to 
the discussion preceding Lemma \ref{lemma:Tl_functor},
there exists a finite dimensional integrable $ U_q(\mathfrak{g}) $-module  
$ V_1 $ and a surjective homomorphism $ V_1 \otimes F \End(M(l)) \rightarrow H $. Moreover, as $ FU_q(\mathfrak{g}) $ is Noetherian by 
Theorem \ref{fuqgnoetherian}, the kernel of the projection map is again in $ \HC_l $. Hence we obtain a short exact sequence of the form
$$
\xymatrix{
V_2 \otimes F \End(M(l)) \ar@{->}[r] & V_1 \otimes F \End(M(l)) \ar@{->}[r] & H \ar@{->}[r] & 0
     }
$$
Since $ \F_l $ is exact by Lemma \ref{Flexact} and $ \T_l $ is right exact we obtain a commutative diagram 
$$
\xymatrix{
V_2 \otimes F \End(M(l)) \ar@{->}[r] \ar@{->}[d]^\cong & V_1 \otimes F \End(M(l)) \ar@{->}[r] \ar@{->}[d]^\cong & H \ar@{->}[r] \ar@{->}[d] & 0 \\
V_2 \otimes F \End(M(l)) \ar@{->}[r] & V_1 \otimes F \End(M(l)) \ar@{->}[r] & \F_l \T_l(H) \ar@{->}[r] & 0
     }
$$
with exact rows. The $ 5 $-lemma shows that the right hand vertical arrow is an isomorphism. \end{proof}

\begin{lemma} \label{simplicitylemma}
Assume that $ l \in \mathfrak{h}^*_q $ is dominant. If $ V \in \O_l $ is simple then $ \F_l(V) = 0 $ or $ \F_l(V) $ is simple. All simple objects 
in $ \HC_l $ are obtained in this way. 
\end{lemma}

\begin{proof}  Assume $ \F_l(V) \neq 0 $ and let $ H \subset \F_l(V) $ be a nonzero submodule. If $ i: H \rightarrow \F_l(V)  $ denotes the embedding map then 
the map $ j: \T_l(H) \rightarrow V $ corresponding to $ i $ is nonzero according to Proposition \ref{fltladjointduality}. Since $ V $ is simple 
it follows that $ j $ is a surjection. Since $ \F_l $ is exact by Lemma \ref{Flexact} we conclude that $ \F_l(j): \F_l \T_l(H) \rightarrow \F_l(V) $ 
is a surjection. Composing the latter with the isomorphism $ H \rightarrow \F_l \T_l(H) $ we reobtain our original map $ i $. Hence $ i $ is surjective, 
which means that $ H = \F_l(V) $. Hence $ \F_l(V) $ is simple. 

Now let $ H \in \HC_l $ be an arbitrary simple object. Then we have $ H \cong \F_l \T_l(H) $ by Proposition \ref{unitisomorphism}, so that $ \T_l(H) $ is nonzero. 
Since every object in $ \O_l $ has finite length we find a simple quotient $ V $ of $ \T_l(H) $. Then $ \F_l(V) $ is simple, 
and the quotient map $ \T_l(H) \rightarrow V $ corresponds to a nonzero homomorphism $ H \rightarrow \F_l(V) $. Since both $ H $ and $ \F_l(V) $ are simple this 
means $ H \cong \F_l(V) $. \end{proof}  

In fact, in Proposition \ref{domregsimpleobjects} below we will see that with additional hypotheses on $l$, the bimodule $\F_l(V)$ is never zero for $V$ simple, so that $\F_l$ maps simple objects to simple objects.
For this we need to introduce the concept of a regular weight and discuss the translation principle.

Let $ \lambda \in \mathfrak{h}^*_q $. Recall from Section \ref{secdefo} 
that $ \hat{W}_{[\lambda]} = \{\hat{w} \in \hat{W} \mid \hat{w} . \lambda - \lambda \in \roots \subset \mathfrak{h}^*_q/\frac{1}{2} i\hbar^{-1} \roots^\vee \} $. 
We adopt the following definition from Section 8.4.9 in \cite{Josephbook}. 

\begin{definition} 
A weight $ \lambda \in \mathfrak{h}^*_q $ is called regular if 
the only element $ \hat{w} \in \hat{W}_{[\lambda]} $ with $ \hat{w} . \lambda = \lambda $ is $\hat{w} = e $.
\end{definition}

Observe that 
any element $\hat{w} \in \hat{W}$ which satisfies $\hat{w}.\lambda = \lambda$ is automatically in $\hat{W}_{[\lambda]}$, so the definition is equivalent to demanding that $e$ is the only element of $\hat{W}$ which fixes $\lambda$.
We conclude that 
the group $ \hat{W} $ acts freely on the set of all regular elements in $ \mathfrak{h}^*_q $. 

\begin{lemma} \label{dominantregularweightlemma}
Let $ \mu \in \weights^+ $. If $ l \in \mathfrak{h}^*_q $ is dominant and regular then 
$$ 
\hat{W} . (l + \mu) \cap (l + \weights(V(\mu))) = \{l + \mu\}, 
$$ 
where $ \weights(V(\mu)) \subset \weights $ denotes the set of all weights of $ V(\mu) $. 
\end{lemma}

\begin{proof}  It is clear that $ l + \mu $ is contained in $ \hat{W} . (l + \mu) \cap (l + \weights(V(\mu))) $. 

Conversely, assume $ \hat{w} \in \hat{W} $ satisfies $ \hat{w} . (l + \mu) \in (l + \weights(V(\mu))) $. Then $ \hat{w} \in \hat{W}_{[l]} $ since 
$$
\hat{w}(l + \mu + \rho) \in l + \mu + \rho + \roots 
$$
and $ \hat{W}_{[l + \mu + \rho]} = \hat{W}_{[l]} $. Since $ l $ is dominant it is maximal in its $ \hat{W}_{[l]} $-orbit according to 
Proposition \ref{dominantchar}, which means that we have $ \hat{w}^{-1} . l = l - \gamma $ for some $ \gamma \in \roots^+ $, or 
equivalently $ \hat{w}(l - \gamma + \rho) = l + \rho $. From the latter relation we conclude
$$ 
\hat{w} . (l + \mu) = \hat{w}(l + \mu + \rho) - \rho = l + \hat{w}(\mu + \gamma), 
$$
which implies in particular $ \hat{w}(\mu + \gamma) \in \weights(V(\mu)) $. 
Moreover, since $ \gamma \in \roots^+ $ and $ \mu \in \weights^+ $ we have $ (\mu, \gamma) \geq 0 $, and we get  
$$ 
(\hat{w}(\mu + \gamma), \hat{w}(\mu + \gamma)) = (\mu + \gamma, \mu + \gamma) = (\mu, \mu) + 2 (\mu, \gamma) + (\gamma, \gamma) \geq (\mu, \mu) 
$$
with a strict inequality iff $ \gamma \neq 0 $. However, for $ \nu \in \weights(V(\mu)) $ we have $ (\nu, \nu) \leq (\mu, \mu) $, see 
Lemma 1.7 in \cite{Duflocomplexlnm}. 
Hence $ \gamma = 0 $ and we conclude $ \hat{w}^{-1} . l = l $. This implies $ \hat{w} = e $ by regularity of $ l $. \end{proof}  

Let $ M \in \O $. For a central character $ \chi: ZU_q(\mathfrak{g}) \rightarrow \mathbb{C} $, let us define the $ \chi $-primary component of $ M $ by 
$$
M^\chi = \{m \in M \mid \text{for all } Z \in ZU_q(\mathfrak{g}) \text{ we have } (Z - \chi(Z))^n \cdot m = 0 \text{ for some } n \in \mathbb{N} \}, 
$$
\nomenclature{$M^\chi$}{$\chi$-primary component of a module $M\in\O$}%
compare Section 1.12 in \cite{HumphreysO}. Then $ M^\chi \subset M $ is a $ U_q(\mathfrak{g}) $-submodule, and as in the classical case 
one checks that $ M $ decomposes into a finite direct sum of its primary components.

Using this fact we obtain a direct sum decomposition of $ \O $ into the full subcategories $ \O^\chi $
\nomenclature[o$O_l$]{$\O^\chi$}{primary component of category $\O$ corresponding to a central character $\chi$}%
of modules for which $ M = M^\chi $. 
We also refer to $ \O^\chi $ as the primary component corresponding to $ \chi $.

Note that according to Proposition \ref{uqgcentralcharacters} every central character is of the 
form $ \chi = \xi_\lambda $ for some $ \lambda \in \mathfrak{h}^*_q $ in our setting. 

\begin{prop} \label{translationprinciple}
Let $ \mu \in \weights^+ $ and let $ l \in \mathfrak{h}^*_q $ be dominant and regular. Then for all $ \hat{w} \in \hat{W} $ the Verma module 
$ M(\hat{w} . (l + \mu)) $ is a direct summand of $ V(\mu) \otimes M(\hat{w} . l) $. 
\end{prop} 

\begin{proof}  
Consider the primary component of $ M = V(\mu) \otimes M(\hat{w} . l) $ corresponding to $ \xi_{l + \mu}: ZU_q(\mathfrak{g}) \rightarrow \mathbb{C} $. 
All Verma modules in the primary component $ \O^{\xi_{l + \mu}} $ have central character $ \xi_{l + \mu} $, and hence highest weights in $ \hat{W} . (l + \mu) $. 
By Lemma \ref{Vermatensorfiltration}, $ M $ admits a Verma flag with highest weights of the form $ \hat{w} . l + \eta $ for $ \eta \in \weights(V(\mu)) $. 
Hence the Verma modules occurring as subquotients of $M$ for this flag
are of the form $ M(\nu) $ for $ \nu \in (\hat{W} . (l + \mu)) \cap (\hat{w} . l + \weights(V(\mu))) $.
Since $ \weights(V(\mu)) $ 
is stable under the (unshifted) action of $ W $ we obtain 
$$
\hat{w} . l + \weights(V(\mu)) = \hat{w}(l + \rho) - \rho + \weights(V(\mu)) = \hat{w}(l + \weights(V(\mu)) + \rho) - \rho = \hat{w} . (l + \weights(V(\mu)))
$$
Since clearly $ \hat{W} . (l + \mu) = \hat{w} . \hat{W} . (l + \mu) $ we deduce 
$$
(\hat{W} . (l + \mu)) \cap (\hat{w} . l + \weights(V(\mu))) = \hat{w} . (\hat{W} . (l + \mu) \cap (l + \weights(V(\mu)))) = \hat{w} .(l + \mu) 
$$ 
by Lemma \ref{dominantregularweightlemma}. Hence the only possibility is $ \nu = \hat{w} .(l + \mu) $. The result follows. \end{proof}  

We shall now show that $ \F_l $ maps simple objects to simple objects in the case that $ l $ is regular. 

\begin{prop} \label{domregsimpleobjects}
Assume $ l \in \mathfrak{h}^*_q $ is dominant and regular. Then $ \F_l(V) $ is simple for each simple module $ V \in \O_l $. 
\end{prop} 

\begin{proof}  We have to show that $ F \Hom(M(l), V) $ is nonzero if $ V $ is a simple module. 
For this we may assume $ V = V(l + \lambda) $ for some $ \lambda \in \weights $. According to 
Proposition \ref{translationprinciple} we have that $ M(l + \mu) $ is a direct summand of $ V(\mu) \otimes M(l) $ for any $ \mu \in \weights^+ $. 
Since $ F\Hom(M(l), V) \otimes V(\mu)^* \cong F\Hom(V(\mu) \otimes M(l), V) $ it therefore 
suffices to show that $ F \Hom(M(l + \mu), V(l + \lambda)) $ is nonzero for some $ \mu $.

Choose $ \mu $ such that $ \nu = \mu - \lambda \in \weights^+ $, and let $V(-\nu)$ denote the simple module of lowest weight $-\nu$.
Then Proposition \ref{principalseriesmultiplicity} and Lemma \ref{locallyfiniteduallemma} show that
$$
[F\Hom(M(l + \mu), M(l + \lambda)^\vee): V(-\nu)] = \dim(V(-\nu)_{-\nu}) = 1 
$$
and 
$$
[F\Hom(M(l + \mu), M(l + \lambda - \gamma)^\vee): V(-\nu)] = \dim(V(-\nu)_{-\nu-\gamma}) = 0 
$$
for all $ \gamma \in \weights^+ \setminus \{0 \} $. Since $ V(l + \lambda) \subset M(l + \lambda)^\vee $ and the corresponding quotient has a filtration with 
subquotients isomorphic to $ V(l + \lambda - \gamma) $ for $ \gamma \in \weights^+ \setminus \{0 \} $ this yields the claim. \end{proof}  

Let us now summarize the results obtained so far. 

\begin{theorem} \label{categoryequivalence}
Let $ l \in \mathfrak{h}^*_q $ be dominant. 
\begin{bnum}
\item[a)] The functor $ \T_l: \HC_l \rightarrow \O_l $ embeds $ \HC_l $ as a full subcategory into $ \O_l $. 

\item[b)] If $ l \in \mathfrak{h}^*_q $ is regular then $ \F_l: \O_l \rightarrow \HC_l $ is an equivalence of $ \mathbb{C} $-linear categories. 
\end{bnum}
\end{theorem} 

\begin{proof}  $ a) $ It follows from Propositions \ref{fltladjointduality} and \ref{unitisomorphism}
that we have 
$$
\Hom_{\HC}(H, K) \cong \Hom_{\HC}(H, \F_l \T_l(K)) \cong \Hom_{U_q(\mathfrak{g})}(\T_l(H), \T_l(K))
$$
for all $ H, K \in \HC_l $, hence $ \T_l $ is fully faithful. 
This means precisely that $ \T_l $ embeds $ \HC_l $ as a full subcategory into $ \O_l $. 

$ b) $ 
Assume first that $ H \in \HC_l $ is simple. If $ \T_l(H) $ is non-simple there is a composition 
series $ 0 \subset V_0 \subset \cdots \subset V_n = \T_l(H) $ in $ \O_l $ with simple subquotients for some $ n > 0 $. 
Since $ \F_l $ is exact by Lemma \ref{Flexact} and maps simple modules to simple 
modules by Proposition \ref{domregsimpleobjects}, 
this induces a composition series $ 0 \subset \F_l(V_0) \subset \cdots \subset \F_l(V_n) = \F_l \T_l(H) \cong H $ with simple subquotients. 
Now our assumption that $ H $ is simple implies $ n = 0 $, which is a contradiction. Hence $ \T_l(H) $ is simple. 

We conclude that $ \T_l $ maps simple objects in $ \HC_l $ to simple objects in $ \O_l $. 
Taking into account Proposition \ref{domregsimpleobjects} it follows that $ \F_l $ and $ \T_l $ induce inverse equivalences on simple objects. 

Let $ 0 \rightarrow K \rightarrow E \rightarrow Q \rightarrow 0 $ be an extension in $ \O_l $, 
and assume that the counit of the adjunction from Proposition \ref{fltladjointduality} induces isomorphisms
$ \T_l \F_l(K) \cong K $ and $ \T_l \F_l(Q) \cong Q $. Consider the commutative diagram 
$$
\xymatrix{
0 \ar@{->}[r] & \T_l \F_l(K) \ar@{->}[r] \ar@{->}[d]^\cong & \T_l \F_l(E) \ar@{->}[d] \ar@{->}[r] & \T_l \F_l(Q) \ar@{->}[d]^\cong \ar@{->}[r] & 0 \\
0 \ar@{->}[r] & K \ar@{->}[r] & E \ar@{->}[r] & Q \ar@{->}[r] & 0
     }
$$
The bottom row is exact by assumption, so the map $ K \cong \T_l \F_l(K) \rightarrow \T_l \F_l(E) $ must be injective. 
Since $ \F_l $ is exact and $ \T_l $ is right exact, it follows that the upper row of the diagram is exact as well. 
This implies $ \T_l \F_l(E) \cong E $ as desired. 

The full subcategory of $ \O_l $ on which the counit $ e: \T_l \F_l \rightarrow \id $ of the adjunction is an isomorphism contains 
all simple objects, and every object in $ \O_l $ has finite length. We conclude $ \T_l \F_l \cong \id $, and this finishes the proof. \end{proof}

\subsection{Irreducible Harish-Chandra modules} 
\label{sec:Irreducible_HC-modules}

We shall now derive consequences for the category of Harish-Chandra modules. Some of these results have been sketched in \cite{Aranocomparison}.

Recall from the comments following Lemma \ref{lem:principal_series_multiplicities} that the principal series representation $\Gamma(\E_{\mu,\lambda})$ has a minimal $K_q$-type
which occurs with multiplicity one.

\begin{definition}
	\label{Def:V_mu_lambda}
	Let $(\mu,\lambda)\in\weights\times\lie{h}_q^*$.  We denote by $V_{\mu,\lambda}$ the unique irreducible subquotient of $\Gamma(\E_{\mu,\lambda})$ with the same minimal $K_q$-type.
	\nomenclature{$V_{\mu,\lambda}$}{irreducible subquotient of $\Gamma(\E_{\mu,\lambda})$ with the same minimal $K_q$-type}%
\end{definition}

To observe that $V_{\mu,\lambda}$ is indeed well-defined, consider the submodule of $\Gamma(\E_{\mu,\lambda})$ generated by any nonzero vector in the minimal $ K_q $-type. This module has a unique maximal submodule, namely the sum of all submodules not containing the minimal $K_q$-type, with corresponding quotient $V_{\mu,\lambda}$. 

Throughout this section, 
given $(\mu,\lambda) \in \weights \times \mathfrak{h}^*_q $, we shall use $(l,r)$ to denote a pair of weights
in $\mathfrak{h}^*_q \times \mathfrak{h}^*_q $ such that
$$ 
 \mu = l-r , \qquad \lambda+2\rho = -l -r,
$$
\label{nom:lr_shifted2}%
compare Corollary \ref{principalseriesversusFHom}.
We recall that the pair $(l,r)$ is well-defined only up to addition by an element of the form $(\half i\hbar^{-1} \alpha^\vee , \half i\hbar^{-1} \alpha^\vee)$ for some $\alpha^\vee \in \roots^\vee$.

We record the following equivalent formulations of a common integrality condition on the parameters $(\mu,\lambda)$. 

\begin{lemma}
	\label{lem:integer_point}
	Let $(\mu,\lambda)\in\weights\times\lie{h}_q^*$ and let $l,r\in\lie{h}_q^*$ be such that $ 
 \mu = l-r$ and $\qquad \lambda+2\rho = -l -r$.  For any positive root $\alpha\in{\bf\Delta}^+$, the following are equivalent:
	\begin{bnum}
		\item[$a)$]
		 $q_\alpha^{(\lambda,\alpha^\vee)} \in q_\alpha^{-|(\mu,\alpha^\vee)|-2\mathbb{N}}$.
		\item[$b)$]
		 Both $q_\alpha^{(l+\rho,\alpha^\vee)} \in \pm q_\alpha^{\mathbb{N}}$ and $q_\alpha^{(r+\rho,\alpha^\vee)} \in \pm q_\alpha^{\mathbb{N}}$.
	\end{bnum}
  Likewise, the following are equivalent:
	\begin{bnum}
	\item[$c)$]
	$q_\alpha^{(\lambda,\alpha^\vee)} \in q_\alpha^{|(\mu,\alpha^\vee)|+2\mathbb{N}}$.
	\item[$d)$]
	Both $q_\alpha^{(l+\rho,\alpha^\vee)} \in \pm q_\alpha^{-\mathbb{N}}$ and $q_\alpha^{(r+\rho,\alpha^\vee)} \in \pm q_\alpha^{-\mathbb{N}}$.
\end{bnum}  
\end{lemma}

\begin{proof}
	Let us put $\hbar_\alpha = \hbar d_\alpha$.  Then	condition $a)$ is equivalent to 
	\[
	 (\lambda,\alpha^\vee) \in -|(\mu,\alpha^\vee)|-2\mathbb{N} 
	  \bmod  i \hbar_\alpha^{-1}\mathbb{Z}.
	\]
	If we define $l'=l+\rho$, $r'=r+\rho$, so that we have
	\[
	\mu = l'-r', \quad \lambda=-l'-r',
	\]
	then the above condition is equivalent to having both
	\begin{align*}
	 -(l' + r',\alpha^\vee)
	  \in -(l' - r',\alpha^\vee) - 2\mathbb{N}
	  	  \bmod  i \hbar_\alpha^{-1}\mathbb{Z}, 
	  	  \\
	 -(l' + r',\alpha^\vee)
    \in -(r'-l',\alpha^\vee) - 2\mathbb{N}
    \bmod  i \hbar_\alpha^{-1}\mathbb{Z}.
	\end{align*}
	These conditions simplify to
	\begin{align*}
	 (l',\alpha^\vee) \in \mathbb{N}
	  \bmod \textstyle\half i \hbar_\alpha^{-1}\mathbb{Z},
	 \\
 	 (r',\alpha^\vee) \in \mathbb{N}
	 \bmod \textstyle\half i \hbar_\alpha^{-1}\mathbb{Z},
	\end{align*}
	which is equivalent to $b)$.  
	
	The equivalence of $c)$ and $d)$ follows by replacing $(\mu,\lambda)$ by $(-\mu,-\lambda)$.
\end{proof}

We now begin the study of the simple modules $V_{\mu,\lambda}$.

\begin{lemma} \label{lem:generic_principa_series}
		If there is no $ \alpha \in {\bf \Delta}^+  $ such that $ q_\alpha^{(l + \rho, \alpha^\vee)} \in \pm q_\alpha^{\mathbb{Z}} $ 
		then $\Gamma(\E_{\mu,\lambda})$ is simple, and we have
		\[
		V_{\mu,\lambda} \cong \F_l(V(r)) = \F_l(M(r)^\vee) \cong \Gamma(\E_{\mu,\lambda}).
		\]
\end{lemma}

\begin{proof}  
	Since $r-l$ is a weight, the hypothesis implies $ q_\alpha^{(r + \rho, \alpha^\vee)} \notin \pm q_\alpha^{\mathbb{Z}} $ for any $\alpha\in{\bf\Delta^+}$.  In particular, $r$ is antidominant, so $M(r)$ is simple by Theorem \ref{Vermacomplexirreducible}.  We get
	$ \F_l(V(r)) = \F_l(M(r)^\vee) \cong \Gamma(\E_{\mu, \lambda}),$
	which is nonzero, and hence irreducible by Lemma \ref{simplicitylemma}. By definition, we have $ V_{\mu, \lambda} = \Gamma(\E_{\mu, \lambda}) $ in this situation. 
\end{proof}

Our next result gives a sufficient condition for the simple modules $V_{\mu,\lambda}$ and $V_{\mu',\lambda'}$ to be isomorphic.  We will shortly show that this condition is also necessary, see Theorem \ref{HCequivalences}.

\begin{theorem}
	\label{thm:V_isomorphisms_sufficient}
	Let $(\mu,\lambda) \in \weights \times \lie{h}_q^*$.  Then $V_{\mu,\lambda} \cong V_{w\mu,w\lambda}$ for all $w\in W$.
\end{theorem}

\begin{proof}
	Fix $w\in W$ and let us write $(\mu',\lambda')=(w\mu,w\lambda)$.
	We begin by proving that $V_{\mu',\lambda'} \cong V_{\lambda,\mu}$ in the case where $\mu$ is fixed and $\lambda$ belongs to the dense subset
	\begin{align*}
	Z & = \{ \lambda\in\lie{h}_q^* \mid
	q_\alpha^{(\lambda-\mu,\alpha^\vee)} \notin  q_\alpha^{2\mathbb{Z}}
	\text{ for all }\alpha \in {\bf\Delta}^+ \} \subset \lie{h}_q^*.	
	\end{align*}
	Using the the relation $2l = \mu-\lambda-2\rho$ we get $q_\alpha^{(l,\alpha^\vee)} \notin  \pm q_\alpha^\mathbb{Z}$ for $\lambda\in Z$, and hence $V_{\mu,\lambda}=\Gamma(\E_{\mu,\lambda})$ by Lemma \ref{lem:generic_principa_series}.   
	The condition $q_\alpha^{(l,\alpha^\vee)} \notin \pm q_\alpha^\mathbb{Z}$ implies that ${\bf \Delta}_{[l]} = \emptyset$, in the notation of Subsection \ref{sec:dominant_and_antidominant}, and thus also ${\bf \Delta}_{[r]} = \emptyset$ since $r$ and $l$ differ by an integral weight.  Note also that $l$ and $r$ are both dominant and antidominant in this case.

	The parameters associated to $(\mu',\lambda')$ can be chosen as $l'=w.l$, $r'=w.r$, which are again both dominant and antidominant.   By Theorem \ref{chilinkage} we have $\xi_{l'} = \xi_l$, and so $\HC_l = \HC_{l'}$.  Therefore, Lemma \ref{simplicitylemma} implies that $\Gamma(\E_{\mu',\lambda'})$ is isomorphic to $\F_{l}(V)$ for some simple module $V\in \O_l$.  Specifically, we have $V\cong V(s)$ for some $s\in l+\weights$.  From the above condition on $l$ we see that $s$ is again dominant and antidominant.  In particular, according to Theorem \ref{Vermacomplexirreducible}, we have $V \cong M(s)^\vee$ and thus $ \F_{l'}(M(r')^\vee) \cong \F_l(M(s)^\vee) $.  
	
	The minimal $K_q$-type of $\F_{l'}(M(r')^\vee)$ is conjugate to $-\mu'=r'-l'$ under the Weyl group action and hence to $-\mu = r-l$.  Likewise, the minimal $K_q$-type of $\F_l(M(s)^\vee) $ is conjugate to $s-l$.  It follows that $r-l \in \weights$ and $s-l\in\weights$ are in the same orbit of the Weyl group action, and hence $s-r \in \roots$.  Moreover, since $\xi_s = \xi_{r'}$, we see that $s$ and $r'$ must be $\hat{W}$-linked, and hence $s$ and $r$ are $\hat{W}$-linked as well.  This implies that $s=\hat{v}.r$ for some $\hat{v}\in\hat{W}_{[r]}$, in the notation of Subsection \ref{sec:dominant_and_antidominant}.
	But we observed above that ${\bf\Delta}_{[r]}=\emptyset$, so $\hat{W}_{[r]}$ is trivial, according to Proposition \ref{prop:root_subsystem}.  
	Therefore $s=r$, and we conclude that $V_{\mu,\lambda} = \Gamma(\E_{\mu,\lambda}) \cong \Gamma(\E_{w\mu,w\lambda}) = V_{w\mu,w\lambda}$ for all $\lambda\in Z$.
	
	Now we extend this to all $\lambda\in\lie{h}_q^*$ as follows.  Recall from Subsection \ref{sec:principal_series_def} that the principal series modules $\Gamma(\E_{\mu,\lambda})$ for fixed $\mu$ all have the same underlying $\DF(K_q)$-module $\Gamma(\E_\mu)$.
	Let us write $\sigma$ for the minimal $K_q$-type in $\Gamma(\E_\mu)$ and denote the subspaces of $K_q$-type $\sigma$ by $E\subset\Gamma(\E_\mu)$ and $E'\subset\Gamma(\E_{w\mu})$, respectively.  Since the minimal $K_q$-type occurs with multiplicity one, there is a unique $K_q$-invariant isomorphism $\varphi:E \to E'$ up to scalar.

	Recall from Section \ref{sec:action_on_K_q_types} that under any of the representations $\pi_{\mu,\lambda}$, the subalgebra $U_q^\mathbb{R}(\lie{g})^{\sigma,\sigma} \subset U_q^\mathbb{R}(\lie{g})$  preserves $E$.    By the above discussion, when $\lambda\in Z$, there is a bijective intertwiner $\Gamma(\E_{\mu,\lambda}) \to \Gamma(\E_{w\mu,w\lambda})$, and after rescaling we may assume that it agrees with $\varphi$ on $E$.  In particular, for any $X\in U_q^\mathbb{R}(\lie{g})^{\sigma,\sigma}$ we have $ \varphi \circ \pi_{\mu,\lambda}(X)|_E =  \pi_{w\mu,w\lambda}(X)|_{E'} \circ \varphi$.  
	From the definition of the principal series representations $\pi_{\mu,\lambda}$ in Subsection \ref{sec:principal_series_def}, it is easy to check that the functions
	\begin{align*}
	\lie{h}_q^*  & \to \End(E) ;  \quad\;  \lambda  \mapsto \pi_{\mu,\lambda}(X), \\
	\lie{h}_q^*  & \to \End(E') ; \quad  \lambda  \mapsto \pi_{w\mu,w\lambda}(X) 
	\end{align*}
	are algebraic.  
	Since $Z$ is dense in $\lie{h}_q^*$ we 
	deduce that $\varphi\circ \pi_{\mu,\lambda}(X)|_E = \pi_{w\mu,w\lambda}(X)|_{E'}\circ \varphi$ for all $\lambda\in\lie{h}_q^*$.  Therefore, $E$ and $E'$ are isomorphic as $U_q^\mathbb{R}(\lie{g})^{\sigma,\sigma}$-modules for all $\lambda\in\lie{h}_q^*$, and so Proposition \ref{prop:simple_modules_vs_simple_K-types} shows that $V_{\mu,\lambda} \cong V_{w\mu,w\lambda}$ as $U_q^\mathbb{R}(\lie{g})$-modules.   
\end{proof}

To obtain more precise information on the structure of $\Gamma(\E_{\mu,\lambda})$, we will make
use of the following technical lemma.  We use the notation $\|\mu\|^2 = (\mu,\mu)$, when $\mu\in\weights$.  
Note that when $\mu,\mu'\in\weights^+$ we have $\mu \leq \mu'$ implies $\|\mu\| \leq \|\mu'\|$.

\begin{lemma} \label{dominantweightlemma}
	Let $\alpha\in{\bf\Delta}^+$ and assume that $ l,r \in \mathfrak{h}^* $ 
  are such that  $l-r\in\weights$ and that $ (l + \rho, \alpha^\vee) \in \mathbb{N}_0 $ and $(r + \rho, \alpha^\vee) \in \mathbb{N}_0 $. Then  
	$$
	\|l - r\| \leq \|l - s_\alpha . r\|.  
	$$
	Moreover equality holds iff the shifted action of $ s_\alpha $ stabilizes $ l $ or $ r $. 
\end{lemma} 
\begin{proof}  Note first that $ s_\alpha . r = r - (r + \rho, \alpha^\vee) \alpha $ and hence 
	$ l - s_\alpha . r = l - r + (r + \rho, \alpha^\vee) \alpha $.  
	We compute 
	\begin{align*}
	\|l - s_\alpha . r\|^2 - \|l - r\|^2 &= (l - s_\alpha . r, l - s_\alpha . r) - (l - r, l - r) \\
	&= 2 (l - r, (r + \rho, \alpha^\vee) \alpha) + (r + \rho, \alpha^\vee)^2(\alpha, \alpha) \\
	&= 2 (l - r, (r + \rho, \alpha^\vee) \alpha) + 2 (r + \rho, \alpha^\vee)(r + \rho, \alpha) \\
	&= 2 (r + \rho + l - r, (r + \rho, \alpha^\vee) \alpha) \\ 
	&= 2 (l + \rho, \alpha)(r + \rho, \alpha^\vee).  
	\end{align*}
	This yields the claim. \end{proof}

We now consider the modules $V_{\mu,\lambda}$ when the associated parameter $l$ is dominant.

\begin{lemma} \label{dominantvermahcproperties}
	Let $(\mu,\lambda) \in \weights \times \lie{h}_q^*$ such that the associated parameter $ l \in \mathfrak{h}^*_q $ is dominant. 
\begin{bnum} 
\item[a)] If there is no $ \alpha \in {\bf \Delta}^+ $ such that both $ q_\alpha^{(l + \rho, \alpha^\vee)} = \pm 1 $ 
and $ q_\alpha^{(r + \rho, \alpha)} \in \pm q_\alpha^{\mathbb{N}} $ 
then 
\[ 
 V_{\mu, \lambda} \cong \F_l(V(r)) \subset \F_l(M(r)^\vee) \cong \Gamma(\E_{\mu,\lambda}).
\]
Otherwise $ \F_l(V(r)) = 0 $.
\item[b)] 
The simple module $V_{\mu,\lambda}$ is a submodule of $\Gamma(\E_{\mu,\lambda})$.  Explicitly, there exists $\hat{u} \in \hat{W}$ with $\hat{u}.l=l$ such that the pair $(l, \hat{u}.r)$ verifies the conditions in $a)$ above, and we have
\[
  V_{\mu,\lambda} \cong \F_l(V(\hat{u}.r)) \subset \F_l(M(\hat{u}.r)^\vee) \cong \F_l(M(r)^\vee) \cong \Gamma(\E_{\mu,\lambda}).
\]
\end{bnum} 
\end{lemma}

\begin{proof}  
$ a) $ 
If $l$ fulfils the hypothesis of Lemma \ref{lem:generic_principa_series} we are done, so we assume henceforth that this is not the case.  Given that $l$ is assumed dominant, this means that $ q_\alpha^{(l + \rho, \alpha^\vee)} \in \pm q_\alpha^{\mathbb{N}_0} $ for some $ \alpha  \in {\bf \Delta}^+$.

Let $I(r)$ denote the kernel of the surjection $M(r) \to V(r)$.
By Lemma \ref{lemascendinghighestweightmodules}, we have a filtration
\[
  0 =M_0 \subset M_1 \subset M_2 \subset \cdots \subset M_{n} =I(r),
\]
where for each $1\leq i \leq n$ the quotient $M_i/M_{i-1}$ is a highest weight module.  By the BGG Theorem \ref{BGGtheorem}, if $t\in\lie{h}_q^*$ is a highest weight for any $M_i/M_{i-1}$ with $1\leq i\leq n$ then $t$ is strongly linked to $r$, meaning that $t= s_{k_1, \alpha_1} \cdots s_{k_m, \alpha_m} . r$ for some nonempty chain of affine reflections $s_{\alpha_j,k_j}$ with $ \alpha_j \in {\bf \Delta}^+ $ and $ k_j \in \mathbb{Z}_2 $, satisfying
\[
r > s_{k_m, \alpha_m} . r > s_{k_{m - 1}, \alpha_{m - 1}} s_{k_m, \alpha_{m}} . r > \cdots > s_{k_1, \alpha_1} \cdots s_{k_m, \alpha_m} . r = t.
\]
Note that this implies 
\[
 q_{\alpha_j}^{(s_{\alpha_{j+1}} \cdots s_{\alpha_m} . r + \rho, \alpha_j^\vee)} \in \pm q_{\alpha_j}^{\mathbb{N}} 
\]
for $1\leq j \leq m$.  
Moreover, since $r - s_{k_{j+1},\alpha_{j+1}} \cdots s_{k_m,\alpha_m} . r \in \roots$ and $l-r \in \weights$ we have $q_{\alpha_{j}}^{(l+\rho,\alpha_j^\vee)} \in \pm q_{\alpha_{j}}^{\mathbb{Z}}$ for each $\alpha_j$ appearing in this chain, and hence $q_{\alpha_{j}}^{(l+\rho,\alpha_j^\vee)} \in \pm q_{\alpha_{j}}^{\mathbb{N}_0}$  by dominance.  

Therefore, for each $1\leq i\leq n$ we have a surjective morphism
$ M(t_i) \to M_i/M_{i-1} $
for some $t_i$ strongly linked to $r$ as above, and therefore an injective 
morphism
\[
  \F_l((M_i/M_{i+1})^\vee) \to \F_l(M(t_i)^\vee).
\]
If there is no $ \alpha \in {\bf \Delta}^+ $ such that $ q_\alpha^{(l + \rho, \alpha^\vee)} = \pm 1 $ 
and $ q_\alpha^{(r + \rho, \alpha^\vee)} \in \pm q_\alpha^{\mathbb{N}} $, then an inductive application 
of Lemma \ref{dominantweightlemma} along the chain of reflections $s_{\alpha_j}$ shows that $ \|\Re(l) - \Re(t_i)\| > \|\Re(l) - \Re(r)\| $, where the real part of an element in $\lie{h}_q^*$ refers to its component in the $\mathbb{R}$-span of $\bf\Delta$, see the remarks before Proposition \ref{prop:root_subsystem}. Note that we have a strict inequality here because the strict inequalities in the strongly linking chain imply that the reflection $s_{\alpha_j}$ does not fix $s_{\alpha_{j+1}} \cdots s_{\alpha_m} . \Re(r) $.
According to Proposition \ref{principalseriesmultiplicity} a) and Proposition \ref{principalseriesversusverma}, the minimal $ K_q $-type in $ \F_l(M(t_i)^\vee) $ is strictly larger than that of $\F_l((M(r)^\vee)$, and therefore $\F_l((M_i/M_{i-1})^\vee)$ does not contain the minimal $K_q$-type of $\F_l(M(r)^\vee)$.  By induction, using the short exact sequences
\[
 0 \to \F_l((M_i/ M_{i-1})^\vee) \to \F_l(M_i^\vee) \to \F_l(M_{i-1}^\vee) \to 0,
\]
we deduce that $\F_l(M_i^\vee)$ does not contain the minimal $K_q$-type of $\F_l(M(r)^\vee)$ for any $1\leq i\leq n$.  

Therefore, in the exact sequence
\[
  0 \to \F_l(V(r)) \to \F_l(M(r)^\vee) \to \F_l(I(r)^\vee) \to 0,
\]
the inclusion on the left is necessarily an isomorphism on the minimal $K_q$-type.  It follows that $ \F_l(V(r)) $ is nonzero, and hence it is an irreducible Harish-Chandra module by Lemma \ref{simplicitylemma}. 
Moreover, since $ \F_l(V(r)) $ contains the minimal $ K_q $-type  of $\F_l(M(r)^\vee)$ we have $ \F_l(V(r)) \cong V_{\mu, \lambda} $ in this case. 

Conversely, assume that there is an $ \alpha \in {\bf \Delta}^+ $ satisfying both $ q_\alpha^{(l + \rho, \alpha^\vee)} = \pm 1 $ 
and $ q_\alpha^{(r + \rho, \alpha^\vee)} \in \pm q_\alpha^{\mathbb{N}} $. 
Then $ M(s_{k, \alpha} . r) $ is a submodule of $ M(r) $ for suitable $ k $, by Theorem \ref{Vermatheorem}, and hence $ V(r) $ is a quotient of $ M(r)/M(s_{k, \alpha} . r) $. 
Equivalently, $ V(r) $ is contained in the kernel of the natural projection $ M(r)^\vee \rightarrow M(s_{k, \alpha} . r)^\vee $. 
Applying the exact functor $ \F_l $ shows that $ \F_l(V(r)) $ is contained in the kernel of the surjective 
map $ \F_l(M(r)^\vee) \rightarrow \F_l(M(s_{k, \alpha} . r)^\vee) $. 
Note that $s_{k, \alpha}$ fixes $l$ up to a translation by an element of $\half i\hbar^{-1}\roots^\vee$, and since $s_{k, \alpha}. r$ and $r$ differ by an element of $\weights$ the same must be true of $s_{k, \alpha}. l$ and $l$, whence $ s_{k, \alpha} . l = l $.  It follows that both $\F_l(M(r)^\vee) \cong \Gamma(\E_{\mu,\lambda})$ and $\F_l(M(s_{k,\alpha}. r)^\vee) \cong \Gamma(\E_{s_\alpha\mu,s_\alpha\lambda})$ have the same $K_q$-type multiplicities, so that the surjection $ \F_l(M(r)^\vee) \rightarrow \F_l(M(s_{k, \alpha} . r)^\vee) $ is in fact an isomorphism.  Therefore $\F_l(V(r))$ is zero.

$ b) $
Under the conditions stated in $a)$ we obtain the assertion with $\hat{u}=1$, so it remains to consider the case where both $ q_\alpha^{(l + \rho, \alpha^\vee)} = \pm 1 $ 
and $ q_\alpha^{(r + \rho, \alpha^\vee)} \in \pm q_\alpha^{\mathbb{N}} $ for some $ \alpha \in {\bf \Delta}^+ $. 
Let $ \beta_1, \dots, \beta_{m} \in {\bf \Delta}^+ $ be the positive roots such that $ q_{\beta_i}^{(l + \rho, \beta_i^\vee)} = \pm 1 $. 
Consider the subgroup $ U $ of $ \hat{W} $ generated by the affine reflections corresponding to 
$ \beta_1, \dots, \beta_{m} $ and pick $ \hat{u}
\in U $ such that $ r' = \hat{u} . r $ 
is minimal. We claim that then $ q_{\beta_i}^{(r' + \rho, \beta_i^\vee)} \notin \pm q_{\beta_i}^{\mathbb{N}} $ for all $ 1\leq i \leq m $. 
Indeed, if $ q_{\beta_i}^{(r' + \rho, \beta_i^\vee)} \in \pm q_{\beta_i}^{\mathbb{N}} $ we 
would get 
$ s_{k, \beta_i} . r' < r' $ for suitable $ k $, contradicting our choice of $ r' $. 
Thus, by our definition of $ r' $ 
there is no $ \alpha \in {\bf \Delta}^+ $ such that $ q_\alpha^{(l + \rho, \alpha^\vee)} = \pm 1 $ 
and $ q_\alpha^{(r' + \rho, \alpha^\vee)} \in \pm q_\alpha^{\mathbb{N}} $. 

We are thus back in a situation as in $ a) $. 
Using Theorem \ref{thm:V_isomorphisms_sufficient}, we obtain an embedding
\[
 V_{\mu,\lambda} \cong V_{u\mu,u\lambda} \cong \F_l(\hat{u}.r) \subseteq \F_l(M(\hat{u}.r)^\vee),
\]
where $u$ denotes the image of $\hat{u}$ under the canonical projection $\hat{W} \to W$.
We may assume that $ \hat{u} $ is of the form $ \hat{u} = s_{k_1, \beta_{i_1}} \cdots s_{k_p, \beta_{i_p}} \in U $ 
such that 
$$ 
s_{k_1, \beta_{i_1}} \cdots s_{k_p, \beta_{i_p}} . r < s_{k_2, \beta_{i_2}} \cdots s_{k_p, \beta_{i_p}} . r < \cdots < s_{k_p, \beta_{i_p}} . r  < r ,
$$
and we have $s_{k_j,\beta_{i_j}} \ldots s_{k_p,\beta_{i_p}}.l=l$ for all $1\leq j \leq p$.
In this situation, as in $ a) $, we see that the corresponding quotient maps 
$$ 
\F_l(M(r)^\vee) \rightarrow \F_l(M(s_{k_p, \beta_{i_p}} . r)^\vee) \rightarrow \cdots \rightarrow \F_l(M(\hat{u}. r)^\vee) 
$$ 
are all isomorphisms.  Hence we get an inclusion $ V_{\mu, \lambda} \subset \F_l(M(\hat{u}.r)^\vee) \cong \F_l(M(r)^\vee) \cong \Gamma(\E_{\mu,\lambda}) $ as claimed.
\end{proof}  

\begin{lemma}
	\label{lem:FlVs_are_nonisomorphic}
	Assume $l\in\lie{h}_q^*$ is dominant, and let $r,r'\in l+\weights$ with $\F_l(V(r))\neq 0$.  Then $\F_l(V(r)) \cong \F_l(V(r'))$ iff $r=r'$.
\end{lemma}

\begin{proof}
	Assume $ r \neq r' $ and $ \F_l(V(r)) \cong \F_l(V(r')) $.  Let us write $H = \F_l(V(r))$.  Since $H\neq 0$, it is simple by Lemma \ref{simplicitylemma}.
	The adjointness relation from Proposition \ref{fltladjointduality} implies that there are nonzero $ U_q(\mathfrak{g}) $-linear maps $ \T_l(H) \rightarrow V(r), \T_l(H) \rightarrow V(r') $. 
	These are necessarily surjective by the simplicity of $ V(r), V(r') $. Consider the direct sum $ \T_l(H) \rightarrow V(r) \oplus V(r') $ of the maps 
	thus obtained. This map is again surjective because $ V(r) $ is not isomorphic to $ V(r') $. 
	By exactness of $ \F_l $ and Proposition \ref{unitisomorphism}, we obtain a surjection $ H \cong \F_l \T_l(H) \rightarrow \F_l(V(r)) \oplus \F_l(V(r')) $. 
	This contradicts the simplicity of $ H $. Hence $ \F_l(V(r)) \not\cong \F_l(V(r')) $. 
\end{proof}

\begin{theorem} \label{HCequivalences}
Two simple Harish-Chandra modules of the form $ V_{\mu, \lambda} $ and $ V_{\mu', \lambda'} $ 
for $ (\mu, \lambda), (\mu', \lambda') \in \weights \times \mathfrak{h}^*_q $ 
are isomorphic iff there exists $ w \in W $ 
such that 
$ (\mu', \lambda') = (w \mu, w \lambda) $. 
\end{theorem} 

\begin{proof}  
Sufficiency was proven in Theorem \ref{thm:V_isomorphisms_sufficient}, so it remains to check that the condition is necessary.  Suppose then that $V_{\mu,\lambda} \cong V_{\mu',\lambda'}$, and let $(l,r)$ and $(l',r')$ be the associated parameters in $\lie{h}_q^*\times\lie{h}_q^*$.  By Theorem \ref{chilinkage}, $l$ and $l'$ are $\hat{W}$-linked.  Therefore, after applying the isomorphisms of Theorem \ref{thm:V_isomorphisms_sufficient}, we may assume that $l=l'$ and moreover that $l$ is dominant.

By Theorem \ref{dominantvermahcproperties} $b)$, we have
$V_{\mu,\lambda} \cong \F_l(V(\hat{u}.r))$ and $V_{\mu',\lambda'} \cong \F_l(V(\hat{v}.r'))$ for some $\hat{u},\hat{v} \in \hat{W}$ such that $\hat{u}.l = \hat{v}.l = l$.  By Lemma \ref{lem:FlVs_are_nonisomorphic}, this implies $\hat{u}.r = \hat{v}.r'$.  Therefore, $(l',r') = (l,\hat{v}^{-1}\hat{u}.r) = (\hat{v}^{-1}\hat{u}.l , \hat{v}^{-1}\hat{u}.r)$ and hence
$(\mu',\lambda') = (v^{-1}u\mu,v^{-1}u\lambda)$, where $u$ and $v$ are the images in $W$ of the canonical quotient map $\hat{W} \to W$.
This completes the proof.
\end{proof}

\begin{theorem} \label{subquotient}
Every simple Harish-Chandra module is of the form $ V_{\mu, \lambda} $ for some $ (\mu, \lambda) \in \weights \times \mathfrak{h}^*_q $.  
\end{theorem} 

\begin{proof}  Let $ V $ be a simple Harish-Chandra module. Then the elements of the algebra 
$ 1 \otimes ZU_q(\mathfrak{g}) \subset FU_q(\mathfrak{g}) \otimes FU_q(\mathfrak{g}) $ 
commute with the action of $ U_q^\mathbb{R}(\mathfrak{g}) $,  
in particular with the action of $ U_q^\mathbb{R}(\mathfrak{k}) $. Hence they preserve the minimal $ K_q $-type of $ V $.  
According to Theorem \ref{subquotienthelp} and Schur's Lemma they act by scalars. 
Hence, the action of $ 1 \otimes ZU_q(\mathfrak{g}) $ is determined by a central character $\xi_l: ZU_q(\mathfrak{g}) \rightarrow \mathbb{C} $ for some $ l \in \mathfrak{h}^*_q $, see Proposition \ref{uqgcentralcharacters}. 
Since the central characters are invariant under the shifted $\hat{W}$-action 
by Theorem \ref{chilinkage} we may assume without loss of generality that $ l $ is dominant. It follows that
$ V $ is contained in $ \HC_l $,
and according to Lemma \ref{simplicitylemma}
$ V $ is isomorphic to $ \F_l(V(r)) $ for some simple module $ V(r) \in \O_l $. 
Since $V \neq 0$, Lemma \ref{dominantvermahcproperties} $a)$ shows that $V\cong \F_l(V(r)) \cong V_{\mu,\lambda}$ where $\mu=l-r$ and $\lambda = -l-r-2\rho$.
\end{proof}

Theorem \ref{subquotient} implies in particular that every irreducible Harish-Chandra module has a unique minimal $K_q$-type, which occurs with multiplicity one.

\begin{theorem} \label{principalirreducible}
Let $ (\mu, \lambda) \in \weights \times \mathfrak{h}^*_q $. 
Then the principal series module $ \Gamma(\E_{\mu, \lambda}) $ is an 
irreducible Yetter-Drinfeld module iff $ q_\alpha^{(\lambda, \alpha^\vee)} \neq q_\alpha^{\pm (|(\mu, \alpha^\vee)| + 2{k})} $ 
for all $ {k} \in \mathbb{N} $ and all $ \alpha \in {\bf \Delta}^+ $. 
\end{theorem} 

\begin{proof}  
The assumption on $ (\mu, \lambda) $
is equivalent to saying that for every $ \alpha \in {\bf \Delta}^+ $ 
we do not have nonzero integers $ m, n \in \mathbb{Z} $ of the same sign 
such that $ q_\alpha^{(l + \rho, \alpha^\vee)} = \pm q_\alpha^m $ and $ q_\alpha^{(r + \rho, \alpha^\vee)} = \pm q_\alpha^n $. 

We begin with the case where $l$ is dominant.  Under the above condition, this means that for each $\alpha\in{\bf\Delta}^+$, either $q_\alpha^{(l+\rho,\alpha^\vee)} \notin \pm q_\alpha^\mathbb{Z}$, in which case $q_\alpha^{(r+\rho,\alpha^\vee)} \notin \pm q_\alpha^\mathbb{Z}$, or $q_\alpha^{(l+\rho,\alpha^\vee)} = \pm q_\alpha^m$ for some $m\in\NN$, in which case $q_\alpha^{(r+\rho,\alpha^\vee)} \notin \pm q_\alpha^\NN$.  In either case, $r$ is antidominant, so that $M(r)$ is irreducible by Theorem \ref{Vermacomplexirreducible} and therefore $\Gamma(\E_{\mu,\lambda}) \cong \F_l(M(r)^\vee)$ is irreducible by Lemma \ref{simplicitylemma}.

Conversely,  if the above condition on $l$ and $r$ does not hold, with $l$ dominant, then there is some $\alpha\in{\bf\Delta}^+$ such that $q_\alpha^{(l+\rho,\alpha^\vee)}$ and $q_\alpha^{(r+\rho,\alpha^\vee)}$ both belong to $\pm q_\alpha^\mathbb{N}$.  In this case, we have an embedding $M(s_{k,\alpha}.r) \subset M(r)$ for some $k\in\mathbb{Z}_2$ by Theorem \ref{Vermatheorem}, and hence a surjective morphism 
\[
 \Gamma(\E_{\mu,\lambda}) \cong \F_l(M(r)^\vee) \to \F_l(M(s_{k,\alpha}.r)^\vee) \cong \Gamma(\E_{\mu',\lambda'}),
\]
where $\mu' = l-s_{k,\alpha}.r$ and $\lambda = -l-s_{k,\alpha}.r-2\rho$.
Moreover, since $(\Re(l)+\rho,\alpha^\vee) \neq0$, Lemma \ref{dominantweightlemma} shows that $\|\Re(l)-\Re(r)\| < \|\Re(l)-s_\alpha. \Re(r)\|$, so that the minimal $K_q$-type of $\F_l(M(s_{k,\alpha}.r)^\vee)$ is strictly larger than that of $\F_l(M(r)^\vee)$.  It follows that $\Gamma(\E_{\mu,\lambda})$ is not irreducible.

Consider now an arbitrary weight $ (\mu, \lambda) $ and let $ \hat{w}  \in \hat{W} $ be such that $ \hat{w}. l $ is dominant.  
Due to Theorem \ref{HCequivalences} we know that $ V_{w \mu, w \lambda} $ has the same $ K_q $-type multiplicities 
as $ V_{\mu, \lambda} $, where $w\in W$ denotes the image of $\hat{w}$ under the canonical projection of $\hat{W}$ onto $W$.   By Lemma \ref{lem:principal_series_multiplicities}, we also know that the principal series modules 
$ \Gamma(\E_{\mu, \lambda}) $ and $ \Gamma(\E_{w \mu, w \lambda}) $ have the same $ K_q $-type multiplicities. 
Since the condition $ q_\alpha^{(\lambda, \alpha^\vee)} \neq q_\alpha^{\pm (|(\mu, \alpha^\vee)| + 2{k})} $ 
for all $ {k} \in \mathbb{N} $ and all $ \alpha \in {\bf \Delta}^+ $ is stable under the action of $ W $, the result follows.
\end{proof}

\subsection{The principal series for $SL_q(2,\mathbb{C})$}
\label{sec:sl2_principal_series}

In this section we will determine the structure of the principal series representations for $SL_q(2,\mathbb{C})$ as well as all intertwiners between them.  As a result we obtain a complete classification of the irreducible Harish-Chandra modules for $SL_q(2,\mathbb{C})$.  These results were first proven by Pusz-Woronowicz \cite{PWquantumlorentzgelfand}.  Here we obtain them as a consequence of the relationship between the categories $\HC_l$ and $\O_l$ described in the previous sections.

Recall that when $\lie{g} = \lie{sl}(2,\mathbb{C})$ we identify $\lie{h}^*$ with $\mathbb{C}$ via the correspondence which sends $\lambda\in\lie{h}^*$ to $\half (\lambda,\alpha^\vee)$.  With this convention, $\weights = \half\mathbb{Z}$,  $\roots^\vee = \roots = \mathbb{Z}$ and $\lie{h}_q^* = \mathbb{C} / i\hbar^{-1}\mathbb{Z}$.  The unique simple root is $\alpha=1$ and we have $\rho=\half$.  It is important to note that the invariant bilinear form on $\lie{h}^*$ is given by $(\lambda_1,\lambda_2) = 2\lambda_1\lambda_2$.  

As in the previous section, for parameters $(\mu,\lambda) \in \half\mathbb{Z} \times \mathbb{C} / i \hbar^{-1} \mathbb{Z}$, we let $(l,r) \in \mathbb{C} / i \hbar^{-1}\mathbb{Z} \times \mathbb{C} / i \hbar^{-1}\mathbb{Z}$ denote a pair such that 
\[
\mu = l-r, \qquad \lambda + 1 = -l-r.
\]
There are always two such choices for $(l,r)$ which differ by $(\half i \hbar^{-1}, \half i \hbar^{-1})$.  
A weight $l\in\mathbb{C}/i\hbar^{-1}\mathbb{Z}$ is dominant in the sense of Definition \ref{defdominantantidominant} iff
\[
\textstyle l + \half \notin -\half\mathbb{N} \bmod \frac14 i \hbar^{-1} \mathbb{Z},
\]
and antidominant iff 
\[
\textstyle l + \half \notin \half\mathbb{N} \bmod \frac14 i \hbar^{-1} \mathbb{Z},
\]
where as usual $\mathbb{N} = \{1,2,3,\ldots\}$.

\begin{theorem}
	\label{thm:SL2_irreducibilities}
	Let $(\mu,\lambda)\in \half\mathbb{Z}\times\CC/i\hbar^{-1}\mathbb{Z}$.  The principal series representation $\Gamma(\E_{\mu,\lambda})$ for $SL_q(2,\mathbb{C})$ is irreducible iff $\lambda\notin \pm(|\mu|+\mathbb{N}) \bmod \half i \hbar^{-1}\mathbb{Z}$.  Moreover, when $\lambda \in \pm(|\mu|+\mathbb{N}) \bmod \half i \hbar^{-1}\mathbb{Z}$ the simple module $V_{\mu,\lambda}$ is finite dimensional and $\Gamma(\E_{\mu,\lambda})$ is a non-split extension of simple modules as follows:
	\begin{align*}
	 &0 \to V_{\mu,\lambda} \to \Gamma(\E_{\mu,\lambda})	\to \Gamma(\E_{\lambda,\mu}) \to 0 
	 && \text{if } \lambda\in-|\mu|-\mathbb{N},
	 \\
	 &0 \to V_{\mu,\lambda} \to \Gamma(\E_{\mu,\lambda})	\to \Gamma(\E_{\lambda+\half i\hbar^{-1},\mu+\half i\hbar^{-1}}) \to 0 
	 && \text{if } \lambda\in-|\mu|-\mathbb{N}+\half i\hbar^{-1}\mathbb{Z},
	 \\
	 &0 \to \Gamma(\E_{\lambda,\mu}) \to \Gamma(\E_{\mu,\lambda}) 
	 \to V_{\mu,\lambda} \to 0
	 && \text{if } \lambda\in|\mu|+\mathbb{N},
	 \\
	 &0 \to \Gamma(\E_{\lambda+\half i\hbar^{-1},\mu+\half i\hbar^{-1}}) \to \Gamma(\E_{\mu,\lambda}) 
	 \to V_{\mu,\lambda} \to 0
	 && \text{if } \lambda\in|\mu|+\mathbb{N}+\half i\hbar^{-1}\mathbb{Z}.
	\end{align*} 
\end{theorem}

\begin{proof}
  The statement about irreducibility follows directly from Theorem \ref{principalirreducible}.  
  
  If $\lambda\in-|\mu|-\mathbb{N}$ then $l+\half, r+\half \in \half\mathbb{N}$, so that $l$ is dominant and regular and $r$ is not antidominant.  We have a short exact sequence in category $\O$,
  \[
   0 \to V(r) \to M(r)^\vee \to M(-r-1)^\vee \to 0.
  \]
  Note that $\F_l(M(-r-1)^\vee) \cong \Gamma(\E_{-\lambda,-\mu})$ is irreducible, and so is isomorphic to $\Gamma(\E_{\lambda,\mu})$ by Theorem \ref{HCequivalences}.  Therefore, applying the functor $\F_l$ to the above short exact sequence gives the first of the stated extensions.  It is non-split by Theorem \ref{categoryequivalence}.  To see that $V_{\mu,\lambda}$ is finite dimensional, it suffices to compare the $K_q$-type multiplicities in the two principal series representations in the extension.  Specifically, Lemma \ref{lem:principal_series_multiplicities} shows that $\Gamma(\E_{\mu,\lambda})$ contains the $K_q$-types with highest weights $\nu\in|\mu|+\mathbb{N}$, each with multiplicity one.
  
  The case $\lambda\in-|\mu|-\mathbb{N}+\half i \hbar^{-1}$ is obtained similarly, and the final two cases follow by duality using Lemma \ref{lem:principal_series_bilinear_pairing}.
\end{proof}

\begin{theorem} \label{thm:slq2_intertwiners}
	The following is an exhaustive list of non-zero intertwining operators between principal series representations of $ SL_q(2, \mathbb{C}) $.
	\begin{enumerate}
		\item[1.] {\bf Trivial intertwiners:} For every $ (\mu, \lambda) \in \weights \times \mathfrak{h}_q^* $, the only self-intertwiners 
		$ \Gamma(\E_{\mu, \lambda}) \rightarrow \Gamma(\E_{\mu, \lambda}) $ are the scalar multiples of the identity.
		\item[2.] {\bf Standard intertwiners:} For every $ (\mu, \lambda) \in \weights \times \mathfrak{h}_q^* $ there is a unique intertwiner, up to 
		scalars, $ \Gamma(\E_{\mu, \lambda}) \rightarrow \Gamma(\E_{-\mu, -\lambda}) $ as follows:
		\begin{bnum}
			\item[a)] If $ \lambda \in -|\mu| - \mathbb{N} \bmod \frac{1}{2} i \hbar^{-1}\mathbb{Z} $ the intertwiner is Fredholm with kernel and cokernel isomorphic to $V_{\mu,\lambda}$.
			\item[b)] If $ \lambda \in |\mu| + \mathbb{N} \bmod \frac{1}{2} i \hbar^{-1}\mathbb{Z} $ the intertwiner is finite-rank and factors through $V_{\mu,\lambda}$.
			\item[c)] Otherwise, the intertwiner is bijective.
		\end{bnum}
		
		\item[3.] \textbf{\u{Z}elobenko intertwiners:}
		With $ \lambda \in -|\mu| - \mathbb{N} $, there are additional intertwiners as follows:
		\[ \label{eq:diamond1}
		\xymatrix{
			& \Gamma(\E_{-\lambda, -\mu}) 
			\ar@{<->}[dd]^(0.25)\cong 
			\ar@{^{(}->}[dr]
			\\
			\Gamma(\E_{\mu,\lambda}) 
			\ar@<0.3ex>[rr]
			\ar@{->>}[ur] 
			\ar@{->>}[dr] 
			&& \Gamma(\E_{-\mu, -\lambda}) \ar@<0.3ex>[ll]
			\\
			& \Gamma(\E_{\lambda, \mu})  
			\ar@{^{(}->}[ur]
		}
		\]
		and
		\[ \label{eq:diamond2}
		\xymatrix{
			& \Gamma(\E_{-\lambda,-\mu + \frac{1}{2} i \hbar^{-1}}) \ar@{<->}[dd]^(0.25)\cong 
			\ar@{^{(}->}[dr]
			\\
			\Gamma(\E_{\mu,\lambda + \frac{1}{2} i \hbar^{-1}})
			\ar@<0.3ex>[rr]
			\ar@{->>}[ur]
			\ar@{->>}[dr]
			&& \Gamma(\E_{-\mu,-\lambda + \frac{1}{2} i \hbar^{-1}}) \ar@<0.3ex>[ll] 
			\\
			& \Gamma(\E_{\lambda,\mu + \frac{1}{2} i \hbar^{-1}})  \ar@{^{(}->}[ur]
		}
		\]
		where the right- and left-pointing horizontal arrows are the standard intertwiners from $2a)$ and $ b)$ respectively, the vertical arrows are the standard intertwiners from $2c)$, and the diagonal arrows are obtained using the extensions in Theorem \ref{thm:SL2_irreducibilities}. 
		These diagrams commute, up to scalar multiple, except for the finite rank intertwiners from right to left which have composition zero with any other. 
	\end{enumerate}
\end{theorem}

\begin{proof}
	By Theorem \ref{thm:SL2_irreducibilities}, $\Gamma(\E_{\mu,\lambda})$ is either simple or a non-split extension of two non-isomorphic simple modules.  In either case, $\Gamma(\E_{\mu,\lambda})$ does not admit any nontrivial self-intertwiners.
	
	Inspecting the extensions in Theorem \ref{thm:SL2_irreducibilities}, we see that the non-simple principal series representations are all mutually non-isomorphic.
	It follows that there exists a unique intertwiner, up to scalars, between two principal series modules $\Gamma(\E_{\mu_1,\lambda_1})$ and $\Gamma(\E_{\mu_2,\lambda_2})$ if and only if $\Gamma(\E_{\mu_1,\lambda_1})$ admits a simple quotient which is isomorphic to a simple submodule of $\Gamma(\E_{\mu_2,\lambda_2})$.
	The simple Harish-Chandra modules are all of the form $V_{\mu,\lambda}$ for some $(\mu,\lambda)\in\weights\times\lie{h}_q^*$, and we have $V_{\mu,\lambda} = \Gamma(\E_{\mu,\lambda})$ if $\lambda\notin\pm(|\mu|+\mathbb{N}) \bmod \half i\hbar^{-1} \mathbb{Z}$.  Also, by Theorem \ref{HCequivalences} we have $V_{\mu,\lambda} \cong V_{\mu',\lambda'}$ if and only if $(\mu',\lambda') = \pm(\mu,\lambda)$.  Therefore, to complete the proof, we can make a case-by-case examination using the structure of the principal series modules given in Theorem \ref{thm:SL2_irreducibilities}.
  	
  Here are the details.  First, consider the case where both $\lambda_1\notin \pm(|\mu_1|+\mathbb{N}) \bmod \half i \hbar^{-1}\mathbb{Z}$ and $\lambda_2\notin \pm(|\mu_2|+\mathbb{N}) \bmod \half i \hbar^{-1}\mathbb{Z}$.  Then $\Gamma(\E_{\mu_1,\lambda_1})$ and $\Gamma(\E_{\mu_2,\lambda_2})$ are both simple, so the only nontrivial intertwiners in this situation are the bijective intertwiners $\Gamma(\E_{\mu,\lambda}) \to \Gamma(\E_{-\mu,-\lambda})$ from case $2c)$.
  
  Next, we consider the cases where either $\lambda_1\in \pm(|\mu_1|+\mathbb{N})$ or $\lambda_2\in \pm(|\mu_2|+\mathbb{N})$ on the nose.  We will obtain precisely the intertwiners in the first diagram from part 3.
	\begin{itemize}
  \item
   	If $\lambda_1\in-|\mu_1|-\mathbb{N}$ then the image of $T$ must be isomorphic to $\Gamma(\E_{\lambda_1,\mu_1}) = V_{\lambda_1,\mu_1}$.
   	 The only principal series modules which contain $\Gamma(\E_{\lambda_1,\mu_1})$ as a submodule are $\Gamma(\E_{\lambda_1,\mu_1})$ itself, $\Gamma(\E_{-\lambda_1,-\mu_1})$ which is simple and isomorphic to $\Gamma(\E_{\lambda_1,\mu_1})$, and $\Gamma(\E_{-\mu_1,-\lambda_1})$ which contains $\Gamma(\E_{-\lambda_1,-\mu_1})$ as a proper submodule according to Theorem \ref{thm:SL2_irreducibilities}.   	 
   	 These correspond to the three intertwiners issuing from the left-hand module in the first diagram. 
   \item 
    If $\lambda_1\in|\mu_1|+\mathbb{N}$ then the image of $T$ must be isomorphic to $V_{\mu_1,\lambda_1}$ and the only principal series module containing this as a submodule is $\Gamma(\E_{-\mu_1,-\lambda_1})$.  This corresponds to the horizontal arrow from right to left.
   \item 
    If $\lambda_2\in-|\mu_2|-\mathbb{N}$ then the image of $T$ must be isomorphic to $V_{\mu_2,\lambda_2}$ and the only principal series admitting this as a quotient module is $\Gamma(\E_{-\mu_2,-\lambda_2})$. This corresponds again to the horizontal arrow from right to left.
   \item
    If $\lambda_2\in|\mu_2|+\mathbb{N}$ then the image of $T$ must be isomorphic to $\Gamma(\E_{\lambda_2,\mu_2})$ and the only principal series modules admitting this as a quotient module are $\Gamma(\E_{\lambda_2,\mu_2}) \cong \Gamma(\E_{-\lambda_2,-\mu_2})$ and the non-simple module $\Gamma(\E_{-\mu_2,-\lambda_2})$.  These correspond to the three intertwiners mapping to the right-hand module.  
	\end{itemize}
	
	Finally, the cases where $\lambda_1 \in \pm(|\mu_1| +\mathbb{N}) +\half i \hbar^{-1}$ or $\lambda_2 \in \pm(|\mu_2| +\mathbb{N}) +\half i \hbar^{-1}$ are treated similarly, and we obtain the intertwiners in the second diagram from part $3$.
\end{proof}

We note that the two diagrams in part $3$ of Theorem \ref{thm:slq2_intertwiners} are related by the action of the group $Z=\weights^\vee/\roots^\vee \cong \mathbb{Z}_2$ which is described in Subsection \ref{sec:action_of_centre}.  
More precisely, the representations in the second diagram become isomorphic to the corresponding representations in the first diagram upon tensoring by the one-dimensional module $\CC_{\half i \hbar^{-1}}$, see Lemma \ref{lem:Z-action_on_principal_series}.  
Therefore, the corresponding principal series representations in these diagrams become identical upon restriction to the ``connected component'' $SL_q(2,\mathbb{C})^0$.  Moreover, under the canonical identification of the vector spaces $\Gamma(\E_{\mu,\lambda}) = \Gamma(\E_{\mu,\lambda+\half i \hbar^{-1}})$ for each pair $(\mu,\lambda)$, the corresponding intertwiners in the two diagrams are given by exactly the same linear operators.

Finally, let us remark that the \u{Z}elobenko intertwiners can be interpreted geometrically.  For instance, from the case $(\mu,\lambda)=(0,1)$ we obtain a system of $\DF(G_q)$-linear operators
 \[ 
		\xymatrix@R=0mm{
			&& \Gamma(\E_{1,0}) 
			\ar[dr]
			\\
			\mathbb{C} \ar[r]&
			\Gamma(\E_{0,-1}) 
			\ar[ur]
			\ar[dr]
			&& \Gamma(\E_{0, 1}) \ar[r]
			& \mathbb{C}			
			\\
			&& \Gamma(\E_{-1, 0})  
			\ar[ur]
		}
\]
where $\mathbb{C}$ carries the trivial representation.
In Subsection \ref{sec:intertwiner_formulas}, we will show that the diagonal intertwiners in this diagram are given by the right regular action of the elements $E$ or $F$ in $U_q^\mathbb{R}(\lie{k})$.

The principal series representation $\Gamma(\E_{0,-1})$ 
is the quantum analogue of the space of polynomial functions on the flag variety $G/B$.  In this vein, the above diagram can be viewed as a quantum analogue of the $(\partial,\bar\partial)$-complex:
 \[ 
	\xymatrix@R=0mm{
		&& \Omega^{0,1}(G/B) 
		\ar@{.}[dd]|\bigoplus
		\ar[dr]^{\partial} 
		\\
		\mathbb{C} \quad \ar[r]&
		\Omega^{0,0}(G/B)  
		\ar[ur]^{\bar\partial}
		\ar[dr]_{\partial}
		&& \Omega^{1,1}(G/B) \ar[r]
		& \mathbb{C}	.		
		\\
		&& \Omega^{1,0}(G/B)  
		\ar[ur]_{\bar\partial}
		}
\]
Similarly, from the other values of $(\mu,\lambda)$ appearing in part $3$ of Theorem \ref{thm:slq2_intertwiners} we obtain quantum analogues of the $(\partial,\bar\partial)$-complex twisted by a $G$-equivariant vector bundle.

\subsection{Intertwining operators in higher rank}

\subsubsection{Intertwiners in the compact picture}
\label{sec:YD_intertwiners}

Fix $(\mu,\lambda) \in \weights \times \lie{h}_q^*$.
Let $ f = \bra v'| \bullet| v \ket \in \CF^\infty(K_q) $ and $ \xi = \bra w'| \bullet |w \ket \in \Gamma(\E_{\mu, \lambda}) $ be matrix coefficients of finite 
dimensional $ \DF(K_q)$-modules $ V $ and $ W $, respectively. That is, we have $ v \in V$, $v' \in V^*$, $w \in W$, $w'\in W^* $, and moreover $ w $ has 
weight $ \mu $. Let $ e_1, \dots, e_n $ be a weight basis for $ V $, with $ e_j $ of weight $ \epsilon_j $ for each $1\leq j\leq n$, and let $ e^1, \dots e^n $ be the dual basis of $ V^* $. 

We shall write $\precon V$ 
\nomenclature[s$V$]{$\precon V$}{pre-contragredient of a $U_q^\mathbb{R}(\lie{k})$-module $V$}%
for the dual space $V^*$ equipped with
the precontragredient representation 
of $ U_q^\mathbb{R}(\mathfrak{k}) $, namely
$$
(X \cdot {v'})(v) = {v'}(\hat{S}^{-1}(X) \cdot v) 
$$
for $ X \in U_q^\mathbb{R}(\mathfrak{k}), v' \in \precon V $ and $ v \in V $. 
This is relevant for the antipode in $\CF^\infty(K_q)$, since
\[
S(\bra v'| \bullet| v \ket) = \bra v| \bullet| v' \ket,
\]
where the latter is a matrix coefficient for the pre-contragredient representation $\precon V$.

\begin{lemma} \label{eq:matrix_coeff_action}
	With the above notation, the Yetter-Drinfeld action of $ f = \bra v'| \bullet| v \ket \in \CF^\infty(K_q) $ on $ \xi = \bra w'| \bullet |w \ket  \in \Gamma(\E_{\mu, \lambda}) $ becomes 
	\begin{align*}
	f \cdot \xi 
	&= \sum_j (K_{2\rho + \lambda}, \bra e^j| \bullet |e_j \ket) \bra v \otimes w' \otimes v'| \bullet |e^j \otimes w \otimes e_j \ket \\
	&= \sum_j q^{(\lambda + 2\rho, \epsilon_j)} \bra v \otimes w' \otimes v'|\bullet |e^j \otimes w \otimes e_j \ket,
	\end{align*}
	where the latter is viewed as a matrix coefficient for $ \precon V \otimes W \otimes V $.
\end{lemma} 

\begin{proof}  This follows immediately from the definitions.  Indeed, we compute 
	\begin{align*}
	(x, f \cdot \xi) &= (x, f_{(1)} \xi S(f_{(3)})) (K_{2 \rho + \lambda}, f_{(2)}) \\
	&= (x_{(3)}, f_{(1)}) (x_{(2)}, \xi) (x_{(1)}, S(f_{(3)})) (K_{2 \rho + \lambda}, f_{(2)}) \\
	&= \sum_{i,j} (K_{2\rho + \lambda}, \bra e^i| \bullet |e_j \ket) 
	(x_{(1)}, \bra v| \bullet |e^j \ket)(x_{(2)}, \bra w'|\bullet |w \ket) (x_{(3)}, \bra v'| \bullet |e_i \ket) \\
	&= \sum_j (x, (K_{2\rho + \lambda}, \bra e^j| \bullet |e_j \ket) \bra v \otimes w' \otimes v'| \bullet |e^j \otimes w \otimes e_j \ket) 
	\end{align*}
	for all $ x \in \DF(K_q) $, as desired. \end{proof}

Let $ T \in \M(\DF(K_q)) $ be an element of weight $ \beta \in \weights $ for the adjoint action, meaning 
$$ 
K_\lambda T K_{-\lambda} = q^{(\lambda, \beta)} T 
$$ 
for all $ \lambda \in \weights $. Then the right regular action of $ T $ on $ \CF^\infty(K_q) $,
\[
 T \hit \xi = (T, \xi_{(2)}) \xi_{(1)},
\]
restricts to a morphism of $ \DF(K_q) $-modules
$ \Gamma(\E_\mu) \rightarrow \Gamma(\E_{\mu + \beta}) $ for any $ \mu \in \weights $. 
We will denote this morphism simply by $ T $.

In general, any $ \DF(K_q) $-linear map $ f: \Gamma(\E_\mu) \rightarrow \Gamma(\E_{\mu + \beta}) $ is of this form for 
some (in fact, many) $ T \in \M(\DF(K_q)) $. To confirm this, use the Peter-Weyl decomposition of $ \CF^\infty(K_q) $ 
to write $ f $ in the form 
$$
f = \bigoplus_{\sigma \in \weights^+} \id \otimes T_\sigma: 
\bigoplus_{\sigma \in \weights^+} V(\sigma)^* \otimes V(\sigma)_\mu \rightarrow 
\bigoplus_{\sigma \in \weights^+} V(\sigma)^* \otimes V(\sigma)_{\mu + \beta}
$$
for some family of linear maps $ T_\sigma: V(\sigma)_\mu \rightarrow V(\sigma)_{\mu + \beta} $. 
Extend these maps to $ T_\sigma: V(\sigma) \rightarrow V(\sigma) $, for instance by zero on all other weight spaces, 
to obtain $ T = \bigoplus_\sigma T_\sigma \in \M(\DF(K_q)) $ which acts as $ f $.

Given this characterization of the $\DF(K_q)$-linear maps between principal series representations, Lemma \ref{eq:matrix_coeff_action} immediately yields the following characterization 
of intertwiners of principal series representations.

\begin{lemma} \label{lem:intertwiner_condition}
	Fix $(\mu_1,\lambda_1), (\mu_2,\lambda_2) \in \weights \times \lie{h}_q^*$ and let $ T \in \M(\DF(K_q)) $ have weight $ \mu_2 - \mu_1$ for the adjoint action. 
	Then the following conditions are equivalent. 
	\begin{bnum} 
		\item[a)] $ T: \Gamma(\E_{\mu_1, \lambda_1}) \rightarrow \Gamma(\E_{\mu_2, \lambda_2}) $ is an intertwiner of $ \DF(G_q) $-modules. 
		\item[b)] For any finite dimensional $ \DF(K_q) $-modules $ V $ and $ W $, 
		for any $ w \in W_{\mu_1} $  we have
		\begin{equation*} \label{eq:intertwiner_condition}
		\sum_j q^{(\lambda_1 + 2 \rho, \epsilon_j)} \; T \cdot (e^j \otimes  w \otimes e_j) = \sum_j q^{(\lambda_2 + 2\rho, \epsilon_j)} \; e^j \otimes T \cdot w \otimes e_j.
		\end{equation*}
		where $ e_1, \dots, e_m $ is a weight basis for $V$ with weights $ \epsilon_1,\ldots, \epsilon_m $, respectively, and $ e^1, \dots, e^m $ is the dual basis
			of $ \precon V $.
	\end{bnum}
\end{lemma}

\subsubsection{Intertwiners associated to simple roots}
\label{sec:simple_intertwiners}

Fix a simple root $\alpha$.  
Consider the Hopf $*$-algebra $U_{q_\alpha}^\mathbb{R}(\lie{su}(2))$ generated by elements $E$, $F$ and $K= K_\varpi$ as usual, where $\varpi$ is the fundamental weight of $\lie{sl}(2,\mathbb{C})$.  We write $U_{q_\alpha}^\mathbb{R}(\lie{psu}(2))$
\nomenclature{$U_{q_\alpha}^\mathbb{R}(\lie{psu}(2))$}{$*$-Hopf subalgebra of $U_{q_\alpha}^\mathbb{R}(\lie{su}(2))$ generated by $E$, $F$ and $K_\alpha=K^2$}%
for the $*$-Hopf subalgebra generated by $E$, $F$ and $K_\alpha=K^2$.

There is a unique Hopf $ * $-algebra morphism 
\[
\iota_\alpha: U_{q_\alpha}^\mathbb{R}(\mathfrak{psu}(2)) \rightarrow U_q^\mathbb{R}(\mathfrak{k})
\]
\label{nom:suq2-embedding1}%
sending $ E , F, K^2 $ to $ E_\alpha, F_\alpha, K_\alpha $, respectively. In this way, every $ U_q^\mathbb{R}(\mathfrak{k}) $-module $ V $ restricts to 
a $ U_{q_\alpha}^\mathbb{R}(\mathfrak{psu}(2)) $-module. If $ V $ is integrable then the restriction extends uniquely to an integrable 
module for $ U_{q_\alpha}^\mathbb{R}(\mathfrak{su}(2)) $, where $ K $ acts by the positive square root of $ K^2 $. This gives a well-defined restriction functor
from integrable $U_q^\mathbb{R}(\lie{k})$-modules to integrable $U_{q_\alpha}^\mathbb{R}(\lie{su}(2))$-modules and thus a map
\begin{align*}
\iota_\alpha: \DF(SU_{q_\alpha}(2)) \rightarrow \M(\DF(K_q)).
\end{align*}
\label{nom:suq2-embedding2}%
We will write $S_\alpha$ 
\nomenclature{$S_\alpha$}{rank one quantum subgroup of $K_q$ associated to a simple root $\alpha$}%
for the quantum subgroup of $K_q$ which is obtained in this way.

Let $ V $ be an integrable $ U_q^\RR(\mathfrak{k}) $-module. If $ v \in V $ is a weight vector of weight $ \mu \in \weights $, then 
$$
\iota_\alpha(K) \cdot v = q^{\half(\alpha, \mu)} v = q_\alpha^{\half(\alpha^\vee, \mu)} v.
$$
In other words, upon restricting $V$ to a $\DF(S_\alpha)$-module, $ v $ has weight $ \frac{1}{2} (\alpha^\vee, \mu) \in\half\mathbb{Z}$.  
Motivated by this, we define
$$
\textstyle \mu_\alpha = \frac{1}{2} (\alpha^\vee,\mu),
$$
\nomenclature{$\mu_\alpha$}{restriction of $\mu\in\weights^*$ with respect to a simple root $\alpha$}%
and call it the restriction of $ \mu $ with respect to $ \alpha $. Note in particular that $ \rho_\alpha = \frac{1}{2} $.

More generally, we write $\lambda_\alpha = \half(\lambda,\alpha^\vee) \in \CC$ for any $\lambda \in \lie{h}^*$.  %
\label{nom:weight_restriction1}%
This yields the orthogonal decomposition
\[
\lambda = \lambda' + \lambda_\alpha\alpha,
\]
where $\lambda' = \lambda - \lambda_\alpha\alpha \in \alpha^\perp$.

We also wish to define $\lambda_\alpha$ when $\lambda \in \lie{h}_q^* = \lie{h}^*/i\hbar^{-1}\roots^\vee$.  We have
\[ \textstyle
\half (i\hbar^{-1} \roots^\vee, \alpha^\vee) = \half i\hbar^{-1} d_\alpha^{-1} (\roots^\vee,\alpha) \subseteq \half i \hbar_\alpha^{-1}\mathbb{Z},
\]
where we are using the notation $\hbar_\alpha = d_\alpha\hbar$.
\nomenclature[o$h_\alpha$]{$\hbar_\alpha$}{$= d_\alpha\hbar$ for a simple root $\alpha$}%
Therefore, we obtain a map
\[
\lie{h}_q^* \to \CC / {\textstyle \half} i\hbar_\alpha^{-1} \mathbb{Z}; \qquad  \lambda \mapsto \lambda_\alpha.
\]
\label{nom:weight_restriction2}%

We point out that that the factor of $\half$ in the above quotient means that the restriction of parameters
\begin{align*}
\weights \times \mathfrak{h}_q^* &\rightarrow \textstyle\frac{1}{2}\mathbb{Z} \times \mathbb{C}/\frac{1}{2} i \hbar_\alpha^{-1} \mathbb{Z} \\
(\mu, \lambda) & \mapsto (\mu_\alpha, \lambda_\alpha)
\end{align*}
sends principal series parameters for $ G_q $ to principal series parameters for the complex quantum group $ SL_{q_\alpha}(2, \mathbb{C})^0 $ rather than $SL_{q_\alpha}(2,\mathbb{C})$,
see the remarks at the end of Subsection \ref{sec:characters_of_UqRb}.   This point will not play a significant role in what follows, thanks to the observation after Theorem \ref{thm:slq2_intertwiners} that the intertwining operators corresponding to different lifts of $\lambda_\alpha$ to $\CC / i\hbar_\alpha^{-1} \mathbb{Z}$ coincide.

Recall from Subsection \ref{sec:YD_intertwiners} that the intertwiners of $ SL_{q_\alpha}(2,\mathbb{C}) $-principal series representations are always given by the right regular action of some 
element $ T \in \M(\DF(SU_{q_\alpha}(2))) $. Since $\alpha$ is a simple root, we can 
map $ T $ to an element $ T_\alpha = \iota_\alpha(T) \in \M(\DF(K_q)) $
\nomenclature[s$T_\alpha$]{$T_\alpha$}{image of an element $T\in\M(\DF(SU_{q_\alpha}(2))) $ under the embedding $\iota_\alpha: \DF(SU_{q_\alpha}(2)) \rightarrow \M(\DF(K_q))$ associated to a simple root $\alpha$}%
\nomenclature{$T_\alpha$}{image of an element $T\in\M(\DF(SU_{q_\alpha}(2))) $ under the embedding $\iota_\alpha: \DF(SU_{q_\alpha}(2)) \rightarrow \M(\DF(K_q))$ associated to a simple root $\alpha$}%
as described above. We will show that the right regular action of $ T_\alpha $ is itself an intertwiner of appropriate $ G_q $-principal series representations.

\begin{lemma} \label{lem:simple_intertwiners}
	Let $\alpha$ be a simple root, and suppose the right regular action of $ T \in \M(\DF(SU_{q_\alpha}(2))) $ defines an intertwiner of $ SL_{q_\alpha}(2, \mathbb{C}) $-principal series representations
	$$
	T: \Gamma(\E_{\mu_1,\lambda_1}) \rightarrow \Gamma(\E_{\mu_2,\lambda_2})
	$$
	for some $ (\mu_1,\lambda_1),~(\mu_2,\lambda_2) \in \frac{1}{2} \mathbb{Z}\times \mathbb{C}/i \hbar_\alpha^{-1} \mathbb{Z} $. 
	Then the right regular action of the corresponding element $ T_\alpha \in \M(\DF(K_q)) $ defines an intertwiner of $ G_q $-principal series representations
	$$
	T_\alpha: \Gamma(\E_{\mu' + \mu_1 \alpha, \lambda' + \lambda_1 \alpha}) \rightarrow \Gamma(\E_{\mu' + \mu_2 \alpha, \lambda' + \lambda_2 \alpha}),
	$$
	for any $ s_\alpha $-fixed points $ \mu' \in \mathfrak{h}^* $ and $ \lambda' \in \mathfrak{h}_q^* $ such that $ \mu' + \mu_i \alpha \in \weights $ for $i=1,2$.  Moreover, the intertwiner $T_\alpha$ is injective or surjective if and only if $T: \Gamma(\E_{\mu_1,\lambda_1}) \rightarrow \Gamma(\E_{\mu_2,\lambda_2})$ is injective or surjective, respectively.  Likewise $T_\alpha$ is injective on the minimal $K_q$-type if and only if $T: \Gamma(\E_{\mu_1,\lambda_1}) \rightarrow \Gamma(\E_{\mu_2,\lambda_2})$ is injective on the minimal $SU_{q_\alpha}(2)$-type.
\end{lemma}

\begin{proof}
	We will appeal to Lemma \ref{lem:intertwiner_condition}. Let $ V $ and $ W $ be finite dimensional essential $ \DF(K_q) $-modules and 
	let $ w \in W_{\mu' + \mu_1 \alpha} $. We need to fix a weight basis for $ V $; let us do so by first decomposing $ V $ into $ \DF(S_\alpha) $-submodules, 
	$ V = \bigoplus_i V_i $, and then fixing a weight basis $ (e_{i,j})_{j = 1}^{N_i} $ for each $ V_i $. Let $ \epsilon_{i,j} $ denote the weight of $ e_{i,j} $. 
	Note that the weights in any fixed $ V_i $ differ by multiples of $ \alpha $, 
	so that we can write
	\[
	\epsilon_{i,j} = \epsilon'_i + k_{i,j} \alpha,
	\]
	where $k_{i,j} = \half(\epsilon_{i,j},\alpha^\vee) \in \half\ZZ$ is the weight of $e_{i,j}$ in $V_i$ as a $ \DF(S_\alpha) $-module,  and the orthogonal component $\epsilon'_i$ does not depend on $j$.  The analogous orthogonal decomposition of $\rho$ is
	\[
	\textstyle\rho = \rho' + \rho_\alpha\alpha = \rho' + \half\alpha,
	\]
	where $\rho' = \rho - \half\alpha$.  Note that $\half \alpha$ is itself the half-sum of positive roots associated to the quantum subgroup $S_\alpha$.

	Now we calculate, as in Lemma \ref{lem:intertwiner_condition},
	\begin{align*}
	\sum_{i,j} & q^{(\lambda' + \lambda_1 \alpha + 2 \rho, \epsilon_{i,j})} T_\alpha (e^{i,j} \otimes  w \otimes e_{i,j}) \\
	&= \sum_i q^{(\lambda' + 2\rho', \epsilon'_i)} \left(\sum_j q^{(\lambda_1 \alpha + \alpha, k_{i,j} \alpha)} T_\alpha (e^{i,j} \otimes  w \otimes e_{i,j}) \right) \\
	&
	= \sum_i q^{(\lambda' + 2\rho', \epsilon'_i)} \left(\sum_j q_\alpha^{2(\lambda_1+1) k_{i,j}} T_\alpha (e^{i,j} \otimes  w \otimes e_{i,j}) \right).
	\end{align*}
	Recall that for the Lie algebra $\mathfrak{sl}(2,\CC)$ with the identification $\lie{h}^* \cong \CC$, the bilinear form on $\lie{h}^*$ is given by $(\nu_1,\nu_2) = 2\nu_1\nu_2$ and the half-sum of positive roots is $\half$.  Therefore, since $T$ defines an $ SL_{q_\alpha}(2,\mathbb{C}) $-intertwiner from $\Gamma(\E_{\mu_1,\lambda_1})$ to $\Gamma(\E_{\mu_2,\lambda_2})$, Lemma \ref{lem:intertwiner_condition} for the group $ SL_{q_\alpha}(2,\mathbb{C}) $ shows that the above sum equals
	\begin{align*}
	\sum_i & q^{(\lambda' + 2\rho', \epsilon'_i)} \left(\sum_j q_\alpha^{2(\lambda_2+1) k_{i,j}} e^{i,j} \otimes T_\alpha(w) \otimes e_{i,j} \right) \\
	&= \sum_{i,j}  q^{(\lambda' + \lambda_2 \alpha + 2\rho, \epsilon_{i,j})} e^{i,j} \otimes T_\alpha(w) \otimes e_{i,j}.
	\end{align*}
	Thus $ T_\alpha $ satisfies the condition of Lemma \ref{lem:intertwiner_condition} for $ G_q $, and so
	$ T_\alpha: \Gamma(\E_{s_\alpha \mu, s_\alpha \lambda}) \rightarrow \Gamma(\E_{\mu,\lambda})$ is an intertwiner of $ G_q $-representations. 
	
	For the statement concerning injectivity, note that $\Gamma(\E_{\mu'+\mu_1\alpha,\lambda'+\lambda_1\alpha})$ is spanned by matrix coefficients of the form $\xi = \bra v'|\bullet| v \ket$ with $v'\in V(\nu)^*$ and $v\in V(\nu)_{\mu'+\mu_1\alpha}$ where $\nu\in\weights^+$, and where we may moreover assume that $v$ belongs to some simple $\DF(S_\alpha)$-submodule of $V(\nu)$.  Note that $v$ has weight $\mu_1$ for the $\DF(S_\alpha)$-action, and so the highest weight of its $\DF(S_\alpha)$-submodule is equal to $|\mu_1|+n$ for some $n\in\mathbb{N}_0$.  Moreover, every such highest weight in $|\mu_1|+\mathbf{N}_0$ can occur for $\nu$ sufficiently large.  
	Since $T_\alpha(\bra v'|\bullet| v\ket ) = \bra v'|\bullet| \iota_\alpha(T)\cdot v\ket$, it follows that $T_\alpha$ annihilates $\xi$ if and only if the action of $T$ annihilates the $\mu_1$-weight space of $V(|\mu_1|+n)$.  Therefore $T_\alpha$ is injective if and only if $T$ acts injectively on the $\mu_1$-weight space of every simple $\DF(SU_{q_\alpha}(2))$-module, and hence if and only if
	$T:\Gamma(\E_{\mu_1,\lambda_1}) \to \Gamma(\E_{\mu_2,\lambda_2})$ is injective.  
	
	The statement about surjectivity is proven similarly.  
	
	Finally, suppose $\xi = \bra v'|\bullet| v \ket$ belongs to the minimal $K_q$-type in $\Gamma(\E_{\mu_1,\lambda_1})$.  Then $v$ is an extremal weight vector in $V(\mu'+\mu_1\alpha)$, meaning that its weight lies in the Weyl orbit of the highest weight $\mu'+\mu_1\alpha$.  In particular, $v$ is annihilated by either $E_\alpha$ or $F_\alpha$, so that it is an extremal weight vector in the $\DF(S_\alpha)$-submodule it generates, which must therefore have highest weight $|\mu_1|$.  
	Since $|\mu_1|$ is the minimal $S_\alpha$-type of $\Gamma(\E_{\mu_1,\lambda_1})$, it follows that $T_\alpha$ annihilates $\xi$ if and only if $T$ annihilates the minimal $S_\alpha$-type of $\Gamma(\E_{\mu_1,\lambda_1})$.
\end{proof}

\begin{theorem} \label{thm:simple_intertwiners}
	Fix a simple root $ \alpha $. For any $ (\mu, \lambda) \in \weights \times \mathfrak{h}_q^* $ there exists a nonzero intertwiner
	\begin{equation*} \label{eq:simple_intertwiner}
	T_\alpha : \Gamma(\E_{\mu,\lambda}) \rightarrow \Gamma(\E_{s_\alpha\mu,s_\alpha\lambda}),
	\end{equation*}
	where $ T \in \DF(SU_{q_\alpha}(2)) $ 
	implements the intertwiner $ T: \Gamma(\E_{\mu_\alpha, \lambda_\alpha}) \rightarrow \Gamma(\E_{-\mu_\alpha, -\lambda_\alpha}) $ 
	of $ SL_{q_\alpha}(2, \mathbb{C})^0 $-principal series. 
	We can identify the following particular cases:
	\begin{bnum}
		\item[a)] If $ \lambda_\alpha \in -|\mu_\alpha| - \mathbb{N} \bmod \frac{1}{2} i \hbar_\alpha^{-1}\mathbb{Z} $, then the intertwiner is zero on the 
		minimal $ K_q $-type of $ \Gamma(\E_{\mu,\lambda}) $.
		\item[b)] If $ \lambda_\alpha \in |\mu_\alpha| + \mathbb{N} \bmod \frac{1}{2} i \hbar_\alpha^{-1}\mathbb{Z} $, then the intertwiner is nonzero on 
		the minimal $ K_q $-type, but is not bijective.
		\item[c)]
		If $\pm \lambda_\alpha \notin |\mu_\alpha| + \mathbb{N} \bmod \frac{1}{2} i \hbar_\alpha^{-1} \mathbb{Z} $, the intertwiner is bijective. 
	\end{bnum}
	
	Moreover, if $ \lambda_\alpha \in -(|\mu_\alpha| + \mathbb{N}) \bmod \frac{1}{2} i \hbar_\alpha^{-1} \mathbb{Z} $, then we have 
	additional \u{Z}elobenko intertwiners as follows:
  \begin{equation*} \label{eq:diamond3}
	\xymatrix{
		& \Gamma(\E_{\mu+m\alpha,\lambda+m\alpha}) \ar@{<->}[dd]^(0.25)\cong \ar@{^{(}->}[dr]
		\\
		\Gamma(\E_{\mu,\lambda}) \ar@<0.3ex>[rr]
		\ar@{->>}[ur]
		\ar@{->>}[dr]
		&& \Gamma(\E_{s_\alpha\mu,s_\alpha\lambda}) 
		\ar@<0.3ex>[ll] 
		\\
		& \Gamma(\E_{\mu-n\alpha,\lambda+n\alpha})  
		\ar@{^{(}->}[ur]
	}
	\end{equation*}
	where $ m = 2\Re(r_\alpha)+1 = -\mu_\alpha - \Re(\lambda_\alpha)$, $ n = 2\Re(l_\alpha)+1 = \mu_\alpha-\Re(\lambda_\alpha)$, and where the horizontal and vertical arrows are those from $a)$, $b)$ and $c)$ above.
\end{theorem}

\begin{proof}
	Let $ (\mu, \lambda) \in \weights \times \mathfrak{h}_q^*$.  As usual, we write $\mu' = \mu - \mu_\alpha\alpha$ for the component of $\mu$ orthogonal to $\alpha$, so that
	$$
	\mu = \mu' + \mu_\alpha \alpha, \qquad s_\alpha \mu = \mu' - \mu_\alpha \alpha.
	$$
	With $ \lambda $ we must be a little careful since in general $ \lambda_\alpha $ is only well-defined 
	modulo $ \frac{1}{2} i \hbar_\alpha^{-1}\mathbb{Z} $ rather than $ i \hbar_\alpha^{-1} \mathbb{Z} $. 
	Let $ \tilde{\lambda} \in \mathfrak{h}^* $ be 
	any lift of $ \lambda \in \mathfrak{h}_q^* $ and put $\tilde{\lambda}' = \tilde{\lambda} - \tilde{\lambda}_\alpha\alpha$.  If we write $\lambda'$ for the projection of $\tilde{\lambda}'$ in $\lie{h}^*_q$, then 
	$$
	\lambda = \lambda' + \tilde{\lambda}_\alpha \alpha, \qquad s_\alpha \lambda = \lambda' - \tilde{\lambda}_\alpha \alpha.
	$$
	Hence applying Lemma \ref{lem:simple_intertwiners} to the standard $ SL_{q_\alpha}(2,\mathbb{C}) $-intertwiner 
	$ T: \Gamma(\E_{\mu_\alpha, \tilde{\lambda}_\alpha}) \rightarrow \Gamma(\E_{-\mu_\alpha, -\tilde{\lambda}_\alpha}) $ from Theorem \ref{thm:slq2_intertwiners} \emph{(2)} yields the intertwiner
	$$
	T_\alpha : \Gamma(\E_{\mu, \lambda}) \rightarrow \Gamma(\E_{s_\alpha\mu, s_\alpha\lambda}).
	$$
	We remark that a different choice of lift $\tilde{\lambda}$ would only alter $\tilde{\lambda}_\alpha$ by an integer multiple of $\half i\hbar_\alpha^{-1}$, and so would result in the same intertwiner $T_\alpha : \Gamma(\E_{\mu, \lambda}) \rightarrow \Gamma(\E_{s_\alpha\mu, s_\alpha\lambda})$ by the remark after Theorem \ref{thm:slq2_intertwiners}.
	
	Using Lemma \ref{lem:simple_intertwiners}, the statements about bijectivity and action on minimal $K_q$-types follow from the analogous statements in Theorem \ref{thm:slq2_intertwiners}.
 	
	The construction of the \u{Z}elobenko intertwiners runs similarly. 
	Suppose that $ \lambda_\alpha \in -(|\mu_\alpha| + \mathbb{N}) \bmod \frac{1}{2} i\hbar_\alpha^{-1}\mathbb{Z} $. 
	As above, we lift $\lambda \in \lie{h}_q^*$ to $\tilde{\lambda}\in\lie{h}^*$ and write 
	\[
	\tilde{\lambda} = \tilde{\lambda}' + \tilde{\lambda}_\alpha\alpha,
	\]
	where $\tilde{\lambda}'$ is fixed by $s_\alpha$ and $\tilde{\lambda}_\alpha \in  -(|\mu_\alpha| + \mathbb{N}) \bmod \frac{1}{2} i\hbar_\alpha^{-1}\mathbb{Z} $.  This means we can apply Lemma \ref{lem:simple_intertwiners} to one of the two diagrams of \v{Z}elobenko intertwiners in Theorem \ref{thm:slq2_intertwiners}.  However, it will simplify the exposition if we arrange to have $\tilde{\lambda}_\alpha \in  -(|\mu_\alpha| + \mathbb{N})$ on the nose.    To this end, suppose that the imaginary part of $\tilde{\lambda}_\alpha$ is equal to $\frac{k}{2}\hbar_\alpha^{-1}$ for some $k\in\mathbb{Z}$.  Note that $\half i \hbar^{-1} \alpha^\vee$ is fixed by $s_\alpha$ modulo $i \hbar^{-1} \roots^\vee$, so that if we decompose $\tilde{\lambda}$ as 
	\[
	 	\tilde{\lambda} = (\tilde{\lambda}' +\frac{k}{2} i \hbar^{-1} \alpha^\vee) + (\tilde{\lambda}_\alpha -\frac{k}{2} i \hbar_\alpha^{-1} )\alpha,
	\]
	then the first term is still fixed by $s_\alpha$ when understood as an element of $\lie{h}_q^*$.  Thus, without loss of generality, we may assume that $ \tilde{\lambda}_\alpha \in -(|\mu_\alpha| + \mathbb{N})$ on the nose.

	In this situation, if we define $m = -\mu_\alpha-\tilde{\lambda}_\alpha$ and $n=\mu_\alpha-\tilde{\lambda}_\alpha$, then
	the various principal series parameters in the diagram in the statement decompose as
		\begin{align*}
		\mu &= \mu' + \mu_\alpha \alpha, &  \lambda &= \tilde{\lambda}' + \tilde{\lambda}_\alpha \alpha, 
		\\
		\mu + m\alpha &= \mu' - \tilde{\lambda}_\alpha \alpha, 
		& \lambda +m \alpha &= \tilde{\lambda}' - \mu_\alpha \alpha, 
		\\
		\mu - n \alpha &= \mu' + \tilde{\lambda}_\alpha \alpha, 
		& \lambda +n \alpha &= \tilde{\lambda}' + \mu_\alpha \alpha, 
		\\
		s_\alpha \mu &= \mu' -\mu_\alpha \alpha, 
		& s_\alpha \lambda &= \tilde{\lambda}' - \tilde{\lambda}_\alpha \alpha.
		\end{align*}
	Applying Lemma \ref{lem:simple_intertwiners} to the intertwiners in part $3a)$ of Theorem \ref{thm:slq2_intertwiners}, we obtain the diagram 
	of \u{Z}elobenko intertwiners as claimed. 
\end{proof}

\subsubsection{Explicit formulas for intertwiners}
\label{sec:intertwiner_formulas}

In this subsection, we will give explicit formulas for some of the intertwining operators between principal series representations.    As usual, if $(\mu,\lambda)\in\weights\times\lie{h}_q^*$ are parameters for the principal series, we will use $(l,r)\in\lie{h}_q^*\times\lie{h}_q^*$ to denote a pair such that
\[
  \mu = l-r, \quad \lambda+2\rho = -l-r.
\]

We begin with the \v{Z}elobenko intertwiners.  Recall from Theorem \ref{Vermatheorem} that if $r'\in\lie{h}_q^*$ is strongly linked to $r$, then we have an inclusion $M(r') \subset M(r)$.  Recall also the algebra involution
$\theta = \tau \hat{S}$ of $U_q(\lie{g})$, which is given on generators by
\[
  \theta(E_i) = -F_i, \quad
  \theta(F_i) = -E_i, \quad
  \theta(K_\nu) = K_{-\nu},
\]
see the remarks before Corollary \ref{cor:Verma_module_iso}.

\begin{theorem}
	\label{thm:Zelobenko_explicit}
  Let $(\mu,\lambda)\in\weights\times\lie{h}_q^*$ with associated parameters $(l,r)\in\lie{h}_q^*\times\lie{h}_q^*$.
  \begin{bnum}
  	\item[$a)$]
  	Suppose that $r'\in\lie{h}_q^*$ is strongly linked to $r$ and fix $Y\in U_q(\lie{n}_-)$ such that $Y\cdot v_r\in M(r)$ is a primitive vector of weight $r'$.  Then the right regular action of $\theta(Y)$ on $\Gamma(\E_\mu)$ defines an intertwiner
  	\[
  	  \theta(Y)\hit \bullet : \Gamma(\E_{\mu,\lambda}) \to \Gamma(\E_{\mu',\lambda'})
  	\]
  	where $\mu'=l-r'$, $\lambda'+2\rho = -l-r'$.
  	
  	\item[$b)$] 
  	Suppose that $l'\in\lie{h}_q^*$ is strongly linked to $l$ and fix $Y\in U_q(\lie{n}_-)$ such that $Y\cdot v_l \in M(l)$ is a primitive vector of weight $l'$.  Then the right regular action of $Y$ on $\Gamma(\E_\mu)$ defines an intertwiner
  	\[
    	Y\hit \bullet : \Gamma(\E_{\mu,\lambda}) \to \Gamma(\E_{\mu',\lambda'})
  	\]
  	where $\mu'=l'-r$, $\lambda'+2\rho = -l'-r$.
  \end{bnum}
\end{theorem}

\begin{proof}
	$a)$  
	It will be convenient to work with the involution $\theta' = \hat{S}\tau$ in place of $\theta=\tau\hat{S}$ in the calculations to follow. 
	Note that $Y$ is of weight $r'-r$ for the adjoint action,
	and hence 
	\[
	 \theta'(Y) = \hat{S}^2(\theta(Y)) = K_{2\rho} \theta(Y) K_{-2\rho} = q^{(2\rho,r-r')} \theta(Y),
	\]
	so that $\theta'(Y)$ differs from $\theta(Y)$ only by a scalar multiple.  Therefore, it is equivalent to prove the claim with $\theta'$ in place of $\theta$.
  Also,	since  $\mu'-\mu = r-r'$, the right regular action
	\[
	 \theta'(Y)\hit \xi = (\theta'(Y),\xi_{(2)})\xi_{(1)}
	\]
	does indeed map $\Gamma(\E_{\mu})$ to $\Gamma(\E_{\mu'})$.  
	
	From Corollary \ref{cor:Verma_module_iso} we have an isomorphism of $U_q^\mathbb{R}(\lie{g})$-modules $M(\mu,\lambda +2\rho) \cong M(l) \otimes M(r)$ which sends the cyclic vector $v_{\mu,\lambda}$ to $v_l\otimes v_r$.  This requires a particular choice of action of $U_q^\mathbb{R}(\lie{g})$ on $M(l) \otimes M(r)$, and in particular for $X\in U_q(\lie{g})$, we have
	\[
	 (X\bowtie 1) \cdot (m\otimes n) = X_{(1)} \cdot m \otimes \theta'(X_{(2)}) \cdot n.
	\]
  Note that $\theta'(Y) \in U_q(\lie{n}_+)$, and so $\hat{\Delta}(\theta'(Y)) \in U_q(\lie{n}_+) \otimes U_q(\lie{b}_+)$.  Since the action of $U_q(\lie{n}_+)$ on $v_l$ factors through $\hat{\epsilon}$, we get
	\begin{align*}
	 (\theta'(Y)\bowtie 1) \cdot (v_l\otimes v_r)	
	 	 &= \theta'(Y)_{(1)}\cdot v_l \otimes \theta'(\theta'(Y)_{(2)}) \cdot v_r \\
	   &= v_l \otimes \theta'(\theta'(Y)) \cdot v_r \\
	   &= v_l \otimes Y\cdot v_r.
	\end{align*}
	Therefore the inclusion $M(l) \otimes M(r') \to M(l) \otimes M(r)$ gives rise to an inclusion of $U_q^\mathbb{R}(\lie{g})$-modules $ j: M(\mu',\lambda'+2\rho) \to M(\mu,\lambda+2\rho)	$ such that $j(v_{\mu',\lambda'+2\rho}) = (\theta'(Y)\bowtie 1)\cdot v_{\mu,\lambda+2\rho}$.
	
	Let $\ext : \Gamma(\E_{\mu,\lambda}) \to \ind_{B_q}^{G_q}(\mathbb{C}_{\mu,\lambda})$ denote the extension from the compact to the induced picture, see Lemma \ref{lem:ext}, so that
	$
	  \ext(\xi)  = \xi\otimes K_{2\rho+\lambda},
  $
	for all $\xi\in\Gamma(\E_{\mu,\lambda})$.  For any $X\bowtie f \in U_q^\mathbb{R}(\lie{g})$, the pairing $M(\mu,\lambda+2\rho) \times \ind_{B_q}^{G_q} (\mathbb{C}_{\mu,\lambda}) \to \mathbb{C}$ from Lemma \ref{lem:Verma-principal-pairing} satisfies
	\begin{align*}
	  (j((X\bowtie f)\cdot v_{\mu',\lambda'+2\rho}), \ext(\xi)) 
	   &= (f,K_{\lambda'+2\rho}) \, (j((X\bowtie 1) \cdot v_{\mu',\lambda'+2\rho}) , \xi\otimes K_{\lambda+2\rho} )
	   \\
	   &= (f,K_{\lambda'+2\rho})\, ((X \theta'(Y) \bowtie 1) \cdot v_{\mu,\lambda+2\rho} , \xi \otimes K_{\lambda+2\rho})
	   \\
	   &= (f,K_{\lambda'+2\rho})\, (X , \xi_{(1)}) \, (\theta'(Y), \xi_{(2)})
	   \\
	   &= ((X \bowtie f)\cdot v_{\mu',\lambda'+2\rho} , \ext(\theta'(Y)\hit \xi)),
	\end{align*}
	for $\xi\in\Gamma(\E_{\mu,\lambda})$.
	Thus $j$ is dual to the right regular action $\theta'(Y) \hit\bullet:\Gamma(\E_{\mu,\lambda}) \to \Gamma(\E_{\mu',\lambda'})$ under this pairing.  By Lemma \ref{lem:Verma-principal-pairing}, we conclude that the right regular action of $\theta'(Y)$ is $U_q^\mathbb{R}(\lie{g})$-linear.
	
	$b)$  The inclusion $M(l') \otimes M(r) \to M(l) \otimes M(r)$ corresponds to a map
	\[
	 j: M(\mu',\lambda'+2\rho) \to M(\mu,\lambda+2\rho)
	\]
	which sends $v_{\mu',\lambda'+2\rho}$ to $(Y\bowtie 1)\cdot v_{\mu,\lambda+2\rho}$, and this map is dual to the right regular action of $Y$ from $\Gamma(\E_{\mu,\lambda})$ to $\Gamma(\E_{\mu',\lambda'})$ under the pairing from Lemma \ref{lem:Verma-principal-pairing}.  The rest of the argument is analogous to part $a)$.
\end{proof}

In particular, let us consider the \v{Z}elobenko intertwiners corresponding to simple roots, that is, the diagonal arrows in the diagram from Theorem \ref{thm:simple_intertwiners} part \emph{(3)}.  Suppose that $(\mu,\lambda) \in \weights \times \lie{h}_q^*$ is such that $\lambda_\alpha \in -|\mu_\alpha|-\mathbb{N} \bmod \textstyle\half i \hbar_\alpha^{-1}\mathbb{Z}$.
This is equivalent to the condition that $l_\alpha$ and $r_\alpha$ belong to $\half\mathbb{N}_0 \bmod \frac{1}{4} i \hbar_\alpha^{-1}\mathbb{Z}$.  Therefore, putting $l' = l - 2\Re(l_\alpha)\alpha$ and $r' = r - 2\Re(r_\alpha)\alpha$ we have inclusions of Verma modules
\[
M(l') \subset M(l) , \qquad M(r')\subset M(r).
\]
Explicitly, $F_\alpha^{2\Re(l_\alpha)}\cdot v_l$ is a primitive vector of weight $l'$ in $M(l)$ and $F_\alpha^{2\Re(r_\alpha)}\cdot v_r$ is a primitive vector of weight $r'$ in $M(r)$.  Therefore, an application of Theorem \ref{thm:Zelobenko_explicit} gives the diagram of intertwiners
\begin{equation} \label{eq:diamond3}
\xymatrix{
	& \Gamma(\E_{\mu+m\alpha,\lambda+m\alpha}) \ar@{<->}[dd]^(0.25)\cong \ar@{^{(}->}[dr]^{F_\alpha^{n}} 
	\\
	\Gamma(\E_{\mu,\lambda}) \ar@<0.3ex>[rr]
	\ar@{->>}[ur]^{E_\alpha^{m}} 
	\ar@{->>}[dr]_{F_\alpha^{n}} 
	&& \Gamma(\E_{s_\alpha\mu,s_\alpha\lambda}) 
	\ar@<0.3ex>[ll]
	\\
	& \Gamma(\E_{\mu-n\alpha,\lambda+n\alpha})  
	\ar@{^{(}->}[ur]_{E_\alpha^{m}}
}
\end{equation}

Next, we give an explicit formula for the standard intertwiners
$
  T_\alpha : \Gamma(\E_{\mu,\lambda}) \to \Gamma(\E_{s_\alpha\mu,s_\alpha,\lambda})
$
corresponding to a simple root, as in Theorem \ref{thm:simple_intertwiners} \emph{(2)}.  

We begin with the case of $G_q=SL_q(2,\mathbb{C})$.  Let $(\mu,\lambda) \in \half\mathbb{Z} \times \mathbb{C}/i\hbar^{-1}\mathbb{Z}$.  Recall that we have 
\[
  \Gamma(\E_{\mu,\lambda}) \cong \bigoplus_{n\in|\mu|+\mathbb{N}_0} V(n)^* \otimes V(n)_\mu.
\]
and therefore any $\DF(K_q)$-linear map $T:\Gamma(\E_{\mu,\lambda}) \to \Gamma(\E_{-\mu,-\lambda})$ is specified uniquely by a sequence of linear maps
\[
 T_n : V(n)_\mu \to V(n)_{-\mu}
\]
for all $n\in|\mu|+\mathbb{N}_0$.  Since the weight spaces of simple $\DF(SU_q(2))$-modules are all one dimensional, the map $T_n$ is specified by a scalar once we have fixed a weight basis of $V(n)$.

For this, we will use the orthonormal weight basis with respect to the invariant inner product on $V(n)$.  In comparison to the weight basis $v_n, v_{n-1}, \ldots, v_{-n}$ from Theorem \ref{vnexplicit}, the orthonormal basis is given by
\[
 e^n_m = q^{\half(n-m)(n+m-1)} 
   \left( \frac{[n+m]_q!}{[2n]_q! [n-m]_q!} \right)^\half v_m.
\]
\nomenclature{$e^n_m$}{orthonormal basis element for the unitary $\DF(SU_q(2))$-representation $V(n)$}%
where $m=-n,-n+1,\ldots,n$.  
To see that it is indeed orthonormal for the invariant inner product, one can first confirm that the action of $U_q(\lie{sl}(2))$ on this basis is given by
\begin{align*}
K \cdot e^n_m &= q^{m} e^n_m, \hspace{2cm} K_\alpha \cdot e^n_m = q^{2m} e^n_m,   \\
E \cdot e^n_m &= q^m [n - m]_q^\frac{1}{2} [n + m + 1]_q^\frac{1}{2} e^n_{m + 1}, \\
F \cdot e^n_m &= q^{-(m-1)} [n + m]_q^\frac{1}{2} [n - m + 1]_q^\frac{1}{2} e^n_{m - 1},
\end{align*}
and then check these formulas are compatible with the involution on  $U_q^\mathbb{R}(\lie{su}(2))$ when they are declared to be orthonormal,
compare Section 3.2.1 in \cite{KS}. We point out that one must be careful with statements about unitary modules in \cite{KS}, since 
there are different versions of the enveloping algebra in play, and their preferred real forms do not always agree under the Hopf algebra morphisms they use.

We can then specify a $\DF(K_q)$-linear operator $ T: \Gamma(\E_{\mu, \lambda}) \rightarrow \Gamma(\E_{-\mu,-\lambda}) $ by 
the scalars $ (T_n)_{n \in |\mu| + \mathbb{N}_0} $ such that $ T(e^n_{\mu}) = T_n e^n_{-\mu} $. 
Our reason for working in this basis is that we will be using the Clebsch-Gordan formulas from \cite{KS}.
Looking at Lemma \ref{lem:intertwiner_condition}, we will need to consider the action of $ T $ on tensor products of the form
$$
e^{\frac{1}{2} \dual}_{\pm \frac{1}{2}} \otimes e^n_\mu \otimes e^\frac{1}{2}_{\pm \frac{1}{2}}
$$
for $ n \in |\mu| + \mathbb{N}_0 $, where $e^{n\dual}_m$ denotes an element of the basis of $\precon V(n)$ dual to $(e^n_m)$.
\nomenclature{$e^{n\dual}_m$}{elements of the basis of $\precon V(n)$ dual to $e^n_m$}%
There is a double multiplicity of the representation $V(n)$ in the tensor product $V(\half) \otimes V(n) \otimes V(\half)$, and to deal with this we will write 
$$
V(\tfrac{1}{2}) \otimes (V(n) \otimes V(\tfrac{1}{2})) \cong V(n + 1) \oplus V(n^+) \oplus V(n^-) \oplus V(n - 1),
$$
where $ V(n^\pm) $ signifies the copy of $ V(n) $ contained in $ V(\tfrac{1}{2}) \otimes V(n \pm \tfrac{1}{2})$, respectively. 

The Clebsch-Gordan coefficients for tensor products with the fundamental representation are as follows:
\begin{align*}
e^n_\mu \otimes e^\frac{1}{2}_\frac{1}{2} = 
\textstyle q^{\frac{1}{2}(n - \mu)} \frac{[n + \mu + 1]^\frac{1}{2}}{[2n + 1]^\frac{1}{2}} e^{n + \frac{1}{2}}_{\mu + \frac{1}{2}}
- q^{-\frac{1}{2}(n + \mu + 1)} \frac{[n - \mu]^\frac{1}{2}}{[2n + 1]^\frac{1}{2}} e^{n - \frac{1}{2}}_{\mu + \frac{1}{2}}, 
\\
e^n_\mu \otimes e^\frac{1}{2}_{-\frac{1}{2}} = 
\textstyle q^{-\frac{1}{2}(n + \mu)} \frac{[n - \mu + 1]^\frac{1}{2}}{[2n + 1]^\frac{1}{2}} e^{n + \frac{1}{2}}_{\mu - \frac{1}{2}}
+ q^{\frac{1}{2}(n - \mu + 1)} \frac{[n + \mu]^\frac{1}{2}}{[2n + 1]^\frac{1}{2}} e^{n - \frac{1}{2}}_{\mu - \frac{1}{2}},  
\\
e^\frac{1}{2}_\frac{1}{2} \otimes e^n_\mu =
\textstyle q^{-\frac{1}{2}(n - \mu)} \frac{[n + \mu + 1]^\frac{1}{2}}{[2n + 1]^\frac{1}{2}} e^{n + \frac{1}{2}}_{\mu + \frac{1}{2}}
+ q^{\frac{1}{2}(n + \mu + 1)} \frac{[n - \mu]^\frac{1}{2}}{[2n + 1]^\frac{1}{2}} e^{n - \frac{1}{2}}_{\mu + \frac{1}{2}},  
\\
e^\frac{1}{2}_{-\frac{1}{2}} \otimes e^n_\mu =
\textstyle q^{\frac{1}{2}(n + \mu)} \frac{[n - \mu + 1]^\frac{1}{2}}{[2n + 1]^\frac{1}{2}} e^{n + \frac{1}{2}}_{\mu - \frac{1}{2}} 
- q^{-\frac{1}{2}(n - \mu + 1)} \frac{[n + \mu]^{\frac{1}{2}}}{[2n + 1]^{\frac{1}{2}}} e^{n - \frac{1}{2}}_{\mu - \frac{1}{2}}.
\end{align*}
For this, see e.g. Equations (3.68)--(3.69) in \cite{KS}, but note the typographical error in Equation (3.68), where the 
factor $ [l \pm j + \frac{1}{2}]^\frac{1}{2} $ should be $ [l \pm j + 1]^\frac{1}{2} $.

For $ SU_q(2) $, the precontragredient representation $ \precon V(n)$ is isomorphic to $ V(n) $. 
In particular, for the fundamental representation the following map is an isomorphism: 
\begin{equation*}
e^{\frac{1}{2}\dual}_\frac{1}{2} \mapsto i q^{-\frac{1}{2}} e^\frac{1}{2}_{-\frac{1}{2}}, \qquad 
e^{\frac{1}{2}\dual}_{-\frac{1}{2}} \mapsto -i q^{\frac{1}{2}} e^\frac{1}{2}_\frac{1}{2}.
\end{equation*}

Combining this with the Clebsch-Gordan formulae above gives
\begin{align*}
e^{\frac{1}{2}\dual}_{\frac{1}{2}} \otimes (e^n_{\mu} \otimes e^\frac{1}{2}_{\frac{1}{2}})
&= iq^{-\frac{1}{2}} \left( e^{\frac{1}{2}}_{-\frac{1}{2}} \otimes (e^n_{\mu} \otimes e^\frac{1}{2}_{\frac{1}{2}}) \right)\\
&= \textstyle i q^{\frac{1}{2}(n-\mu-1)} \frac{[n+\mu+1]^\frac{1}{2}}{[2n+1]^\frac{1}{2}} 
\left(e^\frac{1}{2}_{-\frac{1}{2}} \otimes e^{n+\frac{1}{2}}_{\mu+\frac{1}{2}} \right) \\
& \textstyle  - i q^{-\frac{1}{2}(n+\mu+2)} \frac{[n+\mu+1]^\frac{1}{2}}{[2n+1]^\frac{1}{2}} 
\left( e^\frac{1}{2}_{-\frac{1}{2}} \otimes e^{n-\frac{1}{2}}_{\mu+\frac{1}{2}} \right) \\
&=\ \textstyle i q^{n} \frac{[n+{\mu}+1]^\frac{1}{2} [n-{\mu}+1]^\frac{1}{2}}{[2n+1]^\frac{1}{2} [2n+2]^\frac{1}{2}} e^{n+1}_{\mu}
- iq^{-1} \frac{[n+{\mu}+1]}{[2n+1]^\frac{1}{2} [2n+2]^\frac{1}{2}} e^{n^+}_{\mu}  \\
& \textstyle - iq^{-1} \frac{[n-{\mu}]}{[2n]^\frac{1}{2} [2n+1]^\frac{1}{2}} e^{n^-}_{\mu}
+ iq^{-n-1} \frac{[n+{\mu}]^\frac{1}{2} [n-{\mu}]^\frac{1}{2}}{[2n]^\frac{1}{2} [2n+1]^\frac{1}{2}} e^{n-1}_{\mu}.  \\
e^{\frac{1}{2}\dual}_{-\frac{1}{2}} \otimes (e^n_{\mu} \otimes e^\frac{1}{2}_{-\frac{1}{2}})
&= -iq^\frac{1}{2} \left(e^\frac{1}{2}_\frac{1}{2} \otimes (e^n_\mu \otimes e^\frac{1}{2}_{-\frac{1}{2}}) \right) \\
&= \textstyle -i q^{\frac{1}{2}(-n-\mu+1)} \frac{[n-\mu+1]^\frac{1}{2}}{[2n+1]^\frac{1}{2}} 
\left( e^\frac{1}{2}_\frac{1}{2} \otimes e^{n+\frac{1}{2}}_{\mu-\frac{1}{2}} \right) \\
&\textstyle -i q^{\frac{1}{2}(n-\mu+2)} \frac{[n+\mu]^\frac{1}{2}}{[2n+1]^\frac{1}{2}} 
\left( e^\frac{1}{2}_\frac{1}{2} \otimes e^{n-\frac{1}{2}}_{\mu-\frac{1}{2}} \right) \\
&=\ \textstyle - iq^{-n} \frac{[n+{\mu}+1]^\frac{1}{2} [n-{\mu}+1]^\frac{1}{2}}{[2n+1]^\frac{1}{2} [2n+2]^\frac{1}{2}} e^{n+1}_{\mu}
- iq \frac{[n-{\mu}+1]}{[2n+1]^\frac{1}{2} [2n+2]^\frac{1}{2}} e^{n^+}_{\mu}  \\
& \textstyle - iq \frac{[n+{\mu}]}{[2n]^\frac{1}{2} [2n+1]^\frac{1}{2}} e^{n^-}_{\mu}
- iq^{n+1} \frac{[n+{\mu}]^\frac{1}{2} [n-{\mu}]^\frac{1}{2}}{[2n]^\frac{1}{2} [2n+1]^\frac{1}{2}} e^{n-1}_{\mu}.
\end{align*}
This allows us to write the sums that appear in Lemma \ref{lem:intertwiner_condition} as
\begin{align*}
\sum_{j = \pm\frac{1}{2}}  &q^{(\lambda + 2\rho, \epsilon_j)} e^{\frac{1}{2}\dual}_j \otimes e^n_{\mu} \otimes e^\frac{1}{2}_j \nonumber\\
&= q^{\lambda+1} e^{\frac{1}{2}\dual}_\frac{1}{2} \otimes e^n_{\mu} \otimes e^\frac{1}{2}_\frac{1}{2}
+ q^{-\lambda-1} e^{\frac{1}{2}\dual}_{-\frac{1}{2}} \otimes e^n_{\mu} \otimes e^\frac{1}{2}_{-\frac{1}{2}} \nonumber\\
&= \textstyle i (q^{n+\lambda+1} - q^{-n-\lambda-1} )\frac{[n+{\mu}+1]^\frac{1}{2} [n-{\mu}+1]^\frac{1}{2}}{[2n+1]^\frac{1}{2} [2n+2]^\frac{1}{2}} e^{n+1}_{\mu} \nonumber \\
&\qquad \textstyle  - i\frac{q^{\lambda}[n+{\mu}+1]+q^{-\lambda}[n-\mu+1]}{[2n+1]^\frac{1}{2} [2n+2]^\frac{1}{2}} e^{n^+}_{\mu}  \nonumber \\
&\qquad \textstyle - i \frac{q^\lambda [n-{\mu}]+q^{-\lambda}[n+\mu]}{[2n]^\frac{1}{2} [2n+1]^\frac{1}{2}} e^{n^-}_{\mu} \nonumber \\
&\qquad \textstyle + i (q^{-n+\lambda} - q^{n-\lambda}) \frac{[n+{\mu}]^\frac{1}{2} 
	[n-{\mu}]^\frac{1}{2}}{[2n]^\frac{1}{2} [2n+1]^\frac{1}{2}} e^{n-1}_{\mu} \nonumber \\ 
&= A e^{n+1}_{\mu} + B e^{n^+}_{\mu} + C e^{n^-}_{\mu} + D e^{n-1}_{\mu},
\end{align*}
with coefficients
\begin{align*}
A = A(\mu,\lambda,n) &= i(q-q^{-1}) \frac{[n+\lambda+1] [n+\mu+1]^\frac{1}{2} [n-\mu+1]^\frac{1}{2}}{[2n+1]^\frac{1}{2} [2n+2]^\frac{1}{2}}, \\
B = B(\mu,\lambda,n) &= -i \frac{q^{\lambda} [n+\mu+1] + q^{-\lambda} [n-\mu+1]}{[2n+1]^\frac{1}{2} [2n+2]^\frac{1}{2}}, \\
C = C(\mu,\lambda,n) &= -i \frac{q^{\lambda} [n-\mu] + q^{-\lambda} [n+\mu]}{[2n]^\frac{1}{2} [2n+1]^\frac{1}{2}}, \\
D = D(\mu,\lambda,n) &= -i(q-q^{-1}) \frac{[n-\lambda] [n+\mu]^\frac{1}{2} [n-\mu]^\frac{1}{2}}{[2n]^\frac{1}{2} [2n+1]^\frac{1}{2}}.
\end{align*}
Therefore, the intertwiner condition in Lemma \ref{lem:intertwiner_condition} reduces to the following four conditions on the scalars $ T_n $ 
for $ n \in |\mu_1| + \mathbb{N}_0 $:
\begin{align*}
T_{n + 1} A(\mu_1, \lambda_1,n) &= T_n A(\mu_2,\lambda_2,n),
\\
T_n B(\mu_1,\lambda_1,n) &= T_n B(\mu_2,\lambda_2,n),
\\
T_n C(\mu_1,\lambda_1,n) &= T_n C(\mu_2,\lambda_2,n),
\\
T_{n - 1} D(\mu_1,\lambda_1,n) &= T_n D(\mu_2,\lambda_2,n). 
\end{align*}
Note that when $n=|\mu_1|$ we have 
$D(\mu_1,\lambda_1,n)=0$, so that the last equation should be interpreted as saying $T_nD(\mu_2,\lambda_2,n) = 0$ in this case.

The following theorem gives explicit formulas for the standard intertwiners from part 2 of Theorem \ref{thm:slq2_intertwiners}.

\begin{theorem} \label{thm:slq2_intertwiners_explicit}
	Let $ (\mu, \lambda) \in \half\mathbb{Z} \times \mathbb{C}/i\hbar^{-1}\mathbb{Z} $.   Up to scalar, the unique intertwiner of $SL_q(2,\mathbb{C})$-principal series $ T: \Gamma(\E_{\mu, \lambda}) \rightarrow \Gamma(\E_{-\mu, -\lambda}) $ acts on matrix coefficients as
	\[
	  T (\bra v' | \bullet | e^n_\mu \ket) = T_n \bra v' | \bullet | e^n_{-\mu} \ket
	\]
	where the coefficients $(T_n)_{n\in|\mu|+\mathbb{N}_0}$ are as follows:
	\begin{bnum}
		\item[$a)$] If $ \lambda \in -|\mu| - \mathbb{N} \bmod \frac{1}{2} i \hbar^{-1}\mathbb{Z} $, we have
		\begin{equation*} \label{eq:I}
		T_n = 
		 \begin{cases}
  		0, & |\mu| \leq n < -\Re(\lambda), \\
  	 \prod_{k = -\Re(\lambda) + 1}^{n} \frac{[k-\lambda]_q}{[k+\lambda]_q}, & n \geq -\Re(\lambda).
		 \end{cases}
		\end{equation*}
		\item[$b)$] Otherwise, we have
		\begin{equation*} \label{eq:J}
		 T_n =  \prod_{k = |\mu| + 1}^{n} \frac{[k - \lambda]_q}{[k + \lambda]_q}.  
		\end{equation*}
  \end{bnum}
\end{theorem}

\begin{proof} 	
	We know that such intertwiners exist from Theorem \ref{thm:slq2_intertwiners} \emph{(2)}, so it only remains to deduce the explicit formulas for the coefficients $T_n$ from the recurrence relations stated before the theorem.
	With $(\mu_1,\lambda_1) = (\mu,\lambda)  = (-\mu_2,-\lambda_2)$, the last equation before the theorem reduces to
	\begin{equation*}
	\label{eq:recurrence1}
	T_n [n + \lambda]_q = T_{n - 1}[n-\lambda]_q
	\end{equation*}
	for all $ n \in |\mu| + \mathbb{N} $.    Note that the $q$-numbers $[n \pm \lambda]_q$ are all nonzero except for at most one value as follows:
	\begin{itemize}
		\item if $\lambda \in -|\mu| -\mathbb{N} \bmod \half i\hbar^{-1}\mathbb{Z}$, then $[n+\lambda]_q =0$ for $n = -\Re(\lambda)$;
		\item if $\lambda \in |\mu| + \mathbb{N} \bmod \half i\hbar^{-1}\mathbb{Z}$, then $[n-\lambda]_q =0$ for $n = \Re(\lambda)$;
		\item otherwise, they are all nonzero.
	\end{itemize}
	In every case, we observe that the recurrence relation will have a unique solution up to an overall scalar multiple, and a direct check shows that the values of $T_n$ stated in the Theorem are indeed solutions.  This completes the proof.
	 \end{proof} 

We note that the intertwiners $T:\Gamma(\E_{\mu,\lambda}) \to \Gamma(\E_{-\mu,-\lambda})$ given in part $b)$ of Theorem \ref{thm:slq2_intertwiners_explicit} form a meromorphic family of operators $\Gamma(\E_{\mu}) \to \Gamma(\E_{-\mu})$ as a function of $\lambda$, with poles at the points where $\lambda\in-|\mu|-\mathbb{N} + \half i \hbar^{-1}\mathbb{Z}$.
We also point out, when $\lambda\in|\mu|+\mathbb{N}+\half i \hbar^{-1}\mathbb{Z}$, the coefficients $T_n$ in $b)$ are zero for all $n \geq \Re(\lambda)$, so that the intertwiner $T$ is indeed finite rank in this case.

Using Theorem \ref{thm:simple_intertwiners}, the results of Theorem \ref{thm:slq2_intertwiners_explicit} immediately give an explicit formula for the intertwiners in higher rank associated to a simple reflection.  To specify this formula, we will make use of elements $\Ph(E_\alpha)$ and $\Ph(F_\alpha)$ in $\M(\DF(K_q))$
\nomenclature[o$Ph(E_\alpha)$]{$\Ph(E_\alpha)$, $\Ph(F_\alpha)$}{operator phases of $E_\alpha$ and $F_\alpha$}%
which are defined as the operator phases in the polar decomposition of $E_\alpha$ and $F_\alpha$, respectively.  Explicitly, if $V\cong V(n)$ is a $\DF(SU_{q_\alpha}(2))$-submodule of highest weight $n\in\half\mathbb{N}_0$ inside a $\DF(K_q)$-module, $\Ph(E_\alpha)$ and $\Ph(F_\alpha)$ act as
\[
 \Ph(E_\alpha)\cdot e^n_m = e^n_{m+1}, \qquad \Ph(F_\alpha)\cdot e^n_m = e^n_{m-1},
\]
where $e^n_m, \ldots, e^n_{-m}$ denotes an orthonormal weight basis of $V$ isomorphic to that described above for $V(n)$.  In order to avoid having to deal with different cases, we will abuse notation and write, for any $k\in\mathbb{Z}$
\[
  \Ph(F_\alpha)^k = \begin{cases}
   \Ph(F_\alpha)^k & \text{if } k>0, \\
   1 & \text{if } k=0,\\
   \Ph(E_\alpha)^{-k} & \text{if } k<0. \\
  \end{cases}
\]

We can now give formulas for the standard intertwiners corresponding to simple reflections from Theorem \ref{thm:simple_intertwiners} $a)$, $b)$ and $c)$.

\begin{cor}
	\label{cor:simple_intertwiners_explicit}
	Let $(\mu,\lambda)\in\weights\times\lie{h}_q^*$ and fix a simple root $\alpha$.  The intertwiner $T:\Gamma(\E_{\mu,\lambda}) \to \Gamma(\E_{s_\alpha\mu,s_\alpha\lambda})$ acts on matrix coefficients by
	\[
	  T(\bra v' | \bullet | v\ket) = T_n \,\bra v' | \bullet | \Ph(F_\alpha)^{2\mu_\alpha} v\ket
	\]
	whenever $v$ belongs to a $\DF(SU{q_\alpha}(2))$-submodule of highest weight $n$, and	where the coefficients $T_n\in\mathbb{C}$ are as follows:
		\begin{itemize}
		\item[$a)$] If $ \lambda_\alpha \in -|\mu_\alpha| - \mathbb{N} \bmod \frac{1}{2} i \hbar_\alpha^{-1}\mathbb{Z} $,
		\begin{equation*} \label{eq:I}
		T_n = 
		\begin{cases}
		0, & |\mu_\alpha| \leq n < -\Re(\lambda_\alpha), \\
		\prod_{k = -\Re(\lambda_\alpha) + 1}^n \frac{[k-\lambda_\alpha]_{q_\alpha}}{[k+\lambda_\alpha]_{q_\alpha}}, & n \geq -\Re(\lambda_\alpha).
		\end{cases}
		\end{equation*}
		\item[$b)$] Otherwise,
		\begin{equation*} \label{eq:J2}
		T_n =  \prod_{k = |\mu_\alpha| + 1}^n \frac{[k - \lambda_\alpha]_{q_\alpha}}{[k + \lambda_\alpha]_{q_\alpha}}.  
		\end{equation*}
	\end{itemize}
\end{cor}

\subsection{Submodules and quotient modules of the principal series}

In this section we will elaborate on the structure of the principal series representation $\Gamma(\E_{\mu,\lambda})$.  We closely follow the proofs of the corresponding classical results in \cite{Duflocomplexlnm}.

For $(\mu,\lambda)\in\weights\times\lie{h}_q^*$ and $\alpha\in{\bf\Delta}^+$ any positive root, we continue to use the notation 
\[
\textstyle \mu_\alpha = \half(\mu,\alpha^\vee) \in\half\mathbb{Z}, \qquad \lambda_\alpha = \half(\lambda,\alpha^\vee) \in\lie{h}^*/\half i \hbar_\alpha^{-1}\mathbb{Z},
\]
as in the previous section.

Recall that for any $w\in W$ we define
	\begin{align*}
	 S(w) &= \{\alpha\in{\bf\Delta}^+ \mid w\alpha \in {\bf\Delta}^- \}
	   = {\bf\Delta}^+ \cap w^{-1}{\bf\Delta}^-,
	\end{align*}
and note that $l(w) = |S(w)|$.

\begin{prop}
	\label{prop:intvertible_intertwiners}
	Let $(\mu,\lambda)\in\weights\times\lie{h}_q^*$ with associated parameters $(l,r)\in\lie{h}_q^* \times \lie{h}_q^*$ such that $\mu=l-r$, $\lambda+2\rho = -l-r$ as usual.  Let $w\in W$ and suppose that for all $\alpha\in S(w)$ we have
	\[
	 \pm \lambda_\alpha \notin |\mu_\alpha|+\mathbb{N} \bmod \textstyle \half i \hbar_\alpha^{-1} \mathbb{Z}.
	\]
	Then $\Gamma(\E_{\mu,\lambda}) \cong \Gamma(\E_{w\mu,w\lambda})$
\end{prop}

\begin{proof}
	We proceed by induction on the length of $w$.  If $w=1$, the result is obvious.  So suppose that $l(w)\geq1$, and write $w=s_\alpha v$ for some simple root $\alpha$ and some $v\in W$ with $l(w) = l(v)+1$.  In this case, we have $S(w) = S(v) \cup \{v^{-1}\alpha\}$.
	In particular, by the inductive hypothesis we have $\Gamma(\E_{\mu,\lambda}) \cong \Gamma(\E_{v\mu,v\lambda})$.
	
  Further, since $v^{-1}\alpha \in S(w)$, the hypothesis gives that
	\[ \textstyle
	 \pm \half (\lambda, {v}^{-1}\alpha^\vee) 
	  \notin |\half(\mu,{v}^{-1}\alpha^\vee)| + \mathbb{N} 
	   \bmod  \half i \hbar_\alpha^{-1}\mathbb{Z},
	\]
	and hence $\pm (v\lambda)_\alpha \notin |(v\mu)|_\alpha+\mathbb{N} \bmod \half i \hbar_\alpha^{-1}\mathbb{Z}$.  Therefore, by Theorem \ref{thm:simple_intertwiners}, we have a bijective intertwiner from  $\Gamma(\E_{v\mu,v\lambda})$ to $\Gamma(\E_{w\mu, w\lambda})$.  The result follows.
\end{proof}

The next two results are quantum analogues of Theorems I.4.3 and I.4.2 in \cite{Duflocomplexlnm}.

\begin{theorem}
	\label{thm:minimal_Kq-type_in_submodules}
	Let $(\mu,\lambda)\in\weights\times\lie{h}_q^*$. Every nonzero submodule of $\Gamma(\E_{\mu,\lambda})$ contains the minimal $K_q$-type if and only if  $q_\alpha^{(\lambda,\alpha^\vee)} \notin  q_\alpha^{|(\mu,\alpha^\vee)|+2\mathbb{N}}$ for every positive root $\alpha\in{\bf\Delta}^+$.
  In this case, the submodule generated by the minimal $K_q$-type is equal to $V_{\mu,\lambda}$, and this is the only simple submodule of $\Gamma(\E_{\mu,\lambda})$.
\end{theorem}

\begin{proof}	
	We fix a pair $(l,r)\in\lie{h}^*_q \times \lie{h}^*_q$ with $\mu=l-r$ and $\lambda+2\rho=-l-r$ as usual. 
	Suppose that $q_\alpha^{(\lambda,\alpha^\vee)} \notin  q_\alpha^{|(\mu,\alpha^\vee)|+2\mathbb{N}}$ for every positive root $\alpha\in{\bf\Delta}^+$.  By Lemma \ref{lem:integer_point}, this is equivalent to saying that for every $\alpha\in{\bf\Delta}^+$, the quantities $q_\alpha^{(l+\rho,\alpha^\vee)}$ and $q_\alpha^{(r+\rho,\alpha^\vee)}$ do not both lie in $\pm q_\alpha^{-\mathbb{N}}$.
	
	Recall that we write $\Re(l)$ for the component of $l$ in the $\mathbb{R}$-span of $\bf\Delta$, see the remarks before Proposition \ref{prop:root_subsystem}.
  To begin with, let us suppose that the real part of $l+\rho$ lies in the closure of the dominant Weyl chamber, meaning that
	\[
	 (\Re(l)+\rho,\alpha^\vee) \geq 0
	\]
	for all $\alpha\in{\bf\Delta}^+$. 
	This means in particular that $l$ is dominant.  According to Lemma \ref{dominantvermahcproperties} $b)$,
	there exists $\hat{u}\in\hat{W}$ such that, if we put $r'=\hat{u}.r$, we have $\Gamma(\E_{\mu,\lambda})\cong\F_l(M(r')^\vee)$ and there is no $\alpha\in{\bf\Delta}^+$ with both $q_\alpha^{(l+\rho,\alpha^\vee)} = \pm1$ and $q_\alpha^{(r'+\rho, \alpha^\vee)} \in \pm q_\alpha^{\mathbb{N}}$.
	
	Let $H$ be a submodule of $\F_l(M(r')^\vee)$.  By Proposition \ref{fltladjointduality}, we have a nonzero $U_q(\lie{g})$-linear map $\varphi: \T_l(H) \to M(r')^\vee$.  Every nontrivial quotient module of $M(r')$ projects to $V(r')$, and consequently every nontrivial submodule of $M(r')^\vee$ contains $V(r')$.  In particular, $V(r') \subset \im(\varphi)$.
	
	Using Lemma \ref{dominantvermahcproperties} $b)$ and applying the exact functor $\F_l$, we obtain an inclusion
	\[
	 V_{\mu,\lambda} \cong \F_l(V(r')) \subset \F_l(\im(\varphi)) \cong \im(\F_l(\varphi)).
	\]
	By Proposition \ref{unitisomorphism},
	the inclusion of $H$ into $\F_l(M(r')^\vee)$ factorizes as
	\[
	 H \cong \F_l(T_l(H)) \stackrel{\F_l(\varphi)}{\longrightarrow} \F_l(M(r')^\vee),
	\]
	so that $\im(\F_l(\varphi)) = H$ and so the inclusion above shows that $V_{\mu,\lambda}$ is isomorphic to a submodule of $H$.  By considering minimal $K_q$-types, it follows that $V_{\mu,\lambda}$ is indeed a submodule of $H$ as claimed.

	Now consider the case of arbitrary $l$.  We can find $w\in W$ such that the real part of $w(l+\rho)$ belongs to the closure of the dominant Weyl chamber.  For any $\alpha\in S(w)$ we have $w\alpha \in {\bf\Delta}^-$, so that
	\[
	 (\Re(l)+\rho,\alpha^\vee) = (w(\Re(l)+\rho),w\alpha^\vee) \leq 0.
	\]  
	In particular, 
	according to Lemma \ref{lem:integer_point}, 
	we do not have $q_\alpha^{(\lambda,\alpha^\vee)} \in q_\alpha^{-|(\mu,\alpha^\vee)|-2\mathbb{N}}$, and by hypothesis we do not have $q_\alpha^{(\lambda,\alpha^\vee)} \in q_\alpha^{|(\mu,\alpha^\vee)|+2\mathbb{N}}$.
	Therefore, Proposition \ref{prop:intvertible_intertwiners} shows that $\Gamma(\E_{\mu,\lambda}) \cong \Gamma(\E_{w\mu,w\lambda})$.  By the previous case, every nonzero submodule of $\Gamma(\E_{w\mu,w\lambda})$ contains $V_{w\mu,w\lambda}$, and the result follows.
	
	In this situation, the intersection of all nonzero submodules of $\Gamma(\E_{\mu,\lambda})$ contains the minimal $K_q$-type, and in particular is nonzero.  It is  necessarily simple, and so must be equal to $V_{\mu,\lambda}$.  Moreover, there can be no other simple submodule, since otherwise the above intersection would be zero.
	
	Conversely, suppose that there exists $\alpha\in{\bf\Delta}^+$ with $q_\alpha^{(\lambda,\alpha^\vee)} \in  q_{\alpha}^{|(\mu,\alpha^\vee)|+2\mathbb{N}}$, and so $q_\alpha^{(l+\rho,\alpha^\vee)} \in \pm q_\alpha^{-\mathbb{N}}$ and $q_\alpha^{(r+\rho,\alpha^\vee)}\in \pm q_\alpha^{-\mathbb{N}}$, by Lemma \ref{lem:integer_point}.  
  In this case, $r$ is strongly linked to $r'=s_{k,\alpha}.r$ for some $k\in\mathbb{Z}_2$, see Definition \ref{def:strongly_linked} and the remark which follows it.
  By Theorem \ref{thm:Zelobenko_explicit} we have an nonzero intertwiner $T:\Gamma(\E_{\mu',\lambda'}) \to \Gamma(\E_{\mu,\lambda})$ where $\mu'=l-r'$ and $\lambda'+2\rho = -l-r'$.  
	
	We have $(\Re(l)+\rho,\alpha^\vee) \in -\mathbb{N}$ and $(\Re(r)+\rho,\alpha^\vee) \in -\mathbb{N}$, so
	$(s_\alpha.\Re(l)+\rho,\alpha^\vee) \in \mathbb{N}$ and $(s_\alpha.\Re(r)+\rho,\alpha^\vee) \in \mathbb{N}$.  Thus, neither $\Re(l)$ nor $\Re(r)$ is fixed by the shifted action of $s_\alpha$, and an 
  application of Lemma \ref{dominantweightlemma} shows that $\|\mu'\| > \|\mu\|$.  Therefore, $\im(T)$, is a submodule of $\Gamma(\E_{\mu,\lambda})$ which does not contain the minimal $K_q$-type.  
	This completes the proof.  
\end{proof}

\begin{cor}
  \label{cor:minimal_Kq-type_generating}
  Let $(\mu,\lambda)\in\weights\times\lie{h}_q^*$.  The principal series representation $\Gamma(\E_{\mu,\lambda})$ is generated by its subspace of minimal $K_q$-type if and only if $q_\alpha^{(\lambda,\alpha^\vee)} \notin  q_\alpha^{-|(\mu,\alpha^\vee)|-2\mathbb{N}}$ for every positive root $\alpha\in{\bf\Delta}^+$.
\end{cor}

\begin{proof}
	First suppose that $q_\alpha^{(\lambda,\alpha^\vee)} \notin   q_\alpha^{-|(\mu,\alpha^\vee)|-2\mathbb{N}}$ for every positive root $\alpha\in{\bf\Delta}^+$.
	Let $V\subset \Gamma(\E_{\mu,\lambda})$ be the submodule generated by the minimal $K_q$-type.  Using the invariant pairing from Lemma \ref{lem:principal_series_bilinear_pairing}, we can define a submodule of $\Gamma(\E_{-\mu,-\lambda})$ by
	\[
    V^\perp = \{ \eta \in \Gamma(\E_{-\mu,-\lambda}) \mid 
     \phi(\xi\eta) = 0 \text{ for all } \xi \in V \} .
	\]
	The pairing decomposes as a sum of nondegenerate pairings between the subspaces of each $K_q$-type.  Since $V$ contains the minimal $K_q$-type,  $V^\perp$ does not contain the minimal $K_q$-type of $\Gamma(\E_{-\mu,-\lambda})$.  Note, however, that $\Gamma(\E_{-\mu,-\lambda})$ fulfils the conditions of  Theorem \ref{thm:minimal_Kq-type_in_submodules}, so we must have $V^\perp=\{0\}$, and hence $V = \Gamma(\E_{\mu,\lambda})$.
	
	Conversely, suppose that $\Gamma(\E_{\mu,\lambda})$ is generated by its subspace of minimal $K_q$-type.  Let $U\subset \Gamma(\E_{-\mu,-\lambda})$ be a submodule which does not contain the minimal $K_q$-type.  Then the module 
	\[
   U^\perp = \{ \xi \in \Gamma(\E_{\mu,\lambda}) \mid 
    \phi(\xi\eta) = 0 \text{ for all } \eta \in \Gamma(\E_{-\mu,-\lambda}) \} \subset \Gamma(\E_{\mu,\lambda})
  \]	
	contains the minimal $K_q$-type, and hence $U^\perp = \Gamma(\E_{\mu,\lambda})$ by hypothesis.  Therefore $U=0$, and we deduce that every nonzero submodule of $\Gamma(\E_{-\mu,-\lambda})$ contains the minimal $K_q$-type.  Theorem \ref{thm:minimal_Kq-type_in_submodules} implies that $q_\alpha^{(-\lambda,\alpha^\vee)} \notin  q_\alpha^{|(-\mu,\alpha^\vee)|+2\mathbb{N}}$ for every positive root $\alpha\in{\bf\Delta}^+$, and the result follows.
\end{proof}

\subsection{Unitary representations}
\label{sec:unitary}

In this section we briefly comment on the question of unitarizability of Harish-Chandra modules. For more information 
we refer to the work of Arano \cite{Aranospherical}, \cite{Aranocomparison}.

Let us write $\lie{a}^* \subset \lie{h}^*$
\nomenclature[o$a^*$]{$\lie{a}^*$}{$\mathbb{R}$-span of the roots $\bf\Delta$}%
for the $\mathbb{R}$-span of the roots $\bf\Delta$, and $\lie{t}_q^* = i\lie{a}^* /  i\hbar^{-1}\roots^\vee$.
\nomenclature[o$t_q^*$]{$\lie{t}_q^*$}{$=i\lie{a}^* /  i\hbar^{-1}\roots^\vee$}%
Recall that for $\lambda\in\lie{h}_q^*$ we have a decomposition $\lambda = \Re(\lambda) + i \Im(\lambda)$ where $\Re(\lambda)\in\lie{a}^*$ and $i\Im(\lambda) \in \lie{t}_q^*$.  Accordingly we write $\overline{\lambda} = \Re(\lambda) - i\Im(\lambda)$. 
\nomenclature[s]{$\overline{\lambda}$}{conjugate of $\lambda\in\lie{h}_q^*$}%
Note that the characters $K_\lambda \in \M(\DF(K_q))$ with $\lambda\in\lie{h}_q^*$, as defined in Subsection \ref{sec:characters_of_UqRb}, satisfy $K_\lambda^* = K_{\overline{\lambda}}$.

\begin{prop} \label{prop:sesquilinear_pairing}
Let $ (\mu, \lambda) \in \weights \times \mathfrak{h}_q^* $. The standard inner product on $\CF^\infty(K_q)$ defined by
$
 \bra \xi, \eta \ket = \phi(\xi^*\eta),
$
where $\phi$ is the Haar functional, restricts to a non-degenerate sesquilinear pairing 
between $ \Gamma(\E_{\mu,-\overline{\lambda}}) $ and $ \Gamma(\E_{\mu, \lambda}) $ which is $G_q$-invariant in the sense that
\begin{equation*} \label{eq:invariant_inner_product}
\bra \xi, \pi_{\mu,\lambda}(x \bowtie f)\eta \ket = \bra \pi_{\mu,-\overline{\lambda}}((x \bowtie f)^*)\xi, \eta \ket
\end{equation*}
for all $ x \bowtie f \in \DF(G_q), \xi \in \Gamma(\E_{\mu,-\overline{\lambda}}), \eta \in \Gamma(\E_{\mu,\lambda}) $.
\end{prop}

\begin{proof} 
Consider the bijective conjugate-linear map $*: \Gamma(\E_{\mu,-\overline{\lambda}}) \to \Gamma(\E_{-\mu,-\lambda})$ which sends $\xi$ to $\xi^*$.  We claim that for any $x\bowtie f\in \DF(G_q)$, we have
\[
 * \circ \pi_{\mu,-\overline{\lambda}}((x\bowtie f)^*)
 =
 \pi_{-\mu,-\lambda}(\hat{S}_{G_q}^{-1}(x\bowtie f)) \circ *.
\]
It suffices to prove this for elements of the form $x\bowtie 1$ and $1\bowtie f$ with $x\in\DF(K_q)$ and $f\in\CF^\infty(K_q)$.  Let $\xi \in \Gamma(\E_{-\mu,-\lambda})$.  Using the compatibility of the $*$-structures, as recorded after Definition \ref{def:CKq}, we calculate
\begin{align*}
   \pi_{-\mu,-\lambda}(\hat{S}_{G_q}^{-1}(x\bowtie 1)) (\xi^*)
   &= (x,\xi_{(1)}^*)\, \xi_{(2)}^* \\
   &= \overline {(\hat{S}(x^*),\xi_{(1)})}\, \xi_{(2)}^* \\   
   &= (\pi_{\mu,-\overline{\lambda}}(x^*\bowtie 1) \xi)^*,
\end{align*}
and similarly,
\begin{align*}
  \pi_{-\mu,-\lambda}(\hat{S}_{G_q}^{-1}(1\bowtie f)) (\xi^*)
   &= S^{-1}(f_{(3)}) \xi^* f_{(1)} \, (K_{2\rho-\lambda}, S^{-1}(f_{(2)})) \\
   &= S^{-1}(f_{(3)}) \xi^* f_{(1)} \, (\hat{S}(K_{2\rho-\overline{\lambda}}^*), f_{(2)}) \\
   &= (f_{(1)}^* \xi S(f_{(3)}^*))^* \, \overline{(K_{2\rho-\overline{\lambda}}, f_{(2)}^*) } \\
   &= (\pi_{\mu,-\overline{\lambda}}(1\bowtie f^*) \xi)^*,
\end{align*}
which proves the claim.

Therefore, using the invariance property of the non-degenerate bilinear pairing
 $(\ ,\ ) : \Gamma(\E_{-\mu,-\lambda}) \times \Gamma(\E_{\mu,\lambda}) \to \mathbb{C}$ from Lemma \ref{lem:principal_series_bilinear_pairing} we have
\begin{align*}
  \bra \xi, \pi_{\mu,\lambda}(x \bowtie f)\eta \ket 
   &= ( \xi^*, \pi_{\mu,\lambda}(x \bowtie f)\eta) \\
   &= ( \pi_{-\mu,-\lambda}(\hat{S}_{G_q}^{-1}(x \bowtie f))(\xi^*), \eta) \\
   &= ( (\pi_{\mu,-\overline{\lambda}}((x \bowtie f)^*)\xi)^*, \eta ) \\
   &= \bra \pi_{\mu,-\overline{\lambda}}((x \bowtie f)^*)\xi, \eta \ket,
\end{align*}
for any $\xi\in\Gamma(\E_{\mu,-\overline{\lambda}})$ and $\eta\in\Gamma(\E_{\mu,\lambda})$ and $x\bowtie f \in \DF(G_q)$, as claimed.
\end{proof}

A nondegenerate $ \DF(G_q)$-module $ V $ is \emph{unitarizable} if it admits a positive definite Hermitian form $ (\;,\;) $ which 
is invariant in the sense that 
$$
(x \cdot v, w) = (v, x^* \cdot w) 
$$
for all $ v, w \in V $ and $ x \in \DF(G_q) $. 
By Proposition \ref{prop:sesquilinear_pairing}, the invariant sesquilinear forms $ (\;,\;)  $ 
on $ \Gamma(\E_{\mu,\lambda}) $ are in one-to-one correspondence with the intertwining 
operators $ T: \Gamma(\E_{\mu,\lambda}) \rightarrow \Gamma(\E_{\mu,-\overline{\lambda}}) $ via the formula
$$
(\eta, \xi) = \bra T(\eta), \xi \ket, \qquad \xi, \eta \in \Gamma(\E_{\mu,\lambda}).
$$ 
Moreover, if $ \Gamma(\E_{\mu,\lambda}) $ admits a nonzero invariant sesquilinear form $ (\;,\;) $ then we can arrange it to be Hermitian.  
Indeed, by polarization, we must have $ (\xi, \xi) \neq 0 $ for some $ \xi $, and after multiplying the form by some scalar we 
can ensure this is strictly positive. Then the Hermitian form defined by
$$
\bra \eta, \xi \ket = \frac{1}{2} \left((\eta, \xi) + \overline{(\xi, \eta)} \right)
$$
is invariant and nonzero since $ \bra \xi,\xi \ket > 0 $. 

On the other hand, there is no guarantee that this Hermitian form will be positive definite. Two particular classes of unitarizable principal series 
are well-known.  As before, we use the notation $\lie{a}^*\subset \lie{h}^*$ for the $\mathbb{R}$-span of the roots $\bf\Delta$ and put $\lie{t}^*_q = i\lie{a}^* / i\hbar^{-1}\roots^\vee \subset \lie{h}_q^*$.
Thus, $\lie{t}_q^*$ is a compact torus of dimension $N$.

\begin{theorem}
Let $ G_q $ be a complex semisimple quantum group. 
\begin{bnum}
\item[a)] (Unitary principal series)
The principal series representation $ \Gamma(\E_{\mu, \lambda}) $ equipped with the standard inner product is unitary if and only if $ \lambda \in \mathfrak{t}_q^*$.
\item[b)] (Complementary series)
Let $ \alpha $ be a simple root. There is an invariant inner product on the principal series representations $ \Gamma(\E_{\mu, \lambda' + t \alpha}) $ 
where $ \mu \in \weights, \lambda' \in \mathfrak{t}_q^* $ are both $ s_\alpha $-fixed and $ -1 < t < 1 $.
\end{bnum}
\noindent All of these representations are irreducible, and they are unitarily equivalent if and only if their parameters lie in the same orbit of the Weyl 
group action on $ \weights \times \mathfrak{h}_q^* $.
\end{theorem}
\begin{proof} We have the following facts. 
\begin{bnum}
\item[a)] The standard inner product is $ G_q $-invariant if and only if $ (\mu, \lambda) = (\mu, -\overline{\lambda}) $, which 
is equivalent to $ \lambda \in \mathfrak{t}_q^*$.
\item[b)] Let $ \mu, \lambda = \lambda' + t \alpha $ be as stated. Note that 
$$
-\overline{\lambda} = \lambda' - t\alpha = s_\alpha(\lambda),
$$
so the standard intertwiner $ \Gamma(\E_{\mu,\lambda}) \rightarrow \Gamma(\E_{s_\alpha \mu, s_\alpha \lambda}) = \Gamma(\E_{\mu,-\overline{\lambda}}) $ from Theorem \ref{thm:simple_intertwiners} combined with the sesquilinear pairing from Proposition \ref{prop:sesquilinear_pairing} yields an invariant Hermitian form 
on $ \Gamma(\E_{\mu,\lambda}) $. Moreover, since $ \lambda_\alpha = t $ and $\mu_\alpha=0$, this intertwiner is bijective for all $ -1 < t < 1 $. 
At $ t = 0 $, the Hermitian form is positive definite by part $a)$. Therefore, by a standard continuity argument, they are positive definite for all $ -1 < t < 1 $.  
\end{bnum}
The irreducibility of these representations follows from Theorem \ref{principalirreducible}. 
The intertwiners associated to the simple reflections are therefore bijective and automatically unitary. 
\end{proof} 

This is far from an exhaustive list of irreducible unitary representations of $ G_q $. In particular, some subquotients of generalized principal series 
are unitarizable. Arano \cite{Aranocomparison} has used continuity arguments in $ q $ to compare the unitary dual 
of $ G_q $ with the classical unitary dual of $ G $.  This yields a complete classification for $ SL_q(n,\mathbb{C}) $, and an almost complete classification for general $G_q$. 

The case of $ SL_q(2,\mathbb{C}) $ was already completed by Pusz \cite{Puszunitaryrepresentations}; see also \cite{PWquantumlorentzgelfand}. 
We state their result without proof.

\begin{theorem}
Up to unitary equivalence, the irreducible unitary representations of $ SL_q(2,\mathbb{C}) $ are the Hilbert space completions of the 
following Harish-Chandra modules. 
\begin{bnum}
\item[a)] The unitary principal series $ \Gamma(\E_{\mu,\lambda}) $ with $ \mu \in \frac{1}{2} \mathbb{Z}, \lambda \in i \mathbb{R}/i \hbar^{-1} \mathbb{Z} $, 
modulo the unitary equivalences $ \Gamma(\E_{\mu,\lambda}) \cong \Gamma(\E_{-\mu,-\lambda}) $,
\item[b)] The two complementary series $ \Gamma(\E_{0,t}) $ and $ \Gamma(\E_{\frac{1}{2} i \hbar^{-1},t}) $ with $ -1 < t < 1 $ equipped with their unitary inner products, modulo the unitary 
equivalences $ \Gamma(\E_{0,t}) \cong \Gamma(\E_{0,-t}) $ and $ \Gamma(\E_{\frac{1}{2} i \hbar^{-1},t}) \cong \Gamma(\E_{\frac{1}{2} i \hbar^{-1},-t}) $,
\item[c)] The trivial representation and the unitary character $ \chi_z $ where $ z $ is the nontrivial element of the centre of $ SU(2) $.  
\end{bnum}
\end{theorem}

As a final remark, we recall that Corollary \ref{cor:simple_intertwiners_explicit} gives explicit formulas for the intertwiners between unitary principal series representations.  In particular, we can obtain a particularly simple formula for the intertwiners between the base of principal series representations. This fact was needed in \cite{VYfredholm}.

\begin{cor}
For any $ \mu \in \weights $ and any simple root $ \alpha $, the operator 
$$
\Ph(F_\alpha)^{2 \mu_\alpha}: \Gamma(\E_{\mu, 0}) \rightarrow \Gamma(\E_{s_\alpha \mu, 0}),
$$
defined by the right regular action, is a unitary intertwiner. 
\end{cor}

We recall that before Corollary \ref{cor:simple_intertwiners_explicit} we introduced the notation $\Ph(F_\alpha)^k \in \M(\DF(K_q))$ for the operator phase of $F_\alpha^{k}$  if $k\geq0$ or $E_\alpha^{-k}$ if $k<0$.

\newpage




\nomenclature[a,1]{$\NN$}{$=\{1,2,\ldots\}$ \nomnorefpage}
\nomenclature[a,2]{$\NN_0$}{$=\{0,1,2,\ldots\}$ \nomnorefpage}
\nomenclature[a,3]{$\KH$}{base field \nomnorefpage} 
\nomenclature[a,4]{$\sigma$, $\Sigma$}{flip maps on tensor products, $a\otimes b \mapsto b\otimes a$ \nomnorefpage}



\nomenclature[g01]{$\hit$}{
	left action of $ U_q(\mathfrak{g}) $ on $ \Poly(G_q) $ given by $ X \hit f = f_{(1)} (X, f_{(2)}) $, p.\pageref{nom:hit_Uqg} \\
	--- of $\DF(K)$ on $\CF^\infty_c(K)$ given by $x \hit h = h_{(1)} (x, h_{(2)}) $, p.\pageref{nom:hit_DK}
	\nomnorefpage 
}

\nomenclature[g02]{$\hitby$}{
	right action of $ U_q(\mathfrak{g}) $ on $ \Poly(G_q) $ given by $ f \hitby X = (X, f_{(1)}) f_{(2)} $, p.\pageref{nom:hit_Uqg} \\
	--- of $\DF(K)$ on $\CF^\infty_c(K)$ given by $h \hitby x = (x, h_{(1)}) h_{(2)} $, p.\pageref{nom:hit_DK}
	\nomnorefpage
}

\nomenclature[g03]{$\rightarrow$}{
	left adjoint action of a (regular multiplier) Hopf algebra, p.\pageref{nom:adjoint_action1} \\
	--- in particular, for $\tilde{U}_q(\lie{g})$ or $U_q(\lie{g})$, p.\pageref{nom:left_adjoint_action1}, p.\pageref{nom:left_adjoint_action2} \\
	--- coadjoint action of $U_q(\lie{g})$ on $U_q(\lie{g})^*$ or $\Poly(G_q)$, p.\pageref{nom:coadjoint_action}
	\nomnorefpage
}

\nomenclature[g041]{$\overset{\gamma}{\rightarrow}$}{
	$\gamma$-twisted adjoint action of $U_q(\lie{g})$, p.\pageref{nom:gamma-twisted}
}%

\nomenclature[g04]{}{
	--- any compatible ``adjoint'' action of $U_q(\lie{g})$ on an $FU_q(\lie{g})$-bimodule, 
	p.\pageref{defadcompatibility}
	--- adjoint action of $U_q(\lie{g})$ on $\Hom(M,N)$, 
	p.\pageref{nom:adjoint_action_on_Hom} \\
	--- adjoint action of $U_q(\lie{g})$ on $(M\otimes N)^*$, 
	p.\pageref{nom:adjoint_action_on_MxN}, p.\pageref{nom:adjoint_action_on_MxN2} \\
	--- adjoint action of $U_q(\lie{g})\cong U_q^\mathbb{R}(\lie{k})$ on a $U_q^\mathbb{R}(\lie{g})$-module,
	p.\pageref{nom:FUq-actions} \\
	--- two versions of adjoint action of $U_q(\lie{g})$ on $F\Hom(M(l),M(r)^\vee)$, 
	p.\pageref{nom:FHom_bimodule}
	\nomnorefpage
}

\nomenclature[g05]{$\leftarrow$}{
	right adjoint action of $\tilde{U}_q(\lie{g})$ or $U_q(\lie{g})$ on itself, p.\pageref{nom:right_adjoint_action1}
	\nomnorefpage
}






\nomenclature[g12]{$w.\lambda$}{
	shifted action of the Weyl group $W$ on $\lie{h}_q^*$ given by $w.\lambda = w(\lambda+\rho)-\rho$ \\
  --- similarly for the extended Weyl group p.\pageref{nom:shifted_extended_Weyl_group_action} 
	\nomnorefpage
}


\nomenclature[o$(\ ,\ )$1]{$(\ ,\ )$}{
	canonical skew-pairing $\hat{H}\times H \to \KH$; $(x,f)=x(f)$, p.\pageref{nom:pairing}\\
	--- reverse skew-pairing $H\times\hat{H}\to \KH$; $(f,x) = (x,S^{-1}(f)) = (\hat{S}(x),f)$, p.\pageref{def:reverse_pairing}
	\nomnorefpage
}

\nomenclature[o$(\ ,\ )$2]{}{
	--- skew-pairing $U_q(\lie{g}) \times \O(G_q) \to \KH$, \pageref{nom:OGq-pairing} \\
	--- skew-pairing $U_q^\mathbb{R}(\lie{k}) \times \CF^\infty(K_q) \to \CH$, p.\pageref{nom:K_q-pairing} \\
	--- reverse skew-pairing $\CF^\infty(K_q) \times U_q^\mathbb{R}(\lie{k}) \to \CH$,
	p.\pageref{nom:K_q-reverse-pairing} \\
	--- reverse skew-pairing $\CF^\infty(K_q) \times \DF(K_q) \to \CH$, p.\pageref{nom:G_q-pairing} \\
	--- skew-pairing $\DF(G_q) \times \CF^\infty_c(G_q) \to \mathbb{C}$, p.\pageref{nom:G_q-pairing} \\
	--- pairing $(U_q(\lie{g})\times U_q(\lie{g})) \times (\Poly(G_q)\times\Poly(G_q)) \to \CH$, p.\pageref{nom:UU-OO-pairing} 
	\nomnorefpage
}

\nomenclature[o$(\ ,\ )$3]{}{
	--- bilinear pairing $M(\mu,2\rho+\lambda) \times \ind_{B_q}^{G_q}(\CH_{\mu,\lambda}) \to \CH$, p.\pageref{lem:Verma-principal-pairing}   
	\nomnorefpage
}

\nomenclature[o$(\ ,\ )$4]{}{
	--- bilinear pairing  $\Gamma(\E_{\mu,\lambda}) \times \Gamma(\E_{-\mu,-\lambda}) \to \mathbb{C}$,
	p.\pageref{lem:principal_series_bilinear_pairing}
	\nomnorefpage
}

\nomenclature[o$(\ ,\ )$5]{}{
	--- rescaling of Killing form on $\lie{g}$ with $(\alpha,\alpha)=2$ for shortest root $\alpha$, p.\pageref{sec:semisimple_Lie_algebras} 
	\nomnorefpage
}


\nomenclature[o$(,)$4]{
	$\bra\ ,\ \ket$}{scalar product on $L^2(G)$ for an algebraic quantum group, defined via $\bra f, g \ket = \phi(f^* g)$, p.\pageref{nom:GNS-ip} \\
	--- invariant scalar product on a simple $U_q^\mathbb{R}(\lie{k})$-module, p.\pageref{nom:V-ip} \\
	--- sesquilinear pairing of principal series representations $\Gamma(\E_{\mu,-\overline{\lambda}}) \times \Gamma(\E_{\mu, \lambda}) \to \CH$, p.\pageref{prop:sesquilinear_pairing}
	\nomnorefpage
}%


\nomenclature[o$?$]{$*$}{convolution of linear functionals on a Hopf algebra, p.\pageref{nom:convolution}
	\nomnorefpage}

\nomenclature[o$?$]{$\bowtie$}{
	Drinfeld double of multiplier Hopf algebras, p.\pageref{def:Drinfeld_double} \\
	--- Drinfeld double of algebraic quantum groups, defined via: \\
	\hspace*{1cm} the group algebra $\DF(K\bowtie L)$, p.\pageref{nom:alg_double} \\
	\hspace*{1cm} the algebra of functions $\CF_c^\infty(K\bowtie L)$, p.\pageref{nom:codouble}
	\nomnorefpage
}

\nomenclature[o$q$1]{$q$}{parameter $q\in\mathbb{K}^\times$, usually with $q=s^L$,  not a root of unity, p.\pageref{secqnumbers}, \pageref{def:no_Serre}, \pageref{defuqg}\\
 in Chapters \ref{chcomplexqg}, \ref{chcato}, \ref{chreptheory}, $q=e^h\in\RR^\times_+$, $q\neq1$
  \nomnorefpage }

\nomenclature[o$q$2]{$q_i$}{$=q^{d_i} = q^{(\alpha_i,\alpha_i)/2}$,
  p.\pageref{def:no_Serre}
  \nomnorefpage}

\nomenclature[o$q$3]{$q_\beta$}{$=q^{(\beta,\beta)/2}$ for any root $\beta\in\bf\Delta$,
  p.\pageref{nom:q_beta_r}, p.\pageref{nom:q_beta} \nomnorefpage}

\nomenclature[o$\Delta$1]{$\Delta$}{
	coproduct of a coalgebra, p.\pageref{nom:coalgebra_coproduct}\\
		 --- of a multiplier Hopf algebra, p.\pageref{nom:multiplier_coproduct}\\
		 --- of a Hopf $C^*$-algebra, p.\pageref{defhopfcstar} \\
		 --- of a locally compact quantum group, p.\pageref{defqg}
		 \nomnorefpage
 }
\nomenclature[o$\Delta$2]{$\hat{\Delta}$}{
	coproduct of a dual multiplier Hopf algebra, p.\pageref{dualmulthopfconstruct} \\
	--- of $U_q(\lie{g})$,  p.\pageref{lemuqhopf}, p.\pageref{nom:Uqg_Hopf} \\
	--- of the group von Neumann algebra of a locally compact quantum group, p.\pageref{nom:dual_vN-coproduct}
	\nomnorefpage
}

\nomenclature[o$S$1]{$S$}{
	antipode of a multiplier Hopf algebra, p.\pageref{nom:S} \\
	--- of a locally compact quantum group, p.\pageref{nom:vN-antipode}
	\nomnorefpage
}

\nomenclature[o$S$2]{$\hat{S}$}{
	antipode of a dual multiplier Hopf algebra, p.\pageref{dualmulthopfconstruct}\\
	--- of $U_q(\lie{g})$, p.\pageref{lemuqhopf}, p.\pageref{nom:Uqg_Hopf}
	\nomnorefpage
}

\nomenclature[o$\epsilon$2]{$\hat{\epsilon}$}{
	counit of a dual multiplier Hopf algebra, p.\pageref{dualmulthopfconstruct},\\
	--- of $U_q(\lie{g})$, p.\pageref{lemuqhopf}, p.\pageref{nom:Uqg_Hopf}
	\nomnorefpage
}

\nomenclature[o$(K_i)$2]{$[K_i;l]_{q_i}$}{
	$=\displaystyle \frac{q_i^lK_i - q_i^{-l}K_i^{-1}}{q_i - q_i^{-1}}$,
	p.\pageref{nom:K:0}, p.\pageref{nom:K:l}
  \nomnorefpage
}

\nomenclature[o$\phi$1]{$\phi$}{
	left invariant integral, p.\pageref{nom:left_integral} \\
	   --- left invariant integral of $\Poly(G_q)\cong\CF^\infty(K_q)$, p.\pageref{nom:OGq_integral} \\
     --- left Haar weight of a locally compact quantum group, p.\pageref{nom:Haar_weights}\\
     --- left invariant integral of an algebraic quantum group, p.\pageref{defaqg} \\
     --- (left) invariant integral of a compact quantum group, p.\pageref{nom:Haar-compact}
     \nomnorefpage
}

\nomenclature[o$\psi$1]{$\psi$}{
	right invariant integral, p.\pageref{nom:right_integral}, \\
	--- right Haar weight of a locally compact quantum group, p.\pageref{nom:Haar_weights}\\
	--- right invariant integral of an algebraic quantum group, p.\pageref{defaqg}\\
	--- (right) invariant integral of a compact quantum group, p.\pageref{nom:Haar-compact}
	\nomnorefpage
}

\nomenclature[o$\phi$2]{$\hat{\phi}$}{
	left integral of	a dual multiplier Hopf algebra, p.\pageref{nom:dual_left_integral},\\
	--- left Haar functional of $\DF(K_q)$, p.\pageref{nom:dual_Haar_weights} \\
	--- left Haar weight of the group von Neumann algebra of a locally compact quantum group, p.\pageref{nom:dual_vN_Haar_weight}	
	\nomnorefpage
}

\nomenclature[o$\psi$2]{$\hat{\psi}$}{
	right Haar functional of $\DF(K_q)$, p.\pageref{nom:dual_Haar_weights} 
	\nomnorefpage
}

\nomenclature[o$ad$]{$\ad$}{adjoint action of $\tilde{U}_q(\lie{g})$ or $U_q(\lie{g})$ on itself, 
	p.\pageref{nom:left_adjoint_action1}, p.\pageref{nom:left_adjoint_action2}
	\nomnorefpage
}

\nomenclature{$K_\lambda$}{
	Cartan generators in $U_q(\lie{g})$ with $\lambda\in\weights$, \pageref{def:no_Serre}, \pageref{defuqg} \\
	--- generalizations in $\M(\DF(K_q))$ with $\lambda\in\lie{h}_q^*$, p.\pageref{nom:K_lambda}, p.\pageref{nom:K-central}
  \nomnorefpage
}

\nomenclature{$K_i$}{
	Cartan element $K_{\alpha_i}\in U_q(\lie{g})$, \pageref{def:no_Serre}, \pageref{defuqg}
	\nomnorefpage
}

\nomenclature{$E_i$}{
	generator of $U_q(\lie{g})$, p.\pageref{def:no_Serre}, p.\pageref{defuqg}
	\nomnorefpage
}

\nomenclature{$F_i$}{
  generator of $U_q(\lie{g})$, \pageref{def:no_Serre}, \pageref{defuqg}
	\nomnorefpage
}

\nomenclature{$\beta$}{
	bar automorphism of $\mathbb{Q}(s)$, p.\pageref{nom:beta_field_automorphism} \\
	--- bar automorphism of the restricted integral form $U_q^\A(\mathfrak{g})$,
	p.\pageref{defbarinvolution}
  \nomnorefpage}%

\nomenclature{$\kappa$}{
  Rosso form, p.\pageref{def:Rosso_form} \\
  --- quantum Killing form, p.\pageref{defquantumkilling}
  \nomnorefpage
}

\nomenclature{$\omega$}{algebra automorphism, coalgebra antihomomorphism of $\tilde{U}_q(\mathfrak{g})$, p.\pageref{lemflippingEandF} \\
	--- of $U_q(\mathfrak{g})$, p.\pageref{defOmega}
	\nomnorefpage}

\nomenclature{$\Omega$}{algebra anti-automorphism of $U_q(\mathfrak{g})$, p.\pageref{defOmega}
	\nomnorefpage}

\nomenclature[o$\Omega$2]{$\Omega$}{
	quantum Casimir element, p.\pageref{nom:Casimir}
  \nomnorefpage
}

\nomenclature[o$Uqg$]{$\tilde{U}_q(\mathfrak{g})$, $\tilde{U}_q(\mathfrak{n_\pm})$, $\tilde{U}_q(\mathfrak{b_\pm})$}{analogues of ${U}_q(\mathfrak{g})$, ${U}_q(\mathfrak{n_\pm})$, ${U}_q(\mathfrak{b_\pm})$ without Serre relations,
 pp.\pageref{def:no_Serre}--\pageref{nom:Uqb_tilde}
 \nomnorefpage}

\nomenclature{$s$}{element in $\mathbb{K}^\times$ such that $q=s^L$, p.\pageref{def:no_Serre}, p.\pageref{defuqg}, p.\pageref{def:integral_form}
	\nomnorefpage}

\nomenclature[o$h_q^*$]{$\lie{h}_q^*$}{character group of $U_q(\lie{h})$ written additively, p.\pageref{def:hq-star} \\
	--- identified with $\lie{h}^*/i\hbar^{-1}\roots^\vee$ when the ground field is $\mathbb{C}$, p.\pageref{nom:hq-star-C} 
 \nomnorefpage}

\nomenclature[o$<$]{$\leq$}{relation on $ \mathfrak{h}^*_q $ defined by $	\lambda \leq \mu$ if $\mu-\lambda \in \roots^+$, 
  p.\pageref{def:hq-star-order}, p.\pageref{sec:O_Artinian} \\
  --- and on the quotient $\mathfrak{h}^*_q / \half i\hbar^{-1}\roots^\vee$, p.\pageref{nom:leq_in_hqstar}
  \nomnorefpage}

\nomenclature{$V(\lambda)$}{simple weight module of highest weight $\lambda\in\lie{h}_q^*$, p.\pageref{nom:simple_quotient} \\
  --- in particular, the simple integrable module of highest weight $\lambda\in\weights^+$, p.\pageref{sl2irreps}, p.\pageref{thmfdirreducible}
  \nomnorefpage}

\nomenclature[o$T_\pm$]{$\T_\pm$}{braid automorphisms of an integrable $U_q(\lie{sl}(2,\mathbb{K}))$-module,
  p.\pageref{nom:T_pm} \nomnorefpage }

\nomenclature[o$T_i$]{$\T_i$}{braid automorphism of an integrable $U_q(\lie{g})$-module,  p.\pageref{nom:T_i}  \\
  --- of $U_q(\lie{g})$, p. \pageref{deftiautomorphism}
  \nomnorefpage }

\nomenclature{$E_{\beta_r}$}{quantum root vector, p.\pageref{nom:quantum_root_vector_w}, p.\pageref{def:quantum_root_vectors} \nomnorefpage}

\nomenclature{$F_{\beta_r}$}{quantum root vector, p.\pageref{nom:quantum_root_vector_w}, p.\pageref{def:quantum_root_vectors} \nomnorefpage}

\nomenclature[o$R$]{$\R$}{
  universal $R$-matrix, p.\pageref{defquasitriangular} \\
  --- in particular for $\D(G_q)$, p.\pageref{Rmatrixformula}
  \nomnorefpage}

\nomenclature{$r$}{universal $r$-form, p.\pageref{defrform} \\
	--- in particular for $\Poly(G_q)$, p.\pageref{rformconstruction}
  \nomnorefpage}%

\nomenclature[o$Y_q$]{$\mathbf{Y}_q$}{
	  set of elements of $\lie{h}_q^*$ of order $2$, p.\pageref{def:order_2} \\
    --- identified with $\mathbf{Y}_q \cong \half i\hbar^{-1} \roots^\vee / i\hbar^{-1} \roots^\vee$ when the ground field is $\mathbb{C}$, p.\pageref{nom:order_2_C}
    \nomnorefpage}%

\nomenclature[o$e_i$]{$\tilde{e}_i$}{
	 Kashiwara operators, p.\pageref{def:crystal}, p.\pageref{nom:Kashiwara_operators}
	 \nomnorefpage
}%

\nomenclature[o$f_i$]{$\tilde{f}_i$}{
	Kashiwara operators, p.\pageref{def:crystal}, p.\pageref{nom:Kashiwara_operators}
	\nomnorefpage
}%



\nomenclature[o$W$11]{$W$}{
	multiplicative unitary of a locally compact quantum group, \pageref{nom:multiplicative_unitary} \\
	--- of an algebraic quantum group, \pageref{nom:W_alg}
  \nomnorefpage}

\nomenclature[o$Ghat$]{$\hat{G}$}{
	Pontrjagin dual of a locally compact quantum group, \pageref{nom:lc-dual} \\
	--- of an algebraic quantum group, \pageref{nom:alg_dual}
  \nomnorefpage}%

\nomenclature{$L^2(G)$}{
  GNS Hilbert space of a locally compact quantum group, \pageref{nom:lc_L2G},\\
  --- of an algebraic quantum group, \pageref{nom:L2_alg}
  \nomnorefpage}%

\nomenclature{$\Lambda$}{
	GNS map $\N_\phi \to L^2(G)$ of a locally compact quantum group, \pageref{nom:lc_L2G} \\
	--- of an algebraic quantum group, \pageref{nom:alg_GNS}
	\nomnorefpage
}

\nomenclature[o$F$]{$\F$}{
	linear isomorphism $H \to \hat{H}$ from a multiplier Hopf algebra to its dual, $\F(f) = \phi(\;\cdot\; f)$, \pageref{nom:F_Hopf}\\
  --- from an algebraic quantum group to its dual, \pageref{Fouriertransform}
  \nomnorefpage}



\nomenclature{$F_\lambda$}{
  matrices instituting $S^2$ on simple representations a compact quantum group, p.\pageref{nom:F_general}
  \nomnorefpage}


\nomenclature{$\omega^\mu_{ij}$}{
	elements of basis for $\D(G_q)$ dual to matrix coefficients, p.\pageref{nom:dual_basis_Gq}\\
	--- and for $\DF(K)$ where $K$ is a compact quantum group, p.\pageref{nom:dual_basis_K}
	\nomnorefpage
}%

\nomenclature{$F$}{
	element of $\M(\DF(K))$ instituting $S^2$, p.\pageref{nom:F_general} \\
	--- in particular, $F=K_{-2\rho}$ in $\M(\DF(K_q))$, p.\pageref{lem:F_is_K2rho}
  \nomnorefpage
}

\nomenclature{$\iota$}{
	diagonal embedding $\iota:\DF(G_q) \to \M(\DF(K_q) \otimes \DF(K_q))$, \pageref{nom:iota-DGq} \\
	---diagonal embedding $\iota: U_q^\mathbb{R}(\lie{g}) \to U_q(\lie{g})\otimes U_q(\lie{g})$, \pageref{uqrgdiagonalembedding}
  \nomnorefpage
}%

\nomenclature{$P$}{
	Kostant partition function, p.\pageref{nom:Kostant_partition_function1}, p.\pageref{nom:Kostant_partition_function2}
	\nomnorefpage
}%

\nomenclature[s$M^\vee$]{$M^\vee$}{
	restricted dual of a weight module $M$, 
	p.\pageref{nom:restricted_dual1}, p.\pageref{nom:restricted_dual2}
	\nomnorefpage
}%

\nomenclature[o$P(V(\mu))$]{$ \weights(V(\mu)) $}{
	set of weights of $V(\mu)$, 
	p.\pageref{nom:weights_of_V}
	\nomnorefpage
}%

\nomenclature[o$C_\lambda$]{$\mathbb{C}_\lambda$}{
	one-dimensional $U_q(\mathfrak{h})$-module or $U_q(\mathfrak{b})$-module,
	p.\pageref{nom:C_lambda}, p.\pageref{nom:C_lambda2}
	\nomnorefpage
}

\nomenclature[o$ann$]{$\ann$}{
	annihilator ideal of a module,
	p.\pageref{annihilatorlemma}, p.\pageref{nom:ann2}
  \nomnorefpage}

\nomenclature{$\theta$}{
	algebra automorphism, coalgebra anti-automorphism $\theta = \tau \hat{S}$ of $U_q(\lie{g})$,
  p.\pageref{nom:theta1}, p.\pageref{nom:theta2}
  \nomnorefpage
}%

\nomenclature{$l$}{component of the decomposition of
	 $(\mu,\lambda)\in\lie{h}_q^*\times \mathfrak{h}^*_q$ determined either by:
	 --- $\mu = l-r$, $\lambda= -l-r$, p.\pageref{nom:lr} \\
	 --- $\mu=l-r$, $\lambda+2\rho = -l-r$,
	 p.\pageref{nom:lr_shifted}, p.\pageref{nom:lr_shifted2}
   \nomnorefpage
}%

\nomenclature{$r$}{component of the decomposition of
	$(\mu,\lambda)\in\lie{h}_q^*\times \mathfrak{h}^*_q$ determined either by:
	--- $\mu = l-r$, $\lambda= -l-r$, p.\pageref{nom:lr} \\
	--- $\mu=l-r$, $\lambda+2\rho = -l-r$,
	p.\pageref{nom:lr_shifted}, p.\pageref{nom:lr_shifted2}
	\nomnorefpage
}%

\nomenclature{$\iota_\alpha$}{
	embedding $\iota_\alpha: U_{q_\alpha}^\mathbb{R}(\mathfrak{psu}(2)) \rightarrow U_q^\mathbb{R}(\mathfrak{k})$ or $\iota_\alpha: \DF(SU_{q_\alpha}(2)) \rightarrow \M(\DF(K_q))$ associated to a simple root,
	p.\pageref{nom:suq2-embedding1}
  \nomnorefpage
}%

\nomenclature{$\lambda_\alpha$}{
	restriction of $\lambda\in\lie{h}^*$ with respect to a simple root $\alpha$,
	p.\pageref{nom:weight_restriction1} \\
	--- restriction of $\lambda\in\lie{h}_q^*$ with respect to a simple root $\alpha$,
  p.\pageref{nom:weight_restriction2}
  \nomnorefpage
}%

\nomenclature{$\sigma$}{flip map on tensor products,
	p.\pageref{nom:flip1}, p.\pageref{nom:flip2}
  \nomnorefpage
}

\nomenclature[o$\beta$2]{$\beta_M$}{bar involution on a $ U_q(\mathfrak{g}) $-module $ M $, 
  p.\pageref{nom:bar_on_V}, p.\pageref{nom:bar_on_M}
  \nomnorefpage
}%

\nomenclature{$e^\mu$}{
	element of the group ring $\mathbb{C}[\weights]$, p.\pageref{nom:e_lambda} \\
	--- function in $\CF^\infty(T)$ corresponding to a character $\chi_\mu$ of $U_q^\mathbb{R}(\lie{t})$, p.\pageref{nom:e_mu}
	\nomnorefpage%
}%

\printnomenclature

\newpage

\bibliographystyle{plain}

\bibliography{cvoigt}

\end{document}